\documentclass[12pt,titlepage]{report}
\usepackage[francais]{babel}
\usepackage[latin1]{inputenc}
\usepackage{amsmath}
\usepackage{amsfonts}
\usepackage{diagrams}
\usepackage{makeidx}
\usepackage{oldgerm}
\font\Goth=yinitas scaled \magstep0

\newcommand{\Gth}[1]{\lower5mm\hbox{\Goth #1}}

\author{Régis PELLISSIER}
\title{CATEGORIES ENRICHIES FAIBLES}

\newtheorem{defin}{Définition}[section]
\newtheorem{lem}[defin]{Lemme}
\newtheorem{prop}[defin]{Proposition}
\newtheorem{cor}[defin]{Corollaire}
\newtheorem{theo}[defin]{Théorème}
\newtheorem{ex}[defin]{Exemple}
\newtheorem{hyp}[defin]{Hypothèses}
\newcommand{\cmf}{catégorie de modèles fermée}
\newcommand{\cmfcof}{catégorie de modèles fermée engendrée de manière
cofibrante}
\newcommand{\cmfmono}{catégorie de modèles fermée engendrée par
monomorphismes}
\newcommand{\precat}{$\mathcal{C}$-pré\-caté\-gorie}
\newcommand{\precats}{$\mathcal{C}$-précaté\-gories}
\newcommand{\cat}{$\mathcal{C}$-caté\-gorie}
\newcommand{\cats}{$\mathcal{C}$-caté\-gories}
\newcommand{\obc}{objet régal}
\newcommand{\obcs}{objets régaux}
\newcommand{\eqc}{alliance}
\newcommand{\eqcs}{alliances}
\newcommand{\eq}{équi\-valence}
\newcommand{\eqs}{équi\-valences}
\newcommand{\discret}{discrétisante}

\newcommand{\prg}{propriété de relèvement à gauche}
\newcommand{\prd}{propriété de relèvement à droite}
\newcommand{\colim}[1]{\operatornamewithlimits{colim}_{\overrightarrow{\ #1 \ }}}
\newcommand{\colimite}[1]{\operatornamewithlimits{colim}_{\ #1 \ }}
\newcommand{\segal}[2]{#1_{#2}\rightarrow #1_1\times_{#1_0}\ldots\times_{#1_0}#1_1}
\newcommand{\sefb}[3]{#1_{#2}(#3_0,\ldots,#3_{#2})\rightarrow
#1_1(#3_0,#3_1)\times\ldots\times #1_1(#3_{#2-1},#3_{#2})}

\makeindex

\begin{document}
\maketitle

\begin{titlepage}
\begin{center}
{\bf
UNIVERSITÉ DE NICE - SOPHIA ANTIPOLIS

UFR SCIENCES

École Doctorale de Sciences Fondamentales et Appliquées

Département de Mathématiques

Laboratoire J. A. Dieudonné, U.M.R. C.N.R.S. 6621}

\vspace{1.2cm}
{\Large \bf THÈSE}

\vspace{1.2cm}

{\bf présentée pour obtenir le titre de

DOCTEUR ÈS SCIENCES

de l'Université de Nice-Sophia Antipolis

Spécialité : Mathématiques

par}

\vspace{5mm}
{\Large \bf Régis PELLISSIER}

\vspace{1cm}
{\LARGE \linespread{1.6} \bf Catégories enrichies faibles. \par}

\vspace{1cm}
Soutenue publiquement le jeudi 27 juin 2002 devant le jury composé de :
\end{center}

$$\begin{array}{ll}
\mbox{M. Clemens BERGER} & \mbox{Maître de conférence à l'Université de Nice}\\
\mbox{M. André HIRSCHOWITZ} & \mbox{Professeur à l'Université de Nice}\\
\mbox{M. André JOYAL} & \mbox{Professeur à l'Université de Montreal}\\
\mbox{M. Bernhard KELLER} & \mbox{Professeur à l'Université de Jussieu}\\
\mbox{M. Jean-Michel LEMAIRE} & \mbox{Professeur à l'Université de Nice}\\
\mbox{M. Georges MALTSINIOTIS} & \mbox{Directeur de recherche CNRS à Jussieu}\\
\mbox{M. Carlos SIMPSON} & \mbox{Directeur de recherche CNRS à Nice}\\
\mbox{M. Bertrand TOEN} & \mbox{Chargé de recherche CNRS à Nice}\\
\end{array}$$

\vspace{1,5mm}
\begin{center}
à 15h30, en salle de conférences du Laboratoire J.A. Dieudonné.
\end{center}
\end{titlepage}

\newpage

\begin{verse}
{\LARGE \textgoth{Deploration sur le trespas de feu Jean Molinet, faicte par
maistre Guillaume Cretin.}}\\
\end{verse}
\begin{flushleft}
\hbox{\Gth{M}\LARGE\textgoth{olinet n'est sans bruyt, ne sans nom, non,}}
{\LARGE \textgoth{Il a son son, et comme tu vois, voix,\\
Son doulx plaid plaist mieulx que ne faict ton ton,\\
Ton vif art ard plus cher que charbon bon,\\
Tes trenchans chantz perchent ses parois roidz,\\
D'entre gent gent ont nobles Franchois choiz,\\
Se ne doibz doigtz doubter en son laict laid,\\
Car souvent vent vient au molinet nect.}}
\end{flushleft}

\vspace{5mm}
\begin{flushright}(Guillaume Crétin, 1507)\end{flushright}

\newpage

......

\newpage

\begin{flushright}
 \emph{En l'honneur des Pellissier Dragon}
\end{flushright}

\newpage
....

\chapter*{Introduction}

\newpage

Il existe à l'heure 
actuelle plusieurs définitions de $n$-catégories. L'une des pionnières 
est celle de Brown et Loday (approche cubique). Street dans 
{\it The algebra of oriented simplexes} (Journal of Pure
and Applied Algebra 1987) a proposé une définition basée sur la condition
de Kan. Baez et Dolan ont développé dans leur prépublication {\it Higher
dimensional algebra III: n-categories and the algebra of opetopes} (1997) une 
théorie de $n$-catégories basée sur les opérades. Tamsamani a défini 
dans sa thèse {\it Sur des notions de n-catégorie et n-groupoïde non
stricts via des ensembles multisimpliciaux} (K-theory 1999) la notion de 
$n$-nerf qui est une généralisation de celle du nerf d'une petite 
catégorie. Enfin Hirschowitz et Simpson se sont basés pour leur
prépublication {\it Descente pour les n-champs} (1998) sur la notion de
$n$-catégorie de Segal, une version simpliciale de la notion de $n$-nerf.\\

En regardant de plus près les travaux de Tamsamani et d'Hirschowitz-Simpson, on
remarque une forte ressemblance à la fois dans la construction des
$n$-catégories faibles, comme on l'a remarqué ci-dessus, et dans la démonstration de l'existence d'une
structure de \cmf~les concernant, la démonstration pour les $n$-nerfs se
trouvant dans la prépublication de Simpson {\it A closed model structure for
n-categories, internal Hom, n-stacks and generalized Seifert-Van Kampen} et
les idées de celle pour les $n$-catégories de Segal dans {\it Descente pour les n-champs}.
Le principe général est de définir récursivement une $n$-catégorie
faible comme un préfaisceau de la catégorie simpliciale dans la catégorie
des $(n-1)$-précatégories qui niveau par niveau est une $(n-1)$-catégorie
faible et
qui vérifie certaines conditions dites conditions de Segal. De même, on démontre qu'il existe une
structure de \cmf~pour les $n$-catégories faibles en supposant montrée celle
pour les $(n-1)$-catégories faibles. On remarquera que la structure de
\cmf~n'est pas mise directement sur la catégorie des $n$-catégories faibles
mais sur celle des $n$-précatégories faibles, c'est-à-dire des
préfaisceaux de la catégorie simpliciale dans la catégorie des
$(n-1)$-précatégories.\\

Pour résumer ce qui précède, on se donne une \cmf~$\mathcal{C}$ et l'on
regarde les préfaisceaux de la catégorie simpliciale vers $\mathcal{C}$, que
l'on appellera \precats. Puis on définit une notion de catégorie sur
$\mathcal{C}$, ce que l'on appellera \cat, en demandant aux \precats~de
vérifier certaines conditions, principalement celles de Segal. Enfin on
démontre qu'il existe sur la catégorie des \precats~une structure de
\cmf~basée sur la notion de \cat. Le but de cette thèse est donc
d'étudier cette forme de généralisation de la notion de catégorie, ce
qu'on appellera catégorie enrichie faible, tout d'abord en définissant en
détail la notion de \cat, puis en cherchant à montrer l'existence d'un bon
procédé de catégorisation des \precats. Après quoi nous démontrerons
l'existence d'une structure de \cmf~sur les \precats~avec pour \eqs~faibles
les morphismes de \precats~qui après catégori\-sation sont des \eqs~de \cats.
Et pour finir, nous montrerons qu'avec cette structure de \cmf~sur les \precats, on retrouve la structure de \cmf~pour les catégories de Segal
de la prépublication {\it Descente pour les n-champs}.\\

Dans un premier chapitre, nous nous intéressons à la notion de catégorie
enrichie faible sur une catégorie quelconque $\mathcal{C}$. Afin de simplifier
les notations, on nommera \cats~les catégories enrichies faibles sur
$\mathcal{C}$. Pour commencer, étant donnée une catégorie $\mathcal{C}$ quelconque, nous
définissons les \precats~comme les préfaisceaux de la catégorie
simpliciale, notée $\Delta$, dans la catégorie $\mathcal{C}$. Toutefois comme nous ne
voulons considérer que les préfaisceaux $A$ dont l'espace des objets $A_0$
est un ensemble, nous demandons que la catégorie $\mathcal{C}$ soit munie
d'objets discrets, ce qui nous conduit à définir la notion de catégorie \discret.
Ainsi les préfaisceaux de la catégorie simpliciale $\Delta$ dans une
catégorie \discret~$\mathcal{C}$ ayant
leurs ensembles d'objets discrets seront les seuls à pouvoir être appelés \precats.\\

La définition de \cat~est beaucoup plus difficile à définir. Pour suivre les exemples
des $n$-catégories faibles de Tamsamani ou de Segal évoquées précédemment, nous aimerions dire
qu'une \cat~est une \precat~$A$ telle que, pour tout $n$ entier, $A_n$ est une
"catégorie de $\mathcal{C}$" et que, pour tout $m>1$, les morphismes dits de
Segal $A_m\rightarrow A_1\times_{A_0}\ldots\times_{A_0} A_1$ sont des "\eqs~de
catégories dans $\mathcal{C}$". De même, nous voudrions dire qu'un morphisme
$f$ de \cats~est une \eq~de \cats~s'il est essentiellement surjectif et
pleinement fidèle. Pour la pleine fidélité, il s'agit simplement de
demander que pour tout couple d'objets $(x,y)$ de la source de $f$, le morphisme
de $\mathcal{C}$ $f_1(x,y)$ soit une "\eq~de catégories dans $\mathcal{C}$".
Pour définir l'essentielle surjectivité, toujours en suivant les exemples des
$n$-catégories faibles de Tamsamani et de Segal, cela nécessite l'existence
d'un foncteur noté $\tau_0$ des \cats~vers les ensembles généralisant la
notion d'\eq~d'objets dans une catégorie. Or ce foncteur se définit très
facilement si pour les "catégories de $\mathcal{C}$" un tel foncteur existe
déjà.\\

Ainsi l'on remarque que les définitions de \cat~et d'\eq~de \cats~nécessitent
trois notions sur $\mathcal{C}$ qui ont vocation à être respectivement
celles de catégorie, d'\eq~de catégorie et d'\eq~d'objets d'une catégorie.
Nous les nommerons respectivement \obc, \eqc~d'\obcs~et $\tau_0$. La catégorie
$\mathcal{C}$ munie de ces trois notions et vérifiant un certain nombre de
propriétés techniques constituera une donnée de Segal. C'est grâce à
cette notion de donnée de Segal que les définitions de \cat~et d'\eq~de
\cats~telles qu'elles ont été données ci-dessus auront un sens.\\

Après avoir défini correctement les notions de \cat~et d'\eq\\ de \cats,
nous allons chercher à construire un procédé de catégorisation qui
associera à toute \precat~une \cat~tout en préservant son type
d'homotopie. En nous inspirant de la prépublication de Simpson {\it A
  closed model structure for n-categories, internal Hom, n-stacks and
  generalized Seifert-Van Kampen} \cite{s}, nous allons chercher un
tel procédé parmi les procédés d'addition transfinie de flèches. C'est
pourquoi nous allons entièrement consacrer notre deuxième chapitre à
l'étude des procédés d'addition de flèches que nous appellerons plan
d'addition de cellules. Un des principaux rôles de ces plans
d'addition de cellules est de résoudre le problème suivant. On se
donne une catégorie $\mathcal{C}$ cocomplète et une famille $I$ de
morphismes de $\mathcal{C}$. On appelle objet $I$-injectif tout objet
ayant la \prd~par rapport aux flèches de $I$. Le problème qui se pose
alors est de trouver un procédé qui transforme tout objet en objet
$I$-injectif. La solution de ce problème se trouve en faisant des
colimites transfinies de sommes amalgamées des flèches de $I$, ce qui
n'est autre qu'appliquer un plan d'addition de cellules pour $I$. Et
c'est un argument du type petit objet qui assurera que tel ou tel
plan d'addition de cellules de $I$ rend effectivement $I$-injectif.\\

Cependant cette astuce de passer par les plans d'addition de cellules
pour résoudre le problème de la catégorisation des \precats~ne prend
pas en compte la préservation de l'homotopie. On doit alors résoudre
un autre problème. Parmi les plans d'addition de cellules $P$
fonctoriels rendant $I$-injectif, y en a-t-il un tel que, pour tout
objet $A$, l'image par $P$ du morphisme naturel de $A$ dans le
résultat du plan $P(A)$ est une \eq~d'objets $I$-injectifs, en
supposant qu'une telle notion existe. Toujours en s'inspirant de
\cite{s}, nous allons définir un plan d'addition de cellules rendant
$I$-injectif ayant la propriété de factoriser de manière unique tout
morphisme à valeur dans un objet $I$-injectif dont les relèvements par
rapport aux morphismes de $I$ sont marqués. Ce plan noté $E_I$ sera le
prototype de la catégorisation $Cat$. Son avantage est que s'il existe
un plan d'addition de cellules rendant $I$-injectif et préservant le
type d'homotopie pour les objets $I$-injectifs, cela entraînera que
notre procédé $E_I$ a bien la propriété voulue de
préservation de l'homotopie pour tous les objets.\\

Comme l'un des atouts de la catégorisation est de pouvoir définir des
colimites pour les \cats, dont les colimites en tant que \precats~ne sont
généralement pas des \cats, il est utile de pouvoir comparer la
catégorisation d'une colimite de \precats~avec la colimite des
catégori\-sées. Pour cela les plans d'addition de cellules nous seront à nouveau d'un
grand secours, car ce problème consiste, en termes de plans d'addition de
cellules, à comparer deux plans d'addition de cellules rajoutant les mêmes
flèches mais pas dans le même ordre. En effet dans les plans d'addition de cellules, on ne fait souvent
pas attention à l'ordre dans lequel on additionne les cellules. Or il existe
certains plans d'addition de cellules, dit rationnels, pour lequel l'ordre
d'addition est bien maîtrisé. On montrera donc que tout plan d'addition
de cellules est rationalisable, ce qui permettra de montrer que deux plans
d'addition de cellules sont équivalents s'ils ont même rationalisation.\\

Cependant pour pouvoir appliquer ces plans d'addition de cellules au problème
de la catégorisation, encore faut-il qu'il existe un ensemble de morphismes
$I$ pour lequel les \precats~$I$-injectives sont des \cats, on dira alors que
$I$ engendre la notion de \cat. Comme il n'existe en
général pas de tel ensemble caractérisant les \cats, nous serons obligés
de passer par une notion plus forte de \cat, que l'on nommera \cat~facile, pour
laquelle il existe un ensemble $I$ l'engendrant. C'est donc l'objet du
troisième chapitre que de bien définir cette notion de \cat~facile et de
trouver son ensemble générateur.\\

Afin de définir cette notion de
\cat~facile, nous définissons tout d'abord les données de Segal
proto-faciles comme des données de Segal pour lesquelles les notions d'\obc~et
d'\eqc~d'\obcs~sont elles-mêmes faciles, c'est-à-dire caractérisées par des propriétés de relèvement à
droite par rapport à des familles de flèches. Les \cats~faciles sont alors
simplement les \cats~correspondant aux données de Segal proto-faciles, qui
forment bien évidemment des données de Segal. Afin de montrer qu'il existe
un ensemble de flèches qui caractérise par relèvement les \cats~faciles
ainsi définies, nous
définissons une construction $\Theta$ qui est un bifoncteur qui à un
ensemble simplicial $X$ et à un objet $C$ de $\mathcal{C}$ associe une
\precat~$X\Theta C$. Cette construction a la propriété fondamentale
suivante : si $X$ est un ensemble simplicial représentable, noté $\Delta[n]$
avec $n$ entier, alors les morphismes de \precats~de $\Delta[n]\Theta C$ vers
une \precat~$A$ correspondent aux morphismes dans $\mathcal{C}$ de
$C$ vers $A_n$. Grâce à cette construction, nous pouvons alors exhiber les
flèches génératrices de \cats~faciles et définir le procédé de
catégorisation $Cat$ comme la plan d'addition de cellules $E_I$ pour l'ensemble
des flèches génératrices de \cats~faciles. En rassemblant toutes les
hypothèses sur la donnée de Segal proto-facile qui permettent d'obtenir
cette construction $Cat$, nous obtenons la notion de donnée de Segal
pré-facile qui non seulement donne lieu à une notion de \cat~mais aussi
assure l'existence du foncteur $Cat$ transformant les \precats~en \cats.\\

Toutefois il reste encore à montrer que la construction $Cat$ préserve le
type d'homotopie des \precats. Ce problème est assez délicat mais, comme on
l'a vu dans le deuxième chapitre, si l'on est capable de construire une
catégorisation préservant le type d'homotopie pour les \cats~faciles, alors la
construction $Cat$ préservera bien le type d'homotopie des \precats. Aussi
allons-nous consacrer le quatrième chapitre de cette thèse à la
construction d'un tel procédé. Mais tout d'abord nous montrerons que les \eqs~de
\cats~vérifient la propriété du "trois pour deux", c'est-à-dire que si
parmi $f,g,g\circ f$ deux morphismes sont des \eqs~de \cats~alors le troisième
aussi. En effet ce résultat est indispensable pour pouvoir appliquer le
résultat du deuxième chapitre. Puis nous allons chercher une
catégorisation préservant l'homotopie pour les \cats~faciles. Ceci va passer par la mise en place de
constructions $Reg$ et $Seg$ permettant de comprendre au niveau de la catégorie $\mathcal{C}$
ce qu'il se passe lorsque l'on fait une somme amalgamée par l'une des
flèches génératrices. Ces constructions vont alors montrer que certaines
de ces sommes amalgamées sont compatibles entre elles, ce qui induira un
plan d'addition de cellules $Bigcat$ pour les flèches génératrices qui
niveau par niveau respecte l'homotopie lorsqu'on l'applique aux \cats. En
appliquant alors le résultat du deuxième chapitre, l'on obtient que sous
certaines hypothèses, qui donneront lieu à la notion de donnée de Segal
facile, la construction $Cat$ préserve bien l'homotopie des \precats~et donc
résout bien comme il faut le problème de la catégorisation.\\

A la fin du quatrième chapitre, nous avons donc obtenu le résultat suivant :
pour une donnée de Segal facile, il existe une notion de \cat~ainsi qu'un
bon procédé de catégorisation $Cat$. Cependant parmi les hypothèses
d'une donnée de Segal facile, on a demandé à la catégorie $\mathcal{C}$
d'être une \cmf, car cette notion permet de bien gérer l'homotopie. Aussi
est-il naturel de demander à la catégorie des \precats~d'avoir une structure
de \cmf~traduisant l'homotopie des \cats. Comme le foncteur de catégorisation
$Cat$ préserve l'homotopie des \precats, on définira les \eqs~faibles de
\precats~comme les morphismes de \precats~dont l'image par $Cat$ est une \eq~de
\cats. Afin de faciliter la démonstration du fait que les \precats~avec leurs
\eqs~faibles forment une \cmf, nous prendrons pour cofibrations les
monomorphismes et pour fibrations les morphismes de \precats~ayant la \prd~par
rapport aux monomorphismes qui sont aussi des \eqs~faibles de \precats.\\

Dans le cinquième chapitre, nous allons donc montrer qu'il existe sur les
\precats~une structure de \cmf~avec les notions de cofibrations, d'\eqs~faibles et de
fibrations définies ci-dessus. Même si, comme dans la prépublication
d'Hirschowitz-Simpson {\it Descente pour les $n$-champs}, nous allons montrer en
fait qu'il existe une structure de \cmfcof, la démonstration suit de près, en
la généralisant, celle de la prépublication de Simpson {\it A closed model structure for
n-categories, internal Hom, n-stacks and generalized Seifert-Van Kampen}. Nous
commençons ce chapitre par des rappels de quelques résultats sur la structure de \cmf, en particulier en rappelant le lemme de
reconnaissance d'une \cmfcof~suivant :\\
\\
{\bf Lemme de reconnaissance d'une \cmfcof~:}\\
\\
{\it Soit $\mathcal{M}$ une catégorie munie d'une classe de morphismes W et de deux
ensembles de morphismes I et J. Supposons que $\mathcal{M}$ vérifie les
propriétés suivantes :\\
\\
0) $\mathcal{M}$ est complète et cocomplète,\\
\\
1) W est stable par rétracts et satisfait à l'axiome "trois pour
deux" : pour tout couple $(f,g)$ de morphismes
composables, si deux morphismes parmi $f,\; g$ et $g\circ f$ sont dans W alors le morphisme
restant est aussi dans W,\\
\\
2) I et J permettent l'argument du petit objet,\\
\\
3) Toute J-cofibration est à la fois une I-cofibration et dans W,\\
\\
4) Tout morphisme I-injectif est à la fois J-injectif et dans W,\\
\\
5) L'une des deux propriétés suivantes est vraie : \\
 - les morphismes à la fois I-cofibrations et dans W sont
des J-cofibrations,\\
 - les morphismes à la fois J-injectifs et dans W sont
I-injectifs.\\
\\
Alors il existe une structure de \cmfcof~sur $\mathcal{M}$ avec W pour classe
des équivalences faibles, I pour l'ensemble des cofibrations génératrices
et J pour l'ensemble des cofibrations triviales génératrices.}\\

Nous allons tout au long du cinquième chapitre montrer que les \precats~avec
les monomorphismes et les \eqs~faibles de \precats\\ vérifient successivement
les six hypothèses de ce lemme de reconnaissance de la structure de \cmfcof.
L'axiome 0) ayant déjà été  montré dans le premier chapitre, nous
commençons en fait par montrer l'axiome 1) qui concerne de la stabilité des
\eqs~faibles de \precats. Or cet axiome se réduit à montrer simplement la stabilité par
rétracts des \eqs~faibles de \precats~car l'axiome "trois pour deux" découle directement du fait que les
\eqs~de \cats~le vérifient ainsi que du fait que la construction $Cat$
préserve l'homotopie, résultats ayant été montrés dans le chapitre précédent. Puis en
supposant que la catégorie $\mathcal{C}$ elle-même est une \cmfcof, nous
définissons l'ensemble $\mathcal{I}$ des cofibrations génératrices des \precats~à
partir des cofibrations génératrices de $\mathcal{C}$ en utilisant la
construction $\Theta$ et nous définissons les cofibrations triviales
génératrices comme les cofibrations triviales à sources et buts petits
pour un cardinal transfini suffisamment grand. De ces définitions découle
directement la vérification de l'axiome 2).\\ 

Etant donné que l'axiome 3) est le plus difficile à vérifier et qu'en
outre pour vérifier l'axiome 5) l'on se sert des résultats de l'axiome 3),
nous allons donc faire la vérification de l'axiome 4) avant celle du 3).
Pour démontrer l'axiome 4), nous
commençons par montrer que l'ensemble $\mathcal{I}$ engendre les
cofibrations par des arguments et des techniques classiques dans la théorie
des catégories de modèles fermées. De ce résultat découle le fait que les morphismes $I$-injectifs se
relèvent par rapport aux cofibrations et donc en particulier par rapport aux
cofibrations triviales, ce qui montre la première partie de l'hypothèse 4). 
Il ne reste donc plus qu'à montrer, pour terminer la vérification de l'axiome 4) que les morphismes $I$-injectifs
sont des \eqs~faibles, ce qui nécessite au préalable quelques rappels sur la
préservation des \eqs~faibles par colimites dans une \cmf.\\

La vérification de l'axiome 3) est bien comme prévu la plus compliquée de
toutes et comme elle est nécessaire pour pouvoir vérifier l'axiome 5), nous
sommes obligés de vérifier d'abord l'axiome 3). Comme les propriétés des cofibrations entraînent
que les morphismes engendrés par les cofibrations triviales génératrices
sont bien des cofibrations, l'axiome 3) revient alors à montrer que ce sont des \eqs~faibles.
Or nous avons déjà vu que ces cofibrations génératrices vérifient
l'argument du petit objet (axiome 2)). Ainsi tout morphisme engendré n'est
autre qu'un rétract d'une colimite séquentielle transfinie de sommes amalgamées de cofibrations
génératrices, qui sont par définition des cofibrations triviales. Comme la
stabilité par rétract découle de l'axiome 1) et des propriétés des
cofibrations, si l'on montre que les cofibrations triviales sont
stables par somme amalgamée et par colimite séquentielle transfinie, on
aura le résultat attendu.\\

Pour démontrer ce résultat, nous allons commencer par le cas particulier des
cofibrations triviales qui sont des isomorphismes sur les objets, que l'on
appellera iso-cofibrations triviales. Grâce aux propriétés de la
construction $Cat$ vis-à-vis des colimites découlant des rationalisations
de plans d'addition de cellules du deuxième chapitre, montrer la stabilité
des iso-cofibrations triviales par colimites revient à montrer celle de leurs
catégorisées. On se contentera donc de regarder les iso-cofibrations
triviales entre \cats. Or justement pour les \cats, les
iso-cofibrations triviales ne sont autres que les morphismes de \cats~qui niveau
par niveau sont des cofibrations triviales. Ainsi la stabilité des
iso-cofibrations triviales de \cats~revient à une stabilité niveau par niveau, ce qui
est assuré par la \cmf~$\mathcal{C}$. Toutefois comme on l'a vu plus haut, les
colimites de \cats~ne sont pas des \cats, aussi allons-nous nous retrouver à
avoir montré que certains morphismes de \precats~sont niveau par niveau des
cofibrations triviales. C'est pourquoi afin de montrer la stabilité des
iso-cofibrations triviales, nous allons dans un premier temps caractériser les
iso-\eqs~de \cats~et dans un second temps montrer qu'un morphisme de \precats~qui niveau par niveau
est une \eq~faible de $\mathcal{C}$ est en fait une \eq~faible de \precats.\\

Afin d'obtenir la stabilité des cofibrations triviales quelconques, nous
devons montrer un théorème sur le produit des \cats~qui après réduction
revient à montrer qu'un morphisme entre \precats~ayant certaines
propriétés est une \eq~faible. Pour montrer ce dernier résultat, nous
allons faire intervenir la notion de \precat~ordonnée. Nous définissons les
\precats~ordonnées comme des \precats~dont l'ensemble des objets est totalement ordonné et
vérifiant certaines propriétés. Un des
avantages de cette notion est qu'elle possède une caractérisation très
simple d'\eq~faible. En effet si un morphisme $f$ entre \precats~ordonnées respectant
l'ordre vérifie la propriété que, pour tout couple d'objets adjacents $(x,y)$ de la source de
$f$, $f_1(x,y)$ est une \eq~faible dans $\mathcal{C}$, alors $f$ est une
\eq~faible de \precats. Ce résultat montre en particulier que les flèches
génératrices de \cats~faciles sont bien des \eqs~faibles. Mais ce n'est pas
sa seule utilité car il est aussi le point clef du théorème sur le produit
des catégories. Ce théorème sur le produit des catégories est comme on
l'a dit un des éléments importants de la démonstration de la stabilité des cofibrations
triviales par somme amalgamée.\\

Avec ce dernier résultat, nous sommes en mesure de montrer le théorème qui
nous donne une \eq~entre un produit de \precats~et
le produit de leurs catégorisations. La stabilité des iso-cofibrations
triviales par sommes amalgamées et colimites transfinies, que l'on a montrée
auparavant, permet en effet de réduire ce théorème au cas des \precats~ordonnées dont la
propriété sur les \eqs~faibles démontre ce théorème. De ce théorème
découle la préservation des \eqs~faibles par produit, résultat dont nous
avons besoin pour montrer la préservation par somme amalgamée des
cofibrations triviales quelconques.\\

Pour démontrer la préservation des cofibrations triviales par somme
amalgamée, nous faisons une récurrence transfinie sur la différence entre
le nombre d'objets du but et celui de la source de la cofibration triviale. Le
cas zéro ayant été déjà démontré, nous nous sommes ramenés à montrer ce qu'il se
passe lorsqu'on rajoute un objet, d'abord isomorphe à un objet existant, puis
seulement équivalent à un objet existant, et ce qu'il se passe lorsque l'on prend une colimite
séquentielle transfinie de cofibrations triviales, ce qui montre en même
temps la stabilité des cofibrations triviales par colimite séquentielle
transfinie. On remarque que pour représenter l'\eq~entre deux objets d'une
\cat, on a dû supposer l'existence d'un bon intervalle, noté $\bar{J}$, qui
est une condition difficile à formuler directement sur $\mathcal{C}$. Ayant donc montré que les cofibrations triviales sont stables par
somme amalgamée et colimite séquentielle transfinie, nous avons terminé la
vérification de l'axiome 3), comme on l'a vu précédemment.\\

Il ne nous reste donc plus qu'à vérifier le dernier axiome, ce qui nous
demande quelques lemmes techniques sur la catégorisation $Cat$ ainsi que des
techniques de cardinaux transfinis. Comme à chaque vérification d'axiome
nous avons ressorti les hypothèses qui nous permettaient de les vérifier,
nous aboutissons à l'énoncé d'une ébauche de notre théorème
central qui, sous vingt-deux hypothèses sur la donnée de Segal proto-facile
portée par $\mathcal{C}$, nous assure l'existence d'une structure de \cmfcof~sur les \precats~avec pour
cofibrations les monomorphismes et pour \eqs~faibles les morphismes qui
deviennent des \eqs~de \cats~par catégorisation, les fibrations étant
définies comme il se doit par la propriété de relèvement à droite par rapport aux
cofibrations triviales. Comme cela fait beaucoup d'hypothèses, dont certaines
modifient l'énoncé de certaines autres quand elles sont prises globalement, nous
allons remanier ces hypothèses en les comparant aux diverses notions de
données de Segal existantes, ce qui donnera de nouvelles définitions qui
permettront d'énoncer de manière plus synthétique le théorème central
et conclura du même coup le cinquième chapitre.\\

Nous terminons cette thèse par un sixième chapitre consacré à
l'application du théorème à la donnée de Segal définissant la notion
de catégorie de Segal afin de retrouver la 
\cmfcof~de la prépublication {\it Descente pour les $n$-champs}
d'Hirschowitz-Simpson. La donnée de Segal qui nous donnera les catégories de
Segal est constituée par la catégorie des ensembles simpliciaux avec pour
\obcs\\ les ensembles simpliciaux quelconques, pour \eqcs~d'\obcs~les \eqs~faibles
d'ensembles simpliciaux, c'est-à-dire les morphismes dont la réalisation
géomé\-trique induit des isomorphismes sur les groupes d'homotopie, et pour
foncteur $\tau_0$ la composée du foncteur réalisation géométrique avec
la foncteur composante connexe.\\

Si les vingt premières hypothèses de l'ébauche du théorème central sont facilement
vérifiées, les deux dernières, consistant en l'existence d'un intervalle
$\bar{J}$ caractérisant l'\eq~d'objets dans une catégorie de Segal muni de certaines
propriétés, nécessitent l'introduction de la notion de groupoïde de
Segal. Les groupoïdes de Segal sont des catégories de Segal dont tous les
morphismes sont inversibles à homotopie près. On peut en outre définir des
groupes d'homotopie pour ces groupoïdes de Segal qui caractérisent les
\eqs~entre groupoïdes de Segal. Toutefois pour notre intervalle, il faut
construire une paire de foncteurs adjoints entre les groupoïdes de Segal et
les espaces topologiques qui préservent l'homotopie. Ce couple est constitué
d'une réalisation géométrique des précatégories de Segal et d'un
foncteur singulier, comme c'est le cas pour la paire d'\eqs~de Quillen entre les
ensembles simpliciaux et les espaces topologiques. Pour montrer la
préservation de l'homotopie par cette paire, nous définissons les notions de
catégorie topologique et de groupoïde topologique de manière très
similaire à celles de catégorie de Segal et groupoïde de Segal. Nous
obtenons finalement le fait qu'un morphisme de groupoïdes de Segal est une
\eq~de catégories de Segal si et seulement si sa réalisation
géométrique est une \eq~d'homotopie faible. Ce résultat nous permet de
montrer que notre intervalle $\bar{J}$ existe et vérifie les deux dernières
hypothèses de l'ébauche du théorème central. Ainsi notre théorème central appliqué à la
donnée de Segal sur les ensembles simpliciaux ci-dessus nous redonne bien la
structure de \cmfcof~qu'Hirschowitz et Simpson ont construite pour les
catégories de Segal dans leur prépublication {\it Descente pour les
$n$-champs}.\\

Avant de commencer la thèse proprement dite, je voudrais souligner le fait
que ce travail nécessite encore quelques améliorations. La principale
consiste à exprimer les deux dernières hypothèses du théorème central
concernant l'intervalle $\bar{J}$ caractérisant l'\eq~d'objets des \cats~en terme
d'hypothèses directes sur la catégorie de base $\mathcal{C}$. En effet dès
que l'on arrivera à faire cela, toutes les hypothèses du théorème
central seront directement vérifiables sur la catégorie de base
$\mathcal{C}$. Ceci non seulement permettra d'appliquer plus facilement le
théorème central à d'autres exemples que celui traité dans cette
thèse, mais aussi à amorcer une autre amélioration qui consiste à
montrer que la \cmfcof~des \precats~vérifie elle-même les hypothèses du
théorème central. Dans ce cas, en effet, nous aurions démontré le pas de
la récurrence permettant de définir les $n$-$\mathcal{C}$-catégories et de
leur construire une \cmfcof. Je termine en remerciant mon directeur de thèse Jean-Michel Lemaire pour
sa gentillesse et son appui, André Hirschowitz et Carlos Simpson pour leur
aide permanente et leur encadrement, André Joyal et Georges Maltsiniotis pour
avoir accepté de rapporter cette thèse, Bernhard Keller et Bertrand Toen pour leur participation à mon jury, Nadjia Hohweiller pour ses précieux conseils concernant le nouveau chapitre de la thèse et pour finir Clemens Berger et mes
collègues, de bureau, de l'étage et de l'équipe, pour leurs encouragements et leur soutient sans faille.\\

\newpage
....

\tableofcontents

\newpage
......

\chapter{Catégories faibles enrichies sur une catégorie}

\newpage

Dans ce chapitre, nous allons définir une notion de catégorie enrichie
faible inspirée par les définitions de $n$-nerfs de \cite{t} et de
$n$-catégories de Segal de \cite{h-s}. Cette définition repose
essentiellement sur le fait que toute petite catégorie correspond
bi-univoquement à un ensemble simplicial dit nerf de catégorie, ayant la
propriété que certains morphismes, dits morphismes de Segal, sont des
isomorphismes. L'idée est donc de regarder non plus seulement des préfaisceaux de
la catégorie simpliciale $\Delta$ vers la catégorie des ensembles mais des
préfaisceaux de la catégorie simpliciale $\Delta$ vers une catégorie quelconque
$\mathcal{C}$. Comme on désire ne travailler qu'avec des catégories
enrichies dont l'espace des objets forme un ensemble, on demandera à ces
préfaisceaux de $\Delta$ vers $\mathcal{C}$ que leur espace d'objets soit un
ensemble et on appellera \precats~de tels préfaisceaux.\\ 

 Puis en supposant que
la catégorie $\mathcal{C}$ est munie d'une notion d'\eq, on demandera que les
morphismes de Segal des \precats~soient des \eqs~de $\mathcal{C}$ pour que ces
\precats~soient des catégories faibles enrichies sur $\mathcal{C}$. Toutefois,
dans les exemples de $n$-catégories faibles de \cite{t} et \cite{h-s}, qui
sont en fait des catégories faibles enrichies sur les $n-1$-précatégories,
on demande aux $n$-catégories faibles non seulement que leurs morphismes de
Segal soient des \eqs~mais aussi que niveau par niveau ce soient des
$n-1$-catégories. Aussi définira-t-on sur $\mathcal{C}$ une notion de
donnée de Segal qui consistera essentiellement à munir $\mathcal{C}$ de
notions de catégories de $\mathcal{C}$ et d'\eqs~entre ces catégories de
$\mathcal{C}$, ce que nous nommerons respectivement \obcs~et \eqcs~d'\obcs. Munie
d'une donnée de Segal sur $\mathcal{C}$, nous pourrons alors définir les
\cats~qui seront donc des \precats~qui niveau par niveau sont des \obcs~et dont
les morphismes de Segal sont des \eqcs~d'\obcs.\\

 Nous débuterons donc ce chapitre en regardant à quelles conditions sur la
catégorie $\mathcal{C}$ la notion de \precat~est bien définie. Ces conditions sont
essentiellement dues au fait de demander aux \precats~d'avoir un espace d'objets
qui est un ensemble, ce qui conduira aux notions d'objet discret et de
catégorie \discret.

\newpage

\section{Objets discrets}

Par objet discret, nous voulons modéliser le fait qu'un objet d'une
catégorie $\mathcal{C}$ ressemble à un ensemble. Comme, dans la catégorie
des ensembles, tout ensemble est isomorphe à un coproduit sur son cardinal de
l'objet final, nous proposons donc la définition suivante.

\begin{defin}\index{objets discrets}
Soit $\mathcal{C}$ une catégorie ayant un objet final. Un objet $X$ de $\mathcal{C}$ est un objet
discret si : $$ X=\coprod_{*\rightarrow X} * $$
\end{defin}

Cependant on aimerait beaucoup que ces objets discrets se comportent comme les
ensembles. On va donc définir un foncteur ensemble sous-jacent de
$\mathcal{C}$ vers la catégorie des ensembles et un foncteur objet discret de
la catégorie des ensembles vers $\mathcal{C}$ qui auront vocation à être
une \eq~de catégories entre la catégorie des ensembles et la
sous-catégorie des objets discrets de $\mathcal{C}$.

\begin{defin}
Soit $\mathcal{C}$ une catégorie ayant un objet final et les coproduits de
l'objet final. On définit le couple suivant de foncteurs entre la catégorie $\mathcal{C}$ et la catégorie des ensembles :
\item on note $Enssj$ le foncteur $\operatorname{Hom}_{\mathcal{C}}(*,-)$ qui
à un objet discret associe l'ensemble des morphismes dans $\mathcal{C}$ de
l'objet final, noté $*$, vers l'objet discret.
\item on note $Discret$ le foncteur qui à un ensemble associe le coproduit de l'objet final
indexé par l'ensemble
\end{defin}

A partir de maintenant on identifiera un objet discret $X$ avec l'ensemble des
morphismes de l'objet final vers $X$.

\begin{lem}
Soit $\mathcal{C}$ une catégorie ayant un objet final et les coproduits de
l'objet final. On a toujours que $Discret\circ Enssj$
est l'identité des objets discrets.
\end{lem}
{\it Preuve :} par définition d'objet discret !\\

En revanche la composée $Enssj\circ
Discret$ n'est pas toujours naturellement isomorphe à l'identité des
ensembles. En effet dans une catégorie pointée, les coproduits de l'objet
final sont tous des objets finaux et donc $Enssj\circ
Discret$ est naturellement isomorphe au foncteur constant à valeur un
singleton. Or nous voudrions que $\mathcal{C}$ possède une sous-catégorie
équivalente à celle des ensembles. En outre, pour manipuler plus facilement
les futures \precats, nous aimerions également que tout morphisme dont le but
s'envoie dans un objet discret soit caractérisé par ses fibres en les points
de l'objet discret.\\
\\

{\bf Notation :} \\
Soit $Z\rightarrow X$ un morphisme vers un objet discret. Soit $x$ un élément de
$X$, on note $Z(x)$ le produit fibré de $Z$ par le morphisme $x:*\rightarrow
X$ et on l'appelle fibre de $Z$ en $x$.
Soit un diagramme :
\begin{diagram}
Z & & \rTo^f & & Z'\\
  & \rdTo & & \ldTo & \\
  & & X & & \\
\end{diagram}
On note $f(x)$ le produit fibré, dans la catégorie des morphismes, de $f$
par le morphisme $x:*\rightarrow
X$ et on l'appelle fibre de $f$ en $x$.

\begin{defin}\label{mod}\index{catégorie!\discret}
Une catégorie $\mathcal{C}$ est dite \discret~si elle admet un objet final, les coproduits
de l'objet final, les fibres et possède les propriétés
suivantes :
\item 1) pour tout morphisme $Z\rightarrow X$ vers un objet discret, on a
$$Z=\coprod_{x\in X} Z(x).$$
\item 2) $Enssj\circ Discret$ est naturellement isomorphe à l'identité de la
catégorie des ensembles.
\item 3) les objets discrets sont stables par limites.
\end{defin}

Le dernière hypothèse est assez naturelle car, l'hypothèse 2) assurant
qu'une catégorie \discret~contient une sous-catégorie pleine équivalente
à la caté\-gorie des ensembles, on peut espérer que les limites et colimites
de la catégorie des ensembles se retrouvent dans la sous-catégorie des
ensembles discrets comme limites et colimites dans la catégorie \discret. Comme en outre, on peut remarquer que le foncteur $Discret$ est un adjoint à gauche du foncteur $Enssj$, les colimites d'objets discrets sont automatiquement des objets discrets, c'est pourquoi l'hypothèse 3) ne requiert que la stabilité des objets discrets pour les limites.\\
\\

La proposition suivante montre que la définition de catégorie \discret\\
nous donne bien sur $\mathcal{C}$ les propriétés désirées, à savoir que sa
sous-catégorie d'objets discrets est équivalente à la catégorie des
ensembles et que les morphismes au-dessus d'un objet discret sont
caractérisés par leurs fibres.

\begin{prop}
Soit $\mathcal{C}$ une catégorie \discret. 
\item - Sa sous-catégorie des objets
discrets est équivalente à la catégorie des ensembles, à travers les
foncteurs $Enssj$ et $Discret$. 
\item - Pour tout diagramme :
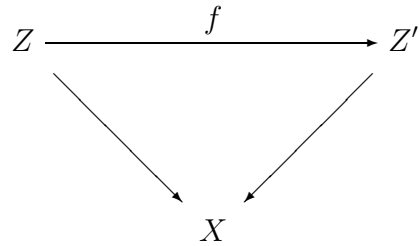
\begin{diagram}
Z & & \rTo^f & & Z'\\
  & \rdTo & & \ldTo & \\
  & & X & & \\
\end{diagram}
on a :
$$f=\coprod_{x\in X} f(x).$$
\end{prop}
{\it Preuve :} la première partie de la proposition découle de la
propriété 2 de la définition de catégorie \discret~et du lemme ci-dessus, la seconde partie est une
vérification triviale qui découle de la propriété 1.

\begin{ex}
La catégorie des ensembles, qu'on notera $\mathcal{ENS}$, est trivialement
\discret.\index{$\mathcal{ENS}$} 
\end{ex}

\begin{ex}
$\mathcal{ENSSIMP}$, la catégorie des ensembles simpliciaux, est\\ 
\discret.\label{enssimpdiscret}\index{$\mathcal{ENSSIMP}$}
\end{ex}

\newpage

\section{\precats}

Munis de la définition de catégorie \discret, nous pouvons maintenant
définir la notion de \precat.

\begin{defin}\index{\precat}
Si $\mathcal{C}$ est une catégorie \discret, on appelle \precat~tout
préfaisceau $A$ de la catégorie simpliciale $\Delta$ vers
$\mathcal{C}$ dont l'espace $A_0$, dit l'ensemble des objets de $A$, est un objet discret de
$\mathcal{C}$. On appelle $\mathcal{C-PC}$\index{$\mathcal{C-PC}$} la sous-catégorie pleine de la
catégorie des préfaisceaux de $\Delta$ vers $\mathcal{C}$ d'objets les
\precats.
\end{defin}

\begin{ex}
Si on prend pour catégorie \discret~$\mathcal{ENS}$, on obtient que la
catégorie des \precats~n'est autre que la catégorie des ensembles
simpliciaux.
\end{ex}

En vue de réitérer la construction des \precats, il est intéressant de voir
si la catégorie des \precats~est elle-même \discret.

\begin{prop}\label{precatdiscret}
Si $\mathcal{C}$ est une catégorie \discret, la catégorie $\mathcal{C-PC}$ des \precats~est 
\discret.
\end{prop}
{\it Preuve :}\\
Par hypothèse, $\mathcal{C}$ est possède l'objet final, les coproduits de
l'objet final et les fibres. Or les catégories de préfaisceaux à valeurs dans une
catégorie possédant l'objet final, les coproduits de
l'objet final et les fibres possèdent également l'objet final, les coproduits de
l'objet final et les fibres, ces derniers étant calculés niveau par niveau.
C'est en particulier le cas pour la catégorie des préfaisceaux de
$\Delta$ vers $\mathcal{C}$. En outre les fibres d'objets discrets dans $\mathcal{C}$ sont encore des objets
discrets et, par définition, l'objet final et ses coproduits sont des objets
discrets. Donc l'objet final, les coproduits de
l'objet final et les fibres niveau par niveau dans les préfaisceaux de $\Delta$ vers
$\mathcal{C}$ préservent la condition que le niveau 0 est discret. D'où
$\mathcal{C-PC}$ possède l'objet final, les coproduits de
l'objet final et les fibres, ces derniers étant calculés niveau par niveau.\\

Soit $Z\rightarrow X$ un morphisme de $\mathcal{C-PC}$ dont le but est un objet
discret. Comme $\mathcal{C-PC}$ est une catégorie de préfaisceaux,
$\coprod_{x\in X} Z(x)$ se calcule niveau par niveau, or pour tout entier $n$, on a
$Z_n=\coprod_{x\in X} Z_n(x)$ car $\mathcal{C}$ est \discret. Ainsi
$\coprod_{x\in X} Z(x)$ n'est autre que $Z$ et $\mathcal{C-PC}$ vérifie la propriété 1 de la
définition~\ref{mod}.\\

On remarque facilement que se donner un morphisme dans $\mathcal{C-PC}$ de l'objet final vers un coproduit de
l'objet final équivaut à se donner un morphisme dans $\mathcal{C}$ de
l'objet final vers le même coproduit de l'objet final, car pour un objet
discret de $\mathcal{C-PC}$, tous les niveaux sont isomorphes au niveau 0. Ainsi
il vient pour tout ensemble $E$ : 
$$\operatorname{Hom}_{\mathcal{C-PC}}\Big(*,\coprod_{E}*\Big)=\operatorname{Hom}_{\mathcal{C}}\Big(*,\coprod_{E}*\Big)\cong
E $$
où l'isomorphisme est naturel en $E$, car $\mathcal{C}$ vérifie la
propriété 2 de la définition\\ \ref{mod}. On a donc montré
la propriété 2 pour $\mathcal{C-PC}$.\\

Pour montrer la stabilité des objets discrets par limite, il
suffit de rappeler que les limites dans une catégorie de
préfaisceaux sont les limites niveau par niveau. Or niveau par
niveau les objets discrets sont stables par limites car $\mathcal{C}$
est \discret. De plus nos objets étant discrets, leurs niveaux sont tous
isomorphes. Donc nos limites d'objets discrets sont niveau par
niveau des objets discrets et tous leurs niveaux sont isomorphes, ce qui en fait
des objets discrets et montre la propriété 3 de la définition~\ref{mod}.\\ 
CQFD.

\begin{ex}
Dans l'exemple précédent, on a vu qu'$\mathcal{ENSSIMP}$ la caté\-gorie des
ensembles simpliciaux n'est autre que $\mathcal{ENS-PC}$. Par la proposition
précé\-dente, on obtient donc que la catégorie $\mathcal{ENSSIMP}$ est
\discret.\\
Avec la catégorie \discret~$\mathcal{ENSSIMP}$, les \precats~sont les
ensembles bi-simpliciaux dont l'espace des objets est discret.
\end{ex}

\begin{prop}\label{precatclos}
Soit $\mathcal{C}$ une catégorie \discret.\\
Si $\mathcal{C}$ est complète (respectivement cocomplète), alors $\mathcal{C-PC}$ est complète
(respectivement cocomplète).
\end{prop}
{\it Preuve :}\\ 
Les catégories de préfaisceaux
vers une catégorie complète (respectivement cocomplète) sont complètes
(respectivement cocomplètes) avec les limites (respectivement les colimites) niveau par
niveau. Et comme les limites (respectivement les colimites) de $\mathcal{C}$ préservent la discrétude des
objets, les limites (respectivement les colimites) niveau par niveau préservent la condition d'être
discrètes au niveau 0. Par conséquent, $\mathcal{C-PC}$ est complète (respectivement cocomplète) avec
les limites (respectivement les colimites) niveau par niveau.\\
CQFD.\\

Munis de ces propositions, on peut définir par récurrence une notion de
$n$-\precat.

\begin{prop}[définition]\index{$n$-\precat}
Si $\mathcal{C}$ est une catégorie \discret, on appelle 1-\precats~les
\precats. Soit $n$ un entier strictement supé\-rieur à un et supposons
définie la catégorie $(n-1)-\mathcal{C-PC}$ des $(n-1)$-\precats, on appelle
alors $n$-\precats~les $(n-1)-\mathcal{C-PC}$-précaté\-gories et
$n-\mathcal{C-PC}$\index{$n-\mathcal{C-PC}$} la sous-catégorie pleine de la
catégorie des préfaisceaux de $\Delta$ vers $(n-1)-\mathcal{C-PC}$ d'objets les
$n$-\precats.\\
Si $\mathcal{C}$ est complète (respectivement cocomplète) et que ses objets
discrets sont stables par
limites (respectivement par colimites), alors, pour tout entier $n$, $n-\mathcal{C-PC}$ est complète
(respectivement cocomplète) et ses objets discrets sont stables par
limites (respectivement par colimites).
\end{prop}
{\it Preuve :} c'est une démonstration par récurrence qui utilise les deux
propositions précédentes à chaque étape. CQFD.

\begin{ex}
En prenant pour catégorie \discret~$\mathcal{ENS}$, on obtient les
$n$-précatégories de Tamsamani de \cite{s}.
\end{ex}

\begin{ex}
En prenant pour catégories \discret~$\mathcal{ENSSIMP}$, on obtient les
$n$-précatégories de Segal de \cite{h-s}. On peut remarquer que comme la
catégorie $\mathcal{ENSSIMP}$ est en fait celle des 1-précatégories de
Tamsamani, pour tout entier $n$, les $n$-précatégories de Segal ne sont
autres que les $n+1$-précatégories de Tamsamani.
\end{ex}

\newpage

\section{\cats}

Nous avons vu que pour définir une \precat, il fallait juste que $\mathcal{C}$ soit
\discret. Mais pour définir une notion de \cat, la seule
donnée de $\mathcal{C}$ ne suffit pas. Il faut pouvoir prendre en compte des
notions de catégories et d'\eqs~de catégories dans $\mathcal{C}$, que nous
appellerons respectivement \obcs~et \eqcs~d'\obcs. Nous allons donc définir la notion de
donnée de Segal qui prendra en compte ces éléments et à partir de
laquelle la notion de \cat~découlera. Comme en outre pour définir par la
suite la notion d'\eq~de \cats~nous avons besoin de savoir si des objets d'une
\cat~sont équivalents, nous allons demander aussi à la donnée de Segal un
foncteur noté $\tau_0$ qui a pour vocation d'associer à un \obc~l'ensemble
des classes d'\eq~de ses objets.

\begin{defin}\index{donnée de Segal}
Une donnée de Segal est un quadruplet
$(\mathcal{C},\mathcal{C}_c,\mathcal{C}_{eq},\tau_0)$ consti\-tué d'une
catégorie $\mathcal{C}$ \discret~et possédant les produits fibrés de deux
objets au-dessus d'un objet discret, d'une sous-catégorie pleine 
$\mathcal{C}_c$ de $\mathcal{C}$, d'une sous-catégorie $\mathcal{C}_{eq}$ de
$\mathcal{C}_c$, ayant les mêmes objets que $\mathcal{C}_c$, et d'un foncteur $\tau_0$ de $\mathcal{C}_c$ vers la
catégorie des ensembles, satisfaisant les propriétés suivantes :
\item 1) $\mathcal{C}_c$ est replète et contient les objets discrets de $\mathcal{C}$.
\item 2) Le produit fibré de deux objets de $\mathcal{C}_c$ au-dessus d'un
objet discret appartient à $\mathcal{C}_c$.
\item 3) $\mathcal{C}_{eq}$ contient les isomorphismes de $\mathcal{C}_c$ et les
produits fibrés, dans la catégories des morphismes de $\mathcal{C}$, de deux
morphismes de $\mathcal{C}_{eq}$ au-dessus d'un objet discret.
\item 4) Etant donné un diagramme :
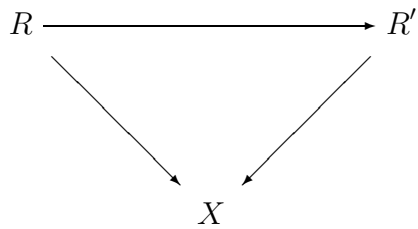
\begin{diagram}
R & & \rTo & & R'\\
  & \rdTo & & \ldTo & \\
  & & X & & \\
\end{diagram}
avec $R$ et $R'$ des \obcs~et $X$ objet discret, $R\rightarrow R'$ est un morphisme de $\mathcal{C}_{eq}$
si et seulement si pour tout élément $x$ de $X$, les morphismes induits
$R(x)\rightarrow R'(x)$ sont des morphismes de  $\mathcal{C}_{eq}$.
\item 5) La restriction du foncteur $\tau_0$ aux objets discrets n'est autre que le foncteur $\operatorname{Hom}_{\mathcal{C}}(*,-)$.
\item 6) $\tau_0$ préserve le produit fibré de deux objets de $\mathcal{C}_c$ au-dessus d'un objet discret.
\item 7) $\tau_0$ envoie les morphismes de $\mathcal{C}_{eq}$ sur les bijections
ensemblistes.\index{$\tau_0$!dans $\mathcal{C}$}
\end{defin}

On appellera \obcs~les objets de $\mathcal{C}_c$ et
\eqcs~d'\obcs~les morphismes de $\mathcal{C}_{eq}$.
\index{objet régal}\index{\eqc~d'\obcs} A partir de maintenant, par abus de notation
nous noterons indifféremment $\mathcal{C}$ la donnée de Segal
$(\mathcal{C},\mathcal{C}_c,\mathcal{C}_{eq},\tau_0)$ et sa catégorie
sous-jacente.\\

Au niveau des objets discrets, on aura pu remarquer qu'on leur a demandé d'être des \obcs~et que le foncteur $\tau_0$ se restreint sur eux en le foncteur $\operatorname{Hom}_{\mathcal{C}}(*,-)$. Ainsi si l'on identifie les objets discrets avec leurs ensembles sous-jacents, on obtient que la restriction aux objets discrets de $\tau_0$ n'est autre que l'identité. En outre les \eqcs~d'objets discrets sont des \eqcs~d'\obcs~et donc s'envoyent par $\tau_0$ sur les bijections d'ensembles. Or comme l'on vient de voir que la restriction du foncteur $\tau_0$ aux objets discrets n'est autre que l'identité, on obtient que les \eqcs~d'objets discrets sont des isomorphismes d'objets discrets. La réciproque est incluse dans la propriété 3) de donnée de Segal. Ainsi les \eqcs~d'\obcs~entre objets discrets sont exactement les isomorphismes d'objets discrets.  

\begin{ex}\label{hypens1}
Pour $\mathcal{ENS}$, on peut prendre la donnée de Segal suivante :
la sous-catégorie $\mathcal{ENS}_c$ est $\mathcal{ENS}$ tout entière,
les morphismes de $\mathcal{ENS}_{eq}$ sont les bijections ensemblistes. Enfin
pour $\tau_0$, on prend l'identité. Il est facile de voir que
cette donnée est bien une donnée de Segal.
\end{ex}

\begin{ex}\label{hypenssimp1}
Pour $\mathcal{ENSSIMP}$, on peut prendre la donnée de Segal suivante :
la sous-catégorie $\mathcal{ENSSIMP}_c$ est $\mathcal{ENSSIMP}$ tout entière,
les morphismes de $\mathcal{ENSSIMP}_{eq}$ sont les équivalences faibles
d'ensembles simpliciaux (i.e. les morphismes dont les réalisations induisent
des équivalences sur les groupes d'homotopie). Enfin pour $\tau_0$, on prend la composée du foncteur réalisation
géométrique par le foncteur composante connexe. Il est facile de voir que
cette donnée est bien une donnée de Segal.
\end{ex}

\begin{ex}
$(\mathcal{C},\mathcal{C}_c,\mathcal{C}_{eq},\tau_0)$ une donnée de Segal, il
est assez facile de voir que si on remplace les \eqcs~d'\obcs~par les
isomorphismes de $\mathcal{C}$ entre \obcs, on obtient une nouvelle donnée de
Segal plus stricte que la précédente. En effet, la seule chose à
vérifier est la préservation des isomorphismes par fibre et par produit fibré au-dessus
d'un objet discret, ce qui est vrai car $\mathcal{C}$ est une catégorie
\discret.
\end{ex}

\begin{defin}\index{\cat}
Soit $\mathcal{C}$ une donnée de Segal, on dit d'une \precat~$A$ qu'elle est une \cat~si, pour tout entier $m$ strictement positif, $A_m$
est un \obc~et si, pour tout entier $m$ supérieur ou
égal à deux, le morphisme de Segal
$$\segal{A}{m}$$ induit par les applications de {\bf 1}
dans {\bf m} qui à 0 et 1 associent i et i+1 (pour i compris entre 0 et m-1),
est une \eqc~d'\obcs. Un morphisme de \cats~est un morphisme
de \precats. 
\end{defin}

{\bf Notations :}\\ 
$A_1(x,y)$ note la fibre en $(x,y)$ de l'application $A_1\rightarrow
A_0\times A_0$ induite par les morphismes source et but. C'est donc l'\obc~des morphismes de $A$ allant de $x$ à $y$. Plus
généralement, $A_m(a_0,\ldots,a_m)$ note la fibre en $(a_0,\ldots,a_m)$ de
l'application $A_m\rightarrow A_0\times\ldots\times A_0$ induite par les
applications sommets de {\bf 0} vers {\bf m} qui envoyent 0 sur i, pour i compris
entre 0 et m. On considère $A_m(a_0,\ldots,a_m)$ comme l'\obc~des $m$-simplexes construits sur la suite de sommets
$(a_0,\ldots,a_m)$.\\ 

\begin{lem}\label{cat2}
Soit $\mathcal{C}$ une donnée de Segal, une \precat~$A$ est une \cat~si et
seulement si, pour tout entier $m$ strictement positif, $A_m$
est un \obc~et si, pour tout entier $m$ supérieur ou
égal à deux, et pour tout $m$-uplet $(a_0,\ldots,a_m)$ d'objets de $A$, le
morphisme $$\sefb{A}{m}{a}$$ induit par le morphisme de Segal, est une \eqc~d'\obcs.
\end{lem}
{\it Preuve :} application directe de la propriété 4 de donnée de Segal au
diagramme :
\begin{diagram}
A_m & & \rTo & & A_1\times_{A_0}\ldots\times_{A_0} A_1\\
  & \rdTo & & \ldTo & \\
  & & A_0\times\ldots\times A_0 & & \\
\end{diagram}
où les morphismes diagonaux sont induits par les applications sommets de {\bf
0} vers {\bf m} et de {\bf 0} vers $\mbox{{\bf 1}}\coprod_{\mbox{{\bf
0}}}\ldots\coprod_{\mbox{{\bf 0}}}\mbox{{\bf 1}}$).\\
CQFD.

\newpage

\section{Equivalences de \cats}

Maintenant que nous avons une notion de \cat, il nous faut aussi définir une
notion d'\eq~de \cats. L'idée est encore une fois de s'inspirer du cas des
petites catégories. Comme une \eq~de catégories est un foncteur
essentiellement surjectif et pleinement fidèle, nous allons aussi demander de
telles propriétés à nos \eqs~de \cats. Si la pleine fidélité peut se
définir de manière très similaire à la pleine fidélité des
foncteurs, à la différence près que nous ne demanderons pas un
isomorphisme mais seulement une \eqc, l'essentielle surjectivité fait quant
à elle intervenir la notion d'\eq~entre objets d'une \cat, notion qui est
gérée par le foncteur $\tau_0$. Ainsi afin de définir la notion d'\eq~de \cats, construisons
maintenant un foncteur $\tau_0$ pour les \cats. Pour cela nous
commençons par définir un foncteur $\tau_1$ des \cats~vers les ensembles
simpliciaux. Ce foncteur aura en outre la propriété de tomber dans les nerfs
de catégories.

\begin{defin}\index{$\tau_1$}
Soit $\mathcal{C}$ une donnée de Segal. On définit le foncteur~$\tau_1$ de
la sous-catégorie des \cats~vers la catégorie $\mathcal{ENSSIMP}$
ainsi :\\
pour tout entier $m$, $\tau_1(A)_m=\tau_0(A_m)$ et, pour tout
$f:\mbox{{\bf m}}\rightarrow\mbox{{\bf n}}$, $\tau_1(A)(f)=\tau_0(A(f))$. Enfin
pour $F:A\rightarrow B$ morphisme de \cats, on pose,
pour tout entier $m$, $\tau_1(F)_m=\tau_0(F_m)$.
\end{defin}

\begin{lem}
Soit $\mathcal{C}$ une donnée de Segal.
Le foncteur $\tau_1$, défini ci-dessus, a valeur dans la sous-catégorie des
nerfs de petites catégories.
\end{lem}
{\it Preuve :}\\
Comme $A$ est une \cat~et que par définition $\tau_0$ préserve les produits
fibrés de deux \obcs~au dessus d'un objet discret et envoie les \eqcs~d'\obcs~sur des
bijections d'ensembles, il vient que $\tau_1(A)$ est un nerf de
catégorie car ses morphismes de Segal sont bijectifs. Donc le foncteur
$\tau_1$ a valeurs dans les nerfs de petites catégories.\\
CQFD.\\

Ce foncteur $\tau_1$, qui transforme les \cats~en nerfs de petites catégories,
permet de définir le foncteur $\tau_0$ des classes d'\eq~d'objets d'une \cat, car un tel
foncteur existe naturellement pour les petites catégories. Par la suite, nous
identifierons par abus de langage les petites catégories et leurs nerfs.

\begin{defin}\index{$\tau_0$!pour \cat}
Soit $\mathcal{C}$ une donnée de Segal. On définit $\tau_0$ pour les \cats~comme la composée de $\tau_1$ par le foncteur qui à
une (petite) catégorie associe l'ensemble des classes d'équivalence de ses
objets (ce foncteur n'est autre que le $\tau_0$ des petites catégories).
\end{defin}

Munis du foncteur $\tau_0$ qui nous donne une notion d'\eq~entre objets d'une 
\cat, nous pouvons désormais définir la notion d'\eq~de \cats.

\begin{defin}\index{equivalence!de \cat}
Soit $\mathcal{C}$ une donnée de Segal. Un morphisme de \cats\\ $f:A\rightarrow B$ est une équivalence de \cats~si l'application d'ensembles 
$\tau_0(f):\tau_0(A)\rightarrow \tau_0(B)$ est surjective (on dit alors que $f$
est essentiellement surjective)\index{\cat!essentiellement surjective} et, si pour tout couple $(x,y)$ d'objets de A, le morphisme $f_1(x,y):A_1(x,y)\rightarrow
B_1(f(x),f(y))$ est une \eqc~d'\obcs~(on dit alors que
$f$ est pleinement fidèle).\index{\cat!pleinement fidèle}
\end{defin}

Donnons maintenant les premiers exemples de données de Segal avec leurs notions de
\cats~et d'\eqs~de \cats~associées.

\begin{ex}
En prenant la donnée de Segal de l'exemple~\ref{hypens1}, on obtient pour
\cats~les petites catégories, pour \eqs~de \cats~les \eqs~de catégories et
pour $\tau_0$ le foncteur qui à une catégorie associe l'ensemble des classes
d'\eq~de ses objets.
\end{ex}

\begin{ex}
En prenant maintenant sur la catégorie des ensembles simpliciaux la donnée
de Segal de l'exemple~\ref{hypenssimp1}, on obtient les notions de \precat,
de \cat~et d'équivalence de \cats~qui ne sont autres que celles de
précatégorie de Segal, de catégorie de Segal et d'équivalence de
catégories de Segal introduite par Hirschowitz et Simpson dans \cite{h-s}.
\end{ex}
\index{catégorie!de Segal}

Il est intéressant de constater, en particulier pour pouvoir itérer les
constructions de \cats~et d'\eqs~de \cats, que la catégorie $\mathcal{C-PC}$
des \precats~forme une donnée de Segal avec les \cats, les \eqs~de
\cats~et le foncteur $\tau_0$ pour les \cats.

\begin{lem}\label{cpcse}
Soit $\mathcal{C}$ une donnée de Segal.
Munie des notions de \cat, d'équivalence de \cats~et du foncteur $\tau_0$
définis ci-dessus, la catégorie $\mathcal{C-PC}$ des \precats~constitue une
donnée de Segal.
\end{lem}
{\it Preuve :}\\
D'après la proposition~\ref{precatdiscret}, la catégorie $\mathcal{C-PC}$ est \discret~et,
par la proposition~\ref{precatclos}, elle possède les produits fibrés de
deux \precats~au-dessus d'un objet discret. La sous-catégorie des \cats~a
été définie comme sous-catégorie pleine de $\mathcal{C-PC}$.
En outre, il est facile de voir que les isomorphismes
de \cats~et les composées d'équivalences de \cats~vérifient la
définition d'équivalence de \cats, car $\mathcal{C}_{eq}$ est une
sous-catégorie contenant les isomorphismes, ce qui montre que la
classe des \eqs~de \cats~donne lieu à une sous-catégorie de la catégorie
des \cats~ayant tous les objets et vérifiant la première partie de la propriété 3 de donnée de
Segal.\\

On remarque que si une \precat~est isomorphe à une \cat, alors
chacun de ses niveaux est isomorphe à un \obc~et,
comme $\mathcal{C}_c$ est replète, chaque niveau de notre \precat~est un 
\obc. En outre il est facile de voir que pour notre
\precat~isomorphe à une \cat, le morphisme de Segal est bien une
\eqc~d'\obcs~car $\mathcal{C}_{eq}$ contient
les isomorphismes et est stable par composition. Ceci montre que toute
\precat~isomorphe à une \cat~est une \cat, ce qui prouve la première partie
de la propriété 1 des données de Segal. On remarque ensuite que les objets discrets de $\mathcal{C-PC}$ sont
clairement des \cats, car niveau par niveau ce sont des objets discrets et que
les applications de Segal sont des isomorphismes. Ce qui prouve la seconde
partie de la propriété 1.\\ 

En outre les notions de \cat~et d'\eq~de \cats~sont stables par produit fibré au-dessus d'un
objet discret, car le produit fibré au-dessus d'un
objet discret est calculé niveau par niveau et que niveau par niveau on
applique la propriété 2 et la seconde partie de la propriété 3 pour
$\mathcal{C}$. Ainsi la propriété 2 et la seconde partie de la propriété 3 sont
vérifiées.\\

On peut remarquer par construction de $\tau_1$ que l'image par $\tau_1$ d'un
objet discret n'est autre que la catégorie discrète d'objets les morphismes
de l'objet final de $\mathcal{C-PC}$ vers l'objet discret et donc $\tau_0$ pour
les \cats~vérifie la propriété 5. Par composition de foncteurs
préservant le produit fibré au-dessus d'un objet discret
, le foncteur $\tau_0$ pour les \cats~a lui aussi la propriété 6.
La définition d'équivalence de \cats~entraîne que $\tau_1(f)$ est une
équivalence de petites catégories. Une des conséquences de la pleine fidélité
de $\tau_1(f)$ est que $\tau_0(f)$ est injectif. D'où le foncteur $\tau_0$ pour les
\cats~envoie bien les équivalences de \cats~sur des bijections et donc
vérifie la propriété 7.\\

Considérons le diagramme suivant dans $\mathcal{C-PC}$ :
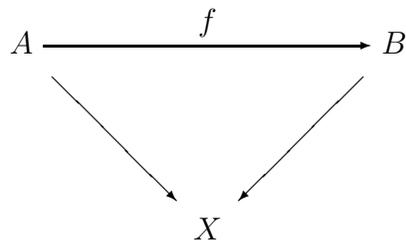
\begin{diagram}
A & & \rTo^{f} & & B\\
  & \rdTo & & \ldTo & \\
  & & X & & \\
\end{diagram}
où X est une \precat~discrète et $f$ un morphisme de \cats. La seconde partie de la propriété 3, que
l'on a montrée pour $\mathcal{C-PC}$, entraîne que si $f$ est une
\eq~de \cats, alors toutes ses fibres $f(x):A(x)\rightarrow B(x)$ le sont.
Supposons maintenant que pour tout élément $x$ de $X$, les morphismes
$f(x):A(x)\rightarrow B(x)$ sont des \eqs~de \cats. Comme la catégorie des
ensembles est \discret~et que l'on a montré que $\tau_0$ commute aux fibres, il vient que
$\tau_0(f)$ est le coproduit des $\tau_0(f(x))$. Or on a montré que $\tau_0$
vérifie la propriété 7, donc tous les $\tau_0(f(x))$ sont bijectifs et
leur coproduit $\tau_0(f)$ aussi, ce qui montre l'essentielle surjectivité de
$f$. Pour montrer la pleine fidélité de $f$, on applique à $f_1(a,a')$,
pour tout couple $(a,a')$ d'objets de A, la propriété 4 de donnée de Segal vérifiée par
$\mathcal{C}$. Soit donc $x$ un élément de $X$. Si $a$ et $a'$ sont des
objets de $A(x)$, alors $f(a)$ et $f(a')$ sont des objets de $B(x)$, et l'on a
que la fibre en $x$ de $f_1(a,a')$ n'est autre que $f(x)_1(a,a')$, qui est une
\eqc. Sinon la fibre en $x$ de $f_1(a,a')$ est l'identité de l'ensemble vide,
qui est bien une \eqc. Donc par la propriété 4 vérifiée par
$\mathcal{C}$, $f_1(a,a')$ est une \eqc, ce
qui montre la pleine fidélité de $f$. On a donc montré que $\mathcal{C-PC}$
vérifie la propriété 4.\\
CQFD.\\

Dans la démonstration du lemme, on a obtenu le résultat important suivant :

\begin{cor}
Soit $\mathcal{C}$ une donnée de Segal.
Si le morphisme de \cats~$f:A\rightarrow B$ est une équivalence de \cats, alors $\tau_1(f)$ est une
équivalence de petites catégories et $\tau_0(f)$ est une bijection
d'ensembles.
\end{cor}

Un autre corollaire de ce lemme est que l'on peut définir les $n$-\cats.

\begin{cor}\index{$n$-\cat}\index{\eq!de $n$-\cats}
Si $\mathcal{C}$ est une donnée de Segal, on appelle 1-\cats~les
\cats~et 1-\eqs~les \eqs~de \cats.\\ Soit $n$ un entier strictement supérieur à un et supposons
définie la donnée de Segal $(n-1)-\mathcal{C-PC}$ des $(n-1)$-\cats, on appelle
alors $n$-\cats~les $(n-1)-\mathcal{C-PC}$-catégories, $n$-\eqs~les \eqs~de
$n$-\cats~et on note $n-\mathcal{C-PC}$ la donnée de
Segal sur la catégorie $n-\mathcal{C-PC}$ des $n$-\precats~formée par les
$n$-\cats, les $n$-\eqs~et le foncteur $\tau_0$ des $n$-\cats.
\end{cor}
{\it Preuve :} c'est une démonstration par récurrence immédiate qui
utilise le lemme précédent à chaque étape. CQFD.

\begin{ex}
En prenant la donnée de Segal de l'exemple~\ref{hypens1} et en lui appliquant
le corollaire ci-dessus, on obtient pour
$n$-\cats~les $n$-nerfs, pour \eqs~de \cats~les $n$-\eqs~entre $n$-nerfs et
pour $\tau_0$ le foncteur de troncature tels que ces trois notions ont été
définies dans \cite{t}.\\ 
En particulier, en appliquant seulement le lemme~\ref{cpcse} à la donnée de
Segal de l'exemple~\ref{hypens1}, on obtient que la
catégorie des ensembles simpliciaux, munie des notions de catégories,
d'équivalences de catégories et du foncteur $\tau_0$ des classes
d'\eq~d'objets, forme une
donnée de Segal. Pour cette dernière, les \precats~sont les ensembles
bi-simpliciaux à espace d'objets discret, les \cats~sont les 2-nerfs et les
\eqs~de \cats~sont les \eqs~de 2-nerfs de \cite{t}. Avec la
donnée de Segal qui engendre les catégories de Segal, cela nous fait deux
données de Segal différentes sur les ensembles simpliciaux.
\end{ex}

\begin{ex}
Reprenons la donnée de Segal de l'exemple~\ref{hypens1}. On obtient par le
lemme~\ref{cpcse}, une donnée de Segal sur les ensembles simpliciaux avec les
catégories et les \eqs~de catégories. Prenons maintenant la donnée de
Segal plus stricte obtenue en remplaçant les \eqs~de catégories par les
isomorphismes, on obtiendra alors pour \cats~les 2-catégories strictes et
pour \eqs~de \cats~les morphismes essentiellement surjectifs et strictement
pleinement fidèles. Par le lemme~\ref{cpcse}, ces notions forment une donnée
de Segal sur les 2-\precats~de Tamsamani. Et comme préce\-demment, on peut
considérer la donnée de Segal plus stricte sur les 2-\precats~en prenant les
isomorphismes au lieu des \eqs~de 2-catégories, on obtiendra alors les
3-catégories avec pour \eqs~les morphismes essentiellement surjectifs et
strictement pleinement fidèles. Et ainsi de suite...\\
Par récurrence, en strictifiant à chaque étape la donnée de Segal
obtenue, c'est-à-dire en remplaçant les \eqs~par les isomorphismes, on
retrouve les $n$-catégories strictes.\\
On remarque au passage que l'on vient de munir la catégorie des ensembles
simpliciaux d'une troisième donnée de Segal qui engendre les 2-catégories
strictes.
\end{ex}

\begin{ex}
En prenant la donnée de Segal de l'exemple~\ref{hypenssimp1} et en lui appliquant
le corollaire ci-dessus, on obtient pour
$n$-\cats~les $n$-catégories de Segal, pour \eqs~de \cats~les \eqs~de
$n$-catégories de Segal et
pour $\tau_0$ le foncteur de troncature des $n$-catégories de Segal tels que ces trois notions ont été
définies dans \cite{h-s}.\\ 
Comme on a déjà remarqué que les $n$-\precats~de Segal sont exactement les
$n+1$-\precats~de Tamsamani, on obtient qu'il existe au moins trois données de
Segal différentes sur la catégories des $n$-\precats~de Segal dont l'une va
engendrer les $n+1$-catégories de Segal, l'autre les $n+2$-catégories
strictes et la troisième les $n+2$-nerfs de Tamsamani.
\end{ex}

\newpage

\section{Problème de la catégorisation}

Les catégories enrichies faibles ainsi définies ne sont pas stables par
colimite. Pour s'en convaincre, considérons deux exemplaires de la catégorie
ayant un seul morphisme et recollons les de telle sorte que les deux morphismes
soient composables. Ce recollement correspond bien à la somme amalgamée dans
les ensembles simpliciaux mais pas à celle des catégories car le résultat
obtenu est une précatégorie avec deux morphismes composables sans
composition correspondante. L'une des idées permettant de remédier à ce
problème est de catégoriser le résultat obtenu. Ainsi on définirait les
colimites de catégories comme les catégorisées des colimites des
précatégories sous-jacentes.\\

On s'est donc ramener à un problème de catégorisation. On cherche donc un
foncteur, notons-le $Cat$, muni d'une transformation naturelle notée $can$ de
l'identité vers $Cat$, qui à toute \precat~associe une \cat. Bien
évidemment on voudrait que la \cat~associée par ce foncteur soit la
minimale, c'est à dire que tout morphisme d'une \precat~dans une \cat~se
factorise par la catégorisation de la \precat.
\begin{diagram}
A & & \rTo &  & B_{\mbox{ catégorie}} \\
& \rdTo^{can_A} & & \ruDotsto_{\exists} & \\
 & & Cat(A) & & \\
\end{diagram}

En outre on aimerait aussi que la catégorisation donne une \cat~associée
équivalente à la \precat~de départ au sens suivant : soit $A$ une \precat,
le morphisme $Cat(can_A):Cat(A)\rightarrow Cat(Cat(A))$ est une \eq~de \cats.
Ceci nous donnerait d'ailleurs une notion d'\eq~de \precats. On pourra en effet
poser qu'un morphisme de \precats~est une \eq~si sa catégorisation l'est.\\

Comme il est assez difficile de trouver une telle catégorisation, en suivant
l'idée de \cite{s}, nous allons ramener ce problème à un problème de
relèvement de diagramme. S'il existe en effet une famille I de flèches telle que les objets ayant la \prd~par rapport à I sont des \cats, une bonne façon d'avoir ces propriétés de relèvement est de construire
la catégorisation comme une colimite séquentielle transfinie de sommes
amalgamées des flèches de I, ce que nous appellerons un plan infini d'addition de cellules de I.\\

 Dans cette optique nous allons tout d'abord traiter de ces
problèmes de relèvement et de ceux d'engendrement qui leur sont liés avant
de chercher par rapport à quelles flèches nos \cats~pourraient se relever.

\newpage
...

\chapter{Engendrements}

\newpage

Dans ce chapitre nous allons tout d'abord rappeler certaines définitions et
certains résultats classiques sur les relèvements de diagrammes. Une des
définitions principales est celle d'objet I-injectif qui est un objet ayant la
\prd~par rapport à l'ensemble de flèches I. En effet le but de ce chapitre
est de donner des constructions rendant I-injectifs les objets qui ne le sont
pas. Pour ce faire, nous allons utiliser certaines techniques apparaissant dans l'argument
du petit objet de Quillen qui permet de répondre en partie à notre
problème de rendre I-injectif. L'un des points clefs de cet argument, outre la
notion de petitesse, est le procédé de colimite séquentielle transfinie de
sommes amalgamées de flèches de I, ce que nous appellerons plan d'addition
de cellules.\\
\\

Nous donnerons alors une formalisation de plan d'addition de
cellules et nous donnerons des conditions pour lesquelles ces plans d'addition
de cellules $P$ I-injectivisent, c'est-à-dire répondent aux trois critères suivants :\\
- pour tout objet $A$, $P(A)$ est I-injectif,\\
- tout morphisme $A\rightarrow B$ vers un objet I-injectif se factorise à
travers le morphisme naturel $A\rightarrow P(A)$,\\
- pour tout objet $A$, l'image par $P$ du morphisme naturel $A\rightarrow P(A)$
est une \eq~d'objets I-injectifs, si l'on possède une notion naturelle
d'\eq~entre objets I-injectifs.\\
\\

Si les deux premiers critères sont simples à traiter, ce n'est pas le cas du
troisième qui nous obligera d'une part à trouver un procédé rendant
I-injectif avec la propriété d'avoir unicité des factorisations des
morphismes à but I-injectifs. Pour cela, nous ferons intervenir la notion de
marquage d'un relèvement.
D'autre part, nous allons montrer que, muni d'un tel procédé et en supposant
l'existence d'un autre procédé d'I-injectivisation qui ne vérifie la
propriété de stabilité homotopique que pour les objets I-injectifs, le
procédé I-injectivant à factorisation unique vérifie bien les trois
critères.\\
\\

Comme l'une des utilités de l'I-injectivisation est aussi de pourvoir les
objets I-injectifs de colimites, nous allons terminer ce chapitre en comparant
les plans d'addition de cellules via un procédé de rationalisation des
plans et appliquer cette comparaison à certains types de plans qui sont
stables lorsqu'on les compose avec l'un de leurs sous-plans. Avec ces plans, on
pourra comparer la colimite des I-injectivés avec l'I-injectivisation de la
colimite, ce qui nous permettra d'obtenir une bonne caractérisation d'\eq~pour
les colimites d'objets quelconques en termes de colimites de leurs
I-injectivés.

\newpage

\section{Rappels sur les relèvements de diagrammes}

Comme l'une des propriétés que l'on demande à la catégorisation est une
propriétés de relèvement, nous allons préciser cette notion dans la
définition suivante prise dans \cite{h}:

\begin{defin}
Soient $i:A\rightarrow B$ et $p:C\rightarrow D$ deux morphismes tels que pour
tous les carrés commutatifs du type suivant :
\begin{diagram}
A & \rTo & C\\
\dTo^{i} & \ruDotsto_{\exists} & \dTo_{p}\\
B & \rTo & D\\
\end{diagram}
la flèche en pointillé existe et fait commuter les deux parties du
diagramme.\\ On dira alors que:
\item -cette flèche est un relèvement de ce diagramme,
\item -$i$ a la \prg~par rapport à $p$,
\item -$p$ a la \prd~par rapport à $i$,
\item -$C$ a la \prd~par rapport à $i$, si $D$ est un objet final.
\end{defin}
\index{relèvement!propriété de}
\index{relèvement!à gauche}
\index{relèvement!à droite}
 
Pour la suite, nous aurons souvent besoin de certaines sommes amalgamées ou de
certains produits fibrés. Pour simplifier les appellations, nous prendrons
toujours en suivant \cite{h} la convention suivante.

\begin{defin}
Soit un carré commutatif :
\begin{diagram}
A & \rTo^a & C \\
\dTo^f & & \dTo_g \\
B & \rTo_b & D
\end{diagram}
\item - Si ce carré est cocartésien, on dit que $g$ est la somme amalgamée de $f$
le long de $a$.
\item - Si ce carré est cartésien, on dit que $f$ est le produit fibré de
$g$ le long de $b$.
\end{defin}

Les flèches ayant des propriétés de relèvement sont stables par
certaines limites et colimites, ce que nous rappelons ci-dessous.

\begin{lem}\label{prst}
\item - Les classes de flèches ayant une \prg~contiennent les isomorphismes et sont
stables par composition, par rétract, par somme amalgamée le long d'un
morphisme, par colimite séquentielle transfinie et par coproduit dans la
catégorie des morphismes.
\item - Les classes de flèches ayant une \prd~contiennent les isomorphismes et sont
stables par composition, par rétract, par produit fibré le long d'un
morphisme et par produit dans la catégorie des morphismes.
\end{lem}
{\it Preuve :}\\ 
Nous devons juste montrer la stabilité par coproduit (resp. produit)
dans la catégorie des morphismes car le reste de la preuve se trouve dans
l'ouvrage d'Hirschhorn~\cite{h}. Soient J un ensemble et $(f_j:A_j\rightarrow
B_j)_{j\in J}$ une famille de morphismes ayant la \prg~par rapport à un
morphisme $p$ fixé. Considérons le diagramme suivant :
\begin{diagram}
\coprod_{j\in J} A_j & \rTo^a & C \\
\dTo^{\coprod_{j\in J} f_j} & & \dTo_p \\
\coprod_{j\in J} B_j & \rTo_b & D
\end{diagram}
Notons $i_j$ le morphisme canonique de $A_j$ dans $\coprod_{j\in J} A_j$ et
$i'_j$ le morphisme canonique de $B_j$ dans $\coprod_{j\in J} B_j$. Comme chaque
$f_j$ se relève par rapport à $p$, il vient que pour tout $j$, il existe un
morphisme $\phi_j:B_j\rightarrow C$ tel que $p\circ\phi_j=b\circ i'_j$ et
$\phi_j\circ f_j=a\circ i_j$. Par propriété universelle du coproduit
$\coprod_{j\in J} B_j$, il existe un unique morphisme $\phi:\coprod_{j\in J}
B_j\rightarrow C$ tel que pour tout $j$, on ait $\phi\circ i'_j=\phi_j$. On
obtient donc les égalités suivantes : $p\circ\phi\circ
i'_j=p\circ\phi_j=b\circ i'_j$ et $\phi\circ\coprod_{j\in J} f_j\circ
i_j=\phi\circ i'_j\circ f_j=\phi_j\circ f_j=a\circ i_j$. Et par propriété
universelle des coproduits $\coprod_{j\in J} B_j$ et $\coprod_{j\in J} A_j$, il
vient que $p\circ\phi=b$ et $\phi\circ\coprod_{j\in J} f_j=a$ et donc que $\phi$
est un relèvement pour le diagramme $(\coprod_{j\in J} f_j,a,b,p)$. On a donc
montré que les morphismes ayant une \prg~sont stables par coproduit dans la
catégorie des morphismes. Dualement on aura que ceux ayant une \prd~sont
stables par produit dans la catégorie des morphismes.\\
CQFD.\\

On remarquera
au passage que pour la somme amalgamée (resp. le produit fibré) dans la catégorie des
morphismes, on n'a généralement pas de morphisme induit par les relèvements de chaque
membre de la somme amalgamée (resp. du produit fibré) car la non-unicité
des relèvements est une obstruction à leur compatibilité avec la somme
amalgamée (resp. le produit fibré). D'où la définition et le lemme
suivants.

\begin{defin}
Soient $i:A\rightarrow B$ et $p:C\rightarrow D$ deux morphismes tels que pour
tous les carrés commutatifs du type suivant :
\begin{diagram}
A & \rTo & C\\
\dTo^{i} & \ruDotsto_{\exists} & \dTo_{p}\\
B & \rTo & D\\
\end{diagram}
la flèche en pointillé existe et fait commuter les deux parties du
diagramme. On dira alors :
\item - si $i$ est muni d'un choix de relèvements pour chaque diagramme du
type ci-dessus avec $p$ fixé, que $i$ est marqué à gauche par rapport à $p$.
\item - si $p$ est muni d'un choix de relèvements pour chaque diagramme du
type ci-dessus avec $i$ fixé, que $p$ est marqué à droite par rapport à
$i$.\\
\\
On obtient alors une catégorie dont les objets sont les morphismes marqués
et les morphismes les carrés commutatifs compatibles avec le marquage.
\end{defin}

\begin{lem}\label{mst}
\item - La catégorie des morphismes marqués à gauche par rapport à un ensemble
de flèches est stables par colimites.
\item - La catégorie des morphismes marqués à droite par rapport à un ensemble
de flèches est stables par limites.
\end{lem}
{\it Preuve :}\\ 
Comme conséquence du lemme précédent, ces catégories sont stables
respectivement par coproduit et produit. Il suffit donc de montrer qu'elles sont
stables respectivement par co-égalisateur et par égalisateur pour avoir le
résultat. Soit le diagramme commutatif suivant :
\begin{diagram}
A_0 & \pile{\rTo^{f_0}\\ \rTo_{f'_0}} & B_0 & \rTo^{\tilde{f}_0} & C_0 & \rTo^x &
X\\
\dTo^a & & \dTo_b & & \dTo_c & & \dTo_p \\
A_1 & \pile{\rTo^{f_1}\\ \rTo_{f'_1}} & B_1 & \rTo^{\tilde{f}_1} & C_1 & \rTo^y &
Y\\
\end{diagram}
où $c$ est le morphisme universel du co-égalisateur $C_0$ de $f_0,f'_0$ vers le co-égalisateur $C_1$ de $f_1,f'_1$. 
Comme $a$ est marqué par rapport à $p$, il vient avec un morphisme
$\phi:A_1\rightarrow X$ tel que $p\circ\phi=y\circ\tilde{f}_1\circ
f_1=y\circ\tilde{f}_1\circ f'_1$ et $\phi\circ a=x\circ\tilde{f}_0\circ
f_0=x\circ\tilde{f}_0\circ f'_0$. De même comme $b$ est marqué par rapport
à $p$, il vient avec un morphisme $\psi:B_1\rightarrow X$ tel que
$p\circ\psi=y\circ\tilde{f}_1$ et $\psi\circ b=x\circ\tilde{f}_0$. Or les
morphismes $f$ et $f'$ sont compatibles avec les marquages et donc il vient
$\psi\circ f_1=\phi=\psi\circ f'_1$. Par la propriété universelle du
co-égalisateur $C_1$, il existe un morphisme $\chi:C_1\rightarrow X$ tel que
$\chi\circ\tilde{f}_1=\psi$. On obtient donc les égalités suivantes :
$p\circ\chi\circ\tilde{f}_1=p\circ\psi=y\circ\tilde{f}_1$ et
$\chi\circ c\circ\tilde{f}_0=\chi\circ\tilde{f}_1\circ b=\psi\circ
b=x\circ\tilde{f}_0$. Par la propriété universelle des co-égalisateurs $C_1$ et
$C_0$, il vient que $p\circ\chi=y$ et $\chi\circ c=x$. Donc $\chi$ est un
relèvement pour le diagramme $(c,x,y,p)$ et c'est celui qui marquera $c$ par
rapport à $p$. On a donc montré que la catégorie des morphismes marqués
à gauche par rapport à un morphisme fixé est stable par co-égalisateur.
Et on aura dualement la stabilité des morphismes marqués à droite par
égalisateur.\\
CQFD.\\

\newpage

\section{Morphismes I-injectifs et I-cofibrations}

Parmi la multitude de flèches ayant des propriétés de relèvement,
certaines nous intéressent plus que d'autres. En effet, souvent nous allons nous
fixer un ensemble $I$ de flèches et regarder la classe des flèches se
relevant à droite par rapport à I et celle des flèches se relevant à
gauche par rapport à ces dernières. C'est pourquoi nous rappelons les
définitions de $I$-injectifs et de I-cofibrations ainsi que certaines de
leurs propriétés, toujours en suivant \cite{h}.

\begin{defin}\index{I-injectif}\index{I-cofibration}
Soient $\mathcal{C}$ une catégorie et I un ensemble de morphismes de
$\mathcal{C}$.
\item 1) Un morphisme est I-injectif s'il a la \prd~par rapport aux morphismes
de I.
\item 2) Un objet est I-injectif si son unique morphisme vers l'objet final est
I-injectif.
\item 3) Un morphisme est une I-cofibration s'il a la \prg~par rapport aux
morphismes I-injectifs.
\end{defin}

\begin{lem}\label{Ist}
Soient $\mathcal{C}$ une catégorie et I un ensemble de morphismes de
$\mathcal{C}$.
\item - Les morphismes I-injectifs forment une sous-catégorie de $\mathcal{C}$
contenant les isomorphismes et stable par rétract, par produit fibré le
long d'un morphisme et par produit dans la catégorie des morphismes. 
\item - Les I-cofibrations forment une sous-catégorie de $\mathcal{C}$
contenant les isomorphismes et les morphismes de I et qui est stable par rétract, par somme amalgamée le
long d'un morphisme, par colimite séquentielle transfinie et par coproduit
dans la catégorie des morphismes.
\end{lem}
{\it Preuve :} application directe du lemme~\ref{prst}

\begin{lem}\label{inj}
Soient $\mathcal{C}$ une catégorie et I et K des ensembles de morphismes de
$\mathcal{C}$. Si les morphismes I-injectifs sont exactement les morphismes
K-injectifs alors les I-cofibrations sont exactement les K-cofibrations
\end{lem}

Nous arrivons maintenant à l'une des définitions les plus importantes de
cette partie, celle de famille génératrice d'une sous-catégorie de
morphismes.

\begin{defin}\index{famille de flèches!génératrice}
Soient $\mathcal{C}$ une catégorie et $\mathcal{C}'$ une sous-catégorie
ayant les mêmes objets. Soit I un ensemble de morphismes de $\mathcal{C}$. On
dit que I engendre $\mathcal{C}'$ si $\mathcal{C}'$ est la sous-catégorie des
I-cofibrations.
\end{defin}

La définition de famille génératrice d'une classe de morphismes que l'on
vient de donner n'est très intuitive. En fait, quand on parle d'une famille de flèches engendrant certains
morphismes, on pense à des morphismes fabriqués comme colimites
séquentielles transfinies de sommes amalgamées des flèches
génératrices. Ceci suggère donc que les I-cofibrations sont des
colimites séquentielles transfinies de sommes amalgamées d'éléments de
I. Mais pour avoir ce résultat, on se sert de l'argument du petit objet de
Quillen.

\newpage

\section{Argument du petit objet}

L'argument du petit objet de Quillen permet de résoudre le problème suivant.
On se donne une catégorie $\mathcal{C}$ cocomplète et une famille I de
morphismes de $\mathcal{C}$ et l'on cherche à factoriser tous morphismes de
$\mathcal{C}$ en une I-cofibration suivie d'un morphisme I-injectif. Un des
corollaires de ce résultat est que I engendre les I-cofibrations au sens
auquel on s'attend. Pour arriver à cette factorisation, l'idée de l'argument
du petit objet consiste à limiter la taille des flèches de I, ce que l'on
modélise par la notion de petitesse que nous allons rappeler. Mais auparavant rappelons les définitions de $\lambda$-séquence et de cardinal régulier.

\begin{defin}
Un cardinal $\alpha$ est régulier si, pour toute famille d'ensembles $(A_i)_{i\in I}$ telle que le cardinal de $I$ et les cardinaux des $A_i$, pour $i\in I$, sont strictement inférieurs à $\alpha$, le cardinal de la réunion $\cup_{i\in I} A_i$ est strictement inférieur à $\alpha$.
\end{defin}

\begin{defin}\index{$\lambda$-séquence}
Soient $\mathcal{A}$ une catégorie admettant les colimites séquentielles
transfinies et $\mathcal{F}$ une famille de morphismes
de $\mathcal{A}$. Soit $\lambda$ un ordinal.\\
Une $\lambda$-séquence de morphismes de $\mathcal{F}$ est un foncteur $X:\lambda\rightarrow \mathcal{C}$ tel que, pour tout ordinal $\beta<\lambda$, le morphisme $X_{\beta}\rightarrow X_{\beta+1}$ appartient à $\mathcal{F}$ et, pour tout ordinal limite $\beta<\lambda$, le morphisme $\colimite{\gamma<\beta} X_{\gamma}\rightarrow X_{\beta}$ est un isomorphisme.
\end{defin}

\begin{defin}\index{$\alpha$-petit}\label{alphapetit}
Soit $\alpha$ un cardinal régulier et $\mathcal{F}$ une famille de morphismes
d'une catégorie $\mathcal{A}$ admettant les colimites séquentielles
transfinies,
un objet $A$ de $\mathcal{A}$ est $\alpha$-petit par rapport à $\mathcal{F}$ si pour tout $\lambda$ cardinal régulier strictement
plus grand qu'$\alpha$ et pour toute $\lambda$-séquence de morphismes de $\mathcal{F}$ $X_0\rightarrow
X_1\rightarrow \ldots$, on a :
$$\colim{\beta<\lambda}
\mathrm{Hom}(A,X_{\beta})\stackrel{\cong}{\longrightarrow}\mathrm{Hom}(A,\colim{\beta<\lambda}
X_{\beta}) $$
On dit que $A$ est $\alpha$-petit s'il est $\alpha$-petit par rapport à tous
les morphismes de $\mathcal{A}$.
\end{defin}

Afin de montrer qu'un objet est petit, voici quelques résultats techniques,
inspirés de \cite{h}, sur la stabilité de la petitesse vis-à-vis de
certaines opérations ainsi que quelques résultats sur la petitesse des
\precats.

\begin{lem}\label{cpdtpetit}
Soient $\alpha$ un cardinal régulier et $\mathcal{A}$ une catégorie
admettant les coproduits et les colimites séquentielles
transfinies. Soit I un ensemble d'indices de cardinal strictement inférieur à $\alpha$ et
soit $(A_i)_{i\in I}$ une famille d'objets de $\mathcal{A}$ $\alpha$-petits alors $\coprod_{i\in
I} A_i$ est $\alpha$-petit.
\end{lem}

\begin{lem}\label{colimpetit}
Soient $\alpha$ un cardinal régulier et $\mathcal{A}$ une catégorie
admettant les coproduits et les colimites séquentielles
transfinies. Soit $F$ un foncteur d'une petite catégorie $\mathcal{I}$ vers la
catégorie $\mathcal{C}$ tel que pour tout objet $i$ de $\mathcal{I}$ l'objet
$F(i)$ de $\mathcal{C}$ est $\alpha$-petit. Alors la colimite de $F$,
$\colimite{i\in \mathcal{I}}F(i)$, est $\alpha'$-petite, avec $\alpha'$ un
cardinal plus grand qu'$\alpha$ et que le cardinal de l'ensemble des morphismes
de $\mathcal{I}$.
\end{lem}

\begin{lem}\label{precatpetit}
Soit $\alpha$ un cardinal régulier strictement plus grand qu'$\aleph_0$. Soit
$A$ un préfaisceau de $\Delta$ vers une catégorie $\mathcal{C}$ admettant
les colimites séquen\-tielles transfinies, si,
pour tout $n\in\mathbb{N}$, $A_n$ est $\alpha$-petit alors $A$ est
$\alpha$-petit.
\end{lem}

On va maintenant énoncer la définition qui précise que, sous certaines
conditions, un ensemble de flèches permet l'argument du petit objet de Quillen

\begin{defin}\index{famille de flèches!permettant l'argument du petit objet}\label{ptobj}
Soit $\mathcal{A}$ une catégorie admettant
les colimites séquentielles transfinies.
Soit E un ensemble de morphismes de $\mathcal{A}$, on dit que E
permet l'argument du petit objet si, pour tout morphisme de E, il existe un
cardinal régulier $\alpha$ tel que la source du morphisme est $\alpha$-petite
par rapport à E.
\end{defin}

Un des cas les plus fréquents d'ensembles de flèches permettant l'argument
du petit objet provient d'ensembles de flèches dont les sources et buts sont
petits. On dira de ces flèches qu'elles sont petites ou limitées. Bien
entendu, de tels ensembles permettent l'argument du petit objet, ce qu'énonce
le lemme suivant.

\begin{lem}\label{enspt}
Soient $\alpha$ un cardinal régulier et $\mathcal{A}$ une catégorie admettant
les colimites séquentielles transfinies.
Soit E un ensemble de morphismes de $\mathcal{A}$ dont les sources et
buts sont $\alpha$-petits, alors E permet l'argument du petit objet.
\end{lem}

Comme on l'a déjà dit plus haut, l'argument du petit objet sert à obtenir
la factorisation de toute flèche en une I-cofibration suivie d'un morphisme
I-injectif pour un ensemble I assez petit. Ce fait est énoncé dans le
théorème suivant dont on trouvera la démonstration dans \cite{h}.

\begin{theo}[L'argument du petit objet]\label{fact}
Soient $\mathcal{C}$ une catégorie cocomplète et I un ensemble de morphismes
de $\mathcal{C}$ permettant l'argument du petit objet. Alors tout morphisme de
$\mathcal{C}$ se factorise de manière fonctorielle en une I-cofibration suivie
d'un morphisme I-injectif.
\end{theo}

L'un des corollaires de ce théorème qui nous intéresse le plus est celui qui
nous donne la caractérisation des I-cofibrations en termes de colimites
séquentielles transfinies de sommes amalgamées de flèches de I lorsque I
est assez petit. Mais en fait, les I-cofibrations ne sont pas toujours de
telles colimites même si I est assez petit, car en général ce sont des
rétracts de telles colimites. Aussi vais-je tout d'abord rappeler la notion de
rétract d'un morphisme avant d'énoncer le corollaire sur l'engendrement.

\begin{defin}
Soient $f:A\rightarrow B$ et $g:C\rightarrow D$ deux morphismes, on dit que $f$
est un rétract de $g$ s'il existe un tel diagramme commutatif :
\begin{diagram}
A & \rTo^{i} & C & \rTo^{r} & A \\
\dTo^{f} & & \dTo^{g} & & \dTo^{f} \\
B & \rTo_{j} & D & \rTo_{s} & B \\
\end{diagram}
avec $r\circ i=Id_A$ et $s\circ j=Id_B$
\end{defin}
\index{retract}

\begin{cor}\label{Icof}
Soient $\mathcal{C}$ une catégorie cocomplète et I un ensemble de morphismes
de $\mathcal{C}$ permettant l'argument du petit objet. Alors les I-cofibrations
sont exactement les rétracts de colimites séquentielles transfinies de sommes
amalgamées de morphismes de I.
\end{cor}
  
Grâce à ce corollaire nous savons que les morphismes engendrés par un bon
ensemble I, ce que l'on appelle des I-cofibrations, s'expriment en termes de
colimites séquentielles transfinies de sommes amalgamées d'éléments de
I. De plus, si l'on regarde la démonstration du théorème~\ref{fact}
donnée dans \cite{h}, on s'aperçoît que pour construire la
factorisation d'une flèche quelconque $f:A\rightarrow B$ en une I-cofibration suivi d'un
morphisme I-injectif, on va en fait rajouter à $A$ toutes les flèches de I
rendant commutatif le diagramme suivant et répéter ce procédé un nombre
transfini de fois.
\begin{diagram}
X & \rTo & A\\
\dTo_{\in I} & & \dTo_f \\
Y & \rTo & B\\
\end{diagram}
On remarque donc que cette construction est encore du type colimite
séquen\-tielle transfinie de sommes amalgamées d'éléments de I. C'est en
outre cette méthode que nous allons utiliser pour catégoriser nos
\precats~une fois explicité l'ensemble I de flèches pour lequel nos
\cats~sont des objets I-injectifs. C'est la raison pour laquelle nous allons
consacrer la prochaine partie à l'étude de ces procédés que nous
appellerons plans d'addition de cellules, du fait de leur similarité d'avec la
construction des complexes cellulaires.

\newpage

\section{Plan d'addition de cellules}

Dans cette partie, nous allons donner un formalisme permettant d'étudier les
procédés d'addition transfinie de flèches. Nous commençons d'abord par
définir les plans simples ou plans à un pas avant de s'attaquer aux plans
finis et infinis.

\begin{defin}\index{plan d'addition de cellules!à un pas}
Soient $\mathcal{C}$ une catégorie cocomplète et $\Phi$ un ensemble de
morphismes de $\mathcal{C}$. Soit $X$ un objet de $\mathcal{C}$. Un plan $P$ 
d'addition de cellules à un pas pour $X$ (ou plus
simplement plan simple pour $X$) est un ensemble de couples $(D,\alpha_D)$ où
$\alpha_D$ est un cardinal et $D=(d,\phi_D)$ est un diagramme comme suit :
\begin{diagram}
S_{\phi_D} & \rTo^{d} & X\\
\dTo^{\phi_D} & & \\
B_{\phi_D} & & \\
\end{diagram}
où $\phi_D$ est un morphisme de $\Phi$.
On note $P(X)$ le résultat du plan d'addition, défini de la manière
suivante :
$$P(X)=X\coprod_{\coprod_{D\in P} \coprod_{\alpha_D} S_{\phi_D}} \Bigg(
\coprod_{D\in P} \coprod_{\alpha_D} B_{\phi_D}\Bigg) $$
\end{defin}

Donnons un premier exemple de plan d'addition de cellules simple qui sera très
utile pour la suite.

\begin{defin}\index{$e_{\Psi,\lambda}$}
Soient $\mathcal{C}$ une catégorie cocomplète et $\Phi$ un ensemble de
morphismes de $\mathcal{C}$. Soit $\lambda$ un cardinal (transfini ou pas). On définit alors
$e_{\Phi,\lambda}$ le plan d'addition de cellules à un
pas qui, pour un objet $X$ quelconque de $\mathcal{C}$, est constitué par
l'ensemble des couples $(D,\lambda)$ où $D=(d,\phi_D)$ avec $\phi_D$
appartient à $\Phi$ et $d$ ayant pour but $X$. 
\end{defin}

Une fois les plans simples d'addition de cellules définis, on peut facilement
définir les plans d'addition de cellules finis comme composées des plans
simples et les plans d'addition de cellules infinis comme colimites de plans
simples. En voici les définitions précises.

\begin{defin}\index{plan d'addition de cellules!fini}
Soient $\mathcal{C}$ une catégorie cocomplète et $\Phi$ un ensemble de
morphismes de $\mathcal{C}$. Soit $X$ un objet de $\mathcal{C}$. 
Un plan $P_.$ à deux pas sur $X$ est une suite
$(P_i)_{1\leq i\leq 2}$ tel que $P_1$ est un plan simple sur $X$ et $P_2$ un
plan simple sur $P_1(X)$. Le résultat $P_.(X)$ n'est autre $P_2(P_1(X))$. On
notera parfois ce plan $P_2\circ P_1$, ce qu'on appellera composition (naïve) de plans simples.\\
Par récurrence, on définit de la même manière les plans finis comme
composition d'un nombre fini de plans simples.
\end{defin}

\begin{defin}\index{plan d'addition de cellules!infini}
Soient $\mathcal{C}$ une catégorie cocomplète et $\Phi$ un ensemble de
morphismes de $\mathcal{C}$. Soient $X$ un objet de $\mathcal{C}$ et $\lambda$
un cardinal transfini. 
Un plan $P_.$ infini d'addition de cellules sur $X$ est une suite des plans simples $(P_{\beta})_{\beta<\lambda}$ telle qu'il existe une $\lambda$-séquence $(X_{\beta})_{\beta\leq\lambda}$ vérifiant les propriétés suivantes :
\item - $X_0$ n'est autre que $X$,
\item - pour tout $\beta<\lambda$, $X_{\beta+1}$ est le résultat
du plan simple $P_{\beta}$ sur $X_{\beta}$,
\item - pour tout $\gamma\leq\lambda$ ordinal limite, $X_{\gamma}$ est la colimite des $X_{\beta}$ pour $\beta<\gamma$.\\
\\
 Le résultat $P_.(X)$ du plan infini $P_.$ est soit
$P_{\lambda-1}(X_{\lambda-1})$, si $\lambda$ n'est pas un cardinal limite, soit la
colimite des $X_{\beta}$, pour $\beta<\lambda$ si $\lambda$ est un cardinal
limite.\\
En outre, on dit que $\lambda$ est la longueur du plan $P$. 
\end{defin}

Pour simplifier les écritures, nous allons prendre les conventions
suivantes : soit $A$ un objet de $\mathcal{C}$ et $P$ un plan d'addition de
cellules sur $A$ de longueur $\lambda$. On pose $P_0(A)=A$ puis, pour tout
ordinal $\beta\leq\lambda$, on note $P_{\beta}(A)$ le résultat du plan
$(P_{\alpha})_{\alpha<\beta}$. Ainsi les diagrammes de $P_{\beta}$ sont des
diagrammes sur $P_{\beta}(A)$.\\

Lors d'un plan fini ou infini d'addition de cellules, il se peut que l'on rajoute la même cellule à plusieurs étapes du plan. Or ceci obscurcit la compréhension de l'opération effectuée par le plan d'addition de cellules. Certains plans d'addition de cellules cependant n'ont pas ce défaut, i.e. les cellules $D$ ne seront ajoutées qu'à un seul moment donné de la
construction et exactement $\alpha_D$ fois dans toute la construction,
 c'est ce que l'on appellera des plans rationnels.

\begin{defin}\index{plan d'addition de cellules!rationnel}
Soient $\mathcal{C}$ une catégorie cocomplète et $\Phi$ un ensemble de
morphismes de $\mathcal{C}$. Soient $X$ un objet de $\mathcal{C}$ et $P_.$ un plan d'addition de cellules sur $X$ de longueur un certain cardinal $\lambda$. On dit que $P_.$ est un plan d'addition de cellules rationnel si, pour tout ordinal $\beta\leq\lambda$ et pour tout diagramme $D$ de $P_{\beta}$, l'extension de $D$ à $P(X)$ n'est pas l'extension à $P(X)$ d'un diagramme de $P_{\alpha}$ avec $\alpha<\beta$. 
\end{defin}

On a vu dans la définition de plan d'addition fini de cellules la composition
des plans simples. Mais on peut composer en fait n'importe quel plan.

\begin{defin}\index{composition de plans d'addition de cellules!naïve}
Soient $\mathcal{C}$ une catégorie cocomplète et $\Phi$ un ensemble de
morphismes de $\mathcal{C}$. Soient $X$ un objet de $\mathcal{C}$, $P_.$ un plan
fini ou infini sur $X$ et $Q_.$ un plan fini ou infini sur $P_.(X)$. La
composition naïve des plans $P_.$ et $Q_.$ est le plan noté $Q_.\circ P_.$
obtenu par concaténation des suites de plans simples de $P_.$ et $Q_.$.
\end{defin}

Cette composition a de bonnes propriétés : elle est associative et admet
comme identités les plans vides ou ceux dont tous les coefficients
$\alpha_D$ sont nuls. Toutefois, on l'appelle "naïve" car lorsque l'on
applique le plan $Q_.$ à $P_.(X)$ on pourra ré-additionner des cellules
déjà additionnées à $X$ par $P_.$. Aussi afin de mieux contrôler
l'addition des cellules lors d'une composition de plans, on va définir une
nouvelle composition dite rationnelle, qui nous donnera une construction simple de plans
rationnels, notion dans laquelle le phénomène décrit ci-dessus n'aura
plus lieu. 

\begin{defin}\index{composition de plans d'addition de cellules!rationnelle}
Soient $\mathcal{C}$ une catégorie cocomplète et $\Phi$ un ensemble de
morphismes de $\mathcal{C}$. Soient $X$ un objet de $\mathcal{C}$, $P_.$ un plan
fini ou infini sur $X$ et $Q_.$ un plan simple sur $P_.(X)$. Notons
$(X_{\beta})_{\beta\leq\lambda}$ la suite correspondant à $P_.$. La
composition rationnelle des plans $P_.$ et $Q_.$ est le plan noté $Q_.*P_.$
obtenu comme composition naïve du plan $P_.$ avec le plan $Q'_.$ défini
ainsi : un diagramme $D$ de $Q_.$ appartient à $Q'_.$ si et seulement si $d$
ne se factorise pas par l'un des $X_{\beta}$, pour $\beta<\lambda$.
\end{defin}

L'un des intérêts de cette composition rationnelle est que le plan $Q_.*P_.$ est moralement un plan rationnel. Toutefois ceci n'est pas vrai en général. Considérons par exemple pour $X$ un ensemble simplicial construit avec trois segments $f,g,h$ tel que $g$ et $h$ soient précomposables par $f$. Considérons le plan simple $P_.$ sur $X$ consistant à rajouter le triangle de faces principales $f,g$ et considérons le plan simple $Q_.$ sur $P(X)$ consistant à rajouter le triangle de faces principales $f,h$ et à identifier $g$ et $h$ en un même segment noté $i$. La composition naïve des plans $P_.$ et $Q_.$ n'est autre que leur composition rationnelle et le résultat de cette composition rationnelle $(Q*P)(X)$ consiste en deux triangles collés le long de leurs faces principales $f,i$. Ainsi au final une même cellule, ici le triangle de face principale $f,i$, a été additionnée deux fois mais à deux moments différents de la construction. Aussi ce plan $Q_.*P_.$ issu d'une composition rationnelle n'est pas un plan rationnel. Le principal empêchement vient du fait que l'une des cellules, celle qui identifie les segments, n'est pas un monomorphisme, ce qui amène au final à identifier deux diagrammes qui au départ étaient différents.\\
Nous allons donc montrer que sous certaines hypothèses sur $\Phi$ concernant les monomorphismes, les plans d'addition de cellules obtenus à partir de compositions rationnelles sont bien des plans d'addition rationnels au sens défini précédemment.

\begin{prop}\label{comprat}\index{plan d'addition de cellules!à compositions rationnelles}
Soient $\mathcal{C}$ une catégorie cocomplète et $\Phi$ un ensemble de
morphismes de $\mathcal{C}$ tel que les $\Phi$-cofibrations soient des monomorphismes. Soient $X$ un objet de $\mathcal{C}$ et $P_.$ un plan
fini ou infini sur $X$ dont toutes les compositions le composant sont des compositions rationnelles. Un tel plan sera appelé plan à compositions rationnelles.\\
Alors $P_.$ est un plan d'addition rationnel.
\end{prop}
{\it Preuve :}\\
Pour montrer que $P_.$ est un plan d'addition de cellules rationnel sur $X$, nous allons montrer par l'absurde qu'il vérifie la définition de plan rationnel. Notons $\lambda$ la longueur du plan $P_.$. Soit $\beta\leq\lambda$ un ordinal et soit $D$ un diagramme de $P_{\beta}$. Supposons par l'absurde que l'extension de ce diagramme à $P(X)$ est égale à l'extension à $P(X)$ d'un diagramme $D'$ de $P_{\alpha}$ pour un certain ordinal $\alpha<\beta$. On remarque que le morphisme naturel de $P_{\beta}(X)$ vers $P(X)$ est une colimite séquentielle finie ou transfinie de sommes amalgamées le long d'un morphisme de flèches de $\Phi$. Par le lemme~\ref{Ist}, il vient que les colimites séquentielles de sommes amalgamées le long d'un morphisme d'éléments de $\Phi$ sont des $\Phi$-cofibrations, donc des monomorphismes par hypothèse sur $\Phi$. On en déduit que le diagramme $D$ sur $P_{\beta}(X)$ n'est autre que l'extension du diagramme $D'$ à $P_{\beta}(X)$. Or par hypothèse sur $P_.$, la composition du plan simple $P_{\beta}$ sur $P_{\beta}(X)$ avec le plan $(P_{\gamma})_{\gamma<\beta}$ sur $X$ est une composition rationnelle, or par définition même de composition rationnelle ceci contredit le fait que le diagramme $D$ sur $P_{\beta}(X)$ est égal à l'extension à $P_{\beta}(X)$ du diagramme $D'$ de $P_{\alpha}$ avec $\alpha<\beta$. On a donc bien montré par l'absurde que $P_.$ est un plan rationnel d'addition de cellules.\\
CQFD. 

\begin{ex}
Soient $\mathcal{C}$ une catégorie cocomplète et $\Phi$ un ensemble de
morphismes de $\mathcal{C}$ tel que les $\Phi$-cofibrations soient des monomorphismes.\\
Considérons le plan d'addition de cellules $e_{\Phi,\lambda}$, si l'on
compose simplement $e_{\Phi,\lambda}$ avec lui-même, on va remettre une deuxième fois toutes
les cellules ajoutées la première fois. Si en revanche, on compose
rationnellement $e_{\Phi,\lambda}$ avec lui-même, on ne mettra à la seconde étape que les
cellules qui n'ont pas déjà été mises à la première.
\end{ex}

Dans ces définitions, la notion de plan dépend de l'objet sur lequel le
plan doit agir. Toutefois on peut définir des plans pour tous les objets de la
catégorie, auquel cas on peut se demander si cette collection de plans est
fonctorielle. Nous allons donc définir la notion de plan fonctoriel.

\begin{defin}\index{plan d'addition de cellules!fonctoriel}
Soient $\mathcal{C}$ une catégorie cocomplète et $\Phi$ un ensemble de
morphismes de $\mathcal{C}$.
Soit $P=(P_A)_{A\in\mathcal{C}}$ une collection de plans d'addition de cellules,
on dit que $P$ est un plan d'addition de cellules fonctoriel si $P$ induit un
foncteur, encore noté $P$, de la catégorie $\mathcal{C}$ dans elle-même tel
que, pour tout objet $A$ de $\mathcal{C}$, $P(A)$ est le résultat $P_A(A)$ du
plan d'addition $P_A$.
\end{defin} 

La fonctorialité des plans d'addition de cellules se comporte bien
vis-à-vis des compositions de plans et peut se reconnaître pour les plans
non simples par la fonctorialité des plans simples qui les composent.

\begin{lem}\label{plsfonc}
Soient $\mathcal{C}$ une catégorie cocomplète et $\Phi$ un ensemble de
morphismes de $\mathcal{C}$. On a les résultats suivants :
\item - si $P$ et $Q$ sont des plans d'addition de cellules fonctoriels alors leurs
composées naïve et rationnelle sont fonctorielles,
\item - si $P=(P_{\beta})_{\beta<\lambda}$ est un plan d'addition transfinie de
cellules dont tous les plans simples sont fonctoriels, alors $P$ est un plan
d'addition de cellules fonctoriel.
\end{lem}
{\it Preuve :} cela découle directement du fait que les foncteurs sont
composables et possèdent les colimites. CQFD.\\

Il est intéressant de donner une condition suffisante simple pour qu'une
collection de plans d'addition de cellules constitue un plan d'addition de
cellules fonctoriel. Pour cette condition, nous utiliserons la notion
d'extension de plans d'addition de cellules le long d'un morphisme que nous
allons tout d'abord définir.

\begin{defin}\index{plan d'addition de cellules!étendu}
Soient $\mathcal{C}$ une catégorie cocomplète et $\Phi$ un ensemble de
morphismes de $\mathcal{C}$.
Soient $f:A\rightarrow B$ un morphisme de $\mathcal{C}$ et $P$ un plan
d'addition de cellules sur $A$. On appelle plan d'addition de cellules étendu
par $f$ le plan d'addition de cellules sur $B$ constitué par l'extension par
$f$ des diagrammes du plan $P$ pris avec la même cardinalité que dans $P$.
\end{defin}

\begin{lem}\label{plfcrit}
Soient $\mathcal{C}$ une catégorie cocomplète et $\Phi$ un ensemble de
morphismes de $\mathcal{C}$.
Soit $P=(P_A)_{A\in\mathcal{C}}$ une collection de plans d'addition de cellules
de même longueur.
Si pour tout morphisme $f:A\rightarrow B$ de $\mathcal{C}$, l'extension par $f$ du plan
d'addition de cellules $P_A$ sur $A$ est un sous-ensemble du plan d'addition de
cellules $P_B$ sur $B$, alors $P$ est un plan d'addition fonctoriel.
\end{lem} 
{\it Preuve :}\\
D'après le lemme~\ref{plsfonc}, pour montrer qu'une collection de plans non
simples de même longueur est fonctorielle,
il suffit de montrer que chaque collection de plans simples qui la compose est
fonctorielle. Donc on va montrer le lemme pour une collection $P=(P_A)_{A\in\mathcal{C}}$ de
plans simples. Soit $f:A\rightarrow B$ un morphisme de $\mathcal{C}$. Le
morphisme $f$ induit une fonction non nécessairement injective, encore notée
$f$, de l'ensemble
$D_A$ des diagrammes du plan $P_A$ vers l'ensemble $D_B$ des diagrammes du plan
$P_B$ obtenue en étendant les diagrammes sur $A$ par $f$. Comme l'extension de
$P_A$ par $f$ est un sous-ensemble de $P_B$, ceci signifie que les diagrammes de
$B$ provenant par $f$ de diagrammes de $A$ sont pris avec la même
cardinalité que leurs antécédents. On en déduit le
carré commutatif suivant :
\begin{diagram}
\coprod_{D\in P_A} \coprod_{\alpha_D} X & \rTo & \coprod_{D'\in P_B}
\coprod_{\alpha_D'} X \\
\dTo^{\coprod_{D\in P_A} \coprod_{\alpha_D} \phi_D} & & \dTo_{\coprod_{D'\in P_B}
\coprod_{\alpha_D'} \phi_{D'}} \\
\coprod_{D\in P_A} \coprod_{\alpha_D} Y & \rTo & \coprod_{D'\in P_B}
\coprod_{\alpha_D'} Y \\ 
\end{diagram}
dans lequel on envoie par l'identité les sources et buts d'une flèche
$\phi_D$ d'un diagramme $D\in D_A$ sur ceux de la flèche $\phi_{f(D)}$ de
l'extension $f(D)\in D_B$ du diagramme $D$ par $f$. On déduit de ce diagramme
le cube commutatif suivant :
\begin{diagram}
\coprod_{D\in P_A} \coprod_{\alpha_D} X & \rTo & & \coprod_{D'\in P_B}
\coprod_{\alpha_{D'}} X & & \\
\dTo~{\coprod_{D\in P_A} \coprod_{\alpha_D} \phi_D}& \rdTo &  & \dTo~{\coprod_{D'\in P_B}
\coprod_{\alpha_{D'}} \phi_{D'}} & \rdTo & \\
 & & A &  & \rTo~{f}& B \\
\coprod_{D\in P_A} \coprod_{\alpha_D} Y & \rTo & & \coprod_{D'\in P_B}
\coprod_{\alpha_{D'}} Y & & \\ 
 & \rdTo & \dTo &  & \rdTo & \dTo\\
 & & P_A(A) & \rTo_{\exists !} & & P_B(B)\\ 
\end{diagram}
où les faces verticales gauche et droite sont des diagrammes cocartésien. Il
est facile de voir que si dans le cube on prend pour $f$ l'identité de $A$, le
morphisme du bas est l'identité. De plus, si l'on a deux morphismes
composables $f$ et $g$, on obtient deux cubes commutatifs qui se juxtaposent rendant le
grand cube total commutatif, ce qui entraîne que la composée des
morphismes induits par $f$ et $g$ est bien le morphisme induit par la composée
$f\circ g$. Posons alors, pour tout objet $A$ de $\mathcal{C}$, $P(A)=P_A(A)$
et, pour tout morphisme $f$, notons $P(f)$ le morphisme induit par le cube.
On a donc que $P$ ainsi défini est bien un foncteur de la catégorie $\mathcal{C}$
dans elle-même, ce qui montre le résultat.\\
CQFD.

\begin{lem}\label{ephif}
Soient $\mathcal{C}$ une catégorie cocomplète, $\Phi$ un ensemble de
morphismes de $\mathcal{C}$ et $\lambda$ un cardinal.
$e_{\Phi,\lambda}$ est un plan d'addition de cellules fonctoriel.\\
De plus tout plan fini ou transfini d'addition de cellules dont tous les plans
simples sont du type $e_{\Phi,\lambda}$ est fonctoriel.
\end{lem}
{\it Preuve :}
Pour montrer le premier résultat il suffit d'appliquer le critère du
lemme~\ref{plfcrit}. Soit $f:A\rightarrow B$ un morphisme de $\mathcal{C}$,
l'extension de $e_{\Phi,\lambda,A}$ n'est autre que le plan simple d'addition de
cellules constitué par les couples $(D,\lambda)$ où $D$ décrit les
extensions par $f$ de diagrammes de $\Phi$ pour $A$. Or ces couples sont bien des
couples du plan simple $e_{\Phi,\lambda,B}$, ce qui montre que
$e_{\Phi,\lambda}$ vérifie la condition suffisante pour être un plan
d'addition de cellules fonctoriel.\\
Le second résultat résulte du premier auquel on applique le
lemme~\ref{plsfonc}.\\
CQFD.\\

Le formalisme du rajout de cellules ayant été précisé, nous sommes
maintenant en mesure de donner un procédé qui permettra de catégoriser les
\precats.

\newpage

\section{Plans d'addition de cellules I-injectivants}

On peut en effet reformuler le problème de la catégorisation en termes de
plans d'addition de cellules et d'objets I-injectifs. Supposons
en effet que nous ayons exhibé une famille I de morphismes dont les objets
I-injectifs sont les \cats. Le problème devient alors de rendre I-injectif un
objet qui ne l'est pas. L'idée la plus répandue pour faire cela est de
rajouter les flèches de I jusqu'à rendre l'objet I-injectif, c'est-à-dire
d'appliquer à l'objet un plan d'addition de cellules. Ainsi l'équivalent en
termes de plan d'addition de cellules et d'objets I-injectifs du problème de
la catégorisation est de trouver un plan d'addition de cellules $P$ fonctoriel
vérifiant :\\
- i) tout morphisme $f:A\rightarrow B$ dont le but est I-injectif se factorise à
travers $P(A)$,\\
- ii) pour tout objet $A$, $P(A)$ est I-injectif,\\
- iii) si l'on possède une notion d'\eq~d'objets I-injectifs, pour tout objet $A$,
l'image par $P$ du morphisme naturel $A\rightarrow P(A)$ est une \eq~d'objets
I-injectifs.\\
Tout plan fonctoriel d'addition de cellules résolvant ce problème sera
appelé plan I-injectivant ou I-injectivisation.\\

Comme on l'a dit plus haut, l'idée que nous allons suivre pour résoudre ce
problème d'I-injectivisation est de construire un plan d'addition de flèches
de I. Or il est assez facile de voir, avec les définitions d'objets I-injectifs et de
I-cofibrations, que tout plan d'addition de flèches de I est une I-cofibration
et donc, par \prg~par rapport aux objets I-injectifs, vérifie le i). C'est ce
que nous donne le lemme suivant.

\begin{lem}\label{plIcof}
Soient $\mathcal{C}$ une catégorie cocomplète et I un ensemble de morphismes
de $\mathcal{C}$. Soient $A$ un objet de $\mathcal{C}$ et $P$ un plan d'addition
de cellules de I sur $A$. Alors le morphisme naturel $A\rightarrow P(A)$ est une
I-cofibration. En particulier pour tout morphisme $f:A\rightarrow B$ vers un
objet I-injectif $B$, il existe un morphisme de $P(A)$ vers $B$ qui
précomposé avec le morphisme naturel de $A$ vers $P(A)$ redonne $f$.
\begin{diagram}
A & \rTo^{f} & B \\
\dTo & \ruTo_{\exists} \\
P(A) & & \\
\end{diagram}
\end{lem}
{\it Preuve :}\\
Tout d'abord il convient de remarquer que le morphisme naturel $A\rightarrow
P(A)$ est une colimite séquentielle transfinie de sommes amalgamées de
flèches de I. Or les flèches de I sont des I-cofibrations et les
I-cofibrations sont stables par somme amalgamée le long d'un morphisme et par
colimite séquentielle transfinie d'après le lemme~\ref{Ist}. Donc le
morphisme naturel $A\rightarrow P(A)$ est bien une I-cofibration. En outre par
définition de I-cofibration, $A\rightarrow P(A)$ possède la \prg~par rapport
aux morphismes I-injectifs et donc en particulier par rapport aux objets
I-injectifs.\\
CQFD.\\

Ce lemme nous a donc permis d'avoir une infinité de solutions pour le
problème i) de l'I-injectivisation. Il suffit que le plan d'addition de
cellules rajoute des flèches de I à la source d'un morphisme vers un objet
I-injectif pour factoriser ce morphisme à travers le résultat du plan
appliqué à la source. Toutefois on peut aussi bien rajouter à chaque
étape du plan une partie quelconque de I ou bien avoir une partition de I 
et un plan qui à chaque étape rajoute les flèches d'un
certain nombre d'éléments de la partition de I. Nous allons donc donner une
définition générale qui permet de modéliser de tels plans.

\begin{defin}\index{plan d'addition de cellules!suivant un ordre}
Soient $\mathcal{C}$ une catégorie cocomplète et $\Phi$ un ensemble de
morphismes de $\mathcal{C}$. Soit $(\Phi_k)_{k\in K}$ une partition de $\Phi$. Donnons nous une fonction $t$ d'un cardinal $T$ quelconque
vers l'ensemble des parties de $K$. On appelle plan d'addition de cellules
suivant l'ordre $t$ tout plan $P$ de longueur $T$ tel que, pour tout ordinal
$\tau<T$, les diagrammes apparaissant dans $P_{\tau}$ ne sont obtenus qu'à
partir des flèches des $\Phi_k$ avec $k$ décrivant $t(\tau)$.\\
On notera $t^{-1}(k)$ pour $k\in K$ l'ensemble des ordinaux $\tau<T$ pour
lesquels $t(\tau)$ contient $k$.
\end{defin}

Pour résoudre le problème ii) de l'I-injectivisation, l'idée est de
rajouter autant de fois que nécessaire les flèches de I pour rendre les
objets quelconques I-injectifs. On se demande alors quelles conditions doit
vérifier un plan d'addition de flèches de I pour que son résultat soit
I-injectif. Comme les plans ont chacun une certaine longueur, si l'ensemble I
est petit par rapport à cette longueur, toutes flèches d'une source d'un
morphisme de I vers le résultat du plan se factorisera à travers l'une des
étapes du plan. Deux cas se présente alors : soit le diagramme à cette
étape possède déjà un relèvement, et dans ce cas c'est bon, soit le
diagramme à cette étape n'en possède pas. Or si l'on veut que le diagramme
se relève pour le résultat du plan, il faudra qu'il se relève pour une
étape ultérieure du plan. Ces considérations nous amènent au lemme
suivant qui nous donne une condition suffisante pour qu'un plan d'addition
suivant un ordre ait son résultat I-injectif.

\begin{lem}
Soient $\mathcal{C}$ une catégorie cocomplète et $\Phi$ un ensemble de
morphismes de $\mathcal{C}$. Soient $(\Phi_k)_{k\in K}$ une partition de $\Phi$ et $\alpha$ un cardinal régulier pour lequel les sources des morphismes de $\Phi$ sont $\alpha$-petites.\\ 
Soit $A$ un objet de $\mathcal{C}$. Donnons nous une fonction $t$ d'un cardinal
régulier $T$ strictement supérieur à $\alpha$
vers l'ensemble des parties de $K$ et un plan d'addition de cellules $P$ sur $A$
suivant l'ordre $t$. Si $P$ vérifie la propriété suivante :\\
pour tout ordinal $\tau<T$, pour tout $k\in K$ et pour tout diagramme $D$ sur
$P_{\tau}(A)$ fabriqué avec une flèche de $\Phi_k$ qui ne se relève pas dans $P_{\tau}(A)$, il existe $\tau'\in t^{-1}(k)$
tel que $\tau'\geq\tau$ et que le diagramme $D$ étendu à $P_{\tau'}(A)$
appartienne au plan $P_{\tau'}$,\\
alors $P(A)$ est $\Phi$-injectif.
\end{lem}
{\it Preuve :}\\
Soit $P$ un plan sur $A$ vérifiant la condition du lemme. Montrons que $P(A)$
est $\Phi$-injectif. Considérons alors le diagramme suivant :
\begin{diagram}
X & \rTo^g & P(A) \\
\dTo^{\phi} & & \\
Y & & \\
\end{diagram}
Comme par hypothèse la source $X$ de $\phi$ est $\alpha$-petite et que $A\rightarrow P(A)$ est une colimite
séquentielle transfinie de longueur $T$ cardinal régulier strictement
supérieur à $\alpha$, il existe un ordinal $\tau<T$ tel que $g$ se factorise
par $P_{\tau}(A)$. Soit le diagramme se relève dans $P_{\tau}(A)$ et, par extension à $P(A)$, on obtient que le diagramme ci-dessus se relève. Soit le diagramme ne se relève pas dans $P_{\tau}(A)$, et dans ce cas, par la propriété de $P$, il vient qu'il existe
$\tau'\geq\tau$ tel que $t(\tau')$ contienne l'indice de la partie de $\Phi$ à
laquelle appartient $\phi$ et que la restriction du diagramme à $P_{\tau'}(A)$
appartienne au plan $P_{\tau'}$. Ainsi $P_{\tau'+1}(A)$ se relève par rapport
à $\phi$, et par extension à $P(A)$,on obtient que $P(A)$ se relève aussi, ce qui montre que le
diagramme considéré se relève.\\
CQFD.

\begin{cor}\label{plIinj}
Soient $\mathcal{C}$ une catégorie cocomplète et $\Phi$ un ensemble de
morphismes de $\mathcal{C}$. Soient $(\Phi_k)_{k\in K}$ une partition de $\Phi$ et $\alpha$ un cardinal régulier pour lequel les sources des flèches de $\Phi$ sont $\alpha$-petites.\\ 
Soit $A$ un objet de $\mathcal{C}$. Donnons nous une fonction $t$ d'un cardinal
régulier $T$ strictement supérieur à $\alpha$
vers l'ensemble des parties de $K$ et un plan d'addition de cellules $P$ sur $A$
suivant l'ordre $t$.\\
Supposons que $t$ vérifie que pour tout $\tau<T$ et pour tout $k\in t(\tau)$
il existe $\tau'\geq\tau$ tel que $\tau'$ appartienne à $t^{-1}(k)$.
Si $P$ vérifie la propriété suivante :\\
pour tout ordinal $\tau<T$, le plan $P_{\tau}$ contient tous les diagrammes
obtenus à partir des morphismes de $\Phi_{k}$, pour $k$ décrivant $t(\tau)$,
qui ne se relèvent pas dans $P_{\tau}(A)$,
alors $P(A)$ est $\Phi$-injectif.
\end{cor}
{\it Preuve :}\\
Montrons que sous les conditions du corollaire le plan $P$ vérifie la
condition du lemme précédent. Soit donc un ordinal $\tau<T$ et un diagramme
$D$ sur $P_{\tau}(A)$ obtenu à partir d'un morphisme d'un $\Phi_k$ avec $k\in t(\tau)$ qui ne se relève pas dans $P_{\tau}(A)$. Par
propriété de $t$, il existe un ordinal $\tau'\geq\tau$ tel que $k$
appartienne à $t(\tau')$. Prenons le plus petit $\tau'$ vérifiant cela. En outre par propriété de $P$, il vient que
$P_{\tau'}$ contient tous les diagrammes obtenus avec $\Phi_k$ (car $k\in
t(\tau')$) ne se relevant pas dans $P_{\tau'}(A)$, c'est en particulier le cas
pour $D$ étendu à $P_{\tau'}(A)$, par minimalité de $\tau'$. Ainsi $P$ vérifie la condition du lemme
et donc $P(A)$ est $\Phi$-injectif.\\
CQFD.\\ 

Donnons deux exemples de plans d'addition de cellules qui vérifient les
propriétés i) et ii) de l'I-injectivisation.

\begin{lem}\label{Ephi1}
Soient $\mathcal{C}$ une catégorie cocomplète et $\Phi$ un ensemble de
morphismes de $\mathcal{C}$ dont les sources sont $\alpha$-petites pour un certain cardinal régulier
$\alpha$. Soient $\alpha'$ le plus petit cardinal régulier supérieur à
$\alpha$ et $\lambda$ un cardinal quelconque.\\
Notons $E_{\Phi,\lambda}$ le plan d'addition de cellules fonctoriel de longueur
$\alpha'$ dont tous les plans simples sont $e_{\Phi,\lambda}$ et notons
$E'_{\Phi,\lambda}$ le plan d'addition de cellules fonctoriel de
longueur $\alpha'$ obtenu à partir des compositions rationnelles et dont tous les plans simples sont $e_{\Phi,\lambda}$.\\
Alors, pour tout objet $A$ de $\mathcal{A}$, $E_{\Phi,\lambda}(A)$ et
$E'_{\Phi,\lambda}(A)$ sont $\Phi$-injectifs et tout morphisme $f:A\rightarrow
B$ dont le but est $\Phi$-injectif se factorise à travers le morphisme naturel
$A\rightarrow E_{\Phi,\lambda}(A)$ et aussi à travers $A\rightarrow
E'_{\Phi,\lambda}(A)$.
\end{lem}
{\it Preuve :}\\
Tout d'abord on remarque que comme $E_{\Phi,\lambda}$ et $E'_{\Phi,\lambda}$
sont des plans d'addition de morphismes de $\Phi$, par le lemme~\ref{plIcof}, les
morphismes naturels $A\rightarrow E_{\Phi,\lambda}(A)$ et $A\rightarrow E'_{\Phi,\lambda}(A)$ sont des $\Phi$-cofibrations qui donc ont la
\prg~par rapport aux objets $\Phi$-injectifs. Ceci prouve la seconde partie du
lemme.\\

Pour montrer que les résultats des plans $E_{\Phi,\lambda}$ et $E'_{\Phi,\lambda}$ sont $\Phi$-injectifs, il
suffit de montrer que ces plans vérifient le corollaire~\ref{plIinj}. Tout
d'abord on remarque que dans notre cas il n'y a pas de partition de $\Phi$, ce
qui revient à dire que la partition n'a qu'un élément. Du coup la fonction
$t$ correspondant à nos plans n'est autre que la fonction constante du
cardinal $\alpha'$ dans l'ensemble des parties du singleton à valeur le
singleton. Ainsi $t$ vérifie de manière triviale la condition du
corollaire~\ref{plIinj}. Il ne reste plus qu'à montrer que pour tout ordinal
$\beta<\alpha'$, les plans simples $(E_{\Phi,\lambda})_{\beta}$ et
$(E'_{\Phi,\lambda})_{\beta}$ contiennent tous les diagrammes obtenus avec des
morphismes de $\Phi$ qui ne se relèvent pas dans
$(E_{\Phi,\lambda})_{\beta}(A)$ et $(E'_{\Phi,\lambda})_{\beta}(A)$. Or les
plans simples $(E_{\Phi,\lambda})_{\beta}$ et $(E'_{\Phi,\lambda})_{\beta}$ sont
en fait les plans $e_{\Phi,\lambda}$ qui par définition contiennent tous les
diagrammes possibles obtenus à partir de $\Phi$. Si dans
$(E_{\Phi,\lambda})_{\beta}$, on a bien tous les diagrammes possibles sur
$(E_{\Phi,\lambda})_{\beta}(A)$, en revanche par rationalité des compositions de
$E'_{\Phi,\lambda}$, dans $(E'_{\Phi,\lambda})_{\beta}$
on ne garde que les diagrammes qui ne relèvent pas dans $(E'_{\Phi,\lambda})_{\beta}$. Dans les deux cas
la propriété du corollaire~\ref{plIinj} est vérifiée, ce qui nous donne
que les résultats des plans $E_{\Phi,\lambda}$ et $E'_{\Phi,\lambda}$ sont bien $\Phi$-injectifs.\\
CQFD.\\

On a donc exhibé deux plans d'addition de cellules qui vérifient les
propriétés i) et ii) des plans I-injectivants. Le problème consiste
maintenant à savoir s'ils vérifient ou non la propriété iii) concernant
les \eqs. En fait cette propriété iii) est la plus difficile à obtenir.
En nous inspirant des idées de \cite{s} pour résoudre ce problème, nous
allons chercher deux plans d'addition de cellules particuliers vérifiant i) et ii). Le premier aura une
propriété d'unicité de relèvement par rapport aux objets I-injectifs et
le second aura la propriété iii) mais uniquement pour les objets
I-injectifs. Nous allons donc construire le premier plan et
montrer que, si l'on suppose l'existence du second, le premier vérifie iii).

\newpage

\section{Plan I-injectivant marqué}

L'idée pour trouver un plan d'addition I-injectivant est de construire un plan
vérifiant les propriétés i) et ii) mais avec pour la propriété de
factorisation i) l'unicité de cette factorisation. Bien évidemment on ne
pourra pas obtenir une telle propriété avec un objet I-injectif quelconque
qui en général n'a pas la propriété de se relever de manière unique
par rapport aux I-cofibrations. Cependant la notion de marquage déjà vu
auparavant va nous permettre d'aboutir à l'unicité des relèvements. Afin
de mieux manipuler cette notion de marquage, nous allons donner une définition
d'objet partiellement I-marqué.

\begin{defin}\index{objet!partiellement I-marqué}
Soit $\mathcal{C}$ une catégorie cocomplète et $I$ un ensemble de morphismes
de $\mathcal{C}$.
Un objet partiellement I-marqué est un couple $(A,\mu)$, où
$A$ est un objet quelconque de $\mathcal{C}$ et $\mu$ une fonction d'un ensemble de diagrammes solides sur $A$ du type
suivant :
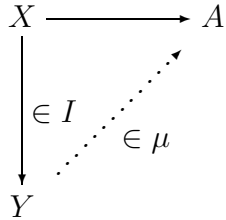
\begin{diagram}
X & \rTo & A \\
\dTo_{\in I} & \ruDotsto_{\in \mu} & \\
Y & & \\
\end{diagram}
qui associe à chaque diagramme solide de cet ensemble un relèvement, que l'on dira marqué. Par abus de langage, on dira d'un diagramme qu'il appartient à $\mu$ pour dire qu'il appartient au domaine de définition de la fonction $\mu$ et donc qu'on lui a choisi un relèvement.
\end{defin}

On remarque que si $A$ est totalement marquée (i.e. $\mu$ contient tous les
diagrammes des flèches de I vers $A$) alors $A$ a
bien évidemment la \prd~par rapport à toutes les flèches de I, ce qui fait
de $A$ un objet I-injectif pour lequel tous les
relèvements sont marqués, c'est ce qu'on va appeler un objet I-injectif marqué.
 
\begin{defin}\index{I-injectif!marqué}
Soit $\mathcal{C}$ une catégorie cocomplète et $I$ un ensemble de morphismes
de $\mathcal{C}$.
Un objet I-injectif marqué est un couple $(A,\mu)$, où
$A$ est un objet I-injectif et $\mu$ une fonction de l'ensemble des diagrammes solides de flèches de I
vers $A$ qui à chaque diagramme solide associe un relèvement, qui sera le relèvement
marqué.
\end{defin}

Les notions d'objet totalement marqué et d'objet I-injectif marqué sont
équivalentes, comme on l'a vu plus haut. Il nous reste donc à définir ce que
sont les morphismes d'objets partiellement marqués.

\begin{defin}\index{morphisme!d'objets partiellement marqués}
\index{morphisme!préservant le marquage}
Soit $\mathcal{C}$ une catégorie cocomplète et $I$ un ensemble de morphismes
de $\mathcal{C}$.
Un morphisme $f:A\rightarrow B$ entre les objets partiellement
marqués $(A,\mu)$ et $(B,\nu)$ est un morphisme d'objets partiellement marqués si, pour tout
diagramme de $\mu$, son prolongement par $f$ est un diagramme de $\nu$. On dit
alors que le morphisme préserve le marquage.\\
Un morphisme d'objets I-injectifs marqués n'est autre qu'un morphisme d'objets
partiellement I-marqués.
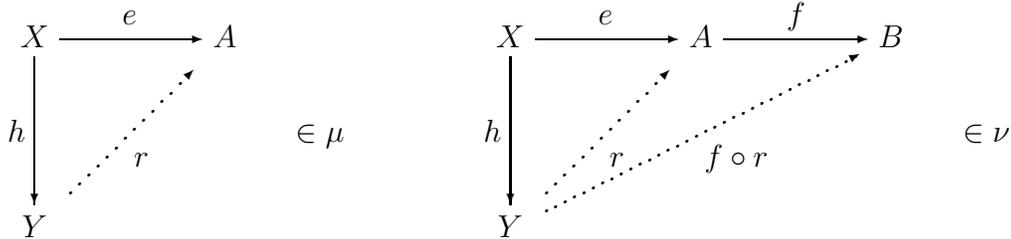
\begin{diagram}
X        & \rTo^{e}      & A &       & & X        & \rTo^{e} & A & \rTo^{f} & B & \\
\dTo^{h} & \ruDotsto_{r} &   &\in\mu & & \dTo^{h} &
\ruDotsto_{r} & &\ruDotsto(4,2)_{f\circ r} & & \in\nu \\
Y        &               &   &       & & Y        &          &   &       &  & \\
\end{diagram}
\end{defin}

On obtient ainsi deux catégories : $\mathcal{C}_m$, la
catégorie des objets partiellement marqués et sa sous-catégorie pleine
$\mathcal{C}_{I-inj m}$ des objets I-injectifs marqués. Tout d'abord, remarquons qu'un même
objet I-injectif peut avoir plusieurs marquages différents et ainsi peut donner lieu à
plusieurs objets I-injectifs marqués, il en va de même pour les objets
quelconques. Ainsi nos deux nouvelles catégories ne sont pas des sous-catégories de $\mathcal{C}$.
En revanche, elles sont toutes deux munies d'un foncteur Oubli vers
$\mathcal{C}$ qui est fidèle mais n'est pas plein car tout morphisme ne
préserve pas nécessairement le marquage. En outre on a vu que la notion
d'objets I-injectifs n'est pas stable par limite mais en revanche celle d'objets
I-injectifs marqués l'est.

\begin{lem}
Soit $\mathcal{C}$ une catégorie complète et $I$ un ensemble de morphismes
de $\mathcal{C}$.
La catégorie $\mathcal{C}_{I-inj m}$ des objets I-injectifs marqués avec les morphismes
préservant les marquages possède les limites.
\end{lem}
{\it Preuve :} application directe du lemme~\ref{mst} !\\

Cette notion de marquage va nous permettre de construire un plan d'addition de
cellules qui à un objet quelconque associera un objet I-injectif marqué et
factorisera de manière unique tout morphisme à but I-injectif marqué. Non seulement
ce plan sera fonctoriel mais il servira d'adjoint au foncteur Oubli des objets
I-injectifs marqués vers les objets quelconques de $\mathcal{C}$.\\

Commençons tout d'abord par remarquer que si un plan d'addition de cellules
rajoute plusieurs fois la même cellule il ne pourra pas en général
factoriser de manière unique un morphisme à valeur dans un objet
I-injectif marqué.\\
Par exemple considérons l'ensemble $I=\{A\rightarrow B\}$ et l'objet
$C=B\coprod_A B$. Ce dernier est bien I-injectif et on peut le marquer de deux
manières selon l'exemplaire $B$ que l'on choisit. Choisissons un marquage.
Considérons maintenant le morphisme naturel $A\rightarrow C$ et le plan
$e_{I,2}$. Alors $e_{I,2}(A)$ n'est autre que $C$ et on le marque en choisissant
un exemplaire de $B$. Sur les quatre morphismes possibles de $C$ dans lui-même
laissant $A$ invariant, il y en a deux qui préservent le marquage : celui qui
envoie l'exemplaire non marqué de $B$ sur le non marqué et celui qui l'envoie
sur le marqué.\\

Cet exemple montre qu'un plan d'addition de cellules rajoutant plusieurs fois la
même flèche et donc en particulier tout plan d'addition de cellules non rationnel n'a pas en
général la propriété de factoriser de manière unique les morphismes
à but I-injectif marqué. Le lemme suivant nous donne de bons candidats pour
cette factorisation unique.

\begin{lem}
Soit $\mathcal{C}$ une catégorie cocomplète et $I$ un ensemble de morphismes
de $\mathcal{C}$.
Soient $A$ un objet de $\mathcal{C}$ et $P$ un plan d'addition de cellules de I sur
$A$ rationnel et dans lequel à chaque étape les diagrammes à rajouter ne sont pris qu'une
fois. Alors $P(A)$ a un marquage partiel naturel et tout morphisme $A\rightarrow B$ vers
un objet I-injectif marqué se factorise de manière unique à travers le
morphisme naturel $A\rightarrow P(A)$ en un morphisme préservant le marquage.
\end{lem}
{\it Preuve :}\\
Soit $f:A\rightarrow (B,\nu)$ un morphisme vers un objet I-injectif marqué. Notons $\lambda$ la longueur
du plan et définissons par récurrence (transfinie) sur $\beta\leq\lambda$
à la fois le marquage $\mu$ de $P(A)$ et le relèvement préservant le marquage de
$A\rightarrow P(A)$ par rapport à $B$.\\ 
Pour $\beta$ nul, $P_0(A)$ n'est autre
que $A$ qu'on laisse sans marquage, i.e. $\mu_0=\emptyset$. Le relèvement est dans ce cas $f$ qui
préserve bien le marquage.\\ 
Supposons maintenant que, pour $\beta<\lambda$,
$P_{\beta}(A)$ est partiellement marqué par $\mu_{\beta}$ et qu'il existe un
morphisme préservant le marquage $(P_{\beta}(A),\mu_{\beta})\rightarrow (B,\nu)$ qui précomposé par le morphisme
naturel $A\rightarrow P_{\beta}(A)$ redonne $f$. Le résultat $P_{\beta+1}(A)$
du plan simple $P_{\beta}$ appliqué à $P_{\beta}(A)$ n'est autre que la
somme amalgamée de $P_{\beta}(A)$ par les diagrammes de $P_{\beta}$. Comme le
plan est rationnel, les extensions à $P(A)$ de ces diagrammes ne sont pas les extensions à $P(A)$ des diagrammes des
plans simples précédents mais peuvent être celles de diagrammes déjà
relevés dans $A=P_0(A)$. Or comme $\mu_0$ est vide, les diagrammes qui se
relèvent déjà dans $A$ ne sont pas marqués. Ainsi aucune des extensions à $P(A)$ des diagrammes
de $P_{\beta}$ n'est l'extension à $P(A)$ d'un diagramme déjà marqué. De plus par
hypothèse sur $P$, les diagrammes de $P_{\beta}$ ne sont pris qu'une fois.
Ainsi tous les diagrammes de $P_{\beta}$ peuvent être marqués sans
ambiguïté ni incohérence dans $P_{\beta+1}(A)$ avec pour relèvements
marqués les morphismes naturels des buts des flèches de I dans la somme
amalgamée $P_{\beta+1}(A)$. On peut donc construire $\mu_{\beta+1}$ comme
réunion de $\mu_{\beta}$ avec les marquages naturels des diagrammes de
$P_{\beta}$. En outre $(B,\nu)$ étant un objet I-injectif marqué, il existe un
relèvement marqué de $B$ par rapport aux extensions par $f$ de tous les
diagrammes de $P_{\beta}$. Par propriété universelle de la colimite
$P_{\beta+1}(A)$, il vient qu'il existe un morphisme de $P_{\beta+1}(A)$ vers
$B$ qui d'une part précomposé par le morphisme naturel de $P_{\beta}(A)\rightarrow P_{\beta+1}(A)$
redonne le morphisme préservant le marquage $P_{\beta}(A)\rightarrow B$ et qui d'autre part
préserve les marquages par construction. Enfin par propriété du morphisme
$P_{\beta}(A)\rightarrow B$, le morphisme de $P_{\beta+1}(A)$ vers $B$ ainsi
obtenu redonne bien $f$ lorsqu'on le précompose par le morphisme naturel
$A\rightarrow P_{\beta+1}(A)$, ce qui finit de montrer l'hypothèse au rang
$\beta+1$. On remarque au passage qu'avec ces marquages le morphisme naturel
$P_{\beta}(A)\rightarrow P_{\beta+1}(A)$ préserve le marquage.\\

Soit maintenant $\beta\leq\lambda$ un ordinal limite (si $\lambda$ est un
cardinal transfini). Supposons donc construits pour tout $\alpha<\beta$ un
marquage $\mu_{\alpha}$ pour $P_{\alpha}(A)$ et un morphisme préservant le
marquage $(P_{\alpha}(A),\mu_{\alpha})\rightarrow (B,\nu)$ qui précomposé
par le morphisme naturel $A\rightarrow P_{\alpha}(A)$ redonne $f$. Comme
$P_{\beta}(A)$ est la colimite séquentielle transfinie des $P_{\alpha}(A)$ pour $\alpha<\beta$, définissons alors
$\mu_{\beta}$ comme la réunion des $\mu_{\alpha}$ pour $\alpha<\beta$. En
outre par propriété universelle de la colimite séquentielle transfinie, il existe un morphisme de
$P_{\beta}(A)$ vers $B$ qui précomposé par les morphismes naturels
$P_{\alpha}(A)\rightarrow P_{\beta}(A)$ redonnent les morphismes préservant le
marquage $(P_{\alpha}(A),\mu_{\alpha})\rightarrow (B,\nu)$. Par construction ce
morphisme préserve le marquage et précomposé par $A\rightarrow P_{\beta}(A)$ redonne $f$, ce qui
montre l'hypothèse de récurrence pour les ordinaux limites. On remarque là
encore que les morphismes naturels $P_{\alpha}(A)\rightarrow P_{\beta}(A)$ préservent les marquages grâce à
notre définition de $\mu_{\beta}$.\\
Par récurrence transfinie, on a bien montré qu'il existe un marquage naturel
$\mu$ sur $P(A)$ et qu'il existe un morphisme préservant le marquage de $P(A)$
vers $B$ qui précomposé par le morphisme naturel $A\rightarrow P(A)$ redonne
$f$.\\

Montrons maintenant par l'absurde que le morphisme préservant le marquage de
$P(A)$ vers $B$ qui précomposé par le morphisme naturel $A\rightarrow P(A)$ redonne
$f$ est unique. Supposons qu'il existe deux tels morphismes. Posons $\beta$
le plus petit ordinal inférieur ou égal à $\lambda$ pour lequel les
précompositions de ces morphismes par le morphisme naturel
$P_{\beta}(A)\rightarrow P_{\lambda}(A)$ sont différentes. Notons-les
$f_1,f_2$. Cet ordinal existe
bien car l'ensemble des ordinaux vérifiant cette propriété est non vide du
fait qu'il contient $\lambda$. Supposons que $\beta$ est un ordinal limite
(uniquement si $\lambda$ est un cardinal transfini), alors pour tout
$\alpha<\beta$ les précompositions de $f_1$ et $f_2$ par le morphisme naturel
$P_{\alpha}(A)\rightarrow P_{\beta}(A)$ sont des morphismes identiques par
minimalité de $\beta$. Mais donc par propriété universelle de la colimite
séquentielle transfinie $P_{\beta}(A)$, les morphismes $f_1$ et $f_2$ sont
identiques, ce qui contredit la définition de $\beta$. Supposons alors que
$\beta$ n'est pas un ordinal limite. Soit $\beta$ est nul, ce qui donnerait que
$f_1=f_2=f$, ce qui est absurde, soit $\beta$ admet un prédécesseur. Par
minimalité de $\beta$, on sait que les précompositions de $f_1$ et $f_2$ par le morphisme naturel
$P_{\beta-1}(A)\rightarrow P_{\beta}(A)$ sont des morphismes identiques. Comme
$P_{\beta}(A)$ est la somme amalgamée de $P_{\beta-1}(A)$ par les diagrammes de
$P_{\beta-1}$, regardons ce qui se passe au niveau des diagrammes. On remarque
en outre que les morphismes $f_1$ et $f_2$ préservent le marquage comme composés
de morphismes préservant le marquages par le morphisme naturel
$P_{\beta}(A)\rightarrow P(A)$ qui préserve le marquage, comme on l'a vu dans
la première partie de la démonstration. Par
définition du marquage des diagrammes de $P_{\beta-1}$, par rationalité de $P$ et du fait que $f_1$ et
$f_2$ préservent le marquage, il vient que les restrictions de $f_1$ et $f_2$
aux buts des flèches de I apparaissant dans les diagrammes de $P_{\beta-1}$
sont identiques et valent le relèvement marqué dans $\nu$ pour ces
flèches. Ainsi par propriété universelle de la somme amalgamée
$P_{\beta}(A)$, il vient que $f_1$ et $f_2$ sont égaux, ce qui contredit la
définition de $\beta$.\\
On a donc montré par l'absurde l'unicité de la factorisation de $f$ à
travers $A\rightarrow P(A)$ par un morphisme préservant le marquage.\\
CQFD.\\

Après avoir donné un critère pour qu'un plan d'addition de cellules
factorise de manière unique à travers un morphisme préservant le marquage
les morphismes à but I-injectif marqué, donnons un exemple de plan
vérifiant ce critère.

\begin{cor}
Soit $\mathcal{C}$ une catégorie cocomplète et $I$ un ensemble de morphismes
de $\mathcal{C}$ tel que les $I$-cofibrations soient des monomorphismes. Tout plan d'addition de cellules $P$ obtenu par des compositions rationnelles et dont tous les plans
simples sont $e_{I,1}$ vérifie les propriétés suivantes :
\item - pour tout objet $A$ de $\mathcal{C}$, $P(A)$ a un marquage partiel
naturel,
\item - tout morphisme $A\rightarrow B$ vers
un objet I-injectif marqué se factorise de manière unique à travers le
morphisme naturel $A\rightarrow P(A)$ en un morphisme préservant le marquage.
\end{cor} 
{\it Preuve :}\\
Comme par la proposition~\ref{comprat}, $P$ est rationnel car obtenu par des compositions rationnelles, il suffit de voir qu'à chaque étape
les diagrammes n'apparaissent qu'avec la cardinalité un, ce qui est bien le
cas ici car à chaque étape le plan simple n'est autre que $e_{I,1}$. On peut
donc appliquer le lemme précédent à $P$, ce qui nous donne le résultat.\\
CQFD.\\

En réunissant ce corollaire avec le lemme~\ref{Ephi1}, on obtient un plan
d'addition de cellules satisfaisant les propriétés i) et ii) de la
I-injectivisation avec en prime l'unicité pour i) quand on considère les
objets I-injectifs marqués.

\begin{prop}\label{Ephi2}
Soient $\mathcal{C}$ une catégorie cocomplète et $\Phi$ un ensemble de
morphismes de $\mathcal{C}$ dont les sources sont $\alpha$-petites pour un certain cardinal régulier
$\alpha$ et tel que les $\Phi$-cofibrations sont des monomorphismes. Soit $\alpha'$ le plus petit cardinal régulier supérieur à
$\alpha$.\\
Notons $E_{\Phi}$ le plan d'addition de cellules fonctoriel de
longueur $\alpha'$ obtenu par des compositions rationnelles et dont tous les plans simples sont $e_{\Phi,1}$.\\
Alors, pour tout objet $A$ de $\mathcal{A}$, 
$E_{\Phi}(A)$ est un objet $\Phi$-injectif marqué et tout morphisme $f:A\rightarrow
B$ dont le but est $\Phi$-injectif marqué se factorise de manière unique à travers le morphisme naturel
$A\rightarrow E_{\Phi}(A)$ en un morphisme préservant le marquage.
\end{prop}
{\it Preuve :}\\
Le plan ici noté $E_{\Phi}$ n'est autre que le plan $E'_{\Phi,1}$ qui
d'après le lemme~\ref{Ephi1} rend I-injectifs les objets auxquels on
l'applique. En outre $E_{\Phi}$ étant obtenu par des compositions rationnelles et uniquement composé de plans simples
$e_{\Phi,1}$, on peut lui appliquer le corollaire précédent, ce qui termine
de montrer la proposition.\\
CQFD.

\begin{cor}
Soient $\mathcal{C}$ une catégorie cocomplète et $\Phi$ un ensemble de
morphismes de $\mathcal{C}$ dont les sources sont $\alpha$-petites pour un certain cardinal régulier
$\alpha$ et tel que les $\Phi$-cofibrations sont des monomorphismes. Soit $\alpha'$ le plus petit cardinal régulier supérieur à
$\alpha$.\\
Notons $E_{\Phi}$ le plan d'addition de cellules fonctoriel de
longueur $\alpha'$ obtenu par des compositions rationnelles et dont tous les plans simples sont $e_{\Phi,1}$.\\
Alors $E_{\Phi}$ est un foncteur de $\mathcal{C}$ vers $\mathcal{C}_{I-inj m}$
qui est adjoint à gauche du foncteur $Oubli$, i.e. on a un isomorphisme
naturel en $A$ et $(B,\nu)$ induit par le morphisme naturel $A\rightarrow
E_{\Phi}(A)$ :
$$\operatorname{Hom}_{\mathcal{C}_{I-inj
m}}(E_{\Phi}(A),(B,\lambda))=\operatorname{Hom}_{\mathcal{C}}(A,Oubli((B,\nu)))$$
\end{cor}
{\it Preuve :} cela découle directement de l'unicité de la factorisation
par un morphisme préservant le marquage de tout morphisme à but un objet
I-injectif marqué. CQFD.\\
\\

Comme on l'a vu plus haut, le point fort de cette construction particulière est que si
l'on exhibe un plan d'addition de cellules vérifiant les propriétés i) et
ii) de l'I-injectivisation mais aussi la propriété iii) pour les objets
I-injectifs, alors notre construction vérifiera aussi iii), ce qui en fera un
procédé d'I-injectivisation.

\newpage

\section{I-injectivisation et \eq~d'objets I-injectifs}

Maintenant que l'on a un plan d'addition de cellules $E_{I}$ fonctoriel,
vérifiant les propriétés i) et ii) de l'I-injectivisation et qui assure
l'unicité de la factorisation du i) pour les objets I-injectifs marqués, il ne
nous reste plus qu'à montrer que $E_{I}$ vérifie la propriété iii), à
savoir que pour tout objet $A$, l'image par $E_{I}$ du morphisme naturel
$A\rightarrow E_{I}(A)$ est une \eq~d'objets I-injectifs, lorsqu'on s'est donné une telle notion. Pour cela on va
supposer que l'on connaît un plan d'addition de cellules fonctoriel $P$ vérifiant i) et ii)
tel que, pour tout objet $A$ I-injectif, le morphisme naturel $A\rightarrow
P(A)$ est une \eq~d'objets I-injectifs.

\begin{prop}\label{Ephi3}
Soient $\mathcal{C}$ une catégorie cocomplète et $\Phi$ un ensemble de
morphismes de $\mathcal{C}$ dont les sources sont $\alpha$-petites pour un certain cardinal régulier
$\alpha$ et tel que les $\Phi$-cofibrations soient des monomorphismes. Soit $\alpha'$ le plus petit cardinal régulier supérieur à
$\alpha$. Considérons sur $\mathcal{C}$ une notion d'\eq~d'objets
$\Phi$-injectifs vérifiant les propriétés suivantes :
\item - pour tout couple $(f,g)$ de morphismes composables entre objets
$\Phi$-injectifs, si parmi $f,g,g\circ f$ deux morphismes sont des \eqs~d'objets
$\Phi$-injectifs, alors le troisième morphisme aussi,
\item - pour tout couple $(f,g)$ de morphismes composables entre objets
$\Phi$-injectifs, si la composée $g\circ f$ est l'identité et que la
composée $f\circ g$ une \eq~d'objets $\Phi$-injectifs, alors $f$ et $g$ sont
des \eqs~d'objets $\Phi$-injectifs,
\item - les isomorphismes entre objets $\Phi$-injectifs sont des \eqs~d'objets
$\Phi$-injectifs.\\
\\
Notons $E_{\Phi}$ le plan d'addition de cellules fonctoriel de
longueur $\alpha'$ obtenu par des compositions rationnelles et dont tous les plans simples sont $e_{\Phi,1}$.\\
Soit $P$ un plan d'addition de cellules de $\Phi$ fonctoriel vérifiant les
propriétés suivantes :
\item - pour tout objet $A$ de $\mathcal{C}$, $P(A)$ est $\Phi$-injectif,
\item - tout morphisme $A\rightarrow B$ dont le but est $\Phi$-injectif se
factorise à travers le morphisme naturel $A\rightarrow P(A)$,
\item - pour tout objet $A$ $\Phi$-injectif, le morphisme naturel $A\rightarrow
P(A)$ est une \eq~d'objets $\Phi$-injectifs.\\
\\
Alors $E_{\Phi}$ est une $\Phi$-injectivisation,
c'est-à-dire vérifie les propriétés suivantes :
\item -i) pour tout objet $A$ de $\mathcal{C}$, $E_{\Phi}(A)$ est $\Phi$-injectif,
\item -ii) tout morphisme $A\rightarrow B$ dont le but est $\Phi$-injectif se
factorise à travers le morphisme naturel $A\rightarrow E_{\Phi}(A)$,
\item -iii) pour tout objet $A$ de $\mathcal{C}$, l'image par $E_{\Phi}$ du morphisme naturel $A\rightarrow
E_{\Phi}(A)$ est une \eq~d'objets $\Phi$-injectifs.\\
\\
De plus pour tout morphisme $f$ de $\mathcal{C}$ entre objets I-injectifs, on
aura l'\eq~suivante :\\
\\
$f$ est une \eq~d'objets $\Phi$-injectifs si et seulement si $E_{\Phi}(f)$ l'est.
\end{prop}
{\it Preuve :}\\
Par propriété du plan $P$, pour tout objet $A$ de $\mathcal{C}$, $P(A)$ est
$\Phi$-injectif. Choisis\-sons-lui un marquage. Il est facile de voir alors que
pour la propriété de factorisation par $P(A)$ de tout morphisme de $A$ vers
un but $\Phi$-injectif marqué, il existe un morphisme marqué de $P(A)$ vers
$B$ préservant le marquage et réalisant la factorisation voulue. En effet,
$P(A)$ étant une colimite séquentielle transfinie de sommes amalgamées de
flèches de $\Phi$, il suffit d'envoyer les flèches marquées de $P(A)$
sur celles marquées de $B$ et d'utiliser la \prd~de $B$ pour les flèches de
$\Phi$ qui ne porte pas le marquage de $P(A)$. On remarque cependant qu'une
telle factorisation par un morphisme marqué n'est pas unique du fait de la
présence dans $P(A)$ de flèches de $\Phi$ ne portant pas le marquage.\\

La construction $E_{\Phi}$ vérifie les propriétés i) et ii) de la
$\Phi$-injectivisation par la proposition~\ref{Ephi2}. De plus par cette
proposition, tout morphisme à valeur dans un objet
$\Phi$-injectif marqué admet une factorisation marquée unique à travers
$E_{\Phi}$. Il ne reste dons plus qu'à montrer iii). Pour cela on va tout
d'abord comparer les plans $P$ et $E_{\Phi}$, ce qui nous permettra de montrer
que $E_{\Phi}$ vérifie lui aussi la troisième propriété de $P$. C'est de
là que découlera la propriété iii) pour $E_{\Phi}$.\\

Pour la suite, choisissons un objet $A$ quelconque de $\mathcal{C}$ et notons
$e_A$ le morphisme naturel de $A$ vers $E_{\Phi}(A)$ et $p_A$ le morphisme
naturel de $A$ vers $P(A)$. Fixons un choix de marquage pour $P(A)$. D'après
ce qui précède, il existe un morphisme $pe_A$ de $P(A)$ vers $E_{\Phi}(A)$
préservant le marquage et tel que $pe_A\circ p_A=e_A$. De plus par
propriété de $E_{\Phi}(A)$, il existe un unique morphisme $ep_A$ de
$E_{\Phi}(A)$ vers $P(A)$ préservant le marquage et tel que $ep_A\circ
e_A=p_A$.
Considérons l'égalité suivante : $pe_A\circ ep_A\circ e_A=pe_A\circ
p_A=e_A$. Or par propriété d'unicité de factorisation marquée par $e_A$, il vient
que $pe_A\circ ep_A$ n'est autre que $Id_{E_{\Phi}(A)}$. Considérons
maintenant l'égalité suivante : $ep_A\circ pe_A\circ p_A=ep_A\circ e_A=p_A$.
Si $A$ est $\Phi$-injectif, alors par propriété de $P$, on a que $p_A$ est une
\eq~d'objets $\Phi$-injectifs. Ici on a que $p_A$ et $(ep_A\circ pe_A)\circ p_A=p_A$
sont des \eqs~d'objets $\Phi$-injectifs. Par la première propriété des \eqs~d'objets
$\Phi$-injectifs, il vient que $ep_A\circ pe_A$ est une \eq~d'objets
$\Phi$-injectifs. On a donc que $pe_A\circ ep_A$ est une identité et que
$ep_A\circ pe_A$ est une \eq~d'objets $\Phi$-injectifs. Alors par la seconde
propriété des \eqs~d'objets $\Phi$-injectifs, il vient que $ep_A$ et $pe_A$ sont des
\eqs~d'objets $\Phi$-injectifs, lorsque $A$ est $\Phi$-injectif.\\
Supposons encore que $A$ est $\Phi$-injectif. Considérons l'égalité
$ep_A\circ e_A=p_A$. Deux des trois morphismes sont des \eqs~d'objets
$\Phi$-injectifs : en effet $p_A$ est une \eq~par propriété de $P$ et $ep_A$
est une \eq~par ce qui précède. Par la première propriété des
\eqs~d'objets $\Phi$-injectifs, il vient que le troisième morphisme $e_A$ est
une \eq~d'objets $\Phi$-injectifs. Ceci montre que pour tout objet $A$
$\Phi$-injectif, le morphisme naturel $e_A:A\rightarrow E_{\Phi}(A)$ est une
\eq~d'objets $\Phi$-injectifs.\\

Supposons maintenant $A$ quelconque. Appliquons la propriété de
factorisation par $e_{E_{\Phi}(A)}$ à l'identité de $E_{\Phi}(A)$. On obtient un unique
morphisme $r:E_{\Phi}(E_{\Phi}(A))\rightarrow E_{\Phi}(A)$ préservant le
marquage et tel que $r\circ e_{E_{\Phi}(A)}=Id_{E_{\Phi}(A)}$. Comme
$E_{\Phi}(A)$ est $\Phi$-injectif, par le résultat précédent, il vient que
$e_{E_{\Phi}(A)}$ est une \eq~d'objets $\Phi$-injectifs. Comme l'identité en
est aussi une, par la première propriété des \eqs~d'objets
$\Phi$-injectifs, on obtient que $r$ est une \eq~d'objets $\Phi$-injectifs.
Considérons alors le diagramme commutatif suivant :
\begin{diagram}
A & \rTo^{e_A} & E_{\Phi}(A) & \rTo^{Id} & E_{\Phi}(A) \\
\dTo^{e_A} & & \dTo^{e_{E_{\Phi}(A)}} & \ruTo_r & \\
E_{\Phi}(A) & \rTo_{E_{\Phi}(e_A)} & E_{\Phi}(E_{\Phi}(A)) & & \\
\end{diagram}
Par ce diagramme, on obtient l'égalité suivante : $r\circ E_{\Phi}(e_A)\circ
e_A=r\circ e_{E_{\Phi}(A)}\circ e_A=Id\circ e_A=e_A$. Par unicité de la
factorisation marquée par $e_A$, il vient que $r\circ E_{\Phi}(e_A)$ n'est
autre que l'identité de $E_{\Phi}(A)$. Or l'identité et $r$ sont des
\eqs~d'objets $\Phi$-injectifs, donc, toujours par la première propriété
des \eqs~d'objets $\Phi$-injectifs, on a bien que $E_{\Phi}(e_A)$ est une
\eq~d'objets $\Phi$-injectifs. Ceci montre la propriété iii) de la
$\Phi$-injectivisation pour $E_{\Phi}$. Comme $E_{\Phi}$ a déjà les propriétés i) et
ii), on a donc montré que $E_{\Phi}$ est bien une $\Phi$-injectivisation.\\

Pour terminer, considérons un morphisme $f:A\rightarrow B$ quelconque entre
objets $\Phi$-injectifs. On a alors le diagramme commutatif suivant :
\begin{diagram}
A & \rTo^f & B\\
\dTo^{e_A} & & \dTo_{e_B}\\
E_{\Phi}(A) & \rTo_{E_{\Phi}(f)} & E_{\Phi}(B)\\
\end{diagram}
Comme $A$ et $B$ sont $\Phi$-injectifs, par ce qui précède, on a que les
morphismes $e_A$ et $e_B$ sont des \eqs~d'objets $\Phi$-injectifs. En appliquant
à ce diagramme la première propriété des \eqs~d'objets $\Phi$-injectifs,
il vient que $f$ est une \eq~d'objets $\Phi$-injectifs si et seulement si
$E_{\Phi}(f)$ l'est, ce qui montre le dernier résultat de la proposition.\\
CQFD.\\

Au passage, on remarque lors de cette démonstration l'utilité de l'unicité
de la factorisation du plan $E_{\Phi}$. En effet, dans la démonstration, on
montre d'abord que $E_{\Phi}$ vérifie que pour tout objet $\Phi$-injectif $A$,
le morphisme naturel $e_A:A\rightarrow E_{\Phi}(A)$ est une \eq~d'objets
$\Phi$-injectifs. A ce moment-là, $E_{\Phi}$ et $P$ ont alors les mêmes
propriétés à l'exception de l'unicité de la factorisation. Or c'est
précisément grâce à cette dernière propriété que l'on montre que $E_{\phi}$
vérifie iii). \\
En outre la dernière propriété de la proposition nous permet d'étendre la notion
d'\eq~pour les objets $\Phi$-injectifs aux objets quelconques.

\begin{cor}\label{Ephi3'}
Sous les hypothèses de la proposition précédente,
définissons une notion de $\Phi$-\eq~comme suit :\\
un morphisme $f$ quelconque de $\mathcal{C}$ est une $\Phi$-\eq~si $E_{\Phi}(f)$
est une \eq~d'objets $\Phi$-injectifs.\\ 
Alors ces $\Phi$-\eqs~vérifient les propriétés suivantes :
\item - pour tout couple de morphismes
composables $(f,g)$, si parmi $f,g,g\circ f$ deux morphismes sont des
$\Phi$-\eqs, alors le morphisme restant aussi. On appellera cette propriété
de la $\Phi$-\eq~la propriété du "trois pour deux".
\item - pour tout couple de morphismes
composables $(f,g)$, si $f\circ g$ est l'identité et $g\circ f$ une
$\Phi$-\eq, alors $f$ et $g$ sont des $\Phi$-\eqs.
\item - les isomorphismes de $\mathcal{C}$ sont des $\Phi$-\eqs.
\item - pour les morphismes d'objets $\Phi$-injectifs, les notions de $\Phi$-\eq~et d'\eq~d'objets
$\Phi$-injectifs se confondent.
\end{cor}
{\it Preuve :}\\
Soient $f$ et $g$ deux morphismes composables. Supposons que parmi $f,g,g\circ
f$ deux morphismes sont des $\Phi$-\eqs, ceci signifie que parmi
$E_{\Phi}(f),E_{\Phi}(g)$ et $E_{\Phi}(g)\circ E_{\Phi}(f)$ les deux morphismes
correspondant sont des \eqs~d'objets $\Phi$-injectifs. Or par propriété de
l'\eq~d'objets $\Phi$-injectifs, ceci nous donne que le troisième morphisme
est une \eq~d'objets $\Phi$-injectifs, donc que le troisième morphisme parmi
$f,g,g\circ f$ est une $\Phi$-\eq.\\
Si cette fois $f\circ g$ est l'identité et $g\circ f$ est une $\Phi$-\eq,
alors comme $E_{\Phi}$ est un foncteur, $E_{\Phi}(f)\circ E_{\Phi}(g)$ est
l'identité et $E_{\Phi}(g)\circ E_{\Phi}(f)$ une \eq~d'objets
$\Phi$-injectifs. Or par propriété de
l'\eq~d'objets $\Phi$-injectifs, ceci nous donne que $E_{\Phi}(f)$ et
$E_{\Phi}(g)$ sont des \eqs~d'objets $\Phi$-injectifs, donc que $f$ et $g$ sont
des $\Phi$-\eqs.\\
Soit $f$ un isomorphisme de $\mathcal{C}$. Comme $E_{\Phi}$ est un foncteur
rendant $\Phi$-injectif, alors $E_{\Phi}(f)$ est un isomorphisme entre objets
$\Phi$-injectifs. Par la troisième propriété vérifiée par les
\eqs~d'objets $\Phi$-injectifs, il vient que $E_{\Phi}(f)$ est une \eq~d'objets
$\Phi$-injectifs, donc que $f$ est une $\Phi$-\eq.\\
La proposition~\ref{Ephi3} montre que, pour tout morphisme $f$ entre objets
$\Phi$-injectifs, $f$ est une \eq~d'objets $\Phi$-injectifs si et seulement si
$E_{\Phi}(f)$ l'est, c'est-à-dire si et seulement si $f$ est une $\Phi$-\eq.
Ainsi pour les morphismes d'objets $\Phi$-injectifs, les notions de $\Phi$-\eq~et
d'\eq~d'objets $\Phi$-injectifs sont confondues.\\
CQFD.\\

Comme pour la catégorisation, un des avantages de l'I-injectivisation est de
pouvoir définir des sortes de colimites pour les objets I-injectifs. Il est
alors intéressant de pouvoir comparer la colimite des I-injectivisés avec la
I-injectivisation de la colimite. Dans les deux cas, il s'agit toujours de plans
d'addition de cellules. Ainsi, il serait utile de pouvoir comparer
les plans d'addition de cellules entre eux, ce que nous allons faire dans les
sections suivantes, par l'intermédiaire de la notion de rationalisation de
plans d'addition de cellules.

\newpage

\section{Rationalisation des plans d'addition de cellules}

Comme la notion de plan rationnel permet de mieux comprendre l'opération d'addition finie ou
infinie de cellules, on va montrer que l'on peut réduire tout plan d'addition
à un plan rationnel, ce qui nous permettra de mieux le manipuler. Pour cela
nous avons besoin d'un lemme technique sur la commutativité des colimites
séquentielles transfinies avec les sommes amalgamées.

\begin{lem}
Soient $\mathcal{C}$ une catégorie cocomplète et $\Phi$ un ensemble de
morphismes de $\mathcal{C}$. Soient $\lambda$ un cardinal transfini et
$A=\colimite{\beta<\lambda}(A_{\beta})$ une $\lambda$-séquence d'objets de
$\mathcal{C}$. Soit $P$ un plan simple sur $A$ dont tous les diagrammes se
factorisent par l'un des $A_{\beta}$. Pour tout $\beta<\lambda$, notons $P_{\beta}$ le plan
simple défini ainsi : un diagramme $D$ de $P$ est dans $P_{\beta}$ si et
seulement si $\beta$ est le plus petit indice tel que $d$ se factorise par
$A_{\beta}$. Et notons $P_{<\beta}$ le plan simple réunion des plans simples
$P_{\alpha}$ pour $\alpha<\beta$.\\
Posons maintenant $A'_0=A_0$ et, pour tout $\beta<\lambda$, posons :
$$A'_{\beta}=A_{\beta}\coprod_{\coprod_{D\in P_{<\beta}} \coprod_{\alpha_D}
S_{\phi_D}} \Bigg(\coprod_{D\in P_{<\beta}} \coprod_{\alpha_D}
B_{\phi_D}\Bigg).$$
On a alors l'égalité suivante :
$$\Bigg(\colim{\beta<\lambda} A_{\beta}\Bigg)\coprod_{\coprod_{D\in P} \coprod_{\alpha_D}
S_{\phi_D}} \Bigg(\coprod_{D\in P} \coprod_{\alpha_D}
B_{\phi_D}\Bigg)\cong\colim{\beta<\lambda} A'_{\beta} $$
\end{lem}
{\it Preuve :}\\
En utilisant les propriétés universelles des sommes
amalgamées, des coproduits et des colimites séquentielles transfinies, on
exhibe un unique morphisme de la somme amalgamée vers la colimite et un unique
morphisme de la colimite vers la somme amalgamée. Toujours en utilisant les
propriétés universelles, on montre que ces morphismes sont inverses l'un de
l'autre.\\
CQFD.

\begin{prop}\label{ratio}
Soient $\mathcal{C}$ une catégorie cocomplète et $\Phi$ un ensemble de
morphismes de $\mathcal{C}$ tel que les $\Phi$-cofibrations sont des monomorphismes. Alors tout plan d'addition de cellules fini ou
infini est rationnalisable dans le sens suivant. Soient $X$ un objet de
$\mathcal{C}$ et $P_.$ un plan quelconque sur $X$. Notons
$(X_{\beta})_{\beta\leq\lambda}$ sa suite associée. Alors il existe un plan
rationnel $Q_.$ sur $X$, dont la suite associée sera notée
$(Z_{\beta})_{\beta\leq\lambda}$ et, pour tout $\beta\leq\lambda$, des morphismes
dans $\mathcal{C}$ de $X_{\beta}$ vers $Z_{\beta}$ tels que le morphisme entre $P_.(X)$ et
$Q_.(X)$ soit un isomorphisme et que le diagramme suivant soit commutatif :
\begin{diagram}
X=X_0 & \rTo & \ldots & \rTo & X_{\beta} & \rTo & \ldots & \rTo & X_{\lambda}=P_.(X)\\
\dTo^{=}&    &        &      & \dTo      &      &        &      & \dTo_{\cong} \\
X=Z_0 & \rTo & \ldots & \rTo & Z_{\beta} & \rTo & \ldots & \rTo & Z_{\lambda}=Q_.(X)\\
\end{diagram}
\end{prop}
{\it Preuve :}\\
Soit $\lambda$ un cardinal transfini quelconque. Montrons que tout plan de
longueur inférieure ou égale à $\lambda$ est rationnalisable en un plan à compositions rationnelles, par
récurrence transfinie sur la longueur du plan. En effet comme par hypothèse les $\Phi$-cofibrations sont des monomorphismes, par la proposition~\ref{comprat}, les plans à compositions rationnelles sont des plans rationnels. On va donc montrer la proposition dans le cas plus fort des plans à compositions rationnelles.
Si
le plan est à un pas alors il est nécessairement à composition rationnelle car non
constitué de composition de plans simples.\\ 
Supposons la proposition vraie pour tout plan à $\beta$ pas où
$\beta<\lambda$ est
un ordinal quelconque. Soit $P_.$ un plan sur $X$ à $\beta+1$ pas. Notons
$(X_{\alpha})_{\alpha\leq\beta+1}$ la $\beta+1$-séquence associée. Par
définition de plan d'addition, on a que $P_.(X)$ est la somme amalgamée de
$X_{\beta}$ par un coproduit de cellules indexé par $P_{\beta}$. Or
$X_{\beta}$ est le résultat d'un plan à $\beta$ pas donc, par hypothèse,
il existe un plan à composition rationnelle à $\beta$ pas $Q_.$, dont la $\beta$-séquence
associée sera notée $(Z_{\alpha})_{\alpha\leq\beta}$, et, pour tout
$\alpha\leq\beta$, des morphismes de
$X_{\alpha}$ vers $Z_{\alpha}$ dans $X/\mathcal{C}$ tel que celui entre
$X_{\beta}$ et $Q_.(X)$ est un isomorphisme. 
Avec les notations du lemme précédent, il vient que la somme amalgamée de
$X_{\beta}$ par un coproduit de cellules indexé par
$\Big(P_{\beta}\Big)_{<\beta}$ est isomorphe dans $X/\mathcal{C}$ à la 
somme amalgamée de $Z_{\beta}$ par un coproduit de cellules indexé par
$P'$, extension de $\Big(P_{\beta}\Big)_{<\beta}$ par
l'isomorphisme de $X_{\beta}$ à $Z_{\beta}$. L'existence de morphismes de
$X_{\alpha}$ vers $Z_{\alpha}$ assure que les diagrammes de $P'$ se factorisent
par les $Z_{\alpha}$ pour $\alpha<\beta$. Par le lemme précédent, cette somme
amalgamée est isomorphe dans $X/\mathcal{C}$ à une colimite
séquentielle transfinie de $Z'_{\alpha}$. Par définition, $Z'_{\alpha+1}$ est
la somme amalgamée de $Z_{\alpha+1}$ par un coproduit de cellules indexé par
$P'_{<\alpha+1}$. Or $Z_{\alpha+1}$ est lui-même la somme
amalgamée de $Z_{\alpha}$ par un coproduit de cellules indexé par
$Q_{\alpha}$. Donc $Z'_{\alpha+1}$ est la somme amalgamée de $Z_{\alpha}$ par
un coproduit de cellules indexé par $P'_{<\alpha+1}$ et
$Q_{\alpha}$. Comme la somme amalgamée de $Z_{\alpha}$ par
un coproduit de cellules indexé par $P'_{<\alpha}$ n'est autre
que $Z'_{\alpha}$, il vient que $Z'_{\alpha+1}$ est la somme amalgamée de $Z'_{\alpha}$ par
un coproduit de cellules indexé par $P'_{\alpha}$ et
$Q_{\alpha}$. Posons donc $Q'_{\alpha}$ le plan simple réunion de $P'_{\alpha}$ et
$Q_{\alpha}$. On vient donc de montrer que la suite des $Q'_{\alpha}$ n'est autre
que le plan associé à la $\beta$-séquence des $Z'_{\alpha}$.\\
\\

 Montrons que pour tout $\alpha<\beta$ les
diagrammes de $Q'_{\alpha}$ ne se factorisent pas par $Z'_{\gamma}$, pour
$\gamma<\alpha$. Soit $D$ un diagramme de $Q'_{\alpha}$ qui se factorise par
$Z'_{\gamma}$, avec $\gamma\leq\alpha$, comme en outre ce diagramme a valeur dans $Z_{\alpha}$, alors il
se factorise par le produit fibré de $Z'_{\gamma}$ par $Z_{\alpha}$ sur
$Z'_{\alpha}$ qui n'est autre que $Z_{\gamma}$. Or si $D$ est un diagramme de
$Q_{\alpha}$, comme $Q_.$ est un plan à compositions rationnelles, alors $\gamma$ n'est autre que $\alpha$
et si $D$ est dans $P'_{\alpha}$, par définition même de $P'_{\alpha}$, on a
donc $\gamma$ égal à $\alpha$. Ceci montre
que la suite des $Q'_{\alpha}$ forme un plan sur $X$ dont toutes les compositions sont rationnelles. En outre on a montré qu'il est isomorphe dans
$X/\mathcal{C}$ à la somme amalgamée de $X_{\beta}$ par un coproduit de
cellules indexé par $\Big(P_{\beta}\Big)_{<\beta}$ et muni de morphismes de
$X_{\alpha}$ vers $Z'_{\alpha}$, pour tout $\alpha\leq\beta$. On termine en posant :
$Q'_{\beta}$ est l'extension par cet isomorphisme de $\Big(P_{\beta}\Big)_{\beta}$ et $Z'_{\beta+1}$ la somme
amalgamée de $Z'_{\beta}$ par un coproduit de cellules indexé par
$Q'_{\beta}$. Par un argument similaire à ce qui précède en utilisant la
définition de $\Big(P_{\beta}\Big)_{\beta}$, on montre que les
diagrammes de $Q'_{\beta}$ ne se factorisent pas par $Z'_{\gamma}$, avec
$\gamma<\beta$, ainsi le plan $Q'_.$ à $\beta+1$ pas est à compositions rationnelles. En outre
comme la somme amalgamée de $X_{\beta}$ par un coproduit de cellules, indexé
par $\Big(P_{\beta}\Big)_{<\beta}$, et $Z'_{\beta}$ sont isomorphes dans $X/\mathcal{C}$, il vient
que $P_.(X)$ et la somme amalgamée de $Z'_{\beta}$ par $Q'_{\beta}$ sont
isomorphes dans $X/\mathcal{C}$. On a donc montré que $P_.$ a pour
rationalisation le plan à compositions rationnelles $Q'_.$, i.e. que $P_.$ est rationnalisable.\\
\\

Supposons la proposition vraie pour tout plan d'addition de cellules de longueur strictement
inférieure à $\beta$, pour $\beta<\lambda$ ordinal limite.
Soit $P_.$ un plan sur $X$ à $\beta$ pas. Notons
$(X_{\alpha})_{\alpha\leq\beta}$ la $\beta$-séquence associée. Par
définition de plan infini, on a que $P_.(X)$ est la colimite des
$X_{\alpha}$ pour $\alpha<\beta$. Construisons par récurrence transfinie sur
$\alpha<\beta$ une suite de couples de plans infinis $(Q_{\alpha,.},P_{\alpha,.})_{\alpha<\beta}$
dont les séquences associées sont notées respectivement
$(Z_{\alpha,\gamma})_{\gamma\leq\alpha}$ et
$(X_{\alpha,\gamma})_{\gamma>\alpha}$ vérifiant les propriétés suivantes :
les plans $Q_{\alpha,.}$ sont à compositions rationnelles, ces couples de plans sont munis de morphismes entre eux dans $X/\mathcal{C}$, ceux de
$X_{\gamma,\alpha}$ vers $Z_{\alpha,\alpha}$ et ceux
de $X_{\gamma,\beta}$ vers $X_{\alpha,\beta}$, pour tout $\gamma<\alpha$, sont
des isomorphismes.\\ 
Pour $\alpha$ nul, $Q_{\alpha,.}$ sera
un plan vide et $P_{\alpha,.}$ le plan $P$. Supposons construits ces couples
jusqu'au rang $\alpha$, alors le plan composé $P_{\alpha,\alpha}\circ
Q_{\alpha,.}$ est un plan à $\alpha+1$ pas. Comme $\alpha<\beta$ et que $\beta$
est un ordinal limite, par hypothèse de récurrence, la proposition
s'applique. Notons alors $Q_{\alpha+1,.}$ le plan à compositions rationnelles rationalisation de $P_{\alpha,\alpha}\circ
Q_{\alpha,.}$, cette rationalisation vient avec des morphismes dans
$X/\mathcal{C}$ entre la séquence de $P_{\alpha,\alpha}\circ
Q_{\alpha,.}$ et celle de $Q_{\alpha+1,.}$ tels que celui de
$X_{\alpha,\alpha+1}$ vers $Z_{\alpha+1,\alpha+1}$ est un
isomorphisme. On pose alors $P_{\alpha+1,.}$ l'extension de $P_{\alpha,.+1}$ par cet
isomorphisme. On obtient donc bien que le morphisme entre $P_{\alpha,.}\circ
Q_{\alpha,.}(X)$ et $P_{\alpha+1,.}\circ Q_{\alpha+1,.}(X)$ est un
isomorphisme. Supposons ces plans construits pour tout $\gamma<\alpha$, ordinal
limite. Alors on définit $Q_{\alpha,\gamma}$ comme la réunion des plans
$Q_{\delta,\gamma}$ pour tout $\gamma<\delta<\alpha$ et $Z_{\alpha,\gamma}$
comme la colimite séquentielle transfinie des 
$Z_{\delta,\gamma}$ pour tout $\gamma<\delta<\alpha$. Il est facile de voir
d'une part que, par un résultat similaire à celui du lemme précédent, 
$Q_{\alpha,.}$ est bien le plan infini associé à la suite des
$(Z_{\alpha,\gamma})_{\gamma<\alpha}$ et d'autre part que, par intervertion des
colimites, le $Z_{\alpha,\alpha}$ ainsi obtenu est isomorphe dans $X/\mathcal{C}$ à tous les
$X_{\gamma,\alpha}$, pour $\gamma<\alpha$. La compatibilité de ces
isomorphismes entraîne que les extensions par ces isomorphismes des plans
$P_{\gamma,.+\alpha-\gamma}$ donnent lieu à un unique plan que l'on notera
$P_{\alpha,.}$. On a donc les morphismes des $X_{\gamma,\beta}$ vers
$X_{\alpha,\beta}$ sont des isomorphismes. Enfin le plan $Q_{\alpha,.}$ est bien
à compositions rationnelles car un diagramme de $Q_{\alpha,\gamma}$ ne peut pas se factoriser par
un $Z_{\alpha,\gamma'}$, avec $\gamma'<\gamma<\alpha$, car il appartient en fait
à l'un des $Q_{\delta,\gamma}$, pour un certain $\delta<\alpha$, et donc s'il
se factorise par $Z_{\alpha,\gamma'}$, il se factoriserait aussi par
$Z_{\delta,\gamma'}$, ce qui contredirait la rationalité des compositions du plan
$Q_{\delta,.}$. On a donc bien le résultat cherché par récurrence
transfinie sur $\alpha<\beta$.\\
\\

Ainsi pour tout ordinal $\gamma<\beta$, on a une suite de plans simples
$\Big(Q_{\alpha,\gamma}\Big)_{\gamma<\alpha<\beta}$ et une $\beta$-séquence
d'objets $\Big(Z_{\alpha,\gamma}\Big)_{\gamma<\alpha<\beta}$. Posons $Q^{\gamma}_.$ le plan simple réunion
de cette suite de plans simples et $Z^{\gamma}$ la colimite de la
$\beta$-séquence des $Z_{\alpha,\gamma}$. Il est facile de voir, en utilisant
un résultat similaire au lemme précédent, que la suite
des $Q^{\gamma}_.$ forme un plan infini dont 
$(Z^{\gamma})_{\gamma\leq\beta}$ est la $\beta$-séquence associée, où
$Z^{\beta}$ note la colimite de la $\beta$-séquence. Nous allons
donc montrer que d'une part le plan infini des $Q^{\gamma}_.$ est à compositions rationnelles et
qu'il existe un isomorphisme dans $X/\mathcal{C}$ entre $P_.(X)$ et
$Z_{\beta}$.\\
Soient $\gamma<\beta$ un ordinal et $D$ un diagramme de $Q^{\gamma}$ tel que $d$
se factorise par un $Z^{\gamma'}$, avec $\gamma'\leq\gamma$. Comme $Q^{\gamma}$ est une réunion de plans
simples, il existe un ordinal $\alpha$ tel que $D$ soit en fait un diagramme de
$Q_{\alpha,\gamma}$. Donc $d$ a pour but $Z_{\alpha,\gamma}$. Or comme $d$ se
factorise par $Z^{\gamma}$, alors $d$ se factorise par $Z_{\alpha,\gamma'}$. Or
$D$ est un diagramme de $Q_{\alpha,\gamma}$ et par hypothèse $Q_{\alpha,.}$
est à compositions rationnelles donc $\gamma'$ est égal à $\gamma$. Ceci montre que le plan
infini des $Q^{\gamma}_.$ est à compositions rationnelles.\\
Par intervertion des colimites, on a que la colimite des $X_{\alpha,\beta}$ pour
$\alpha<\beta$ n'est autre que celle des $Z^{\gamma}$ pour $\gamma<\beta$. Or le
morphisme de $X_{\beta}$ dans cette colimite est une colimite séquentielle
transfinie d'isomorphismes par construction, donc on a bien que le morphisme
dans $X/\mathcal{C}$ de $X_{\beta}$ vers $Z^{\beta}$ est un isomorphisme.\\
Donc le plan $P_.$ à $\beta$ pas est bien rationnalisable. Donc par
récurrence transfinie sur $\beta\leq\lambda$, on a bien montré que tout plan
infini d'au plus $\lambda$ pas est rationnalisable.\\
CQFD.\\

Il est intéressant mais aussi utile pour la suite d'énoncer un corollaire,
tiré directement de la démonstration de rationalisation des plans
d'addition de cellules, qui donne une construction explicite de rationalisation des plans
d'addition infinie de cellules.

\begin{cor}\label{plrat}
Soient $\mathcal{C}$ une catégorie cocomplète et $\Phi$ un ensemble de
morphismes de $\mathcal{C}$ tel que les $\Phi$-cofibrations sont des monomorphismes. Soient $X$ un objet de
$\mathcal{C}$ et $P_.$ un plan quelconque sur $X$. Notons
$(X_{\beta})_{\beta\leq\lambda}$ sa suite associée. 
Alors il possède comme
rationalisation le plan à compositions rationnelles $Q_.$ sur $X$, dont la suite associée sera notée
$(Z_{\beta})_{\beta\leq\lambda}$, défini de la manière suivante.\\ 
Pour tout couple d'ordinaux $\alpha,\beta$ inférieur à $\lambda$, notons
$(P_{\alpha})_{\beta}$ l'ensemble des
couples $(D',\alpha_{D'})$ tel que $D'$ se factorise par $Z_{\beta}$ mais pas
par $Z_{\gamma}$ pour $\gamma<\beta$, où les $D'$ sont les extensions par
$X_{\alpha}\rightarrow Z_{\alpha}$ des diagrammes $D$ de $P_{\alpha}$ et où
$\alpha_{D'}$ est la cardinalité de $D$ dans $P_{\alpha}$.\\
On définit alors $Q_{\beta}$ comme la colimite des plans
$(P_{\alpha})_{\beta}$ pour $\alpha<\lambda$.
\end{cor}
{\it Preuve :} C'est le plan à compositions rationnelles associé à $P_.$ construit par
récurrence transfinie dans la démonstration de la proposition
précédente. CQFD.\\

Cette notion de rationalisation des plans d'addition de cellules nous permettra de
comparer les plans d'addition de cellules, car il découle de la proposition de
rationalisation des plans d'addition de cellules que deux plans d'addition de
cellules sont isomorphes si et seulement si leurs rationalisations le sont.

\newpage

\section{Plan rationnel $HCat$}

Dans cette partie, nous allons construire un plan d'addition de cellules rationnel $HCat_{\lambda}$ qui
sera le modèle de certains plans d'addition de cellules rationnels. En fait,
il s'agit d'un plan du type $E'_{\Phi,\lambda}$ au sens du lemme~\ref{Ephi1},
mais dont le nombre d'exemplaires d'une même cellule qu'il rajoute à chaque
fois n'est autre que le cardinal de sa longueur. Un tel plan d'addition de
cellules a la bonne propriété que composé avec un de ses sous-plans, il
donne un nouveau plan équivalent à lui-même. Donnons d'abord sa
définition explicite.

\begin{defin}\index{$HCat_{\lambda}$}
Soient $\mathcal{C}$ une catégorie cocomplète et $\Phi$ un ensemble de
morphismes de $\mathcal{C}$. Soit $\lambda$ un ordinal transfini limite.
On définit $HCat_{\lambda}$ par récurrence
transfinie. Posons $e^0_{\lambda}$ le plan à un pas $e_{\Phi,Card(\lambda)}$. Supposons
construit le plan $e^{\beta}_{\lambda}$, pour $\beta<\lambda$, alors on pose
$e^{\beta+1}_{\lambda}=e_{\Phi,Card(\lambda)}*e_{\lambda}^{\beta}$. Supposons construits
les plans $e^{\alpha}_{\lambda}$ pour tout $\alpha<\beta$ avec $\beta$ cardinal
limite, on définit alors le plan $e^{\beta}_{\lambda}$ comme colimite des
$e^{\alpha}_{\lambda}$, pour $\alpha<\beta$. Ainsi par récurrence transfinie,
on a donc défini pour tout $\beta<\lambda$ les plans $e^{\beta}_{\lambda}$.\\
On définit donc le plan $HCat_{\lambda}$ comme colimite des
$e^{\beta}_{\lambda}$ pour $\beta<\lambda$. Par définition même,
$HCat_{\lambda}$ est un plan à compositions rationnelles, et donc un plan rationnel si les $\Phi$-cofibrations sont des monomorphismes par proposition~\ref{comprat}.
\end{defin}

Ces plans rationnels d'addition de cellules $HCat_{\lambda}$ sont différents uniquement en fonction de leur
longueur. Toutefois pour certaines longueurs, ces plans d'addition de cellules
sont quand même identiques.

\begin{lem}\label{hcat}
Soient $\mathcal{C}$ une catégorie cocomplète et $\Phi$ un ensemble de
morphismes de $\mathcal{C}$ dont les sources sont $\mu$-petites, pour un certain ordinal $\mu$. Soit
$\lambda$ un ordinal transfini limite dont le cardinal est supérieur à $\mu$. Alors pour tout
ordinal limite $\lambda'$ ayant même cardinalité que $\lambda$, on a un
isomorphisme dans $X/\mathcal{C}$ entre $HCat_{\lambda}$ et $HCat_{\lambda'}$
\end{lem}
{\it Preuve :}\\ 
On remarque tout d'abord que, par définition de composition rationnelle, si un
plan est vide alors les seuls plans composables rationnellement avec lui
sont également vides. On remarque aussi que les deux plans $HCat_{\lambda}$ et
$HCat_{\lambda'}$ ne diffèrent que par leur longueur. Comparons les au plan
$HCat_{Card_{\lambda}}$ qui ne diffère des deux autres lui aussi que par sa longueur.
On remarque qu'à l'étape $Card(\lambda)$ des trois plans $HCat(X)$, par petitesse des sources des morphismes de
$\Phi$, tous les diagrammes de $\Phi$ à valeur dans $X_{Card(\lambda)}$ se
factorisent par des $X_{\beta}$ pour $\beta<Card(\lambda)$. Ceci entraîne
que le plan simple des trois $HCat(X)$ au niveau $Card(\lambda)$ est vide et donc de
même pour les suivants, ce qui montre que les plans $HCat_{\lambda}$ et
$HCat_{\lambda'}$ sont isomorphes au plan $HCat_{Card(\lambda)}$.\\
CQFD.\\

Avant de voir les propriétés de stabilité de notre plan rationnel
d'addition de cellules $HCat$ vis-à-vis de ses sous-plans, définissons ce
que l'on entend par sous-plan d'un plan d'addition de cellules.

\begin{defin}\index{sous-plan d'un plan d'addition de cellules}
Soient $\mathcal{C}$ une catégorie cocomplète et $\Phi$ un ensemble de
morphismes de $\mathcal{C}$ tel que les $\Phi$-cofibrations sont des monomorphismes. Soient $X$ un objet de $\mathcal{C}$ et $P_.$ et $R_.$ des plans d'addition de cellules sur $X$ non nécessairement de même longueur. Notons $X_{\beta}$ la $\beta$-ième étape du plan $P_.$ appliqué à $X$.\\
 On dit que le plan $P_.$ est un sous-plan de $R_.$ si :
\item - pour tout $\beta$, il existe un monomorphisme de $X_{\beta}$ vers $R(X)$ faisant commuter le diagramme suivant :
\begin{diagram}
X & \rTo & R(X) \\
\dTo & \ruTo & \\
X_{\beta} & & \\
\end{diagram}
Notons alors $P'_{\beta}$ (respectivement $R'_{\gamma}$) les extensions à $R(X)$ des plans simples $P_{\beta}$ (respectivement $R_{\gamma}$).
\item - la réunion des $P'_{\beta}$ est un
sous-ensemble de la réunion des $R'_{\gamma}$.
\end{defin}

Cette définition quelque peu compliquée de sous-plan traduit le fait que l'on considère comme sous-plan d'un plan d'addition de cellules donné $R_.$ tout plan prenant dans un ordre quelconque et avec un cardinal plus petit certaines des cellules du plan $R_.$. Cependant la longueur du sous-plan peut être plus grande que la longueur du plan mais il ne faut pas que le sous-plan ajoute plus de cellules que le plan $R_.$ ou qu'il ajoute des cellules différentes. Un moyen simple de vérifier cela est de regarder les extensions des cellules à $R(X)$ et de comparer alors l'ensemble des cellules à rajouter du sous-plan avec celui de $R_.$. Montrons maintenant la stabilité des plans de type $HCat_{\lambda}$ par rapport à leurs sous-plan.  

\begin{prop}
Soient $\mathcal{C}$ une catégorie cocomplète et $\Phi$ un ensemble de
morphismes de $\mathcal{C}$ dont les sources sont $\nu$-petites, pour un certain ordinal $\nu$, et tel que les $\Phi$-cofibrations sont des monomorphismes. Soient
$\lambda$ un ordinal limite transfini dont le cardinal est supérieur à celui
de $\nu$ et $P_.$
un sous-plan à compositions rationnelles de $HCat_{\lambda}$ de longueur $\mu$ tel que
$Card(\mu)\leq Card(\lambda)$.
Alors le plan composé $HCat_{\lambda}\circ P_.$ se rationalise en
$HCat_{\lambda}$.
\end{prop}
{\it Preuve : }\\
Afin de simplifier les notations, notons $(X_{\beta})_{\beta<\mu+\lambda}$ la
$\mu+\lambda$-séquence associée à $HCat_{\lambda}\circ P_.$, notons
$(e_{\beta})_{\beta<\mu}$ les plans simples composant $P_.$ et
$(e_{\beta})_{\mu\leq\beta<\mu+\lambda}$ ceux composant $HCat_{\lambda}$. Notons
enfin $Q_.$ la rationalisation de $HCat_{\lambda}\circ P_.$ obtenue par le
corollaire~\ref{plrat} et $(Z_{\beta})_{\beta<\mu+\lambda}$ sa
$\mu+\lambda$-séquence associée.\\
Montrons dans un premier temps que $Q_{\beta}$
est isomorphe à la $\beta$-ième étape de $HCat_{\mu+\lambda}$. Notons donc $A$,
respectivement $B$, leurs ensembles sous-jacents des diagrammes. Comme pour tout
$\alpha<\mu+\lambda$, les diagrammes de $(e_{\alpha})_{\beta}$ sont des diagrammes de
$\Phi$ à valeur dans $Z_{\beta}$ ne se factorisant pas par $Z_{\delta}$ pour
$\delta<\beta$, on a bien l'inclusion de $A$ dans $B$. Soit $D$ un diagramme
quelconque de $\Phi$ à valeur dans $Z_{\beta}$ ne se factorisant pas par $Z_{\delta}$ pour
$\delta<\beta$. On peut donc étendre $D$ à $Z_{\mu+\lambda}$ qui est
isomorphe à $X_{\mu+\lambda}$. Comme par hypothèse, les sources des morphismes de $\Phi$ est $\nu$-petites
et que $\lambda$ est un ordinal limite supérieur à $\nu$, l'extension de $D$
à $X_{\mu+\lambda}$ se factorise par un certain $X_{\gamma}$, avec $\gamma$
minimal. Or, si $\gamma\neq\mu$, $e_{\gamma}$ contient tous les diagrammes de $\Phi$ à valeur dans
$X_{\gamma}$ ne se factorisant pas par $X_{\delta}$ pour $\delta<\gamma$ et donc
en particulier l'extension de $D$ à $X_{\mu+\lambda}$. Et si $\gamma=\mu$,
alors comme $e_{\mu}$ contient tous les diagrammes de $\Phi$ à valeur dans
$X_{\mu}$, il contient en particulier l'extension de $D$. Ceci montre bien dans
les deux cas que $D$ appartient à $(e_{\gamma})_{\beta}$ donc à $Q_{\beta}$.
Et donc les ensembles $A$ et $B$ sont isomorphes.\\
\\
 Pour avoir l'isomorphisme
entre $Q_{\beta}$ et la $\beta$-ième étape de $HCat_{\mu+\lambda}$, il
suffit donc de montrer que tous les diagrammes ont dans $Q_{\beta}$ le cardinal
$Card(\mu+\lambda)=Card(\lambda)$.\\
Soit $D$ un diagramme de $Q_{\beta}$, donc il appartient à un certain
$(e_{\alpha})_{\beta}$, avec $\alpha$ minimal. Il est donc l'extension d'un diagramme $D'$ de $\Phi$
à valeur dans $X_{\alpha}$. Considérons deux cas.\\ 
Si $\alpha<\mu$, comme
$P_.$ est à compositions rationnelles, $D'$ n'apparaît dans $P_.$ qu'à l'étape $\alpha$. Or
à l'étape $\mu$ du plan composé, on s'autorise à reprendre des
diagrammes existants, donc $D'$ apparaît aussi dans $e_{\mu}$ mais il
n'apparaît plus au-delà par rationalité des compositions de $HCat_{\lambda}$. Or le
cardinal associé à $D'$ dans $e_{\mu}$ est $Card(\lambda)$ et celui dans
$e_{\alpha}$ est inférieur à $Card(\lambda)$ car $P_.$ est un sous-plan de
$HCat_{\lambda}$. Ainsi le cardinal associé à $D$ dans $Q_{\beta}$ est bien
$Card(\lambda)$.\\
Si $\alpha\geq\mu$, comme $\alpha$ est minimal, cela signifie que $D'$ ne se
factorise pas par $X_{\gamma}$ pour $\gamma<\alpha$. De plus comme $HCat_{\lambda}$
est à compositions rationnelles, $D'$ n'apparaît pas dans les $e_{\gamma}$ pour $\gamma>\alpha$.
Donc $D'$ n'apparaît que dans $e_{\alpha}$ où son cardinal est
$Card(\lambda)$. Ainsi le cardinal associé à $D$ dans $Q_{\beta}$ est bien
$Card(\lambda)$.\\
Donc dans les deux cas, le cardinal associé à $D$ dans $Q_{\beta}$ est bien
$Card(\lambda)$, ce qui montre que $Q_{\beta}$ est isomorphe à la
$\beta$-ième étape de $HCat_{\mu+\lambda}$, pour tout $\beta<\mu+\lambda$.
Il s'en déduit que $Q_.$ est isomorphe à $HCat_{\mu+\lambda}$ dans
$X/\mathcal{C}$.\\
\\
On conclut en remarquant que comme les sources des morphismes de $\Phi$ sont $\nu$-petites, que $\lambda$ est
un ordinal limite dont le cardinal est supérieur à celui de $\nu$ et que $Card(\mu+\lambda)=Card(\lambda)$
par hypothèse, le lemme~\ref{hcat} nous donne un
isomorphisme dans $X/\mathcal{C}$ entre $HCat_{\mu+\lambda}$ et
$HCat_{\lambda}$. Ceci montre bien que $HCat_{\lambda}$ est la rationalisation
de $HCat_{\lambda}\circ P_.$.\\
CQFD.\\

Bien évidemment ce résultat est intéressant pour le plan $HCat$ mais ce
résultat est en fait plus général qu'il n'y paraît. En effet tout plan
d'addition de cellules dont la rationalisation est de la forme $HCat$ vérifie
aussi un tel résultat de stabilité vis-à-vis de ses sous-plans. C'est ce
que l'on va énoncer dans le corollaire suivant.

\begin{cor}\label{ratplcomp}
Soient $\mathcal{C}$ une catégorie cocomplète et $\Phi$ un ensemble de
morphismes de $\mathcal{C}$ dont les sources sont $\nu$-petites, pour un certain ordinal $\nu$, et tel que les $\Phi$-cofibrations sont des monomorphismes. Soient
$\lambda$ un ordinal limite transfini dont le cardinal est supérieur à celui
de $\nu$, $R_.$ un plan infini dont la rationalisation est $HCat_{\lambda}$ et
$P_.$ un sous-plan de $R_.$ de longueur $\mu$ tel que
$Card(\mu)\leq Card(\lambda)$.
Alors le plan composé $R_.\circ P_.$ se rationalise en
$HCat_{\lambda}$.\\ 
On en déduit donc que, sous ces hypothèses, le plan composé $R_.\circ P_.$ et le plan $R_.$
sont isomorphes dans $X/\mathcal{C}$ car ils ont même rationalisation.
\end{cor}
{\it Preuve :}\\ 
Comme tout plan se rationalise en plan à compositions rationnelles, il vient que $P_.$ est isomorphe
dans $X/\mathcal{C}$ à un plan à compositions rationnelles $P'_.$ de même longueur $\mu$. En
outre comme $P_.$ est un sous-plan de $R_.$, c'est aussi un sous-plan de sa
rationalisation $HCat_{\lambda}$ et donc la rationalisation $P'_.$ de $P_.$
est bien un sous-plan de $HCat_{\lambda}$. D'où il vient que le plan composé
$R_.\circ P_.$ est isomorphe dans $X/\mathcal{C}$ au plan composé
$HCat_{\lambda}\circ P'_.$ qui vérifie bien les hypothèses de la proposition
précédente, ce qui permet ainsi de conclure.\\ 
CQFD.\\

Grâce à ce corollaire, on obtient le résultat fondamental suivant : pour tout plan
infini d'addition de cellules dont la rationalisation est un $HCat$, on peut
faire d'abord une partie (éventuellement tout) du plan puis faire ensuite le
plan tout entier sur ce que l'on a obtenu et on aboutit au final à la même chose que
si l'on avait fait directement tout le plan dès le début.\\
Nous allons voir des conséquences de ce résultat pour les plans d'addition
de cellules que nous utiliserons dans la suite de cette thèse pour les
\precats.

\newpage

\section{Exemples particuliers de plans d'addition de cellules à usage pour
les \precats}

Dans cette partie nous allons nous intéresser aux plans d'addition de cellules
qui formalisent ceux qui interviennent dans la construction $Bigcat$ des
\precats. En effet, nous avons déjà construit toute une machinerie pour la
catégorisation. Mais on a vu que pour avoir l'invariance d'homotopie par
catégorisation, il nous fallait exhiber un autre procédé de
catégorisation ayant la propriété de préserver le type d'homotopie des
\cats. C'est ce que fera $Bigcat$. L'idée de cette construction peut
s'exprimer en termes de plans d'addition de cellules de la manière suivante.
On a une partition de l'ensemble des cellules $\Phi$. Pour certains groupements
des éléments de la partition qu'on notera $\Psi$, on a un plan d'addition de
cellules qui a pour vocation de vérifier la propriété iii) pour les objets $\Phi$-injectifs tout en étant une étape de la $\Psi$-injectivisation. On construit
alors un nouveau plan d'addition qui combine ces plans d'addition de cellules
pour chaque groupement ayant cette propriété. Ce nouveau plan aura alors les
propriétés voulues. Donnons les formalisations de ces plans d'addition de
cellules.

\begin{defin}\index{$Raj_{\Psi}$}\index{$Cat_{\Psi,\lambda}$}
Soient $\mathcal{C}$ une catégorie cocomplète et $\Phi$ un ensemble de
morphismes de $\mathcal{C}$. Soient $\Psi$ un sous-ensemble de $\Phi$ et
$\lambda$ un ordinal transfini limite. Soit $Raj_{\Psi}$ un sous-plan (donc fonctoriel) de $e_{\Psi,1}$.
On définit $Cat_{\Psi,\lambda}$ par récurrence
transfinie. Posons $Cat^0_{\Psi,\lambda}$ le plan simple $e_{\Psi,1}$. Supposons
construit le plan $Cat^{\beta}_{\Psi,\lambda}$, pour $\beta<\lambda$, alors on pose
$Cat^{\beta+1}_{\Psi,\lambda}=Raj_{\Psi}\circ Cat_{\Psi,\lambda}^{\beta}$. Supposons construits
les plans $Cat^{\alpha}_{\Psi,\lambda}$ pour tout $\alpha<\beta$ avec $\beta$ cardinal
limite, on définit alors le plan $Cat^{\beta}_{\Psi,\lambda}$ comme colimite des
$Cat^{\alpha}_{\Psi,\lambda}$, pour $\alpha<\beta$. Ainsi par récurrence transfinie,
on a donc défini pour tout $\beta<\lambda$ les plans $Cat^{\beta}_{\Psi,\lambda}$.\\
On définit donc le plan $Cat_{\Psi,\lambda}$ comme colimite des
$Cat^{\beta}_{\Psi,\lambda}$ pour $\beta<\lambda$. En général,
$Cat_{\Psi,\lambda}$ n'est pas un plan rationnel.
\end{defin}

Il est important de remarquer que si le plan d'addition de cellules $Cat_{\Psi,\lambda}$ commence par un plan simple tel $e_{\Psi,1}$ qui prend en compte toutes les flèches de $\Psi$, il peut se poursuivre par des plans simples qui ne prennent pas toutes les flèches de $\Psi$ et qu'on a appelés $Raj_{\Psi}$.

\begin{defin}\index{$BCat(t)$}
Soient $\mathcal{C}$ une catégorie cocomplète et $\Phi$ un ensemble de
morphismes de $\mathcal{C}$. Soit $(\Phi_{k})_{i\in K}$ une partition de $\Phi$. 
Donnons nous une fonction $t$ d'un ordinal limite $T$ quelconque
vers l'ensemble des parties de $K$ telle
que, pour tout $k$ dans $K$, on ait les propriétés suivantes :
\item - pour tout ordinal $\tau<T$, il existe un ordinal
$\tau\leq\tau'<T$ tel que $t(\tau')$ contienne $k$,
\item - le cardinal de l'ensemble $t^{-1}(k)$
est égal au cardinal de $T$.\\
\\ 
On définit $BCat(t)$ comme le plan d'addition de cellules suivant l'ordre $t$
qui pour tout ordinal $\tau<T$ a pour plan d'addition de cellules
$Cat_{\cup_{k\in t(\tau)}\Phi_k,T}$.
Le plan $BCat(t)$ est de longueur
$T^2$ et, en général, ce n'est pas un plan rationnel.
\end{defin}

Tout d'abord vérifions que ces plans d'addition de cellules ont bien les
propriétés attendues.

\begin{lem}\label{bcat}
Soient $\mathcal{C}$ une catégorie cocomplète et $\Phi$ un ensemble de
morphismes de $\mathcal{C}$. Soient $(\Phi_k)_{k\in K}$ une partition de $\Phi$  et $\alpha$ un cardinal régulier pour lequel les sources des flèches de $\Phi$ sont $\alpha$-petites.\\ 
Donnons nous une fonction $t$ d'un ordinal
limite $T$ de cardinal un cardinal régulier strictement supérieur à $\alpha$
vers l'ensemble des parties de $K$ telle que, pour tout $k$ dans $K$, on ait les propriétés suivantes : 
\item - pour tout ordinal $\tau<T$, il existe un ordinal
$\tau\leq\tau'<T$ tel que $t(\tau)$ contienne $k$,
\item - le cardinal de l'ensemble $t^{-1}(k)$ est égal au cardinal de $T$.\\
\\  
Alors le plan d'addition de cellules $BCat(t)$ vérifie les propriétés suivantes :
\item - pour tout objet $A$ de $\mathcal{C}$, $BCat(t)(A)$ est $\Phi$-injectif,
\item - tout morphisme $A\rightarrow B$ de $\mathcal{C}$ à but $\Phi$-injectif
se factorise par le morphisme naturel $A\rightarrow BCat(t)(A)$.\\
\end{lem}
{\it Preuve :}\\
Le morphisme naturel $A\rightarrow BCat(t)$ se construit comme colimite
séquentielle transfinie de sommes amalgamées de flèches de $\Phi$, donc
c'est une $\Phi$-cofibration. Ainsi le morphisme naturel $A\rightarrow BCat(t)$
a la \prg~par rapport aux objets $\Phi$-injectifs, ce qui montre la deuxième
propriété.\\
Pour la première propriété, on va montrer que $BCat(t)$ vérifie les
hypothèses du corollaire~\ref{plIinj}. On remarque tout d'abord que par
hypothèse sur $t$, la fonction $t$ vérifie bien les hypothèses du
lemme. En outre à chaque étape $\tau<T$, on applique de manière simple un
plan du type $Cat_{\Psi,T}$
qui par définition commence avec le plan simple $e_{\Psi,1}$. Or ce dernier contient tous les diagrammes obtenus à partir de $\Psi$,
où ici $\Psi$ n'est autre que la réunion des $\Phi_k$ pour $k$ décrivant
$t(\tau)$. Ainsi $BCat(t)$ vérifie bien les hypothèses du
corollaire~\ref{plIinj}, ce qui montre la première propriété.\\
CQFD.\\

En plus de formaliser le plan Bigcat qui catégorise tout en préservant l'homotopie des \cats, les plans de type $BCat(t)$ ont la propriété de se rationaliser en $HCat$,
ce qui aura pour conséquence de les rendre stables par précomposition avec un
de leurs sous-plans.

\begin{prop}\label{bcatrat}
Soient $\mathcal{C}$ une catégorie cocomplète et $\Phi$ un ensemble de
morphismes de $\mathcal{C}$ dont les sources sont $\mu$-petites, pour un certain cardinal $\mu$, et tel que les $\Phi$-cofibrations sont des monomorphismes. Soit $(\Phi_{k})_{i\in K}$ une partition de $\Phi$. 
Soient $T$ un ordinal transfini limite dont le cardinal est un cardinal régulier
strictement supérieur à $\mu$ et $t:T\rightarrow \mathcal{P}(K)$
une fonction telle que, pour tout $k$ dans $K$, on ait les propriétés suivantes : 
\item - pour tout ordinal $\tau<T$, il existe un ordinal
$\tau\leq\tau'<T$ tel que $t(\tau)$ contienne $k$,
\item - le cardinal de l'ensemble $t^{-1}(k)$ est égal au cardinal de $T$.\\
\\ 
Alors le plan infini $BCat(t)$ a pour rationalisation
$HCat_{T^2}$.
\end{prop}
{\it Preuve :}\\
D'après le corollaire~\ref{plrat}, BCat(t), dont la suite associée sera
notée $(X_{\beta})_{\beta\leq T^2}$, a pour rationalisation le plan
$Q_.$ dont la suite associée sera notée $(Z_{\beta})_{\beta\leq T^2}$
et tel que $Q_{\beta}$ est la colimite des plans $(BCat_{\alpha})_{\beta}$,
où $(BCat_{\alpha})_{\beta}$ est l'ensemble des
couples $(D',\alpha_{D'})$ tel que $D'$ se factorise par $Z_{\beta}$ mais pas
par $Z_{\gamma}$ pour $\gamma<\beta$, où les $D'$ sont les extensions par
$X_{\alpha}\rightarrow Z_{\alpha}$ des diagrammes $D$ de $BCat(t)_{\alpha}$ et où
$\alpha_{D'}$ est la cardinalité de $D$ dans $BCat(t)_{\alpha}$.\\ 
\\
Montrons que, pour tout $\beta<T^2$, $Q_{\beta}$ est isomorphe à la
$\beta$-ième étape de $HCat_{T^2}$ dans $X/\mathcal{C}$.
Pour cela, explicitons $Q_{\beta}$. On remarque
tout de suite que, par définition, $(BCat_{\alpha})_{\beta}$ est vide pour
$\alpha<\beta$.
Prenons sur $T^2$ l'ordre lexicographique, la première
coordonnée représentant la source de $t$. On remarque que si
$\alpha=(\gamma,\gamma')\geq\beta$, $BCat(t)_{\alpha}$
n'est autre qu'un sous-plan de $e_{\cup_{k\in t(\gamma)}\Phi_{k},1}$. Ainsi dans ce
cas, $(BCat_{\alpha})_{\beta}$ est l'ensemble des
couples $(D',1)$ tel que $D'$ se factorise par $Z_{\beta}$ mais pas
par $Z_{\delta}$ pour $\delta<\beta$, où les $D'$ sont les extensions par
$X_{\alpha}\rightarrow Z_{\alpha}$ des diagrammes $D$ de
$\cup_{k\in t(\gamma)}\Phi_{k}$ pris à l'étape $BCat(t)_{\alpha}$.
Ainsi l'ensemble sous-jacent des diagrammes de $Q_{\beta}$, ensemble qu'on
notera $A$, est un ensemble de diagrammes de $\Phi$ à valeur dans $Z_{\beta}$
mais pas dans $Z_{\delta}$, pour $\delta<\beta$, donc un sous-ensemble de
l'ensemble $B$ de tous les diagrammes de $\Phi$ à valeur dans $Z_{\beta}$
mais pas dans $Z_{\delta}$, pour $\delta<\beta$.\\
\\
Montrons qu'en fait $A$ et $B$
sont en bijection. Considérons donc un diagramme $D$ quelconque de $\Phi$ à
valeur dans $Z_{\beta}$ mais ne se factorisant pas par $Z_{\delta}$, pour
$\delta<\beta$. Alors le diagramme $D$ peut s'étendre à $Z_{T^2}$ qui
est isomorphe à $BCat(t)(X)=X_{T^2}$. Or par hypothèse, les sources des morphismes de $\Phi$ sont
$\mu$-petites et $T$ est un ordinal limite de cardinal un cardinal régulier strictement supérieur à $\mu$, alors il
existe $\alpha<T^2$ tel que l'extension de $D$ à $X_{T^2}$ se
factorise par $X_{\alpha}$. Comme $\Phi$ admet une partition en $\Phi_k$ pour
$k\in K$, il existe un $k\in K$ tel que $D$ soit un diagramme de $\Phi_k$. 
En outre on peut écrire $\alpha$ sous la forme $(\gamma,\gamma')$. Or par
hypothèse sur $t$, il existe $\gamma\leq\epsilon<T$ tel que
$t(\epsilon)$ contienne $k$. Considérons alors
$BCat(t)_{(\epsilon, 0)}$, ce plan étant de la forme $e_{\Psi,1}$ contient tous les diagrammes de
$\Phi_k$ à valeur dans $X_{(\epsilon, 0)}$ donc en particulier ceux à
valeur dans $X_{\alpha}$, car $\alpha\leq(\epsilon,0)$ par ordre lexicographique, et par conséquent ce plan contient l'extension de
$D$ à $X_{T^2}$. Ceci montre bien que $D$ appartient à
$(BCat(t)_{(\epsilon, 0)})_{\beta}$ et donc à $Q_{\beta}$. On a donc
montré que les ensembles $A$ et $B$ sont bien en bijection.\\
\\
Comme ce sont les ensembles sous-jacents des diagrammes respectivement de
$Q_{\beta}$ et de la $\beta$-ième étape de $HCat_{T^2}$ qui sont en bijection, pour montrer
que $Q_{\beta}$ et la $\beta$-ième étape de $HCat_{T^2}$ sont
isomorphes, il suffit de montrer que les cardinaux des diagrammes sont les
mêmes. Comme le cardinal de tout diagramme de la $\beta$-ième étape de
$HCat_{T^2}$ est $Card(T^2)=Card(T)$, il ne reste plus qu'à
montrer que le cardinal de chaque diagramme de $Q_{\beta}$ est
$Card(T)$.\\
\\
Soit donc $D$ un diagramme de $Q_{\beta}$, ce diagramme appartient à un
certain $(BCat(t)_{\alpha})_{\beta}$, avec $\alpha$ minimal. Si on écrit $\alpha$ sous la forme
$(\gamma,\gamma')$, ce diagramme est l'extension par $X_{\alpha}\rightarrow
Z_{\alpha}$ d'un diagramme $D'$ appartenant à $\cup_{k\in t(\gamma)}\Phi_{k}$
à valeur dans $X_{\alpha}$. Notons $k$ l'indice de la
partie de $\Phi$ à laquelle appartient $D'$. Alors $D'$ appartient au moins à tous les plans $BCat(t)_{\alpha'}$ tel que $\alpha'\geq\alpha$ et
$\alpha'=(\epsilon,0)$ avec $\epsilon$ décrivant $t^{-1}(k)$ et au plus à tous
les plans $BCat(t)_{\alpha'}$ tel que $\alpha'\geq\alpha$ et
$\alpha'=(\epsilon,\epsilon')$ avec $\epsilon$ décrivant $t^{-1}(k)$ et
$\epsilon'$ quelconque. En effet ces plans sont des sous-plans quelconques d'un plan simple de type $e_{\Psi,1}$, plan simple prenant une seule fois tous les
diagrammes de $\Phi_k$ à valeur dans $X_{\alpha'}$ donc en particulier à
valeur dans $X_{\alpha}$ car $\alpha\leq\alpha'$. Ainsi le cardinal associé
à $D$ dans $Q_{\beta}$ est compris entre celui des ordinaux $\epsilon>\gamma$ de $t^{-1}(k)$, qui par hypothèse sur $t$ n'est autre que
$Card(T)$, et celui des ordinaux $\alpha'=(\epsilon,\epsilon')$ avec $\epsilon\geq \gamma$ décrivant $t^{-1}(k)$ et $\epsilon'$ quelconque, qui par hypothèse sur $t$ n'est autre que $Card(T^2)=Card(T)$.\\
CQFD.\\

Comme on l'a dit plus haut, la stabilité de ces plans $BCat(t)$ vis-à-vis de
la composition avec leurs sous-plans nous intéressent particulièrement lorsqu'on veut faire
des colimites d'objets $\Phi$-injectifs et comparer la colimite des
$\Phi$-injectivisés avec la $\Phi$-injectivisation de la colimite. Nous allons
donc traiter de ce problème dans la section suivante.

\newpage

\section{Application de l'I-injectivisation aux colimites d'objets I-injectifs.}

Nous nous intéressons maintenant à comparer l'I-injectivisation des colimites
avec les colimites des objets I-injectivés. En fait nous aimerions en
particulier savoir s'il y a une I-\eq~au sens du corollaire~\ref{Ephi3'} entre
les deux. Pour cela, nous allons supposer que le plan d'addition de cellules
vérifiant les propriétés i) et ii) de la I-injectivisation ainsi que la
propriété iii) pour les objets I-injectifs est de la forme $BCat(t)$. En
effet, dans ce cas le morphisme naturel d'un objet quelconque $A$ dans le résultat du
plan sera une I-\eq. De là et du fait que $BCat(t)$ se rationalise en $HCat$,
il découlera que tout morphisme naturel d'un objet $A$ quelconque dans le
résultat d'un sous-plan de $BCat(t)$ est encore une I-\eq. Les colimites
des I-injectivisés étant des sous-plans de l'I-injectivisation des
colimites, on pourra leur appliquer ce résultat. C'est pourquoi dans toute
cette section, nous allons supposer les hypothèses suivantes.

\begin{hyp}\label{colimIinj}
Soient $\mathcal{C}$ une catégorie cocomplète et $\Phi$ un ensemble de
morphismes de $\mathcal{C}$ dont les sources sont $\alpha$-petites pour un certain cardinal régulier
$\alpha$, et tel que les $\Phi$-cofibrations sont des monomorphismes. Soit $\alpha'$ le plus petit cardinal régulier supérieur à
$\alpha$. Donnons-nous une notion d'\eq~d'objets
$\Phi$-injectifs vérifiant les propriétés suivantes :
\item - pour tout couple $(f,g)$ de morphismes composables entre objets
$\Phi$-injectifs, si parmi $f,g,g\circ f$ deux morphismes sont des \eqs~d'objets
$\Phi$-injectifs, alors le troisième morphisme aussi,
\item - pour tout couple $(f,g)$ de morphismes composables entre objets
$\Phi$-injectifs, si la composée $g\circ f$ est l'identité et que la
composée $f\circ g$ une \eq~d'objets $\Phi$-injectifs, alors $f$ et $g$ sont
des \eqs~d'objets $\Phi$-injectifs.
\item - les isomorphismes entre objets $\Phi$-injectifs sont des \eqs~d'objets
$\Phi$-injectifs.\\
\\
Notons $E_{\Phi}$ le plan d'addition de cellules fonctoriel à compositions rationnelles de
longueur $\alpha'$ dont tous les plans simples sont $e_{\Phi,1}$.\\
\\
Soit $(\Phi_{k})_{i\in K}$ une partition de $\Phi$. 
Soient $T$ un ordinal transfini limite dont le cardinal est
strictement supérieur à $\alpha$ et $t:T\rightarrow \mathcal{P}(K)$
une fonction telle que, pour tout $k$ dans $K$, on ait les propriétés suivantes : 
\item - pour tout ordinal $\tau<T$, il existe un ordinal
$\tau\leq\tau'<T$ tel que $t(\tau)$ contienne $k$,
\item - le cardinal de l'ensemble $t^{-1}(k)$ est égal au cardinal de $T$.\\
\\ 
Supposons que le plan d'addition de cellules fonctoriel $BCat(t)$ suivant
l'ordre $t$ vérifie la propriété suivante :\\
pour tout objet $A$ $\Phi$-injectif, le morphisme naturel $A\rightarrow
BCat(t)(A)$ est une \eq~d'objets $\Phi$-injectifs. 
\end{hyp}

Sous ces hypothèses, le plan d'addition de cellules $E_{\Phi}$ est bien une $\Phi$-injectivisation par la
proposition~\ref{Ephi3}, ce qui permet de définir la notion de $\Phi$-\eq~qui
vérifie la propriété du "trois pour deux" d'après le
corollaire~\ref{Ephi3'}. De plus sous ces hypothèses, par le lemme~\ref{bcat}, le plan d'addition de
cellules $BCat(t)$ vérifie bien les propriétés i) et ii) de la
$\Phi$-injectivisation.\\
Montrons tout d'abord que sous ces hypothèses, pour tout objet $A$ de la
catégorie $\mathcal{C}$, le morphisme naturel $A\rightarrow BCat(t)(A)$ est
une $\Phi$-\eq.

\begin{lem}\label{bcateq}
Sous les hypothèses~\ref{colimIinj}, pour tout objet $A$ de la catégorie
$\mathcal{C}$, le morphisme naturel $A\rightarrow BCat(t)(A)$ est une
$\Phi$-\eq, c'est-à-dire que l'image par $E_{\Phi}$ de ce morphisme est une
\eq~d'objets $\Phi$-injectifs.
\end{lem}
{\it Preuve :}\\
Tout d'abord on remarque que sous les hypothèses~\ref{colimIinj}, on peut
appliquer la proposition~\ref{bcatrat} à $BCat(t)$, ce qui nous permet
d'obtenir que $BCat(t)$ se rationalise en $HCat_{T^2}$. Or le plan à compositions rationnelles d'addition de
cellules $E_{\Phi}$ est un sous-plan de $HCat_{T^2}$ donc de $BCat(t)$ de
longueur $\alpha'$ avec $Card(\alpha')\leq Card(T^2)$. On peut alors appliquer le
corollaire~\ref{ratplcomp}, qui nous donne pour tout objet $A$ de la catégorie
$\mathcal{C}$ le diagramme commutatif suivant :
\begin{diagram}
A & \rTo & E_{\Phi}(A)\\
\dTo & & \dTo\\
BCat(t)(A) & \rTo^{\cong} & BCat(t)(E_{\Phi}(A))\\
\end{diagram} 
Le corollaire~\ref{ratplcomp} nous certifie la commutativité de ce diagramme
ainsi que l'isomorphisme entre $BCat(t)(A)$ et $BCat(t)(E_{\Phi}(A))$. En outre
par la proposition~\ref{Ephi3}, le morphisme naturel $A\rightarrow E_{\Phi}(A)$
est une $\Phi$-\eq~et $E_{\Phi}(A)$ est $\Phi$-injectif. Par hypothèse sur
$BCat(t)$, il vient que le morphisme naturel $E_{\Phi}(A)\rightarrow
BCat(t)(E_{\Phi}(A))$ est une \eq~d'objets $\Phi$-injectifs, ce qui par le
corollaire~\ref{Ephi3'} revient à dire que c'est une $\Phi$-\eq. Comme par
corollaire~\ref{Ephi3'}, les isomorphismes sont des $\Phi$-\eqs~et les
$\Phi$-\eqs~vérifient la propriété du "trois pour deux", il vient que le
morphisme naturel $A\rightarrow BCat(t)(A)$ est une $\Phi$-\eq.\\
CQFD.\\

Munis de cette propriété de stabilité homotopique de $BCat(t)$, montrons
que tout sous-plan de $BCat(t)$ conserve cette stabilité homotopique.

\begin{lem}\label{splbcat}
Supposons vraies les hypothèses~\ref{colimIinj}. Soit $P$ un sous-plan de $BCat(t)$
dont le cardinal de la longueur est inférieur au cardinal de la longueur de
$BCat(t)$. Alors pour tout objet $A$ de la catégorie $\mathcal{C}$, le
morphisme naturel $A\rightarrow P(A)$ est une $\Phi$-\eq.
\end{lem}
{\it Preuve :}\\
Comme on l'a vu dans la démonstration précédente, le plan d'addition de
cellules $BCat(t)$ se rationalise en un $HCat$. Comme en outre, par hypothèse
sur le sous-plan $P$ de $BCat(t)$, le cardinal de la longueur du sous-plan $P$ est
inférieur à celui de la longueur de $BCat(t)$, le corollaire~\ref{ratplcomp}
s'applique et nous donne un
isomorphisme dans $X/\mathcal{C}$ entre $BCat(t)(A)$ et $BCat(t)(P(A))$ ainsi que le
diagramme commutatif suivant :
\begin{diagram}
A & \rTo & P(A)\\
\dTo_{\sim} & & \dTo^{\sim}\\
BCat(t)(A) & \rTo^{\cong} & BCat(t)(P(A))\\
\end{diagram}
En outre, par le lemme précédent, les
morphismes canoniques $A\rightarrow BCat(t)(A)$ et $P(A)\rightarrow BCat(t)(P(A))$ sont des
$\Phi$-\eqs. Comme, par le corollaire~\ref{Ephi3'}, les $\Phi$-\eqs~vérifient
la propriété du "trois pour deux" et que les isomorphismes sont des
$\Phi$-\eqs, on obtient que le morphisme $A\rightarrow P(A)$ est aussi une
$\Phi$-\eq.\\
CQFD.\\

Comme le fait de prendre la colimite d'objets $\Phi$-injectivisés est en fait
un sous-plan de la $\Phi$-injectivisation de la colimite, on peut donc appliquer
le résultat que l'on vient de montrer afin de donner un critère utile pour
montrer que les morphismes naturels des colimites sont des $\Phi$-\eqs.

\begin{lem}\label{colimrec}
Supposons vraies les hypothèses~\ref{colimIinj}.
\item Soit $F:\mathcal{I}\rightarrow \mathcal{C}$ un foncteur d'une petite catégorie
d'indices $\mathcal{I}$ vers la catégorie $\mathcal{C}$. Soit $i$ un objet de
$\mathcal{I}$. Alors le morphisme naturel $F(i)\rightarrow
\colimite{i\in\mathcal{I}} F(i)$ est une $\Phi$-\eq~si et seulement si le morphisme
naturel $E_{\Phi}(F(i))\rightarrow \colimite{i\in\mathcal{I}} E_{\Phi}(F(i))$ en
est une.
\item Soit $h:F\rightarrow G$ une transformation naturelle entre foncteurs d'une
petite catégorie d'indices $\mathcal{I}$ vers la catégorie $\mathcal{C}$.
Alors le morphisme universel, induit par $h$, $\colimite{i\in\mathcal{I}}
F(i)\rightarrow \colimite{i\in\mathcal{I}}G(i)$, est une $\Phi$-\eq~si et seulement
si le morphisme universel, induit par $E_{\Phi}(h)$, $\colimite{i\in\mathcal{I}}
E_{\phi}(F(i))\rightarrow \colimite{i\in\mathcal{I}}E_{\Phi}(G(i))$ en est une.
\end{lem}
{\it Preuve :}\\
Considérons le diagramme commutatif suivant :
\begin{diagram}
F(i) & \rTo & \colimite{i\in\mathcal{I}} F(i)\\
\dTo^{\sim} & & \dTo_{\sim}\\
E_{\Phi}(F(i)) & \rTo & \colimite{i\in\mathcal{I}} E_{\Phi}(F(i))\\
\end{diagram}
Par la proposition~\ref{Ephi3}, la flèche verticale de gauche est une
$\Phi$-\eq. En outre la flèche verticale de droite consiste à
additionner des cellules de $\Phi$ sur chacun des $F(i)$, or ceci est un
sous-plan de $BCat(t)$ appliqué à $\colimite{i\in\mathcal{I}} F(i)$. Comme
en outre le cardinal de la longueur de $E_{\Phi}$ est inférieur à celui de
$BCat(t)$, on peut appliquer le lemme~\ref{splbcat} qui nous certifie que la
flèche verticale de droite est une $\Phi$-\eq. Par le corollaire~\ref{Ephi3'},
les $\Phi$-\eqs~vérifient la propriété du "trois pour deux", ce qui nous
donne que la flèche horizontale du haut est une $\Phi$-\eq~si et seulement si
la flèche horizontale du bas l'est.\\
\\
Considérons maintenant le diagramme commutatif suivant :
\begin{diagram}
\colimite{i\in\mathcal{I}} F(i) & \rTo^{\colimite{} h} & \colimite{i\in\mathcal{I}} G(i)\\
\dTo^{\sim} & & \dTo_{\sim}\\
\colimite{i\in\mathcal{I}} E_{\Phi}(F(i)) & \rTo_{\colimite{} E_{\Phi}(h)} & \colimite{i\in\mathcal{I}}
E_{\Phi}(G(i))\\
\end{diagram}
Par le même argument de sous-plan que dans la première partie de la preuve, on montre que
les deux flèches verticales sont des $\Phi$-\eqs, et on en conclut, avec l'aide
de la propriété du "trois pour deux" des $\Phi$-\eqs, que la flèche horizontale du haut est une $\Phi$-\eq~si et seulement si
la flèche horizontale du bas l'est.\\
CQFD.\\

Ce lemme de reconnaissance du caractère $\Phi$-\eq~des morphismes
naturels des colimites nous sera très utile lorsque dans la suite de cette
thèse, nous montrerons la stabilité des cofibrations triviales par somme
amalgamée le long d'un morphisme et par colimite séquentielle transfinie.
Nous allons donc terminer ce chapitre par quelques petits lemmes techniques sur la
$\Phi$-injectivisation qui seront également nécessaires pour la suite de
cette thèse.

\newpage

\section{Lemmes techniques sur l'I-injectivisation}

Dans la suite de cette thèse, nous aurons à montrer des résultats qui
nécessitent de connaître les propriétés de stabilité de
l'I-injectivisation par rapport aux monomorphismes, aux intersections de sous
objets et à la petitesse.

\begin{lem}\label{monoadd}
Soient $\mathcal{C}$ une catégorie cocomplète et $\Phi$ un ensemble de
morphismes de $\mathcal{C}$ tel que les $\Phi$-cofibrations sont les monomorphismes. Soient $\Psi$ un sous-ensemble de $\Phi$ et
$\lambda$ un cardinal quelconque (transfini ou non). Soit enfin $X\rightarrow Y$ un
monomorphisme de $\mathcal{C}$. Alors $e_{\Psi,\lambda}(X\rightarrow Y)$ est
encore un monomorphisme.\\ 
De même, pour tout plan d'addition $P_.$ dont les
plans simples sont du type $e_{\Psi,\lambda}$, $P_.(X\rightarrow Y)$ est un
monomorphisme.
\end{lem}
{\it Preuve :}\\
Tout d'abord, il est facile de voir que comme $X\rightarrow Y$ est un
monomorphisme, l'ensemble des diagrammes de $\Psi$ à valeur
dans $X$ s'injecte par extension à travers $X\rightarrow Y$ dans l'ensemble
des diagrammes de $\Psi$ à valeur dans $Y$. De plus, comme par hypothèse les $\Phi$-cofibrations sont les monomorphismes, on obtient d'une part que les flèches de $\Phi$ sont des monomorphismes et que d'autre part les monomorphismes sont stables par somme amalgamée le long d'un morphisme et colimite séquentielle transfinie.\\
Considérons maintenant la somme amalgamée, notée $Z$, de $e_{\Psi,\lambda}(X)$ avec $Y$ au-dessus
de $X$. Tout d'abord par commutation des sommes amalgamées, il vient que $Z$
est isomorphe à la somme amalgamée de $Y$ par l'ensemble des extensions le
long de $X\rightarrow Y$ des diagrammes de $\Psi$ à valeur dans $X$ pris
chacun $\lambda$ fois. Le morphisme canonique de $e_{\Psi,\lambda}(X)$ vers la somme
amalgamée $Z$ est un monomorphisme comme somme amalgamée du monomorphisme
$X\rightarrow Y$ le long du morphisme $X\rightarrow e_{\Psi,\lambda}(X)$. En
utilisant les remarques précédentes, on obtient que le morphisme universel de la
somme amalgamée $Z$ dans $e_{\Psi,\lambda}(Y)$ est en fait le résultat d'un
plan d'addition par les diagrammes de $\Psi$ à valeur dans $Y$ pris
$\lambda$-fois qui ne sont pas extensions par $X\rightarrow Y$ de diagrammes de
$\Psi$ à valeur dans $X$. Comme $\Phi$ est un ensemble de monomorphismes, le morphisme cellulaire
$Z\rightarrow e_{\Psi,\lambda}(Y)$ est un monomorphisme comme somme amalgamée
de monomorphismes. On en déduit que le morphisme
$e_{\Psi,\lambda}(X\rightarrow Y)$, composée des monomorphismes
$e_{\Psi,\lambda}(X)\rightarrow Z$ et $Z\rightarrow e_{\Psi,\lambda}(Y)$, est
bien un monomorphisme.\\
En utilisant ce résultat ainsi que la préservation des monomorphismes par
colimite séquentielle transfinie le long de colimites séquentielles transfinies de monomorphismes, on
obtient le résultat pour $P_.$.\\
CQFD.

\begin{lem}\label{interpl}
Soit $\mathcal{C}$ une catégorie cocomplète telle que les réunions de deux
sous-objets sont exactement les sommes amalgamées de ces sous-objets au-dessus
de leur intersection. Soit $\Phi$ un ensemble de
morphismes de $\mathcal{C}$ tel que les $\Phi$-cofibrations sont des monomorphismes. Soient $\Psi$ un sous-ensemble de $\Phi$ et
$\lambda$ un cardinal quelconque (transfini ou non). Soient $A$ et $B$ deux objets de
$\mathcal{C}$.
Alors pour tout plan d'addition $P_.$ dont les
plans simples sont du type $e_{\Psi,\lambda}$, on a : 
$$P(A\cap B)=P(A)\cap P(B)$$
\end{lem}
{\it Preuve :}\\
Par le lemme~\ref{monoadd}, les monomorphismes sont stables par tout plan $P_.$
d'addition de cellules composé de plans simples du type $e_{\Psi,\lambda}$.
Donc $P(A)$ et $P(B)$ sont des sous-objets de $P(C)$ et donc leur
intersection a un sens. De plus $P(A\cap B)$ est un sous-objet à la fois de
$P(A)$ et de $P(B)$, car $A\cap B$ est inclus dans $A$ et $B$ et que le plan
$P_.$ préserve les monomorphismes. Ainsi $P(A\cap B)$ est un
sous-objet de $P(A)\cap P(B)$. Il ne reste donc plus qu'à montrer l'autre
sens. Comme le plan $P_.$ est la colimite des plans simples $e_{\Psi,\lambda}$, il suffit
de montrer que $e_{\Psi,\lambda}(A\cap B)$ contient $e_{\Psi,\lambda}(A)\cap
e_{\Psi,\lambda}(B)$. La construction $e_{\Psi,\lambda}$ est un
procédé d'addition de cellules. Donc pour
montrer que $e_{\Psi,\lambda}(A)\cap e_{\Psi,\lambda}(B)$ est inclus dans
$e_{\Psi,\lambda}(A\cap B)$, il suffit de
montrer que les cellules du premier sont en fait des cellules du second. Or cela
découle directement du fait que les diagrammes de flèches de $\Phi$
apparaissant pour $A$ n'apparaissent pour $B$ que si les sources des
flèches de $\Phi$ s'envoyent dans $A\cap B$. Ainsi on a montré que la
construction $e_{\Psi,\lambda}$ préserve l'intersection et donc par récurrence transfinie
on pourra montrer que le plan $P_.$ aussi préserve l'intersection.\\
CQFD.

\begin{lem}\label{petitpl}
Soient $\mathcal{C}$ une catégorie cocomplète et $\Phi$ un ensemble de
flèches de $\mathcal{C}$. Supposons qu'il existe un cardinal transfini
$\alpha$ pour lequel les sources et buts de toutes les flèches de $\Phi$ sont $\alpha$-petites.
Soient $\lambda$ et $\lambda'$ deux ordinaux quelconques et $P_.$ un plan d'addition de cellules
fonctoriel de longueur $\lambda$ et uniquement composé de plans simples du
type $e_{\Psi,\lambda'}$ avec $\Psi$ un sous-ensemble de $\Phi$.
Soit $\alpha'$ un cardinal régulier strictement supérieur à
$2^{\lambda}$ et $2^{\alpha}$. Alors il existe un cardinal régulier $\alpha''$ strictement
supérieur à $\alpha'$ tel que, pour tout
objet $\alpha'$-petit $A$ de $\mathcal{C}$, $P(A)$ est $\alpha''$-petit.
\end{lem}
{\it Preuve :}\\
On rappelle que $P(A)$ est la colimite d'une $\lambda$-séquence
$(A^{\beta})_{\beta\leq\lambda}$ où $A^{\beta+1}=e_{\Psi,\lambda'}(A^{\beta})$. Montrons par
récurrence transfinie sur $\beta\leq\lambda$ que $A^{\beta}$ est
$\alpha^{(\beta)}$-petit pour un cardinal régulier $\alpha^{(\beta)}>\alpha'$
qui ne dépend pas de $A$. Le cas $\beta=0$ est vrai car $A^0$ n'est autre que $A$ qui est
par hypothèse $\alpha'$-petit. Supposons le cas $\beta$ montré et montrons
le cas $\beta+1$. Comme $A^{\beta+1}=e_{\Psi,\lambda}(A^{\beta})$ et que, par hypothèse de
récurrence, $A^{\beta}$ est $\alpha^{(\beta)}$-petit, il suffit de montrer que
le plan simple $e_{\psi,\lambda'}$ préserve la petitesse. Or
$e_{\Psi,\lambda'}(A^{\beta})$ consiste en une somme
amalgamée de $A^{\beta}$ par un coproduit de flèches de $\Psi$,
sous-ensemble de $\Phi$, chacune de ces flèches étant prise $\lambda'$ fois. Ce coproduit est
petit car la sous-famille $\Psi$ de la famille $\Phi$ est un ensemble par
hypothèse et que la collection des morphismes à but fixé et à
source indexée par un ensemble est un ensemble. Par hypothèse sur $\Phi$,
les sources et buts des flèches de $\Phi$
sont $\alpha$-petits, donc en particulier $\alpha^{(\beta)}$-petits car
$\alpha^{(\beta)}>\alpha'>\alpha$. Ainsi $e_{\Psi,\lambda'}(t^{\beta})$ est la colimite d'un petit diagramme
dont tous les objets sont $\alpha^{(\beta)}$-petits. D'après le
lemme~\ref{colimpetit}, une telle colimite est $\alpha^{(\beta+1)}$-petite pour un cardinal
régulier $\alpha^{(\beta+1)}$ strictement plus grand que $\alpha^{(\beta)}$ et
que le cardinal de l'ensemble des morphismes de la petite catégorie indexant
le diagramme. Or ce cardinal dépend de $\alpha^{(\beta)}$, du cardinal de
l'ensemble $\Psi$ et du cardinal de $\lambda'$, ces trois cardinaux étant indépendants de $A$, ce qui montre que
$e_{\Psi,\lambda'}(A^{\beta})=A^{\beta+1}$ est bien $\alpha^{(\beta+1)}$-petit avec
$\alpha^{(\beta+1)}>\alpha'$ cardinal régulier indépendant de $A$. Ainsi est prouvée
l'hypothèse de récurrence au rang $\beta+1$.\\ 
Supposons maintenant
l'hypothèse vraie pour tout ordinal $\gamma<\beta$ avec $\beta\leq\lambda$
ordinal limite et
montrons qu'elle est aussi vraie pour $\beta$. Posons $\alpha^{(\beta)}$ le plus
petit cardinal régulier majorant l'ensemble $\{\alpha^{(\gamma)},
\gamma<\beta\}$, qui existe par propriété des cardinaux et qui est
indépendant de $A$ car $\beta$ et les $\alpha^{(\gamma)}$ le sont. Comme $A^{\beta}$ est la colimite
des $\alpha^{(\gamma)}$-petits $A^{\gamma}$ pour $\gamma<\beta$ et que
$\beta\leq\lambda$
est strictement inférieur à $\alpha^{(\beta)}$ car
$\alpha^{(\beta)}>\alpha^{(0)}=\alpha'>\lambda$, par
régularité de $\alpha^{(\beta)}$, il vient que la colimite $A^{\beta}$ est
$\alpha^{(\beta)}$ petite avec $\alpha^{(\beta)}>\alpha'$ cardinal régulier
indépendant de $A$, ce qui montre le cas $\beta$ ordinal limite. Par récurrence
transfinie, on a donc montré que, pour tout $\beta\leq\lambda$, $A^{\beta}$ est
$\alpha^{(\beta)}$-petit, avec $\alpha^{(\beta)}>\alpha'$ cardinal régulier
indépendant de $A$. Donc ceci est vrai en particulier pour $A^{\lambda}=P(A)$, ce
qui est le résultat attendu en prenant pour $\alpha''$ le cardinal
$\alpha^{(\lambda)}$.\\
CQFD.\\

Nous avons dans ce chapitre fait une mise au point sur la théorie des plans
d'addition de cellules qui nous a fourni les outils nécessaires pour
I-injectiviser. Si donc nous trouvons un ensemble I de flèches pour lequel les
objets I-injectifs sont des \cats, en appliquant les constructions
d'I-injectivisation exhibées dans ce chapitre, nous serons à même de
catégoriser nos \precats. Aussi le but du prochain chapitre est de trouver un
ensemble de morphismes de \precats~I pour lequel les \precats~I-injectives sont
des \cats.

\chapter{Engendrement des \cats}

\newpage

A la fin du premier chapitre, nous avons expliqué l'intérêt d'avoir un
procédé de catégorisation des \cats, mais nous avons également parlé
des difficultés d'en obtenir un. L'idée principale pour construire une
caté\-gorisation est de caractériser certaines \cats~en terme
d'objets I-injectifs. Autrement dit, on cherche à exhiber une famille de
morphismes de \precats~telle que les objets ayant la \prd~par rapport aux
morphismes de cette famille sont des \cats. Comme cette caractérisation des
\cats~par une propriété de relèvement par rapport à une famille de
morphismes n'est pas vraie en général, on va devoir s'intéresser à une
notion un peu plus forte de \cats~qui aura le bon goût de se relever par
rapport à certaines flèches. Ce faisant, nous allons faciliter la
construction d'un procédé de catégorisation.\\

Comme les \cats~sont définies à partir d'une donnée de Segal, on peut
aisément imaginer que si dans cette donnée de Segal, on a déjà des
notions d'\obcs~et d'\eqcs~d'\obcs~caractérisées par une propriété de
relèvement alors ce sera aussi le cas au niveau des \cats~et des \eqs~de
\cats. Aussi allons-nous dans un premier temps définir les données de Segal
proto-faciles pour lesquelles \obcs~et \eqcs~d'\obcs~sont caractérisés par
des propriétés de relèvement, et définir les notions de \cats~faciles et
d'\eqs~faciles à partir de telles données de Segal.\\

Afin de remonter au niveau des \precats~les propriétés de relèvements des
\obcs~et des \eqcs~d'\obcs~et ainsi montrer que nos \cats~faciles et nos \eqs~faciles
sont caractérisées par des propriétés de relèvement, nous allons dans
un second temps définir une construction $\Theta$ permettant de traduire les
diagrammes des \precats~en diagrammes dans $\mathcal{C}$ et réciproquement.
Dans un troisième temps, nous utiliserons cette construction pour montrer que
les \cats~faciles sont bien caractérisées par une propriété de
relèvement, ce qui nous conduira à la construction $Cat$ qui transformera
toute \precat~en \cat~facile.

\newpage

\section{\cats~faciles et équivalences faciles}

Comme on l'a écrit plus haut, les notions de \cat~et d'équivalence de \cats~n'ont
pas de caractérisation simple par propriété de relèvement. Nous allons donc introduire une notion de
\cat~facile et d'équivalence facile qui admettront une caractérisation par
relèvement vis-à-vis de certaines flèches. Mais comme on l'a vu précédemment
nous devons d'abord demander à la donnée de Segal $\mathcal{C}$ d'avoir une
famille d'\obcs~et une famille d'\eqcs~d'\obcs~caractérisées par des
propriétés de relèvement. Plus précisément, on va rajouter à la
donnée de Segal deux familles de morphismes, notées $\mathcal{F}_1$ et
$\mathcal{F}_2$, telles que les objets $\mathcal{F}_1$-injectifs soient des
\obcs~et les morphismes $\mathcal{F}_2$-injectifs des \eqcs~d'\obcs. C'est ce
qu'on appellera une donnée de Segal proto-facile. 

\begin{defin}\index{donnée de
Segal!proto-facile}\index{$\mathcal{F}_1$}\index{$\mathcal{F}_2$}
Une donnée de Segal proto-facile est une triplet
$(\mathcal{C},\mathcal{F}_1,\mathcal{F}_2)$ constitué
d'une donnée de Segal $\mathcal{C}$ et de deux familles $\mathcal{F}_1$ et
$\mathcal{F}_2$ de morphismes de $\mathcal{C}$ qui sont des ensembles satisfaisant les propriétés
suivantes :
\item 8) Les objets de $\mathcal{C}$ ayant la \prd~par rapport à
$\mathcal{F}_1$ sont des \obcs~et sont stables par produit fibré au-dessus
d'un objet discret.
\item 9) Les objets discrets de $\mathcal{C}$ ont la \prd~par rapport à
$\mathcal{F}_1$.
\item 10) Les morphismes qui ont la \prd~par rapport à
$\mathcal{F}_2$ et dont la source et le but sont des
\obcs~ayant la \prd~par rapport à $\mathcal{F}_1$ sont des \eqcs~d'\obcs~et le produit fibré, dans
la catégorie des morphismes, de deux tels morphismes au-dessus d'un objet
discret a encore la \prd~par rapport à $\mathcal{F}_2$.
\end{defin}

Les objets de $\mathcal{C}$ ayant la \prd~par rapport à
$\mathcal{F}_1$ sont appelés \obcs~faciles de $\mathcal{C}$ et les
morphismes entre \obcs~faciles ayant la \prd~par rapport à
$\mathcal{F}_2$ sont appelés \eqcs~faciles d'\obcs~faciles.
\index{objet régal!facile}
\index{\eqc~d'\obcs!facile}\\

On pourra remarquer au passage qu'il y a plusieurs choix possibles pour les
familles $\mathcal{F}_1$ et $\mathcal{F}_2$. Ainsi une même donnée de Segal
pourra donner lieu à plusieurs données de Segal proto-faciles.\\
\\

\begin{ex}\label{hypenssimp2}
Avec la donnée de Segal pour $\mathcal{ENSSIMP}$ de
l'exemple~\ref{hypenssimp1}, on obtient une donnée de Segal proto-facile en prenant la
famille vide pour $\mathcal{F}_1$ et,
pour $\mathcal{F}_2$, celle des inclusions des bords des simplexes standards dans ces
derniers, qui forme bien un ensemble. On a alors que tout ensemble simplicial est un \obc~facile
et que les \eqcs~faciles sont les équivalences faibles
qui sont aussi des fibrations de Kan.
\end{ex}

\begin{ex}
Toujours avec la donnée de Segal pour $\mathcal{ENSSIMP}$ de
l'exem\-ple~\ref{hypenssimp1}, on obtient une donnée de Segal proto-facile en prenant la
famille vide pour $\mathcal{F}_1$ et,
pour $\mathcal{F}_2$, la famille suivante : pour tout entier $n$, on prend d'une
part l'unique morphisme du vide vers le $n$-simplexe $\Delta[n]$ et
d'autre part l'unique morphisme du coproduit de deux exemplaires de $\Delta[n]$
vers $\Delta[n]$ induit par les identités. On a alors que tout ensemble simplicial est un \obc~facile
et que les \eqcs~faciles sont uniquement les isomorphismes d'ensembles
simpliciaux.\\
Si au lieu de prendre le vide pour $\mathcal{F}_1$, on prend la famille des
inclusions des cornes des simplexes dans ces derniers, alors on aura une
troisième donnée de Segal proto-facile avec pour \obcs~les ensembles simpliciaux fibrants de Kan.
\end{ex}

Ces notions d'\obc~facile et d'\eqc~facile d'\obcs~faciles
d'une part vont constituer pour $\mathcal{C}$ une nouvelle donnée de Segal,
comme le montre le lemme ci-dessous, mais d'autre part vont nous permettre
de définir les notions de \cat~facile et d'\eq~facile de \cats~faciles. 

\begin{lem}
Soit $(\mathcal{C},\mathcal{F}_1,\mathcal{F}_2)$ une donnée de Segal
proto-facile.\\
Alors la catégorie sous-jacente $\mathcal{C}$ avec
sa sous-catégorie pleine d'\obcs~faciles, sa classe d'\eqcs~faciles et son
foncteur $\tau_0$ vérifie toutes les propriétés d'une donnée de Segal.
\end{lem}
{\it Preuve}\\
La propriété 9 de $(\mathcal{C},\mathcal{F}_1,\mathcal{F}_2)$ assure que les objets discrets sont des
\obcs~faciles de $\mathcal{C}$ et le lemme~\ref{prst} prouve que tout
objet isomorphe à un \obc~facile est un \obc~facile, ce qui
nous donne la propriété 1. La propriété 8 de
$(\mathcal{C},\mathcal{F}_1,\mathcal{F}_2)$ assure que les \obcs~faciles
sont stables par produit fibré au-dessus d'un objet discret, ce qui montre la
propriété 2. De même la propriété 10 de
$(\mathcal{C},\mathcal{F}_1,\mathcal{F}_2)$ assure que
les \eqcs~faciles sont stables par produit fibré au-dessus d'un
objet discret dans la catégorie des morphismes, ce qui montre la seconde partie de la
propriété 3. Par le lemme~\ref{prst}, les
\eqcs~faciles sont stables par composition et les isomorphismes sont
des \eqcs~faciles, ce qui montre que la classe des \eqcs~faciles est
stable par composition et vérifie la première partie de la propriété 3.
Par propriété 10 de $(\mathcal{C},\mathcal{F}_1,\mathcal{F}_2)$, les \eqcs~faciles sont stables par fibre, ce
qui nous donne le premier sens de l'équivalence de la propriété 4. Considérons un
morphisme d'\obcs~faciles dont toutes les fibres sont des \eqcs~faciles. 
Le fait que $\mathcal{C}$ soit \discret~entraîne que tout diagramme de
relèvement de ce morphisme vis-à-vis de la famille $\mathcal{F}_2$ n'est
autre que le
coproduit de ses fibres. On obtient alors le relèvement attendu en prenant le
coproduit dans la catégorie des morphismes des relèvements de chaque fibre
(cf preuve du lemme~\ref{prst}), ce qui montre l'autre sens
de la propriété 4. La propriété 8 de $(\mathcal{C},\mathcal{F}_1,\mathcal{F}_2)$ assure que le foncteur $\tau_0$ est
bien défini sur la sous-catégorie pleine des \obcs~faciles de
$\mathcal{C}$, qui est une sous-catégorie pleine des \obcs~de
$\mathcal{C}$, et donc que $\tau_0$ se restreint à la catégorie des
\obcs~faciles en conservant les propriétés 5 et 6. Comme
par propriété 10 de $(\mathcal{C},\mathcal{F}_1,\mathcal{F}_2)$, les \eqcs~faciles sont des
\eqcs~d'\obcs, $\tau_0$ garde aussi la
propriété 7.\\
CQFD.\\

Comme on vient de le voir une donnée de Segal proto-facile engendre une
donnée de Segal avec les \obcs~faciles et les \eqcs~faciles. Pour cette
donnée de Segal, on a bien évidemment l'existence de notions de \cats~et d'\eqs~de
\cats, que l'on appellera respectivement \cats~faciles et \eqs~faciles de
\cats~faciles.

\begin{defin}
Soit $(\mathcal{C},\mathcal{F}_1,\mathcal{F}_2)$ une donnée de Segal proto-facile.\\
Une \precat~$A$ est une \cat~facile si pour tout entier $m$ strictement positif $A_m$
est un \obc~facile dans $\mathcal{C}$ et si pour tout entier $m$ supérieur ou
égal à deux le morphisme de Segal $\segal{A}{m}$ est une \eqc~facile
d'\obcs~faciles.
 Un morphisme de \cats~faciles est un morphisme de \precats.
\end{defin}
\index{\cat!facile}

\begin{defin}
Soit $(\mathcal{C},\mathcal{F}_1,\mathcal{F}_2)$ une donnée de Segal proto-facile.\\
Un morphisme de \cats~faciles $f:A\rightarrow B$ est une
équivalence facile de
\cats~faciles si $f_0:A_0\rightarrow B_0$ est surjective et, si pour tout couple $(x,y)$ d'objets de A, le morphisme $f_1(x,y):A_1(x,y)\rightarrow
B_1(f(x),f(y))$ est une \eqc~facile d'\obcs~faciles.
\end{defin}
\index{equivalence!facile de \cat}

Si ces \cats~faciles et ces \eqs~faciles de \cats~faciles ont vocation à être
caractérisées par des propriétés de relèvement, il faut au minimum que les
\cats~faciles soient en particulier des \cats~et que les \eqs~faciles soient en
particulier des \eqs~de \cats.

\begin{lem}\label{fac}
Soit $(\mathcal{C},\mathcal{F}_1,\mathcal{F}_2)$ une donnée de Segal proto-facile.
Les \cats~faciles sont des \cats~et les équivalences faciles
de \cats~faciles sont des équivalences de \cats. 
\end{lem}
{\it Preuve :} immédiat en utilisant les propriétés 8 et 10.\\

La notion de \cat~facile a été introduite essentiellement afin de
permettre de caractériser les \precats~qui sont des \cats~faciles en terme de
relèvement par rapport à une certaine famille de flèches, ce qui permettra
de construire une opération de catégorisation des \precats~basée sur ces
relèvements. En vue d'exhiber cette famille de flèches, nous allons
introduire une construction fonctorielle $\Theta$ qui permet d'exprimer les
diagrammes de \precats~en diagrammes de $\mathcal{C}$ et réciproquement. Ainsi
les caractérisations par relèvement des \obcs~faciles et des \eqcs~faciles
entraîneront de telles caractérisations pour les \cats~faciles et les
\eqs~faciles.

\newpage

\section{La construction $\Theta$}

Par le lemme de Yoneda, on sait que les morphismes entre ensembles simpliciaux
du $n$-simplexe $\Delta[n]$ vers un ensemble simplicial $A$ correspondent exactement aux
éléments de l'ensemble $A_n$. Comme nos \precats~sont des objets simpliciaux
sur la catégorie $\mathcal{C}$, on cherche une construction bifonctorielle
$\Theta$ respectant les colimites qui à un ensemble simplicial $E$ et à un objet
$C$ de $\mathcal{C}$ associe une \precat~et telle que les morphismes de
$\Delta[n]\Theta C$ vers une \precat~$A$ correspondent aux morphismes dans
$\mathcal{C}$ de $C$ vers $A_m$. Si $\mathcal{C}$ est la catégorie des
ensembles, en prenant pour $C$ le point $*$, on retrouve le résultat
découlant de Yoneda sur les ensembles simpliciaux.\\

L'idée pour construire $\Delta[n]\Theta C$ muni de cette propriété est de
considérer $\Delta[n]$ comme une \precat~discrète et de remplacer son $n$-simplexe non dégénéré par $C$ et de répercuter ce changement à toute
la \precat~par les applications de faces et de dégénérescences. Comme en
outre on désire que $*\Theta C$ commute aux colimites, nous allons définir
$\Theta$ ainsi : à un ensemble simplicial X et à un objet Y de
$\mathcal{C}$, on associe une
\precat~$X\Theta Y$ formée sur la \precat~discrète X avec les mêmes objets que X mais dans laquelle
chaque simplexe de X non dégénérescence d'un point est remplacé par Y.

\begin{defin}
Soit $\mathcal{C}$ une catégorie \discret~possédant les coproduits.
Soit $X$ un ensemble simplicial et $Y$ un objet de $\mathcal{C}$, on définit la
\precat~$X\Theta Y$ ainsi :
$$(X\Theta Y)_0=X_0 \mbox{ et, pour }n>0, \; (X\Theta Y)_n= X_0\coprod\Bigg(\coprod_{X_n\backslash
s^n(X_0)}Y\Bigg)$$
où $s^n:X_0\rightarrow X_n$ est induite par l'unique application de {\bf n}
vers {\bf 0}.\\
On notera $Y^x$ l'exemplaire de $Y$ associe à l'élément
$x\in X_n\backslash s^n(X_0)$.\\
Soit $f$ un morphisme de $\Delta$ de {\bf m} vers {\bf n}. Définissons le
morphisme induit\\ $f^*:(X\Theta Y)_n\rightarrow (X\Theta Y)_m$. Sur la composante
$X_0$ de $(X\Theta Y)_n$, $f^*$ induit l'identité. Pour la composante $Y^x$,
il y a deux cas. Si $X(f)(x)$ n'est pas dégénérescence d'un point, alors $f^*$ induit
l'identité entre $Y^x$ et $Y^{X(f)(x)}$. Sinon $X(f)(x)$ est la
dégénérescence d'un point $a\in X_0$ et alors $f^*$ induit l'application
constante à valeur $a$ de $Y^x$ vers $X_0$. On vérifie aisément que
$X\Theta Y$ ainsi défini est une \precat.
\end{defin}
\index{$*\Theta *$}

Nous allons maintenant faire de la construction $\Theta$ définie ci-dessus un
bifoncteur $*\Theta *:\mathcal{ENSSIMP}\times(\mathcal{C-PC})\rightarrow
\mathcal{C-PC}$.

\begin{defin}
Soit $\mathcal{C}$ une catégorie \discret~possédant les coproduits.
\item - Soient $X$ un ensemble simplicial et $g:Y\rightarrow Y'$ un morphisme de
\precats, on définit $X\Theta g:X\Theta Y\rightarrow X\Theta Y'$ de la manière suivante :
$$(X\Theta g)_0=Id_{X_0} \mbox{ et, pour }n>0, \; (X\Theta g)_n= Id_{X_0}\coprod\Bigg(\coprod_{X_n\backslash
s^n(X_0)}g\Bigg) .$$
\item - Soient $Y$ une \precat~et $f:X\rightarrow X'$ un morphisme d'ensembles
simpliciaux, on définit $f\Theta Y:X\Theta Y\rightarrow X'\Theta Y$ de la
manière suivante :\\
$(f\Theta Y)_0$ est l'application 
$f_0:X_0\rightarrow X'_0$ et, pour n strictement positif, on
définit $(f\Theta Y)_n$ sur chacune des composantes de $(X\Theta Y)_n$. Sur
$X_0$, $(f\Theta Y)_n$ est l'application $f_0$. Pour $x\in X_n\backslash
s^n(X_0)$ tel que $f_n(x)$ n'est pas dégénérescence d'un point, $(f\Theta
Y)_n$ est l'identité de $Y^x$ vers $Y^{f_n(x)}$. Pour $x\in X_n\backslash
s^n(X_0)$ tel que $f_n(x)$ est dégénérescence d'un point $a'\in X'_o$, $(f\Theta
Y)_n$ est l'application constante à valeur $a'$ de $Y^x$ dans $X'_0$.
\end{defin}

\begin{lem}
Soit $\mathcal{C}$ une catégorie \discret~possédant les coproduits.\\
$*\Theta *:\mathcal{ENSSIMP}\times(\mathcal{C-PC})\rightarrow
\mathcal{C-PC}$ est un bifoncteur préser\-vant les sommes amalgamées en
chacune de ses composantes.
\end{lem} 
{\it Preuve :}\\
Il est facile de vérifier que, pour tout ensemble simplicial X et pour tout
objet Y de $\mathcal{C}$, $X\Theta *$ et $*\Theta Y$ sont des foncteurs et que
l'on a la propriété suivante : pour tout $f:X\rightarrow X'$ morphisme
d'ensembles simpliciaux et pour tout $g:Y\rightarrow Y'$ morphisme de
$\mathcal{C}$, on a l'égalité $f\Theta Y'\circ X\Theta g=X'\Theta g\circ
f\Theta Y$.\\
La vérification de la préservation de la somme amalgamée en chacune des variables est immédiate
à partir de la définition, car la somme amalgamée des \precats~est la
somme amalgamée niveau par niveau.\\
CQFD.\\

A partir de maintenant nous allons essentiellement appliquer cette construction
$\Theta$ à deux types d'ensembles simpliciaux.

\begin{defin} 
\item - Pour tout entier $m$, on définit $\Delta[m]$ comme l'ensemble
simplicial représenté par {\bf m}.\index{$\Delta[m]$}
\item - Pour tout entier $m$ strictement positif, on définit $\Upsilon(m)$ comme
la somme amalgamée de $m$ exemplaires de $\Delta[1]$ au-dessus de $\Delta[0]$
en alternant but et source.\index{$\Upsilon(m)$} On note $i_m:\Upsilon(m)\rightarrow
\Delta[m]$ l'inclusion canonique.
\end{defin}

On peut se représenter $\Delta[m]$ comme un $m$-simplexe mais aussi comme une
catégorie engendrée par $m$ morphismes composables entre $m+1$ points
distincts. Et l'on peut se représenter $\Upsilon(m)$ comme une suite de m
morphismes composables entre m+1 points distincts. Ainsi $\Delta[m]$ est la
catégorie associée à $\Upsilon(m)$.\\

Nous allons maintenant montrer que le bifoncteur $\Theta$ a bien la propriété que l'on
voulait, à savoir que les morphismes de $\Delta[m]\Theta C$ vers une
\precat~$A$ correspondent aux morphismes dans $\mathcal{C}$ de $C$ vers $A_m$.

\begin{prop}\label{theta}
Soit $\mathcal{C}$ une catégorie \discret~possédant les coproduits et les
produits fibrés au-dessus d'un objet discret.
\item 1) Il existe un isomorphisme naturel en $X$ et en $A$ entre l'ensemble des morphismes de
\precats~de $\Delta[m]\Theta X$ vers $A$ et l'ensemble des morphismes dans
$\mathcal{C}$ de X vers $A_m$ tels que le morphisme composé\\ $X\rightarrow
A_m\rightarrow A_0\times\ldots\times A_0$ est constant (i.e. se factorise par
l'objet final de $\mathcal{C}$), où le morphisme $A_m\rightarrow
A_0\times\ldots\times A_0$ est induit par les m+1 applications sommets de {\bf
0} vers {\bf m}.
\item 2) Il existe un isomorphisme naturel en $X$ et en $A$ entre l'ensemble des morphismes de
\precats~de $\Upsilon[m]\Theta X$ vers $A$ et l'ensemble des morphismes dans
$\mathcal{C}$ de X vers $A_1\times_{A_0}\ldots\times_{A_0} A_1$ tels que le morphisme composé $X\rightarrow
A_1\times_{A_0}\ldots\times_{A_0} A_1\rightarrow A_0\times\ldots\times A_0$ est constant, où le morphisme $A_1\times_{A_0}\ldots\times_{A_0} A_1\rightarrow
A_0\times\ldots\times A_0$ est induit par les m+1 applications sommets de {\bf
0} vers $\mbox{\bf 1}\coprod_{\mbox{\bf 0}}\ldots\coprod_{\mbox{\bf 0}}\mbox{\bf 1}$.
\item 3) Les deux isomorphismes ci-dessus sont compatibles avec l'inclusion $i_m:\Upsilon(m)\rightarrow
\Delta[m]$.
\end{prop}
{\it Preuve :}\\
1) Tout d'abord, on remarque que $(\Delta[m]\Theta X)_n$ est le coproduit de
$\Delta[m]_0$ avec un exemplaire de X pour chaque morphisme de {\bf n} dans {\bf
m} ne se factorisant pas par {\bf 0}. En outre, pour se donner un morphisme de
$\Delta[m]\Theta X$ vers $A$, il faut se donner pour chaque entier n des
morphismes de $(\Delta[m]\Theta X)_n$ vers $A_n$ compatibles avec la structure
simpliciale. Or pour se donner un morphisme de $(\Delta[m]\Theta X)_n$ vers
$A_n$, il faut et il suffit de se donner un morphisme de $\Delta[m]_0$ vers
$A_n$ et un morphisme de $X^f$ dans $A_n$ pour tout morphisme $f$ de {\bf n} dans {\bf
m} ne se factorisant pas par {\bf 0}. Par respect de la structure simpliciale,
on obtient les carrés commutatifs suivants :
\begin{diagram}
X^{Id_{\mbox{m}}} & \rTo & A_m         & & & \Delta[m]_0             & \rTo & A_0        & & X^{Id_{\mbox{m}}} & \rTo & A_m\\
\dTo^{Id_X}           &      & \dTo_{A(f)} & & & \dTo^{Id_{\Delta[m]_0}} &      & \dTo_{s^n} & & \dTo^{*_i}            &      & \dTo_{A(\delta^i)}\\
X^f                   & \rTo & A_n         & & & \Delta[m]_0             & \rTo & A_n        & & \Delta[m]_0           & \rTo & A_0\\
\end{diagram}
où $*_i$ désigne l'application constante à valeur i et $\delta^i$
l'application i-ème sommet de {\bf 0} vers {\bf m} envoyant 0 sur i.\\
Ainsi le premier diagramme montre qu'un morphisme de $X^f$ dans $A_n$ est entièrement déterminé par $f$
et la donnée d'un morphisme de $X^{Id_{\mbox{m}}}$ dans $A_m$. Quant au
morphisme de $\Delta[m]_0$ vers $A_n$, le deuxième diagramme montre qu'il est entièrement déterminé par la
donnée d'un morphisme de $\Delta[m]_0$ vers $A_0$, c'est-à-dire par la
donnée de m+1 objets de A $(a_0,\ldots,a_m)$. En outre le troisième
diagramme montre le morphisme $X^{Id_{\mbox{m}}}$ dans $A_m$ composé avec
les morphismes induits par les m+1 sommets de {\bf m} est constant. Ainsi se donner un
morphisme de $\Delta[m]\Theta X$ vers $A$, c'est se donner un morphisme de $X$
dans $A_m$ tel que la composée avec le morphisme $A_m\rightarrow
A_0\times\ldots\times A_0$ est constant et a pour valeur le m+1-uplet d'objets
de $A$ images du morphisme de $\Delta[m]_0$ vers $A_0$. Un autre manière de le
formuler est de dire que la donnée d'un morphisme de $\Delta[m]\Theta X$ vers
$A$ n'est autre que la donnée d'un morphisme de $X$ vers $A_m(a_0,\ldots,a_m)$
pour un certain $(a_0,\ldots,a_m)$ m+1-uplet d'objets de $A$.\\
Avec ce qui précède, on peut aisément construire l'isomorphisme cherché
et montrer qu'il est naturel en $X$ et $A$.\\
\\
2) La démonstration est presque identique à celle du 1). Cette fois-ci, l'objet
$(\Upsilon(m)\Theta X)_n$ est le coproduit de
$\Upsilon(m)_0$ avec un exemplaire de X pour chaque morphisme de {\bf n} dans {\bf
m} se factorisant par l'un des morphismes de faces principales de {\bf 1} dans
{\bf m} qui à 0 et 1 associe i et i+i, pour i compris entre 0 et m-1. Pour des
raisons similaires à celles du 1), se donner un morphisme de
$\Upsilon(m)\Theta X$ vers $A$, c'est se donner $m$ morphismes de $X$ vers $A_1$,
un pour chaque face principale, et tel que leurs composées avec le morphisme
$A_1\rightarrow A_0\times A_0$, induit par source et but, sont constantes et ont
pour valeur les deux objets correspondants par l'application de $\Upsilon(m)_0$
vers $A_0$. Là encore, l'isomorphisme s'en déduit aisément et il est
immédiat de vérifier sa naturalité en $X$ et $A$.\\
\\
3) Ici aussi, il s'agit d'une simple vérification. Il faut juste remarquer que
l'on passe de la donnée de $X\rightarrow A_m$ pour $\Delta[m]\Theta
X\rightarrow A$ à celle de\\ $X\rightarrow A_1\times_{A_0}\ldots\times_{A_0}A_1$
pour $\Upsilon(m)\Theta X\rightarrow A$ en composant la première donnée par
le morphisme de Segal de A au cran m.\\
{\it CQFD.}\\

On peut reformuler la proposition ci-dessus en disant que, pour
tout objet $C$ de $\mathcal{C}$ et pour tout entier $n$, le foncteur des
\precats~vers les ensembles qui à une \precat~$A$ associe le coproduit, sur les
$m+1$-uplets $(a_0,\ldots,a_m)$ d'objets de $A$, des ensembles de morphismes de
$C$ vers $A_m(a_0,\ldots,a_m)$ est représenté par $\Delta[m]\Theta C$. De
même, celui qui à $A$ associe le coproduit, sur les
$m+1$-uplets $(a_0,\ldots,a_m)$ d'objets de $A$, des ensembles de morphismes de
$C$ vers $A_1(a_0,a_1)\times\ldots\times A_1(a_{m-1},a_m)$ est représenté
par $\Upsilon(m)\Theta C$.\\

Comme on l'a vu, le fait que les \precats~aient leurs ensembles d'objets discrets
entraîne une condition de constance sur les morphismes de $X$ vers $A_m$
pour qu'ils proviennent d'un morphisme de $\Delta[m]\Theta X$ vers $A$. Il en
est de même avec $\Upsilon(m)$. Afin de simplifier cette condition de
constance, donnons une condition sur $X$ pour laquelle tout morphisme de $X$ vers $A_m$ ou
$A_1\times_{A_0}\ldots\times_{A_0}A_1$ vérifie la constance sur le produit des
$A_0$ de la proposition~\ref{theta}. Cette condition sera une
généralisation de la connexité.

\begin{defin}\label{connex}\index{connexité}
Soit $\mathcal{C}$ une catégorie admettant un objet final. Un objet $X$ de
$\mathcal{C}$ est connexe s'il vérifie la propriété que tout morphisme de
$X$ vers un objet discret se factorise par l'objet final.
\end{defin}

\begin{lem}\label{cons}
Soit $\mathcal{C}$ une catégorie \discret~possédant les coproduits et les
produits fibrés au-dessus d'un objet discret. Soient $X$ un objet quelconque de
$\mathcal{C}$ et $A$ une \precat. Tout morphisme dans $\mathcal{C}$ de $X$ vers $A_m$ ou
$A_1\times_{A_0}\ldots\times_{A_0}A_1$ se factorisant à travers un objet
connexe vérifie la condition de la
proposition~\ref{theta} et donc correspond univoquement à un morphisme de $\Delta[m]\Theta
X$ ou de $\Upsilon(m)\Theta X$ vers
$A$.
\end{lem}
{\it Preuve :}\\
On étend le morphisme de $X$ vers $A_m$ ou
$A_1\times_{A_0}\ldots\times_{A_0}A_1$ à l'objet discret
$A_0\times\ldots\times A_0$ en composant par les morphismes induits par les
applications sommets. Comme le morphisme de $X$ vers $A_m$ ou
$A_1\times_{A_0}\ldots\times_{A_0}A_1$ se factorise à travers un objet
connexe $X'$, son extension se factorise en un morphisme de $X$ vers $X'$ suivi
d'un morphisme de $X'$ vers $A_0\times\ldots\times A_0$. Par définition
d'objet connexe, le morphisme de l'objet connexe $X'$ vers l'objet discret
$A_0\times\ldots\times A_0$ se factorise à travers l'objet final et donc l'extension à $A_0\times\ldots\times A_0$ du
morphisme de $X$ vers $A_m$ ou $A_1\times_{A_0}\ldots\times_{A_0}A_1$ aussi,
comme composée d'un morphisme se factorisant à travers l'objet final. Ceci
montre que le morphisme de $X$ vers $A_m$ ou
$A_1\times_{A_0}\ldots\times_{A_0}A_1$ vérifie la condition de la
proposition~\ref{theta}.\\
CQFD.\\

Montrons maintenant que, dans $\mathcal{C-PC}$, la notion de connexité possède une
caractérisation simple.

\begin{lem}\index{connexité!des \precats}
Soient $\mathcal{C}$ une catégorie \discret, 
Une \precat~$A$ est connexe si et seulement si, pour tout couple d'objets $(a,b)$ de $A$, il
existe un entier m et une suite d'objets $(a_0,\ldots,a_m)$ de $A$ allant de $a$
à $b$ tel que, pour i allant de 0 à m-1, $A_1(a_i,a_{i+1})$ est non vide.
(On dira alors que $a$ et $b$ sont reliés.)
\end{lem}
{\bf Preuve :}\\
Soit $A$ une \precat~vérifiant la caractérisation ci-dessus et montrons que
$A$ est connexe. Soit $f:A\rightarrow B$ un morphisme de \precats~avec $B$ discret. 
Soient $a$ et $a'$ deux objets de $A$.
Par hypothèse sur $A$, alors il existe
une suite $(a_0,\ldots,a_m)$ d'objets de $A$ allant de $a$ à $a'$, pour un
certain m, et tel que, pour i allant de 0 à m-1, $A_1(a_i,a_{i+1})$ est non
vide. Donc son image par $f$,
$B_1(f(a_i),f(a_{i+1}))$, est non vide. Or $B$ est discret donc il vient que
$f(a_i)$ et $f(a_{i+1})$ sont égaux. Comme les $f(a_i)$ sont égaux
deux-à-deux, alors ils sont tous égaux à un certain objet $b$ de $B$. En particulier,
$a$ et $a'$ ont même image par $f$. Donc $f_0$ est constante (à valeur $b$).
De ce fait, tous les $A_n(a_0,\ldots,a_n)$ ont pour image $B_n(b,\dots,b)$ et
donc $A_n$ qui en est le coproduit (car $\mathcal{C}$ est \discret) aussi et ceci quelque soit n. Donc $A$ a pour
image le sous-ensemble $\{b\}$ de $B$ et donc $f$ est constant à valeur
$\{b\}$. Ce qui montre que $A$ est connexe.\\

Soit $A$ une \precat~ne vérifiant pas la caractérisation ci-dessus et montrons que
$A$ n'est pas connexe. Par hypothèse, il existe deux objets $a$ et $a'$ de $A$
tels que pour tout entier m et pour toute suite d'objets $(a_0,\ldots,a_m)$
allant de $a$ à $a'$, on a qu'il existe i entre 0 et m-1 tel que 
$A_1(a_i,a_{i+1})$ soit vide. En particulier, ceci entraîne que tout objet
relié à $a$ n'est pas relié à $a'$ et réciproquement. On va construire
un morphisme de $A$ vers $*\coprod *$ de la manière suivante. On envoie tous
les objets reliés à $a$ sur un exemplaire de l'objet final, notons-le $*_1$, et les autres
objets sur l'autre exemplaire de l'objet final, notons-le $*_2$. Comme dans un
$m$-simplexe tous les sommets sont reliés entre eux par les arêtes, par ce qui
précède, il vient que les $m$-simplexes sont de deux types uniquement : ceux
dont tous les sommets sont reliés à $a$ et ceux dont aucun sommet n'est
relié à $a$. On envoie donc les premiers sur $*_1$ et les seconds sur $*_2$.
Ceci nous définit bien un morphisme de \precats~de $A$ vers $*\coprod *$ qui
par construction ne se factorise pas par l'objet final, ce qui montre que $A$
n'est pas connexe.\\ 
CQFD.

\begin{ex}\label{enssimpconnex}
Dans $\mathcal{ENSSIMP}$, la notion de connexité est caractérisé par le
lemme ci-dessus car les ensembles simpliciaux sont des \precats~avec pour
$\mathcal{C}$ la catégorie des ensembles.
\end{ex}

Grâce à cette notion de connexité et à sa caractérisation simple, nous
allons pouvoir utiliser plus facilement les propriétés de la construction
$\Theta$ qui nous permettra de caractériser les \cats~faciles et les
équivalences faciles de \cats~faciles par des propriétés de relèvements
vis-à-vis de flèches construites à partir des $\Delta[m]\Theta
X$, des $\Upsilon(m)\Theta X$ et des familles de flèches caractérisant
\obcs~faciles et \eqcs~faciles.

\newpage

\section{Flèches génératrices des \cats~faciles}

Par définition, une \cat~facile a ses niveaux qui sont des \obcs~faciles et
ses morphismes de Segal qui sont des \eqcs~faciles d'\obcs~faciles. Ainsi si $A$
est une \cat~facile, pour tout entier $m$, $A_m$ est un \obc~facile donc possède la
\prd~par rapport à la famille $\mathcal{F}_1$. Par la construction $\Theta$
précédente (et sous certaines conditions de connexité), ceci revient à
dire que $A$ a la \prd~par rapport aux flèches de la forme $\Delta[m]\Theta f$
où $f$ décrit $\mathcal{F}_1$. Il reste donc à formuler par propriété
de relèvement que les morphismes de Segal sont des \eqcs~faciles. C'est ce que
va faire la construction $Boit_m(g)$ que nous donnons ci-dessous.

\begin{defin}\index{$Boit_m(g)$}
Soit $\mathcal{C}$ une catégorie \discret~possédant les coproduits et les
sommes amalgamées. 
Soient $m$ un entier et $g:E\rightarrow F$ un morphisme de $\mathcal{C}$. On
définit la flèche $Boit_m(g)$ de la manière suivante.
Considérons le diagramme suivant :
\begin{diagram}
\Upsilon(m)\Theta E        & \rTo^{i_m\Theta E} & \Delta[m]\Theta E &  & \\
\dTo^{\Upsilon(m)\Theta g} &  & \dTo  & \rdTo(2,4)^{\Delta[m]\Theta g} & \\
\Upsilon(m)\Theta F        & \rTo     & B(m,g)   &     & \\
  & \rdTo(4,2)_{i_m\Theta F}  &  & \rdDotsto~{\exists ! Boit_m(g)} & \\
  &    &     &    & \Delta[m]\Theta F \\                          
\end{diagram}
On note $B(m,g)$ la somme amalgamée de $\Upsilon(m)\Theta F $ avec
$\Delta[m]\Theta E$ au-dessus de $\Upsilon(m)\Theta E$. 
L'extérieur du diagramme commute car $*\Theta *$ est un bifoncteur, donc il
existe un unique morphisme de $B(m,g)$ vers $\Delta[m]\Theta F$ que l'on a noté
$Boit_m(g)$.
\end{defin}

\begin{lem}
Soit $\mathcal{C}$ une catégorie \discret~possédant les coproduits, les
produits fibrés au-dessus d'un objet discret et les
sommes amalgamées. 
Soient $A$ une \precat, $m$ un entier et $g:E\rightarrow F$ un morphisme de
$\mathcal{C}$. Considérons le diagramme commutatif suivant :
 \begin{diagram}
 & E & \rTo & & A_m\\
(II) & \dTo^{g} & &\ruDotsto(3,2)_{\exists}& \dTo_{Segal}\\
 & F & \rTo & & A_1\times_{A_0}\ldots\times_{A_0}A_1\\
\end{diagram}
Si le but de $g$ est connexe, alors ce diagramme équivaut au diagramme suivant
:
\begin{diagram}
B(m,g) & \rTo & A & \\
\dTo^{Boit_m(g)} & \ruDotsto_{\exists} &  & (II')\\
\Delta[m]\Theta F & & & \\
\end{diagram}
\end{lem}
{\it Preuve :}\\
Comme le but de $g$ est connexe, en appliquant le lemme~\ref{cons} et
la proposition~\ref{theta}, la partie solide du diagramme (II) devient :
\begin{diagram}
\Upsilon(m)\Theta E        & \rTo^{i_m\Theta E} & \Delta[m]\Theta E &            & \\
\dTo^{\Upsilon(m)\Theta g} &                    & \dTo              & \rdTo(2,4) & \\
\Upsilon(m)\Theta F        & \rTo               & B(m,g)            &            & \\
                           & \rdTo(4,2)         &                   & \rdDotsto~{\exists !} & \\
                           &                    &                   &            & A \\                   
\end{diagram}
Comme la partie solide
du diagramme (II) équivaut à la commutativité de la partie extérieure du
diagramme ci-dessus, il existe bien un unique morphisme de $B(m,g)$ vers $A$. C'est ce
morphisme qui traduit donc la partie solide du diagramme (II). Comme $F$ est
connexe, en appliquant le lemme~\ref{cons} et
la proposition~\ref{theta} au relèvement du diagramme (II), on obtient le
relèvement du diagramme (II').\\
CQFD.\\
\\

Avec cette nouvelle construction $Boit_m(g)$ qui représente au niveau des
\precats~ce qui se passe au niveau des morphismes de Segal et avec la
construction $\Delta[m]\Theta *$ qui représente au niveau des \precats~ce
qui se passe au niveau $m$ des \precats, nous sommes en mesure de donner les
flèches génératrices des \cats~faciles.

\begin{defin}\index{$\mathcal{FG}1$}\index{$\mathcal{FG}2$}\index{$\mathcal{FG}_1$}
Soit $(\mathcal{C},\mathcal{F}_1,\mathcal{F}_2)$ une donnée de Segal proto-facile
dont la catégorie sous-jacente $\mathcal{C}$ possède les coproduits et les
sommes amalgamées.\\ 
On définit la famille $\mathcal{FG}1$ des flèches
génératrices de \cats~faciles de type 1 comme la famille des morphismes de la forme
$\Delta[m]\Theta f$, avec $m$ un entier strictement positif et $f$ une flèche de
$\mathcal{F}_1$.\\
On définit la famille $\mathcal{FG}2$ des flèches
génératrices de \cats~faciles de type 2 comme la famille des morphismes de la forme
$Boit_m(g)$, avec m un entier supérieur ou égal à deux et g une flèche de
$\mathcal{F}_2$.
On note $\mathcal{FG}_1$ la réunion des deux familles
$\mathcal{FG}1$ et $\mathcal{FG}2$. 
\end{defin}

On appelle les morphismes de $\mathcal{FG}1$ les flèches génératrices de
type 1,\index{flèches génératrices!de type 1} car ce sont celles qui caractérisent la condition 1) de la notion de
\cat~facile. Et les morphismes de $\mathcal{FG}2$ sont nommés les flèches génératrices de
type 2\index{flèches génératrices!de type 2} car ce sont celles qui caractérisent la condition 2) de la notion de
\cat~facile. La famille $\mathcal{FG}_1$, réunion des familles $\mathcal{FG}1$
et $\mathcal{FG}2$, est nommée famille des flèches
génératrices des \cats~faciles, car ce sont bien par rapport à ces
flèches que les \cats~faciles se relèvent, comme le montre la proposition
ci-dessous.\index{flèches génératrices!des \cats}

\begin{prop}\label{catfac}
Soit $(\mathcal{C},\mathcal{F}_1,\mathcal{F}_2)$ une donnée de Segal proto-facile
dont la catégorie sous-jacente $\mathcal{C}$ possède les coproduits et les
sommes amalgamées et
telle que les sources et buts des flèches de la famille $\mathcal{F}_1$ soient
connexes ainsi que les buts de la famille $\mathcal{F}_2$.\\
Une \precat~$A$ est une \cat~facile si et seulement si elle a la \prd~par
rapport à la famille $\mathcal{FG}1$ et par rapport à la famille
$\mathcal{FG}2$. 
\end{prop}
{\it Preuve :}\\
Tout d'abord, explicitons la définition de \cat~facile.
Une \precat~$A$ est une \cat~facile si:
\begin{itemize}
\item 1)pour tout entier m, $A_m$ est une catégorie facile dans $\mathcal{C}$,
c'est-à-dire possède la \prd~par rapport aux flèches de la famille
$\mathcal{F}_1$.
\item 2)pour tout entier m supérieur ou égal à deux, le morphisme de Segal
$A_m\rightarrow A_1\times_{A_0}\ldots\times_{A_0}A_1$ est une équivalence
facile d'\obcs~faciles, c'est-à-dire possède la \prd~par rapport aux flèches de la famille
$\mathcal{F}_2$. 
\end{itemize}
La condition 1) se traduit par le diagramme suivant :
\begin{diagram}
 & C & \rTo & A_m & \\
(I) & \dTo^{f} & \ruDotsto_{\exists} & & f\in\mathcal{F}_1\\
 & D & & & \\
\end{diagram}
La condition 2) se traduit par le diagramme suivant :
\begin{diagram}
 & E & \rTo & & A_m & \\
(II) & \dTo^{g} & &\ruDotsto(3,2)_{\exists}& \dTo_{Segal} & g\in\mathcal{F}_2\\
 & F & \rTo & & A_1\times_{A_0}\ldots\times_{A_0}A_1 & \\
\end{diagram}
Comme les sources et buts des flèches de $\mathcal{F}_1$ sont connexes, en
appliquant successivement le lemme~\ref{cons} et
la proposition~\ref{theta}, le diagramme (I) équivaut au diagramme suivant :
\begin{diagram}
\Delta[m]\Theta C & \rTo & A & \\
\dTo^{\Delta[m]\Theta f} & \ruDotsto_{\exists} &  & (I')\\
\Delta[m]\Theta D & & & \\
\end{diagram}
Ceci montre que la condition 1) est vérifiée si et seulement si le
diagramme (I') se relève.\\
Comme en outre les buts des flèches de $\mathcal{F}_2$ sont connexes alors par
le lemme précédent, il vient que le diagramme (II) équivaut au diagramme
(II') du lemme précédent. Ceci montre que la condition 2) est vérifiée
si et seulement si le diagramme (II') se relève.\\
CQFD.

\begin{ex}\label{enssimpconnex2}
Pour les familles de flèches prises dans
l'exemple~\ref{hypenssimp2} pour $\mathcal{ENSSIMP}$, on a que leurs buts sont
connexes car les $m$-simplexes standards sont connexes au sens de
l'exemple~\ref{enssimpconnex}.
\end{ex}

Passons maintenant à la caractérisation par relèvement des équivalences
faciles de \cats~faciles.

\begin{defin}\index{$\mathcal{FG}_2$}
Soit $(\mathcal{C},\mathcal{F}_1,\mathcal{F}_2)$ une donnée de Segal proto-facile
dont la catégorie sous-jacente $\mathcal{C}$ possède les coproduits.
On définit la famille $\mathcal{FG}_2$ des flèches génératrices
d'\eqs~faciles comme la famille constituée de l'inclusion $\emptyset\rightarrow *$ et des
morphismes $\Delta[1]\Theta g$, avec $g$ une flèche de $\mathcal{F}_2$.
\index{flèches génératrices!des équivalences faciles}
\end{defin}

\begin{prop}\label{eqfac}
Soit $(\mathcal{C},\mathcal{F}_1,\mathcal{F}_2)$ une donnée de Segal proto-facile
dont la catégorie sous-jacente $\mathcal{C}$ possède les coproduits et
telle que les buts des flèches de la famille $\mathcal{F}_2$ soient
connexes. 
Un morphisme $f:A\rightarrow B$ de \cats~faciles est une équivalence facile de
\cats~faciles si et seulement si $f$ a la \prd~par rapport aux flèches de la famille
$\mathcal{FG}_2$. 
\end{prop}
{\it Preuve :}\\
Un morphisme $f:A\rightarrow B$ de \cats~faciles est une équivalence facile de
\cats~faciles si :
\begin{itemize}
\item 1') $f_0$ est surjective,
\item 2') pour tout couple $(x,y)$ d'objets de $A$, le morphisme dans $\mathcal{C}$
$f_1(x,y):A_1(x,y)\rightarrow B_1(f(x),f(y))$ est une équivalence facile
d'\obcs~faciles dans
$\mathcal{C}$, c'est-à-dire se relève à droite par rapport aux flèches
de la famille $\mathcal{F}_2$. 
\end{itemize}
La condition 1') se traduit simplement par le diagramme suivant :
\begin{diagram}
 & \emptyset & \rTo & A \\
(III) & \dTo & \ruDotsto_{\exists} & \dTo_{f}\\
 & * & \rTo & B\\
\end{diagram}
La condition 2') se traduit par le diagramme suivant :
\begin{diagram}
 & E & \rTo & A_1(x,y) &  & \\
(IV) & \dTo^{g} & \ruDotsto_{\exists} & \dTo_{f_1(x,y)} & & g\in\mathcal{F}_2\\
 & F & \rTo & B_1(f(x),f(y))& &\\
\end{diagram}
Comme les buts des flèches de la famille $\mathcal{F}_2$ sont connexes, en utilisant le
lemme~\ref{cons} et la
proposition~\ref{theta}, le diagramme (IV)
équivaut au diagramme suivant : 
\begin{diagram}
\Delta[1]\Theta E & \rTo & A & \\
\dTo^{\Delta[1]\Theta g} & \ruDotsto_{\exists} & \dTo_{f} & (IV')\\
\Delta[1]\Theta F & \rTo & B & \\
\end{diagram}
CQFD.\\

Les deux propositions précédentes nous assurent donc que les \cats\\ faciles
et les \eqs~faciles sont caractérisées par des propriétés de
relève\-ment. Ainsi on a muni la donnée de Segal des \cats~de deux familles
$\mathcal{FG}_1$ et $\mathcal{FG}_2$ définissant par relèvement les
\cats~faciles et les \eqs~faciles de \cats~faciles. Il est alors tout naturel de
se demander si avec ces familles la donnée de Segal des \cats~n'est pas une
donnée de Segal proto-facile. Pour cela on va d'abord montrer que les
propriétés de connexité que l'on demande aux familles $\mathcal{F}_1$ et
$\mathcal{F}_2$ sont aussi valables sur $\mathcal{FG}_1$ et $\mathcal{FG}_2$.

\begin{lem}\label{fgcon}
Soit $(\mathcal{C},\mathcal{F}_1,\mathcal{F}_2)$ une donnée de Segal proto-facile
dont la caté\-gorie sous-jacente $\mathcal{C}$ possède les coproduits et les
sommes amalgamées et telle que les sources et buts des flèches de la famille $\mathcal{F}_1$ soient
connexes non vides ainsi que les buts de la famille $\mathcal{F}_2$. Alors les
sources et buts de $\mathcal{FG}_1$ sont connexes non vides ainsi que les but de
$\mathcal{FG}_2$.
\end{lem}
{\it Preuve :}\\
Tout d'abord remarquons qu'en dehors de
l'inclusion $\emptyset\rightarrow *$, dont le but est bien connexe et non vide,
la famille $\mathcal{FG}_2$ est de la forme $\Delta[1]\Theta g$.
Remarquons ensuite que $(\Delta[m]\Theta X)_1(x,y)$ et $(\Upsilon(m)\Theta X)_1(x,x+1)$ ne
sont autres qu'un exemplaire de $X$ (celui qui est associé à l'application
de {\bf 1} dans {\bf m} qui à 0 et 1 associe x et y, respectivement x et x+1).
Donc si $X$ est non vide, il vient que $\Delta[m]\Theta X$ et $\Upsilon(m)\Theta
X$ sont connexes. Leurs ensembles d'objets étant non vides, ils sont eux-mêmes
non vides. Comme par hypothèse sur la donnée de Segal proto-facile, pour
tout morphisme $f:C\rightarrow D$ dans $\mathcal{F}_1$, C et D sont connexes non
vides, alors $\Delta[m]\Theta f:\Delta[m]\Theta C\rightarrow \Delta[m]\Theta D$
a bien sa source et son but connexes et non vides. Ceci montre que les flèches
de $\mathcal{FG}1$ ont leurs sources et buts connexes non vides. De même, les buts des
flèches $g:E\rightarrow F$ de $\mathcal{F}_2$ sont connexes et non vides, donc les buts des
flèches $\Delta[1]\Theta g:\Delta[1]\Theta E\rightarrow \Delta[1]\Theta F$ et $Boit_m(g):B(m,g)\rightarrow \Delta[m]\Theta F$ ont aussi leurs buts
connexes et non vides. Ainsi les buts de flèches de $\mathcal{FG}_2$ sont
connexes non vides. Enfin on remarque que $B(m,g)_1(x_i,x_{i+1})$ n'est autre
qu'un exemplaire de $F$ et donc $B(m,g)$ est connexe non vide. Ainsi les sources
et buts des flèches de $\mathcal{FG}2$ sont connexes non vides.\\
CQFD.\\

\begin{lem}\label{FGens}
Soit $(\mathcal{C},\mathcal{F}_1,\mathcal{F}_2)$ une donnée de Segal proto-facile
dont la caté\-gorie sous-jacente $\mathcal{C}$ possède les coproduits et les
sommes amalgamées et telle que les sources et buts des flèches de la famille $\mathcal{F}_1$ soient
connexes non vides ainsi que les buts de la famille $\mathcal{F}_2$. Alors la
donnée de Segal $\mathcal{C-PC}$ définie dans le lemme~\ref{cpcse} munie des
familles $\mathcal{FG}_1$ et $\mathcal{FG}_2$ constitue une donnée de Segal
proto-facile.
\end{lem}
{\it Preuve :}\\
Montrons tout d'abord que les familles $\mathcal{FG}_1$ et $\mathcal{FG}_2$ sont
des ensembles. Comme $(\mathcal{C},\mathcal{F}_1,\mathcal{F}_2)$ est une
donnée de Segal proto-facile, les familles $\mathcal{F}_1$ et $\mathcal{F}_2$ sont par définition des
ensembles. Or $\mathcal{FG}_1$ est constituée des flèches de type
$\Delta[n]\Theta f$ et $Boit_m(g)$ avec $f$ décrivant $\mathcal{F}_1$, $g$
décrivant $\mathcal{F}_2$ et $m$ et $n$ décrivant $\mathbb{N}$. Donc
$\mathcal{FG}_1$ est donc bien un ensemble comme réunion dénombrable
d'ensembles. De même, comme $\mathcal{FG}_2$ est constituée des flèches de type
$\Delta[1]\Theta g$ avec $g$ décrivant $\mathcal{F}_2$ et de la flèche
$\emptyset\rightarrow *$, c'est un ensemble comme réunion
d'ensembles. Donc les familles $\mathcal{FG}_1$ et $\mathcal{FG}_2$ sont des ensembles.\\
Par la proposition~\ref{catfac}, les \precats~ayant la \prd~par rapport à la famille
$\mathcal{FG}_1$ sont des \cats~faciles et donc des \cat~par le
lemme~\ref{fac}, ce qui donne la première partie de la propriété 8. De
même par la proposition~\ref{eqfac}, les morphismes de \cats~faciles ayant la \prd~par rapport à la famille
$\mathcal{FG}_2$ sont des \eqs~faciles de \cats~faciles et donc des \eqs~de
\cats~par le lemme~\ref{fac},
ce qui donne la première partie de la propriété 10. Comme les produits
fibrés au-dessus d'un objet discret dans $\mathcal{C-PC}$ ou dans sa
catégorie de morphismes se calculent niveau par niveau et que les
\obcs~faciles et les \eqs~faciles d'\obcs~faciles sont stables par produit
fibré au-dessus d'un objet discret, par les propriétés 8) et 10) de la
donnée de Segal proto-facile
$(\mathcal{C},\mathcal{F}_1,\mathcal{F}_2)$, on obtient que les \cats~faciles et
les \eqs~faciles de \cats~faciles sont stables par produit fibré au-dessus
d'un objet discret. Comme par les lemmes~\ref{catfac} et~\ref{eqfac}, les
\cats~faciles sont les objets ayant la \prd~par rapport à $\mathcal{FG}_1$ et
les \eqs~faciles de \cats~faciles sont les morphismes de \cats~faciles ayant la \prd~par rapport
à $\mathcal{FG}_2$, les secondes parties des propriétés 8 et 10 sont donc
bien vérifiées.
Par le lemme~\ref{fgcon}, les sources des
morphismes de la famille $\mathcal{FG}_1$ sont connexes et donc tout morphisme
de la source d'un morphisme de $\mathcal{FG}_1$ vers un objet discret se
factorise par l'objet final. De fait les \precats~discrètes ont la \prd~par
rapport à $\mathcal{FG}_1$, ce qui donne la propriété 9.\\
CQFD.\\

Nos \cats~faciles étant caractérisées par une propriété de
relèvement par rapport à une famille de flèches, on va donc pouvoir
appliquer toute la machinerie développée au chapitre précédent pour
obtenir le procédé de catégorisation tant recherché.

\newpage

\section{Construction Cat}

Nous avons montré dans la section précédente que la notion de \cat~facile est simple à manipuler car elle est
caractérisée par des propriétés de relèvements. En effet, on a
montré que les \precats~$\mathcal{FG}_1$-injectives correspondent exactement
aux \cats~faciles. L'idée pour
catégoriser est donc de prendre la $\mathcal{FG}_1$-injectivisation des
\precats. Ainsi on obtiendra une \cat~facile associée à travers laquelle
tout morphisme de notre \precat~vers une \cat~facile se factorisera. La seule
ombre à ce tableau idyllique est que l'on veut que ce procédé préserve
l'homotopie, ce que l'on est incapable de montrer directement. Aussi a-t-on vu
au chapitre précédent l'intérêt d'une I-injectivisation possédant la
propriété d'unicité de la factorisation des morphismes à but
I-injectifs. Comme au chapitre précédent, on va construire une telle
$\mathcal{FG}_1$-injectivisation en se servant de la notion de marquage partiel des
\precats, qui n'est autre que la notion d'objet partiellement $\mathcal{FG}_1$-marquée.

\begin{defin}\index{\precat!partiellement marquée}
\index{relèvement marqué}
Soit $(\mathcal{C},\mathcal{F}_1,\mathcal{F}_2)$ une donnée de Segal proto-facile
dont la catégorie sous-jacente $\mathcal{C}$ possède les coproduits et les
sommes amalgamées et
telle que les sources et buts des flèches de la famille $\mathcal{F}_1$ soient
connexes ainsi que les buts de la famille $\mathcal{F}_2$.
Une \precat~partiellement marquée est la donnée d'un couple $(A,\mu)$, où
$A$ est une \precat~et $\mu$ une fonction d'un ensemble de diagrammes solides de type (I') et (II')
pour $A$ qui associe à chaque diagramme solide de cet ensemble un relèvement, que l'on dira marqué. Par abus de langage, on dit qu'un diagramme appartient à $\mu$ pour dire qu'il appartient au domaine de définition de la fonction $\mu$ et donc qu'on lui a choisi un relèvement.
\end{defin}

On remarque que si $A$ est totalement marquée (i.e. $\mu$ contient tous les
diagrammes des flèches génératrices de \cats~faciles vers A) alors $A$ a
bien évidemment la \prd~par rapport à toutes les flèches génératrices
de \cats~faciles, ce qui fait de $A$ une \cat~facile pour laquelle tous les
relèvements sont marqués, c'est ce qu'on appelle une \cat~marquée.
 
\begin{defin}\index{\cat!marquée}
Soit $(\mathcal{C},\mathcal{F}_1,\mathcal{F}_2)$ une donnée de Segal proto-facile
dont la catégorie sous-jacente $\mathcal{C}$ possède les coproduits et les
sommes amalgamées et
telle que les sources et buts des flèches de la famille $\mathcal{F}_1$ soient
connexes ainsi que les buts de la famille $\mathcal{F}_2$.
Une \cat~marquée est la donnée d'un couple $(A,\mu)$, où
$A$ est une \cat~facile et $\mu$ une fonction de l'ensemble des diagrammes solides de type (I') et (II')
pour $A$ qui à chaque diagramme solide associe un relèvement.
\end{defin}

Les notions de \precat~totalement marquée et de \cat~marquée sont
équivalentes, comme on l'a vu plus haut. Il nous reste donc à définir ce que
sont les morphismes de \precats~partiellement marquées.

\begin{defin}\index{morphisme!de \precats~partiellement marquées}
\index{morphisme!préservant le marquage}
Soit $(\mathcal{C},\mathcal{F}_1,\mathcal{F}_2)$ une donnée de Segal proto-facile
dont la catégorie sous-jacente $\mathcal{C}$ possède les coproduits et les
sommes amalgamées et
telle que les sources et buts des flèches de la famille $\mathcal{F}_1$ soient
connexes ainsi que les buts de la famille $\mathcal{F}_2$.
Un morphisme $f:A\rightarrow B$ entre les \precats~partiellement
marquées $(A,\mu)$ et $(B,\nu)$ est un morphisme de \precats~partiellement marquées si, pour tout
diagramme de $\mu$, son prolongement par $f$ est un diagramme de $\nu$. On dit
que le morphisme préserve le marquage.\\
Un morphisme de \cats~marquées n'est autre qu'un morphisme de
\precats~partiellement marquées.
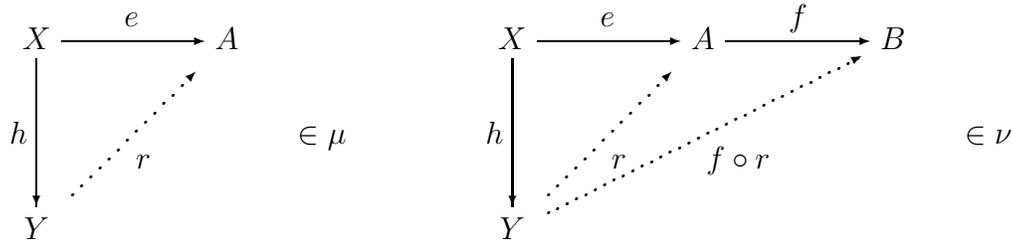
\begin{diagram}
X        & \rTo^{e}      & A &       & & X        & \rTo^{e} & A & \rTo^{f} & B & \\
\dTo^{h} & \ruDotsto_{r} &   &\in\mu & & \dTo^{h} &
\ruDotsto_{r} & &\ruDotsto(4,2)_{f\circ r} & & \in\nu \\
Y        &               &   &       & & Y        &          &   &       &  & \\
\end{diagram}
\end{defin}

On obtient ainsi deux catégories : $\mathcal{C-PC}_m$, la
catégories des \precats~partiellement marquées et sa sous-catégorie pleine
$\mathcal{C-C}_m$ des \cats~mar\-quées. Tout d'abord, remarquons qu'une même
\cat~facile peut avoir plusieurs marquages différents et ainsi peut donner lieu à
plusieurs \cats~mar\-quées, il en va de même pour les \precats. Ainsi nos deux
nouvelles catégories ne sont pas des sous-catégories de $\mathcal{C-PC}$.
En revanche, elles sont toutes deux munies d'un foncteur Oubli vers
$\mathcal{C-PC}$ qui est fidèle mais n'est pas plein car tout morphisme ne
préserve pas le marquage. En outre on a vu que la notion de \cat~facile n'est
pas stable par limite mais en revanche celle de \cats~marquées l'est.

\begin{lem}
Soit $(\mathcal{C},\mathcal{F}_1,\mathcal{F}_2)$ une donnée de Segal proto-facile
dont la caté\-gorie sous-jacente $\mathcal{C}$ possède les coproduits, les limites et les
sommes amalgamées et
telle que les sources et buts des flèches de la famille $\mathcal{F}_1$ soient
connexes ainsi que les buts de la famille $\mathcal{F}_2$.
La catégorie $\mathcal{C-C}_m$ des \cats~marquées avec les morphismes
préservant les marquages possède les limites.
\end{lem}
{\it Preuve :} application directe du lemme~\ref{mst} !\\

Cette notion de marquage va nous permettre de définir un procédé de caté\-gorisation qui
à une \precat~va associer une \cat~marquée, procédé qui non seulement
sera fonctoriel mais servira d'adjoint au foncteur Oubli des \cats~marquées
vers les \precats. On nommera ce procédé Cat. Pour définir Cat, on va se
servir de la construction $E_{\Phi}$ définie au chapitre précédent avec
pour $\Phi$ la famille $\mathcal{FG}_1$. Toutefois si l'on veut que le
procédé ainsi défini $\mathcal{FG}_1$-injectivise bien, c'est-à-dire
catégorise bien, il faut que les flèches de $\Phi$, ici de $\mathcal{FG}_1$,
ait des propriétés de petitesse et que les $\Phi$-cofibrations soient des monomorphismes. Le lemme suivant nous montre que si les
familles $\mathcal{F}_1$ et $\mathcal{F}_2$ de la donnée de Segal proto-facile
ont leurs sources et buts petits alors $\mathcal{FG}_1$ aussi.

\begin{lem}\label{flgenpet}
Soit $(\mathcal{C},\mathcal{F}_1,\mathcal{F}_2)$ une donnée de Segal proto-facile
dont la caté\-gorie sous-jacente possède les coproduits et les sommes
amalgamées et vérifie les deux propriétés suivantes :
\item -les sources et buts des flèches de la famille $\mathcal{F}_1$ sont
connexes ainsi que les buts de la famille $\mathcal{F}_2$,
\item -il existe un cardinal transfini régulier strictement supérieur à
$\aleph_0$, que l'on notera $\alpha$, pour lequel tout morphisme de
$\mathcal{F}_1$ et $\mathcal{F}_2$ a sa source et son but $\alpha$-petits au
sens~\ref{alphapetit}.\\ 
\\
Alors tout morphisme de $\mathcal{FG}_1$ et $\mathcal{FG}_2$ a aussi sa source et son but
$\alpha$-petits.
\end{lem}
{\it Preuve :}\\
Par construction, $(\Delta[m]\Theta X)_n$ et $(\Upsilon(m)\Theta X)_n$ sont des
sommes amalgamées d'un ensemble fini avec un coproduit fini d'exemplaire de X.
Si X est $\alpha$-petit alors $(\Delta[m]\Theta X)_n$ et $(\Upsilon(m)\Theta
X)_n$ aussi, par le lemme~\ref{cpdtpetit} et, comme on a pris $\alpha$
régulier et strictement supérieur à $\aleph_0$, par le
lemme~\ref{precatpetit}, $\Delta[m]\Theta X$ et $\Upsilon(m)\Theta X$ sont
$\alpha$-petits ainsi que $B(m,g)$, somme amalgamée des précédents, si $g$
est $\alpha$-petite. Comme, par ailleurs, l'inclusion $\emptyset\rightarrow *$ a sa source et son but
$\aleph_0$-petits, les ensembles $\mathcal{FG}_1$ et $\mathcal{FG}_2$ ont bien leurs morphismes
ayant les sources et buts $\alpha$-petits.\\
CQFD.\\

Un des corollaires de ce lemme est que la famille $\mathcal{FG}_1$ permet l'argument du petit objet. Ainsi les $\mathcal{FG}_1$-cofibrations sont des rétracts de colimites séquentielles transfinies de sommes amalgamées de flèches de $\mathcal{FG}_1$. Montrons alors que si les rétracts de colimites séquentielles transfinies de sommes amalgamées de flèches de $\mathcal{F}_1$ et de $\mathcal{F}_2$ sont des monomorphismes, alors les $\mathcal{FG}_1$-cofibrations sont des monomorphismes.

\begin{lem}\label{flgenmono}
Soit $(\mathcal{C},\mathcal{F}_1,\mathcal{F}_2)$ une donnée de Segal proto-facile
dont la caté\-gorie sous-jacente possède les colimites et vérifie les trois propriétés suivantes :
\item -les sources et buts des flèches de la famille $\mathcal{F}_1$ sont
connexes ainsi que les buts de la famille $\mathcal{F}_2$,
\item -il existe un cardinal transfini régulier strictement supérieur à
$\aleph_0$, que l'on notera $\alpha$, pour lequel tout morphisme de
$\mathcal{F}_1$ et $\mathcal{F}_2$ a sa source et son but $\alpha$-petits au
sens~\ref{alphapetit},
\item -les $\mathcal{F}_1\cup\mathcal{F}_2\cup\{\emptyset\rightarrow *\}$-cofibrations sont des monomorphismes.\\ 
\\
Alors les $\mathcal{FG}_1\cup\mathcal{FG}_2$-cofibrations sont des monomorphismes.
\end{lem}
{\it Preuve :}\\
Par le lemme précédent, les familles $\mathcal{FG}_1$ et $\mathcal{FG}_2$ permettent l'argument du petit objet et donc la famille réunion $\mathcal{FG}_1\cup\mathcal{FG}_2$ aussi. Ceci permet de caractériser les $\mathcal{FG}_1\cup\mathcal{FG}_2$-cofibrations comme les rétracts de colimites séquen\-tielles transfinies de sommes amalgamées de flèches de $\mathcal{FG}_1\cup\mathcal{FG}_2$. Or niveau par niveau, les flèches de $\mathcal{FG}_1\cup\mathcal{FG}_2$ sont des coproduits de flèches de $\mathcal{F}_1\cup\mathcal{F}_2\cup\{\emptyset\rightarrow *\}$. Ainsi nos $\mathcal{FG}_1\cup\mathcal{FG}_2$-cofibrations sont donc niveau par niveau des rétracts de colimites séquentielles transfinies de sommes amalgamées de flèches de $\mathcal{F}_1\cup\mathcal{F}_2\cup\{\emptyset\rightarrow *\}$, donc ce sont niveau par niveau des $\mathcal{F}_1\cup\mathcal{F}_2\cup\{\emptyset\rightarrow *\}$-cofibrations. Or par hypothèse, les $\mathcal{F}_1\cup\mathcal{F}_2\cup\{\emptyset\rightarrow *\}$-cofibrations sont des monomorphismes. Donc les $\mathcal{FG}_1\cup\mathcal{FG}_2$-cofibrations sont niveau par niveau des monomorphismes. Comme $\mathcal{C-PC}$ est une catégorie de préfaisceaux sur $\mathcal{C}$, les monomorphismes de $\mathcal{C-PC}$ sont les monomorphismes niveau par niveau. Ainsi les $\mathcal{FG}_1\cup\mathcal{FG}_2$-cofibrations sont bien des monomorphismes dans $\mathcal{C-PC}$.\\
CQFD.

\begin{ex}\label{enssimppet}
Les $m$-simplexes standards et leurs bords ayant un nombre fini de simplexes non
dégénérés sont $\aleph_0$-petits et donc, avec les familles de
flèches choisies dans l'exemple~\ref{hypenssimp2}, $\mathcal{ENSSIMP}$
vérifie l'hypothèse de petitesse de la proposition-définition de Cat
ci-dessous.\\
 En outre les familles de flèches de l'exemple~\ref{hypenssimp2} sont des monomorphismes. De plus dans $\mathcal{ENSSIMP}$, $\emptyset\rightarrow *$ est un monomorphisme et les monomorphismes sont stables par rétract, somme amalgamée le long d'un morphisme, coproduit et comilite séquentielle transfinie. Ainsi l'hypothèse concernant les monomorphismes de la proposition-définition de Cat
ci-dessous est vérifiée également.  
\end{ex} 

Nous avons maintenant tous les éléments en main pour pouvoir donner la proposition-définition
du procédé de catégorisation Cat.

\begin{prop}[-définition]\index{Cat}
Soit $(\mathcal{C},\mathcal{F}_1,\mathcal{F}_2)$ une donnée de Segal proto-facile
dont la catégorie sous-jacente est cocomplète et vérifie les trois propriétés suivantes :
\item -les sources et buts des flèches de la famille $\mathcal{F}_1$ sont
connexes ainsi que les buts de la famille $\mathcal{F}_2$,
\item -il existe un cardinal transfini régulier strictement supérieur à
$\aleph_0$, que l'on notera $\alpha$, pour lequel tout morphisme de
$\mathcal{F}_1$ et $\mathcal{F}_2$ a sa source et son but $\alpha$-petits au
sens~\ref{alphapetit},
\item -les $\mathcal{F}_1\cup\mathcal{F}_2\cup\{\emptyset\rightarrow *\}$-cofibrations sont des monomorphismes.\\  
\\
Il existe un foncteur Cat des \precats~vers les \cats~faciles mar\-quées muni
d'une transformation naturelle $can$ entre l'identité des \precats~et la
composée de Cat avec le foncteur Oubli des \cats~faciles marquées vers les
\precats~ayant la propriété universelle suivante :\\
\\
pour tout morphisme $f:A\rightarrow B$ d'une \precat~quelconque $A$ vers une
\cat~marquée $(B,\lambda)$, il existe un unique morphisme $\tilde{f}$ de $Cat(A)$ vers $B$
préservant le marquage et dont la précomposition par $can_A$ est $f$.
\begin{diagram}
A & \rTo^{f} & B \\
\dTo{can_A} & \ruDotsto_{\exists !\tilde{f}} & \\
Cat(A) & & \\
\end{diagram}
\end{prop}
{\it Preuve :}\\
Nous allons prendre pour foncteur Cat le plan d'addition de cellules $E_{\Phi}$
avec pour $\Phi$ la famille $\mathcal{FG}_1$. Les hypothèses sur la donnée
de Segal proto-facile permettent de vérifier les hypothèses des lemmes~\ref{flgenpet} et~\ref{flgenmono}, ce qui nous donne que les morphismes de $\Phi$ ont leurs
sources et buts $\alpha$-petit et que les $\Phi$-cofibrations sont des monomorphismes. On peut alors appliquer la
proposition~\ref{Ephi2}, qui nous dit que $E_{\Phi}$ est un foncteur à valeur
dans les objets $\Phi$-injectifs marqués et que tout morphisme $A\rightarrow B$ à but
$\Phi$-injectif marqué se factorise de manière unique à travers le
morphisme naturel $A\rightarrow E_{\Phi}(A)$ en un morphisme
préservant le marquage. Or les objets $\Phi$-injectifs marqués ne sont
autres que les \cats~faciles marquées, ce qui nous donne que Cat va dans les
\cats~faciles marquées et vérifie la propriété universelle. Enfin la
transformation naturelle $can$ est induite par le morphisme naturel
$A\rightarrow E_{\Phi}(A)$.\\
CQFD.\\

Comme on a pu le faire remarquer au sujet du plan d'addition de cellules
$E_{\Phi}$, la factorisation ci-dessus s'applique aussi bien pour les
\cats~faciles mais il n'y aura pas d'unicité de la factorisation à travers $Cat(A)$.

\begin{cor}
Soit $(\mathcal{C},\mathcal{F}_1,\mathcal{F}_2)$ une donnée de Segal proto-facile
dont la catégorie sous-jacente est cocomplète et vérifie les trois propriétés suivantes :
\item -les sources et buts des flèches de la famille $\mathcal{F}_1$ sont
connexes ainsi que les buts de la famille $\mathcal{F}_2$,
\item -il existe un cardinal transfini régulier strictement supérieur à
$\aleph_0$, que l'on notera $\alpha$, pour lequel tout morphisme de
$\mathcal{F}_1$ et $\mathcal{F}_2$ a sa source et son but $\alpha$-petits au
sens~\ref{alphapetit},
\item -les $\mathcal{F}_1\cup\mathcal{F}_2\cup\{\emptyset\rightarrow *\}$-cofibrations sont des monomorphismes.\\ 
\\
La construction Cat est un foncteur de $\mathcal{C-PC}$ vers
$\mathcal{C-C}_m$ qui est adjoint à gauche du foncteur Oubli, i.e. on a un
isomorphisme naturel en $A$ et $(B,\lambda)$, induit par le morphisme $can_A$ :
$$\operatorname{Hom}_{\mathcal{C-C}_m}(Cat(A),(B,\lambda))=\operatorname{Hom}_{\mathcal{C-PC}}(A,Oubli((B,\lambda)))$$
De ce fait, le foncteur Cat préserve les colimites.
\end{cor}
{\it Preuve :}\\
De la propriété universelle découle l'unicité à isomorphisme près de
$Cat(A)$, ce qui entraîne que $Cat$ est en
fait un adjoint à gauche du foncteur Oubli des \cats~marquées vers les
\precats. Enfin comme foncteur adjoint à gauche, Cat préserve les
colimites.\\ 
CQFD.\\

On prendra bien garde que, contrairement aux limites, les colimites dans
$\mathcal{C-C}_m$ ne sont pas les colimites niveau par niveau dans $\mathcal{C-PC}$ des
\precats~sous-jacentes avec pour marquage la colimites des marquages. Par
exem\-ple, pour les (petites) catégories, une colimite niveau par niveau de
nerfs de catégories n'est pas un nerf de catégorie en général.\\

Nous avons désormais un procédé de catégorisation Cat ayant de bonnes
propriétés mais pas encore celle de préserver le type d'homotopie.
Toutefois, il est intéressant de récapituler les hypothèses faites sur la
donnée de Segal proto-facile qui ont permis la construction Cat. Cela donnera
lieu à la définition de donnée de Segal pré-facile et à une
proposition montrant que, pour toute donnée de Segal pré-facile, il existe
une catégorisation de type Cat.

\begin{defin}\index{donnée de Segal!pré-facile}
Une donnée de Segal pré-facile est une triplet
$(\mathcal{C},\mathcal{F}_1,\mathcal{F}_2)$ constitué
d'une donnée de Segal $\mathcal{C}$, dont la catégorie sous-jacente est
cocomplète, et de deux familles $\mathcal{F}_1$ et
$\mathcal{F}_2$ de morphismes de $\mathcal{C}$ qui sont des ensembles satisfaisant les propriétés
suivantes :
\item 8) Les objets de $\mathcal{C}$ ayant la \prd~par rapport à
$\mathcal{F}_1$ sont des \obcs~et sont stables par produit fibré au-dessus
d'un objet discret.
\item 9) Les objets discrets de $\mathcal{C}$ ont la \prd~par rapport à
$\mathcal{F}_1$.
\item 10) Les morphismes qui ont la \prd~par rapport à
$\mathcal{F}_2$ et dont la source et le but sont des
\obcs~ayant la \prd~par rapport à $\mathcal{F}_1$ sont des \eqcs~d'\obcs~et le produit fibré, dans
la catégorie des morphismes, de deux tels morphismes au-dessus d'un objet
discret a encore la \prd~par rapport à $\mathcal{F}_2$.
\item 11) Les sources et buts des flèches de la famille $\mathcal{F}_1$ sont
connexes non vides ainsi que les buts de la famille $\mathcal{F}_2$. 
\item 12) Il existe un cardinal transfini régulier strictement supérieur à
$\aleph_0$, que l'on notera $\alpha$, pour lequel tout morphisme de
$\mathcal{F}_1$ et $\mathcal{F}_2$ a sa source et son but $\alpha$-petits au
sens~\ref{alphapetit}.
\item 13) Les $\mathcal{F}_1\cup\mathcal{F}_2\cup\{\emptyset\rightarrow *\}$-cofibrations sont des monomorphismes.
\end{defin} 

On remarque qu'ici nous avons demandé en hypothèse 11) en plus de la connexité des sources et buts des flèches génératrices leur non vacuité. Si la non vacuité des sources et buts des flèches génératrices est superflue pour assurer l'existence de Cat, en revanche c'est l'une des hypothèses permettant d'assurer que la donnée de Segal des \cats~munie des familles $\mathcal{FG}_1$ et $\mathcal{FG}_2$ est une donnée de Segal proto-facile. Montrons maintenant qu'une donnée de Segal pré-facile donne bien lieu à un foncteur Cat.

\begin{prop}\label{catadj}
Soit $(\mathcal{C},\mathcal{F}_1,\mathcal{F}_2)$ une donnée de Segal
pré-facile.
Alors il existe un foncteur Cat de $\mathcal{C-PC}$ vers
$\mathcal{C-C}_m$ qui est adjoint à gauche du foncteur Oubli, i.e. on a un
isomorphisme naturel en $A$ et $(B,\lambda)$ :
$$\operatorname{Hom}_{\mathcal{C-C}_m}(Cat(A),(B,\lambda))=\operatorname{Hom}_{\mathcal{C-PC}}(A,Oubli((B,\lambda)))$$
De ce fait, le foncteur Cat préserve les colimites.
\end{prop}
{\it Preuve :} c'est juste une reformulation du corollaire précédent.\\

On a déjà vu que la catégorie $\mathcal{C-PC}$ munie des \cats, des
\eqs~de \cats~et des familles $\mathcal{FG}_1$ et $\mathcal{FG}_2$ forme une
donnée de Segal proto-facile. Il est donc intéressant de savoir si cette
donnée de Segal proto-facile est aussi une donnée de Segal pré-facile.

\begin{lem}\label{FGpre}
Soit $(\mathcal{C},\mathcal{F}_1,\mathcal{F}_2)$ une donnée de Segal
pré-facile.\\
Alors la donnée de Segal $\mathcal{C-PC}$ définie dans le lemme~\ref{cpcse} munie des
familles $\mathcal{FG}_1$ et $\mathcal{FG}_2$ constitue une donnée de Segal
pré-facile.
\end{lem}
{\it Preuve :}\\
On a déjà montré dans le lemme~\ref{FGens} qu'avec les familles
$\mathcal{FG}_1$ et $\mathcal{FG}_2$, la donnée de Segal sur $\mathcal{C-PC}$
des \cats~et \eqs~de \cats~forme une donnée de Segal proto-facile.
Il ne reste donc plus qu'à montrer que $\mathcal{C-PC}$ est cocomplète et
que les familles $\mathcal{FG}_1$ et $\mathcal{FG}_2$ vérifient les
propriétés 11), 12) et 13).\\
Tout d'abord, comme $\mathcal{C-PC}$ est une catégorie de préfaisceaux sur
la catégorie $\mathcal{C}$ qui est cocomplète par hypothèse, alors
$\mathcal{C-PC}$ est cocomplète. En outre, par hypothèse, les familles $\mathcal{F}_1$
et $\mathcal{F}_2$ vérifient la propriété 11) de connexité, donc , par
le lemme~\ref{fgcon}, les familles $\mathcal{FG}_1$ et $\mathcal{FG}_2$ aussi. De même, par hypothèse, les familles $\mathcal{F}_1$
et $\mathcal{F}_2$ vérifient la propriété 11) ainsi que la propriété
12) sur la petitesse, donc , par
le lemme~\ref{flgenpet}, les familles $\mathcal{FG}_1$ et $\mathcal{FG}_2$
vérifient aussi la propriété 12). Enfin les familles $\mathcal{F}_1$
et $\mathcal{F}_2$ vérifient par hypothèse les propriétés 11), 12) et 13), donc, par le lemme~\ref{flgenmono}, les familles $\mathcal{FG}_1$ et $\mathcal{FG}_2$
vérifient bien la propriété 13) concernant les monomorphismes.\\
CQFD.\\

Grâce à notre nouvelle notion de donnée de Segal pré-facile, nous avons
assurer l'existence d'un foncteur de catégorisation Cat
adjoint au foncteur Oubli des \cats~faciles marquées vers les \precats.
Cependant nous ne savons toujours pas si ce procédé préserve le type
d'homotopie, c'est-à-dire si pour toute \precat~$A$ l'image par Cat du
morphisme naturel $can_A:A\rightarrow Cat(A)$ est une \eq~de \cats. C'est
l'objet du chapitre suivant que de résoudre ce problème.

\newpage
...

\chapter{Préservation du type d'homotopie par catégorisation}

\newpage

Nous avons construit le procédé de catégorisation Cat qui est un adjoint
au foncteur Oubli des \cats~faciles marquées vers les \precats. Cependant,
pour que cette construction soit une bonne catégorisation, nous lui demandons
de préserver le type d'homotopie des \precats. Plus précisément, nous voulons que, pour
toute \precat~$A$, l'image par Cat du morphisme naturel $can_A:A\rightarrow
Cat(A)$ soit une \eq~de \cats. Or il est très difficile de montrer directement
ce résultat. En nous appuyant sur la proposition~\ref{Ephi3}, nous savons que
si les \eqs~de \cats~vérifient certaines propriétés, dont celle de "trois
pour deux", et si nous trouvons un procédé de catégorisation préservant
l'homotopie des \cats, alors notre construction Cat préservera le type
d'homotopie des \precats.\\

C'est pourquoi nous allons débuter ce chapitre par montrer les propriétés
que les \eqs~de \cats~doivent vérifier pour satisfaire les hypothèses de la
proposition~\ref{Ephi3}. Puis nous poursuivrons cette partie par la recherche
d'une nouvelle catégorisation que l'on nommera Bigcat qui aura justement la
propriété de préserver le type d'homotopie des \cats. Comme cette
recherche est assez compliquée, nous allons tout d'abord donner des
constructions simples qui non seulement se comportent bien pour la composition
finie et transfinie mais aussi formalisent la somme amalgamée avec les
flèches génératrices de \cats~faciles. La bonne compréhension de ces
constructions nous mènera à des plans d'addition de certaines flèches
génératrices de \cats~faciles qui ont la faculté de préserver
l'homotopie des \cats. Nous pourrons alors définir Bigcat comme un plan
d'addition formé de tels plans. Nous terminerons ce chapitre par l'application
de la proposition~\ref{Ephi3} au triplet Cat, Bigcat, \eqs~de \cats.

\newpage

\section{Propriétés des \eqs~de \cats}

Commençons cette section en rappelant brièvement la définition d'une
\eq~de \cats. Un morphisme $f:A\rightarrow B$ entre \cats~est une \eq~de
\cats~si :\\
- l'application d'ensembles $\tau_0(f):\tau_0(A)\rightarrow \tau_0(B)$
est surjective,\\ 
- pour tout couple $(x,y)$ d'objets de $A$, le morphisme $f_1(x,y):A_1(x,y)\rightarrow
B_1(f(x),f(y))$ est une \eqc~d'\obcs.\\ 
Parmi les propriétés demandées aux \eqs~de \cats~par la
proposition~\ref{Ephi3}, deux concernent des propriétés du genre "trois pour
deux" et une demande que les isomorphismes soient des \eqs~de \cats. Nous avons
déjà cette propriété par le lemme~\ref{cpcse}. Toutefois on va ici
montrer un résultat plus général qui nous assure que les \eqs~niveau par
niveau, c'est-à-dire les morphismes de \cats~qui niveau par niveau sont des
\eqcs~d'\obcs, sont bien des \eqs~de \cats. Ce résultat sera souvent très utile car
dans de nombreux cas où l'on aura à montrer qu'un morphisme est une \eq~de
\cats, en particulier pour montrer que Bigcat préserve l'homotopie des \cats,
on pourra se réduire à montrer que le morphisme est niveau par niveau une
\eqc~d'\obcs.

\begin{defin}\index{equivalence!de \cats~niveau par niveau}
Soit $\mathcal{C}$ une donnée de Segal.\\
On appelle \eq~de \cats~niveau par niveau tout morphisme de \cats~$f$ tel que :
\item - $f_0$ est une bijection d'ensembles,
\item - pour tout entier $m>0$, $f_m$ est une \eqc~d'\obcs.
\end{defin}

\begin{lem}\label{eqcniv}
Soit $\mathcal{C}$ une donnée de Segal.
Alors les \eqs~de \cats~niveau par niveau sont des \eqs~de \cats.
\end{lem}
{\it Preuve :}\\
Soit $f:A\rightarrow B$ une \eq~de \cats~niveau par niveau. Comme $f_0$ est une
bijection et que, pour tout entier $m>0$, $f_m$ est une \eqc~d'\obcs, alors, pour
tout entier $n\geq 0$, $\tau_0(f_n)$ est une bijection car $\tau_0$ envoie les
\eqcs~d'\obcs~et les bijections d'objets discrets sur des bijections
ensemblistes. Ainsi $\tau_1(f)$ est un isomorphisme de catégories, ce qui
entraîne que $\tau_0(f)$ est une bijection. Donc $f$ est bien
essentiellement surjective.\\ 
Soit $(a,a')$ un couple d'objets de
$A$. Comme $f_1(a,a')$ n'est autre que le produit fibré de $f_1$ et de
l'identité de l'objet final au-dessus de l'isomorphisme $f_0\times f_0$ et que
les \eqcs~d'\obcs~sont stables par produit fibré au-dessus d'un objet discret,
il vient que $f_1(a,a')$ est une \eqc~d'\obcs, et donc $f$ est bien pleinement
fidèle.\\
CQFD.\\

Un des cas particuliers de ce lemme concerne les morphismes de \cats~qui sont niveau par niveau des
isomorphismes de $\mathcal{C}$. Comme les isomorphismes d'\obcs~sont des
\eqcs~d'\obcs, les isomorphismes niveau par niveau sont bien des \eqs~de
\cats~niveau par niveau, donc par le lemme ce sont des \eqs~de \cats. Ceci est
une autre manière de montrer que les isomorphismes de \cats~sont des \eqs~de
\cats, ce qui est l'une des trois hypothèses sur les \eqs~de \cats~de la
proposition~\ref{Ephi3}. Montrons maintenant les deux autres qui sont du type
"trois pour deux".

\begin{lem}\label{3pour2c}
Soit $\mathcal{C}$ une donnée de Segal vérifiant les propriétés suivantes :
\item 1) La catégorie sous-jacente $\mathcal{C}$ est une \cmf~dont tous les
objets sont cofibrants.
\item 2) Les \eqcs~d'\obcs~sont exactement les
équivalences faibles de la \cmf~$\mathcal{C}$ entre \obcs.
\item 3) Le foncteur $\tau_0$ est tel que, pour tout objet $C$ de
$\mathcal{C}$, $\tau_0(C)$ est un quotient de l'ensemble des morphismes dans
$\mathcal{C}$ de l'objet final vers $C$.
\item 4) Il existe un \obc~contractile $K$ muni de deux
morphismes de l'objet final vers $K$ notés 0 et 1, et ayant la propriété
suivante : pour tout couple $(f,g)$ de morphismes de l'objet final vers un \obc~$C$ tel que leurs images par $\tau_0$ soient égales, il existe un
morphisme de $K$ vers $C$ envoyant 0 sur $f$ et 1 sur $g$.
\item 5) Les équivalences faibles de $\mathcal{C}$ sont stables par produit dans la catégorie
des morphismes de $\mathcal{C}$.
\item 6) Le produit fibré d'une équivalence faible entre
\obcs~le long d'une fibration est une équivalence faible.\\
\\
Alors les \eqs~de \cats~vérifient
la propriété de "trois pour deux", i.e. pour tout couple $(f,g)$ de morphismes de
\cats~composables, si deux morphismes parmi $f,\; g \mbox{ et }g\circ f$ sont
des \eqs~de \cats, alors le morphisme restant est aussi une
\eq~de \cats.
\end{lem}
{\it Preuve :}\\
Soient $f:A\rightarrow B$ et $g:B\rightarrow C$ deux morphismes de
\cats. Supposons d'abord que f et g sont des équivalences de \cats.
On a donc par propriété du foncteur $\tau_0$ que $\tau_0(f)$ et $\tau_0(g)$
sont des bijections et donc $\tau_0(g\circ f)=\tau_0(g)\circ \tau_0(f)$ est une
bijection comme composée de bijections, ce qui rend $g\circ f$ essentiellement
surjective. En outre, pour tout couple $(a,a')$ d'objets de $A$, on a $f_1(a,a')$
et $g_1(f(a),f(a'))$ \eqcs~d'\obcs~et donc, par
hypothèse 2), ce sont des équivalences faibles.
Comme $\mathcal{C}$ est une \cmf~et donc vérifie l'axiome "trois pour deux",
il vient que $(g\circ f)_1(a,a')$ est une
équivalence faible entre \obcs, et par 2), une \eqc~d'\obcs. D'où $g\circ f$ est pleinement fidèle. On a donc
montré que si $f$ et $g$ sont des équivalences de \cats, alors $g\circ f$
aussi.\\

Supposons maintenant que $g\circ f$ et $g$ sont des équivalences de \cats.
Alors on a $\tau_0(g\circ f)$ et $\tau_0(g)$ bijectifs, ce qui entraîne que
$\tau_0(f)$ aussi. De même, pour tout couple $(a,a')$ d'objets de $A$,
$(g\circ f)_1(a,a')$ et $g_1(f(a),f(a'))$ sont des \eqcs~d'\obcs~et donc par 2) des équivalences faibles, d'où, par l'axiome "trois pour deux" vérifié par
la \cmf~$\mathcal{C}$, $f_1(a,a')$ est aussi une équivalence faible entre
\obcs, donc par 2) une \eqc~d'\obcs.
Ceci montre que si $g\circ f$ et $g$ sont des équivalences de \cats, alors $f$
aussi.\\

Supposons enfin que $g\circ f$ et $f$ sont des équivalences de \cats. Comme
pour le cas précédent, on a $\tau_0(g\circ f)$ et $\tau_0(f)$ bijectifs, ce qui entraîne que
$\tau_0(g)$ aussi et donc que $g$ est essentiellement surjective. On veut
maintenant montrer que pour tout couple $(b,b')$ d'objets de $B$, $g_1(b,b')$
est une \eqc~d'\obcs. On remarque
tout d'abord que comme $f$ est essentiellement surjective, il existe $a$ et $a'$
objets de $A$ tels que $f(a)$ soit équivalent à $b$ et $f(a')$ équivalent
à $b'$. Or $f_1(a,a')$ et $(g\circ f)_1(a,a')$ sont des \eqcs~d'\obcs~et donc
par 2) des équivalences faibles, d'où, par l'axiome "trois pour deux" vérifié par
la \cmf~$\mathcal{C}$, $g_1(f(a),f(a'))$ est aussi une équivalence faible
entre \obcs, donc par 2) une \eqc~d'\obcs. Nous cherchons donc à transmettre cette propriété à
$g_1(b,b')$ à travers les équivalences d'objets $f(a)\sim b$ et $f(a')\sim
b'$.\\

Pour ce faire, nous allons construire, pour tout objet $x$ de $A$, un foncteur
$F$
de la catégorie comma $(x,\tau_1(A))$
vers $Ho-\mathcal{C}$ la catégorie homotopique de $\mathcal{C}$, qui à un
objet $y$ de $(x,A)$ associe le remplacement fibrant de $A_1(x,y)$, que l'on notera
$A^f_1(x,y)$. Soit $f$ un morphisme de l'objet final vers $A_1(y,z)$ représentant
d'une classe de $(x,\tau_1(A))$, nous voulons définir $F(f)$.\\ 

Considérons le
morphisme de Segal $A_2(x,y,z)\rightarrow A_1(x,y)\times A_1(y,z)$, c'est une
équivalence faible du fait que $A$ est une \cat~et en utilisant le
lemme~\ref{cat2} et 2). 
Comme  les morphismes de $A_1(x,y)$ et de $A_1(y,z)$ vers leurs remplacements fibrants sont
des équivalences faibles, par hypothèse 5), le morphisme produit
$A_1(x,y)\times A_1(y,z)\rightarrow A^f_1(x,y)\times A^f_1(y,z)$ est une
équivalence faible et, par l'axiome "trois pour deux" dans la
\cmf~$\mathcal{C}$, sa précomposition par le morphisme de Segal est aussi une
équivalence faible. Comme $\mathcal{C}$ est une \cmf, on peut décomposer
l'équivalence faible $A_2(x,y,z)\rightarrow A^f_1(x,y)\times A^f_1(y,z)$ en
une cofibration triviale suivie d'une fibration, qui par l'axiome "trois pour
deux" sera aussi triviale. Notons $A'_2(x,y,z)$ l'objet par lequel
l'équivalence se factorise. Posons $A_2(x,y,f)$ le produit fibré de la
fibration $A'_2(x,y,z)\rightarrow A^f_1(y,z)$ le long du morphisme $*\rightarrow
A^f_1(y,z)$ induit par $f$. Par produit fibré le long d'un morphisme d'une fibration, $A_2(x,y,f)$
est fibrant. Considérons le diagramme suivant :
\begin{diagram}
A_2(x,y,f) & \rTo & A'_2(x,y,z)\\
\dTo & & \dTo \\
A^f_1(x,y) & \rTo & A^f_1(x,y)\times A^f_1(y,z)\\
\dTo & & \dTo \\
* & \rTo & A^f_1(y,z)\\
\end{diagram}
Le diagramme du bas est cartésien et par définition de $A_2(x,y,f)$ le
diagramme total est cartésien. Donc le diagramme du haut est aussi cartésien
et donc $A_2(x,y,f)\rightarrow A^f_1(x,y)$ est une fibration triviale comme
produit fibré de la fibration triviale $A'_2(x,y,z)\rightarrow
A^f_1(x,y)\times A^f_1(y,z)$ le long d'un morphisme.\\ 
En outre $A_2(x,y,z)\rightarrow A'_2(x,y,z)$ étant une cofibration triviale,
le morphisme $A_2(x,y,z)\rightarrow A^f_1(x,z)$  s'étend à $A'_2(x,y,z)$, et
l'on obtient un morphisme\\ $A_2(x,y,f)\rightarrow A^f_1(x,z)$. Considérons donc le diagramme suivant :
\begin{diagram}
A_2(x,y,f) & \rTo^{\sim} & A^f_1(x,y)\\
\dTo & \ldDotsto_{\exists F(f)} & \\
A^f_1(x,z) & & \\
\end{diagram}
Comme par hypothèse 1), dans $\mathcal{C}$ tous les objets sont cofibrants et que les objets $A_2(x,y,f)$,
$A^f_1(x,y)$ et $A^f_1(x,z)$ sont
fibrants, ce diagramme existe bien dans la catégorie homotopique
$Ho-\mathcal{C}$. En outre $A_2(x,y,f)\rightarrow A^f_1(x,y)$ est une
équivalence faible dans $\mathcal{C}$, donc sa classe dans $Ho-\mathcal{C}$
est un isomorphisme. En composant sa classe inverse par la classe de
$A_2(x,y,f)\rightarrow A^f_1(x,z)$, on obtient une classe qu'on notera $F(f)$ de
$A^f_1(x,y)$ vers $A^f_1(x,z)$, que l'on peut concevoir comme la composition
homotopique au but par $f$.\\

Montrons maintenant que $F(f)$ est indépendant du choix du représentant pour
la classe de $f$. Soit donc $g$ un autre morphisme de l'objet final vers
$A_1(y,z)$ ayant même image que $f$ par $\tau_0$, donc représentant la
même classe que $f$ dans $(x,\tau_1(A))$. Comme $A_1(y,z)$ est un \obc~comme fibre
de l'\obc~$A_1$, par hypothèse 4), il existe un morphisme de l'\obc~contractile
$K$ vers $A_1(y,z)$ envoyant 0 sur $f$ et 1 sur $g$. Soit $K'$ l'objet fibrant par
lequel $K\rightarrow A_1(y,z)\rightarrow A^f_1(y,z)$ se factorise en cofibration
triviale suivie d'une fibration. Posons $P$ l'objet fibrant produit fibré de
$A'_2(x,y,z)\rightarrow A^f_1(y,z)$ le long de $K'\rightarrow A^f_1(y,z)$. Par
produit fibré d'une fibration le long d'un morphisme , $P\rightarrow K'$ est
une fibration. Considérons le diagramme :
\begin{diagram}
A_2(x,y,f) & \rTo & * \\
\dTo & & \dTo \\
P & \rTo & K'\\
\dTo & & \dTo \\
A'_2(x,y,z) & \rTo & A^f_1(y,z)\\
\end{diagram}
Le diagramme du bas et le diagramme total sont cartésiens par définition de
$P$ et $A_2(x,y,f)$. Donc le diagramme du haut est cartésien et par suite 
le morphisme $A_2(x,y,f)\rightarrow P$
est le produit fibré de $*\rightarrow K'$ induit par 0 le long de
$P\rightarrow K'$. Comme $K$ et donc $K'$
est contractile, par hypothèse 6), $A_2(x,y,f)\rightarrow P$ est une
équivalence faible comme produit fibré d'une équivalence faible le long
d'une fibration. il en va de même pour $A_2(x,y,g)\rightarrow P$ (induit
par 1). On obtient ainsi dans $Ho-\mathcal{C}$ un isomorphisme entre
$A_2(x,y,f)$ et $A_2(x,y,g)$ qui entraîne l'égalité de $F(f)$ et $F(g)$.
Ainsi $F(f)$ est indépendant du choix d'un représentant pour la classe de
$f$.\\

Pour montrer que $F$ est un foncteur, on remarque tout d'abord que, pour tout $f:y\rightarrow z$
et $g:z\rightarrow w$, il existe $"g\circ f":y\rightarrow w$ tel que la classe
de $"g\circ f"$ soit la composée de la classe de $f$ avec celle de $g$, car
$(x,\tau_1(A))$ est
une catégorie. On définit alors une équivalence faible
$A_3(x,y,f,g)\rightarrow A^f_1(x,y)$ comme on l'a fait pour
$A_2(x,y,f)\rightarrow A^f_1(x,y)$ qui permet de montrer, comme pour l'indépendance du
choix de $F(f)$, que $F(g)\circ F(f)$ et $F("g\circ f")$ sont égales. De
même, on construira une équivalence faible $A_2(x,y,Id_y)\rightarrow
A^f_1(x,y)$, dont l'inverse homotopique est induit par l'application
de dégénéres\-cence, et qui permet de montrer que $F$ préserve
l'identité.\\

On a ainsi montré que $F$ est un foncteur. On en déduit que si
$f:*\rightarrow A_1(y,z)$ est une équivalence entre $y$ et $z$ alors il existe
$g:*\rightarrow A_1(z,y)$ tel que la composée des classes de $f$ et de $g$
ainsi que celle des classes de $g$ et de $f$ sont les identités et donc $F(f)$ et
$F(g)$ induisent des isomorphismes entre $A^f_1(x,y)$ et $A^f_1(x,z)$ dans $Ho-\mathcal{C}$. On montre
de manière symétrique que si $x$ et $x'$ sont deux objets équivalents de
$A$ et $y$ un objet fixé de $A$, alors $A^f_1(x,y)$ et $A^f_1(x',y)$ sont
isomorphes dans $Ho-\mathcal{C}$.\\

Or nous avons que $b\sim f(a)$ et $b'\sim f(a')$, donc par ce qui précède on
a dans $Ho-\mathcal{C}$ des isomorphismes entre $B^f_1(b,b')$ et
$B^f_1(f(a),f(a'))$ et entre $C^f_1(g(b),g(b'))$ et $C^f_1(g\circ f(a),g\circ
f(a'))$. Or on a vu que $g_1(f(a),f(a'))$ est une équivalence faible, donc par
l'axiome "trois pour deux" son remplacement fibrant $g^f_1(f(a),f(a')):
B^f_1(f(a),f(a'))\rightarrow C^f_1(g\circ f(a),g\circ f(a'))$ est aussi une
équivalence faible, et donc sa classe dans $Ho-\mathcal{C}$ est un
isomorphisme. On obtient donc que la classe de $g^f_1(b,b'):
B^f_1(b,b')\rightarrow C^f_1(g(b),g(b'))$ dans $Ho-\mathcal{C}$ est un
isomorphisme, donc que le remplacement fibrant de $g_1(b,b')$ est une
équivalence faible et par l'axiome "trois pour deux" que $g_1(b,b')$
elle-même est une équivalence faible entre \obcs,
ce qui par 2) en fait une \eqc~d'\obcs. $g$ est donc pleinement fidèle et comme on a déjà montré
son essentielle surjectivité, $g$ est une équivalence de \cats. On a ainsi
montré que si $f$ et $g\circ f$ sont des équivalences de \cats~alors $g$ aussi.\\
CQFD.\\

\begin{cor}\label{3pour2ccor}
Sous les hypothèse du lemme~\ref{3pour2c}, soient $f$ et $g$ deux morphismes
de \cats~tels que $f\circ g$ soit une \eq~de \cats~et $g\circ f$ soit
l'identité, alors $f$ et $g$ sont des \eqs~de \cats.
\end{cor}
{\it Preuve :}
Par fonctorialité de $\tau_0$, on a que $\tau_0(f)\circ \tau_0(g)$ est une
bijection donc que $\tau_0(f)$ est surjectif.\\
On a en outre que, pour tout couple $(x,y)$ d'objets de la source de $f$,
$f_1(x,y)\circ g_1(f(x),f(y))$ est une \eqc~d'\obcs, donc, par hypothèse 2), une équivalence faible, et que
$g_1(f(x),f(y))\circ f_1(x,y)$ est l'identité. Comme par hypothèse,
$\mathcal{C}$ est une \cmf, il vient que $f_1(x,y)$ est une équivalence
faible entre \obcs, et par 2) une \eqc~d'\obcs.\\
Donc $f$ est une équivalence de \cats~et, par le lemme~\ref{3pour2c} que
l'on vient de montrer, $g$ aussi est une équivalence de \cats.\\
CQFD.\\

Nous venons donc de montrer que sous des hypothèses assez fortes sur la
catégorie $\mathcal{C}$, comme la structure de \cmf~par exemple, les \eqs~de
\cats~vérifient les hypothèses de la proposition~\ref{Ephi3}. Il ne nous
reste donc plus qu'à trouver une construction catégorisante préservant
l'homotopie des \cats~pour pouvoir appliquer cette proposition. Pour cela nous
allons définir de nouvelles constructions formalisant différemment les
sommes amalgamées par les flèches génératrices de \cats~faciles.

\newpage

\section{Constructions Reg et Seg}

Afin de construire une catégorisation dont on saura montrer la propriété
de préservation d'homotopie des \cats, nous allons essayer de comprendre
comment les plans d'addition de flèches génératrices de \cats~faciles agissent
niveau par niveau sur les \precats. Les flèches génératrices se
répartissent en deux types : le type 1 a pour fonction de forcer les niveaux
des \precats~à être des \obcs, on dira que ces flèches régalisent, et le
type 2 a pour fonction de forcer les morphismes de Segal des \precats~à être
des \eqcs~d'\obcs, on dira que ces flèches ségalisent. Aussi va-t-on donner
deux constructions Reg et Seg qui permettront de comprendre comment la
régalisation et la ségalisation se passent au niveau de la catégorie
$\mathcal{C}$.\\

Commençons par la construction Reg. Comme le but est de régaliser un
niveau particulier d'une \precat, il va donc falloir remplacer ce niveau par son
remplacement régal. Aussi la construction Reg va-t-elle construire une
nouvelle \precat~à partir d'une \precat~dont on aura modifier l'un des
niveaux.

\begin{defin}\index{Reg}
Soit $\mathcal{C}$ une catégorie \discret, finiment cocomplète et ayant les
produits d'objets discrets. Soient $A$
une \precat, $p$ un entier strictement positif, $B$ un objet de $\mathcal{C}$ et
$f:A_p\rightarrow B$ et $g:B\rightarrow A_0\times\ldots\times A_0$ des
morphismes de $\mathcal{C}$ tels que $g\circ f:A_p\rightarrow
A_0\times\ldots\times A_0$ soit induit par les applications sommets. On
définit la \precat~$Reg(A,f,g)$ de la manière suivante :
$$Reg(A,f,g)_0=A_0 \mbox{ et, pour $q$ strictement positif, }$$
$$Reg(A,f,g)_q=A_q\coprod_{\coprod_{\Delta(q,p)^0}
A_p} \Bigg(\coprod_{\Delta(q,p)^0}B\Bigg)$$ où
$\Delta(q,p)^0$ est l'ensemble des applications de {\bf q} vers {\bf p} qui ne se
factorisent pas par {\bf 0}. On notera $B^x$ l'exemplaire de $B$ associé à
l'application $x:\mbox{{\bf q}}\rightarrow \mbox{{\bf p}}$.
Soit $y:\mbox{{\bf r}}\rightarrow \mbox{{\bf q}}$ un morphisme de $\Delta$. On
définit le morphisme induit $y^*:Reg(A,f,g)_q\rightarrow Reg(A,f,g)_r$ ainsi.
Sur la composante $A_q$ de $Reg(A,f,g)_q$, $y^*$ est le morphisme
$A(y):A_q\rightarrow A_r$. Pour la composante $B^x$, il y a deux cas. Si $x\circ
y:\mbox{{\bf r}}\rightarrow \mbox{{\bf p}}$ ne se factorise pas par 0, alors $y^*$ induit l'identité de $B^x$ à
$B^{x\circ y}$. Sinon $x\circ y$ se factorise en $\delta_i\circ \sigma^r$, où 
$\delta_i:\mbox{{\bf 0}}\rightarrow \mbox{{\bf p}}$ est l'application $i$-ème
sommet et $\sigma^r:\mbox{{\bf r}}\rightarrow \mbox{{\bf 0}}$ l'unique
application de {\bf r} vers {\bf 0}. Alors $y^*$ est la composée suivante :
\begin{diagram}
B^x & \rTo^g & A_0\times\ldots\times A_0 & \rTo^{\pi_i} & A_0 &
\rTo^{A(\sigma^r)} & A_r\\
\end{diagram}
où $\pi_i$ est la $i$-ème projection. On pourra aisément vérifier que
ceci définit bien le morphisme $y^*$ et que $Reg(A,f,g)$ est bien une
\precat~munie d'un morphisme naturel $A\rightarrow Reg(A,f,g)$.
\end{defin}

Le principe de la construction Seg est semblable à celui de la construction
Reg. Cette fois il s'agit de ségaliser un morphisme de Segal
particulier d'une \precat. Disons que l'on veut ségaliser le morphisme de
Segal de niveau $m$. Pour cela on va remplacer le niveau $m$ par un objet qui
sera en \eqc~avec le but du morphisme de Segal de niveau $m$.
La construction Seg pour $m$ va donc construire une
nouvelle \precat~à partir d'une \precat~dont on aura modifier le
niveau $m$ de manière à le mettre en \eqc~avec le but du morphisme de Segal
de niveau $m$.

\begin{defin}\index{Seg}
Soit $\mathcal{C}$ une catégorie \discret, finiment cocomplète et ayant les
produits fibrés au-dessus d'objets discrets. Soient $A$
une \precat, $m$ un entier strictement supérieur à un, $Q$ un objet de $\mathcal{C}$ et
$f:A_m\rightarrow Q$ et $g:Q\rightarrow A_1\times_{A_0}\ldots\times_{A_0}A_1$ des
morphismes de $\mathcal{C}$ tels que $g\circ f:A_m\rightarrow
A_1\times_{A_0}\ldots\times_{A_0}A_1$ soit le morphisme de Segal. On
définit la \precat~$Seg(A,f,g)$ de la manière suivante :
$$Seg(A,f,g)_0=A_0 \mbox{ et, pour $q$ strictement positif, }$$
$$Seg(A,f,g)_q=A_q\coprod_{\coprod_{\Delta(q,m)^1}
A_m} \Bigg(\coprod_{\Delta(q,m)^1}Q\Bigg)$$ où
$\Delta(q,m)^1$ est l'ensemble des applications de {\bf q} vers {\bf m} qui ne se
factorisent pas par les applications de faces principales, i.e. les applications
de {\bf 1} vers {\bf m} qui à 0 et 1 associe respectivement $i$ et $i+1$, pour
$i$ compris entre 0 et $m-1$. On notera $Q^x$ l'exemplaire de $Q$ associé à
l'application $x:\mbox{{\bf q}}\rightarrow \mbox{{\bf m}}$.
Soit $y:\mbox{{\bf r}}\rightarrow \mbox{{\bf q}}$ un morphisme de $\Delta$. On
définit le morphisme induit $y^*:Seg(A,f,g)_q\rightarrow Seg(A,f,g)_r$ ainsi.
Sur la composante $A_q$ de $Seg(A,f,g)_q$, $y^*$ est le morphisme
$A(y):A_q\rightarrow A_r$. Pour la composante $Q^x$, il y a deux cas. Si $x\circ
y:\mbox{{\bf r}}\rightarrow \mbox{{\bf m}}$ ne se factorise pas par les applications de faces principales,
alors $y^*$ induit l'identité de $Q^x$ à
$Q^{x\circ y}$. Sinon $x\circ y$ se factorise en $fp_i\circ \tau$, où 
$fp_i:\mbox{{\bf 1}}\rightarrow \mbox{{\bf m}}$ est l'application $i$-ème
face principale et $\tau:\mbox{{\bf r}}\rightarrow \mbox{{\bf 1}}$ une
application dégénérée de {\bf r} vers {\bf 1}. Alors $y^*$ est la composée suivante :
\begin{diagram}
Q^x & \rTo^g & A_1\times_{A_0}\ldots\times_{A_0}A_1 & \rTo^{\pi_{i+1}} & A_1 &
\rTo^{A(\tau)} & A_r\\
\end{diagram}
où $\pi_{i+1}$ est la $i+1$-ème projection. On pourra aisément vérifier que
ceci définit bien le morphisme $y^*$ et que $Seg(A,f,g)$ est bien une
\precat~munie d'un morphisme naturel de $A\rightarrow Seg(A,f,g)$.
\end{defin}

Le but des constructions Reg et Seg est essentiellement de formaliser la
régalisation et la ségalisation des \precats. Aussi s'attend-on à ce que
les sommes amalgamées par les flèches génératrices des \cats~faciles
s'expriment au moyen de ces constructions. 

\begin{lem}\label{flrs}
Soit $(\mathcal{C},\mathcal{F}_1,\mathcal{F}_2)$ une donnée de Segal
pré-facile.
Soient $m$ un entier
strictement positif, $A$ une \precat,  $f$ un morphisme de $\mathcal{F}_1$ et $g$ un morphisme de
$\mathcal{F}_2$. Alors il existe des morphismes $\eta:A_m\rightarrow B$,
$\nu:B\rightarrow A_0\times\ldots\times A_0$, $\phi:A_m\rightarrow P$ et $\psi:P\rightarrow
A_1\times_{A_0}\ldots\times_{A_0} A_1$ de $\mathcal{C}$ tels que $\nu\circ\eta$ soit
le morphisme induit par les applications sommets, $\psi\circ\phi$ soit le
morphisme de Segal de $A$ et tels que la somme
amalgamée de $A$ par $\Delta[m]\Theta f$ soit égale à $Reg(A,\eta,\nu)$ et
la somme amalgamée de $A$ par $Boit_m(g)$ soit égale à $Seg(A,\phi,\psi)$.
\end{lem}
{\it Preuve :}\\
Par le lemme~\ref{cons}, il vient que la somme amalgamée de $A$ par
$\Delta[m]\Theta f$ équivaut dans $\mathcal{C}$ à la somme amalgamée de
$A_m(x_0,\ldots,x_m)$ par $f$, pour un certain $m$-uplet $(x_0,\ldots,x_m)$. Prenons donc pour $\eta$ le morphisme naturel de
$A_m$ dans la somme amalgamée de $A_m$ par $f$ et pour $\nu$ le morphisme
universel de la somme amalgamée vers $A_0\times\ldots\times A_0$ induit par
les applications sommets de $A$ et le morphisme constant du but de $f$ vers $A_0\times\ldots\times A_0$ à valeur
$(x_0,\ldots,x_m)$. Par construction, $\nu\circ\eta$ est bien induit par les
applications sommets. En utilisant la formule de simplification des sommes
amalgamées suivante :
$$A\coprod_{D}\Bigg(D\coprod_{C}B\Bigg)=A\coprod_{C}B $$
on trouve facilement que, pour $q$ entier strictement positif, le
niveau $q$ de la somme amalgamée de $A$ par $\Delta[m]\Theta f$ et $Reg
(A,\eta,\nu)_q$ sont tous deux égaux à la somme amalgamée de $A_q$ par un coproduit de $f$
indexé par $\Delta(q,m)^0$.\\

Par le lemme~\ref{cons}, un morphisme de la source de $Boit_m(g)$ vers $A$
équivaut dans $\mathcal{C}$ au diagramme commutatif suivant :
\begin{diagram}
X & \rTo & A_m\\
\dTo^{g} & & \dTo_{Segal}\\
Y & \rTo & A_1\times_{A_0}\ldots\times_{A_0} A_1\\
\end{diagram}
Considérons la somme amalgamée de ce diagramme, posons alors $\phi$ le
morphisme naturel de $A_m$ dans la somme amalgamée et $\psi$ le morphisme
universel de la somme amalgamée vers $A_1\times_{A_0}\ldots\times_{A_0} A_1$.
Par construction, $\psi\circ\phi$ est bien le morphisme de Segal de $A$. En 
utilisant la formule de simplification des sommes
amalgamées rappelée ci-dessus, on trouve facilement que, pour $q$ entier strictement positif, le
niveau $q$ de la somme amalgamée de $A$ par $Boit_m(f)$ et $Seg
(A,\phi,\psi)_q$ sont tous deux égaux à la somme amalgamée de $A_q$ par un
coproduit de $g$ indexé par $\Delta(q,m)^1$.\\ 
CQFD.\\

Il est intéressant de constater que ces deux opérations Reg et Seg sont reliées entre elles : en effet faire une
opération de type Reg sur $A_m$ revient à faire une opération de type Reg
sur $A_1$ suivie d'une opération de type Seg sur $A_m$. Remarquons au passage
que du fait qu'il y ait des applications de {\bf m} vers {\bf 1} qui ne se
factorisent pas par une face principale, l'opération Seg ne laisse pas $A_1$ invariant.

\begin{lem}\label{regseg}
Soit $\mathcal{C}$ une catégorie \discret, finiment cocomplète et ayant les
produits fibrés au-dessus d'objets discrets. Soient $A$
une \precat, $m$ un entier strictement supérieur à un, $P$ un objet de $\mathcal{C}$ et
$\phi:A_m\rightarrow P$ et $\psi:P\rightarrow A_0\times\ldots\times A_0$ des
morphismes de $\mathcal{C}$ tels que $\psi\circ \phi:A_m\rightarrow
A_0\times\ldots\times A_0$ soit induit par les applications sommets. Alors il
existe des objets $B$ et $Q$ de $\mathcal{C}$, des morphismes
$\eta:A_1\rightarrow B$, $\nu:B\rightarrow A_0\times A_0$,
$\phi':Reg(A,\eta,\nu)_m\rightarrow Q$ et $\psi':Q\rightarrow
Reg(A,\eta,\nu)_1\times_{Reg(A,\eta,\nu)_0}\ldots\times_{Reg(A,\eta,\nu)_0} Reg(A,\eta,\nu)_1$
de $\mathcal{C}$ tels que $\nu\circ\eta:A_1\rightarrow A_0\times A_0$ soit
l'application induite par les morphismes source et but,
$\psi'\circ\phi':Reg(A,\eta,\nu)_m\rightarrow
Reg(A,\eta,\nu)_1\times_{Reg(A,\eta,\nu)_0}\ldots\times_{Reg(A,\eta,\nu)_0}
Reg(A,\eta,\nu)_1$ soit le morphisme de Segal et que l'on ait la relation suivante :
$$Reg(A,\phi,\psi)=Seg(Reg(A,\eta,\nu),\phi',\psi')$$
\end{lem}
{\it Preuve:}\\
Prenons pour $B$ la somme amalgamée de $A_1$ par le coproduit de
$\phi:A_m\rightarrow P$ indexé par les faces principales (qui induiront le
morphisme du coproduit de $A_m$ vers $A_1$). On posera alors $\eta$ le morphisme
naturel de $A_1$ dans la somme amalgamée $B$ et $\nu$ le morphisme universel
de la somme amalgamée $B$ vers $A_0\times A_0$ induit par les morphismes
source et but de $A_1$ vers $A_0$ et les projections du morphisme $\psi$. Par
construction, on a bien que $\nu\circ\eta$ n'est autre que le morphisme induit
par source et but.\\

Considérons maintenant la \precat~$Reg(A,\eta,\nu)$. On remarque facilement
que $Reg(A,\eta,\nu)_1$ n'est autre que $B$ et que le morphisme naturel du
coproduit de $P$ vers la somme amalgamée $B$ induit un unique morphisme de $P$ vers
$B\times_{A_0}\ldots\times_{A_0}B$, qui n'est autre que $Reg(A,\eta,\nu)_1\times_{Reg(A,\eta,\nu)_0}\ldots\times_{Reg(A,\eta,\nu)_0}
Reg(A,\eta,\nu)_1$. Prenons pour $Q$ la somme amalgamée de $P$ avec
$Reg(A,\eta,\nu)_m$ au-dessus de $A_m$, pour $\phi'$ le morphisme naturel de $Reg(A,\eta,\nu)_m$
dans la somme amalgamée $Q$ et pour $\psi'$ le morphisme universel de la somme
amalgamée $Q$ vers $B\times_{A_0}\ldots\times_{A_0}B$ induit par le morphisme
de $P$ vers $B\times_{A_0}\ldots\times_{A_0}B$ et par le morphisme de Segal de
$Reg(A,\eta,\nu)$. Par construction, on a bien que $\psi'\circ\phi'$ n'est autre
que le morphisme de Segal de $Reg(A,\eta,\nu)$.\\

Pour montrer la relation entre Reg et Seg, on utilise la propriété suivante
de simplification des sommes amalgamées : $$A\coprod_{D}\Bigg(D\coprod_{C}B\Bigg)=A\coprod_{C}B $$
Par cette formule, $Reg(A,\eta,\nu)_q$, qui est la somme amalgamée de $A_q$ par
un coproduit de $\eta:A_1\rightarrow B$ indexé par $\Delta(q,1)^0$, devient la
somme amalgamée de $A_q$ par un double coproduit de $\phi:A_m\rightarrow P$.
Ce double coproduit, étant
indexé par $\Delta(q,1)^0$ et les faces principales, n'est autre qu'un
coproduit indexé par les applications de {\bf q} vers {\bf m} se factorisant par les
faces principales mais par par {\bf 0}. Ceci montre que $Reg(A,\eta,\nu)_q$ n'est autre que la somme
amalgamée de $A_q$ avec un coproduit de $\phi:A_m\rightarrow P$ indexé par les applications de {\bf q} vers {\bf m} se factorisant par les
faces principales mais pas par {\bf 0}.\\

Toujours par la formule de simplification des sommes amalgamées écrite ci-dessus, 
$Seg(Reg(A,\eta,\nu),\phi',\psi')_q$, qui est la somme amalgamée de $Reg(A,\eta,\nu)_q$ par
un coproduit de $\phi':Reg(A,\eta,\nu)_m\rightarrow Q$ indexé par $\Delta(q,m)^1$, devient la somme amalgamée de $Reg(A,\eta,\nu)_q$ par
un coproduit de $\phi:A_m\rightarrow P$ indexé par $\Delta(q,m)^1$ . En utilisant l'expression simplifiée de $Reg(A,\eta,\nu)_q$, on
obtient que $Seg(Reg(A,\eta,\nu),\phi',\psi')_q$ est la somme amalgamée de $A_q$
 avec le coproduit d'un coproduit de $\phi:A_m\rightarrow P$ indexé par $\Delta(q,m)^1$ avec un
 coproduit de $\phi:A_m\rightarrow P$ indexé par les applications de {\bf q} vers {\bf m} se factorisant par les
faces principales mais pas par {\bf 0}. Or ce coproduit de
coproduit de $\phi$ se simplifie en un simple coproduit indexé par
l'ensemble des applications de {\bf q} vers {\bf m} qui soit se factorisent par
les faces principales mais pas par {\bf 0} soit ne se factorisent pas par les
faces principales (éléments de $\Delta(q,m)^1$) donc pas par {\bf 0} non
plus. Or cet ensemble n'est autre que $\Delta(q,m)^0$, i.e. celui des
applications de {\bf q} vers {\bf m} ne se factorisant pas par {\bf 0}.
Ainsi $Seg(Reg(A,\eta,\nu),\phi',\psi')_q$ est donc la somme amalgamée de $A_q$
par un coproduit de $\phi:A_m\rightarrow P$ indexé par $\Delta(q,m)^0$, ce qui
est par définition $Reg(A,\phi,\psi)_q$.\\
CQFD.\\

Cette remarque fondamentale sur les constructions Reg et Seg nous invite à
partir de maintenant à lier fortement les opérations Reg de niveau 1 et $m$
avec l'opération Seg de niveau $m$. Nous allons donc dans la section suivante
envisager une construction prenant en compte ces trois opérations.

\newpage

\section{Construction RS}

Nous avons défini dans la section précédente deux constructions Reg et Seg dont on a vu les relations entre elles.
Nous allons donc maintenant définir une nouvelle construction RS qui est en même
temps une opération de type Reg sur le niveau 1 et de type Seg sur le niveau
$m$. Cette opération aura la bonne propriété d'être stable par composition
mais également de garder la trace de certains sous-objets des niveaux 1 et
$m$. Ce dernier point est assez crucial car c'est lui qui va permettre de suivre
les variations de l'homotopie des \cats~dans les plans d'addition de flèches
génératrices de \cats~faciles. C'est
pourquoi nous allons définir tout d'abord la notion de \precat~(1,m)-peinte qui
permet de prendre en compte des parties, dites peintes, des niveaux 1 et $m$.

\begin{defin}\index{\precat!peinte}\index{morphisme de \precats!peintes}
Soit $\mathcal{C}$ une catégorie \discret~ayant les
produits fibrés au-dessus d'objets discrets et soit $m$ un entier strictement
supérieur à un. Une \precat~(1,m)-peinte est la donnée d'un triplet $(A,i,j)$ constitué d'une
\precat~$A$ et de deux monomorphismes $i:A_1^*\rightarrow A_1$ et
$j:A_m^*\rightarrow A_m$ tels que le morphisme de Segal de $A$ se restreigne en
un morphisme de $A_m^*$ vers $A_1^*\times_{A_0}\ldots\times_{A_0} A_1^*$.\\
Un morphisme de \precats~peinte $f:(A,i,j)\rightarrow (B,i',j')$ est un
morphisme de \precats~tel que $f(i)$ se factorise par $i'$ et $f(j)$ se
factorise par $j'$.
\end{defin}

\begin{defin}\index{RS}
Soit $\mathcal{C}$ une catégorie \discret, finiment cocomplète et ayant les
produits fibrés au-dessus d'objets discrets. Soient $m$ un entier strictement
supérieur à un, $(A,i,j)$
une \precat~(1,m)-peinte, $B$ et $P$ des objets de $\mathcal{C}$ et
$\eta:A_1^*\rightarrow B$, $\nu:B\rightarrow A_0\times A_0$,
$\phi:A_m^*\rightarrow P$ et $\psi:P\rightarrow B\times_{A_0}\ldots\times_{A_0}
B$ des morphismes de $\mathcal{C}$ tels que $\nu\circ\eta:A_1^*\rightarrow
A_0\times A_0$ soit le morphisme induit par source et but et $\psi\circ
\phi:A_m^*\rightarrow
B\times_{A_0}\ldots\times_{A_0}B$ soit la composée de la restriction à
$A_m^*$ du morphisme de Segal par le produit m fois du morphisme $\eta$. 
Posons $\eta':A_1\rightarrow B'$ la somme amalgamée d'$\eta$ le long du
monomorphisme $i$, $\nu':B'\rightarrow A_0\times A_0$ le morphisme universel
induit par $\nu$ et les morphismes source et but, $\phi':A_m\rightarrow P'$ la
somme amalgamée de $\phi$ le long du monomorphisme $j$, $\psi':P'\rightarrow
B'\times_{A_0}\ldots\times_{A_0}B'$ le morphisme universel induit par $\psi$,
le morphisme de Segal et un produit m fois de $\eta'$,
$\phi'':Reg(A,\eta',\nu')_m\rightarrow P''$ la somme amalgamée de $\phi'$ le
long du morphisme naturel $A_m\rightarrow Reg(A,\eta',\nu')_m$ et
$\psi'':P''\rightarrow B'\times_{A_0}\ldots\times_{A_0}B'$ le morphisme
universel induit par $\psi'$ et le morphisme de Segal de $Reg(A,\eta',\nu')$.
On définit la \precat~(1,m)-peinte $(RS(A,\eta,\nu,\phi,\psi),k,l)$ de la manière suivante :
$$RS(A,\eta,\nu,\phi,\psi)=Seg(Reg(A,\eta',\nu'),\phi'',\psi'')$$
avec $k$ et $l$ les monomorphismes naturels de $B$ dans
$RS(A,\eta,\nu,\phi,\psi)_1$, respectivement de $P$
dans $RS(A,\eta,\nu,\phi,\psi)_m$. Il est facile de voir que ceci définit bien une \precat~(1,m)-peinte
munie d'un morphisme naturel provenant de $A$.
\end{defin}

Cette définition de la construction RS étant un peu complexe, donnons donc
de la construction RS une expression plus simple.

\begin{lem}
Sous les hypothèses et les notations de la définition précédente,
$RS(A,\eta,\nu,\phi,\psi)_0$ n'est autre que $A_0$ et, pour tout entier $q$
strictement positif, $RS(A,\eta,\nu,\phi,\psi)_q$ est la somme amalgamée de
$A_q$ par le coproduit d'un coproduit de $\eta$ indexé par $\Delta(q,1)^0$
avec un coproduit de $\phi$ indexé par $\Delta(q,m)^1$. Les morphismes sont
définis de manière similaire à ceux des constructions Reg et Seg.
\end{lem}
{\it Preuve :} directe après un petit calcul utilisant la formule de simplification
des sommes amalgamées:
$$A\coprod_{D}\Bigg(D\coprod_{C}B\Bigg)=A\coprod_{C}B .$$\\

Comme on l'a annoncé plus haut, cette construction RS est stable par composition mais
aussi par colimite séquentielle transfinie comme le montrent les deux lemmes
suivants.

\begin{lem}
Soit $\mathcal{C}$ une catégorie \discret, finiment cocomplète et ayant les
produits fibrés au-dessus d'objets discrets. Soient $m$ un entier strictement
supérieur à un, $(A,i,j)$
une \precat~(1,m)-peinte, $B$, $C$, $P$ et $Q$ des objets de $\mathcal{C}$ et
$\eta:A_1^*\rightarrow B$, $\nu:B\rightarrow A_0\times A_0$,
$\tilde{\eta}:B\rightarrow C$, $\tilde{\nu}:C\rightarrow A_0\times A_0$, 
$\phi:A_m^*\rightarrow P$, $\psi:P\rightarrow B\times_{A_0}\ldots\times_{A_0}
B$, $\tilde{\phi}:P\rightarrow Q$ et $\tilde{\psi}:Q\rightarrow C\times_{A_0}\ldots\times_{A_0}
C$ des morphismes de $\mathcal{C}$ tels que $\nu\circ\eta:A_1^*\rightarrow
A_0\times A_0$ soit le morphisme induit par source et but, $\tilde{\nu}\circ\tilde{\eta}:B\rightarrow
A_0\times A_0$ soit le morphisme $\nu$, $\psi\circ\phi:A_m^*\rightarrow
B\times_{A_0}\ldots\times_{A_0}B$ soit la composée de la restriction à
$A_m^*$ du morphisme de Segal par le produit m fois du morphisme $\eta$ et 
$\tilde{\psi}\circ\tilde{\phi}:P\rightarrow
C\times_{A_0}\ldots\times_{A_0}C$ soit la composée du morphisme $\psi$ par le
produit m fois du morphisme $\tilde{\eta}$. Alors on a la relation suivante :
$$RS(RS(A,\eta,\nu,\phi,\psi),\tilde{\eta},\tilde{\nu},\tilde{\phi},\tilde{\psi})=RS(A,\tilde{\eta}\circ\eta,
\tilde{\nu}\circ\nu,\tilde{\phi}\circ\phi,\tilde{\psi}\circ\psi) $$
\end{lem}
{\it Preuve :} application directe de la formule de simplification des sommes
amalgamées et du lemme décrivant explicitement $RS(A,\eta,\nu,\phi,\psi)_q$.

\begin{lem}\label{RScolim}
Soit $\mathcal{C}$ une catégorie \discret, cocomplète et ayant les
produits fibrés au-dessus d'objets discrets. Soient $m$ un entier strictement
supérieur à un, $(A,i,j)$
une \precat~(1,m)-peinte, $\lambda$ un cardinal transfini,
$\eta:A_1^*\rightarrow B$ la composée d'une $\lambda$-séquence dans
$\mathcal{C}$ de $\eta_i:B^i\rightarrow B^{i+1}$ et $\phi:A_m^*\rightarrow P$ la composée d'une $\lambda$-séquence dans
$\mathcal{C}$ de $\phi_i:P^i\rightarrow P^{i+1}$. Soient pour tout $i$ inférieur ou
égal à $\lambda$ des morphismes $\nu_i:B^i\rightarrow A_0\times A_0$ et
$\psi_i:P^i\rightarrow B^i\times_{A_0}\ldots\times_{A_0}B^i$ de $\mathcal{C}$
tels que, pour tout $i$ cardinal non limite inférieur ou égal à $\lambda$,
$\nu_i\circ\colimite{j\leq i} \eta_j$ soit le morphisme induit par source et but et
$\psi_i\circ\colimite{j\leq i} \phi_j$ soit le morphisme induit par le morphisme de
Segal de $A$ et par $\colimite{j\leq i} \eta_j$, et tels que, pour tout $i$
cardinal limite inférieur ou égal à $\lambda$, d'une part les morphismes
$\nu_i$ et $\psi_i$ sont les morphismes universels induits respectivement par
les $\nu_j$ et les $\psi_j$ (avec $j$ strictement inférieur à $i$), et d'autre
part,  $\nu_i\circ\colimite{j<i} \eta_j$ soit le morphisme induit par source et but et
$\psi_i\circ\colimite{j<i} \phi_j$ soit le morphisme induit par le morphisme de
Segal de $A$ et par $\colimite{j<i} \eta_j$. 
Posons $RS^0=A$. Pour $i$ cardinal inférieur strictement à
$\lambda$, posons $RS^{i+1}=RS(RS^i,\eta_i,\nu_i,\phi_i,\psi_i)$ avec $B^{i+1}$
et $P^{i+1}$ pour partie peinte. Pour $i$ cardinal limite inférieur ou égal
à $\lambda$, posons $RS^i=\colimite{j<i}RS^j$ avec pour partie peinte $B^i$ et
$P^i$. Alors on a les relations suivantes :
\item -si $\lambda$ est un cardinal limite,  
$$ \colim{i<\lambda} RS^i=RS(A,\eta,\nu,\phi,\psi),$$
\item -si $\lambda$ n'est pas un cardinal limite,  
$$ \colim{i\leq\lambda} RS^i=RS(A,\eta,\nu,\phi,\psi).$$  
\end{lem}
{\it Preuve :} directe à partir de la forme explicite de
$RS(A,\eta,\nu,\phi,\psi)$ en utilisant la commutativité des colimites.\\

La construction RS a été voulue comme une sorte de généralisation des opérations Reg et
Seg. Nous allons montrer que les constructions Reg et Seg peuvent donc bien
s'exprimer au moyen de la construction RS.

\begin{lem}
Soit $\mathcal{C}$ une catégorie \discret, finiment cocomplète et ayant les
produits fibrés au-dessus d'objets discrets. Soient $m$ un entier strictement
supérieur à un et $A$ une \precat.
On peut exprimer les opérations Reg et Seg en terme d'opération RS de la
manière suivante :
\item - Soient $f:A_1\rightarrow B$ et $g:B\rightarrow A_0\times A_0$ tels que
$g\circ f$ soit induit par source et but, $Reg(A,f,g)$ n'est autre que
$RS(A,\eta,\nu,\phi,\psi)$ avec $f$ pour $\eta$, $g$ pour $\nu$, des identités
pour $i$, $j$ et $\phi$, et la composée du morphisme de Segal avec un produit m fois de $f$
pour $\psi$.
\item - Soient $f:A_m\rightarrow P$ et $g:P\rightarrow
A_1\times_{A_0}\ldots\times_{A_0}A_1$ tels que
$g\circ f$ soit le morphisme de Segal, $Seg(A,f,g)$ n'est autre que
$RS(A,\eta,\nu,\phi,\psi)$ avec $f$ pour $\phi$, $g$ pour $\psi$, des identités
pour $i$, $j$ et $\eta$, le morphisme induit par source et but pour $\nu$. 
\end{lem}
{\it Preuve :} directe à partir de la forme explicite de
$RS(A,\eta,\nu,\phi,\psi)_q$.\\

Comme on l'a déjà vu l'opération Reg au niveau $m$ est la composée d'une
opération Reg au niveau 1 avec une opération Seg au niveau $m$. Comme toutes
deux sont des opérations RS et que l'opération RS est stable par
composition, l'opération Reg au niveau $m$ est aussi une opération RS.\\

Notre construction RS permet donc d'exprimer les constructions Reg et Seg qui
expriment quant à elles les sommes amalgamées par les flèches
génératrices de \cats~faciles. En outre cette construction RS permet de suivre à la trace
certaines parties des niveaux 1 et $m$. C'est pourquoi nous allons donc regarder certains plans
d'addition de flèches génératrices de \cats~faciles qui s'expriment au
moyen de la construction RS, ce qui nous permettra de montrer qu'ils
préservent l'homotopie des \cats.

\newpage

\section{Construction Cat(1,m)}

Comme on l'a vu dans les sections précédentes, les sommes amalgamées par
des flèches génératrices de \cats~faciles de type 1 pour les niveaux 1 et
$m$ et celles avec des flèches de type 2 pour le niveau $m$ sont compatibles
et s'expriment au moyen de la construction RS qui est stable par colimite
séquentielle transfinie et garde la trace de parties peintes des niveaux 1 et
$m$. Aussi va-t-on considérer le plan d'addition de cellules
$e_{\Phi,\lambda}$ avec pour $\Phi$ la famille des flèches ci-dessus et
montrer que ce plan n'est autre qu'une construction RS.

\begin{defin}\index{$Raj(1,m)$}
Soit $(\mathcal{C},\mathcal{F}_1,\mathcal{F}_2)$ une donnée de Segal
pré-facile. Soit $m$ un entier strictement supérieur à un. 
Notons $\Phi$ la famille de morphismes de \precats~constituée des flèches de
type $\Delta[1]\Theta f$, $\Delta[m]\Theta f$ et $Boit_m(g)$, avec
$f$ dans $\mathcal{F}_1$ et $g$ dans $\mathcal{F}_2$.\\
Soit $(A,i,j)$ une \precat~(1,m)-peinte.
Définissons le plan d'addition de cellules simple $Raj(1,m)$ pour $(A,i,j)$
comme le plan simple $e_{\Phi,1}$ appliqué à $A$ mais dans lequel on ne
garde que les diagrammes dont les sources des flèches de $\Phi$ s'envoyent
dans les parties peintes de $A$. 
\end{defin}

\begin{lem}\label{raj1m}
Soit $(\mathcal{C},\mathcal{F}_1,\mathcal{F}_2)$ une donnée de Segal
pré-facile. Soit $m$ un entier strictement supérieur à un.\\
Pour toute $(A,i,j)$ une \precat~(1,m)-peinte, il existe des morphismes
$\eta:A_1\rightarrow B$, $\nu:B\rightarrow A_0\times A_0$, $\phi:A_m\rightarrow
P$ et
$\psi:P\rightarrow B\times_{A_0}\ldots\times_{A_0}B$ de $\mathcal{C}$ tels que
$\nu\circ\eta$ soit induit par les morphismes source et but, $\psi\circ\phi$
soit induit par le morphisme de Segal de $A$ et par $\eta$ et tels qu'on ait :
$$Raj(1,m)(A,i,j)=RS(A,\eta,\nu,\phi,\psi)$$
Cette égalité fait de $Raj(1,m)(A,i,j)$ une \precat~(1,m)-peinte dont les parties peintes
sont $B$ et $P$. En outre $\eta$ est une colimite
séquentielle transfinie de sommes amalgamées de $A_1$ par les flèches de
$\mathcal{F}_1$ et $\phi$ une colimite séquentielle transfinie de
sommes amalgamées de $A_m$ par les flèches de $\mathcal{F}_1$ et
$\mathcal{F}_2$.
\end{lem}
{\it Preuve :}\\
Par définition, $Raj(1,m)(A,i,j)$ est une somme amalgamée de $A$ par un
coproduit de flèches génératrices de \cats~faciles, ce qui peut se mettre sous la forme d'une colimite
séquentielle transfinie de sommes amalgamées par une flèche
génératrice de \cats~faciles. Or on a vu au lemme~\ref{flrs} que les sommes amalgamée par des flèches
génératrices de \cats~faciles ne sont autres que des opérations de type RS. Comme en outre
les flèches génératrices de \cats~faciles par lesquelles on fait les sommes amalgamées
ont leurs sources s'envoyant dans les parties peintes de $A$, les opérations
RS sont bien composables et $Raj(1,m)(A,i,j)$ forme ainsi une colimite
séquentielle transfinie d'opérations RS au sens du lemme~\ref{RScolim} et
donc, par ce lemme, $Raj(1,m)(A,i,j)$ est bien une opération de type RS. La
démonstration du lemme montre que le morphisme $\eta$ (respectivement $\phi$) est une colimite
séquentielle transfinie des morphismes $\eta_i$ (respectivement $\phi_i$). Ces
morphismes $\eta_i$ et $\phi_i$
proviennent de l'expression des sommes amalgamées par les flèches
génératrices de \cats~faciles. Une somme amalgamée par une flèche de type
$\Delta[1]\Theta f$, où $f$ appartient à la famille $\mathcal{F}_1$, est une
opération RS avec pour $\eta$ le morphisme naturel de $A_1^*$ dans la somme
amalgamée de $A_1^*$ par $f$ et pour $\phi$ l'identité. Une somme amalgamée par une flèche de type
$\Delta[m]\Theta f$, où $f$ appartient à la famille $\mathcal{F}_1$, est une
opération RS avec pour $\eta$ le morphisme naturel de $A_1^*$ dans la somme amalgamée
de $A_1^*$ par un coproduit de $f$ et pour $\phi$ le morphisme naturel de $A_m^*$ dans la somme
amalgamée de $A_m^*$ par $f$. Enfin une somme amalgamée par une flèche de type
$Boit_m(g)$, où $g$ appartient à la famille $\mathcal{F}_2$, est une
opération RS avec pour $\eta$ l'identité et pour $\phi$ le morphisme naturel de $A_m^*$ dans la somme
amalgamée de $A_m^*$ par $g$. Ainsi $\eta$ est bien une colimite
séquentielle transfinie de sommes amalgamées par des flèches de
$\mathcal{F}_1$ et $\phi$ une colimite
séquentielle transfinie de sommes amalgamées par des flèches de
$\mathcal{F}_1$ et $\mathcal{F}_2$.\\
CQFD.\\

Ainsi le plan simple $Raj(1,m)$ transforme des \precats~peintes en
\precats~peintes. Il serait très pratique que cette transformation soit
fonctorielle.

\begin{lem}\label{raj1mf}
Soit $(\mathcal{C},\mathcal{F}_1,\mathcal{F}_2)$ une donnée de Segal
pré-facile. Soit $m$ un entier strictement supérieur à un.\\
La construction $Raj(1,m)$ est en fait un foncteur de la catégorie des
\precats~peintes vers elle-même. De même, tout plan d'addition de cellules
constitué de plans simples de type $Raj(1,m)$ est un foncteur de la catégorie des
\precats~peintes vers elle-même.
\end{lem}
{\it Preuve :}\\
Pour montrer le premier résultat il suffit d'appliquer le critère du
lemme~\ref{plfcrit}. Soit $f:(A,i,j)\rightarrow (B,i',j')$ un morphisme de
\precats~peintes, l'extension de $Raj(1,m)_A$ n'est autre que le plan simple d'addition de
cellules constitué par les couples $(D,1)$ où $D$ est décrit les
extensions par $f$ de diagrammes de $\Phi$ s'envoyant dans les parties peintes
de $A$. Comme $f$ est un morphisme de \precats~peintes, les extensions par $f$
de diagrammes à valeur dans les parties peintes de $A$ sont justement des
diagrammes à valeur dans les parties peintes de $B$. Ceci montre que ces
couples $(D,1)$ sont bien des
couples du plan simple $Raj(1,m)_B$, d'où il vient que
$Raj(1,m)$ vérifie la condition suffisante pour être un plan
d'addition de cellules fonctoriel.\\
Le second résultat résulte du premier auquel on applique le
lemme~\ref{plsfonc}.\\
CQFD.\\

Nous avons maintenant un plan d'addition de cellules fonctoriel $Raj(1,m)$ qui
est en fait une opération de type RS. Si l'on veut une $\Phi$-injectivisation
qui soit aussi une opération du type RS, il suffit de prendre un plan
d'addition de cellules composé uniquement de $Raj(1,m)$ et qui soit de longueur
suffisamment grande. Comme on a vu que l'opération RS est stable par colimite
séquentielle transfinie, on obtiendra qu'un tel plan est encore de type RS.
 
\begin{defin}\index{$Cat(1,m)$}
Soit $(\mathcal{C},\mathcal{F}_1,\mathcal{F}_2)$ une donnée de Segal
pré-facile. Soit $\alpha$ un cardinal transfini régulier strictement
supérieur à celui qui assure la petitesse des sources et buts des flèches des familles $\mathcal{F}_1$ et
$\mathcal{F}_2$. Soient $m$ un entier strictement supérieur à un et
$A$ une \precat. Nous allons définir $Cat(1,m)(A)$ par récurrence
transfinie.\\
Posons $A^0=A$ avec pour partie peinte $A_1$ et $A_m$.
Supposons $(A^{\beta},i_{\beta},j_{\beta})$ défini pour un ordinal
$\beta<\alpha$, alors on pose :
$$A^{\beta+1}=Raj(1,m)(A^{\beta},i_{\beta},j_{\beta}),$$
ses parties peintes étant celle du lemme~\ref{raj1m}.
Soit un ordinal limite $\gamma<\alpha$, supposons $(A^{\beta},i_{\beta},j_{\beta})$ défini pour tout ordinal
$\beta<\gamma$, alors on pose :
$$(A^{\gamma},i_{\gamma},j_{\gamma})=\Bigg(\colim{\beta<\gamma}
A^{\beta},\colim{\beta<\gamma}i^{\beta},\colim{\beta<\gamma}j^{\beta}\Bigg)$$
Par récurrence transfinie, on a bien défini pour tout ordinal $\beta<\alpha$
une \precat~(1,m)-peinte $(A^{\beta},i_{\beta},j_{\beta})$.
Posons donc :
$$ Cat(1,m)(A)=\colim{\beta<\alpha} A^{\beta}$$.
\end{defin}

\begin{lem}\label{cat1m}
Soit $(\mathcal{C},\mathcal{F}_1,\mathcal{F}_2)$ une donnée de Segal
pré-facile.
Soient $m$ un entier strictement supérieur à un et
$A$ une \precat.
Alors il existe des morphismes
$\eta:A_1\rightarrow B$, $\nu:B\rightarrow A_0\times A_0$, $\phi:A_m\rightarrow
P$ et
$\psi:P\rightarrow B\times_{A_0}\ldots\times_{A_0}B$ de $\mathcal{C}$ tels que
$\nu\circ\eta$ soit induit par les morphismes source et but, $\psi\circ\phi$
soit induit par le morphisme de Segal de $A$ et par $\eta$ et tels qu'on ait :
$$Cat(1,m)(A)=RS(A,\eta,\nu,\phi,\psi).$$
En outre $B$ et $P$ sont des \obcs~faciles, $\psi$ une \eqc~facile, $\eta$ est une colimite
séquentielle transfinie de sommes amalgamées de $A_1$ par les flèches de
$\mathcal{F}_1$ et $\phi$ une colimite séquentielle transfinie de
sommes amalgamées de $A_m$ par les flèches de $\mathcal{F}_1$ et
$\mathcal{F}_2$.
\end{lem}
{\it Preuve :}\\
Par le lemme~\ref{raj1m}, $Cat(1,m)(A)$, étant une colimite séquentielle transfinie
d'opérations $Raj(1,m)$, s'exprime comme une colimite séquentielle transfinie
d'opérations RS au sens du lemme~\ref{RScolim}. Donc d'après ce lemme,
$Cat(1,m)(A)$ est une opération de type RS avec pour $\eta:A_1\rightarrow B$ (respectivement
$\phi:A_m\rightarrow P$) une colimite
séquentielle transfinie des morphismes $\eta$ (respectivement $\phi$) apparaissant dans $Raj(1,m)$,
eux-mêmes colimites séquentielles transfinies de sommes amalgamées par des
flèches de $\mathcal{F}_1$ (respectivement par des flèches de
$\mathcal{F}_1$ et $\mathcal{F}_2$), ce qui nous donne la caractérisation des
morphismes $\eta$ et $\phi$.\\

Montrons que $B$ est un \obc~facile, i.e. se relève par rapport à la famille
$\mathcal{F}_1$. Soit $f:X\rightarrow Y$ une flèche de $\mathcal{F}_1$ et
$x:X\rightarrow B$ un  morphisme quelconque de $\mathcal{C}$. Comme par hypothèse $X$ est
petit pour un cardinal régulier strictement inférieur à $\alpha$ et que $B$ est une colimite séquentielle transfinie de $B^i$,
pour $i<\alpha$, $x$ se factorise à travers $B^{\beta}$ pour un certain
$\beta<\alpha$. Ceci équivaut à un morphisme de $\Delta[1]\Theta X$ vers
$A^{\beta}$ ne concernant que les parties peintes. Par définition de
$A^{\beta+1}$ comme $Raj(1,m)(A^{\beta},i_{\beta},j_{\beta})$, $A^{\beta+1}$ est la
somme amalgamée de $A^{\beta}$ par toutes les flèches de type
$\Delta[1]\Theta h$, $\Delta[m]\Theta h$ et $Boit_m(g)$ dont les sources
s'envoient dans les parties peintes, avec $h$ dans $\mathcal{F}_1$ et $g$ dans
$\mathcal{F}_2$. Donc, par le lemme~\ref{raj1m}, sa partie peinte $B^{\beta+1}$ est une colimite
séquentielle transfinie de sommes amalgamées par des flèches de $\mathcal{F}_1$. En particulier, la somme
amalgamée de $A^{\beta}$ par $\Delta[1]\Theta f$ correspond à un $\eta$ qui n'est autre que
le morphisme naturel de $B^{\beta}$ dans la somme amalgamée de $B^{\beta}$ par
$f$. Donc on obtient le diagramme commutatif suivant :
\begin{diagram}
X & \rTo^{x} & B^{\beta}\\
\dTo^{f} & & \dTo\\
Y & \rTo_{\exists \tilde{x}} & B^{\beta+1}\\
\end{diagram}
L'extension de $\tilde{x}$ à $B$ est donc le relèvement attendu. Donc $B$ se
relève par rapport à la famille $\mathcal{F}_1$, ce qui en fait un
\obc~facile. On démontrera de même, en utilisant $\Delta[m]\Theta f$, pour
$f$ dans $\mathcal{F}_1$, que $P$ est un \obc~facile.\\

Montrons maintenant que $\psi:P\rightarrow B\times_{A_0}\ldots\times_{A_0}B$ est
une \eqc~facile d'\obcs, i.e. se relève par rapport à la famille
$\mathcal{F}_2$. On remarque que $B\times_{A_0}\ldots\times_{A_0}B$ est
bien un \obc~comme produit fibré d'\obcs~au-dessus d'objets discrets.
Considérons le diagramme commutatif suivant :
\begin{diagram}
X & \rTo^x & P\\
\dTo^{g} & & \dTo_{\psi}\\
Y & \rTo_y & B\times_{A_0}\ldots\times_{A_0}B\\
\end{diagram} 
avec $g$ dans $\mathcal{F}_2$, $x$ et $y$ quelconques.
Comme $X$ et $Y$, par hypothèse, sont petits pour un cardinal régulier
strictement inférieur à $\alpha$ et que $\psi$ est le
morphisme universel induit par les $\psi_i$ pour $i<\alpha$, ce diagramme se
factorise à travers $\psi_{\beta}$, pour un certain $\beta<\alpha$. Ce nouveau
diagramme correspond à un morphisme de la source du morphisme $Boit_m(g)$ vers
$A^{\beta}$. Par un argument similaire au précédent, $P^{\beta+1}$ est une
colimite séquentielle transfinie de sommes amalgamées de flèches de
$\mathcal{F}_1$ et $\mathcal{F}_2$ et, en particulier, de $g$ qui provient de la
somme amalgamée de $A^{\beta}$ par $Boit_m(g)$. De ce fait, l'extension du
diagramme avec $\psi_{\beta}$ en diagramme avec $\psi_{\beta+1}$ possède un
relèvement, i.e. une flèche de $Y$ vers $P^{\beta+1}$ faisant commuter les
deux parties du diagramme, et son extension à $P$ est bien un relèvement du
diagramme de départ. Donc $\psi$ se relève par rapport à la famille
$\mathcal{F}_2$, ce qui en fait une \eqc~facile.\\ 
CQFD.\\

Ce plan d'addition de cellules $Cat(1,m)$ s'exprime donc en terme de
cons\-truction RS mais avec de bonnes propriétés : en effet les parties
peintes de $Cat(1,m)$ sont des \obcs~et la restriction du morphisme de Segal aux
parties peintes est une \eqc~d'\obcs. Ce sont justement ces bonnes
propriétés qui permettent à $Cat(1,m)$ de préserver l'homotopie des
\cats.

\begin{lem}\label{eqcat1m}
Soit $(\mathcal{C},\mathcal{F}_1,\mathcal{F}_2)$ une donnée de Segal
pré-facile vérifiant les propriétés suivantes :
\item 1) la catégorie sous-jacente $\mathcal{C}$ est une \cmf,
\item 2) les \eqcs~d'\obcs~sont exactement les \eqs~faibles de la
\cmf~$\mathcal{C}$ entre \obcs,
\item 3) la famille $\mathcal{F}_1$ est incluse dans la classe des cofibrations
triviales de la \cmf~$\mathcal{C}$ et la famille $\mathcal{F}_2$ dans celle des
cofibrations,
\item 4) les \eqs~faibles sont stables par produit fibré au-dessus d'un objet
discret.\\
\\
Alors pour toute \precat~$C$ dont les morphismes de Segal sont des
équi\-valences faibles de $\mathcal{C}$, on a que pour tout entier $m\geq 2$, le morphisme canonique $C\rightarrow Cat(1,m)(C)$ est
niveau par niveau une cofibration triviale de $\mathcal{C}$ et $Cat(1,m)(C)$ conserve la propriété que ses
morphismes de Segal sont des \eqs~faibles de $\mathcal{C}$.
\end{lem}
{\it Preuve :}\\
Soit $C$ une \precat~dont tous les morphismes de Segal sont des
\eqs~fai\-bles de $\mathcal{C}$. On a tout d'abord que $Cat(1,m)(C)_0$ n'est
autre que $C_0$ car les flèches génératrices de \cats~faciles sont des
isomorphismes au niveau des objets.  
Par le lemme~\ref{cat1m}, l'opération $Cat(1,m)(C)$, pour un entier $m\geq 2$,
n'est autre qu'une opération du type RS avec $B$ et $P$ des \obcs~faciles,
$\eta$ une colimite séquentielle transfinie de sommes amalgamées de
flèches de la famille $\mathcal{F}_1$, $\psi$ une \eqc~d'\obcs~et $\phi$ une colimite
séquentielle transfinie de sommes amalgamées de flèches des familles
$\mathcal{F}_1$ et $\mathcal{F}_2$.\\ 

Comme par hypothèse la famille
$\mathcal{F}_1$ est incluse dans la classe des cofibrations triviales, classe
qui, dans la \cmf~$\mathcal{C}$, est 
stable par colimite séquentielle transfinie et par somme amalgamée le long
d'un morphisme, il vient que $\eta$ est une cofibration triviale. De même les
familles $\mathcal{F}_1$ et $\mathcal{F}_2$ étant incluse dans la classe des
cofibrations, classe stable dans $\mathcal{C}$ par colimite séquentielle transfinie et par somme amalgamée le long
d'un morphisme, il vient que $\phi$ est une cofibration. Par hypothèse
sur $C$, son morphisme de Segal $\segal{C}{m}$ est une \eq~faible. En outre, les
\eqs~faibles sont stables par produit fibré au-dessus d'un objet discret,
donc le produit fibré $m$ fois de $\eta$ au-dessus de $C_0$ est une
\eq~faible. Par axiome "trois pour deux" dans la \cmf~$\mathcal{C}$, la
composée du morphisme de Segal de $C$ par ce produit fibré est encore une
\eq~faible. Or par le lemme~\ref{cat1m}, cette composée est égale à
$\psi\circ\phi$ et comme $\psi$ est une \eqc~facile, toujours par "trois pour
deux", $\phi$ est aussi une \eq~faible.\\

Pour tout entier $q>0$, le morphisme
canonique $C_q\rightarrow Cat(1,m)(C)_q$ est la somme amalgamée de $C_q$ par
un coproduit multiple de $\eta$ et de $\phi$, donc ce morphisme est une colimite
séquentielle transfinie de sommes amalgamées de $\eta$ et de $\phi$. Comme
$\eta$ et $\phi$ sont des cofibrations triviales et que les cofibrations
triviales sont stables par sommes amalgamées le long d'un morphisme et par
colimite séquentielle transfinie, il vient que $C_q\rightarrow Cat(1,m)(C)_q$
est une cofibration triviale.\\ 

Considérons pour finir le diagramme commutatif suivant :
\begin{diagram}
C_q & \rTo^{Segal} & C_1\times_{C_0}\ldots\times_{C_0}C_1\\
\dTo & & \dTo\\
Cat(1,m)(C)_q & \rTo_{Segal} & Cat(1,m)(C)_1\times_{C_0}\ldots\times_{C_0}Cat(1,m)(C)_1\\
\end{diagram}
On vient de prouver que le morphisme vertical de gauche est une \eq~fai\-ble. En outre, par stabilité des \eqs~faibles par produit fibré au-dessus d'un objet
discret, le morphisme vertical de droite est aussi une \eq~faible. Comme par hypothèse les morphismes de Segal de $C$ sont des \eqs~faibles, alors par "trois pour deux"
le morphisme de Segal de $Cat(1,m)(C)$ aussi.\\
CQFD.\\

Pour que notre construction $Cat(1,m)$ soit parfaite, nous n'avons plus qu'à
lui demander d'être fonctorielle, ce qui est l'objet du lemme suivant.

\begin{lem}\label{cat1mf}
Soit $(\mathcal{C},\mathcal{F}_1,\mathcal{F}_2)$ une donnée de Segal
pré-facile. Soit $m$ un entier strictement supérieur à un.\\
La construction $Cat(1,m)$ est en fait un foncteur de la catégorie des
\precats~vers elle-même. De même tout plan d'addition de cellules
constitué de plans d'addition de cellules de type $Cat(1,m)$ est un foncteur de la catégorie des
\precats~vers elle-même.
\end{lem}
{\it Preuve :}\\
D'après le lemme~\ref{raj1mf}, tout plan d'addition de cellules constitué de
plans simples de type $Raj(1,m)$ est un foncteur de la catégorie des
\precats~peintes vers elle-même. C'est donc le cas de $Cat(1,m)$. Toutefois
le plan d'addition de cellules $Cat(1,m)$ part des \precats~quelconques et tombe dans les \precats~quelconques.
Or le fait d'associer à une \precat~quelconque la peinture des niveaux 1 et
$m$ tout entier est un procédé fonctoriel de la catégorie des
\precats~vers celle des \precats~peintes. De même, l'oubli de la peinture est
un foncteur des \precats~peintes vers les \precats. Ainsi $Cat(1,m)$ est bien un
foncteur des \precats~vers elles-mêmes.\\
Le second résultat résulte du premier auquel on applique le
lemme~\ref{plsfonc}.\\
CQFD.\\

Nous avons donc trouvé un plan d'addition fonctoriel de flèches génératrices de
\cats~faciles $Cat(1,m)$ qui préserve l'homotopie des \cats. Afin d'obtenir un
procédé de catégorisation préservant l'homotopie des \cats, il suffit de
construire un plan d'addition qui prendra en compte tous les $Cat(1,m)$. C'est
l'objet de la prochaine section.

\newpage

\section{Construction Bigcat}

Grâce aux sections précédentes, nous avons pu exhiber un plan d'addition
de flèches génératrices $Cat(1,m)$ qui préserve le type d'homotopie des
\cats. Toutefois ce plan n'ajoute que certaines flèches génératrices de
\cats~faciles et pas toutes. En revanche, si l'on considère l'ensemble des
$Cat(1,m)$ avec $m$ décrivant les entiers strictement supérieurs à un, on
obtient une famille de plans d'addition qui ensembles prennent toutes les
flèches génératrices. Nous allons donc définir la catégorisation
Bigcat comme un plan d'addition de cellules suivant un ordre et construit à partir des
$Cat(1,m)$, c'est-à-dire comme un plan d'addition de cellules de la forme
$BCat(t)$, où $t$ va parcourir tous les entiers strictement supérieurs à un et ce un nombre suffisamment grand de fois pour assurer que le plan Bigcat catégorise.

\begin{defin}\index{Bigcat}
Soit $(\mathcal{C},\mathcal{F}_1,\mathcal{F}_2)$ une donnée de Segal
pré-facile. Soit $\alpha$ un cardinal régulier strictement supérieur au
cardinal régulier pour lequel les sources et buts des flèches des familles $\mathcal{F}_1$ et $\mathcal{F}_2$
sont petites.\\
Notons $\Phi$ la famille $\mathcal{FG}_1$ des flèches
génératrices de \cats~faciles. Notons $\Phi_1$ la famille constituée des flèches génératrices de la forme
$\Delta[1]\Theta f$ avec $f\in \mathcal{F}_1$ et, pour tout entier $m\geq 2$, notons $\Phi_m$
la famille constituée des flèches génératrices de la forme
$\Delta[m]\Theta f$ et $Boit_m(g)$ avec $f\in \mathcal{F}_1$ et $g\in
\mathcal{F}_2$. On a bien une partition de $\Phi$ en $\Phi_m$ pour $m\geq 1$.\\
Posons $T$ l'ordinal transfini $\alpha\times\aleph_0$ et $t$ la fonction qui à
tout ordinal $(\beta,n)<\alpha\times\aleph_0$ associe le couple d'entier
$(1,n+2)$.\\
On définit alors le plan d'addition de cellules $Bigcat$ comme le
plan d'addition de cellules suivant l'ordre $t$ qui, à chaque étape
$\tau<T$, a pour plan d'addition $Cat(t(\tau))$.
\end{defin}
 
Après avoir défini $Bigcat$, nous allons montrer qu'il s'agit bien d'une
catégorisation. 
 
\begin{prop}\label{bigcat}
Soit $(\mathcal{C},\mathcal{F}_1,\mathcal{F}_2)$ une donnée de Segal
pré-facile.\\
$Bigcat$ est un foncteur de la catégorie des \precats~vers la catégories des
\cats~faciles vérifiant la propriété suivante :\\
pour tout morphisme $f$ d'une \precat~quelconque $A$ vers une \cat~facile $C$,
il existe un morphisme de \precats~non nécessairement unique de $Bigcat(A)$
vers $C$ dont la précomposée avec le morphisme naturel de $A$ vers
$Bigcat(A)$ est le morphisme $f$.
\end{prop}
{\it Preuve :}\\
Tout d'abord, il est facile de voir que $Bigcat$ est un foncteur car c'est une
application directe du lemme~\ref{cat1mf}. Pour avoir le reste de la
proposition, il suffit de montrer que les hypothèses du lemme~\ref{bcat} sont
vérifiées : celles sur $\Phi$ et celle sur $t$. Avant tout, on remarque que les plans du type $Cat(1,m)$ sont des plans d'addition de cellules composés de plans simples de type $Raj(1,m)$ qui sont des sous-plans de $e_{\Phi_1\cup\Phi_m,1}$. Or comme à la première étape de $Cat(1,m)$ on prend pour parties peintes la totalité des niveaux 1 et m, le sous-plan $Raj(1,m)$ à cette étape n'est autre que $e_{\Phi_1\cup\Phi_m,1}$ tout entier. Ainsi les plans $Cat(1,m)$ sont bien du type $Cat_{\Psi,\lambda}$ et donc Bigcat est bien du type $BCat(t)$, ce qui justifie le fait d'utiliser le lemme~\ref{bcat}.\\
\\
Premièrement, par hypothèse $(\mathcal{C},\mathcal{F}_1,\mathcal{F}_2)$ est
une donnée de Segal pré-facile. Donc les sources et buts des flèches des familles $\mathcal{F}_1$ et
$\mathcal{F}_2$ sont petites pour un certain cardinal régulier strictement
inférieur à $\alpha$. Par le lemme~\ref{flgenpet}, ceci entraîne que les sources des flèches de la
famille $\mathcal{FG}_1=\Phi$ sont petites par rapport à ce même cardinal. Comme en
outre $\Phi$ admet une partition, on a bien vérifié les hypothèses sur
$\Phi$.\\
\\
Comme $T$ est le produit $\alpha\times\aleph_0$ dont le cardinal est $\alpha$,
on a bien que le cardinal de $T$ est un cardinal régulier strictement supérieur au cardinal
régulier pour lequel les sources des flèches de la famille $\Phi$ sont petites. Par la suite, on mettra sur
$T$ l'ordre lexicographique.\\
\\ 
Soit $m\geq 1$ un entier. Montrons que pour tout ordinal $(\beta, n)<T$,
il existe un ordinal $\tau<T$ plus grand pour lequel $t(\tau)$ contient $m$. 
Deux cas se présentent. Si $m\geq n+2$, alors l'ordinal $(\beta,m-2)$
vérifie d'une part $(\beta,n)\leq (\beta,m-2)<T$ et d'autre part
$t(\beta,m-2)=(1,m)$. Sinon $m<n+2$. Comme
$\alpha$ est un cardinal transfini et que $\beta$ est strictement inférieur
à $\alpha$ alors $\beta+1$ est encore strictement inférieur à $\alpha$.
Ainsi $\tau=(\beta+1,m-2)$ est plus petit strictement que $T$ mais aussi plus
grand que $(\beta,n)$, car l'ordre est lexicographique. De plus ce $\tau$
vérifie bien que $t(\tau)=(1,m)$. Ceci montre donc la première condition sur $t$.\\
\\
Il est facile de voir avec notre $t$ que, pour $m\geq 2$, $t^{-1}(m)$ n'est
autre que l'ensemble des couples $(\beta,m-2)$, avec $\beta<\alpha$, dont le
cardinal est $\alpha=Card(T)$. Pour $m=1$, $t^{-1}(1)$ est en fait $T$ tout
entier, donc son cardinal est bien celui de $T$. Ceci montre la seconde
condition sur $t$.\\
\\
On peut donc appliquer le lemme~\ref{bcat} qui nous assure que le résultat de
$Bigcat$ est $\Phi$-injectif, c'est-à-dire est une \cat~facile, et que $Bigcat$
a bien la propriété de factorisation de la proposition.\\
CQFD.

\begin{cor}
Soit $(\mathcal{C},\mathcal{F}_1,\mathcal{F}_2)$ une donnée de Segal
pré-facile.
Soit $A$ une \precat. Comme, par la proposition précédente, $Bigcat(A)$ est une \cat~facile, on peut choisir de la marquer en prenant
pour relèvements marqués les premiers qui apparaissent dans l'ordre de la
construction. On obtient alors le résultat suivant :\\
tout morphisme de $A$ vers une \cat~marquée $C$ se factorise à travers
un morphisme préservant le marquage de $Bigcat(A)$ vers $C$, cette
factorisation n'étant pas nécessairement unique.
\end{cor}
{\it Preuve :} cela découle de la démonstration de la proposition. CQFD.\\

Notre nouvelle construction $Bigcat$ est une catégorisation composée de
plans $Cat(1,m)$ dont on a vu qu'ils préservaient l'homotopie des \cats. Il
serait donc normal que $Bigcat$ préserve cette propriété.

\begin{prop}\label{eqbig}
Soit $(\mathcal{C},\mathcal{F}_1,\mathcal{F}_2)$ une donnée de Segal pré-facile
vérifiant les propriétés suivantes :
\item 1) la catégorie sous-jacente $\mathcal{C}$ est une \cmf,
\item 2) les \eqcs~d'\obcs~sont exactement les \eqs~faibles de la
\cmf~$\mathcal{C}$ entre \obcs,
\item 3) la famille $\mathcal{F}_1$ est incluse dans la classe des cofibrations
triviales de la \cmf~$\mathcal{C}$ et la famille $\mathcal{F}_2$ dans celle des
cofibrations,
\item 4) les \eqs~faibles sont stables par produit fibré au-dessus d'un objet
discret.\\
\\
Alors pour toute \cat~$A$, le morphisme canonique $A\rightarrow Bigcat(A)$ est
une \eq~de \cats~niveau par niveau et une cofibration niveau par niveau. En
particulier $A\rightarrow Bigcat(A)$ est une \eq~de \cats.
\end{prop}
{\it Preuve :}\\
Tout d'abord on observe que toutes les constructions
apparaissant dans la construction Bigcat laissent invariant l'ensemble des
objets. Ainsi le morphisme canonique induit l'identité sur l'ensemble des
objets. Montrons par récurrence transfinie que pour tout ordinal $\tau\leq T$,
$C\rightarrow Bigcat_{\tau}(C)$ est niveau par niveau une cofibration
triviale et que les morphismes de Segal de $Bigcat_{\tau}(C)$ sont des
\eqs~faibles dans $\mathcal{C}$. Pour $\tau=0$, c'est vrai car $Bigcat_0(C)$ n'est autre que
$C$ qui est une \cat, donc ses morphismes de Segal sont des \eqcs~d'\obcs, ce
qui par l'hypothèse 2) montre que ses morphismes de Segal sont des
\eqs~faibles dans $\mathcal{C}$.\\ 

Supposons la propriété vraie pour $\tau<T$, alors les morphismes de Segal de
$Bigcat_{\tau}(C)$ sont des \eqs~faibles dans $\mathcal{C}$. Donc $Bigcat_{\tau}(C)$ vérifie les hypothèses du
lemme~\ref{eqcat1m} et donc $Bigcat_{\tau}(C)\rightarrow
Cat(t(\tau))(Bigcat_{\tau}(C))=Bigcat_{\tau+1}(C)$ est niveau par niveau une cofibration
triviale et les
morphismes de Segal de $Cat(t(\tau))(Bigcat_{\tau}(C))=Bigcat_{\tau+1}(C)$ sont des
\eqs~fai\-bles. Comme par hypothèse de récurrence sur $\tau<T$ le morphisme naturel
$C\rightarrow Bigcat_{\tau}(C)$ est niveau par niveau une cofibration triviale
et que dans la \cmf~$\mathcal{C}$ les cofibrations triviales sont stables par
composition, il vient que le morphisme $C\rightarrow Bigcat_{\tau+1}(C)$ est
niveau par niveau une cofibration triviale. On a donc montré l'hypothèse de
récurrence pour $\tau+1$.\\

Soit $\tau\leq T$ un ordinal limite et supposons l'hypothèse de récurrence
vraie pour tout $\tau'<\tau$. Comme $C\rightarrow
Bigcat_{\tau}(C)$ est niveau par niveau une colimite séquentielle transfinie de
cofibrations triviales $C\rightarrow Bigcat_{\tau'}(C)$ et que les cofibrations
triviales de la \cmf~$\mathcal{C}$ sont
stables par colimite séquentielle transfinie, il vient que niveau par niveau
$C\rightarrow Bigcat_{\tau}(C)$ est une cofibration triviale. Et par un raisonnement
identique à celui de la démonstration du lemme~\ref{eqcat1m}, on en déduit
que les morphismes de Segal de $Bigcat_{\tau}(C)$ sont des \eqs~faibles, ce qui
montre l'hypothèse de récurrence pour l'ordinal limite $\tau$.\\

Par récurrence transfinie, on a donc montré que pour tout ordinal $\tau\leq
T$ le morphisme naturel $C\rightarrow Bigcat_{\tau}(C)$ est une cofibration
triviale niveau par niveau et que les morphismes de Segal de $Bigcat_{\tau}(C)$
sont des \eqs~faibles de $\mathcal{C}$. C'est en particulier le cas pour
$\tau=T$, et comme $C$ et $Bigcat(C)$ sont des \cats, $C\rightarrow Bigcat(C)$
est niveau par niveau une \eqc~d'\obcs~par l'hypothèse 2), donc 
c'est bien une \eq~de \cats~niveau par niveau, et par le lemme~\ref{eqcniv},
c'est aussi une \eq~de \cats, ce qui montre la proposition.\\
CQFD.\\

La construction $Bigcat$ est donc non seulement un procédé de
catégorisation mais aussi une opération qui préserve le type d'homotopie
des \cats. Comme nous avons déjà vu que les \eqs~de \cats~vérifient la
propriété de "trois pour deux", nous sommes maintenant en mesure d'appliquer
la proposition~\ref{Ephi3}, ce qui est l'objet de la prochaine section.

\newpage

\section{Préservation du type d'homotopie par la construction Cat}

Afin de résoudre la problème de la catégorisation, nous cherchons à
montrer que la construction Cat adjointe du foncteur oubli des \cats~faciles
marquées vers les \precats~préserve le type d'homotopie des \precats. Or la
proposition~\ref{Ephi3} nous certifie que Cat a bien cette propriété si les
\eqs~de \cats~vérifient certaines propriétés et s'il existe une
catégorisation préservant le type d'homotopie des \cats. Dans les sections
précédentes, nous avons justement montré que, sous certaines hypothèses
sur la donnée de Segal pré-facile, les propriétés voulues pour les
\eqs~de \cats~sont vérifiées et que la construction Bigcat catégorise en
préservant le type d'homotopie des \cats. Rassemblons donc ces hypothèses
qui permettent d'appliquer la proposition~\ref{Ephi3}, ce qui donnera la
définition de donnée de Segal facile.

\begin{defin}\index{donnée de Segal!facile}
Une donnée de Segal pré-facile $(\mathcal{C},\mathcal{F}_1,\mathcal{F}_2)$ est dite facile si elle vérifie les
hypothèses suivantes :
\item 1) La catégorie sous-jacente $\mathcal{C}$ est une \cmf~dont tous les
objets sont cofibrants.
\item 2) Les \eqcs~d'\obcs~sont exactement les
équivalences faibles de la \cmf~$\mathcal{C}$ entre \obcs.
\item 3) Le foncteur $\tau_0$ est tel que, pour tout objet $C$ de
$\mathcal{C}$, $\tau_0(C)$ est un quotient de l'ensemble des morphismes dans
$\mathcal{C}$ de l'objet final vers $C$.
\item 4) Il existe un \obc~contractile $K$ muni de deux
morphismes de l'objet final vers $K$ notés 0 et 1, et ayant la propriété
suivante : pour tout couple $(f,g)$ de morphismes de l'objet final vers un \obc~$C$ tel que leurs images par $\tau_0$ soient égales, il existe un
morphisme de $K$ vers $C$ envoyant 0 sur $f$ et 1 sur $g$.
\item 5) Les \eqs~faibles sont stables par produit fibré au-dessus d'un objet
discret.
\item 6) Le produit fibré d'une équivalence faible entre
\obcs~le long d'une fibration est une équivalence faible.
\item 7) La famille $\mathcal{F}_1$ est incluse dans la classe des cofibrations
triviales de la \cmf~$\mathcal{C}$ et la famille $\mathcal{F}_2$ dans celle des
cofibrations.
\end{defin}

Appliquons maintenant la proposition~\ref{Ephi3} au triplet Cat, Bigcat et
\eq~de \cats~lorsque la donnée de Segal est facile.

\begin{theo}\label{boncat}
Soit $(\mathcal{C},\mathcal{F}_1,\mathcal{F}_2)$ une donnée de Segal facile.
Alors il existe un foncteur $Cat$ de la catégorie des \precats~vers la
catégorie des \cats~faciles marquées qui est adjoint à gauche du foncteur Oubli des
\cats~faciles marquées vers les \precats~et vérifiant que pour toute
\precat~$A$ l'image par $Cat$ du morphisme naturel $can_A:A\rightarrow Cat(A)$
est une \eq~de \cats.
\end{theo}
{\it Preuve :}\\
Les hypothèses de la donnée de Segal assurent, par le lemme~\ref{3pour2c} et
le corollaire~\ref{3pour2ccor}, que les \eqs~de
\cats~vérifient les deux premières hypothèses de la
proposition~\ref{Ephi3}. En outre le fait que les isomorphismes de \cats~sont
des \eqs~de \cats~est une conséquence directe du lemme~\ref{eqcniv}. Ainsi
toutes les hypothèses sur les \eqs~de \cats~sont vérifiées.\\
\\
La proposition~\ref{catfac} a montré que les \cats~faciles sont exactement les
objets $\mathcal{FG}_1$-injectifs et la construction Cat définie comme le
plan $E_{\mathcal{FG}_1}$ est bien un adjoint à gauche du foncteur Oubli, car la
donnée de Segal est facile, donc pré-facile, ce qui permet d'utiliser la
proposition~\ref{catadj}.\\ 
\\
Enfin la proposition~\ref{eqbig}, dont les
hypothèses sont vérifiées par la donnée de Segal facile, montre que le
plan fonctoriel d'addition de flèches de $\mathcal{FG}_1$ nommé Bigcat
vérifie les deux premières propriétés de la $\mathcal{FG}_1$-injectivisation mais
aussi la troisième pour les objets $\mathcal{FG}_1$-injectifs, c'est-à-dire pour les
\cats~faciles. Ceci termine de montrer que les hypothèses de la
proposition~\ref{Ephi3} sont vérifiées.\\ 
\\
On peut alors appliquer cette
proposition qui assure que Cat est une $\mathcal{FG}_1$-injectivisation, en
particulier vérifie la stabilité homotopique des \precats.\\
CQFD.\\

En fait l'application de la proposition~\ref{Ephi3} nous donne plus car elle
nous dit que les morphismes de \cats~faciles sont des \eqs~de
\cats~si et seulement si leurs images par Cat sont des \eqs~de \cats. Cette
remarque va nous permettre d'étendre la notion d'\eq~de \cats~aux \precats~en
disant qu'un morphisme de \precat~est une \eq~si son image par Cat est une
\eq~de \cats.

\begin{defin}\index{equivalence faible!de \precats}
Soit $(\mathcal{C},\mathcal{F}_1,\mathcal{F}_2)$ une donnée de Segal facile.
On dit d'un morphisme de \precats~que c'est une \eq~faible si son image par Cat
est une \eq~de \cats.
\end{defin}

Par cette définition, nous allons pouvoir rapatrier sur les \eqs~faibles les
propriétés des \eqs~de \cats, comme le "trois pour deux" par exemple.

\begin{lem}\label{3pour2}
Soit $(\mathcal{C},\mathcal{F}_1,\mathcal{F}_2)$ une donnée de Segal facile.\\
Les \eqs~faibles de \precats~vérifient les propriétés suivantes :
\item - les isomorphismes de \precats~sont des \eqs~faibles,
\item - les morphismes de \cats~faciles sont des \eqs~faibles si et seulement si
ce sont des \eqs~de \cats,
\item - soient $f$ et $g$ deux morphismes de \precats~composables, si parmi $f,g,g\circ f$
deux morphismes sont des \eqs~faibles alors le troisième aussi,
\item - soient $f$ et $g$ deux morphismes de \precats~composables, si $f\circ g$
est l'identité et que $g\circ f$ est une \eq~faible, alors $f$ et $g$ sont des
\eqs~faibles.
\end{lem}
{\it Preuve :}\\
Comme pour la proposition~\ref{boncat}, on démontre que les hypothèses du
corollaire~\ref{Ephi3'}, qui sont les mêmes que celles de la
proposition~\ref{Ephi3}, sont vérifiées, ce lemme étant une conséquence
directe du corollaire~\ref{Ephi3'} car ici les $\Phi$-\eqs~ne sont autres que les
\eqs~faibles.\\
CQFD.\\

En fait, si l'on regarde de plus près, Bigcat vérifie un résultat plus
fort que celui demandé par la proposition~\ref{Ephi3}. En effet, la
proposition demande seulement la préservation de l'homotopie pour les
\cats~faciles alors que Bigcat préserve l'homotopie pour toutes les \cats.
Aussi si l'on refait la démonstration de la proposition~\ref{Ephi3} pour les
\cats~quelconques, on obtient les résultats suivants.

\begin{lem}\label{eqceqf}
Soit $(\mathcal{C},\mathcal{F}_1,\mathcal{F}_2)$ une donnée de Segal facile.
\item - Pour toute \cat~quelconque $A$, le morphisme naturel $A\rightarrow Cat(A)$ est
une \eq~de \cats.
\item - Pour tout morphisme de \cats~quelconque $f:A\rightarrow B$, on a :\\
$f$ est une \eq~de \cats~si et seulement si $f$ est une \eq~faible.
\end{lem}
{\it Preuve :} c'est la même démonstration que la proposition~\ref{Ephi3} en
utilisant le fait que Bigcat préserve l'homotopie de toutes les \cats~et pas
seulement des \cats~faciles.\\

Il serait également intéressant de comparer la notion d'\eq~faible
définie par Cat avec la notion d'\eq~faible qu'on aurait pu définir avec
Bigcat.

\begin{lem}\label{bigcateq}
Soit $(\mathcal{C},\mathcal{F}_1,\mathcal{F}_2)$ une donnée de Segal facile.
\item - Pour toute \precat~$A$, le morphisme naturel $A\rightarrow Bigcat(A)$ est
une \eq~faible de \precats, c'est-à-dire son image par Cat est une \eq~de
\cats.
\item - Pour tout morphisme de \precats~$f:A\rightarrow B$, on a :\\
$f$ est une \eq~faible de \precats~si et seulement si $Bigcat(f)$ est une \eq~de
\cats.
Ce qui signifie :\\
$Cat(f)$ est une \eq~de \cats~si et seulement si $Bigcat(f)$ est une \eq~de
\cats.
\end{lem}
{\it Preuve :}\\
La première partie du lemme est l'application directe du lemme~\ref{bcateq}
dont les hypothèses nommées hypothèses~\ref{colimIinj} ne sont autres que
celles de la proposition~\ref{Ephi3} dont on a montré qu'elles sont
vérifiées par la donnée de Segal facile.\\
\\
Pour la seconde partie, considérons le diagramme commutatif suivant :
\begin{diagram}
A & \rTo^f & B \\
\dTo & & \dTo \\
Bigcat(A) & \rTo_{Bigcat(f)} & Bigcat(B)\\
\end{diagram} 
Par la première partie du lemme, les flèches verticales sont des
\eqs~fai\-bles de \precats. En appliquant le lemme~\ref{3pour2}, il vient que $f$
est une \eq~faible de \precats~si et seulement si $Bigcat(f)$ est une \eq~faible
de \precats. Or $Bigcat(f)$ est un morphisme de \cats~faciles, donc toujours par
le lemme~\ref{3pour2}, $Bigcat(f)$ est une \eq~faible de \precats~si et
seulement si c'est une \eq~de \cats, ce qui montre la seconde partie du lemme.\\
CQFD.\\

Nous avons finalement résolu notre problème de la catégorisation, à
savoir trouver une procédé qui à une \precat~associe une \cat~qui soit de
même type homotopique qu'elle. Pour cela, on a facilité la donnée de Segal
en lui adjoignant deux familles engendrant certains \obcs~et certaines
\eqcs~d'\obcs. Bien sûr, comme on l'a fait remarquer à l'époque ce choix
est arbitraire et rien ne nous dit que la notion de \cat~facile n'est pas
trop forte. C'est pourquoi pour éviter ce phénomène on a requis pour les
catégorisations de préserver le type d'homotopie.\\ 

On a par la suite trouvé
sous quelles hypothèses on a une telle catégorisation, c'est ce que l'on a
appelé donnée de Segal facile, car elle facilite l'obtention d'une bonne
catégorisation. Les hypothèses de cette donnée de Segal facile montre la
marge que l'on a dans le choix des familles génératrices d'\obcs~faciles et
d'\eqcs~faciles d'\obcs~faciles qui permettront une bonne catégorisation.\\
 
On remarque aussi
que les hypothèses des données de Segal faciles sont assez contraignantes
car elles demandent à la catégorie $\mathcal{C}$ de posséder une structure de
\cmf, ce qui est tout à fait normal puisque nous voulons manipuler des sortes
d'\eqs~faibles que sont les \eqcs~d'\obcs~et nous assurer de leur stabilité
par colimites. Par ailleurs, toutes les constructions catégorisantes sont en
fait des I-injectivisations qui sont des processus caractéristiques de la
structure de \cmf.\\
 
Aussi allons nous dans le prochain chapitre nous intéresser à
montrer que la catégorie $\mathcal{C-PC}$ avec les monomorphismes et les
\eqs~faibles de \precats~forme une \cmf. Ceci nous donnera un bon cadre
théorique aux manipulations homotopiques sur les \precats.

\chapter{Structure de \cmfcof~sur les \precats}

\newpage

Après avoir présenté la notion de catégorie enrichie faible sur une
catégorie $\mathcal{C}$ dans un premier chapitre, nous nous sommes posés le
problème de la caté\-gorisation : trouver une construction fonctorielle qui
à toute \precat~associe une \cat~en préservant le type d'homotopie. Les
trois chapitres précédents ont été entièrement consacrés à la
résolution de ce problème. A la fin du chapitre précédent, nous avons pu
établir à quelles conditions sur $\mathcal{C}$ ce problème possède une
solution. Parmi ces conditions se trouve la structure de \cmf~et ce n'est pas un
hasard. En effet la structure de \cmf~est un bon cadre pour traiter de
l'homotopie, des remplacements fibrants et autres. Et c'est ce que l'on faisait
en cherchant à résoudre le problème de la catégorisation, puisque l'on
manipulait les notions d'\eqcs~et d'\obcs~faciles, entre autres.\\

C'est la raison pour laquelle il paraît tout à fait indiquer de chercher
à mettre sur la catégorie $\mathcal{C-PC}$ des \precats~une telle structure
de \cmf. L'existence de cette structure est en fait une propriété de la catégorie et de
trois classes de ses morphismes que l'on nomme cofibrations, équivalences
faibles et fibrations. Comme nous l'avons vu au chapitre précédent, la
notion d'\eq~de \cats~s'étend aux \precats~par l'intermédiaire de la
catégorisation Cat et nous dirons d'un morphisme de \precats~que c'est une
\eq~faible si son image par Cat est une \eq~de \cats. Afin de simplifier les
démonstrations, nous allons prendre pour cofibrations les monomorphismes de
\precats, c'est-à-dire les morphismes de \precats~qui niveau par niveau sont
des monomorphismes de $\mathcal{C}$. Enfin comme
c'est habituellement le cas dans ce genre de problème, nous définirons les
fibrations de \precats~comme les morphismes ayant la \prd~par rapport aux
monomorphismes qui sont des \eqs~faibles.\\

Le but de ce chapitre est donc de montrer que la catégorie $\mathcal{C-PC}$
des \precats~avec les cofibrations, \eqs~faibles et fibrations ainsi définies
est bien une \cmf. Pour des raisons de commodité, nous démontrerons que
c'est en fait une \cmfcof. Afin d'avoir en tête les axiomes de ces structures,
nous débuterons ce chapitre par un rappel des définitions de ces structures
ainsi que par l'énoncé du lemme de reconnaissance de la structure de
\cmfcof. Puis nous allons montrer au fur et à mesure de ce chapitre que la
catégorie $\mathcal{C-PC}$ vérifie les différentes hypothèses de ce
lemme. Ce qui nous amènera au théorème central de la thèse qui énonce
à quelles conditions sur la catégorie $\mathcal{C}$ la catégorie des
\precats~$\mathcal{C-PC}$ est une \cmfcof~avec les trois classes définies
ci-dessus.

\newpage

\setcounter{section}{-1}
\section{Rappels sur la structure de \cmfcof}

Dans cette partie, nous allons surtout faire des rappels concernant la
structure de \cmf~due à Quillen dans~\cite{q} et celle de \cmfcof. En particulier nous allons énoncer le lemme de
reconnaissance de la structure de \cmfcof~de \cite{h}.\\

Commençons tout d'abord par rappeler la définition de \cmf.

\begin{defin}
Une \cmf~est une catégorie $\mathcal{M}$ munie de trois classes de morphismes,
nommés respectivement les cofibrations, les équivalences faibles et les
fibrations, et satisfaisant les axiomes suivants :
\item CM1 : $\mathcal{M}$ est complète et cocomplète.
\item CM2 (axiome "trois pour deux") : si $f$ et $g$ sont deux morphismes
composables dans $\mathcal{M}$ et si deux morphismes parmi $f,\; g \mbox{ et }
g\circ f$ sont des équivalences faibles alors le morphisme restant est aussi
une équivalence faible.
\item CM3 : Les cofibrations, les équivalences faibles et les fibrations sont
stables par rétracts.
\item CM4 : Les cofibrations ont la \prg~par rapport aux fibrations triviales
(morphismes à la fois fibrations et équivalences faibles). Les fibrations
ont la \prd~par rapport aux cofibrations triviales
(morphismes à la fois cofibrations et équivalences faibles).
\item CM5 : Tout morphisme de $\mathcal{M}$ possède une factorisation
(fonctorielle) en une cofibration triviale suivie d'une fibration ainsi qu'une factorisation
(fonctorielle) en une cofibration suivie d'une fibration triviale.
\end{defin}
\index{catégorie de modèles fermée}

Parmi les catégories ayant une structure de \cmf, certaines ont la
propriété que leurs cofibrations et leurs cofibrations triviales sont
engendrées par des ensembles de morphismes, au sens du deuxième chapitre. On
dira alors que ces catégories ont une structure de \cmfcof.

\begin{defin}\index{catégorie de modèles fermée!engendrée de manière cofibrante}
Une \cmfcof~est une \cmf~$\mathcal{M}$ pour laquelle :
\item - il existe un ensemble de morphismes $\mathcal{I}$ (appelé ensemble des
cofibrations généra\-trices) permettant l'argument du petit objet et tel qu'un
morphisme est une fibration triviale si et seulement s'il a la \prd~par rapport
à $\mathcal{I}$.
\item - il existe un ensemble de morphismes $\mathcal{J}$ (appelé ensemble des
cofibrations triviales génératrices) permettant l'argument du petit objet et tel qu'un
morphisme est une fibration si et seulement s'il a la \prd~par rapport
à $\mathcal{J}$.
\end{defin}

Il est assez courant parmi les catégories ayant la structure de \cmfcof~que
leurs cofibrations soient en fait les monomorphismes. On dira, dans ce cas, que
leur \cmf~est engendrée par monomorphismes.

\begin{defin}\index{catégorie de modèles fermée!engendrée par
monomorphismes}
Une \cmfmono~est une \cmfcof\\ pour laquelle les cofibrations sont exactement les
monomorphismes.
\end{defin}

Comme nous l'avons annoncé plus haut, nous allons utiliser un lemme de
reconnaissance de la structure de \cmfcof, que l'on trouve dans le papier
d'Hirschhorn~\cite{h}, afin de montrer que la catégorie des \precats, munie des
\eqs~faibles définies précédemment, est une \cmf~engendrée par
monomorphismes.

\begin{lem}[de reconnaissance]\label{reco}
Soit $\mathcal{M}$ une catégorie munie d'une clas\-se de morphismes W (ses
\eqs~faibles) et de deux
ensembles de morphismes I et J. Supposons que $\mathcal{M}$ vérifie les
propriétés suivantes :
\item 0) $\mathcal{M}$ est complète et cocomplète,
\item 1) W est stable par rétracts et satisfait à l'axiome "trois pour
deux" : pour tout couple $(f,g)$ de morphismes
composables, si deux morphismes parmi $f,\; g$ et $g\circ f$ sont dans W alors le morphisme
restant est aussi dans W,
\item 2) I et J permettent l'argument du petit objet,
\item 3) Toute J-cofibration est à la fois une I-cofibration et dans W,
\item 4) Tout morphisme I-injectif est à la fois J-injectif et dans W,
\item 5) L'une des deux propriétés suivantes est vraie : 
\item - les morphismes à la fois I-cofibrations et dans W sont
des J-cofibrations,
\item - les morphismes à la fois J-injectifs et dans W sont
I-injectifs.\\
\\
Alors il existe une structure de \cmfcof~sur $\mathcal{M}$ avec W pour classe
des équivalences faibles, I pour ensemble des cofibrations génératrices
et J pour ensemble des cofibrations triviales génératrices.
\end{lem}

Munis de ce lemme de reconnaissance, nous allons maintenant pouvoir montrer que
la catégorie $\mathcal{C-PC}$ avec ses \eqs~faibles de \precats~est une
\cmf~engendrée par monomorphismes en montrant qu'elle vérifie les six
hypothèses du lemme de reconnaissance~\ref{reco}.

\newpage

\section{Stabilité des équivalences faibles par rétracts}

Avant de commencer la vérification des hypothèses du lemme de
reconnaissance~\ref{reco} de la structure de \cmfcof, récapitulons les
définitions des cofibrations, \eqs~fai\-bles et fibrations des \precats.

\begin{defin}[rappel des \eqs~faibles de \precats]
Soit $(\mathcal{C},\mathcal{F}_1,\mathcal{F}_2)$ une donnée de Segal facile.
On dit d'un morphisme de \precats~que c'est une \eq~faible si son image par Cat
est une \eq~de \cats.
\end{defin}

\begin{defin}\index{cofibration de \precats}
\index{cofibration triviale de \precats}
Soit $(\mathcal{C},\mathcal{F}_1,\mathcal{F}_2)$ une donnée de Segal facile.
On appelle cofibrations de \precats~les monomorphismes de \precats,
c'est-à-dire les morphismes de \precats~qui niveau par niveau sont des
monomorphismes de $\mathcal{C}$.
On appelle cofibrations triviales de \precats~les monomorphismes de \precats~qui
sont également des \eqs~faibles de \precats.
\end{defin}

\begin{defin}\index{fibration de \precats}
\index{fibration triviale de \precats}
Soit $(\mathcal{C},\mathcal{F}_1,\mathcal{F}_2)$ une donnée de Segal facile.
On appelle fibrations de \precats~les morphismes de \precats~ayant la \prd~par
rapport aux cofibrations triviales de \precats.
On appelle fibrations triviales de \precats~les fibrations de \precats~qui
sont également des \eqs~faibles de \precats.
\end{defin}

Après avoir mis au clair les définitions des classes de cofibrations,
\eqs~faibles et fibrations pour les \precats, nous sommes à même de
débuter la vérification des hypothèses du lemme de
reconnaissance~\ref{reco} de la structure de \cmfcof. Pour cela, nous allons
toujours nous mettre dans le cas où $\mathcal{C}$ est munie d'une donnée de
Segal facile.\\ 

L'hypothèse 0) consiste
en la complétude et la cocomplétude de $\mathcal{C-PC}$. Or comme on ne considère que les données de Segal faciles, la
catégorie $\mathcal{C}$ sous-jacente est toujours une
\cmf, donc elle est complète et cocomplète, ce qui par la
proposition~\ref{precatclos}
assure que $\mathcal{C-PC}$ est complète et cocomplète.\\

L'hypothèse 1) demande aux \eqs~faibles de \precats~d'être stables par
rétracts et par "trois pour deux". Comme le lemme~\ref{3pour2} nous assure la
propriété de "trois pour deux" des \eqs~faibles de \precats, il ne reste
plus qu'à montrer leur stabilité par rétract.

\begin{lem}\label{eqfre}
Soit $(\mathcal{C},\mathcal{F}_1,\mathcal{F}_2)$ une donnée de Segal
pré-facile vérifiant la propriété suivante :\\
les \eqcs~d'\obcs~sont stables par rétracts.\\
Alors les équivalences faibles de $\mathcal{C-PC}$ sont stables par
rétracts.
\end{lem}
{\it Preuve :}\\
Par fonctorialité de Cat, montrer que les équivalences faibles sont stables
par rétracts revient à montrer que les équivalences de \cats~entre
\cats~marquées sont stables par rétracts.\\ 
\\
Soient donc $f:A\rightarrow B$ et
$g:C\rightarrow D$ deux morphismes de \cats~mar\-quées tels que $f$ soit un
rétract de $g$ et $g$ soit une équivalence de \cats. Le fait que $f$ est rétract
de $g$ et que $\tau_0$ est un foncteur entraîne que $\tau_0(f)$ est un
rétract de $\tau_0(g)$. Or on a vu que les équivalences de \cats~induisent
des bijections par $\tau_0$. Donc $\tau_0(g)$ est une bijection et comme les
bijections ensemblistes sont stables par rétracts, il vient que $\tau_o(f)$
est bijectif et donc que $f$ est essentiellement surjectif.\\
\\
Notons maintenant $i,\; j,\; r,\mbox{ et } s$ les morphismes intervenant dans la
définition de rétracts. On a donc les égalités suivantes:
$$g\circ i=j\circ f,\; f\circ r=s\circ g,\; r\circ i=Id_A \mbox{ et } s\circ
j=Id_B.$$
On remarque que pour tout couple $(a,a')$ d'objets de $A$, $f_1(a,a'):A_1(a,a')\rightarrow B_1(f(a),f(a'))$ est un
rétract de $g_1(i(a),i(a')):C_1(i(a),i(a'))\rightarrow D_1(g\circ i(a),g\circ
i(a'))$ car on a les égalités suivantes :
$$g_1(i(a),i(a'))\circ i_1(a,a')=j_1(f(a),f(a'))\circ f_1(a,a')$$
$$ f_1(a,a')\circ
r_1(i(a),i(a'))=s_1(j\circ f(a),j\circ
f(a'))\circ g_1(i(a),i(a'))$$
$$ r_1(i(a),i(a'))\circ i_1(a,a')=Id_{A_1(a,a')}$$
$$ s_1(j\circ f(a),j\circ
f(a'))\circ j_1(f(a),f(a'))=Id_{B_1(f(a),f(a'))}. $$
Or par hypothèse, les \eqcs~d'\obcs~sont stables par rétracts. Comme $g$ est une équivalence de \cats,
on a que $g_1(i(a),i(a'))$ est une \eqc~d'\obcs~et, par conséquent,
$f_1(a,a')$ aussi comme rétract du précédent. Ceci montre que $f$ est
pleinement fidèle et comme on a déjà montré que $f$ est essentiellement
surjective, il vient que $f$ est une équivalence de \cats.\\
CQFD.\\
\\

Bien évidemment si notre donnée de Segal pré-facile est facile, alors les
\eqcs~d'\obcs~sont des \eqs~faibles de la \cmf~$\mathcal{C}$ et sont donc
stables par rétracts, ainsi toute donnée de Segal facile vérifie les
hypothèses de ce lemme, ce qui finit de montrer l'hypothèse 1) du
lemme~\ref{reco}.\\

A partir de l'hypothèse 2), interviennent les ensembles générateurs
des cofibrations et des cofibrations triviales,
notés $\mathcal{I}$ et $\mathcal{J}$. Nous allons donc les définir dans la
section suivante.

\newpage

\section{Ensembles générateurs $\mathcal{I}$ et $\mathcal{J}$}

La structure de \cmfcof~est basée sur le fait que les cofibrations et les
cofibrations triviales sont engendrées par certains ensembles notés
$\mathcal{I}$ et $\mathcal{J}$. Comme nous avons déjà vérifié les deux
premières hypothèses du lemme de reconnaissance~\ref{reco}, il nous reste à
vérifier les quatre dernières qui traitent justement de ces ensembles
générateurs, que nous devons donc définir. Bien que l'on ait pris comme
cofibrations de \precats~les monomorphismes niveau par niveau, nous allons tant
que faire se peut nous placer dans le cadre plus général où les
cofibrations des \precats~sont les cofibrations niveau par niveau.\\

Pour trouver l'ensemble $\mathcal{I}$ des cofibrations génératrices, on remarque tout
d'abord que par définition les cofibrations sont niveau par niveau des
cofibrations de $\mathcal{C}$. Si $\mathcal{C}$ est une
\cmfcof, alors les cofibrations de $\mathcal{C-PC}$ sont niveau par niveau des colimites de sommes amalgamées par
les cofibrations génératrices de $\mathcal{C}$.\\
Soit $f:A\rightarrow B$ une cofibration, tout d'abord c'est une injection au
niveau des objets et donc, pour rajouter à $A$ les éventuels objets que $B$ a en
plus, on prend la somme amalgamée de $A$ par $\emptyset\rightarrow *$.\\
Supposons maintenant que jusqu'au rang $n-1$, $f$ induise des isomorphismes. On
veut alors rajouter à $A_n$ ce qui lui manque pour devenir $B_n$,
c'est-à-dire faire une colimite de somme amalgamée par les cofibrations
génératrices de $\mathcal{C}$. On remarque alors que quand on veut rajouter
une cofibration génératrice $g:X\rightarrow Y$ de $\mathcal{C}$ à $A_n$,
l'image de $Y$ dans le bord de $B_n$ est déjà dans le bord de $A_n$ car $f$
induit des isomorphismes jusqu'au rang $n-1$. Ainsi si l'on suppose $Y$ connexe,
on obtient le diagramme suivant commutatif:
\begin{diagram}
\Delta[n]\Theta X\coprod_{\partial\Delta[n]\Theta X}\partial\Delta[n]\Theta
Y & \rTo & A\\
\dTo^{Attach_n(g)} & & \dTo_{f}\\
\Delta[n]\Theta Y & \rTo & B\\
\end{diagram}
Ces considérations nous conduisent à proposer la définition suivante pour
l'ensemble générateurs des cofibrations $\mathcal{I}$.

\begin{lem}[-définition]\label{I}
Soit $\mathcal{C}$ une catégorie \discret~munie d'une structure de \cmfcof~dont les cofibrations génératrices ont leurs
sources et buts $\alpha$-petits, pour un cardinal régulier $\alpha$ plus grand
qu' $\aleph_0$. Notons
$\mathcal{I}$ l'ensemble constitué de $\emptyset\rightarrow *$ et des
flèches de type $Attach_n(g):\Delta[n]\Theta X\coprod_{\partial\Delta[n]\Theta X}\partial\Delta[n]\Theta
Y\rightarrow \Delta[n]\Theta Y $ pour $n>0$ et $g:X\rightarrow Y$ décrivant
l'ensemble des cofibrations génératrices de
$\mathcal{C}$.\index{$\mathcal{I}$}\index{$Attach_n(g)$} Alors $\mathcal{I}$ permet l'argument du petit
objet au sens~\ref{ptobj}.
\end{lem}
{\bf Preuve :}\\
On remarque que, niveau par niveau, la source d'une flèche de type $Attach_n(g)$ est la somme
disjointe d'un ensemble fini avec un coproduit fini de $X$ et un coproduit fini de
$Y$ et le but est la somme disjointe d'un ensemble fini avec un coproduit fini
de $Y$. Par hypothèse, $X$ et $Y$ sont $\alpha$-petits
alors, par le lemme~\ref{cpdtpetit}, chaque
niveau des extrémités de $Attach_n(g)$ aussi. Et par
lemme~\ref{precatpetit}, car $\alpha$ est régulier et strictement supérieur
à $\aleph_0$, il vient que source et but de $Attach_n(g)$ sont
$\alpha$-petits. On a déjà vu que $\emptyset$ et $*$ sont aussi
$\alpha$-petits. Donc par le lemme~\ref{enspt}, on a bien que $\mathcal{I}$ permet
l'argument du petit objet.\\
CQFD.\\

Pour l'ensemble générateur des cofibrations triviales $\mathcal{J}$, nous
allons suivre la démarche habituelle qui consiste à ne garder que les
classes d'\eq~de cofibrations triviales petites par rapport à un cardinal
fixé. L'avantage de cette démarche est que par définition même
l'ensemble $\mathcal{J}$ vérifie l'argument du petit objet.

\begin{lem}[-définition]\label{J}
Sous les hypothèses et notations du lemme pré\-cédent, choisissons un
cardinal $\alpha'$ régulier et strictement supérieur à $2^{\alpha}$. Notons $\mathcal{J}$
un ensemble de représentants de classe d'isomorphismes 
de morphismes à la fois cofibrations niveau par niveau et équivalences faibles dont les sources et buts sont
$\widetilde{\alpha'}$-petits, au sens de la proposition~\ref{alphatilde}. Alors $\mathcal{J}$ permet l'argument du petit objet.\index{$\mathcal{J}$}
\end{lem}
{\bf Preuve :} application directe du lemme~\ref{enspt}.\\

A part pour la partie concernant l'engendrement des cofibrations triviales par
$\mathcal{J}$, il n'est pas nécessaire de savoir quel cardinal limite les
cofibrations triviales de $\mathcal{J}$ mais il suffit de savoir que les
cofibrations triviales de $\mathcal{J}$ sont limitées, c'est-à-dire que
leurs sources et buts sont $\beta$-petits pour un cardinal $\beta$ régulier et
strictement supérieur à $2^{\aleph_0}$, afin que de savoir que $\mathcal{J}$
permet l'argument du petit objet.\\

Dans cette partie, nous avons non seulement défini les ensembles
générateurs des cofibrations et des cofibrations triviales mais aussi
montrer qu'ils permettent l'argument du petit objet sous des hypothèses plus
contraignante que celles de donnée de Segal facile car il est ici demandé à
la catégorie $\mathcal{C}$ d'être une \cmfcof. Néanmoins nous avons
réglé le cas de l'hypothèse 2) du lemme de reconnaissance~\ref{reco}. Il ne
reste donc plus que les hypothèses 3), 4) et 5). Or l'hypothèse 3) est assez
difficile à démontrer et l'hypothèse 5) découle en partie de la 3).
C'est pourquoi nous allons démontrer d'abord l'hypothèse 4), ce qui sera
l'objet de la prochaine section.

\newpage

\section{$\mathcal{I}$ engendre les cofibrations de \precats}

L'hypothèse 4) du lemme de reconnaissance~\ref{reco} de la structure de
\cmfcof~est en fait constituée de deux parties. Il faut en effet d'une part
montrer que les morphismes de \precats~$\mathcal{I}$-injectifs sont aussi
$\mathcal{J}$-injectifs et d'autre part que ce sont des \eqs~faibles. Pour
montrer le premier point, il suffit de remarquer que l'ensemble $\mathcal{J}$
est un sous-ensemble de la classe des cofibrations. Donc si l'on montre que
$\mathcal{I}$ engendre bien les cofibrations, c'est-à-dire que les
cofibrations sont exactement les $\mathcal{I}$-cofibrations, alors tout
morphisme $\mathcal{I}$-injectif se relèvera par rapport aux cofibrations, donc
en particulier par rapport à $\mathcal{J}$, ce qui en fera un morphisme
$\mathcal{J}$-injectif.\\

Comme dans la section précédente, nous allons faire ces démonstrations le
plus possible dans le cas où les cofibrations des \precats~sont les
cofibrations niveau par niveau et pas uniquement les monomorphismes. C'est
pourquoi nous commencerons par rappeler des résultats de stabilité des
cofibrations.

\subsection{Stabilité des cofibrations dans une \cmf}

Pour la suite, il est important d'avoir à l'esprit quels résultats de stabilité les
cofibrations vérifient dans une \cmf, c'est ce qu'énonce le lemme suivant.

\begin{lem}\label{cofst}
Soit $\mathcal{C}$ une \cmf, alors les cofibrations de $\mathcal{C}$ forment une sous-catégorie contenant les isomorphismes et
stable par rétract, somme amalgamée le long d'un morphisme, colimite
séquentielle transfinie et coproduit dans la catégorie des morphismes.
\end{lem}
{\it Preuve :}\\
Comme $\mathcal{C}$ est une \cmf, les cofibrations sont exactement les morphismes qui
ont la \prg~par rapport aux fibrations triviales et par le lemme~\ref{prst} on a
le résultat. CQFD.\\

Bien entendu cette stabilité des cofibrations dans la \cmf~$\mathcal{C}$ est
encore valable pour les cofibrations niveau par niveau de $\mathcal{C-PC}$.

\begin{cor}\label{cofnpnst}
Soit $\mathcal{C}$ une catégorie \discret~munie d'une structure de \cmf, alors les cofibrations niveau par niveau de
$\mathcal{C-PC}$ vérifient les propriétés du lemme précédent.
\end{cor}
{\it Preuve :}\\
Dans $\mathcal{C-PC}$, les rétracts et les colimites sont calculées niveau
par niveau. En outre les cofibrations de $\mathcal{C-PC}$ sont niveau par niveau
des cofibrations de $\mathcal{C}$. D'où le résultat par le lemme
précédent. CQFD.\\

Un autre résultat particulièrement important et utile est celui de la
préserva\-tion des cofibrations par somme amalgamée dans la catégorie des
morphismes. C'est ce qu'on appelle le lemme du cube pour les cofibrations. Nous
allons le rappeler pour une \cmf~$\mathcal{C}$ puis pour les cofibrations niveau
par niveau de $\mathcal{C-PC}$.

\begin{lem}\label{cofcub}
Soit $\mathcal{C}$ une \cmf, considérons le diagramme suivant :
\begin{diagram}
A  & & \rTo & & B & & \\
   & \rdTo & & & \dLine & \rdTo & \\
\dTo & & C & & \rTo &  & D \\
   & &  & & \dTo & & \\
A' & \rLine & \dTo & \rTo & B' & & \dTo \\
   & \rdTo & & & & \rdTo & \\
   & & C' & & \rTo &  & D' \\ 
\end{diagram}
o\'u $D$ et $D'$ sont les sommes amalgamées $B\coprod_A C$ et $B'\coprod_A'
C'$.\\ Si $C\rightarrow C'$ et $A'\coprod_A B\rightarrow B'$ sont des
cofibrations, alors $D\rightarrow D'$ est aussi une cofibration.
\end{lem}
{\it Preuve :} dans \cite{h}.

\begin{cor}\label{cofnpncub}
Soit $\mathcal{C}$ une catégorie \discret~munie d'une structure de \cmf, alors les cofibrations niveau par niveau de
$\mathcal{C-PC}$ vérifient les propriétés du lemme précédent.
\end{cor}
{\it Preuve :} niveau par niveau on applique le lemme.

\subsection{Astuce ensembliste}

Après ces quelques rappels concernant la stabilité des cofibrations, nous
allons démontrer deux lemmes techniques généraux qui servent lorsque l'on
veut montrer qu'une famille de morphismes est bien engendrée par un ensemble
donné. 

\begin{lem}\label{ast1}
Soient $\mathcal{C}$ une catégorie et $\mathcal{C}'$ une sous-catégorie de
$\mathcal{C}$ ayant les mêmes objets, contenant les isomorphismes et telle que
les isomorphismes ne se factorisent qu'en isomorphismes dans $\mathcal{C}'$.
Supposons que pour tout morphisme $j:A\rightarrow B$ de
$\mathcal{C}'$, la classe contenant un représentant de chaque classe d'isomorphismes d'objets $C$ par lequel $j$ se
factorise dans $\mathcal{C}'$ est un ensemble. On a alors les résultats
suivants :
\item 1) La catégorie des préfaisceaux sur $\mathcal{C}$ avec la
sous-catégorie des morphismes qui niveau par niveau sont des morphismes de
$\mathcal{C}'$ possède également cette propriété.
\item 2) Etant donné le carré commutatif suivant, noté $(j,a,b,p)$ :
\begin{diagram}
A & \rTo^a & X\\
\dTo^{j} & & \dTo_{p}\\
B & \rTo_{b} & Y \\
\end{diagram}
avec $j$ morphisme de $\mathcal{C}'$. Est aussi un ensemble 
la classe contenant un représen\-tant de chaque classe d'isomorphismes de
couples $(C,c)$ avec $C$ objet factorisant $j$ dans
$\mathcal{C}'$ et tels que la flèche pointillée $c$ du diagramme ci-dessous
existe et fait commuter le diagramme :
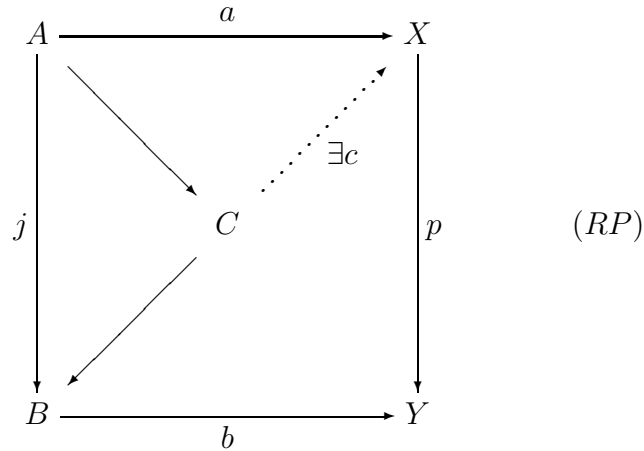
\begin{diagram}
A & & \rTo^a & & X & & \\
  & \rdTo & & \ruDotsto_{\exists c} &  & & \\
\dTo^{j} & & C & & \dTo_{p} & & (RP)\\
  & \ldTo & & & & & \\
B & & \rTo_{b} & & Y & & \\   
\end{diagram}
On dira que $C$ réalise avec $c$ un relèvement partiel du diagramme $(j,a,b,p)$.
\end{lem}
{\it Preuve :}\\
Pour la première partie du lemme, il suffit de voir que la classe de représentants des classes
d'isomorphismes de préfaisceaux factorisant un morphisme de pré\-faisceaux
fixé dans la bonne sous-catégorie est une sous-classe de la réunion
dénombrable des classes de représentants des classes
d'isomorphismes d'objets de $\mathcal{C}$ factorisant dans $\mathcal{C}'$ un des
niveaux du morphisme de préfaisceaux fixé. On obtient le résultat 1) en utilisant
le fait qu'une réunion dénombrable
d'ensembles est un ensemble et qu'une sous-classe d'un ensemble est un
ensemble.\\
\\
Pour un objet $C$ factorisant $j$ dans $\mathcal{C'}$, la classe des morphismes
$c$ de $C$ vers $X$ est un ensemble, par définition de catégorie. En outre
par hypothèse, la classe des représentants des classes d'isomorphismes d'objets $C$ factorisant
$j$ dans $\mathcal{C'}$ est un ensemble. Donc l'union disjointe sur cette
ensemble des ensembles de morphismes $c$ est encore un ensemble. Or la classe des
représentants des classes d'isomorphismes de couples réalisant le
relèvement partiel s'identifie à une sous-classe de cette union
disjointe, ce qui en fait un ensemble.\\ 
CQFD.

\begin{lem}\label{ast2}
Soient $\mathcal{C}$ une catégorie et $\mathcal{C}'$ une sous-catégorie de
$\mathcal{C}$ vérifiant les hypothèses du lemme 
précédent. Si $\mathcal{C}$ admet les colimites séquentielles transfinies,
 que $\mathcal{C}'$ admet les colimites séquentielles transfinies de
$\mathcal{C}$ pour colimites séquentielles transfinies et si, pour tout carré
commutatif $(j,a,b,p)$ avec $j$ morphisme de $\mathcal{C}'$ et $p$ fixé,
l'ensemble contenant un représentant de chaque 
classe d'isomorphismes de couples $(C,c)$ réalisant un relèvement partiel du
diagramme $(j,a,b,p)$ et tels que l'objet $C$ est non isomorphe à la source de
$j$ et factorise $j$ dans $\mathcal{C}'$, est non vide, alors
$(j,a,b,p)$ admet un relèvement global.
\end{lem}
{\it Preuve :}\\
Par propriété 2) du lemme précédent, les représentants des classes d'isomorphismes
de couples $(C,c)$ forment bien un ensemble. On ordonne cet ensemble de la manière suivante :
$(C,c)$ est inférieur à $(C',c')$ s'il existe un morphisme $g$ de $\mathcal{C}'$ allant
de $C$ à $C'$ tel que $c'\circ g=c$. Montrons que cet ensemble est inductif. Considérons une
chaîne, i.e. un sous-ensemble totalement ordonné, $(C_i,c_i)_{I\in I}$. Posons
$$C=\colim{i\in I}C_i \mbox{ et } c=\colim{i\in I}c_i.$$ Comme chaque $C_i$ factorise $j$ dans $\mathcal{C}'$
alors, par propriété de colimite et par hypothèse sur $\mathcal{C}'$, $C$
factorise $j$ dans $\mathcal{C}'$. De même, comme chaque $(C_i,c_i)$ réalise un
relèvement partiel du diagramme $(j,a,b,p)$, alors, par propriété de
colimite, $(C,c)$ aussi. En outre $C$ n'est pas isomorphe à la source de $j$ sinon
comme les isomorphismes ne se factorisent dans $\mathcal{C}'$ qu'en
isomorphismes alors chaque $C_i$ serait isomorphe à la source de $j$, ce qui
n'est pas le cas. Donc $(C,c)$ appartient à l'ensemble et, comme colimite dans
$\mathcal{C}'$, est plus grand que tous les $(C_i,c_i)$, il est donc la borne
supérieure de la chaîne, ce qui montre que l'ensemble est inductif.\\
\\
Nous pouvons donc appliquer le lemme de Zorn à l'ensemble et l'on obtient
qu'il existe un élément maximal $(C',c')$. Supposons que $C'$ ne soit pas isomorphe au but de $j$, alors la partie basse du diagramme
(RP) pour $C'$ est un carré commutatif dont le morphisme de $C'$ vers le but
de $j$ est dans $\mathcal{C}'$. Par hypothèse, il existe un couple $(C'',c'')$
réalisant un relèvement partiel pour la partie basse de (RP), avec $C''$ non
isomorphe à $C'$ et 
factorisant ce morphisme dans $\mathcal{C}'$. Mais $C''$ factorise donc $j$ dans
$\mathcal{C'}$, n'est pas isomorphe à la source de $j$, sans quoi il serait
isomorphe à $C'$, et avec $c''$ réalise en fait un relèvement partiel de $(j,a,b,p)$.
Donc $(C'',c'')$ est dans l'ensemble et strictement supérieur à l'élément
maximal $(C',c')$, ce qui est absurde. Donc l'élément maximal est bien isomorphe
au but de $j$ et donc $(j,a,b,p)$ possède un relèvement global.\\
CQFD.\\

\subsection{Caractérisation des $\mathcal{I}$-cofibrations}

Munis de ces lemmes techniques, nous pouvons désormais démontrer que
$\mathcal{I}$ engendre les cofibrations niveau par niveau. Pour cela nous allons
d'abord montrer que la classe des cofibrations vérifie la propriété du
lemme~\ref{ast2}. Toutefois, nous ne savons montrer ce résultat que dans le
cas où les cofibrations de $\mathcal{C}$ sont les monomorphismes.

\begin{lem}
Soit $\mathcal{C}$ une catégorie \discret~munie d'une structure de
\cmf~engendrée par monomorphismes vérifiant les propriétés suivantes :
\item -les cofibrations génératrices ont leurs
sources et buts $\alpha$-petits, pour un cardinal régulier $\alpha$ plus grand
qu' $\aleph_0$,
\item -les cofibrations génératrices ont leurs buts connexes,
\item -les réunions de deux
sous-objets sont exactement les sommes amalgamées de ces sous-objets au-dessus
de leur intersection.\\
\\ 
Soit $\mathcal{I}$ l'ensemble défini dans le
lemme~\ref{I}. Alors pour toute cofibration niveau par niveau $j:A\rightarrow B$ de
$\mathcal{C-PC}$, qui n'est pas un isomorphisme, il existe un diagramme :
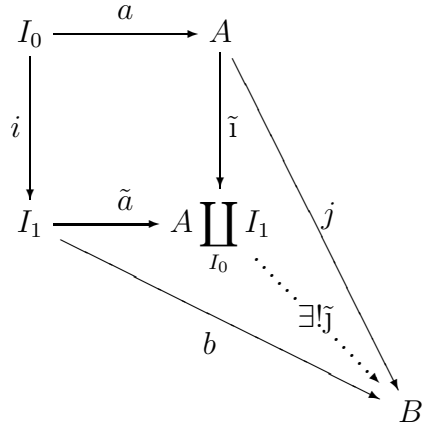
\begin{diagram}
I_0 & \rTo^a & A & & \\
\dTo^{i} & & \dTo_{\tilde{\i}}& \rdTo(2,4)^{j} & \\
I_1 & \rTo^{\tilde{a}} & A\coprod_{I_0}I_1 & & \\
 & \rdTo(4,2)_{b} & & \rdDotsto~{\exists ! \tilde{\j}} & \\
 & & & & B \\
\end{diagram}
avec $i$ appartenant à $\mathcal{I}$ et tel que $\tilde{\i}$ et $\tilde{\j}$ sont des
cofibrations niveau par niveau et que la somme amalgamée $A\coprod_{I_0}I_1$ ne soit pas isomorphe à $A$.
\end{lem}
{\it Preuve :}\\
Si $j$ n'est pas injective au niveau des objets, il suffit de prendre
$\emptyset\rightarrow *$ pour $i$ et $b$ induit par un objet de $B$ qui n'est
pas dans l'image de $A$ par $j$ pour avoir le résultat. Sinon soit $n$
l'entier tel que, pour tout $m<n$, $j_m$ est un isomorphisme mais $j_n$ n'est
pas un isomorphisme. Comme $j_n$ est une cofibration et que $\mathcal{C}$ est
une \cmfcof, alors il existe dans
$\mathcal{C}$ un diagramme :
\begin{diagram}
I_0' & \rTo^{a'} & A_n & & \\
\dTo^{i'} & & \dTo_{\tilde{\i}'}& \rdTo(2,4)^{j_n} & \\
I_1' & \rTo^{\tilde{a}'} & A_n\coprod_{I_0'}I_1' & & \\
 & \rdTo(4,2)_{b'} & & \rdDotsto~{\exists ! \widetilde{\j_n}} & \\
 & & & & B_n \\
\end{diagram}
avec $i'$ cofibration génératrice, $\tilde{\i}'$ et $\widetilde{\j_n}$ des
cofibrations et la somme amalgamée du diagramme non
isomorphe à $A_n$. On étend ce diagramme aux bords de $A_n$ et $B_n$.
$f$ étant un isomorphisme jusqu'au rang $n-1$,
il induit un isomorphisme sur les bords de $A_n$ et $B_n$ et donc il existe un
morphisme naturel de $I'_1$ vers le bord de $A_n$ faisant commuter le diagramme
étendu, ce qui nous donne, par définition des flèches de type $Attach_n(i')$ et par
connexité des buts des cofibrations génératrices de $\mathcal{C}$, le diagramme
suivant :
\begin{diagram}
\Delta[n]\Theta I'_0\coprod_{\partial\Delta[n]\Theta I'_0}\partial\Delta[n]\Theta
I'_1 & & \rTo^a & & A & & \\
\dTo^{Attach_n(i')} & & & & \dTo^{\widetilde{Attach_n(i')}}& \rdTo(2,4)^{j} & \\
\Delta[n]\Theta I'_1 & & \rTo^{\tilde{a}} & & \tilde{A} & & \\
 & \rdTo(6,2)_{b} & & & & \rdDotsto~{\exists ! \tilde{\j}} & \\
 & & & & & & B \\
\end{diagram} 
où $a$ et $b$ sont induits respectivement par $a'$ et $b'$. Comme la somme
amalgamée de $A_n$ par $i'$ n'est pas isomorphe à $A_n$, il vient que la
somme amalgamée $\tilde{A}$ de $A$ par $Attach_n(i')$ n'est pas isomorphe à
$A$. Comme niveau par niveau $Attach_n(i')$ est un coproduit dans la catégorie
des morphismes de cofibrations de $\mathcal{C}$ ($i'$ et
l'identité) alors $Attach_n(i')$ est bien
une cofibration niveau par niveau par lemme~\ref{cofst} et, par le
corollaire~\ref{cofnpnst}, sa somme amalgamée le long de $a$ aussi. Regardons de plus
près le morphisme $\tilde{\j}$. Pour $m<n$, on a $\tilde{\j}_m$ n'est autre
que $j_m$ donc une cofibration de $\mathcal{C}$. Au rang n, on a que
$\tilde{\j}_n$ est $\widetilde{\j_n}$ qui est une cofibration. Pour $m>n$, on va
utiliser la notion de réunion de sous-objets. Comme $j$ est une cofibration
niveau par niveau et que l'on a supposé que les cofibrations de $\mathcal{C}$
sont les monomorphismes, alors $A_m$ est un sous-objet de $B_m$. On a montré
en outre que $\widetilde{\j_n}$ est une cofibration, donc de même il vient
$A_n\coprod_{I'_0} I'_1$ est un sous-objet de $B_n$. Or les relations entre
faces et dégénérescences d'un objet simplicial entraînent que les
dégénérescences sont des monomorphismes et donc, par chaque
dégénérescence de $B_n$ vers $B_m$, $A_n\coprod_{I'_0} I'_1$
est un sous-objet de $B_m$. Pour montrer que $\tilde{\j}_m$ est une cofibration, nous
allons montrer que $\tilde{A}_m$ est une réunion de sous-objets de $B_m$.
Comme on a supposé que les réunions de sous-objets dans $\mathcal{C}$ sont
exactement les sommes amalgamées des sous-objets au-dessus de leur
intersection, il suffit de montrer que $\tilde{A}_m$ est une telle somme
amalgamée. Or $\tilde{A}_m$ est la somme amalgamée du sous-objet $A_m$ par un
coproduit de $A_n\coprod_{I'_0} I'_1$ au-dessus d'un coproduit de $A_n$, ces
coproduits étant indexés par les applications dégénérées de $A_n$
dans $A_m$. On décompose cette multiple somme amalgamée en une suite de
sommes amalgamées par un seul exemplaire de $A_n\coprod_{I'_0} I'_1$ à
chaque étape. Or à chaque étape, on vérifie que c'est bien une réunion
de sous-objets en utilisant en particulier le fait que l'intersection dans $B_m$ de $A_m$ par $A_n\coprod_{I'_0} I'_1$ n'est
autre que $A_n$ et celle entre deux exemplaires de $A_n\coprod_{I'_0} I'_1$ est
$A_{n-1}$. On obtient ainsi que $\tilde{A}_m$ est une réunion de sous-objets de
$B_m$, donc est lui-même un sous-objet de $B_m$, et ainsi $\tilde{\j}_m$, étant l'inclusion du
sous-objet, est une cofibration. Comme par ailleurs $Attach_n(i')$ est
élément de $\mathcal{I}$, on a montré le lemme dans
le second cas.\\
CQFD.\\

Grâce à ce résultat, nous pouvons appliquer le lemme~\ref{ast2}, qui nous
permettra de montrer que les morphismes $\mathcal{I}$-injectifs
sont exactement ceux ayant la \prd~par rapport aux monomorphismes, ce qui
montrera la première partie de l'hypothèse 4) car les morphismes de
$\mathcal{J}$ sont des monomorphismes.

\begin{lem}\label{p4a}
Si $\mathcal{C}$ est une catégorie \discret~munie d'une structure de
\cmf~engendrée par monomorphismes qui avec la sous-catégorie des
monomorphismes vérifie les hypothèses suivantes :
\item -les cofibrations génératrices ont leurs
sources et buts $\alpha$-petits, pour un cardinal régulier $\alpha$ plus grand
qu' $\aleph_0$,
\item -les cofibrations génératrices ont leurs buts connexes,
\item -les réunions de deux
sous-objets sont exactement les sommes amalgamées de ces sous-objets au-dessus
de leur intersection,
\item -pour tout monomorphisme $j:A\rightarrow B$, la classe contenant un représentant de chaque classe d'isomorphismes d'objets $C$ par lequel $j$ se
factorise dans la sous-catégorie des monomorphismes est un ensemble,
\item -les colimites séquentielles transfinies de la sous-catégorie des
monomorphismes de $\mathcal{C}$ existent et sont les colimites séquentielles
transfinies de $\mathcal{C}$.\\
\\
Alors les morphismes
$\mathcal{I}$-injectifs sont exactement ceux qui ont la \prd~par rapport aux
monomorphismes. En particulier tout morphisme
$\mathcal{I}$-injectif est $\mathcal{J}$-injectif. 
\end{lem}
{\it Preuve :}\\
Les monomorphismes vérifient trivialement que les isomorphismes ne se factorisent qu'en isomorphismes dans la sous-catégorie des monomorphismes. Avec cette propriété et les hypothèses de ce lemme, on obtient que la catégorie $\mathcal{C}$ vérifie les hypothèses du lemme~\ref{ast1}. En appliquant la propriété 1) du lemme~\ref{ast1}, on obtient que
$\mathcal{C-PC}$ et sa sous-catégorie des cofibrations niveau par niveau
vérifie les hypothèses du lemme~\ref{ast1}. De plus, les colimites de
$\mathcal{C-PC}$ étant calculées niveau par niveau, on conserve sur
$\mathcal{C-PC}$ la propriété que les colimites
séquentielles transfinies pour les cofibrations niveau par niveau sont celles
de $\mathcal{C-PC}$.\\
\\
Considérons maintenant le diagramme suivant :
\begin{diagram}
A & \rTo^{a} & X\\
\dTo^{j} & & \dTo_{p}\\
B & \rTo_{b} & Y \\
\end{diagram}
avec $j$ une cofibration niveau par niveau et $p$ un morphisme
$\mathcal{I}$-injectif. Si $j$ est un isomorphisme, le diagramme admet un
relèvement naturel. Sinon on applique le lemme précédent et l'on obtient
le diagramme suivant :
\begin{diagram}
I_0 & \rTo^a & A & & \\
\dTo^{i} & & \dTo_{\tilde{\i}}& \rdTo(2,4)^{j} & \\
I_1 & \rTo^{\tilde{a}} & A\coprod_{I_0}I_1 & & \\
 & \rdTo(4,2)_{b} & & \rdDotsto~{\exists ! \tilde{\j}} & \\
 & & & & B \\
\end{diagram}
avec $i$ appartenant à $\mathcal{I}$, $\tilde{\i}$ et $\tilde{\j}$ des cofibrations niveau par
niveau et la somme amalgamée du diagramme, que
l'on notera $C$, non isomorphe à $A$. Comme $p$ est $\mathcal{I}$-injectif, il a la \prd~par rapport à
$i$ ainsi que par rapport à la somme amalgamée $\tilde{\i}$ de $i$ le long de $a$, par le lemme~\ref{prst}.
Donc $C$ réalise un relèvement partiel du diagramme $(j,a,b,p)$. En
appliquant le lemme~\ref{ast2}, il vient que  $(j,a,b,p)$ a un relèvement
global. Donc les morphismes $\mathcal{I}$-injectifs ont la \prd~par rapport aux
cofibrations niveau par niveau.\\
\\
En outre les morphismes de $\mathcal{I}$ sont niveau par niveau des coproduits
de cofibrations élémentaires, donc ce sont des cofibrations niveau par
niveau par lemme~\ref{cofst}. De là il vient que les morphismes ayant la \prd~par rapport aux
cofibrations niveau par niveau se relèvent en particulier par rapport aux
morphismes de $\mathcal{I}$ et sont donc $\mathcal{I}$-injectifs.\\
Enfin comme par définition $\mathcal{J}$ est constitué de monomorphismes niveau par niveau et
que les morphismes $\mathcal{I}$-injectifs se relèvent par rapport aux
monomorphismes, alors ils se relèvent en particulier par rapport aux
morphismes de $\mathcal{J}$, ce qui en fait des morphismes $\mathcal{J}$-injectifs.\\
CQFD.\\

De ce résultat découle le fait que l'ensemble $\mathcal{I}$ engendre bien
les monomorphismes, ce qui justifie son nom d'ensemble générateur des
cofibrations.

\begin{cor}\label{Icof2}
Sous les hypothèses du lemme~\ref{p4a}, $\mathcal{I}$ engendre les
monomorphismes, c'est-à-dire que les monomorphismes sont exactement les
$\mathcal{I}$-cofibrations. En particulier, les monomorphismes sont des
rétracts de colimites séquentielles transfinies de sommes
amalgamées de morphismes de $\mathcal{I}$.
\end{cor}
{\it Preuve :} 
Par le lemme~\ref{I}, $\mathcal{I}$ permet
l'argument du petit objet et, par corollaire~\ref{Icof}, il vient que les
$\mathcal{I}$-cofibrations sont exactement les rétracts de colimites transfinies de sommes
amalgamées de morphismes de $\mathcal{I}$. Or les morphismes de $\mathcal{I}$
sont des cofibrations niveau par niveau et par le corollaire~\ref{cofnpnst}, les cofibrations niveau
par niveau sont stables par rétract, somme amalgamée le long d'un morphisme
quelconque et colimite
séquentielle transfinie. D'où les $\mathcal{I}$-cofibrations sont des cofibrations niveau par niveau. En
outre le lemme précédent montre que les cofibrations niveau par niveau ont la \prg~par
rapport aux morphismes $\mathcal{I}$-injectifs, ce qui en fait des
$\mathcal{I}$-cofibrations, par définition même de
$\mathcal{I}$-cofibrations. Donc les cofibrations niveau par niveau sont exactement les
$\mathcal{I}$-cofibrations et donc $\mathcal{I}$ engendre bien les cofibrations niveau par niveau.\\
CQFD.\\

La première partie de l'hypothèse 4) étant traitée, il ne reste plus
qu'à traiter la seconde. De même que la première partie portant sur les
cofibrations avait nécessité des rappels, la seconde partie concernant les
\eqs~faibles demande elle aussi quelques rappels sur la stabilité des
\eqs~faibles par colimite dans une \cmf.

\subsection{Stabilité des \eqs~faibles dans une \cmf}

La stabilité des \eqs~faibles dans une \cmf~est un sujet assez délicat.
Très souvent elle est couplée à des propriétés des cofibrations et à
la notion de \cmf~propre à gauche.

\begin{defin}\index{catégorie de modèles fermée!propre à gauche}
Une \cmf~est propre à gauche si toute somme amalgamée d'une \eq~faible le long
d'une cofibration est une \eq~faible.
\end{defin}

Une façon simple pour une \cmf~d'être propre à gauche est d'avoir tous
ses objets cofibrants, ce qui est le cas de nos données de Segal faciles.

\begin{lem}
Une \cmf~dont tous les objets sont cofibrants est propre à gauche.
\end{lem}
{\it Preuve :} dans \cite{h}\\

Le premier résultat que nous allons énoncer est le lemme du cube qui
concerne la stabilité des \eqs~faibles par somme amalgamée dans la
catégorie des morphismes. Ce lemme nécessite l'usage d'une \cmf~propre à
gauche.

\begin{prop}\label{eqcube}
Soit $\mathcal{C}$ une \cmf~propre à gauche, considérons le diagramme suivant :
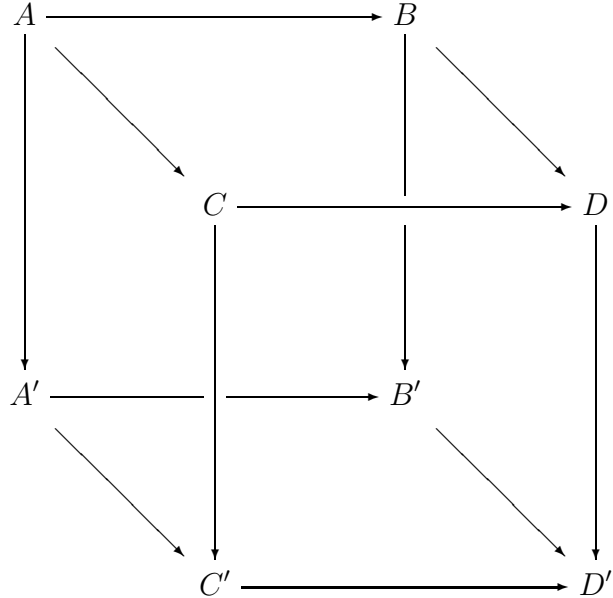
\begin{diagram}
A  & & \rTo & & B & & \\
   & \rdTo & & & \dLine & \rdTo & \\
\dTo & & C & & \rTo &  & D \\
   & &  & & \dTo & & \\
A' & \rLine & \dTo & \rTo & B' & & \dTo \\
   & \rdTo & & & & \rdTo & \\
   & & C' & & \rTo &  & D' \\ 
\end{diagram}
o\'u $D$ et $D'$ sont les sommes amalgamées $B\coprod_A C$ et $B'\coprod_A'
C'$.\\ Si $A\rightarrow A'$, $B\rightarrow B'$ et $C\rightarrow C'$ sont des
\eqs~faibles et si $A\rightarrow C$ et $A'\rightarrow C'$ sont des
cofibrations, alors $D\rightarrow D'$ est aussi une \eq~faible.
\end{prop}
{\it Preuve :} dans \cite{h}.\\

Un autre résultat important est la stabilité des \eqs~faibles par colimite
séquentielle transfinie dans la catégorie des morphismes. Cette fois même
le fait d'être propre à gauche ne suffit pas. Il faudra que la \cmf~soit, en
plus d'être propre à gauche, simpliciale !

\begin{prop}\label{eqcolim}
Soit $\mathcal{C}$ une \cmf~simpliciale propre à gauche. Soit $\lambda$ un
ordinal transfini et $g_.:X_.\rightarrow Y_.$ un morphisme de
$\lambda$-séquences tel que pour tout ordinal $\alpha<\lambda$,
$g_{\alpha}:X_{\alpha}\rightarrow Y_{\alpha}$ soit une \eq~faible et que les morphismes
$X_{\alpha}\rightarrow X_{\alpha+1}$ et $Y_{\alpha}\rightarrow Y_{\alpha+1}$
soient des cofibrations, alors le morphisme induit $\colimite{\alpha<\lambda}
X_{\alpha}\rightarrow \colimite{\alpha<\lambda}
Y_{\alpha}$ est une \eq~faible.
\end{prop}
{\it Preuve :} dans \cite{h}.\\

De ce résultat découle la stabilité des \eqs~faibles par coproduit dans la
catégories des morphismes, résultat qui lui aussi nécessite la
simplicialité et la propreté à gauche de la \cmf.

\begin{cor}\label{eqcoprod}
Soit $\mathcal{C}$ une \cmf~simpliciale dont tous les objets sont cofibrants.
Alors un petit coproduit d'\eqs~faibles est une \eq~faible.
\end{cor}
{\it Preuve :}\\
Le coproduit étant petit, il est indexé par un ensemble, que
l'on peut donc bien ordonner par le théorème de Zermelo. Ainsi notre
coproduit est indexé par un ordinal que l'on notera $\lambda$ et donc il donne
lieu à un morphisme $g_.:X_.\rightarrow Y_.$ de $\lambda$-séquences, chacune d'elles étant une
colimite séquentielle transfinie de coproduits de deux flèches, i.e. de
sommes amalgamées par des morphismes de source l'objet initial. Comme par
hypothèse, dans $\mathcal{C}$ tous les objets sont cofibrants, ces morphismes
sont des cofibrations et, par stabilité des cofibrations par somme amalgamée le long d'un
morphisme dans une \cmf, $g$ vérifie la deuxième condition de la proposition
précédente. Montrons par récurrence transfinie que $g$ vérifie aussi la
première condition. Par hypothèse, $g_0$ est une \eq~faible. Soit $\alpha<\lambda$
un ordinal, supposons que l'on ait montré que $g_{\alpha}$ est une \eq~faible. On a
alors le diagramme suivant :
\begin{diagram}
\emptyset  & & \rTo & & B & & \\
   & \rdTo & & & \dLine & \rdTo & \\
\dTo & & X_{\alpha} & & \rTo &  & X_{\alpha+1} \\
   & &  & & \dTo & & \\
\emptyset & \rLine & \dTo^{g_{\alpha}} & \rTo & B' & & \dTo_{g_{\alpha+1}} \\
   & \rdTo & & & & \rdTo & \\
   & & Y_{\alpha} & & \rTo &  & Y_{\alpha+1} \\ 
\end{diagram}
avec $X_{\alpha+1}$ et $Y_{\alpha+1}$ sommes amalgamées des faces horizontales
du cube et $B\rightarrow B'$ l'une des \eqs~faibles du coproduit. Comme dans
$\mathcal{C}$, tous les objets sont cofibrants, alors $\mathcal{C}$ est une
\cmf~propre à gauche et, par la proposition~\ref{eqcube}, il vient que
$g_{\alpha+1}$ est une \eq~faible. Soit $\gamma<\lambda$ un ordinal limite. Supposons
que pour tout ordinal $\alpha<\gamma$, on ait montré que $g_{\alpha}$ est une
\eq~faible. Alors $g_{\gamma}$ est le morphisme induit par le morphisme de
$\gamma$-séquence $g_.$ qui vérifie bien les hypothèses de la
proposition~\ref{eqcolim}. Comme $\mathcal{C}$ est une \cmf~simpliciale propre
à gauche, il vient que $g_{\gamma}$ est bien une \eq~faible. On a donc montré par
récurrence transfinie que $g$ vérifie aussi l'hypothèse 1) du lemme~\ref{eqcolim}, d'où il vient que $g_{\lambda}$ est une \eq~faible. Or $g_{\lambda}$ est
justement le coproduit d'\eqs.\\
CQFD.

\subsection{Les morphismes $\mathcal{I}$-injectifs sont des \eqs\\ faibles}

Comme l'on vient de montrer la première partie de l'hypothèse 4), il ne nous
reste plus qu'à montrer la seconde partie à savoir que les morphismes
$\mathcal{I}$-injectifs sont des \eqs~faibles de \precats. Pour cette
démonstra\-tion, nous avons besoin du résultat intermédiaire suivant : si
$f:A\rightarrow C$ est un morphisme de \precats~surjectif sur les objets et tel
que, pour tout entier $q>0$ et pour tout $q+1$-uplet $(a_0,\ldots,a_q)$ d'objets
de $A$, le morphisme $f_q(a_0,\ldots,a_q)$ est une \eq~faible dans
$\mathcal{C}$, alors le morphisme $f$ est une \eq~faible de \precats. Pour
montrer ce résultat, nous allons utiliser la construction Bigcat et suivre à
chaque étape ce que devient le morphisme $f$.

\begin{lem}
Soit $(\mathcal{C},\mathcal{F}_1,\mathcal{F}_2)$ une donnée de Segal facile
vérifiant la propriété suivante :\\
les cofibrations sont stables par produit fibré au-dessus d'un objet
discret.\\
Alors pour tout morphisme de \precats~$f:A\rightarrow C$ surjectif sur les
objets et tel que, pour tout entier
$q>0$ et pour tout $q$-uplet $(a_0,\ldots,a_q)$ d'objets de $A$, 
$f_q(a_0,\ldots,a_q):A_q(a_0,\ldots,a_q)\rightarrow C_q(f(a_0),\ldots,f(a_q))$
est une \eq~faible, on a que, pour tout entier $m\geq 2$,
$Cat(1,m)(A)\rightarrow Cat(1,m)(C)$ possède la même propriété que $f$.
\end{lem}
{\it Preuve :}\\
Tout d'abord on remarque que, comme $\mathcal{C}$ est une catégorie \discret,
tout morphisme et tout diagramme équivaut au coproduit de ses fibres. De ce
fait, $Cat(1,m)(A)_q(a_0,\ldots,a_q)$ n'est autre que la somme amalgamée de
$A_q(a_0,\ldots,a_q)$ par un coproduit de $\eta_A(a_i,a_j)$ et de
$\phi_A(a_{i_0},\dots,a_{i_m})$, où on a $0\leq i<j\leq q$ et $0\leq
i_0\leq\ldots\leq i_m\leq q$ avec $i_m\neq i_0+1$. En outre $\eta_A(a_i,a_j)$ et
$\psi_A(a_{i_0},\dots,a_{i_m})$ sont des \eqs~faibles comme produits fibrés au-dessus d'un objet
discret d'\eqs~faibles (respectivement $\eta_A$ et $\psi_A$, d'après la
démonstration du lemme~\ref{eqcat1m}). D'où pour la
même raison, il en est de même pour le produit des
$\eta_A(a_{i_k},a_{i_{k+1}})$
avec $k$ entre 0 et $m-1$. Ce raisonnement reste valable pour la \precat~$C$.\\
\\
Ensuite on remarque que si la source d'une flèche génératrice
s'envoie dans $A$ alors en composant par $f:A\rightarrow C$ elle s'envoie aussi
dans $C$. Ainsi toute somme amalgamée de $A_q$, pour un certain entier $q$,
par un coproduit de flèches génératrices s'envoyant dans $A_q$,
possède un morphisme vers le même type de sommes amalgamées pour
$C_q$. De ce fait il existe un morphisme $\eta_f$ de $B_A$ vers $B_C$ tel que
l'on ait $\eta_f\circ\eta_A=\eta_C\circ f_1$ et $\nu_C\circ\eta_f=f_0\times
f_0\circ\nu_A$. De même il existe un morphisme $\phi_f$ de $P_A$ vers $P_C$
tel que $\phi_f\circ\phi_A=\phi_C\circ f_m$ et
$\psi_C\circ\phi_f=\eta_f\times_{f_0}\ldots\times_{f_0}\eta_f\circ\psi_A$. 
Par les égalités $\eta_C(f(a_i),f(a_j))\circ
f_1(a_i,a_j)=\eta_f(a_i,a_j)\circ \eta_A(a_i,a_j)$ et
$\psi_C(f(a_{i_0}),\dots,f(a_{i_m}))\circ\phi_f(a_{i_0},\dots,a_{i_m})=\eta_f(a_{i_0},a_{i_1})\times\ldots\times\eta_f(a_{i_{m-1}},a_{i_m})\circ\psi_A(a_{i_0},\dots,a_{i_m})$
ainsi que par "trois pour deux" dans la
\cmf~$\mathcal{C}$, on a que $\eta_f(a_i,a_j)$ et
$\phi_f(a_{i_0},\dots,a_{i_m})$ sont des \eqs~faibles. Comme les morphismes
$\eta_A(a_i,a_j)$, $\phi_A(a_{i_0},\dots,a_{i_m})$, $\eta_C(a_i,a_j)$ et 
$\phi_C(a_{i_0},\dots,a_{i_m})$ ainsi que leurs coproduits sont des
cofibrations, alors par la proposition~\ref{eqcube}, il vient que
$Cat(1,m)(f)_q(a_0,\ldots,a_q)$ est une \eq~faible. En outre, comme la
construction $Cat(1,m)$ laisse les objets invariants, $Cat(1,m)(f)$ est bien
surjective sur les objets.\\
CQFD.

\begin{prop}
Soit $(\mathcal{C},\mathcal{F}_1,\mathcal{F}_2)$ une donnée de Segal facile
vérifiant les propriétés suivantes :
\item 1) la \cmf~sous-jacente $\mathcal{C}$ est une \cmf~simpliciale,
\item 2) les cofibrations sont stables par produit fibré au-dessus d'un objet
discret.\\
\\
Alors pour tout morphisme de \precats~$f:A\rightarrow C$ surjectif sur les
objets et tel que, pour tout entier
$q>0$ et pour tout $q$-uplet $(a_0,\ldots,a_q)$ d'objets de $A$, 
$f_q(a_0,\ldots,a_q):A_q(a_0,\ldots,a_q)\rightarrow C_q(f(a_0),\ldots,f(a_q))$
est une \eq~faible, on a que 
$A\rightarrow C$ est une \eq~faible dans $\mathcal{C-PC}$.
\end{prop}
{\it Preuve :}\\ 
Par récurrence transfinie, on montre que pour tout ordinal
$\beta\leq T$ le morphisme $f^{\beta}:A^{\beta}\rightarrow C^{\beta}$ a la même
propriété que $f$, où $f^{\beta}$ dénote la $\beta$-ième étape de
$Bigcat(f)$. Le passage d'un ordinal $\beta<T$ à son successeur $\beta+1$ est
l'application directe du lemme précédent.\\
\\
Supposons maintenant l'hypothèse de
récurrence vraie pour tout ordinal $\gamma$ inférieur à un ordinal limite
$\beta\leq T$. Comme les fibres des morphismes $A_q\rightarrow A^{\gamma}_q$
et $C_q\rightarrow C^{\gamma}_q$ sont des cofibrations (démonstration de la
proposition~\ref{eqbig}), par la proposition~\ref{eqcolim}, il vient que
les fibres de $f^{\beta}_q:A^{\beta}_q\rightarrow C^{\beta}_q$ sont des
\eqs~faibles. Comme en outre les constructions $Cat(1,m)$ préservent les objets,
le morphisme $f^{\beta}$ reste surjectif, ce qui montre l'hypothèse de
récurrence pour un ordinal limite.\\ 
\\
Par récurrence transfinie, on a montré que pour tout ordinal
$\beta\leq T$ le morphisme $f^{\beta}:A^{\beta}\rightarrow C^{\beta}$ a la même
propriété que $f$. C'est donc en particulier le cas pour $Bigcat(f)=f^{T}$.
Ainsi les fibres de $Bigcat(f)_q:Bigcat(A)_q\rightarrow Bigcat(C)_q$ sont des
\eqs~faibles entre \obcs, i.e. niveau par niveau ce sont des \eqcs~d'\obcs. C'est donc en particulier
le cas pour les fibres de $Bigcat(f)_1$ et comme $Bigcat(A)$ et
$Bigcat(C)$ sont des \cats, il vient que
$Bigcat(f):Bigcat(A)\rightarrow Bigcat(C)$ est pleinement fidèle. En outre
comme $Bigcat(f)$ est comme $f$, en particulier surjective sur les objets, alors
$Bigcat(f)$ est essentiellement surjective. Ainsi $Bigcat(f)$
est une \eq~de \cats, donc, par le lemme~\ref{bigcateq}, ceci nous donne que $f$
est une \eq~faible.\\
CQFD.\\

Ce résultat intermédiaire montré, nous sommes à même de terminer la
vérification de l'hypothèse 4) du lemme de reconnaissance~\ref{reco}. 

\begin{lem}\label{p4}
Soit $(\mathcal{C},\mathcal{F}_1,\mathcal{F}_2)$ une donnée de Segal facile
vérifiant les propriétés suivantes :
\item 1) la \cmf~sous-jacente $\mathcal{C}$ est une \cmf~simpliciale engendrée
par monomorphismes,
\item 2) les cofibrations génératrices ont leurs
sources et buts $\alpha$-petits,
\item 3) les cofibrations génératrices ont leurs buts connexes,
\item 4) les réunions de deux
sous-objets sont exactement les sommes amalgamées de ces sous-objets au-dessus
de leur intersection,
\item 5) pour tout monomorphisme $j:A\rightarrow B$, la classe contenant un représentant de chaque classe d'isomorphismes d'objets $C$ par lequel $j$ se
factorise dans la sous-catégorie des monomorphismes est un ensemble,
\item 6) les colimites séquentielles transfinies de la sous-catégorie des
monomorphismes de $\mathcal{C}$ existent et sont les colimites séquentielles
transfinies de $\mathcal{C}$.\\
\\
Alors tout morphisme $\mathcal{I}$-injectif est à la fois $\mathcal{J}$-injectif et
une \eq~faible.
\end{lem} 
{\it Preuve :}\\
Les hypothèses que l'on fait sur la donnée de Segal facile permettent de
vérifier les hypothèses du lemme~\ref{p4a}. On obtient donc que tout morphisme $\mathcal{I}$-injectif est
$\mathcal{J}$-injectif.\\
Toujours d'après ce lemme, les morphismes $\mathcal{I}$-injectifs ont la \prd~par rapport
aux monomorphismes. Soit $f:A\rightarrow B$ un morphisme $\mathcal{I}$-injectif.
Ce morphisme est surjectif sur les objets. En effet soit $b$ un objet de $B$, on
obtient le diagramme commutatif suivant :
\begin{diagram}
\emptyset & \rTo & A\\
\dTo & \ruDotsto_{\exists a} & \dTo_{f}\\
* & \rTo_{b} & B\\
\end{diagram}
Comme $\emptyset\rightarrow *$ est un monomorphisme et que $f$ est
$\mathcal{I}$-injectif, le relèvement existe, ce qui signifie qu'il existe un
objet $a$ de $A$ tel que $f(a)$ est égal à $b$. Donc $f$ est bien surjectif
sur les objets.\\
\\
Montrons maintenant que pour tout entier $q>0$ et pour tout $q$-uplet
$(a_0,\ldots,a_q)$ d'objets de $A$, $f_q(a_0,\ldots,a_q)$ est une fibration
triviale. Considérons donc un diagramme commutatif du type suivant :
\begin{diagram}
X & \rTo & A_q(a_0,\ldots,a_q)\\
\dTo^i & \ruDotsto_{\exists } & \dTo_{f_q(a_0,\ldots,a_q)}\\
Y & \rTo & B_q(f(a_0),\ldots,f(a_q))\\
\end{diagram}
avec $i:X\rightarrow Y$ un monomorphisme de $\mathcal{C}$. Alors ce diagramme équivaut, par
propriété de la construction $\Theta$, au diagramme suivant :
\begin{diagram}
\Delta[q]\Theta X & \rTo & A\\
\dTo^{\Delta[q]\Theta i} & \ruDotsto_{\exists } & \dTo_{f}\\
\Delta[q]\Theta Y & \rTo & B_q\\
\end{diagram}
Comme $\Delta[q]\Theta i$ est niveau par niveau un coproduit du monomorphisme
$i$, il vient que $\Delta[q]\Theta i$ est un monomorphisme niveau par niveau. Or $f$ est
$\mathcal{I}$-injectif, donc le relèvement existe, ce qui montre par
équivalence des diagrammes que $f_q(a_0,\ldots,a_q)$ a la \prd~par rapport aux
monomorphismes, qui sont justement les cofibrations, donc c'est une fibration triviale dans la \cmf~$\mathcal{C}$.\\
\\
On vient donc de montrer que tout morphisme $\mathcal{I}$-injectif est surjectif
sur les objets et que pour tout entier $q>0$ et pour tout $q$-uplet
$(a_0,\ldots,a_q)$ d'objets de $A$, $f_q(a_0,\ldots,a_q)$ est une fibration
triviale donc en particulier une \eq~faible. Les monomorphismes étant stables par produit fibré au-dessus d'un objet discret, toutes les hypothèses de la proposition
précédente sont donc vérifiées et l'on obtient ainsi que $f$ est lui-même une \eq~faible pour
$\mathcal{C-PC}$.\\
CQFD.\\

L'hypothèse 4) du lemme de reconnaissance~\ref{reco} est donc vérifiée
sous des hypothèses très fortes (et très lourdes) comme le fait que $\mathcal{C}$ doit être
une \cmf~simpliciale. Cependant quatre des six hypothèses du lemme de
reconnaissance de la structure de \cmfcof~ont été vérifiées. On n'a donc
plus le choix que de traiter l'hypothèse 3).\\

Cette hypothèse 3) est de loin la plus difficile à vérifier. On doit en
effet montrer que toute $\mathcal{J}$-cofibration est à la fois une
$\mathcal{I}$-cofibration et une \eq~faible, c'est-à-dire une cofibration
triviale, vu que l'on a montré que les $\mathcal{I}$-cofibrations sont les
monomorphismes. Comme on a vu que $\mathcal{J}$
permet l'argument du petit objet, les $\mathcal{J}$-cofibrations sont en fait
des rétracts de colimites séquentielles transfinies de sommes amalgamées
de flèches de $\mathcal{J}$. Il suffit donc de montrer que les cofibrations
triviales sont stables par somme amalgamée et par colimite séquentielle
transfinie pour montrer que les $\mathcal{J}$-cofibrations sont des cofibrations
triviales, la stabilité par rétract ayant déjà été
montrée.\\

Les sections suivantes auront donc pour but de montrer cette stabilité. En
général cette stabilité des cofibrations triviales est difficile à
montrer. Toutefois dans le cas particulier des cofibrations triviales qui
induisent un isomorphisme sur les objets, ce que nous appellerons
iso-cofibrations triviales, cette stabilité est assez simple à
montrer. C'est pourquoi nous allons tout d'abord nous intéresser à ces
iso-cofibrations triviales.

\newpage

\section{Stabilité des iso-cofibrations triviales}

\subsection{Caractérisation des iso-\eqs~de \cats}

Nous cherchons à montrer que les cofibrations triviales sont stables par somme
amalgamée et colimite séquentielle transfinie afin de montrer l'hypothèse
3) du lemme~\ref{reco} de reconnaissance de la structure de \cmfcof. Dans le cas
particulier où ces cofibrations triviales induisent des isomorphismes sur les
objets, cette stabilité est
plus simple à montrer. En effet pour montrer ce résultat, par
l'intermédiaire de la construction Cat, nous allons nous ramener à le
montrer pour les iso-cofibrations triviales de \cats. Or les iso-cofibrations triviales entre \cats~sont exactement les
morphismes de \cats~qui niveau par niveau sont des cofibrations triviales de
$\mathcal{C}$. Cette constatation permettra de remonter la stabilité des
cofibrations triviales dans la \cmf~$\mathcal{C}$ aux iso-cofibrations
triviales de \cats.\\
\\
Mais commençons d'abord par donner les définitions d'iso-\eq~de \cats~et
d'iso-cofibration triviale de \cats.

\begin{defin}\index{iso-\eq!de \cats}\index{iso-cofibration triviale!de \cats}
Soit $\mathcal{C}$ une donnée de Segal.\\
Une \eq~de \cats~qui induit un isomorphisme sur les objets sera nommée
iso-\eq~de \cats. 
Une iso-\eq~de \cats~qui est également un monomorphisme de \precats~sera
nommée iso-cofibration triviale de \cats.
\end{defin}

Donnons aussi les définitions d'iso-\eq~faible de \precats~et
d'iso-cofibration triviale de \precats.

\begin{defin}\index{iso-\eq!faible}\index{iso-cofibration triviale!de \precats}
Soit $\mathcal{C}$ une donnée de Segal facile.\\
Une \eq~faible de \precats~qui induit un isomorphisme sur les objets sera nommée
iso-\eq~faible. 
Une iso-\eq~faible qui est également un monomorphisme de \precats~sera
nommée iso-cofibration triviale de \precats.
\end{defin}

Nous allons maintenant montrer que les iso-\eqs~de \cats~sont exactement les
morphismes de \cats~qui niveau par niveau sont des \eqcs~d'\obcs.

\begin{prop}\label{isoeqc}
Soit $\mathcal{C}$ une donnée de Segal
telle que :
\item -la catégorie sous-jacente $\mathcal{C}$ est une \cmf~simpliciale dont tous les objets sont
cofibrants,
\item -les \eqcs~d'\obcs~sont exactement les \eqs~faibles de la
\cmf~$\mathcal{C}$ entre \obcs,\\
\\
Alors les trois propositions
ci-dessous sont équivalentes :
\item i) $A\rightarrow B$ est un morphisme entre \cats~tel que, pour tout entier
$m$, $A_m\rightarrow B_m$ est une \eqc~d'\obcs.
\item ii) $A\rightarrow B$ est un morphisme entre \cats~tel que $A_0\rightarrow B_0$
est un isomorphisme et $A_1\rightarrow B_1$ est une \eqc~d'\obcs.
\item iii) $A\rightarrow B$ est une iso-\eq~de \cats.
\end{prop}
{\it Preuve :}\\
Montrons tout d'abord que iii) implique ii). Soit $A\rightarrow B$ une \eq~de
\cats~telle que $A_0\rightarrow B_0$ soit un isomorphisme. Il suffit de montrer
que $A_1\rightarrow B_1$ est une \eqc~d'\obcs. Comme $A_0\rightarrow B_0$ est un
isomorphisme, on identifiera les objets de $A$ avec leurs images dans $B$. Par
définition de catégorie \discret, le morphisme $A_1\rightarrow B_1$ n'est autre que le
coproduit des morphismes $A_1(a,a')\rightarrow B_1(a,a')$ pour $a$ et $a'$
décrivant $A_0$. Par hypothèse, comme $A\rightarrow B$ une \eq~de
\cats, alors tous les morphismes $A_1(a,a')\rightarrow B_1(a,a')$ sont des
\eqcs~d'\obcs. En outre, d'après les hypothèses, on peut identifier les
\eqcs~d'\obcs~avec les \eqs~faibles de la \cmf~$\mathcal{C}$ entre \obcs. Comme $\mathcal{C}$ est une \cmf~simpliciale dont tous les objets
sont cofibrants, par corollaire~\ref{eqcoprod}, $A_1\rightarrow B_1$ est une
\eqc~d'\obcs~comme coproduit d'\eqcs~d'\obcs.\\
\\
Montrons maintenant que ii) implique i). Soit $f:A\rightarrow B$ un morphisme entre \cats~tel que $A_0\rightarrow B_0$
est un isomorphisme et $A_1\rightarrow B_1$ est une \eqc~d'\obcs. Soit $m\geq 2$
un entier. Considérons le diagramme commutatif suivant :
\begin{diagram}
A_m & \rTo^{f_m} & B_m \\
\dTo & & \dTo \\
A_1\times_{A_0}\ldots\times_{A_0} A_1 & \rTo & B_1\times_{B_0}\ldots\times_{B_0}
B_1 \\
\end{diagram}
où les flèches verticales sont les morphismes de Segal et la flèche du bas
le produit fibré multiple de $f_1$ au dessus de l'isomorphisme $f_0$. Par
stabilité des \eqcs~d'\obcs~par produit fibré au-dessus d'un objet discret,
il vient que la flèche du bas est une \eqc~d'\obcs. En outre $A$ et $B$ sont
des \cats, donc leurs morphismes de Segal sont aussi des \eqcs~d'\obcs. Et par
l'axiome "trois pour deux" dans la \cmf~$\mathcal{C}$, il vient que $f_m$ est
une \eqc~d'\obcs.\\
Enfin le fait que i) implique iii) découle directement du lemme~\ref{eqcniv}.\\
CQFD.\\

On remarque, d'après la démonstration, que pour toute \cmf~\discret, une
\eq~niveau par niveau est une \eq~de \cats~(voir lemme~\ref{eqcniv}). En revanche une iso-\eq~de
\cats~n'est une \eq~niveau par niveau que dans une
\cmf~\discret~simpliciale dont tous les objets sont cofibrants.\\

Comme les monomorphismes de \precats~sont les monomorphismes niveau par niveau,
il découle directement de cette proposition que les iso-cofibrations triviales
de \cats~sont exactement les morphismes de \cats~qui sont niveau par niveau des
monomorphismes et des \eqcs~d'\obcs. Ceci va nous permettre de montrer la
stabilité des iso-cofibrations triviales de \cats~par colimite en regardant
uniquement ce qui se passe niveau par niveau. Cependant les colimites de
\cats~ne sont pas des \cats. Ainsi on obtiendra seulement à l'arrivée que
les colimites d'iso-cofibrations triviales de \cats~sont des
morphismes de \precats~qui seront des cofibrations triviales niveau par niveau.
Le but de la prochaine section est donc de montrer que ces cofibrations
triviales niveau par niveau sont des cofibration triviales de \precats.

\subsection{Equivalence faible niveau par niveau}

Comme les iso-cofibrations triviales de \cats~sont niveau par niveau des
cofibrations triviales de $\mathcal{C}$ et que, dans la \cmf~$\mathcal{C}$, les
cofibrations triviales sont stables par somme amalgamée et colimite
séquentielle transfinie, on va obtenir que les sommes amalgamées et les
colimites séquentielles transfinies d'iso-cofibrations triviales de \cats~sont
niveau par niveau des cofibrations triviales. Or du fait que les colimites de
\cats~ne sont pas en général des \cats, on a donc à faire à des
morphismes de \precats~qui niveau par niveau sont des cofibrations triviales, ce
qui ne rentre pas dans le cadre de la proposition~\ref{eqcniv}. Aussi
allons-nous montrer dans cette section que les morphismes de \precats~qui niveau
par niveau sont des \eqs~faibles de $\mathcal{C}$ sont en fait des \eqs~faibles
de \precats. Pour cela nous allons utiliser la construction Bigcat qui permet de
suivre la trace de ces morphismes lors des étapes de la catégorisation.

\begin{lem}
Soit $(\mathcal{C},\mathcal{F}_1,\mathcal{F}_2)$ une donnée de Segal facile telle que la \cmf~sous-jacente $\mathcal{C}$ est une \cmf~simpliciale.\\
Alors pour tout morphisme de \precats~$f:A\rightarrow C$ tel que pour tout entier
$q\geq 0$ $f_q:A_q\rightarrow C_q$ est une \eq~faible et pour tout entier $m\geq 2$,
$Cat(1,m)(A)\rightarrow Cat(1,m)(C)$ est niveau par niveau une \eq~faible.
\end{lem}
{\it Preuve :}\\
Tout d'abord on remarque que si la source d'une flèche génératrice
s'envoie dans $A$ alors en composant par $f:A\rightarrow C$ elle s'envoie aussi
dans $C$. Ainsi toute somme amalgamée de $A_q$, pour un certain entier $q$,
par un coproduit de flèches génératrices s'envoyant dans $A_q$,
possède un morphisme vers le même type de sommes amalgamées pour
$C_q$. De ce fait il existe un morphisme $\eta_f$ de $B_A$ vers $B_C$ tel que
l'on ait $\eta_f\circ\eta_A=\eta_C\circ f_1$ et $\nu_C\circ\eta_f=f_0\times
f_0\circ\nu_A$. De même il existe un morphisme $\phi_f$ de $P_A$ vers $P_C$
tel que $\phi_f\circ\phi_A=\phi_C\circ f_m$ et
$\psi_C\circ\phi_f=\eta_f\times_{f_0}\ldots\times_{f_0}\eta_f\circ\psi_A$.\\ 
\\
Comme on l'a vu dans la démonstration du lemme~\ref{eqcat1m}, les morphismes
$\eta$ sont des cofibrations triviales. Comme par hypothèse $f_1$ est une
\eq~faible, par "trois pour deux" dans la \cmf~$\mathcal{C}$, il vient que
$\eta_f$ est une \eq~faible. Et donc
son produit fibré $m$-fois au-dessus de l'isomorphisme $f_0$ aussi. Toujours
d'après la démonstration du lemme~\ref{eqcat1m}, les morphismes $\psi$ sont
des \eqcs~faciles, donc des \eqs~faibles et, par "trois pour deux", il vient que $\phi_f$ est aussi une \eq~faible.
Comme $\mathcal{C}$ est une \cmf~simpliciale dont tous les objets sont
cofibrants, il vient par le corollaire~\ref{eqcoprod} que les coproduits de $\eta_f$ et de $\phi_f$ sont des
\eqs~faibles, ainsi que les coproduits de $f_1$ et $f_m$. Encore d'après la démonstration du lemme~\ref{eqcat1m}, les
morphismes $\eta$ et $\phi$ sont des cofibrations donc leurs coproduits aussi.
On a que $Cat(1,m)(f)_q$ est la somme amalgamée dans la catégorie des
morphismes de $f_q$ par un coproduit de $\eta_f$ et de $\phi_f$ au-dessus d'un
même coproduit de $f_1$ et de $f_m$. En outre les coproduits de $\eta_A$ et de
$\phi_A$ (respectivement $\eta_C$ et $\phi_C$) sont des cofibrations donc, par
la proposition~\ref{eqcube}, $Cat(1,m)(f)_q$ est une \eq~faible.\\
CQFD.

\begin{prop}\label{eqfniv}
Soit $(\mathcal{C},\mathcal{F}_1,\mathcal{F}_2)$ une donnée de Segal facile telle que la \cmf~sous-jacente $\mathcal{C}$ est une \cmf\\ simpliciale.\\
Alors pour tout morphisme de \precats~$f:A\rightarrow C$ tel que pour tout entier
$q\geq 0$ $f_q:A_q\rightarrow C_q$ est une \eq~faible,
$A\rightarrow C$ est une \eq~faible dans $\mathcal{C-PC}$.
\end{prop}
{\it Preuve :}\\ 
Par récurrence transfinie, on montre que pour tout ordinal
$\beta<T$ le morphisme $f^{\beta}:A^{\beta}\rightarrow C^{\beta}$ est niveau par
niveau une \eq~faible, où $f^{\beta}$ est la $\beta$-ième étape de
$Bigcat(f)$. Le pas de la récurrence d'un ordinal $\beta<T$ à l'ordinal
$\beta+1$ est exactement le lemme précédent. Il ne reste qu'à comprendre le
pas de la récurrence transfinie pour un ordinal limite $\beta\leq T$.
Supposons donc que pour tout ordinal $\gamma<\beta$, $f^{\gamma}$ est niveau par
niveau une \eq~faible. Comme en outre les morphismes $A\rightarrow A^{\gamma}$
et $C\rightarrow C^{\gamma}$ sont des cofibrations (démonstration de la
proposition~\ref{eqbig}), par la proposition~\ref{eqcolim}, il vient que
$f^{\beta}$ est niveau par niveau une \eq~faible, ce qui montre l'hypothèse de
récurrence pour l'ordinal limite $\beta$.\\
\\
On a donc bien montré par récurrence transfinie que pour tout ordinal
$\beta<T$ le morphisme $f^{\beta}:A^{\beta}\rightarrow C^{\beta}$ est niveau par
niveau une \eq~faible. C'est donc en particulier le cas pour
$f^{T}=Bigcat(f):Bigcat(A)\rightarrow Bigcat(C)$. 
Comme $Bigcat(f)$ est niveau par niveau une \eq~faible entre
\obcs, i.e. niveau par niveau c'est une \eqc~d'\obcs, et que $Bigcat(A)$ et
$Bigcat(C)$ sont des \cats, alors par le lemme~\ref{eqcniv}, il vient que
$Bigcat(f):Bigcat(A)\rightarrow Bigcat(C)$ est une \eq~de \cats, donc par
le lemme~\ref{bigcateq}, $f$ est une \eq~faible.\\
CQFD.\\

Maintenant que nous savons que les \eqs~faibles niveau par niveau sont des
\eqs~faibles de \precats, nous pouvons nous lancer dans les démonstrations de
stabilité des iso-cofibrations triviales.

\subsection{Sommes amalgamées d'iso-cofibrations triviales}

Avant de montrer le lemme de stabilité des iso-cofibrations triviales par
somme amalgamée le long d'un morphisme, on va montrer un petit lemme technique
qui s'avèrera utile pour toutes ces questions de somme amalgamée de
cofibrations triviales.

\begin{lem}\label{pratique}
Soit $(\mathcal{C},\mathcal{F}_1,\mathcal{F}_2)$ une donnée de Segal facile.
Considérons la somme amalgamée $B\coprod_{A}C$ de \precats. Pour montrer que
le morphisme canonique $C\rightarrow B\coprod_{A}C$ est une \eq~faible, il faut
et il suffit de montrer que le morphisme canonique $Cat(C)\rightarrow
Cat(B)\coprod_{Cat(A)}Cat(C)$ est une \eq~faible. 
\end{lem}
{\it Preuve :}\\
C'est l'application directe du lemme~\ref{colimrec} dont les
hypothèses~\ref{colimIinj}, qui sont en fait celles de la
proposition~\ref{Ephi3}, sont vérifiées par la donnée de Segal facile.\\
CQFD.\\

Ce lemme technique nous permet donc de démontrer la stabilité des
iso-cofibrations triviales par somme amalgamée le long d'un morphisme
uniquement en utilisant les iso-cofibrations triviales de \cats. Or si Cat
transporte bien les \eqs~faibles sur les \eqs~de \cats, par définition même
des \eqs~faibles de \precats, en revanche il faut s'assurer que Cat envoie bien
les cofibrations sur les cofibrations, ce qui est l'objet du lemme suivant.

\begin{lem}\label{cofcat}
Soit $(\mathcal{C},\mathcal{F}_1,\mathcal{F}_2)$ une donnée de Segal facile.
Alors les cofibrations sont préservées par la catégorisation Cat.
\end{lem}
{\it Preuve :}\\
Tout d'abord on remarque que, la catégorisation Cat n'étant autre que le plan
d'addition de cellules $E_{\Phi}$, ses plans simples sont du type $e_{\Psi,\lambda}$, avec pour $\Phi$ la famille $\mathcal{FG}_1$ des flèches génératrices de \cats~faciles.
Ensuite on rappelle que, comme la donnée de Segal est facile, alors par le lemme~\ref{flgenmono}, les $\mathcal{FG}_1$-cofibrations sont des monomorphismes niveau par niveau donc
des monomorphismes de \precats. En prenant donc pour ensemble $\Phi$ la famille
de monomorphismes $\mathcal{FG}_1$,
il vient d'après le lemme~\ref{monoadd} que les monomorphismes sont
préservés par la catégorisation Cat.\\
CQFD.

\begin{lem}\label{ictpush}
Soit $(\mathcal{C},\mathcal{F}_1,\mathcal{F}_2)$ une donnée de Segal facile
telle que la \cmf~sous-jacente $\mathcal{C}$ est une \cmf~simpliciale engendrée
par monomorphismes.\\
Soit $A\rightarrow B$ un morphisme de \precats~à la fois cofibration et
iso-\eq~faible. Alors sa somme amalgamée le long d'un morphisme $A\rightarrow
C$ quelconque est aussi une cofibration et une iso-\eq~faible.
\end{lem}
{\it Preuve :}\\
Tout d'abord, on remarque que, comme les isomorphismes sur
les objets sont stables par somme amalgamée le long d'un morphisme, le
morphisme canonique $C\rightarrow B\coprod_{A} C$ est un
isomorphisme sur les objets. En outre les cofibrations sont des cofibrations
niveau par niveau et, dans la \cmf~$\mathcal{C}$, les cofibrations sont stables
par somme amalgamée le long d'un morphisme. Donc $C\rightarrow B\coprod_{A} C$ est aussi une 
cofibration comme somme amalgamée niveau par niveau le long d'un morphisme
d'une cofibration niveau par niveau. Il reste donc à montrer que c'est une \eq~faible,
ce qui d'après le lemme~\ref{pratique} revient à montrer que $Cat(C)\rightarrow
Cat(B)\coprod_{Cat(A)}Cat(C)$ est une \eq~faible. Or ce morphisme est
lui-même la somme amalgamée de $Cat(A)\rightarrow Cat(B)$ le long de
$Cat(A)\rightarrow Cat(C)$.\\
\\
Comme par hypothèse les cofibrations sont les
monomorphismes,
par le lemme \ref{cofcat}, il vient que $Cat(A)\rightarrow Cat(B)$ est une
cofibration. En outre $Cat$ préserve les objets et donc $Cat(A)\rightarrow
Cat(B)$ est un isomorphisme sur les objets ainsi que sa somme amalgamée $Cat(C)\rightarrow
Cat(B)\coprod_{Cat(A)}Cat(C)$. Enfin comme par hypothèse,
$A\rightarrow B$ est une \eq~faible alors $Cat(A)\rightarrow Cat(B)$ est une
\eq~de \cats. Finalement on obtient que $Cat(A)\rightarrow Cat(B)$ est une
cofibration et une iso-équivalence, et donc par la proposition~\ref{isoeqc}, c'est
une cofibration triviale niveau par niveau. Or la somme amalgamée dans
$\mathcal{C-PC}$ est la somme amalgamée niveau par niveau. Comme
$\mathcal{C}$ est une \cmf, les cofibrations triviales de $\mathcal{C}$ sont
stables par somme amalgamées le long d'un morphisme. Ainsi, on obtient que $Cat(C)\rightarrow
Cat(B)\coprod_{Cat(A)}Cat(C)$ est une cofibration triviale niveau par niveau,
comme somme amalgamée niveau par niveau d'une cofibration triviale niveau par
niveau. Comme c'est aussi un isomorphisme sur les objets, par la proposition~\ref{eqfniv}, il vient que $Cat(C)\rightarrow
Cat(B)\coprod_{Cat(A)}Cat(C)$ est donc une
\eq~faible, ce que l'on voulait démontrer.\\
CQFD.

\subsection{Colimites séquentielles d'iso-cofibrations triviales}

Après avoir montré la stabilité des iso-cofibrations triviales par somme
amalgamée le long d'un morphisme, nous allons regarder la stabilité des
iso-cofibrations triviales par colimite séquentielle transfinie. Mais comme
précédem\-ment, nous allons tout d'abord donner un lemme pratique pour ces
histoires de colimites séquentielles transfinies de cofibrations triviales.

\begin{lem}\label{pratique2}
Soit $(\mathcal{C},\mathcal{F}_1,\mathcal{F}_2)$ une donnée de Segal facile.
Considérons une $\lambda$-séquence $(A_{\alpha})_{\alpha<\lambda}$ de \precats. Pour montrer que
le morphisme canonique $A_0\rightarrow \colimite{\alpha<\lambda}A_{\alpha}$ est une \eq~faible, il faut
et il suffit de montrer que le morphisme canonique $Cat(A_0)\rightarrow
\colimite{\alpha<\lambda} Cat(A_{\alpha})$ est une \eq~faible. 
\end{lem}
{\it Preuve :}\\
Il s'agit à nouveau de l'application directe du lemme~\ref{colimrec} dont les
hypothèses~\ref{colimIinj} sont vérifiées par la donnée de Segal facile.\\
CQFD.

\begin{lem}\label{ictcolim1}
Soit $(\mathcal{C},\mathcal{F}_1,\mathcal{F}_2)$ une donnée de Segal facile
telle que la \cmf~sous-jacente $\mathcal{C}$ est une \cmf~simpliciale engendrée
par monomorphismes.\\
Soit $(A_{\alpha})_{\alpha<\lambda}$ une $\lambda$-séquence de \precats~telle
que pour tout $\alpha<\lambda$ le morphisme $A_{\alpha}\rightarrow A_{\alpha+1}$
est à la fois une cofibration et une iso-\eq~faible. Alors le morphisme
canonique $A_0\rightarrow \colimite{\alpha<\lambda}A_{\alpha}$ est aussi une cofibration et une iso-\eq~faible.
\end{lem}
{\it Preuve :}\\
Tout d'abord, on remarque que, les isomorphismes sur
les objets étant stables par colimite séquentielle transfinie, le
morphisme $A_0\rightarrow \colimite{\alpha<\lambda}A_{\alpha}$ est un
isomorphisme sur les objets. En outre les cofibrations ne sont autres que les
cofibrations niveau par niveau donc elles sont aussi stables par colimite
séquentielle transfinie car cette colimite se calcule niveau par niveau et
que, $\mathcal{C}$ étant une \cmf, les cofibrations de $\mathcal{C}$ sont stables
par colimite séquentielle transfinie. Donc $A_0\rightarrow
\colimite{\alpha<\lambda}A_{\alpha}$ est aussi une cofibration. Il reste donc à
montrer que c'est une \eq~faible, ce qui d'après le lemme précédent
revient à montrer que le morphisme $Cat(A_0)\rightarrow
\colimite{\alpha<\lambda} Cat(A_{\alpha})$ est une \eq~faible. Or ce morphisme est
lui-même la colimite séquentielle transfinie des $Cat(A_{\alpha})\rightarrow
Cat(A_{\alpha+1})$.\\
\\
Comme par hypothèse les cofibrations sont les
monomorphismes,
par le lemme \ref{cofcat}, il vient que $Cat(A_{\alpha})\rightarrow
Cat(A_{\alpha+1})$ est une cofibration. En outre $Cat$ préserve les objets et
donc $Cat(A_{\alpha})\rightarrow Cat(A_{\alpha+1})$ est un isomorphisme sur les
objets, donc, par colimite séquentielle transfinie, $Cat(A_0)\rightarrow
\colimite{\alpha<\lambda} Cat(A_{\alpha})$ aussi. Enfin comme par hypothèse,
$A_{\alpha}\rightarrow A_{\alpha+1}$ est une \eq~faible alors
$Cat(A_{\alpha})\rightarrow Cat(A_{\alpha+1})$ est une
\eq~de \cats. Finalement on obtient que $Cat(A_{\alpha})\rightarrow
Cat(A_{\alpha+1})$ est une
cofibration et une iso-équivalence, et donc par la proposition~\ref{isoeqc}, c'est
une cofibration triviale niveau par niveau. Or les colimites séquentielles
transfinies dans $\mathcal{C-PC}$ sont les colimites séquentielles transfinies niveau par niveau. Comme
$\mathcal{C}$ est une \cmf, les cofibrations triviales de $\mathcal{C}$ sont
stables par colimite séquentielle transfinie. Ainsi, on obtient que $Cat(A_0)\rightarrow
\colimite{\alpha<\lambda} Cat(A_{\alpha})$ est une cofibration triviale niveau par niveau,
comme colimite séquentielle transfinie niveau par niveau de cofibrations
triviales niveau par niveau. Comme c'est aussi un isomorphisme sur les objets, par la proposition~\ref{eqfniv}, il vient que $Cat(A_0)\rightarrow
\colimite{\alpha<\lambda} Cat(A_{\alpha})$ est donc une \eq~faible, ce que l'on voulait démontrer.\\
CQFD.\\

Pour clore cette section, nous allons aussi donner un résultat de stabilité
des iso-cofibrations triviales par colimite séquentielle transfinie dans la
catégorie des morphismes, en commençant comme toujours par un petit lemme
pratique.

\begin{lem}
Soit $(\mathcal{C},\mathcal{F}_1,\mathcal{F}_2)$ une donnée de Segal facile.
Considérons un morphisme de $\lambda$-séquences
$(A_{\alpha})_{\alpha<\lambda}\rightarrow (B_{\alpha})_{\alpha<\lambda}$. Pour montrer que
le morphisme $\colimite{\alpha<\lambda}A_{\alpha}\rightarrow
\colimite{\alpha<\lambda}B_{\alpha}$ est une \eq~faible, il faut
et il suffit de montrer que le morphisme\\ $\colimite{\alpha<\lambda} Cat(A_{\alpha})\rightarrow
\colimite{\alpha<\lambda} Cat(B_{\alpha})$ est une \eq~faible. 
\end{lem}
{\it Preuve :}\\
Il s'agit encore une fois de l'application directe du lemme~\ref{colimrec} dont
la donnée de Segal facile vérifie les hypothèses.\\
CQFD.

\begin{lem}\label{ictcolim2}
Soit $(\mathcal{C},\mathcal{F}_1,\mathcal{F}_2)$ une donnée de Segal facile
vérifiant les propriétés suivantes :
\item 1) la \cmf~sous-jacente $\mathcal{C}$ est une \cmf~simpliciale engendrée
par monomorphismes,
\item 2) les colimites séquentielles transfinies de la sous-catégorie des
monomorphismes de $\mathcal{C}$ existent et sont les colimites séquentielles
transfinies de $\mathcal{C}$.\\  
\\
Soit $(A_{\alpha})_{\alpha<\lambda}\rightarrow (B_{\alpha})_{\alpha<\lambda}$ un
morphisme de $\lambda$-séquences de \precats~tel
que pour tout $\alpha<\lambda$ les morphismes $A_{\alpha}\rightarrow
A_{\alpha+1}$ et $B_{\alpha}\rightarrow B_{\alpha+1}$ sont des cofibrations et
le morphisme $A_{\alpha}\rightarrow B_{\alpha}$ est à la fois une cofibration
et une iso-\eq~faible. Alors le morphisme
$\colimite{\alpha<\lambda}A_{\alpha}\rightarrow
\colimite{\alpha<\lambda}B_{\alpha}$ est aussi une cofibration et une iso-\eq~faible.
\end{lem}
{\it Preuve :}\\
Tout d'abord, on remarque que, les isomorphismes sur
les objets étant stables par colimite séquentielle transfinie dans la
catégorie des morphismes, le
morphisme $\colimite{\alpha<\lambda}A_{\alpha}\rightarrow
\colimite{\alpha<\lambda}B_{\alpha}$ est un
isomorphisme sur les objets. En outre les cofibrations ne sont autres que les
monomorphismes niveau par niveau donc elles sont aussi stables par colimite
séquentielle transfinie dans la catégorie des cofibrations car cette colimite se calcule niveau par niveau et
que, par hypothèse sur $\mathcal{C}$, les monomorphismes de $\mathcal{C}$ sont stables
par colimite séquentielle transfinie dans la catégorie des monomorphismes. Donc $\colimite{\alpha<\lambda}A_{\alpha}\rightarrow
\colimite{\alpha<\lambda}B_{\alpha}$ est aussi une cofibration. Il reste donc à
montrer que c'est une \eq~faible, ce qui d'après le lemme précédent
revient à montrer que le morphisme $\colimite{\alpha<\lambda} Cat(A_{\alpha})\rightarrow
\colimite{\alpha<\lambda} Cat(B_{\alpha})$ est une \eq~faible. Or ce morphisme est
lui-même la colimite séquentielle transfinie dans la catégorie des
morphismes des $Cat(A_{\alpha})\rightarrow
Cat(B_{\alpha})$. Comme par hypothèse les cofibrations sont les
monomorphismes,
par le lemme~\ref{cofcat}, il vient que les morphismes
$Cat(A_{\alpha})\rightarrow Cat(A_{\alpha+1})$, $Cat(B_{\alpha})\rightarrow
Cat(B_{\alpha+1})$ et $Cat(A_{\alpha})\rightarrow
Cat(B_{\alpha})$ sont des cofibrations.\\
Ainsi le morphisme $\colimite{\alpha<\lambda} Cat(A_{\alpha})\rightarrow
\colimite{\alpha<\lambda} Cat(B_{\alpha})$ est la colimite séquentielle transfinie dans la catégorie des
cofibrations des $Cat(A_{\alpha})\rightarrow
Cat(B_{\alpha})$. En outre $Cat$ préserve les objets et
donc $Cat(A_{\alpha})\rightarrow Cat(B_{\alpha})$ est un isomorphisme sur les
objets, et, par colimite séquentielle transfinie dans la catégorie des
morphismes, $\colimite{\alpha<\lambda} Cat(A_{\alpha})\rightarrow
\colimite{\alpha<\lambda} Cat(B_{\alpha})$ aussi. Enfin comme par hypothèse,
$A_{\alpha}\rightarrow B_{\alpha}$ est une \eq~faible alors
$Cat(A_{\alpha})\rightarrow Cat(B_{\alpha})$ est une
\eq~de \cats. Finalement on obtient que $Cat(A_{\alpha})\rightarrow
Cat(B_{\alpha})$ est une
cofibration et une iso-équivalence, et donc par la proposition~\ref{isoeqc}, c'est
une cofibration triviale niveau par niveau. Or les colimites séquentielles
transfinies dans la catégorie des cofibrations de $\mathcal{C-PC}$ sont les
colimites séquentielles transfinies niveau par niveau dans la catégorie des
cofibrations de $\mathcal{C}$. Comme
$\mathcal{C}$ est une \cmf~simpliciale propre à gauche, par la
proposition~\ref{eqcolim}, les cofibrations triviales de $\mathcal{C}$ sont
stables par colimite séquentielle transfinie dans la catégorie des
cofibrations. Ainsi, on obtient que $\colimite{\alpha<\lambda} Cat(A_{\alpha})\rightarrow
\colimite{\alpha<\lambda} Cat(B_{\alpha})$ est une cofibration triviale niveau par niveau,
comme colimite séquentielle transfinie niveau par niveau dans la catégorie
des cofibrations niveau par niveau de cofibrations
triviales niveau par niveau. Comme c'est aussi un isomorphisme sur les objets, par la proposition~\ref{eqfniv}, il vient que $\colimite{\alpha<\lambda} Cat(A_{\alpha})\rightarrow
\colimite{\alpha<\lambda} Cat(B_{\alpha})$ est donc une \eq~faible, ce que l'on voulait démontrer.\\
CQFD.\\

Jusqu'à présent, la stabilité par somme amalgamée le long d'un morphisme
et par colimite séquentielle transfinie des cofibrations triviales n'a été
démontrée que pour les iso-cofibrations triviales. Afin de la montrer pour
toutes les cofibrations triviales, nous avons besoin d'un gros résultat sur le
produit cartésien des \cats, ce qui sera l'objet de la prochaine section.

\newpage

\section{Produit cartésien de \cats}

Le but de cette section est de montrer que le produit cartésien de deux
\precats~est équivalent au produit cartésien de leurs catégorisations. De
ce résultat va découler la stabilité des \eqs~faibles de \precats~par produit
cartésien qui permettra de terminer la démonstration de stabilité des
\eqs~faibles de \precats~par somme amalgamée le long d'un morphisme et par
colimite séquentielle transfinie. Lorsque l'on démontre que le produit cartésien de deux
\precats~est équivalent au produit cartésien de leurs catégorisations, on
se ramène à montrer un résultat particulier d'\eq~entre des \precats. Ce
résultat fait intervenir la notion de \precat~librement ordonnées.

\subsection{\precats~librement ordonnées}

L'idée principale d'une \precat~librement ordonnée est que son ensemble
d'objet est muni d'un ordre total pour lequel il n'existe de morphisme entre
deux objets que s'il existe une chaîne ordonnée d'objets adjacents
reliés par un morphisme. Donnons en une définition
plus précise.

\begin{defin}\index{\precat!librement ordonnée}\index{\precat!librement
ordonnée stricte}
Soit $\mathcal{C}$ une \cmf~\discret.
Une \precat~$A$ est librement ordonnée si son ensemble d'objets $A_0$ est fini
et muni d'un ordre total vérifiant les trois propriétés suivantes :
\item i) pour toute suite d'objets $(a_0,\ldots,a_m)$ non ordonnée dans $A$,
$A_m(a_0,\ldots,a_m)$ est l'ensemble vide,
\item ii) pour toute suite ordonnée d'objets $(a_0,\ldots,a_m)$ dans $A$, le
morphisme\\
$A_m(a_0,\ldots,a_m)\rightarrow A_1(a_0,a_m)$, induit par l'application de {\bf
1} dans {\bf m} qui à 0 et 1 associe 0 et m, est une \eq~faible dans la
\cmf~$\mathcal{C}$,
\item iii) pour tout objet $a$ de $A$, $A_1(a,a)$ est un point.\\
\\
On dira que $A$ est librement ordonnée stricte, si on demande pour la
propriété ii) un isomorphisme au lieu d'une équivalence faible.
\end{defin}

Voici quelques exemples importants de \precats~librement ordonnées.

\begin{lem}\label{FGordre}
Soit $\mathcal{C}$ une catégorie \discret~possédant les coproduits et les
sommes amalgamées. 
\item -Pour tout objet $X$ de
$\mathcal{C}$ et pour tout entier $m$, les \precats~$\Delta[m]\Theta X$ et
$\Upsilon[m]\Theta X$ sont librement ordonnées strictes.
\item -Pour tout morphisme $f$ de $\mathcal{C}$ et pour tout entier
$m\geq 2$, la source du morphisme $Boit_m(f)$ est une \precats~librement
ordonnée stricte et les morphismes $\Delta[m]\Theta f$ et $Boit_m(f)$ sont des
isomorphismes sur les objets préser\-vant strictement l'ordre.
\end{lem}
{\it Preuve:}\\
$\Delta[m]\Theta X$ a pour ensemble d'objets
$\operatorname{Hom}_{\Delta}(\mbox{\bf 0,m)}$, qui est un ensemble fini à
$m+1$ éléments. En identifiant les applications
de {\bf 0} dans {\bf m} avec leur image, on obtient ainsi un ordre total naturel
sur $\Delta[m]\Theta X_0$. On obtient le même ordre sur les objets de
$\Upsilon[m]\Theta X$ et par suite sur ceux de la source de $Boit_m(f)$.\\ 
Soient $q>0$ un entier et $(x_0,\dots,x_q)$ une suite 
d'objets de $\Delta[m]\Theta X_0$. 
On remarque que la suite $(x_0,\dots,x_q)$
s'identifie à une application de {\bf q}
dans {\bf m} si et seulement si cette suite est ordonnée et qu'elle
s'identifie à une application de {\bf q} dans {\bf m} se factorisant par {\bf
0} si et seulement si tous les $x_i$ sont identiques. 
En utilisant cette
remarque ainsi que la définition de $\Delta[m]\Theta X$, il vient que
$\Delta[m]\Theta X_q(x_0,\dots,x_q)$ est l'ensemble vide si la suite n'est pas
ordonnée, un point si elle est constante et $X$ si la suite est ordonnée non
constante.
Dans ces deux derniers cas, le morphisme $\Delta[m]\Theta X_q(x_0,\dots,x_q)\rightarrow
\Delta[m]\Theta X_1(x_0,x_q)$ n'est autre que l'identité du point,
respectivement de $X$. Ceci montre
que $\Delta[m]\Theta X$ est une \precat~librement ordonnée stricte.\\
\\
Pour $\Upsilon[m]\Theta X$, on remarque qu'une suite $(x_0,\dots,x_q)$
d'objets de $\Upsilon[m]\Theta X_0$ s'identifie à une application de {\bf q}
dans {\bf m} se factorisant par une face principale si et seulement si la suite
est ordonnée et que l'on ait : $x_q=x_0+1$. Avec cette remarque et la
définition de $\Upsilon[m]\Theta X$, on obtient que $\Upsilon[m]\Theta
X_q(x_0,\dots,x_q)$ est le vide si la suite n'est pas ordonnée, un point si la
suite est constante, $X$ si la suite est ordonnée à extrémités
adjacentes et le vide sinon. Dans ces trois derniers cas, le morphisme $\Upsilon[m]\Theta X_q(x_0,\dots,x_q)\rightarrow
\Upsilon[m]\Theta X_1(x_0,x_q)$ n'est autre que l'identité du point,
respectivement de $X$ et de l'ensemble vide. Ceci montre que $\Upsilon[m]\Theta
X$ est une \precat~librement ordonnée stricte.\\
\\
Enfin avec les remarques ci-dessus et la définition de $Boit_m(f)$, on obtient
que $B(m,f)_q(x_0,\ldots,x_q)$ est le vide si la suite n'est pas ordonnée, un point si la
suite est constante, $Y$ si la suite est ordonnée à extrémités
adjacentes et $X$ sinon. Dans ces trois derniers cas, le morphisme $B(m,f)_q(x_0,\dots,x_q)\rightarrow
B(m,f)_1(x_0,x_q)$ n'est autre que l'identité du point,
respectivement de $Y$ et de $X$. Ceci montre que la source de $Boit_m(f)$ est
une \precat~librement ordonnée stricte.\\
Pour finir, on termine en constatant que $\Delta[m]\Theta f$ et $Boit_m(f)$ induisent par définition
l'identité sur les objets.\\
CQFD.\\

Comme il n'est pas toujours facile de maîtriser les \precats~librement
ordonnées et encore moins de vérifier qu'un morphisme entre elles préserve
l'ordre, nous allons donner une condition suffisante simple pour qu'un morphisme
entre \precats~librement ordonnées préserve l'ordre.

\begin{lem}\label{mpo}
Soit $\mathcal{C}$ une \cmf~\discret~tel qu'il n'existe aucun morphisme de but
l'ensemble vide qui ne soit pas un isomorphisme.
Soit $f:A\rightarrow B$ un morphisme entre \precats~librement ordonnées. Si
$A$ a la propriété que, pour tout couple ordonné d'objets adjacents $(a,a')$ dans
$A$, $A_1(a,a')$ n'est pas isomorphe à l'ensemble vide, alors $f$ est un morphisme
préservant l'ordre des objets.
\end{lem}
{\it Preuve :}\\
Par l'absurde, si $f$ ne préserve pas l'ordre des objets, il
existe un couple ordonné d'objets adjacents $(a,a')$ tel que le couple
$(f(a),f(a'))$ ne soit pas ordonné. Alors, comme $B$ est librement ordonné,
on a que $B_1(f(a),f(a'))$ est l'ensemble vide, tandis que par hypothèse
$A_1(a,a')$ n'est pas isomorphe à l'ensemble vide. Ainsi $f_1(a,a'):A_1(a,a')\rightarrow
B_1(f(a),f(a'))$ envoie un objet non isomorphe à l'ensemble vide dans
l'ensemble vide, ce qui contredit l'hypothèse sur $\mathcal{C}$.\\
CQFD.\\

La notion de \precats~librement ordonnées est stable par certaines
opérations. Nous allons dans le lemme suivant le montrer pour les opérations
dont nous nous servirons par la suite.

\begin{lem}\label{stordre}
Soit $\mathcal{C}$ une \cmf~\discret~simpliciale dont tous les objets sont
cofibrants et dont les cofibrations et les \eqs~faibles sont stables par
produit fibré au-dessus d'un objet discret.
\item 1) Si $A$ et $B$ est des \precats~librement ordonnées alors $A\times B$
aussi et si elles sont strictes alors $A\times B$ également.
\item 2) Soient $f:A\rightarrow B$ et $g:A\rightarrow C$ deux morphismes entre
\precats~librement ordonnées préservant l'ordre des objets. Si $f$ préserve
strictement l'ordre des objets et que $g$ est une cofibration induisant un
isomorphisme sur les objets, alors la
somme amalgamée $B\coprod_A C$ est une \precat~librement ordonnée. Si $A$,
$B$ et $C$ sont strictes, leur somme amalgamée aussi.
\item 3) Soient $\beta$ un ordinal quelconque et $A$ la colimite d'une
$\beta$-séquence de \precats~librement ordonnées
$(A_{\alpha})_{\alpha<\beta}$ telle que, pour tout $\alpha<\beta$, les morphismes
$A_{\alpha}\rightarrow A_{\alpha+1}$ soient des cofibrations niveau par niveau
induisant un isomorphisme sur les ensembles d'objets préservant leur ordre.
Alors $A$ est une \precat~librement ordonnée.
\end{lem}
{\it Preuve :}\\
1) Tout d'abord on remarque que le produit des \precats~n'est autre que le
produit niveau par niveau. Ainsi l'ensemble des objets de $A\times B$ est fini
comme produit d'ensembles finis et on lui donnera comme ordre total l'ordre
lexicographique : $$(a,b)\leq (a',b') \mbox{ si } a<a' \mbox{ ou } (a=a' \mbox{
et } b\leq b').$$
Soient $m>0$ un entier et $((a_0,b_0),\ldots,(a_m,b_m))$ une suite de points de
$A\times B$. On a l'égalité suivante :
$$\big(A\times B\big)_m\big((a_0,b_0),\ldots,(a_m,b_m)\big)=A_m(a_0,\dots,a_m)\times
B_m(b_0,\ldots,b_m).$$
Si la suite n'est pas ordonnée, alors il existe un indice $i$ entre 0 et $m-1$
tel que $(a_i,a_{i+1})$ n'est pas ordonné ou bien $a_i$ et $a_{i+1}$ sont
égaux et $(b_i,b_{i+1})$ n'est pas ordonné. Comme $A$ et $B$ sont librement
ordonnés, l'un des deux termes du produit est l'ensemble vide, ce qui entraîne
que le produit aussi est l'ensemble vide, ce qui prouve i) pour $A\times B$.\\
\\
Si la suite est ordonnée, comme $A$ et $B$ sont librement
ordonnés, les morphismes $A_m(a_0,\dots,a_m)\rightarrow A_1(a_0,a_m)$ et
$B_m(b_0,\ldots,b_m)\rightarrow B_1(b_0,b_m)$ sont des \eqs~faibles (respectivement des
isomorphismes) donc leur produit aussi par stabilité des \eqs~faibles
(respectivement des isomorphismes) par produit, ce qui prouve ii)
(respectivement ii)-strict) pour $A\times B$.\\
Si la suite est constante, comme $A$ et $B$ sont librement
ordonnés, $A_m(a_0,\dots,a_m)$ et $B_m(b_0,\ldots,b_m)$ sont des points et donc leur produit aussi, ce qui montre
iii) pour $A\times B$. Donc $A\times B$ est librement ordonné (respectivement
librement ordonné strict) si $A$ et $B$ le sont.\\
\\
2) Tout d'abord on remarque que la somme amalgamée des \precats~n'est autre
que la somme amalgamée niveau par niveau. Comme $g_0$ est un isomorphisme préservant l'ordre, on identifiera par la
suite les objets de $A$ avec leurs images dans $C$. En outre, l'isomorphisme
$g_0$ entraîne un isomorphisme entre $B_0$ et $B\coprod_{A} C_0$, ce qui, d'une
part, montre que $B\coprod_{A} C_0$ est un ensemble fini et, d'autre part, munit
ainsi $B\coprod_{A} C_0$ de l'ordre total induit par celui de $B_0$. Encore une
fois, on identifiera les objets de $B_0$ et de $B\coprod_{A} C_0$ via cet
isomorphisme. Soient $m>0$ un entier et $(b_0,\ldots,b_m)$ une suite d'objets de
$B\coprod_{A} C$. Si $b_0,\ldots,b_m$ sont dans l'image de $f_0$, on a
l'égalité suivante :
$$\Big(B\coprod_{A}
C\Big)_m(b_0,\ldots,b_m)=B_m(b_0,\ldots,b_m)\coprod_{\big(\coprod_{a\in f^{-1}(b)}
A_m(a)\big)} \Bigg(\coprod_{a\in f^{-1}(b)}
C_m(a)\Bigg),$$
où $a$ est une notation pour une suite d'objets $(a_0,\ldots,a_m)$ de $A$ et
où $a\in f^{-1}(b)$ signifie que les images de $a_0,\ldots,a_m$ par $f$ sont
respectivement $b_0,\ldots,b_m$.
Sinon on a : $\big(B\coprod_{A} C\big)_m(b_0,\ldots,b_m)=B_m(b_0,\ldots,b_m)$,
auquel cas i), ii) et iii) sont vérifiées car $B$ les vérifie.
Il reste donc à montrer que les trois propriétés sont vraies pour
$b_0,\ldots,b_m$ dans l'image de $f_0$.\\
\\
Tout d'abord on remarque que le
fait que $f$ préserve strictement l'ordre entraîne que $f$ est injective et
donc pour toute suite $(b_0,\ldots,b_m)$ dans l'image de $f_0$, il existe une
unique suite $(a_0,\ldots,a_m)$ telle que son image par $f$ soit
$(b_0,\ldots,b_m)$. Ainsi on obtient l'égalité suivante :
$$\Big(B\coprod_{A}
C\Big)_m(b_0,\ldots,b_m)=B_m(b_0,\ldots,b_m)\coprod_{A_m(a_0,\ldots,a_m)}
C_m(a_0,\ldots,a_m).$$\\
Si la suite n'est pas ordonnée, comme $f$ préserve l'ordre, la suite 
$(a_0,\ldots,a_m)$ ayant pour image $(b_0,\ldots,b_m)$ par $f$ ne sera pas non plus
ordonnée. Comme $A$, $B$ et $C$ sont librement ordonnées, chaque terme de la
somme amalgamée est vide, donc la somme amalgamée aussi, ce qui montre i).
Si la suite est constante, comme $f$ préserve strictement l'ordre, la suite
$(a_0,\ldots,a_m)$ ayant pour image $(b_0,\ldots,b_m)$ par $f$ l'est aussi.
Comme $A$, $B$ et $C$ sont librement ordonnées, chaque terme de la
somme amalgamée est un point, donc la somme amalgamée aussi, ce qui montre
iii).\\
\\
Si la suite est ordonnée mais non constante, comme $f$ préserve l'ordre, la suite 
$(a_0,\ldots,a_m)$ ayant pour image $(b_0,\ldots,b_m)$ par $f$ le sera aussi.
Comme $A$, $B$ et $C$ sont librement ordonnées (respectivement librement
ordonnées strictes), les morphismes
$A_m(a_0,\ldots,a_m)\rightarrow A_1(a_0,a_m)$, $B_m(b_0,\ldots,b_m)\rightarrow
B_1(b_0,b_m)$ et $C_m(a_0,\ldots,a_m)\rightarrow C_1(a_0,a_m)$
sont des \eqs~faibles (respectivement des isomorphismes). Comme par hypothèse $g$ est une cofibration niveau par
niveau donc $g_1$ et $g_m$ sont des cofibrations de la \cmf~$\mathcal{C}$. Comme
en outre, par hypothèse, les cofibrations sont stables par produit fibré
au-dessus d'un objet discret, $g_1(a_0,a_m)$ et $g_m(a_0,\ldots,a_m)$ sont des
cofibrations. Enfin, comme on a supposé que tous les objets de $\mathcal{C}$
sont cofibrants, en appliquant la proposition~\ref{eqcube}, on obtient que
$B\coprod_{A}C_m(b_0,\ldots,b_m)\rightarrow
B\coprod_{A}C_1(b_0,b_m)$
est une \eq~faible, ce qui montre ii). Le cas strict est plus facile car la
somme amalgamée d'isomorphismes dans la catégorie des morphismes est un
isomorphisme.\\
\\
3) Tout d'abord on remarque que la colimite séquentielle transfinie des \precats~n'est autre
que la colimite séquentielle transfinie niveau par niveau. 
Comme pour tout ordinal $\alpha<\beta$, le morphisme $A_{\alpha}\rightarrow
A_{\alpha+1}$ induit un isomorphisme au niveau des objets préservant l'ordre, 
on identifiera donc via ces isomorphismes les ensembles d'objets des
$A_{\alpha}$. Ainsi l'ensemble des objets de $A$ sera donc isomorphe à
ceux-ci, ce qui permet de montrer qu'il est fini et muni d'un ordre total. Soit
$m>0$ un entier et $(a_0,\ldots,a_m)$ une suite d'objets de $A$. On a
l'égalité suivante : $$A_m(a_0,\ldots,a_m)=\colim{\alpha<\beta}
(A_{\alpha})_m(a_0,\ldots,a_m). $$
Si la suite n'est pas ordonnée, comme les $A_{\alpha}$ sont librement
ordonnés, chaque terme de la colimite est l'ensemble vide et donc la colimite
aussi, ce qui montre i) pour $A$. Si la suite est constante, pour la même
raison, chaque terme de la colimite est le point et donc la colimite aussi, ce
qui montre iii) pour $A$. Enfin si la suite est ordonnée mais non constante,
alors comme les $A_{\alpha}$ sont librement ordonnés, pour tout
$\alpha<\beta$, les morphismes
$(A_{\alpha})_m(a_0,\ldots,a_m)\rightarrow (A_{\alpha})_1(a_0,a_m)$ sont des
\eqs~faibles. Or par hypothèse, pour tout $\alpha<\beta$, les morphismes $A_{\alpha}\rightarrow
A_{\alpha+1}$ sont des cofibrations niveau par niveau donc les morphismes
$(A_{\alpha})_1(a_0,a_m)\rightarrow
(A_{\alpha+1})_1(a_0,a_m)$ et
$(A_{\alpha})_m(a_0,\ldots,a_m)\rightarrow$\\
$(A_{\alpha+1})_m(a_0,\ldots,a_m)$ sont des cofibrations, car les cofibrations sont stables par produit fibré
au-dessus d'un objet discret. Comme, par hypothèse,
$\mathcal{C}$ est une \cmf~simpliciale dont tous les objets sont cofibrants, par
la proposition~\ref{eqcolim}, on obtient que $A_m(a_0,\ldots,a_m)\rightarrow A_1(a_0,a_m)$ est une \eq~faible, ce qui
montre ii) pour $A$.\\ 
CQFD.\\

On remarque dans la démonstration que pour qu'une somme amalgamée de
\precats~librement ordonnées, le long d'une cofibration induisant un
isomorphisme sur les objets, vérifie ii) et iii), il faut que le second
morphisme préserve l'ordre strictement. Or lors des constructions Cat ou
Bigcat, les sommes amalgamées par les flèches génératrices, qui sont des
cofibrations induisant des isomorphismes sur les objets, n'ont pas cette
propriété. De ce fait nos deux procédés de catégorisation ne
préservent pas la notion de \precat~librement ordonnée. Dans la section
suivante, nous allons construire une autre catégorisation, dite réduite, qui
préservera les \precats~librement ordonnées. Pour cela, nous regarderons de
plus près les marquages des \cats~essentiellement au niveau des
dégénérescences.

\subsection{Catégorisation réduite}

Lorsqu'on a défini le marquage des \precats, on a fait des choix
arbitraires parmi les relèvements possibles d'un même
diagramme. Mais on n'a jamais demandé que ces choix soient compatibles entre
eux. C'est ce que l'on va demander à un marquage dit "canonique". La seule
compatibilité entre les diagrammes qui pose problème dans la catégorisation d'une \cat~librement ordonnée est celle due aux
applications de dégénérescence de la structure simpliciale. Aussi avant
de définir la notion de \precats~marquées canoniquement, va-t-on définir tout d'abord celle de diagrammes dégénérés.

\begin{defin}\index{diagramme dégénéré}
Soit $(\mathcal{C},\mathcal{F}_1,\mathcal{F}_2)$ une donnée de Segal facile.
Soit $(d,h)$ un diagramme de $\mathcal{C-PC}$ se factorisant par un diagramme $(d',h')$ comme dans
le diagramme commutatif suivant :
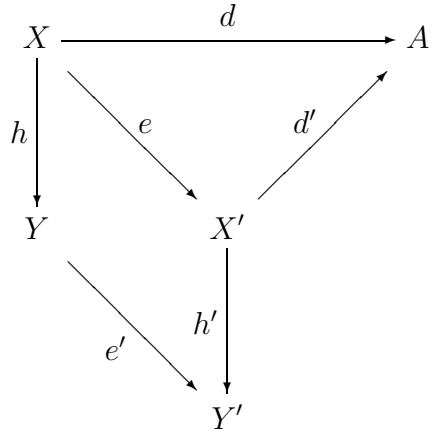
\begin{diagram}
X & & \rTo^{d} & & A \\
\dTo^h & \rdTo^e & & \ruTo^{d'} & \\
Y & & X' & & \\
 & \rdTo_{e'} & \dTo^{h'} & & \\
 & & Y' & & \\ 
\end{diagram}
On dira que le diagramme $(d,h)$ est une dégénérescence du diagramme
$(d',h')$ le long de $(e,e')$ dans les cas suivants :
\item -si $h$ est de la forme $\Delta[m]\Theta f$, avec $f$ dans $\mathcal{F}_1$, et $h'$ est soit un
$\Delta[k]\Theta f$, avec $0<k<m$, soit l'identité du point,
\item -si $h$ est de la forme $Boit_m(g)$, avec $g$ dans $\mathcal{F}_2$, et $h'$ est soit un $Boit_k(g)$, avec
$1<k<m$, soit l'identité de $\Delta[1]\Theta b(g)$, soit l'identité du
point.
\end{defin}

\begin{defin}\index{\precat!partiellement marquée canoniquement}
Soit $(\mathcal{C},\mathcal{F}_1,\mathcal{F}_2)$ une donnée de Segal facile.
Une \precat~partiellement marquée $(A,\mu)$ est marquée canoniquement si, 
pour tout diagramme $(d,h)$ de $\mu$ se factorisant par un diagramme $(d',h')$
le long de $(e,e')$, on a :
\item -si $h'$ est une flèche génératrice, alors le diagramme $(d',h')$ est lui aussi
dans $\mu$ et le relèvement marqué de $(d,h)$ est la composée de celui de
$(d',h')$ par $e'$,
\item -si $h'$ est une identité, alors le relèvement marqué de $(d,h)$ est
la composé $d'\circ e'$.
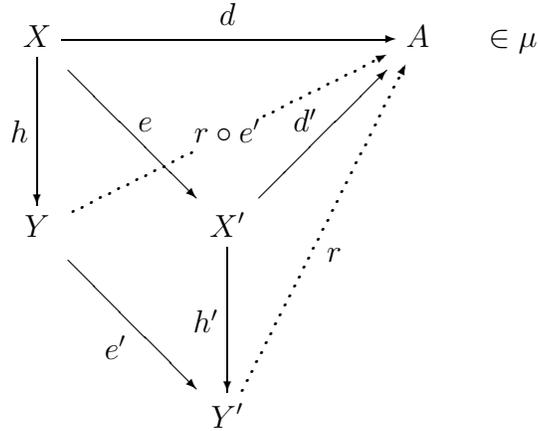
\begin{diagram}
X & & \rTo^{d} & & A & \in\mu \\
\dTo^h & \rdTo^e & & \ruDotsto(4,2)~{r\circ e'}\ruTo^{d'}\ruDotsto(2,4)_{r} & &\\
Y & & X' & & & \\
 & \rdTo_{e'} & \dTo^{h'} & & & \\
 & & Y' & & &\\ 
\end{diagram}
\end{defin}

De même que la notion de \precat~partiellement marquée a donné lieu à
celle de \cat~marquée, on a aussi une notion naturelle de \cat~canoniquement
marquée.

\begin{defin}\index{\cat!marquée canoniquement}
Soit $(\mathcal{C},\mathcal{F}_1,\mathcal{F}_2)$ une donnée de Segal facile.
Une \cat~marquée $(A,\mu)$ est marquée canoniquement si, en tant que
\precat, elle est marquée canoniquement au sens ci-dessus.
\end{defin}

Il est important de voir que si un morphisme d'une \precat~partiellement marquée
canoniquement vers une \precat~partiellement marquée préserve le marquage
alors la restriction du marquage de la \precat~but aux relèvements marqués
induits par le marquage de la \precat~source est canonique.\\
 
On obtient avec ces définitions deux sous-catégories pleine de $\mathcal{C-PC}_{m}$ : $\mathcal{C-PC}_{mc}$, la
sous-catégorie pleine des \precats~partiellement marquées canoniquement et la sous-catégorie pleine
$\mathcal{C-C}_{mc}$ des \cats~mar\-quées canoniquement. Tout d'abord, remarquons que par
construction $Cat$ et $Bigcat$ ne sont pas canoniquement marqués car on y
prend des sommes amalgamées pour des diagrammes dégénérés et que ces
sommes amalgamées serviront de relèvements marqués pour ces diagrammes. 
Ensuite, $\mathcal{C-PC}_{mc}$ et $\mathcal{C-C}_{mc}$ sont toutes deux munies d'un foncteur Oubli vers
$\mathcal{C-PC}$ qui est fidèle mais n'est pas plein car tout morphisme ne
préserve pas le marquage. En outre la notion de \cat~marquée canoniquement
est stable par limite.

\begin{lem}
Soit $(\mathcal{C},\mathcal{F}_1,\mathcal{F}_2)$ une donnée de Segal facile.
La catégorie $\mathcal{C-C}_{mc}$ des \cats~marquées canoniquement avec les
morphismes préservant les marquages est stable par limites.
\end{lem}
{\it Preuve :} application directe du lemme~\ref{mst} !\\

Cette notion de marquage canonique va nous permettre de définir un procédé
de catégorisation réduite qui à une \precat~va associer une \cat~mar\-quée
canoniquement, procédé qui non seulement
sera fonctoriel mais servira d'adjoint au foncteur Oubli des \cats~marquées
canoniquement vers les \precats.\\ 

Commençons par présenter une construction qui à une
\precat~partiellement marquée canoniquement associe une \precat~partiellement
marquée canoniquement ayant la propriété de se relever par rapport aux flèches qui ne se
relevaient pas dans la première.

\begin{prop}[-définition]\index{$Raj_c$}
Soit $(\mathcal{C},\mathcal{F}_1,\mathcal{F}_2)$ une donnée de Segal facile.
Soit $(A,\mu)$ une \precat~partiellement marquée canoniquement, il existe une
\precat~partiellement marquée canoniquement $(Raj_c(A,\mu),\nu)$ munie d'un morphisme
$raj_c:(A,\mu)\rightarrow (Raj_c(A,\mu),\nu)$ préservant le marquage, telle que :
\item 1) pour tout diagramme des flèches de $\mathcal{FG}_1$ vers $A$ n'appartenant pas
à $\mu$, l'extension du diagramme par $raj_c$ appartient à $\nu$,
\item 2) pour toute \precat~partiellement marquée canoniquement $(B,\lambda)$ et pour tout
morphisme $f$ de $A$ vers $B$ préservant le marquage tels que, pour tout diagramme des flèches de $\mathcal{FG}_1$ vers $A$ n'appartenant pas
à $\mu$, l'extension du diagramme par $f$ appartient à $\lambda$, il existe
un unique morphisme $\tilde{f}$ de $Raj_c(A,\mu)$ vers $B$
préservant le marquage et dont la précomposition par $raj_c$ est $f$.
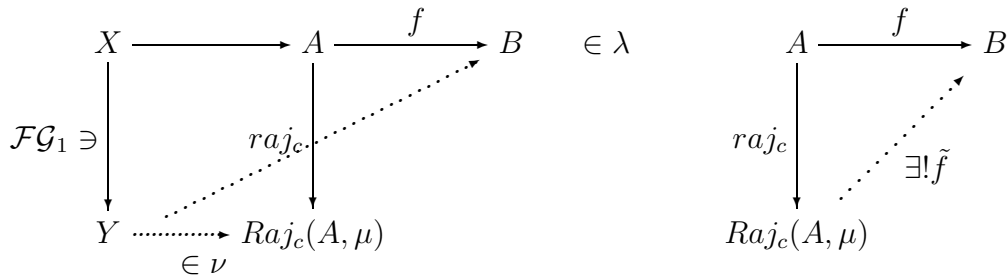
\begin{diagram}
X & \rTo & A & \rTo^{f} & B & \in\lambda& & A & \rTo^{f} & B \\
\dTo^{\mathcal{FG}_1\ni} & &\dTo^{raj_c} & \ruDotsto(4,2)&  & & & \dTo{raj_c} & \ruDotsto_{\exists !\tilde{f}} & \\
Y & \rDotsto_{\in\nu} & Raj_c(A,\mu) & &  & &  & Raj_c(A,\mu) & & \\
\end{diagram}
\end{prop}
{\it Preuve :}\\
On construit $Raj_c(A,\mu)$ de la manière suivante : on prend la somme
amalgamée de $A$ avec les flèches génératrices appartenant à des
diagrammes non dégénérés et non inclus dans $\mu$. Quant à $raj_c$, il
s'agit du morphisme naturel de A vers la somme amalgamée $Raj_c(A,\mu)$. 
Le marquage $\nu$ de $Raj_c(A,\mu)$ est défini ainsi : pour les diagrammes pour
$A$ appartenant à $\mu$, leurs extensions par $raj_c$ appartiennent à $\nu$, et donc
par construction $raj_c$ préserve le marquage. Pour les diagrammes non
dégénérés pour $A$
n'appartenant pas à $\mu$, leurs extensions par $raj_c$ ont pour relèvements
marqués dans $\nu$ les morphismes naturels des buts des flèches
génératrices vers la somme amalgamée $Raj_c(A,\mu)$. Enfin pour les
diagrammes dégénérés pour $A$
n'appartenant pas à $\mu$, on a deux cas. S'ils sont dégénérescences
de diagrammes de flèches génératrices, leurs extensions par $raj_c$ ont pour relèvements
marqués dans $\nu$ ceux induits par les relèvements marqués des diagrammes
non dégénérés dont ils sont les dégénérescences. Pour que cette dernière définition
ait un sens, il faut que tout diagramme dégénéré ne soit la
dégénérescence que d'un diagramme non dégénéré, ce qui est le
cas ici. Si les diagrammes sont dégénérescences de diagrammes avec
l'identité, les relèvements marqués dans $\nu$ de leurs extensions par
$raj_c$ seront les extensions par $raj_c$ des morphismes des buts des flèches
génératrices vers $A$ via l'objet de l'identité. Et donc
par construction le marquage de $(Raj_c(A,\mu),\nu)$ est canonique et
$(Raj_c(A,\mu),\nu)$ munie de $raj_c$ vérifie bien la propriété 1). En
outre, cette construction est bien définie car $\mathcal{C-PC}$ est une
catégorie admettant les coproduits et les sommes amalgamées calculées niveau par niveau.\\
\\
Pour $(B,\lambda)$ et $f:A\rightarrow B$ vérifiant la condition de la
propriété 2), le
diagramme ci-dessus, avec $f$ et $B$ à la place de $raj_c$ et $Raj_c(A,\mu)$ et
où le morphisme $\coprod_{\substack{X\rightarrow
Y\in\mathcal{FG}_1\\ X\rightarrow A\notin\mu}}Y\rightarrow B$ est celui induit
par la condition de la 
propriété 2) pour $B$ et $f$, commute et donc il existe un unique morphisme de
la somme amalgamée $Raj_c(A,\mu)$ vers $B$. Ce morphisme précomposé par $raj_c$ donne $f$ et il préserve le marquage, car $f$
et $raj_c$ préservent le marquage de $A$, que précomposé par un relèvement de
$Raj_c(A,\mu)$ non induit par $\mu$ (i.e. provenant de $\coprod_{\substack{X\rightarrow
Y\in\mathcal{FG}_1\\ X\rightarrow A\notin\mu}}Y\rightarrow Raj_c(A,\mu)$) et non dégénéré, il devient le relèvement correspondant pour
$B$ (qui est une des composantes de $\coprod_{\substack{X\rightarrow
Y\in\mathcal{FG}_1\\ X\rightarrow A\notin\mu}}Y\rightarrow B$) et enfin car $B$ est marqué canoniquement . En outre, la
préservation du marquage imposée au morphisme de $Raj_c(A,\mu)$ dans $B$, le
fait que $Raj_c(A,\mu)$ est une somme amalgamée et que les $\mathcal{FG}_1$-cofibrations sont des monomorphismes entraînent comme pour $E_{\Phi}$ que ce morphisme
est uniquement déterminé par la condition de la propriété 2) pour $B$ et donc est unique, ce
qui démontre la proposition-définition.\\
CQFD.\\

Tout d'abord remarquons que pour toute \cat~facile $C$ munie d'un morphisme
$g:A\rightarrow C$, la construction
ci-dessus fournit un morphisme de $Raj_c(A,\mu)$ vers $C$ dont la précomposition
par $raj_c$ redonne $g$. Toutefois il existe autant de tels morphismes que de
marquages canoniques possibles sur $C$. Mais pour un marquage canonique fixé et à supposé que $g$
préserve alors le marquage, il existe un unique morphisme de $Raj_c(A,\mu)$ vers $C$ dont la précomposition
par $raj_c$ redonne $g$ et qui préserve le marquage (c'est une application
directe de la proposition précédente aux \cats~marquées canoniquement).\\
Par ailleurs la propriété 2) de la proposition-définition assure l'unicité
à isomorphisme près de la construction $Raj_c$ et entraîne sa
fonctorialité en $(A,\mu)$.

\begin{prop}[-définition]\index{$Cat_c$}
Soit $(\mathcal{C},\mathcal{F}_1,\mathcal{F}_2)$ une donnée de Segal facile.
Soit $A$ une \precat~quelconque, il existe une
\cat~marquée canoniquement $Cat_c(A)$ munie d'un morphisme
$can^c_A:A\rightarrow Cat_c(A)$ et ayant la propriété universelle suivante :\\
pour toute \cat~marquée canoniquement $(B,\lambda)$ munie d'un morphisme
$f$ issu de $A$, il existe un unique morphisme $\tilde{f}$ de $Cat_c(A)$ vers $B$
préservant le marquage et dont la précomposition par $can^c_A$ est $f$.
\begin{diagram}
A & \rTo^{f} & B \\
\dTo{can^c_A} & \ruDotsto_{\exists !\tilde{f}} & \\
Cat_c(A) & & \\
\end{diagram}
\end{prop}
{\it Preuve :}\\
Soit $\alpha$ un cardinal régulier strictement plus grand que le plus petit cardinal régulier pour lequel les sources et buts des flèches de $\mathcal{F}_1$ et $\mathcal{F}_2$ sont petites.\\ 
Nous allons construire $Cat_c(A)$ par récurrence transfinie.
Posons : 
$$(A^0,\mu_0)=(A,\emptyset) \mbox{ et } can^c_0=Id_A$$
Supposons $(A^{\beta},\mu_{\beta})$ et $can^c_{\beta}:A\rightarrow A^{\beta}$ définis pour un ordinal
$\beta<\alpha$, alors on pose :
$$(A^{\beta+1},\mu_{\beta+1})=(Raj_c(A^{\beta},\mu_{\beta}),\nu) \mbox{ et }
can^c_{\beta+1}=raj_c\circ can^c_{\beta}$$
Soit un ordinal limite $\gamma<\alpha$, supposons $(A^{\beta},\mu_{\beta})$ et
$can^c_{\beta}:A\rightarrow A^{\beta}$ définis pour tout ordinal
$\beta<\gamma$, alors on pose :
$$(A^{\gamma},\mu_{\gamma})=\Big(\colim{\beta<\gamma}
A^{\beta},\colim{\beta<\gamma} \mu_{\beta}\Big) \mbox{ et }
can^c_{\gamma}=\colim{\beta<\gamma} can^c_{\beta}$$
Par récurrence transfinie, on a bien défini pour tout ordinal $\beta<\alpha$
une \precat~partiellement marquée canoniquement $(A^{\beta},\mu_{\beta})$ munie d'un
morphisme $can^c_{\beta}:A\rightarrow A^{\beta}$.
Posons donc :
$$ Cat_c(A)=\Big(\colim{\beta<\alpha} A^{\beta},\colim{\beta<\alpha}
\mu_{\beta}\Big) \mbox{ et }
can^c_A=\colim{\beta<\alpha} can^c_{\beta}$$
Cette construction est bien définie car $\mathcal{C-PC}$ admet les coproduits,
les sommes amalgamées et les colimites séquentielles transfinies par
proposition~\ref{precatclos} et
par construction $Cat_c(A)$ est partiellement marquée canoniquement.\\
\\
Montrons que $Cat_c(A)$ est
totalement marquée. Soit $X\rightarrow Y$ une flèche de $\mathcal{FG}_1$ et
$X\rightarrow Cat_c(A)$ un morphisme quelconque. Par le lemme~\ref{flgenpet}
, $X$ source d'un morphisme de $\mathcal{FG}_1$ est
$\alpha$-petit, alors, par définition~\ref{alphapetit} de $\alpha$-petit, le morphisme
$X\rightarrow Cat_c(A)$ se factorise à travers un $A^{\beta}$ pour un certain
$\beta<\alpha$. On a alors deux cas. Si $X\rightarrow Y$ possède un relèvement
par rapport à $A^{\beta}$ marqué dans $\mu_{\beta}$ alors son extension par
$A^{\beta}\rightarrow Cat_c(A)$ est dans la colimite des $\mu_{\beta}$ et donc est
marqué pour $Cat_c(A)$. Si en revanche $X\rightarrow Y$ ne se relève pas par
rapport à $A^{\beta}$, alors par la propriété 1) de la construction
$Raj_c$, son extension par
$A^{\beta}\rightarrow Raj_c(A^{\beta},\mu_{\beta})=A^{\beta+1}$ appartient à
$\mu_{\beta+1}$ et donc son extension par
$A^{\beta+1}\rightarrow Cat_c(A)$ est dans la colimite des $\mu_{\beta}$ et donc est
marqué pour $Cat_c(A)$. Ainsi $Cat_c(A)$ est totalement marquée, ce qui fait de
$Cat_c(A)$ une \cat~marquée.\\
\\
Soit $(B,\lambda)$ une \cat~marquée canoniquement munie d'un morphisme $f:A\rightarrow B$.
Construisons par récurrence transfinie le morphime $\tilde{f}:
Cat_c(A)\rightarrow B$ faisant commuter le diagramme de la proposition.
Au cran 0, on a bien $f^0:A^0\rightarrow B$ qui n'est autre que $f$ préservant
le marquage et tel que
$f^0\circ can_0=f$. Supposons construit $f^{\beta}:A^{\beta}\rightarrow B$
préservant le marquage tel que $f^{\beta}\circ can_{\beta}=f$ pour
$\beta<\alpha$, alors par la propriété 2) de la construction $Raj_c$, il existe
un morphisme $f^{\beta+1}:Raj_c(A^{\beta},\mu_{\beta})=A^{\beta+1}\rightarrow B$ préservant
le marquage et dont la précomposition par $raj_c$ est $f^{\beta}$ d'où il vient
que :
$$f^{\beta+1}\circ can_{\beta+1}=f^{\beta+1}\circ raj_c \circ can_{\beta}= f^{\beta}\circ can_{\beta}=f$$
ce que l'on voulait.\\
\\
Soit maintenant un ordinal limite $\gamma<\alpha$ et
supposons construit, pour tout $\beta<\gamma$, $f^{\beta}:A^{\beta}\rightarrow B$
préservant le marquage tel que $f^{\beta}\circ can_{\beta}=f$. Par propriété universelle de la colimite, il existe un
morphisme $f^{\gamma}:\mbox{colim} A^{\beta}=A^{\gamma}\rightarrow B$ préservant le
marquage, colimite des marquages $\mu_{\beta}$, et telle que $f^{\gamma}\circ
can_{\gamma}=f$. Par récurrence transfinie, on a donc pour tout $\beta<\alpha$
des morphismes $f^{\beta}:A^{\beta}\rightarrow B$
préservant le marquage tel que $f^{\beta}\circ can_{\beta}=f$. On définit
alors $\tilde{f}:Cat_c(A)\rightarrow B$ comme colimite des $f^{\beta}$ pour
$\beta<\alpha$ et, par propriété de colimite, $\tilde{f}$ préserve le
marquage et précomposé par $(can_c)_A$ redonne bien $f$. Cette construction est
bien définie car $\mathcal{C-PC}$ possède les coproduits, les sommes
amalgamées et les colimites séquentielles transfinies. En outre, le fait
d'imposer à $\tilde{f}$ de préserver le marquage, allié au fait que les $\mathcal{FG}_1$-cofibrations sont des monomorphismes, entraîne l'unicité de
$\tilde{f}$ et ainsi $Cat_c(A)$ vérifie la propriété universelle de la
proposition.\\
CQFD.\\

On peut remarquer que la construction de la factorisation ci-dessus s'applique aussi bien pour les
\cats~faciles mais, comme on l'a déjà remarqué pour la construction $Raj_c$, il n'y aura pas d'unicité
de la factorisation à travers $Cat_c(A)$.\\

Comme la catégorisation Cat, la catégorisation réduite $Cat_c$ est un
adjoint à gauche du foncteur oubli des \cats~canoniquement marquées vers les
\precats.

\begin{cor}
Soit $(\mathcal{C},\mathcal{F}_1,\mathcal{F}_2)$ une donnée de Segal facile.
La construction $Cat_c$ est un foncteur de $\mathcal{C-PC}$ vers
$\mathcal{C-C}_{mc}$ qui est adjoint à gauche du foncteur Oubli, i.e. on a un
isomorphisme naturel en $A$ et $(B,\lambda)$, induit par le morphisme $can^c_A$ :
$$\operatorname{Hom}_{\mathcal{C-C}_{mc}}(Cat_c(A),(B,\lambda))=\operatorname{Hom}_{\mathcal{C-PC}}(A,Oubli((B,\lambda)))$$
De ce fait, le foncteur $Cat_c$ préserve les colimites.
\end{cor}
{\it Preuve :}\\
De la propriété universelle de
$Cat_c(A)$ découle son unicité à isomorphisme près, ce qui entraîne la fonctorialité de la construction $Cat_c$ et la
naturalité de $can^c_A$ en $A$, mais aussi que $Cat_c$ est en
fait un adjoint à gauche du foncteur Oubli des \cats~marquées canoniquement vers les
\precats. Enfin comme foncteur adjoint à gauche, $Cat_c$ préserve les
colimites.\\ 
CQFD.\\

Comme nous venons de construire un troisième procédé de catégorisation,
il est très utile de le comparer avec Cat notamment pour vérifier si les
notions d'\eqs~faibles induites par ces deux catégorisations sont les mêmes.

\begin{lem}\label{catred}
Soit $(\mathcal{C},\mathcal{F}_1,\mathcal{F}_2)$ une donnée de Segal facile.
Alors pour toute \precat~$A$, le morphisme canonique
$can^c_A:A\rightarrow Cat_c(A)$ est une \eq~faible de $\mathcal{C-PC}$ et, pour
tout morphisme de \precats~$f:A\rightarrow B$, $f$ est une \eq~faible si et
seulement si $Cat_c(f)$ est une \eq~de \cats.
\end{lem}
{\it Preuve :}\\ 
On peut modifier le marquage de $Bigcat(A)$ afin qu'il devienne
canonique. Dès lors il existe un unique morphisme préservant le marquage de $Cat_c(A)$ vers
$Bigcat(A)$ factorisant $bcan_A$. De plus, $Cat_c(A)$ étant une
\cat~marquée, il existe un morphisme (non nécessairement unique)
préservant le marquage de $Bigcat(A)$ vers $Cat_c(A)$. Comme pour $Cat(A)$
(c'est-à-dire comme dans la démonstration de la proposition~\ref{Ephi3}), on
montre que ces morphismes sont des \eqs~de \cats~si $A$ est une \cat~facile, on en
déduit par "trois pour deux" pour les \eqs~de \cats~que $can^c_A:A\rightarrow Cat_c(A)$ est une
\eq~de \cats~si $A$ est une \cat~facile.\\
\\ 
Puis toujours de la même manière que pour $Cat(A)$,
en utilisant la propriété universelle de $Cat_c(A)$ et "trois pour deux"
pour les \eqs~de \cats, on
montre que pour toute \precat~$A$ l'image par $Cat_c$ du morphisme canonique $can^c_A$ est
une \eq~de \cats. Ceci montre que $Cat_c$ est une bonne catégorisation.\\ 
\\
En utilisant alors le lemme~\ref{bcateq} avec la catégorisation $Cat_c$, on
obtient que, pour tout morphisme de \precats~$f$, $Cat_c(f)$ est une \eq~de
\cats~si et seulement si $Bigcat(f)$ l'est. Or par le lemme~\ref{bigcateq},
$Bigcat(f)$ est une \eq~de \cats~si et seulement si $Cat(f)$ l'est, ceci montre
la seconde partie du lemme.\\ 
\\
Comme en outre on a déjà montré que, pour
toute \precat~$A$, l'image par $Cat_c$ du morphisme canonique $can^c_A$ est une
\eq~de \cats, la seconde partie du lemme que l'on vient de prouver montre que
$can^c_A$ est bien une \eq~faible de \precats.\\
CQFD.\\

Nous avons donc bien une nouvelle catégorisation $Cat_c$ qui induit la même
notion d'\eq~faible que Cat. Comme nous l'avons introduite afin de préserver
la notion de \precats~librement ordonnées, nous allons maintenant montrer que
$Cat_c$ conserve bien cette notion.

\subsection{Equivalence de \precats~librement ordonnées}

L'utilité d'avoir considéré les marquages canoniques et d'avoir construit
une catégorisation ne prenant plus de sommes amalgamées de diagrammes
dégénérés est que cette catégorisation réduite préserve la notion
de \precats~librement ordonnées et permettra ainsi de donner un critère pratique
pour montrer qu'un morphisme de \precats~librement ordonnées est une
\eq~faible, notion difficile à manipuler par ailleurs. Commençons donc par
montrer étape par étape que la catégorisation réduite $Cat_c$ préserve
bien la notion de \precats~librement ordonnées. Comme les hypothèses sur la
donnée de Segal facile qui assurent ce résultat sont valables pour tous les
autres résultats de cette partie, nous allons d'abord énoncer ces
hypothèses.

\begin{hyp}\label{plo}
Soit $(\mathcal{C},\mathcal{F}_1,\mathcal{F}_2)$ une donnée de Segal facile
vérifiant les propriétés suivantes :
\item 1) la \cmf~sous-jacente $\mathcal{C}$ est une \cmf~simpliciale, 
\item 2) il n'existe aucun morphisme de but l'ensemble vide qui ne soit pas un
isomorphisme,
\item 3) les cofibrations sont stables par produit fibré au-dessus d'un objet
discret.
\end{hyp} 

\begin{lem}
Supposons vraies les hypothèses~\ref{plo}.\\
Si $(A,\mu)$ est une \precat~librement ordonnée 
et partiellement marquée canoniquement, alors
$Raj_c(A,\mu)$ est librement ordonnée. De plus pour tout couple ordonné d'objets de $A$
adjacents $(a,a')$, le morphisme $\big(raj^c_A\big)_1(a,a'):A_1(a,a')\rightarrow
Raj_c(A)_1(a,a')$ est une cofibration triviale de la \cmf~$\mathcal{C}$.
\end{lem}
{\it Preuve :}\\
Par construction, $Raj_c(A)$ est une somme amalgamée de $A$ par un coproduit
de flèches génératrices qui sont toutes des isomorphismes sur les objets.
Ainsi la construction $Raj_c$ laisse les objets invariants et donc $Raj_c(A)_0$
qui n'est autre que $A_0$ est bien un ensemble fini totalement ordonné.\\
\\
Comme $Raj_c(A)$ est une somme amalgamée de $A$ par un coproduit
de flèches génératrices indexé par les diagrammes non dégénérés
et non inclus dans $\mu$, on peut le développer comme une colimite
séquentielle transfinie de sommes amalgamées par une flèche
génératrice d'un diagramme non dégénéré
et non inclus dans $\mu$. Non seulement ces flèches génératrices sont des
cofibrations entre \precats~librement ordonnées strictes induisant
l'identité sur les objets, d'après le lemme~\ref{FGordre}, mais toute
flèche génératrice $h:X\rightarrow Y$ vérifie la propriété que, pour
tout couple ordonnée $(x,x')$ d'objets adjacents de $X$, $X_1(x,x')$ n'est pas
isomorphe à l'ensemble vide, car d'après la démonstration du
lemme~\ref{FGordre}, $X_1(x,x')$ n'est autre que la source d'une flèche de
$\mathcal{F}_1$
ou le but d'une flèche de $\mathcal{F}_2$ qui sont non isomorphes à l'ensemble vide par
hypothèse. Ainsi par le lemme~\ref{mpo}, tout morphisme de la source d'une
flèche génératrice vers $A$ préserve l'ordre, et même strictement car
n'apparaissent dans $Raj_c(A,\mu)$ que celles des diagrammes non
dégénérés. Ainsi la somme amalgamée par une flèche génératrice
d'un diagramme non dégénéré et non inclus dans $\mu$ vérifie les
hypothèses du 2) du lemme~\ref{stordre}, ce qui en fait une \precat~librement
ordonnée. Ainsi $Raj_c(A,\mu)$ devient une colimite séquentielle transfinie
de \precats~librement ordonnées vérifiant les hypothèses du 3) du
lemme~\ref{stordre}, ce qui fait de $Raj_c(A,\mu)$ une \cat~librement
ordonnée.\\
\\
Soit $(a,a')$ un couple ordonné d'objets adjacents de $A$.
$Raj_c(A,\mu)_1(a,a')$ n'est autre que la somme amalgamée de $A_1(a,a')$ par
un coproduit de $h_1(x,x')$, avec $h$ une flèche génératrice d'un
diagramme $(d,h)$ non dégénéré et non inclus dans $\mu$ et $x$ et $x'$ des
objets de la source de $h$ s'envoyant par $d$ sur $a$ et $a'$. Comme on l'a vu
ci-dessus, $d$ préserve strictement l'ordre, donc est une injection au niveau
des objets. Ceci entraîne que le couple $(x,x')$ est unique et ordonné mais
surtout que $x$ et $x'$ sont adjacents. Si $h$ est de type $\Delta[m]\Theta f$
alors $h_1(x,x')$ n'est autre que $f$, qui étant dans $\mathcal{F}_1$ est une
cofibration triviale. Si $h$ est de type $Boit_m(g)$ alors $h_1(x,x')$ n'est autre que
l'identité du but de $g$, donc une cofibration triviale. Comme par hypothèse,
$\mathcal{C}$ est une \cmf~simpliciale dont tous les objets sont cofibrants, par
corollaire~\ref{eqcoprod}, le coproduit des $h_1(x,x')$ est une \eq~faible et
aussi une cofibration comme coproduit de cofibrations. Par stabilité des
cofibrations triviales par somme amalgamée le long d'un morphisme dans la
\cmf~$\mathcal{C}$, il vient que $raj^c_1(a,a')$ est une cofibration triviale.\\
CQFD.

\begin{lem}\label{catordre}
Supposons vraies les hypothèses~\ref{plo}.\\
Si $A$ est une \precat~librement ordonnée alors
$Cat_c(A)$ est librement ordonnée. De plus pour tout couple ordonné d'objets de $A$
adjacents $(a,a')$, le morphisme $\big(can^c_A\big)_1(a,a'):A_1(a,a')\rightarrow
Cat_c(A)_1(a,a')$ est une cofibration triviale de la \cmf~$\mathcal{C}$.
\end{lem}
{\it Preuve :}\\
Il est facile de voir que la construction $Cat_c$ laisse les objets invariants.
Ainsi $Cat_c(A)_0$ n'est autre que $A_0$ ensemble fini totalement ordonné.\\
Par construction, $Raj_c(A)$ est une somme amalgamée de $A$ par un coproduit
de cofibrations induisant l'identité sur les objets donc le morphisme naturel $A\rightarrow Raj_c(A)$ est une
cofibration induisant l'identité sur les objets. Par cette remarque et par le
lemme précédent, il vient que $Cat_c(A)$ est une colimite séquentielle
transfinie de \precats~librement ordonnées vérifiant les hypothèses du 3)
du lemme~\ref{stordre}, ce qui en fait une \cat~librement ordonnée.\\
\\
On conclut simplement en remarquant que $\big(can^c_A\big)_1(a,a')$ est une
colimite séquen\-tielle transfinie de cofibrations triviales, d'après le lemme
précédent, ce qui en fait une cofibration triviale par stabilité des
cofibrations triviales par colimite séquentielle transfinie dans la
\cmf~$\mathcal{C}$.\\
CQFD.\\

Nous avons donc montré que la catégorisation réduite préserve les
\precats~librement marquées mais aussi l'homotopie de leurs espaces de
morphismes entre objets adjacents. C'est cette dernière propriété qui nous
permettra de montrer qu'un morphisme de \precats~librement ordonnées
préservant l'ordre et induisant sur chaque espace de morphismes entre objets
adjacents une \eq~faible est lui-même une \eq~faible de \precats. Afin de
montrer cette proposition, donnons un lemme qui permet de simplifier la
démonstration de celle-ci.

\begin{lem}
Supposons vraies les hypothèses~\ref{plo}.\\
Soient $A$ une \precat~librement ordonnée et $(x_0,\ldots,x_m)$ une
suite ordonnée d'objets de $A$ deux-à-deux adjacents. Alors les morphismes
\item $Cat_c(A)_m(x_0,\ldots,x_m)\rightarrow Cat_c(A)_1(x_0,x_m),$
\item $Cat_c(A)_m(x_0,\ldots,x_m)\rightarrow Cat_c(A)_1(x_0,x_1)\times\ldots\times
Cat_c(A)_1(x_{m-1},x_m) \mbox{ et}$ 
\item $A_1(x_0,x_1)\times\ldots\times
A_1(x_{m-1},x_m)\rightarrow Cat_c(A)_1(x_0,x_1)\times\ldots\times
Cat_c(A)_1(x_{m-1},x_m)$ sont des \eqs~faibles de $\mathcal{C}$.
\end{lem}
{\it Preuve :}\\
$A$ étant une \precat~librement ordonnée, par le lemme précédent, on a
que $Cat_c(A)$ est librement ordonnée. D'où par le ii) des
\precats~librement ordonnées, on a que
$Cat_c(A)_m(x_0,\ldots,x_m)\rightarrow Cat_c(A)_1(x_0,x_m)$ est une
\eq~faible. Comme $Cat_c(A)$ est une \cat,
$Cat_c(A)_m(x_0,\ldots,x_m)\rightarrow Cat_c(A)_1(x_0,x_1)\times\ldots\times
Cat_c(A)_1(x_{m-1},x_m)$ est une \eqc~d'\obcs, donc c'est une \eq~faible. Enfin comme les $x_i$ sont
deux-à-deux adjacents, par le lemme précédent, on obtient que les
morphismes $A_1(x_i,x_{i+1})\rightarrow Cat_c(A)_1(x_i,x_{i+1})$ sont des
\eqs~faibles et donc leur produit\\ $A_1(x_0,x_1)\times\ldots\times
A_1(x_{m-1},x_m)\rightarrow Cat_c(A)_1(x_0,x_1)\times\ldots\times
Cat_c(A)_1(x_{m-1},x_m)$ aussi.\\
CQFD.

\begin{prop}\label{eqlord}
Supposons vraies les hypothèses~\ref{plo}.\\
Si $f:A\rightarrow B$ est un
morphisme entre \precats~librement ordonnées qui induit un
isomorphisme sur les objets préservant l'ordre et tel que, pour tout couple
ordonné d'objets de $A$ adjacents $(a,a')$, $f_1(a,a'):A_1(a,a')\rightarrow
B_1(a,a')$ est une \eq~faible de $\mathcal{C}$, alors $f$ est une \eq~faible
dans $\mathcal{C-PC}$.
\end{prop}
{\it Preuve :}\\
Pour montrer que $f$ est une \eq~faible, d'après le lemme~\ref{catred}, il
suffit de montrer que $Cat_c(f)$ est une \eq~de \cats. Comme la construction
$Cat_c$ laisse invariant les objets et que $f$ est un isomorphisme sur les
objets alors $Cat_c(f)$ aussi, donc $Cat_c(f)$ est essentiellement surjectif.\\
Montrons maintenant que pour tout couple d'objets $(a,a')$ de $A$, le morphisme
$Cat_c(f)_1(a,a')$ est une \eq~faible dans $\mathcal{C}$. Pour cela nous allons
raisonner par disjonction des cas.\\
Si $a$ et $a'$ sont égaux, alors, comme $Cat_c(A)$ et $Cat_c(B)$ sont librement
ordonnées, car $A$ et $B$ le sont, il vient par le iii) des \precats~librement
ordonnées que $Cat_c(f)_1(a,a')$ n'est autre
que l'identité du point.\\
Si $a$ et $a'$ sont distincts, comme $A_0$ est totalement ordonné, on a soit
$a<a'$ soit $a>a'$. Dans ce deuxième cas, comme $Cat_c(A)$ et $Cat_c(B)$ sont librement
ordonnées, car $A$ et $B$ le sont, il vient par le i) des \precats~librement
ordonnées que $Cat_c(f)_1(a,a')$ n'est autre
que l'identité de l'ensemble vide.\\
\\
Il ne reste donc plus qu'à traiter le cas où $a<a'$. Nous diviserons ce cas
en deux sous-cas. Tout d'abord si $a$ et $a'$ sont adjacents, alors
$(can^c_A)_1(a,a')$ et $(can^c_B)_1(a,a')$ sont des \eqs~faibles par le
lemme~\ref{catordre} et $f_1(a,a')$ est une \eq~faible par hypothèse sur $f$.
Comme on a l'égalité $Cat_c(f)_1(a,a')\circ
(can^c_A)_1(a,a')=(can^c_B)_1(a,a')\circ f_1(a,a')$, par "trois pour deux" dans
la \cmf~$\mathcal{C}$, il vient que $Cat_c(f)_1(a,a')$ est une \eq~faible.\\
\\
Terminons donc par le cas $a<a'$ avec $a$ et $a'$ non adjacents. Comme $A_0$ est
un ensemble fini et totalement ordonné, il existe une suite ordonnée
$(a_0,\ldots,a_m)$ d'objets de $A$ deux-à-deux adjacents telle que $a_0=a$ et
$a_m=a'$. En appliquant le lemme précédent, il vient que $Cat_c(A)_m(x_0,\ldots,x_m)\rightarrow Cat_c(A)_1(x_0,x_m)$ est une
\eq~faible, de même pour $B$, et par "trois pour deux" dans la
\cmf~$\mathcal{C}$, on obtient que $Cat_c(f)_1(a,a')$ est une \eq~faible si et seulement si
$Cat_c(f)_m(a_0,\ldots,a_m)$ l'est. De la
même manière, le fait que $Cat_c(A)_m(x_0,\ldots,x_m)\rightarrow Cat_c(A)_1(x_0,x_1)\times\ldots\times
Cat_c(A)_1(x_{m-1},x_m)$ et $A_1(x_0,x_1)\times\ldots\times
A_1(x_{m-1},x_m)\rightarrow Cat_c(A)_1(x_0,x_1)\times\ldots\times
Cat_c(A)_1(x_{m-1},x_m)$ soient des \eqs~faibles, d'après le lemme
précédent, montre, par "trois pour deux" dans $\mathcal{C}$, que $Cat_c(f)_m(a_0,\ldots,a_m)$ est une
\eq~faible si et seulement si $\prod Cat_c(f)_1(a_i,a_{i+1})$ l'est et que $\prod
Cat_c(f)_1(a_i,a_{i+1})$ est une \eq~faible si et seulement si $\prod f_1(a_i,a_{i+1})$ l'est.
Donc $Cat_c(f)_1(a,a')$ est une \eq~faible si et seulement si $\prod f_1(a_i,a_{i+1})$ l'est,
par transitivité de l'\eq~logique. Or comme les
$a_i$ sont deux-à-deux adjacents, par hypothèse, les morphismes
$f_1(a_i,a_{i+1})$ sont des équi\-valences faibles ainsi que leur produit, dans la \cmf~$\mathcal{C}$.\\
CQFD.\\

Cette proposition sera utile pour démontrer la proposition sur le produit
des \eqs~faibles dans $\mathcal{C}$ mais c'est également un critère très
simple pour démontrer que certains morphismes sont des \eqs~faibles. Ce
critère va nous permettre en particulier de montrer que les flèches génératrices sont des
cofibrations triviales.

\begin{cor}
Supposons vraies les hypothèses~\ref{plo}.\\
Les flèches géné\-ratrices de
la famille $\mathcal{FG}_1$ sont des cofibrations triviales de $\mathcal{C-PC}$.
\end{cor}
{\it Preuve}\\
Comme à chaque niveau les flèches génératrices sont des coproduits de
flèches de $\mathcal{F}_1$ ou $\mathcal{F}_2$ et que ces dernières sont des
cofibrations par définition de la donnée de Segal facile, alors par
stabilité des cofibrations par coproduit dans la \cmf~$\mathcal{C}$, il vient
que les flèches génératrices sont des cofibrations niveau par niveau.\\
Par le lemme~\ref{FGordre}, les flèches génératrices sont des morphismes
de \precats~librement ordonnées strictes qui induisent l'identité sur les
objets. Donc par la proposition précédente, il suffit de montrer que, pour
toute  flèche génératrice $h$ et pour tout couple ordonné $(x,y)$
d'objets adjacents de la source de $h$, $h_1(x,y)$ est une \eq~faible.\\
Si $h$ est du type $\Delta[m]\Theta f$, avec $f$ dans $\mathcal{F}_1$, alors
$h_1(x,y)$ n'est autre que $f$, d'après la démonstration du
lemme~\ref{FGordre}. Or $f$ est une cofibration triviale, donc une
\eq~faible.\\
Si $h$ est du type $Boit_m(g)$, avec $g$ dans $\mathcal{F}_2$, alors
$h_1(x,y)$ n'est autre que l'identité du but de $g$, d'après la démonstration du
lemme~\ref{FGordre}. Donc c'est une \eq~faible.\\
CQFD.\\

Comme nous avons maintenant en main tous les outils nécessaires pour montrer
le résultat cherché sur le produit cartésien de \cats, nous allons
consacrer la prochaine partie à la démonstration de ce résultat.

\subsection{Produit d'\eq~faibles des \precats}

Le but de cette section est d'étudier le produit cartésien des \precats~et
son comportement vis-à-vis de la catégorisation. Avant d'énoncer le
théorème de cette section, intéressons-nous au lemme technique suivant.

\begin{lem}\label{ictcub}
Soit $(\mathcal{C},\mathcal{F}_1,\mathcal{F}_2)$ une donnée de Segal facile
telle que la \cmf~sous-jacente $\mathcal{C}$ est une \cmf~simpliciale engendrée
par monomorphismes.\\  
\\
Soient $A\rightarrow A'$, $B\rightarrow B'$ et
$C\rightarrow C'$ des iso-cofibrations triviales. Supposons que le morphisme $A'\coprod_{A} C\rightarrow C'$
soit une cofibration. Alors leur somme amalgamée
$B\coprod_{A} C\rightarrow B'\coprod_{A'} C'$ est aussi une iso-cofibration triviale.
\end{lem}
{\it Preuve :}\\
Par le lemme~\ref{ictpush}, le morphisme $C\rightarrow A'\coprod_{A} C$ est une
iso-cofibration triviale comme somme amalgamée le long de $A\rightarrow C$ de
l'iso-cofibration triviale $A\rightarrow A'$. Par propriété universelle de
la somme amalgamée $A'\coprod_{A} C$, il existe un morphisme $A'\coprod_{A}
C\rightarrow C'$ dont la précomposition par l'iso-cofibration triviale $C\rightarrow A'\coprod_{A} C$
n'est autre que l'iso-cofibration triviale $C\rightarrow C'$. Par l'axiome
"trois pour deux" dans $\mathcal{C-PC}$ pour les isomorphismes sur les objets et pour les
\eqs~faibles, on obtient que $A'\coprod_{A} C\rightarrow C'$ est une
iso-\eq~faible. Or par hypothèse, c'est aussi une cofibration. Donc c'est une
iso-cofibration triviale. En outre le morphisme $B\coprod_{A} C\rightarrow
B'\coprod_{B}(B\coprod_{A} C)$ est une iso-cofibration
triviale comme somme amalgamée le long de $B\rightarrow B\coprod_{A} C$ de
l'iso-cofibration triviale $B\rightarrow B'$, par lemme~\ref{ictpush}.
Or $B'\coprod_{B}(B\coprod_{A} C)$
n'est autre que $B'\coprod_{A} C$. Donc on a obtenu que le morphisme
$B\coprod_{A} C\rightarrow B'\coprod_{A} C$ est une iso-cofibration triviale.
Toujours par le lemme~\ref{ictpush}, le morphisme $B'\coprod_{A} C\rightarrow
(B'\coprod_{A} C)\coprod_{A'\coprod_{A} C} C$ est une iso-cofibration
triviale comme somme amalgamée le long de $A'\coprod_{A} C\rightarrow
B'\coprod_{A} C$ de l'iso-cofibration triviale $B\coprod_{A} C\rightarrow B'\coprod_{A'} C'$.
Or $(B'\coprod_{A} C)\coprod_{A'\coprod_{A} C} C$ n'est
autre que $B'\coprod_{A'} C'$. Donc le morphisme $B'\coprod_{A} C\rightarrow B'\coprod_{A'} C'$ est une
iso-cofibration triviale. Enfin on remarque que le morphisme universel
$B\coprod_{A} C\rightarrow B'\coprod_{A'} C'$ est une iso-cofibration triviale
comme composée des iso-cofibrations triviales $B\coprod_{A} C\rightarrow
B'\coprod_{A} C$ et $B'\coprod_{A} C\rightarrow B'\coprod_{A'} C'$.\\
CQFD.

\begin{theo}\label{catprod}
Soit $(\mathcal{C},\mathcal{F}_1,\mathcal{F}_2)$ une donnée de Segal facile
vérifiant les propriétés suivantes :
\item 1) la catégorie sous-jacente $\mathcal{C}$ est une \cmf~simpliciale
engendrée par monomorphismes,
\item 2) les cofibrations génératrices de la \cmfcof~$\mathcal{C}$ ont leurs
sources et buts $\alpha$-petits, pour un cardinal régulier $\alpha$ plus grand
qu' $\aleph_0$,
\item 3) les cofibrations génératrices de $\mathcal{C}$ ont leurs buts connexes,
\item 4) les réunions de deux
sous-objets sont exactement les sommes amalgamées de ces sous-objets au-dessus
de leur intersection,
\item 5) pour tout monomorphisme $j:A\rightarrow B$ de $\mathcal{C}$, la classe contenant un représentant de chaque classe d'isomorphismes d'objets $C$ par lequel $j$ se
factorise dans la sous-catégorie des monomorphismes est un ensemble,
\item 6) les colimites séquentielles transfinies de la sous-catégorie des
monomorphismes de $\mathcal{C}$ existent et sont les colimites séquentielles
transfinies de $\mathcal{C}$,
\item 7) il n'existe aucun morphisme de but l'ensemble vide qui ne soit pas un
isomorphisme.\\
\\
Soient $A$ et $B$ deux \precats. Alors le morphisme $A\times B\rightarrow
Cat(A)\times Cat(B)$ est une iso-cofibration
triviale.
\end{theo}

Afin de démontrer ce théorème, nous allons tout d'abord faire les
réductions suivantes.

\begin{lem}
Supposons vraies les hypothèses du théorème~\ref{catprod}.\\
Le théorème~\ref{catprod} équivaut au résultat suivant, qu'on notera (I) :\\
pour toute \precat~$B$ et pour toute flèche $C\rightarrow D$ génératrice 
de \cats~faciles, le morphisme $C\times B\rightarrow D\times B$ est une iso-cofibration triviale.
\end{lem}
{\it Preuve :}\\
Le résultat (I) est un corollaire évident du théorème, donc il suffit de
montrer que le théorème découle du résultat (I) pour avoir l'équivalence
cherchée. Supposons donc le résultat (I) vrai. Soit $A$ une \precat, pour
tout morphisme de la source $C$ d'une flèche génératrice de \cats~faciles vers $A$, on obtient par lemme~\ref{ictpush} que le morphisme $A\times
B\rightarrow (A\coprod_{C} D)\times B$ est une iso-cofibration triviale comme
somme amalgamée le long de $C\times B\rightarrow A\times B$ de
l'iso-cofibration triviale $C\times B\rightarrow D\times B$. Or le morphisme
$A\times B\rightarrow Cat(A)\times B$ s'obtient comme une colimite
séquentielle transfinie de tels morphismes et donc, par lemme~\ref{ictcolim1}, c'est
une iso-cofibration triviale. De même, on peut montrer que $Cat(A)\times
B\rightarrow Cat(A)\times Cat(B)$ est une iso-cofibration triviale. Ainsi par
composition, $A\times B\rightarrow Cat(A)\times Cat(B)$ est encore une
iso-cofibration triviale, ce qui montre le théorème.\\
CQFD.

\begin{lem}
Supposons vraies les hypothèses du théorème~\ref{catprod}.\\
Le résultat (I) est équivalent au résultat (II) suivant :\\
pour tout objet $X$ de $\mathcal{C}$, pour tout entier $n$ et pour toute flèche $E\rightarrow
F$ génératrice de \cats~faciles de type 2, le morphisme
$E\times\Delta[n]\Theta X\rightarrow F\times\Delta[n]\Theta X$ est une
iso-cofibration triviale.
\end{lem}
{\it Preuve}\\
Il est facile de voir que le résultat (II) est un cas particulier du
résultat (I). Il suffit donc de montrer que (II) implique (I). Tout d'abord,
on peut remarquer que le résultat (I) est toujours vrai si la flèche
génératrice $C\rightarrow D$ est de type 1. En effet, dans ce cas, cette
flèche est de la forme $\Delta[m]\Theta C'\rightarrow \Delta[m]\Theta D'$,
avec $C'\rightarrow D'$ un élément de $\mathcal{F}_1$, donc une cofibration triviale dans $\mathcal{C}$. Donc pour
tout entier $n$, on a que $C_n\rightarrow D_n$ est une cofibration triviale dans
$\mathcal{C}$ comme coproduit de cofibrations triviales dans la
\cmf~$\mathcal{C}$. Par hypothèse sur la \cmf~$\mathcal{C}$, il vient que
$(C\times B)_n\rightarrow (D\times B)_n$ est aussi une cofibration triviale
comme produit des cofibrations triviales $C_n\rightarrow D_n$ et $Id_{B_n}$.
Comme en outre $C\rightarrow D$ est un isomorphisme sur les objets alors
$C\times B\rightarrow D\times B$ aussi et, par la proposition~\ref{eqfniv}, il vient
que c'est aussi une cofibration triviale. 
Comme le résultat (I) est toujours vrai pour les flèches génératrices de type 1, on doit
donc seulement montrer que (II) implique (I) pour les flèches génératrices
de type 2.  Par le corollaire~\ref{Icof2}, on a que les
$\mathcal{I}$-cofibrations sont exactement les cofibrations. Donc les objets,
étant tous cofibrants, sont des rétracts d'un objet qui est une colimite séquentielle transfinie de
sommes amalgamées d'éléments de $\mathcal{I}$ à partir de l'objet
initial. Si $B$ est un rétract de $B'$, alors $E\times B\rightarrow F\times B$ est un rétract de
$E\times B'\rightarrow F\times B'$ de manière évidente . Comme les
\eqs~faibles de $\mathcal{C-PC}$ sont stables par
rétracts par lemme~\ref{eqfre}, il suffit de montrer (I) pour les objets $B$
colimites séquentielles transfinies de sommes amalgamées d'éléments de $\mathcal{I}$.\\
\\
Montrons maintenant le résultat (*) suivant : soient $E\times A\rightarrow F\times
A$, $E\times B\rightarrow F\times B$ et $E\times C\rightarrow F\times C$ des
iso-cofibrations triviales, alors leur somme amalgamée $E\times (B\coprod_{A} C)\rightarrow
F\times (B\coprod_{A} C)$ aussi. Ce résultat découle du lemme~\ref{ictcub}
si l'on montre que $(F\times A\coprod_{E\times A} E\times C)\rightarrow F\times
C$ est une cofibration. Or on peut remarquer que le produit fibré de $F\times
A$ et $E\times C$
au-dessus de $F\times C$ n'est autre que $E\times A$. Donc par hypothèse sur
$\mathcal{C}$, on obtient que niveau par niveau le morphisme de la somme
amalgamée dans $F\times C$ est une cofibration. Ce qui montre le résultat (*).\\
\\
Les \precats~de la forme $\partial\Delta[n]\Theta X$ sont des sommes
amalgamées de \precats~de la forme $\Delta[m]\Theta X$. Donc comme on suppose
(II) vrai, on obtient par le résultat (*) que (I) est vraie pour les \precats~de la forme
$\Delta[m]\Theta X$ et $\partial\Delta[n]\Theta X$. En outre les sources et but des flèches de $\mathcal{I}$
sont des sommes amalgamées de \precats~de la forme $\partial\Delta[n]\Theta X$
et $\Delta[m]\Theta X$, donc toujours par résultat (*), on obtient que (I) est
vrai pour les sources et but des flèches de $\mathcal{I}$. Donc toujours par
(*), on obtient que (I) est encore vrai pour les sommes amalgamées de sources
et buts de flèches de $\mathcal{I}$. Par lemme~\ref{ictcolim2}, les
colimites séquentielles transfinies d'iso-cofibrations triviales le long de
cofibrations sont des iso-cofibrations triviales. Comme $\mathcal{I}$ est
contenu dans les cofibrations, les morphismes obtenus comme sommes amalgamées
de $E\times B$ ou de $F\times B$ par les flèches de $\mathcal{I}$ sont des
cofibrations, donc par le lemme~\ref{ictcolim2}, (I) est vrai pour les colimites séquentielles transfinies de
sommes amalgamées d'éléments de $\mathcal{I}$, donc, comme on l'a vu
ci-dessus, pour leurs rétracts et finalement pour toute \precat, car ce sont
toutes des rétracts de tels objets.\\
CQFD. 

\begin{lem}
Supposons vraies les hypothèses du théorème~\ref{catprod}.\\
Le résultat (II) est équivalent au résultat (III) suivant :\\
pour tous objets $X$ et $Y$ de $\mathcal{C}$ et pour tout entier $m$ et $n$, le morphisme
$(\Upsilon(m)\Theta X)\times(\Delta[n]\Theta Y)\rightarrow (\Delta[m]\Theta
X)\times(\Delta[n]\Theta Y)$ est une iso-cofibration triviale.
\end{lem}
{\it Preuve :}\\
Encore une fois (III) découle clairement de (II). Il suffit donc de montrer
que (III) implique (II). Supposons donc (III) vrai. Soit $E\rightarrow F$ une
flèche génératrice de type 2, elle s'écrit donc
$(\Upsilon(m)\Theta F')\coprod_{\Upsilon(m)\Theta E'}(\Delta[m]\Theta
E')\rightarrow \Delta[m]\Theta F'$, avec $E'\rightarrow F'$ une cofibration de
$\mathcal{C}$. On va appliquer le lemme~\ref{ictcub} aux données suivantes :
$A=(\Upsilon(m)\Theta E')\times (\Delta[n]\Theta Y)$, $A'=(\Delta[m]\Theta
E')\times (\Delta[n]\Theta Y)$, $B=B'=(\Delta[m]\Theta
E')\times (\Delta[n]\Theta Y)$,
 $C=(\Upsilon(m)\Theta F')\times (\Delta[n]\Theta Y)$, $C'=(\Delta[m]\Theta
F')\times (\Delta[n]\Theta Y)$. On a bien par (III) que $A\rightarrow A'$ et
$C\rightarrow C'$ sont des iso-cofibrations triviales, de même que
$B\rightarrow B'$ qui n'est autre que l'identité. On peut facilement remarquer
que $A'\times_{C'} C$ n'est autre que $A$, ce qui par hypothèse sur
$\mathcal{C}$, nous donne que $A'\coprod_{A} C\rightarrow C'$ est une
cofibration. Donc par lemme~\ref{ictcub}, on obtient que $B\coprod_{A}
C\rightarrow B'\coprod_{A'} C'$ est une
iso-cofibration triviale, i.e. que $E\times \Delta[n]\Theta Y\rightarrow F\times
\Delta[n]\Theta Y$ en est une.\\
CQFD.\\

Après toutes ces réductions, nous aboutissons au fait que pour démontrer
le théorème~\ref{catprod}, il faut et il suffit de montrer le résultat
(III). C'est donc ce que nous allons faire.\\

{\it Preuve du théorème~\ref{catprod}}\\
Les lemmes précédents nous montrent que, pour démontrer le théorème, il
faut et il suffit de démontrer le résultat (III). Tout d'abord on remarque
que $\Upsilon(m)\Theta X\rightarrow \Delta[m]\Theta X$ est une cofibration et un
isomorphisme sur les objets donc son produit par l'identité de
$\Delta[n]\Theta Y$ aussi. Montrons donc que ce produit est une \eq~faible.\\
Posons $A=(\Upsilon(m)\Theta X)\times(\Delta[n]\Theta Y)$ et $B=(\Delta[m]\Theta
X)\times(\Delta[n]\Theta Y)$. Les ensembles d'objets de $A$ et de $B$ sont
identiques et consistent en l'ensemble $Ob$ des couples $(i,j)$ avec $0\leq i\leq m$
et $0\leq j\leq n$. Soit $S$ un sous-ensemble de $Ob$, on note
$A\{S\}$ et $B\{S\}$ les sous-\precats~pleines de $A$ et de $B$ ayant $S$ pour
ensemble d'objets.\\
On considère divers sous-ensembles de $Ob$. Les rectangles
$R$ de taille $m'\times n'$ sont définis comme l'ensemble des couples $(i,j)$ tel que $0\leq i\leq m'$ et $0\leq
j\leq n'$. Un rectangle de taille $m'\times n'$ sera dit rogné de $m''$ si les couples
$(i,j)$ tel que $0\leq i\leq m''$, avec $m''<m'$ et $j=n'$ n'y sont plus. Une
queue est une suite de couples $(i_k,j_k)_{0\leq k\leq r}$ telle que
$i_{k-1}\leq i_{k} \leq i_{k+1}$ et $j_{k-1}\leq j_{k} \leq j_{k+1}$. Enfin un
rectangle à queue de taille $m'\times n'$ est la réunion d'un rectangle de
taille $m'\times n'$ avec une queue de premier terme $(m',n')$ et de dernier
terme $(m,n)$. De même on peut définir un rectangle rogné à queue.\\
Mettons sur $Ob$ l'ordre lexicographique, ce qui fait de $Ob$ un ensemble
totalement ordonné. Montrons par récurrence sur $(m',n')$ que pour tout
rectangle $R$ de taille $m'\times n'$ à queue, le morphisme $A\{R\}\rightarrow
B\{R\}$ est une \eq~faible.\\ 
\\
Le rang initial $(0,0)$ correspond aux queues allant
de $(0,0)$ à $(m,n)$. Nous allons en fait montrer le résultat pour toutes
les queues. Soit $S$ une queue, montrons d'abord que $A\{S\}$ et $B\{S\}$ sont des
\precats~librement ordonnées strictes. Tout d'abord, d'après le
lemme~\ref{FGordre}, les \precats~$\Upsilon(m)\Theta X$, $\Delta[m]\Theta X$ et
$\Delta[n]\Theta Y$ sont des \precats~librement ordonnées strictes. Il en est
donc de même pour leurs produits $A$ et $B$ par le lemme~\ref{stordre}, car on a munit
$Ob$ de l'ordre lexicographique. Comme $A\{S\}$ et $B\{S\}$ sont des
sous-\precats~pleines des \precats~librement ordonnées strictes $A$ et $B$ et
que $S$ est muni de l'ordre induit par celui de $Ob$, alors ce sont aussi des
\precats~librement ordonnées strictes. On remarque en outre que pour $x<y$
adjacents dans $S$, par définition de $S$, ceci signifie que si $x=(i,j)$
alors $y$ vaut soit $(i,j+1)$ soit $(i+1,j)$ soit $(i+1,j+1)$. Dans le premier
cas, $A_1(x,y)$ et $B_1(x,y)$ valent tous les deux $*\times Y$, dans le
deuxième $X\times *$ et dans le troisième $X\times Y$, car
$\Upsilon(m)\Theta X$ et $\Delta[m]\Theta X$ ont les mêmes composantes sur les
morphismes principaux. Ceci montre que, pour tout couple ordonné $(x,y)$ d'objets
adjacents de $S$, on a que $A\{S\}_1(x,y)$ et $B\{S\}_1(x,y)$ sont égaux. On
peut alors appliquer la proposition~\ref{eqlord} au morphisme de \precats~librement
ordonnées strictes $A\{S\}\rightarrow B\{S\}$ qui est l'identité sur l'ensemble des
objets $S$, ce qui entraîne que $A\{S\}\rightarrow B\{S\}$ est une
\eq~faible et donc montre le rang $(0,0)$ de la récurrence.\\
\\
Supposons montrée l'hypothèse de récurrence pour les rectangles à queue
de taille $m''\times n''$, avec $(m'',n'')$ inférieur à $(m',n')$.
Remarquons tout d'abord que si $(m',n')$ est de la forme $(0,n')$ ou $(m',0)$,
les rectangles à queue correspondant sont des queues allant de $(0,0)$ à
$(m,n)$, et donc l'hypothèse de récurrence est vraie par le cas $(0,0)$. Il
faut donc montrer l'hypothèse de récurrence pour $(m',n')$ avec $m'>0$ et
$n'>0$. Pour cela, nous allons nous servir des rectangles rognés.
On remarque en effet que les rectangles à queue de taille $m'\times n'$
rognés de $m'-1$ ne sont autres que les rectangles à queue de taille
$m'\times (n'-1)$ pour lesquels l'hypothèse de récurrence est vérifiée. On
remarque en outre que les rectangles à queue de taille $m'\times n'$ rognés de $0$
ne sont autres que les rectangles à queue de taille $m'\times n'$ pour
lesquels on veut montrer que l'hypothèse de récurrence est vérifiée. Pour cela, on
va montrer par récurrence sur $1\leq i\leq m'$ que
l'hypothèse est vraie pour les rectangles à queue rognés de $m'-i$ de taille
$m'\times n'$. Le cran initial $1$ a déjà été vu. Supposons
l'hypothèse vraie pour une rognure de $m'-i$ et montrons-la pour une rognure
de $m'-(i+1)$. Notons $S$ un rectangle à queue rogné de $m'-i-1$. Notons $x$
l'extrémité de la rognure, i.e. le point $(m'-i-1,n')$. Notons d'abord
$\bar{S}$ le sous-rectangle à queue de $S$ de taille $m'\times n'$ rogné de
$m'-i$. Notons $S^x$ le
sous-rectangle à queue de $S$ de taille $(m'-i-1)\times (n'-1)$ dont la queue passe
par $x$ et notons $\bar{S}^x$ le sous-rectangle à queue de $S$ de la même taille
mais dont la queue évite $x$ en passant directement de $(m'-i-1,n'-1)$ à
$(m'-i,n')$. On a donc que $S$ n'est autre que la somme amalgamée de $\bar{S}$ avec
$S^x$ au-dessus de $\bar{S}^x$. On en déduit que $A\{S\}$ n'est autre que la
somme amalgamée de $A\{\bar{S}\}$ avec $A\{S^x\}$ au-dessus de
$A\{\bar{S}^x\}$, de même pour $B\{S\}$. Par hypothèse de récurrence, on a
que les morphismes $A\{\bar{S}^x\}\rightarrow B\{\bar{S}^x\}$,
$A\{\bar{S}\}\rightarrow B\{\bar{S}\}$ et $A\{S^x\}\rightarrow B\{S^x\}$ sont
des iso-cofibrations triviales, en outre $A\{\bar{S}\}\times_{B\{\bar{S}\}}
B\{\bar{S}^x\}$ n'est autre
que $A\{\bar{S}^x\}$, ce qui par
hypothèse sur $\mathcal{C}$, nous donne que $A\{\bar{S}\}\coprod_{A\{\bar{S}^x\}}
B\{\bar{S}^x\}\rightarrow B\{\bar{S}\}$ est une cofibration. Donc on
peut appliquer le lemme~\ref{ictcub} et on obtient que $A\{S\}\rightarrow B\{S\}$
est une iso-cofibration triviale, ce qui montre l'hypothèse pour les rognures
de taille $m'-i-1$. Par récurrence, on obtient l'hypothèse vraie pour les
rognures nulles donc pour les rectangles à queue de taille $m'\times n'$, ce
qui montre l'hypothèse de récurrence au rang $(m',n')$. Et donc par
récurrence, on obtient que $A\{S\}\rightarrow B\{S\}$ est une \eq~faible pour les
rectangles à queue de taille $m\times n$, donc pour $S=Ob$, ce qui signifie
que $A\rightarrow B$ est une \eq~faible. Comme en outre on a vu que c'est aussi
une iso-cofibration, on a montré que c'est une iso-cofibration triviale.\\
CQFD.\\

Le théorème sur le produit des catégorisations montré, nous pouvons
maintenant en déduire la stabilité des \eqs~faibles de \precats~par produit
cartésien.  

\begin{cor}\label{eqfprod}
Supposons vraies les hypothèses du théorème~\ref{catprod}.\\
Soient $A\rightarrow A'$ et $B\rightarrow B'$ des \eqs~faibles alors leur
produit $A\times B\rightarrow A'\times B'$ est aussi une \eq~faible.
\end{cor}
{\it Preuve :}\\
Considérons le diagramme suivant :
\begin{diagram}
A\times B & \rTo & A'\times B'\\
\dTo & & \dTo\\
Cat(A)\times Cat(B) & \rTo & Cat(A')\times Cat(B')\\
\end{diagram}
Par le théorème~\ref{catprod}, les flèches verticales de ce diagramme sont
des \eqs~faibles. Donc par l'axiome "trois pour deux" dans $\mathcal{C-PC}$, montrer que $A\times B\rightarrow A'\times B'$ est une
\eq~faible équivaut à montrer que $Cat(A)\times Cat(B)\rightarrow
Cat(A')\times Cat(B')$, qu'on notera $f$, est une \eq~faible.\\ 
\\
Comme on l'a vu à la proposition~\ref{catfac}, une
\cat~facile est un objet ayant une \prd. Or
par le lemme~\ref{prst}, les objets ayant une \prd~sont stables par produit
cartésien. Donc les produits $Cat(A)\times Cat(B)$ et $Cat(A')\times Cat(B')$
sont des \cats~faciles. Ainsi $f$ est une \eq~faible si et seulement si c'est
une \eq~de \cats, par le lemme~\ref{eqceqf}. Or $f$ est le produit des morphismes $Cat(A)\rightarrow
Cat(A')$ et $Cat(B)\rightarrow Cat(B')$, qu'on notera
respectivement $f_1$ et $f_2$. Par hypothèse, $f_1$ et $f_2$ sont des \eqs~de
\cats. Comme $\tau_0$ est un foncteur préservant les produits, on a que
$\tau_0(f)=\tau_0(f_1)\times\tau_0(f_2)$. Or $\tau_0(f_1)$ et $\tau_0(f_2)$ sont
des isomorphismes car $f_1$ et $f_2$ sont des \eqs~de \cats, donc $\tau_0(f)$ est un
isomorphisme, i.e. $f$ est essentiellement surjective. De même $f(x,y)$ est le
produit de $f_1(x,y)$ et de $f_2(x,y)$ qui sont des \eqs~faibles de
$\mathcal{C}$ car $f_1$ et $f_2$ sont des \eqs~de \cats. Or par hypothèse sur
$\mathcal{C}$, le produit d'\eqs~faibles est une \eq~faible, donc $f(x,y)$ est
une \eq~faible, i.e. $f$ est pleinement fidèle. On a donc montré que $f$ est
une \eq~de \cats.\\
CQFD.\\

Ce résultat sur la stabilité des \eqs~faibles par produit cartésien est le
dernier ingrédient dont on avait besoin pour pouvoir démontrer la
stabilité des cofibrations triviales par somme amalgamée le long d'un
morphisme et par colimite séquentielle transfinie, ce que nous allons faire
dans la prochaine section.

\newpage 

\section{Stabilité des cofibrations triviales}

Afin de montrer l'hypothèse 3) du lemme de reconnaissance~\ref{reco} de la
structure de \cmfcof, nous devons montrer la stabilité des cofibrations
triviales par somme amalgamée le long d'un morphisme et par colimite
séquentielle transfinie. Par fonctorialité de Cat, cela revient à le montrer
pour les cofibrations triviales de \cats. Or nous avons déjà un tel résultat pour les
iso-cofibrations triviales. L'idée pour montrer que les sommes amalgamées le
long d'un morphisme des cofibrations triviales sont des cofibrations triviales
est de faire une récurrence transfinie sur la différence entre le nombre
d'objets du but de la cofibration triviale et celui de sa source. Comme le
résultat sur les iso-cofibrations triviales correspond au cas 0, cas de
départ de la récurrence transfinie, nous allons tout d'abord chercher à
montrer le pas de la récurrence qui fait passer d'un ordinal à son
successeur.

\subsection{Stabilité par somme amalgamée de certaines cofibrations
triviales}

Nous allons donc tenter de montrer que toute cofibration triviale, dont le but a
un objet de plus que la source, est stable par somme amalgamée. Comme ce
résultat va utiliser le théorème~\ref{catprod}, nous allons énoncer ses
hypothèses.

\begin{hyp}\label{coftriv1}
Soit $(\mathcal{C},\mathcal{F}_1,\mathcal{F}_2)$ une donnée de Segal facile
vérifiant les propriétés suivantes :
\item 1) la \cmf~sous-jacente $\mathcal{C}$ est une \cmf~simpliciale engendrée
par monomorphismes,
\item 2) les cofibrations génératrices de la \cmfcof~$\mathcal{C}$ ont leurs
sources et buts $\alpha$-petits,
\item 3) les cofibrations génératrices de $\mathcal{C}$ ont leurs buts connexes,
\item 4) les réunions de deux
sous-objets sont exactement les sommes amalgamées de ces sous-objets au-dessus
de leur intersection,
\item 5) pour tout monomorphisme $j:A\rightarrow B$ de $\mathcal{C}$, la classe contenant un représentant de chaque classe d'isomorphismes d'objets $C$ par lequel $j$ se
factorise dans la sous-catégorie des monomorphismes est un ensemble,
\item 6) les colimites séquentielles transfinies de la sous-catégorie des
monomorphismes de $\mathcal{C}$ existent et sont les colimites séquentielles
transfinies de $\mathcal{C}$,
\item 7) il n'existe aucun morphisme de but l'ensemble vide qui ne soit pas un
isomorphisme.
\end{hyp}

Le fait qu'une cofibration triviale de \cats~ait son but possédant un objet de plus que
sa source se produit en particulier si l'objet du but qui n'est pas dans l'image
de la source est isomorphe à un objet dans l'image de la source. C'est en
effet un cas particulier d'essentielle surjectivité. Un tel couple d'objet est
représenté par l'intervalle des ensembles simpliciaux. Nous allons donc
commencer par montrer que le morphisme naturel d'une \cat~vers sa somme
amalgamée par l'intervalle des ensembles simpliciaux est bien une cofibration
triviale.

\begin{lem}
Supposons vraies les hypothèses~\ref{coftriv1}.\\
Soit $\bar{I}$ le groupoïde constitué de deux objets notés 0 et 1 et d'un
isomorphisme entre eux noté $u$. Pour toute \cat~$A$ et pour tout objet $a$ de $A$, le
morphisme $A\rightarrow A\coprod_{0} \bar{I}$, induit par l'identification de $a\in A$ et de $0\in \bar{I}$, est une cofibration triviale.
\end{lem}
{\it Preuve :}\\
Tout d'abord on remarque que l'objet $a$ de $A$ définit un unique morphisme de {\bf 0} vers $A$, de
même l'objet $0$ de $\bar{I}$ définit un unique morphisme de {\bf 0} vers
$\bar{I}$. Notons alors $B$ la somme amalgamée de $A$ avec $\bar{I}$ au-dessus
de ces morphismes et notons $i$, respectivement $j$, le morphisme canonique de
$\bar{I}$, respectivement de $A$, vers $B$. Ensuite on remarque que $A\rightarrow A\coprod_{0} \bar{I}$ est une
cofibration comme somme amalgamée de la cofibration $0\rightarrow \bar{I}$ le
long du morphisme $0\rightarrow A$.\\
\\
Posons $h:\bar{I}\times\bar{I}\rightarrow \bar{I}\times\bar{I}$ le foncteur
entre groupoïdes défini ainsi : $h(0,0)=(0,0)$, $h(0,u)=(0,u)$,
$h(0,u^{-1})=(0,u^{-1})$, $h(0,1)=(0,1)$, $h(u,0)=(u,0)$, $h(u,u)=(0,u)$,
$h(u,u^{-1})=(u,u^{-1})$, $h(u,1)=(0,1)$, $h(u^{-1},0)=(u^{-1},0)$,
$h(u^{-1},u)=(u^{-1},u)$, $h(u^{-1},u^{-1})=(0,u^{-1})$, $h(u^{-1},1)=(0,1)$,
$h(1,0)=(1,0)$, $h(1,u)=(u^{-1},u)$, $h(1,u^{-1})=(u,u^{-1})$, $h(1,1)=(0,1)$.
On remarque que $B\times\bar{I}$ n'est autre que la somme
amalgamée de $\bar{I}\times\bar{I}$ et de $A\times\bar{I}$ au-dessus de
$0\times\bar{I}$. Considérons le diagramme suivant :
\begin{diagram}
0\times\bar{I} & \rTo^{a\times Id} & A\times\bar{I} & & \\
\dTo^{0\times Id} & & \dTo^{j\times Id} & \rdTo(2,4)^{j\times Id} & \\
\bar{I}\times\bar{I} & \rTo^{i\times Id} & B\times\bar{I} & & \\
 & \rdTo(4,2)_{(i\times Id)\circ h} & & \rdTo~{\exists ! f} & \\
 & & & & B\times\bar{I}\\
\end{diagram}
Le diagramme étant commutatif, par universalité de la somme amalgamée
$B\times\bar{I}$, il existe un unique morphisme $f:B\times\bar{I}\rightarrow
B\times\bar{I}$ tel que $f_{|\bar{I}\times\bar{I}}$ soit la composée $(i\times
Id)\circ h$ et $f_{|A\times\bar{I}}$ soit le morphisme $j\times Id$. En
identifiant $B$ avec $B\times 0$, notons $i_0$, respectivement $i_1$, les
inclusions $Id\times 0:B\rightarrow B\times\bar{I}$ et $Id\times 1:B\rightarrow
B\times\bar{I}$. On calcule aisément que $f\circ i_0$ n'est autre que $i_0$
mais que $f\circ i_1$ est le morphisme $(j\circ r)\times 1$, où $r$ est la
rétraction de $B$ dans $A$ définie ainsi. Considérons le diagramme suivant :
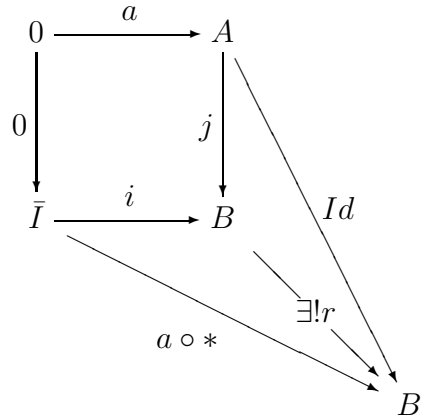
\begin{diagram}
0 & \rTo^{a} & A & & \\
\dTo^{0} & & \dTo^j & \rdTo(2,4)^{Id} & \\
\bar{I} & \rTo^{i} & B & & \\
 & \rdTo(4,2)_{a\circ *} & & \rdTo~{\exists ! r} & \\
 & & & & B\\
\end{diagram}
où $*$ est l'unique morphisme de $\bar{I}$ vers {\bf 0}.
Le diagramme étant commutatif, par universalité de la somme amalgamée $B$,
il existe un unique morphisme $r:B\rightarrow A$ tel que $r\circ j=Id_A$ et
$r\circ i=a\circ *$.\\
\\
Posons maintenant $g:B\times\bar{I}\rightarrow B$ comme étant la composée de $f$ avec la première
projection. On obtient facilement que $g\circ i_0=Id_B$ et que $g\circ
i_1=j\circ r$. Or les morphismes $0:0\rightarrow\bar{I}$ et $1:0\rightarrow\bar{I}$ sont
évidemment des \eqs~de \cats. Par le corollaire~\ref{eqfprod}, il vient que
$i_0=Id\times 0$ et $i_1=Id\times 1$ sont aussi des \eqs~faibles. Comme on a
l'égalité $g\circ i_0=Id_B$, par l'axiome "trois pour deux" dans
$\mathcal{C-PC}$, le morphisme
$g$ est aussi une \eq~faible. Toujours par cet axiome, on obtient que $j\circ
r=g\circ i_1$ est une \eq~faible. Ainsi les morphismes $j$ et $r$ ont la
propriété suivante : $r\circ j=Id$ et $j\circ r$ est une \eq~faible, donc
par le lemme~\ref{3pour2}, $j$ et $r$ sont des \eqs~faibles.\\
CQFD.\\

Bien évidemment ce cas n'est pas le cas général pour l'essentielle
surjectivité. Pour représenter bien comme il faut l'essentielle surjectivité
d'un morphisme de \cats, il faudrait un intervalle, qu'on notera $\bar{J}$, qui
représentât un couple d'objets équivalents dans une catégorie. En supposant
que l'on possède un tel intervalle $\bar{J}$, montrons alors que le morphisme naturel d'une \cat~vers sa somme amalgamée avec cet intervalle $\bar{J}$ est une
cofibration triviale. 

\begin{lem}
Supposons vraies les hypothèses~\ref{coftriv1} et supposons
qu'il existe une \cat~$\bar{J}$ munie de deux objets notés $0$ et
$1$ et ayant les propriétés suivantes :
\item - pour
toute \cat~facile $A$ et pour tout couple d'objets équivalents $a$ et $b$ de $A$, il
existe un morphisme de $\bar{J}$ vers $A$ envoyant $0$ sur $a$ et $1$ sur $b$,
\item - pour toute \cat~$A$ et pour tout morphisme de $\bar{J}$ vers $A$,
les images de $0$ et $1$ par ce morphisme sont des objets équivalents dans
$A$,
\item - notons $\bar{L}$, la sous-\cat~pleine de $\bar{J}$ d'objet $0$, alors il
existe un morphisme $p:\bar{J}\rightarrow \bar{L}$ envoyant $0$ et $1$ sur $0$
qui est une \eq~de \cats.\\
\\
Pour toute \cat~$A$ et pour tout morphisme de $\bar{L}$ vers $A$, le
morphisme $A\rightarrow A\coprod_{\bar{L}}\bar{J}$ est une cofibration triviale.
\end{lem}
{\it Preuve :}\\
Tout d'abord on remarque que $A\rightarrow A\coprod_{\bar{L}}\bar{J}$ est une
cofibration comme somme amalgamée de la cofibration $\bar{L}\subset\bar{J}$ le
long de $\bar{L}\rightarrow A$.\\
Comme $\bar{I}$ est une \cat~facile dans laquelle $0$ et $1$ sont isomorphes, il existe
un morphisme $q:\bar{J}\rightarrow\bar{I}$ envoyant $0$ sur $0$ et $1$ sur $1$.
Il existe donc un morphisme $r:\bar{J}\rightarrow\bar{L}\times\bar{I}$ induit
par $p$ et $q$. Par corollaire~\ref{eqfprod}, le morphisme produit $Id\times
0:\bar{L}\times 0\rightarrow\bar{l}\times\bar{I}$ est une \eq~faible. Par la
suite, on confondra $\bar{L}\times 0$ et $\bar{L}$. Comme on a l'égalité
$pr_1\circ (Id\times 0)=Id$, par l'axiome "trois pour deux" dans $\mathcal{C-PC}$, la première
projection $pr_1$ est une \eq~faible. En outre comme on a $pr_1\circ r=p$ et que
par hypothèse $p$ est une \eq~faible, par l'axiome "trois pour deux", il vient
que $r$ est une \eq~faible. En outre il est facile de voir que $r$ est un
isomorphisme sur les objets, donc $r$ est une iso-\eq~de \cats.\\
Posons $B$ comme étant la somme amalgamée de $A$ avec $\bar{L}\times\bar{I}$
au-dessus de $\bar{L}$ et $C$ la somme amalgamée de $A$ avec $\bar{J}$
au-dessus de $\bar{L}$. On applique le lemme~\ref{ictcub} aux iso-\eqs~faibles
$Id_{\bar{L}}$, $Id_A$ et $r$ qui en vérifient les hypothèses et on obtient
que $C\rightarrow B$ est une iso-\eq~faible.\\
\\
Considérons maintenant le diagramme  commutatif suivant :
\begin{diagram}
0 & \rTo^0 & \bar{I} & &\\
\dTo & & \dTo & \rdTo(2,4)^{(0\circ *)\times Id} & \\
\bar{L} & \rTo & \bar{L}\coprod_{0}\bar{I} & & \\
  & \rdTo(4,2)_{Id\times 0} & & \rdTo~{\exists !} & \\
  & & & & \bar{L}\times\bar{I} \\
\end{diagram}
Par le lemme précédent, le morphisme $\bar{L}\rightarrow
\bar{L}\coprod_{0}\bar{I}$ est une \eq~faible et par le
corollaire~\ref{eqfprod}, le morphisme $Id\times 0$ aussi. Par l'axiome "trois
pour deux" pour $\mathcal{C-PC}$, il vient que le morphisme $\bar{L}\coprod_{0}\bar{I}\rightarrow
\bar{L}\times\bar{I}$ est une \eq~faible. En outre, ce morphisme est aussi une
iso-cofibration, donc c'est une iso-cofibration triviale. Or sa somme
amalgamée le long de $\bar{L}\coprod_{0}\bar{I}\rightarrow
A\coprod_{0}\bar{I}$ n'est autre que le morphisme
$A\coprod_{0}\bar{I}\rightarrow B$ qui, par le lemme~\ref{ictpush}, est aussi
une iso-cofibration triviale. Or par le lemme précédent, le morphisme
$A\rightarrow A\coprod_{0}\bar{I}$ est une cofibration triviale, donc par "trois
pour deux" dans $\mathcal{C-PC}$, la composée $A\rightarrow B$ aussi. Enfin on remarque que
$A\rightarrow B$ est la composée de $A\rightarrow C$ avec l'\eq~faible
$C\rightarrow B$, ce qui par "trois pour deux" nous montre que $A\rightarrow C$
est une \eq~faible.\\
CQFD.\\

Puisque la somme amalgamée d'une \cat~par l'intervalle $\bar{J}$ revient à
rajouter un objet à la \cat~qui soit équivalent à l'un de ses objets et
qu'en outre cette opération est une cofibration triviale, il en découle que
prendre la somme amalgamée d'une \cat~quelconque par une cofibration triviale
de \cats~dont le but a un objet de plus que la source va encore donner une
cofibration triviale.

\begin{lem}\label{sa1}
Supposons vraies les hypothèses~\ref{coftriv1} et supposons
qu'il existe une \cat~$\bar{J}$ munie de deux objets notés $0$ et
$1$ et ayant les propriétés suivantes :
\item - pour
toute \cat~facile $A$ et pour tout couple d'objets équivalents $a$ et $b$ de $A$, il
existe un morphisme de $\bar{J}$ vers $A$ envoyant $0$ sur $a$ et $1$ sur $b$,
\item - pour toute \cat~$A$ et pour tout morphisme de $\bar{J}$ vers $A$,
les images de $0$ et $1$ par ce morphisme sont des objets équivalents dans
$A$,
\item - notons $\bar{L}$, la sous-\cat~pleine de $\bar{J}$ d'objet $0$, alors il
existe un morphisme $p:\bar{J}\rightarrow \bar{L}$ envoyant $0$ et $1$ sur $0$ est
une \eq~de \cats.\\
\\
Soit $A\rightarrow B$ une cofibration triviale entre \cats~faciles telle que $B$ ait
exactement un objet de plus que $A$, alors pour toute \cat~$C$ et pour tout
morphisme de $A$ vers $C$, le morphisme $C\rightarrow C\coprod_{A} B$ est une
cofibration triviale.
\end{lem}
{\it Preuve :}\\
Soit $A'$ la sous-\cat~pleine de $B$ ayant pour objets les images de ceux de $A$
par le morphisme $A\rightarrow B$. Ainsi on obtient que $A\rightarrow A'$ est
une iso-cofibration triviale. Par l'axiome "trois pour deux" dans $\mathcal{C-PC}$, on obtient que
$A'\subset B$ est aussi une cofibration triviale avec $B$ ayant un objet de
plus que $A'$. Or $C\coprod_{A} B$ est aussi la somme amalgamée de
$C\coprod_{A} A'$ et de $B$ au-dessus de $A'$. Si on montre le lemme pour
l'inclusion $A'\subset B$, on obtient que $C\coprod_{A} A'\rightarrow
C\coprod_{A} B$ est une \eq~faible. Or par le lemme~\ref{ictpush}, $C\rightarrow
C\coprod_{A} A'$ est aussi une \eq~faible car $A\rightarrow A'$ est une
iso-cofibration triviale. Et par "trois pour deux", leur composée $C\rightarrow
C\coprod_{A} B$ est une \eq~faible, ce qu'il fallait démontrer.\\
\\
Montrons donc le lemme pour l'inclusion $A'\subset B$. Soit $A'\rightarrow D$ un
morphisme quelconque. On remarque tout d'abord que $D\rightarrow D\coprod_{A'} B$ est une cofibration
comme somme amalgamée de la cofibration $A'\subset B$ le long d'un
morphisme.\\
Notons $b$ l'objet de
$B_0\setminus A'_0$. Comme $A'\subset B$ est une \eq~de \cats, il existe un objet
$a$ de $A'$ tel que $a$ soit équivalent à $b$ dans $B$. Donc, comme $B$ est
une \cat~facile, il existe un
morphisme de $\bar{J}$ vers $B$ envoyant $0$ et $1$ respectivement sur $a$ et
$b$. La sous-\cat~pleine $\bar{L}$ de $\bar{J}$ s'envoie par ce morphisme dans
$A'$ car $A'$ est pleine et contient $a$ l'image de $0$. Soit $E$ la somme
amalgamée $A'\coprod_{\bar{L}}\bar{J}$. Par propriété universelle de $E$,
il existe un morphisme de $E$ vers $B$ factorisant celui de $\bar{J}$ vers $B$
car ce dernier envoie $\bar{L}$ dans $A'$. Or par le lemme précédent,
$A'\rightarrow E$ est une cofibration triviale et comme $A'\subset B$ est une
\eq~faible, par "trois pour deux", $E\rightarrow B$ aussi. En outre on peut
remarquer que c'est un isomorphisme sur les objets. Comme $B$ est une \cat,
l'iso-\eq~faible $E\rightarrow B$ se factorise par $Cat(E)$. Par la
proposition~\ref{boncat}, le morphisme canonique $E\rightarrow
Cat(E)$ est une iso-cofibration triviale. Et donc par "trois pour deux",
$Cat(E)\rightarrow B$ est une iso-\eq~faible entre \cats. Par la
proposition~\ref{isoeqc}, c'est niveau par niveau une \eq~faible. Soit $m$ un
entier. En appliquant dans la \cmf~$\mathcal{C}$ propre à gauche le lemme du
cube pour les \eq~faibles à $Id_{A'_m}$, $Id_{D_m}$ et $Cat(E)_m\rightarrow B_m$,
et les cofibrations $A'_m\rightarrow Cat(E)_m$ et $A'_m\rightarrow B'_m$
on obtient que $(D\coprod_{A'} Cat(E))_m \rightarrow (D\coprod_{A'} B)_m$ est une
\eq~faible. Comme en outre $Cat(E)\rightarrow B$ est un isomorphisme sur les objets
alors $D\coprod_{A'} Cat(E) \rightarrow D\coprod_{A'} B$ aussi et donc c'est une
iso-\eq~faible, par la proposition~\ref{eqfniv}. Par le lemme précédent, le morphisme $D\rightarrow
D\coprod_{\bar{L}}\bar{J}$ est une cofibration triviale, or
$D\coprod_{\bar{L}}\bar{J}$ n'est autre que $D\coprod_{A'} E$. En outre, le
morphisme $D\coprod_{A'} E\rightarrow D\coprod_{A'} Cat(E)$ est une
iso-cofibration triviale, comme somme amalgamée le long d'un morphisme de
l'iso-cofibration triviale $E\rightarrow Cat(E)$ par le lemme~\ref{ictpush}. Donc par "trois
pour deux", les composées successives $D\rightarrow D\coprod_{A'} Cat(E)$ puis
$D\rightarrow D\coprod_{A'} B$ sont des \eqs~faibles, ce
qui montre le lemme pour l'inclusion $A'\subset B$.\\
CQFD.\\

Ce résultat montre donc la stabilité par somme amalgamée le long d'un
morphisme des cofibrations triviales dont le but a un objet de plus que la
source, ce que l'on peut considérer comme le pas de la récurrence transfinie
permettant de passer d'un ordinal à son successeur. Nous devons donc
maintenant regarder ce qu'il se passe pour les ordinaux limites, ce qui va
correspondre à montrer la stabilité par colimite séquentielle transfinie
des cofibrations triviales.

\subsection{Stabilité des cofibrations triviales par colimite\\ séquentielle
transfinie}

L'argument clef pour montrer la stabilité par colimite séquentielle
transfinie des cofibrations triviales de \cats~est le fait que les colimites
séquentielles transfinies de \cats~le long de cofibrations restent des \cats.

\begin{lem}
Soit $\mathcal{C}$ une donnée de Segal telle que :
\item - la catégorie sous-jacente $\mathcal{C}$ est une \cmf~simpliciale
propre à gauche,
\item - les colimites séquentielles transfinies
d'\obcs~le long de cofibrations sont des \obcs,
\item - les \eqcs~d'\obcs~sont exactement les \eqs~faibles de la
\cmf~$\mathcal{C}$ entre \obcs.\\
\\
Alors les colimites séquentielles transfinies de \cats~le long de cofibrations
niveau par niveau sont des \cats.
\end{lem}
{\it Preuve :}
Soit $(A^{\alpha})_{\alpha<\lambda}$ une $\lambda$-séquence de
\cats~telle que pour tout $\alpha$ le morphisme $A^{\alpha}\rightarrow
A^{\alpha+1}$ soit une cofibration. Notons $A^{\lambda}$ la colimite
séquentielle transfinie des
$(A^{\alpha})_{\alpha<\lambda}$. On remarque tout d'abord que pour tout entier
$m$, $A^{\lambda}_m$ est la colimite séquentielle transfinie des
\obcs~$A^{\alpha}_m$ où les morphismes $A^{\alpha}_m\rightarrow
A^{\alpha+1}_m$ sont des cofibrations. Par hypothèse sur $\mathcal{C}$, ceci
nous donne que, pour tout $m$, $A^{\lambda}_m$ est un \obc. Par stabilité des
\obcs~par produit fibré au-dessus des objets discrets,
les $A^{\lambda}_1\times\ldots\times A^{\lambda}_1$ sont aussi des \obcs, mais
également les colimites des $A^{\alpha}_1\times\ldots\times A^{\alpha}_1$. On
obtient donc le diagramme commutatif suivant :
\begin{diagram}
A^0_m & \rTo & A^1_m & \rTo & \ldots & \rTo & A^{\lambda}_m \\
\dTo & & \dTo & & & & \dTo \\
A^0_1\times\ldots\times A^0_1 & \rTo & A^1_1\times\ldots\times A^1_1 & \rTo &
\ldots & \rTo & A^{\lambda}_1\times\ldots\times A^{\lambda}_1 \\
\end{diagram}
où les flèches horizontales sont des cofibrations et les flèches
verticales sont les morphismes de Segal. Comme pour tout $\alpha<\lambda$, les
$A^{\alpha}$ sont des \cats, les morphismes de Segal sont des \eqcs~d'\obcs,
donc des \eqs~faibles. Et
comme $\mathcal{C}$ est une \cmf~simpliciale propre à gauche, par la proposition~\ref{eqcolim}, la colimite du diagramme ci-dessus est aussi
une \eq~faible. Donc
les morphismes de Segal de $A^{\lambda}$ sont des \eqs~faibles entre \obcs, i.e
des \eqcs~d'\obcs, et donc
$A^{\lambda}$ est bien une \cat.\\
CQFD.\\

Cette propriété assez inattendue des \cats~va nous permettre de démontrer
assez aisément la stabilité des cofibrations triviales par colimite
séquentielle transfinie.

\begin{prop}\label{ctcol}
Soit $(\mathcal{C},\mathcal{F}_1,\mathcal{F}_2)$ une donnée de Segal facile
vérifiant les propriétés suivantes :
\item 1) la \cmf~sous-jacente $\mathcal{C}$ est une \cmf~simpliciale engendrée
par monomorphismes,
\item 2) les colimites séquentielles transfinies
d'\obcs~le long de cofibrations sont des \obcs.\\
\\
Les cofibrations triviales sont stables par colimite séquentielle transfinie.
\end{prop}
{\it Preuve :}
Soit $(A^{\alpha})_{\alpha<\lambda}$ une $\lambda$-séquence
telle que pour tout $\alpha$ le morphisme $A^{\alpha}\rightarrow
A^{\alpha+1}$ soit une cofibration triviale. Notons $A^{\lambda}$ la colimite
séquentielle transfinie des
$(A^{\alpha})_{\alpha<\lambda}$. On remarque tout d'abord que $A^0\rightarrow
A^{\lambda}$ est une cofibration comme colimite séquentielle transfinie de
cofibrations. Par le lemme pratique~\ref{pratique2}, pour
montrer que $A^0\rightarrow A^{\lambda}$ est une \eq~faible, il faut et il
suffit de montrer que $Cat(A^0)\rightarrow \colimite{\alpha<\lambda}
Cat(A^{\alpha})$ est une \eq~faible. D'après le lemme~\ref{cofcat}, les
morphismes $Cat(A^{\alpha})\rightarrow Cat(A^{\alpha+1})$ sont des cofibrations.
On peut donc appliquer le lemme précédent aux $Cat(A^{\alpha})$, ce qui nous
donne que $\colimite{\alpha<\lambda} Cat(A^{\alpha})$ est une \cat. Notons cette
\cat~$B$. Ainsi on doit montrer que $Cat(A^0)\rightarrow B$ est une \eq~faible
entre \cats, i.e. est une \eq~de \cats.\\ 
\\
Comme $\tau_0(Cat(A^0)\rightarrow B)$ est la colimite des
$\tau_0(Cat(A^{\alpha})\rightarrow Cat(A^{\alpha+1}))$ qui sont des bijections,
car les $Cat(A^{\alpha})\rightarrow Cat(A^{\alpha+1})$ sont des \eqs~de \cats,
et que les bijections sont stables par colimite séquentielle transfinie, il vient que
$\tau_0(Cat(A^0)\rightarrow B)$ est une bijection et donc que
$Cat(A^0)\rightarrow B$ est essentiellement surjective.\\
\\
Soit $(x,y)$ un couple d'objets de $Cat(A^0)$. Comme les morphismes $Cat(A^{\alpha})\rightarrow Cat(A^{\alpha+1})$ sont des
injections sur les objets, car ce sont des cofibrations, on notera encore
$(x,y)$ l'image du couple $(x,y)$ dans les $Cat(A^{\alpha})$. Or le morphisme $Cat(A^0)_1(x,y)\rightarrow
B_1(x,y)$ n'est autre que la colimite séquentielle transfinie des
$Cat(A^{\alpha})_1(x,y)\rightarrow Cat(A^{\alpha+1})(x,y)$ qui sont des
cofibrations triviales car les $Cat(A^{\alpha})\rightarrow Cat(A^{\alpha+1})$
sont des cofibrations et des \eqs~de \cats. Comme $\mathcal{C}$ est une \cmf,
les cofibrations triviales sont stables par colimite séquentielle transfinie
et donc $Cat(A^0)_1(x,y)\rightarrow B_1(x,y)$ est une cofibration triviale, ce
qui montre que $Cat(A^0)\rightarrow B$ est pleinement fidèle.\\
CQFD.\\

Maintenant que nous avons montré le pas de la récurrence transfinie pour
passer d'un ordinal à son successeur et aussi celui pour les ordinaux limites,
il ne nous reste plus qu'à conclure la récurrence transfinie, autrement dit
à prouver la stabilité par somme amalgamée le long d'un morphisme des
cofibrations triviales.

\subsection{Stabilité des cofibrations triviales par somme amalgamée}

Comme de coutume, nous allons tout d'abord montrer la stabilité par somme
amalgamée pour les cofibrations triviales de \cats~de laquelle va découler celle
pour les \precats. Comme elle provient des deux pas de récurrence déjà
montrés, nous allons recueillir les hypothèses faisant marcher ces deux pas.

\begin{hyp}\label{coftriv}
Soit $(\mathcal{C},\mathcal{F}_1,\mathcal{F}_2)$ une donnée de Segal facile
vérifiant les propriétés suivantes :
\item 1) la \cmf~sous-jacente $\mathcal{C}$ est une \cmf~simpliciale engendrée
par monomorphismes,
\item 2) les cofibrations génératrices de la \cmfcof~$\mathcal{C}$ ont leurs
sources et buts $\alpha$-petits,
\item 3) les cofibrations génératrices de $\mathcal{C}$ ont leurs buts connexes,
\item 4) les réunions de deux
sous-objets sont exactement les sommes amalgamées de ces sous-objets au-dessus
de leur intersection,
\item 5) pour tout monomorphisme $j:A\rightarrow B$ de $\mathcal{C}$, la classe contenant un représentant de chaque classe d'isomorphismes d'objets $C$ par lequel $j$ se
factorise dans la sous-catégorie des monomorphismes est un ensemble,
\item 6) les colimites séquentielles transfinies de la sous-catégorie des
monomorphismes de $\mathcal{C}$ existent et sont les colimites séquentielles
transfinies de $\mathcal{C}$,
\item 7) il n'existe aucun morphisme de but l'ensemble vide qui ne soit pas un
isomorphisme,
\item 8) les colimites séquentielles transfinies
d'\obcs~le long de cofibrations sont des \obcs,
\item 9) il existe une \cat~$\bar{J}$ munie de deux objets notés $0$ et
$1$ et ayant les propriétés suivantes :
\item - pour
toute \cat~facile $A$ et pour tout couple d'objets équivalents $a$ et $b$ de $A$, il
existe un morphisme de $\bar{J}$ vers $A$ envoyant $0$ sur $a$ et $1$ sur $b$,
\item - pour toute \cat~$A$ et pour tout morphisme de $\bar{J}$ vers $A$,
les images de $0$ et $1$ par ce morphisme sont des objets équivalents dans
$A$,
\item - notons $\bar{L}$, la sous-\cat~pleine de $\bar{J}$ d'objet $0$, alors il
existe un morphisme $p:\bar{J}\rightarrow \bar{L}$ envoyant $0$ et $1$ sur $0$ qui soit
une \eq~de \cats.
\end{hyp}

\begin{lem}
Supposons vraies les hypothèses~\ref{coftriv}.\\
Les cofibrations triviales entre \cats~faciles sont stables par somme amalgamée le long d'un
morphisme de \cats~quelconque.
\end{lem}
{\it Preuve :}\\
Soit $\lambda$ un cardinal transfini.
Soit $A\rightarrow B$ une cofibration triviale de \cats~faciles telle que le cardinal de l'ensemble
$B_0\setminus A_0$ est inférieur à $\lambda$ et soit $A\rightarrow C$ un
morphisme de \cats~quelconque. Montrons par récurrence transfinie sur
$\beta\leq \lambda$ que
$C\rightarrow C\coprod_A B$ est une cofibration triviale. Le cas $\beta=0$ est
celui des iso-cofibrations triviales, cas traité par le lemme~\ref{ictpush}.
Supposons la propriété vraie au rang $\beta$ et montrons-la pour $\beta+1$.
Soit $B'_0$ un sous-ensemble de $B_0$ contenant $A_0$ et ayant un objet de moins
que $B_0$ et soit $B'$
la sous-\cat~pleine de $B$ d'objets $B'_0$, en tant que sous-catégorie pleine
d'une \cat~facile c'est encore une \cat~facile. Ainsi le morphisme $A\rightarrow
B'$ est une cofibration triviale de \cats~faciles telle que $Card(B'_0\setminus
A_0)=\beta$ et donc par hypothèse de récurrence, $C\rightarrow C\coprod_A
B'$ est une cofibration triviale. En outre $B'\subset B$ est une cofibration
triviale de \cats~faciles telle que $Card(B_0\setminus B'_0)=1$, donc par le
lemme~\ref{sa1}, pour tout morphisme de \cats~$B'\rightarrow D$, le morphisme
$D\rightarrow D\coprod_{B'} B$ est une cofibration triviale. C'est vrai en
particulier pour $D=C\coprod_A B'$, ce qui donne que $C\coprod_A B'\rightarrow
C\coprod_A B$ est une cofibration triviale. Et par l'axiome "trois pour deux"
dans $\mathcal{C-PC}$,
la composée $C\rightarrow C\coprod_A B$ est une cofibration triviale, ce qui
montre la propriété au rang $\beta+1$.\\
\\
Supposons la propriété vraie pour tous $\gamma<\beta$, avec $\beta$ ordinal
limite, et montrons-la pour $\beta$. On peut construire par récurrence
transfinie une suite croissante de
sous-ensembles $(B^{\gamma}_0)_{\gamma\leq\beta}$ telle que $B^0_0$ n'est autre
que $A_0$, que $B^{\gamma+1}_0=B^{\gamma}_0\coprod *$, que
$B^{\gamma}_0=\bigcup_{\delta<\gamma}B^{\delta}_0$ pour $\gamma$ ordinal limite et
que $B^{\beta}_0$ soit $B_0$. Ceci donne lieu à une suite croissante $(B^{\gamma})_{\gamma\leq\beta}$ de
sous-\cats~pleines de $B$, donc de \cats~faciles. On a pour tout $\gamma<\beta$ que
$B^{\gamma}\rightarrow B^{\gamma+1}$ est une cofibration triviale de
\cats~faciles telle
que $Card(B^{\gamma+1}_0\setminus B^{\gamma}_0)=1$, donc par le lemme~\ref{sa1},
 il vient que $C\coprod_{A}B^{\gamma}\rightarrow
(C\coprod_{A}B^{\gamma})\coprod_{B^{\gamma}}
B^{\gamma+1}=C\coprod_{A}B^{\gamma+1}$ est une cofibration triviale. On obtient
ainsi que $C\coprod_A B^0\rightarrow C\coprod_A B$ est la colimite séquentielle transfinie des
cofibrations triviales\\
 $C\coprod_{A}B^{\gamma}\rightarrow
C\coprod_{A}B^{\gamma+1}$, ce qui, par la proposition~\ref{ctcol}, est bien une
cofibration triviale. Comme par propriété au rang 0, le morphisme
$C\rightarrow C\coprod_A B^0$ est une cofibration triviale car $A\rightarrow
B^0$ est une iso-cofibration triviale de \cats, la composée $C\rightarrow
C\coprod_A B$ est une cofibration triviale par "trois pour deux", ce qui montre
la propriété au rang $\beta$.\\
CQFD.\\

Maintenant que la stabilité par somme amalgamée le long d'un morphisme des
cofibrations triviales est prouvée pour les \cats, ce n'est plus qu'un jeu d'enfant de
l'avoir pour les \precats.
 
\begin{prop}\label{ctpush}
Supposons vraies les hypothèses~\ref{coftriv}.\\
Les cofibrations triviales sont stables par somme amalgamée le long d'un
morphisme quelconque.
\end{prop}
{\it Preuve :}\\
Soient $A\rightarrow B$ une cofibration triviale et $A\rightarrow C$ un
morphisme quelconque. Tout d'abord on remarque que $C\rightarrow C\coprod_A B$
est une cofibration comme somme amalgamée de la cofibration $A\rightarrow B$ le long d'un
morphisme. Il ne reste donc plus qu'à montrer que c'est une \eq~faible.
Or par le lemme pratique~\ref{pratique}, il vient que pour montrer que
$C\rightarrow C\coprod_A B$ est une \eq~faible, il faut et il suffit de montrer
que $Cat(C)\rightarrow Cat(C)\coprod_{Cat(A)} Cat(B)$ est une \eq~faible. Comme
par hypothèse, $A\rightarrow B$ est une cofibration, alors par le
lemme~\ref{cofcat}, $Cat(A)\rightarrow Cat(B)$ est une cofibration de
\cats~faciles, qui en outre
est triviale car $A\rightarrow B$ est une \eq~faible. Par le lemme précédent, on obtient que $Cat(C)\rightarrow
Cat(C)\coprod_{Cat(A)} Cat(B)$ est une \eq~faible.\\
CQFD.\\

Si l'on tire le bilan de cette section, nous avons prouvé que les cofibrations
triviales sont stables par somme amalgamée le long d'un morphisme et par
colimite séquentielle transfinie. Par ailleurs, nous avions au préalable
montrer que les cofibrations et les \eqs~faibles sont stables par rétracts.
Ceci nous donne que les rétracts de colimites séquentielles transfinies de
cofibrations triviales sont des cofibrations triviales. C'est donc en
particulier le cas pour les flèches de $\mathcal{J}$, ce qui montre que les
$\mathcal{J}$-cofibrations sont des cofibrations triviales. Or ce résultat
n'est autre que l'hypothèse 3) du lemme de reconnaissance~\ref{reco} de la
structure de \cmfcof.

\begin{cor}\label{p3}
Supposons vraies les hypothèses~\ref{coftriv}.\\
Les $\mathcal{J}$-cofibrations sont à la fois des $\mathcal{I}$-cofibrations et des
\eqs~faibles.
\end{cor}
{\it Preuve :}\\
Par définition, $\mathcal{J}$ est contenu dans l'intersection de la classe des
$\mathcal{I}$-cofibrations et de celle des \eqs~faibles. Or cette
intersection de classes est stable par rétracts, par somme
amalgamée le long d'un morphisme d'après la proposition précédente, et
par colimite séquentielle transfinie d'après la proposition~\ref{ctcol}.
Comme la classe des $\mathcal{J}$-cofibrations est constituée des rétracts
de colimites séquentielles transfinies de sommes amalgamées d'éléments
de $\mathcal{J}$, car $\mathcal{J}$ permet l'argument du petit objet, on obtient
finalement que la classe des $\mathcal{J}$-cofibrations est incluse dans
l'intersection de celle des $\mathcal{I}$-cofibrations avec celle des
\eqs~faibles.\\
CQFD.\\

Cinq sur les six hypothèses du lemme~\ref{reco} ont été montrées. Il ne
reste plus qu'à montrer la dernière, ce qui est l'objet de la section
suivante.

\newpage

\section{$\mathcal{J}$ engendre les cofibrations triviales}

La dernière hypothèse du lemme~\ref{reco} de reconnaissance de la structure
de \cmfcof, propose un choix. Soit on caractérise les morphismes
$\mathcal{I}$-injectifs, soit l'on montre que $\mathcal{J}$ engendre les
cofibrations triviales. C'est cette deuxième possibilité que nous avons
choisie. Comme nous avons montré dans la section précédente que les
$\mathcal{J}$-cofibrations sont des cofibrations triviales, nous devons donc
montrer l'autre sens. Ceci équivaut à montrer que les cofibrations
triviales ont la \prg~par rapport aux morphismes $\mathcal{J}$-injectifs. Pour
cela nous allons utiliser le même type de démonstration que lorsqu'on a
montré que $\mathcal{I}$ engendre les monomorphismes. Le seul problème qui
va se poser dans cette méthode est de montrer les hypothèses du
lemme technique~\ref{ast2}. Mais tout d'abord commençons par donner quelques
propriétés sur la catégorisation Cat qui seront utiles pour la
résolution de ce problème.

\subsection{Lemmes techniques sur Cat}

Le premier résultat concerne la préservation par la catégorisation Cat de
l'intersection des sous-objets.

\begin{lem}\label{intercat}
Soit $(\mathcal{C},\mathcal{F}_1,\mathcal{F}_2)$ une donnée de Segal facile
telle que les réunions de deux
sous-objets sont exactement les sommes amalgamées de ces sous-objets au-dessus
de leur intersection.\\
Soient $A$ et $B$ deux sous-\precats~d'une \precat~$C$. Alors on a : $Cat(A\cap
B)=Cat(A)\cap Cat(B)$
\end{lem}
{\it Preuve :}\\
Tout d'abord remarquons que Cat n'est autre que la construction $E_{\Phi}$, composée des plans
simples $e_{\Phi,1}$, avec pour ensemble $\Phi$ la famille $\mathcal{FG}_1$ des
flèches génératrices de \cats~faciles. D'autre part, comme on a supposé la donnée de Segal facile, le lemme~\ref{flgenmono} nous assure que les $\mathcal{FG}_1$-cofibrations sont des monomorphismes. Enfin avec l'hypothèse sur les réunions de sous-objets, toutes les hypothèses du
lemme~\ref{interpl} sont vérifiées, ce qui nous donne le résultat voulu.\\
CQFD.\\

Le second résultat montre que la catégorisation Cat préserve la petitesse.
Plus précisément, si l'on se donne un cardinal $\alpha$, il existe un
cardinal $\alpha'$ plus grand qu'$\alpha$ tel que la catégorisation par Cat de
toute \precat~$\alpha$-petite est $\alpha'$-petite. il est important de
remarquer que le cardinal $\alpha'$ ne dépend que d'$\alpha$ et non pas de la
\precat~$\alpha$-petite que l'on veut catégoriser.

\begin{lem}\label{petitcat}
Soit $(\mathcal{C},\mathcal{F}_1,\mathcal{F}_2)$ une donnée de Segal
pré-facile.\\
Soit $\alpha'$ un cardinal régulier strictement supérieur à
$2^{\alpha}$, où l'on note $\alpha$ un cardinal régulier strictement plus
grand que celui pour lequel les
familles $\mathcal{F}_1$ et $\mathcal{F}_2$ sont petites. Alors il existe un cardinal régulier $\alpha''$ strictement
supérieur à $\alpha'$ tel que, pour toute
\precat~$\alpha'$-petite $A$, $Cat(A)$ est $\alpha''$-petite.
\end{lem}
{\it Preuve :}\\
Comme dans la démonstration précédente, on rappelle que Cat n'est autre
que le plan $E_{\Phi}$, avec pour $\Phi$ la famille $\mathcal{FG}_1$ des flèches
génératrices de \cats~faciles. Notons $\beta$ le plus petit cardinal régulier pour lequel les
familles $\mathcal{F}_1$ et $\mathcal{F}_2$ sont $\beta$-petites. Par le
lemme~\ref{flgenpet}, il vient que la famille $\mathcal{FG}_1$ des flèches
génératrices de \cats~faciles est encore $\beta$-petite. Alors le plan
$E_{\Phi}$ composé des plans simples $e_{\Phi,1}$ est de longueur le plus
petit cardinal régulier supérieur à $\beta$. Or par hypothèse celui-ci
est inférieur ou égal à $\alpha$. Ainsi on a que $\alpha'$, étant par
hypothèse strictement supérieur à $2^{\alpha}$, est bien strictement
supérieur à deux puissance la longueur de Cat et aussi strictement
supérieur à deux puissance le cardinal régulier pour lequel
$\mathcal{FG}_1$ est petit. Les hypothèses du lemme~\ref{petitpl} étant
vérifiées, on peut appliquer ce lemme qui nous donne le résultat
désiré.\\
CQFD.\\

Le troisième résultat que nous allons donner concerne encore la petitesse
mais cette fois il s'agit d'un résultat plus subtil. En effet on se donne une
\precat~$A$ quelconque et l'on cherche à montrer qui si l'on connaît une
sous-\precat~petite, notons-la $C$, de la catégorisée $Cat(A)$, alors il
existe une sous-\precat~petite $B$ de $A$ dont la catégorisée $Cat(B)$
contient $C$.

\begin{lem}\label{spap}
Soit $(\mathcal{C},\mathcal{F}_1,\mathcal{F}_2)$ une donnée de Segal facile
dont la catégorie sous-jacente est une \cmf~engendrée par monomorphismes.\\
Soit $\alpha'$ un cardinal régulier strictement supérieur à
$2^{\alpha}$, où l'on note $\alpha$ le plus petit cardinal régulier strictement plus
grand que celui pour lequel les
familles $\mathcal{F}_1$ et $\mathcal{F}_2$ sont petites.\\
Pour tout \precat~$t$ et pour tout $p$ sous-\precat~$\alpha'$-petite de $Cat(t)$,
il existe une sous-\precat~$\alpha'$-petite $p'$ de $t$ telle qu'on ait :\\
$p\subset Cat(p')\subset Cat(t)$.
\end{lem}
{\it Preuve :}\\
Par construction, $Cat(t)$ est la colimite séquentielle transfinie d'une
$\alpha$-séquence de \precats~partiellement marquées notée
$(t^{\beta},\mu_{\beta})_{\beta\leq\alpha}$ telle que\\
$(t^{\beta+1},\mu_{\beta+1})=e_{\Phi,1}(t^{\beta},\mu_{\beta})$, où l'on note
$\Phi$ la famille $\mathcal{FG}_1$ des flèches génératrices de
\cats~faciles. Montrons par
récurrence transfinie sur $\beta\leq\alpha$ que, pour tout $\beta\leq\alpha$, pour
toute sous-\precat~partiellement marquée $\alpha'$-petite $(p,\nu)$ de
$(t^{\beta},\mu_{\beta})$, il existe $p'\subset t$ $\alpha'$-petit tel que $(p,\nu)\subset
Cat_{|\beta}(p')\subset (t^{\beta},\mu_{\beta})$, où $Cat_{|\beta}$ est une notation pour la
construction $Cat$ stoppée à l'étape $\beta$. 
Si $\beta=0$, comme $t^0$ n'est autre que $t$, pour $p'$ il suffit
de prendre $p$.\\ 
\\
Supposons le résultat vrai au rang $\beta<\alpha$, montrons-le pour le rang
$\beta+1$. Soit donc $(p,\nu)\subset (t^{\beta+1},\mu_{\beta+1})$ $\alpha'$-petit.
On remarque que $(t^{\beta+1},\mu_{\beta+1})=e_{\Phi,1}(t^{\beta},\mu_{\beta})$. Or
par construction, $e_{\Phi,1}(t^{\beta},\mu_{\beta})$ est une somme amalgamée
multiple de $t^{\beta}$ par les
diagrammes de flèches génératrices de \cats~faciles à valeurs dans
$t^{\beta}$ non contenus dans $\mu_{\beta}$. $p$ étant une sous-\precat~$\alpha'$-petite de
$e_{\Phi,1}(t^{\beta},\mu_{\beta})$, $p$ ne peut rencontrer qu'un cardinal
$\gamma<\alpha'$ de buts de flèches
génératrices. Notons $p''$ la réunion dans $t^{\beta}$ de
$t^{\beta}\cap p$ et
des images dans $t^{\beta}$ des sources des flèches génératrices dont le but est
rencontré par $p$. Par le lemme~\ref{flgenpet}, comme les sources et buts des
flèches de $\mathcal{F}_1$ et $\mathcal{F}_2$ sont petites pour un cardinal
régulier strictement inférieur à $\alpha$ alors ceux des flèches
génératrices aussi. On obtient donc bien une sous-\precat~de $t^{\beta}$ qui est bien
$\alpha'$-petite comme réunion sur $\gamma+1<\alpha'$ de
sous-\precats~$\alpha'$-petites. Prenons pour marquage de $p''$ la restriction de
celui de $t^{\beta}$. Par la démonstration du lemme~\ref{cofcat},
$e_{\Phi,1}(p'',\mu'')$ est bien une
sous-\precat~de $e_{\Phi,1}(t^{\beta},\mu_{\beta})$. En outre, elle contient
$t^{\beta}\cap p$ ainsi que les
buts de toutes les flèches génératrices dont la source s'envoie dans $p''$,
ce qui nous redonne bien tous les buts que rencontrent $p$. Donc
$e_{\Phi,1}(p'',\mu'')$
contient bien $p$. Quant à son marquage, il s'agit de la restriction de celui
de $e_{\Phi,1}(t^{\beta},\mu_{\beta})$ à $e_{\Phi,1}(p'',\mu'')$ qui contient la restriction du marquage de
$e_{\Phi,1}(t^{\beta},\mu_{\beta})$ à $p$, lui-même plus gros que $\nu$.
On a donc montré qu'il existe $(p'',\mu'')\subset
(t^{\beta},\mu_{\beta})$ $\alpha'$-petit tel que $(p,\nu)\subset
e_{\Phi,1}(p'',\mu'')\subset
(t^{\beta+1},\mu_{\beta+1})$. En appliquant l'hypothèse de récurrence à
$(p'',\mu'')\subset
(t^{\beta},\mu_{\beta})$, il vient qu'il existe $p'\subset t$ $\alpha'$-petit tel
que $(p'',\mu'')\subset Cat_{|\beta}(p')\subset (t^{\beta},\mu_{\beta})$. En appliquant la
construction $e_{\Phi,1}$ à ce résultat et en utilisant la préservation des
cofibrations par catégorisation du lemme~\ref{cofcat}, il vient
$e_{\Phi,1}(p'',\mu'')\subset Cat_{|\beta+1}(p')\subset (t^{\beta+1},\mu_{\beta+1})$. Or on a vu plus
haut que $(p,\nu)\subset e_{\Phi,1}(p'',\mu'')$. D'où il vient que $(p,\nu)\subset
e_{\Phi,1}(p'',\mu'')\subset Cat_{|\beta+1}(p')\subset (t^{\beta+1},\mu_{\beta+1})$, ce qui montre
l'hypothèse au rang $\beta+1$.\\
\\
Supposons le résultat vrai pour tout $\gamma<\beta$, avec $\beta\leq\alpha$ ordinal
limite, et montrons le pour le rang $\beta$. Soit donc $(p,\nu)\subset
(t^{\beta},\mu_{\beta})$ $\alpha'$-petit. Comme $(t^{\beta},\mu_{\beta})$ est la
colimite séquentielle transfinie des $(t^{\gamma},\mu_{\gamma})$, on obtient
que $(p,\nu)$ est la colimite séquentielle transfinie des $(p\cap
t^{\gamma},\nu\cap\mu_{\gamma})$. En appliquant l'hypothèse de récurrence
à chaque $(p\cap
t^{\gamma},\nu\cap\mu_{\gamma})$ $\alpha'$-petits car $p$ l'est, on obtient une $\beta$-séquence de
$p_{\gamma}\subset t$ $\alpha'$-petits tels que $(p\cap
t^{\gamma},\nu\cap\mu_{\gamma})\subset Cat_{|\gamma}(p_{\gamma})\subset
(t^{\gamma},\mu_{\gamma})$. Posons $p'$ la colimite séquentielle transfinie
des $p_{\gamma}$. Comme chaque $p_{\gamma}$ est inclus dans $t$, alors leur
colimite $p'$ aussi. En outre $p'$ est $\alpha'$-petit comme colimite indexée
par $\beta\leq\alpha<\alpha'$ de \precats~$\alpha'$-petites. Comme la catégorisation
préserve les cofibrations, par le lemme~\ref{cofcat}, pour tout $\gamma<\beta$, on obtient que
$Cat_{|\gamma}(p_{\gamma})\subset Cat_{|\gamma}(p')\subset
(t^{\gamma},\mu_{\gamma})$. Donc on en déduit que pour tout $\gamma<\beta$, 
 $(p\cap t^{\gamma},\nu\cap\mu_{\gamma})\subset Cat_{|\gamma}(p')\subset
(t^{\gamma},\mu_{\gamma})$. En outre $A\rightarrow e_{\Phi,1}(A)$ est une
cofibration comme somme amalgamée de cofibrations (les flèches
génératrices) le long d'un morphisme par la démonstration du
lemme~\ref{cofcat}, donc pour tout $\gamma<\beta$, on a que $Cat_{|\gamma}(p')\subset
Cat_{|\beta}(p')$. Ainsi pour tout $\gamma<\beta$, on a que $(p\cap
t^{\gamma},\nu\cap\mu_{\gamma})\subset Cat_{|\beta}(p')\subset
(t^{\beta},\mu_{\beta})$, d'où, comme $(p,\nu)$ est la colimite des $(p\cap
t^{\gamma},\nu\cap\mu_{\gamma})$, il vient que $(p,\nu)\subset
Cat_{|\beta}(p')\subset (t^{\beta},\mu_{\beta})$. Ceci montre l'hypothèse de
récurrence au rang $\beta$.\\
\\
Par récurrence transfinie, on a donc montré que l'hypothèse de
récurrence est vraie pour tout $\beta\leq\alpha$. En particulier, elle est
vraie au rang $\alpha$, ce qui donne le résultat suivant : pour toute
sous-\precat~partiellement marquée $\alpha'$-petite $(p,\nu)$ de $Cat(t)$, il
existe $p'\subset t$ $\alpha'$-petit tel que $(p,\nu)\subset Cat(p')\subset
Cat(t)$. Pour avoir le résultat du lemme, il suffit de remarquer que toute
sous-\precat~$p$ de $Cat(t)$ est naturellement une sous-\precat~partiellement
marquée de $Cat(t)$ avec pour marquage la restriction à $p$ de celui de
$Cat(t)$.\\
CQFD.\\

Ce dernier résultat est très important car il assure que si l'on a des
informations sur la catégorisée d'une \precat~qui sont contenues dans une
sous-\precat~petite de cette catégorisée, alors nous sommes capable
d'extraire de la \precat~de départ une sous-\precat~petite dont la
catégorisée contienne elle-aussi ces informations. Ceci est en particulier
très utile si l'on veut montrer qu'un morphisme de \precats~est une
cofibration triviale, car il faut en fait montrer que le catégorisée du
morphisme est une \eq~de \cats. Or le résultat que l'on vient de montrer nous
permettra de travailler directement sur la cofibration triviale et non sur sa
catégorisée.\\

Un corollaire immédiat de ce résultat est que si deux objets d'une
\precat~sont des objets équivalents dans la catégorisée, il existe une
sous-\precat~petite de la \precat~de départ dont la catégorisée assure
l'\eq~des deux objets. Pour cela, nous allons supposer de nouveau l'existence de
l'intervalle $\bar{J}$ représentant l'\eq~d'objets d'une \cat, mais nous lui
demanderons aussi une propriété de petitesse afin de pouvoir lui appliquer
le résultat que l'on vient de montrer. 

\begin{cor}
Soit $(\mathcal{C},\mathcal{F}_1,\mathcal{F}_2)$ une donnée de Segal facile
dont la catégorie sous-jacente est une \cmf~engendrée par monomorphismes
vérifiant les propriétés suivantes :
\item 1) il existe dans $\mathcal{C-PC}$ une \cat~$\bar{J}$ munie de deux objets notés $0$ et
$1$ et ayant les propriétés suivantes :
   \item a) pour
toute \cat~facile $A$ et pour tout couple d'objets équivalents $a$ et $b$ de $A$, il
existe un morphisme de $\bar{J}$ vers $A$ envoyant $0$ sur $a$ et $1$ sur $b$,
   \item b) pour toute \cat~$A$ et pour tout morphisme de $\bar{J}$ vers $A$,
les images de $0$ et $1$ par ce morphisme sont des objets équivalents dans
$A$,
   \item c) notons $\bar{L}$, la sous-\cat~pleine de $\bar{J}$ d'objet $0$,
   alors il existe un morphisme $p:\bar{J}\rightarrow \bar{L}$ envoyant $0$ et $1$ sur $0$ qui soit
une \eq~de \cats,
\item 2) il existe dans $\mathcal{C-PC}$ une
\precat~$\bar{J}^{pre}$ $\alpha$-petite dont la catégorisation est $\bar{J}$
et vérifiant également la propriété a) de $\bar{J}$, où l'on note $\alpha$ le plus petit cardinal régulier strictement plus
grand que celui pour lequel les
familles $\mathcal{F}_1$ et $\mathcal{F}_2$ sont petites.\\
\\
Soit $\alpha'$ un cardinal régulier strictement supérieur à
$2^{\alpha}$.\\
Soit $A$ une \precat~et soient $(x,y)$ un couple d'objets équivalents de
$Cat(A)$.
Alors il existe une sous-\precat~$\alpha'$-petite $A'$ de $A$ telle que $(x,y)$
soient équivalents dans la sous-\cat~$Cat(A')$ de $Cat(A)$.
\end{cor}
{\it Preuve :}\\
Comme $x$ et $y$ sont équivalents dans $Cat(A)$, par propriété de
$\bar{J}^{pre}$, il existe un morphisme de $\bar{J}^{pre}$ dans $Cat(A)$
envoyant $0$ et $1$ sur $x$ et $y$. Comme par hypothèse, $\bar{J}^{pre}$ est
$\alpha$-petit, donc son image $J$ dans $Cat(A)$ aussi et on peut lui appliquer le
lemme précédent, ce qui nous donne qu'il existe une sous-\precat~$A'$ de $A$
qui d'une part est $\alpha'$-petite et d'autre part vérifie $J\subset
Cat(A')\subset Cat(A)$. Ainsi le morphisme de $\bar{J}^{pre}$ dans $Cat(A)$ se
factorise par $Cat(A')$. Par propriété de Cat, on obtient donc un morphisme
de $\bar{J}$, catégorisation de $\bar{J}^{pre}$, dans $Cat(A')$ envoyant $0$
et $1$ sur $x$ et $y$, ce qui signifie que $x$ et $y$ sont équivalents dans
$Cat(A')$.\\
CQFD.\\

Après tous ces lemmes techniques sur la catégorisation Cat, nous pouvons
maintenant nous attaquer à la vérification des hypothèses du lemme
technique~\ref{ast2} concernant les engendrements.

\subsection{Sous-cofibrations triviales à sources et buts $\alpha'$-petits}

 Le principal problème
qui survient lorsque l'on veut montrer les hypothèses du lemme~\ref{ast2}
pour $\mathcal{J}$ est justement de pouvoir extraire d'une cofibration triviale
quelconque une sous-cofibration triviale qui est dans $\mathcal{J}$,
c'est-à-dire une sous-cofibration triviale petite. L'idée pour faire cela
est de prendre une sous-\precat~petite du but et de montrer que la restriction
de la cofibration triviale à cette sous-\precat~petite est encore une
cofibration triviale. Malheureusement ce n'est pas aussi simple car la
restriction d'une cofibration triviale n'est en général pas une cofibration
triviale. Il va donc falloir forcer cette propriété, ce que nous allons
faire en deux temps. Tout d'abord, nous allons forcer l'essentielle
surjectivité par l'intermédiaire d'une opération $a(.)$ qui conserve la
petitesse et qui est essentiellement
basée sur l'application du corollaire que l'on vient de montrer au sujet
d'objets équivalents dans une \precat.

\begin{lem}
Soit $(\mathcal{C},\mathcal{F}_1,\mathcal{F}_2)$ une donnée de Segal facile
dont la catégorie sous-jacente est une \cmf~engendrée par monomorphismes
vérifiant les propriétés suivantes :
\item 1) il existe dans $\mathcal{C-PC}$ une \cat~$\bar{J}$ munie de deux objets notés $0$ et
$1$ et ayant les propriétés suivantes :
   \item a) pour
toute \cat~facile $A$ et pour tout couple d'objets équivalents $a$ et $b$ de $A$, il
existe un morphisme de $\bar{J}$ vers $A$ envoyant $0$ sur $a$ et $1$ sur $b$,
   \item b) pour toute \cat~$A$ et pour tout morphisme de $\bar{J}$ vers $A$,
les images de $0$ et $1$ par ce morphisme sont des objets équivalents dans
$A$,
   \item c) notons $\bar{L}$, la sous-\cat~pleine de $\bar{J}$ d'objet $0$,
   alors il existe un morphisme $p:\bar{J}\rightarrow \bar{L}$ envoyant $0$ et $1$ sur $0$ qui soit
une \eq~de \cats,
\item 2) il existe dans $\mathcal{C-PC}$ une
\precat~$\bar{J}^{pre}$ $\alpha$-petite dont la catégorisation est $\bar{J}$
et vérifiant également la propriété a) de $\bar{J}$, où l'on note $\alpha$ le plus petit cardinal régulier strictement plus
grand que celui pour lequel les
familles $\mathcal{F}_1$ et $\mathcal{F}_2$ sont petites.\\
\\
Soit $\alpha'$ un cardinal régulier strictement supérieur à
$2^{\alpha}$.\\
Pour toute cofibration triviale $f:s\rightarrow t$, il existe une construction
$a(.)$ qui à toute
sous-\precat~$\alpha'$-petite $p$ de $t$ associe une
sous-\precat~$\alpha'$-petite $a(p)$ comprise entre $p$ et $t$ telle que le
morphisme $f_{|a(p)}:a(p)\cap s\rightarrow a(p)$ vérifie que pour tout objet $y$ de $p$,
il existe un objet $x$ de $a(p)\cap s$ tel que $y$ soit équivalent à $x$
dans $Cat(a(p))$.
\end{lem}
{\it Preuve :}\\
Soit $p$ une sous-\precat~$\alpha'$-petite de $t$ et soit $y$ un objet de $p$.
Comme $f$ est une \eq~faible, $Cat(f)$ est essentiellement surjective. Donc il
existe un objet $x$ de $s$ tel que $y$ et $x$ soient équivalents dans
$Cat(t)$. Par le corollaire précédent, il vient qu'il existe une
sous-\precat~$\alpha'$-petite $d(y)$ de $t$ telle que $x$ et $y$ soient
équivalents dans $Cat(d(y))\subset Cat(t)$.\\ 
\\
Posons maintenant $a_y(p)=p\cup d(y)$. Comme réunion de deux sous-\precats~de
$t$ $\alpha'$-petites, $a_y(p)$ est bien une
sous-\precat~de $t$ $\alpha'$-petite et qui par construction contient $p$. En
outre comme par le lemme~\ref{cofcat}, les cofibrations sont stables par
catégorisation, $Cat(a_y(p))$ est comprise entre $Cat(d(y))$ et $Cat(t)$, donc
$x$ et $y$ sont équivalents dans $Cat(a_y(p))\subset Cat(t)$. Enfin on
remarque que $x$ est un objet à la fois de $s$ et de $a_y(p)$, car $x$
appartient à $d(y)$. Donc le morphisme $a_y(p)\cap s\rightarrow a_y(p)$
vérifie la propriété du lemme pour $y$.\\
\\
Définissons enfin $a(p)$ comme la réunion des $a_y(p)$ pour $y$ décrivant
l'ensemble des objets de $p$. Comme $p$ est $\alpha'$-petit, cet ensemble est de
cardinal strictement inférieur à $\alpha'$. En outre comme chaque $a_y(p)$ est une
sous-\precat~$\alpha'$-petite de $t$ contenant $p$, par hypothèse sur
$\alpha'$, il vient que $a(p)$ est une sous-\precat~$\alpha'$-petite de $t$
contenant $p$. Montrons que $a(p)\cap s\rightarrow a(p)$ vérifie la
propriété du lemme. Soit $y$ un objet de $p$. Par construction de $a_y(p)$,
il existe $x$ objet de $a_y(p)\cap s$ donc de $a(p)\cap s$ tel que $x$ soit
équivalent à $y$ dans $Cat(a_y(p))\subset Cat(a(p))\subset Cat(t)$.\\
CQFD.\\

Après avoir traité de l'essentielle surjectivité, donnons une construction
$b(.)$ qui va forcer la pleine fidélité des restrictions de cofibrations
triviales tout en gardant stable la petitesse.

\begin{lem}
Soit $(\mathcal{C},\mathcal{F}_1,\mathcal{F}_2)$ une donnée de Segal facile
dont la catégorie sous-jacente est une \cmf~engendrée par monomorphismes telle que les réunions de deux
sous-objets sont exactement les sommes amalgamées de ces sous-objets au-dessus
de leur intersection.\\
Soit $\alpha'$ un cardinal régulier strictement supérieur à
$2^{\alpha}$, où l'on note $\alpha$ le plus petit cardinal régulier strictement plus
grand que celui pour lequel les
familles $\mathcal{F}_1$ et $\mathcal{F}_2$ sont petites.\\
Alors il existe $\bar{\alpha'}$ un cardinal régulier strictement supérieur à
$\alpha'$ tel que 
pour toute cofibration triviale $f:s\rightarrow t$, il existe une construction
$b(.)$ qui à toute
sous-\precat~$\alpha'$-petite $p$ de $t$ associe une
sous-\precat~$\bar{\alpha'}$-petite $b(p)$ comprise entre $p$ et $t$ telle que le
morphisme $f_{|b(p)}:b(p)\cap s\rightarrow b(p)$ vérifie la propriété
suivante : pour tout couple $(x,y)$ d'objets de $s\cap p$, il existe un objet
$\bar{\alpha'}$-petit $D$ compris entre $Cat(p)_1(x,y)$ et $Cat(t)_1(x,y)$ tel que
$D\cap Cat(s)_1(x,y)\rightarrow D$ est une sous-cofibration triviale à source
et but $\bar{\alpha'}$-petits de la cofibration $Cat(f_{|b(p)})_1(x,y)$.
\end{lem}
{\it Preuve :}\\
Soit $p$ une sous-\precat~$\alpha'$-petite de $t$. Pour tout couple de points
$(x,y)$ de $s\cap p$, comme $f$ est une cofibration triviale, $Cat(f)_1(x,y)$ est une
cofibration triviale dans $\mathcal{C}$. Comme $p$ est
$\alpha'$-petit, alors par le lemme~\ref{petitcat}, $Cat(p)$ est
$\alpha''$-petite, avec $\alpha''>\alpha'$ cardinal régulier indépendant de
$p$, et donc $Cat(p)_1(x,y)$ aussi. En outre comme les
cofibrations sont stables par catégorisation par le lemme~\ref{cofcat},
$Cat(p)_1(x,y)$ est un sous-\obc~de $Cat(t)_1(x,y)$. Comme $\mathcal{C}$ est une
\cmfcof, elle vérifie la proposition suivante. On applique alors la
proposition suivante à la cofibration triviale de $\mathcal{C}$ $Cat(f)_1(x,y):
Cat(s)_1(x,y)\rightarrow Cat(t)_1(x,y)$ et au sous-objet $\alpha''$-petit
$Cat(p)_1(x,y)$ de $Cat(t)_1(x,y)$. Il existe donc un objet
$\widetilde{\alpha''}$-petit $D_{x,y}$, avec $\widetilde{\alpha''}>\alpha''$ un cardinal
régulier indépendant de $f$ et $p$, tel que $D_{x,y}$ est compris entre $Cat(p)_1(x,y)$ et
$Cat(t)_1(x,y)$, tel que $D_{x,y}\cap Cat(s)_1(x,y)\rightarrow D_{x,y}$
soit une sous-cofibration triviale de $Cat(f)_1(x,y)$. Consi\-dérons $D'_{x,y}$
la somme amalgamée de $Cat(p)$ par $\Delta[1]\Theta D_{x,y}$ au-dessus de
$\Delta[1]\Theta Cat(p)_1(x,y)$. Comme cette somme amalgamée n'est autre
que la sous-\precat~de $Cat(t)$ engendrée par les $\widetilde{\alpha''}$-petits $Cat(p)$ et
$D_{x,y}$, on obtient que $D'_{x,y}$ est une sous-\precat~$\widetilde{\alpha''}$-petite de
$Cat(t)$. En appliquant le lemme~\ref{spap} à la
sous-\precat~$\widetilde{\alpha''}$-petite $D'_{x,y}$ de $Cat(t)$, il vient qu'il existe une
sous-\precat~$\widetilde{\alpha''}$-petite $b'_{x,y}(p)$ de $t$ telle que $D'_{x,y}\subset
Cat(b'_{x,y}(p)\subset Cat(t))$.\\ 
\\
On pose maintenant $b_{x,y}(p)=p\cup
b'_{x,y}(p)$ qui est bien une sous-\precat~$\widetilde{\alpha''}$-petite de $t$ comme
réunion de telles \precats. En outre par
construction, $b_{x,y}(p)$ contient $p$. Par stabilité des cofibrations par
catégorisation du lemme~\ref{cofcat}, il vient que $Cat(b_{x,y}(p))$ contient $Cat(b'_{x,y}(p))$ qui
elle-même contenait $D'_{x,y}$. En outre $(D'_{x,y})_1(x,y)$ contient
$D_{x,y}$. Donc au final, on a que $D_{x,y}\subset
Cat(b_{x,y}(p))_1(x,y)\subset Cat(t)_1(x,y)$.\\
\\
On définit maintenant $b(p)$ comme la réunion des $b_{x,y}(p)$ pour $x$ et
$y$ décrivant l'ensemble des objets de $p$ dont le cardinal est strictement
inférieur à $\alpha'$ car $p$ est $\alpha'$-petit. Donc par hypothèse sur
$\widetilde{\alpha''}$, on obtient que $b(p)$ est une
sous-\precat~$\widetilde{\alpha''}$-petite de $t$ contenant $p$. Montrons que $f_{|b(p)}:b(p)\cap
s\rightarrow b(p)$ a la propriété du lemme. Soit $(x,y)$ un couple d'objets
de $s\cap p$, alors comme on l'a vu ci-dessus, il existe un objet $D_{x,y}$
$\widetilde{\alpha''}$-petit compris entre $Cat(p)_1(x,y)$ et $Cat(t)_1(x,y)$ et
tel que $D_{x,y}\cap Cat(s)_1(x,y)\rightarrow D_{x,y}$ est une cofibration
triviale. On remarque que la source et le but de cette dernière sont bien
$\widetilde{\alpha''}$-petits car
$D_{x,y}$ et donc $D_{x,y}\cap Cat(s)_1(x,y)$ sont $\widetilde{\alpha''}$-petits. Or $D_{x,y}$ est inclus dans $Cat(b_{x,y}(p))_1(x,y)$ donc dans
$Cat(b(p))_1(x,y)$. Ainsi $D_{x,y}\cap Cat(s)_1(x,y)\rightarrow D_{x,y}$
est bien une sous-cofibration triviale à source et but $\widetilde{\alpha''}$-petits de la
cofibration $Cat(b(p))_1(x,y)\cap
Cat(s)_1(x,y)\rightarrow Cat(b(p))_1(x,y)$ qui n'est autre que
$Cat(f_{|b(p)})_1(x,y)$ car $Cat$ commute aux intersections par le
lemme~\ref{intercat}. Comme $\widetilde{\alpha''}$ est indépendant de $p$ et de $f$
mais seulement de $\alpha''$ qui lui même ne dépend que de $\alpha'$, en
posant $\bar{\alpha'}=\widetilde{\alpha''}$, on a le résultat attendu.\\
CQFD.\\

En appliquant un nombre transfini de fois les opérations $a(.)$ et $b(.)$,
nous allons transformer toute restriction de cofibration triviale en cofibration
triviale et ce tout en préservant la petitesse. Si donc pour une cofibration
triviale quelconque, on extrait de son but une sous-\precat~petite, en appliquant un
tel procédé à la restriction de la cofibration triviale à cette
sous-\precat~petite, on obtiendra une sous-cofibration triviale petite de notre
cofibration triviale de départ, résultat auquel nous voulions arriver. 

\begin{prop}\label{alphatilde}
Soit $(\mathcal{C},\mathcal{F}_1,\mathcal{F}_2)$ une donnée de Segal facile
vérifiant les propriétés suivantes :
\item 1) la \cmf~sous-jacente $\mathcal{C}$ est une \cmf~simpliciale
engendrée par monomorphismes,
\item 2) les cofibrations génératrices de la \cmfcof~$\mathcal{C}$ ont leurs
sources et buts $\alpha$-petits, où l'on note $\alpha$ le plus petit cardinal régulier strictement plus
grand que celui pour lequel les
familles $\mathcal{F}_1$ et $\mathcal{F}_2$ sont petites,
\item 3) les cofibrations génératrices de $\mathcal{C}$ ont leurs buts connexes,
\item 4) les réunions de deux
sous-objets sont exactement les sommes amalgamées de ces sous-objets au-dessus
de leur intersection,
\item 5) pour tout monomorphisme $j:A\rightarrow B$ de $\mathcal{C}$, la classe contenant un représentant de chaque classe d'isomorphismes d'objets $C$ par lequel $j$ se
factorise dans la sous-catégorie des monomorphismes est un ensemble,
\item 6) les colimites séquentielles transfinies de la sous-catégorie des
monomorphismes de $\mathcal{C}$ existent et sont les colimites séquentielles
transfinies de $\mathcal{C}$,
\item 7) il n'existe aucun morphisme de but l'ensemble vide qui ne soit pas un
isomorphisme,
\item 8) les colimites séquentielles transfinies
d'\obcs~le long de cofibrations sont des \obcs,
\item 9) il existe une \cat~$\bar{J}$ munie de deux objets notés $0$ et
$1$ et ayant les propriétés suivantes :
\item a) pour
toute \cat~facile $A$ et pour tout couple d'objets équivalents $a$ et $b$ de $A$, il
existe un morphisme de $\bar{J}$ vers $A$ envoyant $0$ sur $a$ et $1$ sur $b$,
\item b) pour toute \cat~$A$ et pour tout morphisme de $\bar{J}$ vers $A$,
les images de $0$ et $1$ par ce morphisme sont des objets équivalents dans
$A$,
\item c) notons $\bar{L}$, la sous-\cat~pleine de $\bar{J}$ d'objet $0$, alors
il existe un morphisme $p:\bar{J}\rightarrow \bar{L}$ envoyant $0$ et $1$ sur $0$ qui soit
une \eq~de \cats,
\item 10) il existe dans $\mathcal{C-PC}$ une
\precat~$\bar{J}^{pre}$ $\alpha$-petite dont la catégorisation est $\bar{J}$
et vérifiant également la propriété a) de $\bar{J}$.\\
\\
Soit $\alpha'$ un cardinal régulier strictement supérieur à
$2^{\alpha}$.\\
Alors il existe $\widetilde{\alpha'}$ un cardinal régulier strictement supérieur à
$\alpha'$ tel que
pour toute cofibration triviale $f:s\rightarrow t$ et pour toute
sous-\precat~$\alpha'$-petite $p$ de $t$, il existe une
sous-\precat~$\widetilde{\alpha'}$-petite $p'$ comprise entre $p$ et $t$ telle qu'on ait le
diagramme suivant :
\begin{diagram}
s\cap p' & \rTo^{\subset} & p' \\
\dTo^{\subset} & & \dTo_{\subset} \\
s & \rTo^{\subset} & t \\
\end{diagram}
où les flèches horizontales sont des \eqs~faibles.
Dans ce cas, l'inclusion de $s\cap p'$ dans $p'$ est dans $\mathcal{J}$, si
$\alpha'$ vaut le cardinal de la définition de $\mathcal{J}$, et la
cofibration triviale $s\rightarrow t$ se factorise à travers $s\cup p'$ en
deux cofibrations triviales.
\end{prop}
{\it Preuve :}\\
Supposons que le cardinal $\widetilde{\alpha'}$ existe et qu'il existe $p'$ $\widetilde{\alpha'}$-petit compris entre $p$ et $t$ vérifiant
le diagramme ci-dessus. Montrons d'abord que $s\cap p'\rightarrow p'$ est dans
$\mathcal{J}$, si $\alpha'$ vaut le cardinal de la définition de $\mathcal{J}$.
On a déjà que $s\cap p'\rightarrow p'$ est une cofibration
triviale. Or par hypothèse, $p'$ et donc $s\cap p'$ sont $\widetilde{\alpha'}$-petits donc
$s\cap p'\rightarrow p'$ est une cofibration triviale à source et but
$\widetilde{\alpha'}$-petits, i.e. un élément de $\mathcal{J}$.\\
\\
Toujours avec la supposition ci-dessus, montrons maintenant que la cofibration triviale $s\rightarrow t$ se factorise à travers $s\cup p'$ en
deux cofibrations triviales. Considérons le diagramme cocartésien suivant :
\begin{diagram}
s\cap p' & \rTo^{\subset} & p' \\
\dTo^{\subset} & & \dTo_{\subset} \\
s & \rTo^{\subset} & s\cup p' \\
\end{diagram}
Comme les cofibrations triviales sont stables par somme amalgamée le long d'un
morphisme par proposition~\ref{ctpush} et que $s\cap p'\rightarrow p'$ est une
cofibration triviale, alors $s\rightarrow s\cup p'$ est une cofibration
triviale. En outre, $s$ et $p'$ sont des sous-\precats~de $t$, donc en utilisant
le fait que $s\cup p'$ soit la somme amalgamée de $s$ et $p'$ au-dessus de
leur intersection, par hypothèse sur $\mathcal{C}$, il vient que le morphisme
universel $s\cup p'\rightarrow t$ est une inclusion donc une cofibration. Comme
enfin la précomposition de ce morphisme par la cofibration triviale $s\rightarrow s\cup p'$ n'est
autre que la cofibration triviale $s\rightarrow t$, par l'axiome "trois pour
deux" dans $\mathcal{C}$, il vient que $s\cup p'\rightarrow t$ est une \eq~faible. Donc $s\cup p'\rightarrow t$ est bien une cofibration
triviale et $s\rightarrow t$ se factorise donc bien par $s\cup p'$ en deux
cofibrations triviales.\\
\\
Définissons $p'$ et $\widetilde{\alpha'}$ par récurrence. Soit $p'_0=p$ et
$\alpha'_0=\alpha'$. Supposons
$p'_n$ et $\alpha'_n$ construits pour $n\geq 0$, alors on pose
$p'_{n+1}=b(a(p'_n))$ et $\alpha'_{n+1}=\bar{\alpha'_n}$, au sens du lemme
précédent.  Enfin on définit $p'$ comme
la colimite des $p'_n$ pour $n\in\mathbb{N}$ et $\widetilde{\alpha'}$ comme le plus
petit cardinal régulier majorant l'ensemble $\{\alpha'_n,n\in\mathbb{N}\}$, ce qui existe par les propriétés des cardinaux. Il est facile de voir qu'à
chaque étape le cardinal défini dépend uniquement de $\alpha'$. On remarque, comme
$\widetilde{\alpha'}>\aleph_0$ est régulier et est supérieur strictement à
chaque $\alpha'_n$, que $p'$ est
$\widetilde{\alpha'}$-petit comme réunion dénombrable de \precats~$\alpha'_n$-petites.  Comme
$p'$ est une réunion croissante de \precats~comprises entre $p$ et $t$ alors
$p'$ aussi. Il ne reste donc qu'à montrer la propriété de $s\cap
p'\rightarrow p'$. Comme c'est clairement une inclusion, montrons que c'est une
\eq~faible, i.e. montrons que $Cat(s\cap p')\rightarrow Cat(p')$ est une \eq~de
\cats.\\
\\
Soit $y$ un objet de $p'$. Comme $p'$ est une colimite, il existe un
entier $n$ tel que $y$ soit un objet de $p'_n$. Appliquons la
construction $a(.)$ à $p'_n$, on obtient alors l'existence d'un objet $x$
de $a(p'_n)\cap s$ tel que $x$ soit équivalent à $y$ dans
$Cat(a(p'_n))$. Comme $b(a(p'_n))$ qui n'est autre que
$p'_{n+1}$ contient $a(p'_n)$, que $p'$ contient $p'_{n+1}$ et que
les cofibrations sont stables par catégorisation par le lemme~\ref{cofcat}, il vient que $x$ appartient
à $p'\cap s$ et que $x$ est équivalent à $y$ dans $Cat(p')$, ce qui
montre l'essentielle surjectivité de $Cat(s\cap p')\rightarrow Cat(p')$.\\
\\
Soit $(x,y)$ un couple d'objets de $p'\cap s$. Il existe donc un entier $n$ tel
que ce couple appartienne à $p'_n$ donc à $a(p'_n)$. On
applique la construction $b(.)$ à tous les $a(p'_m)$ pour $m\geq
n$ et, en utilisant la commutation de Cat avec l'intersection du
lemme~\ref{intercat}, l'on obtient le diagramme suivant : 
\begin{diagram}
Cat(p'_m\cap s)_1(x,y) & \rTo^{\subset}  & Cat(a(p'_m)\cap s)_1(x,y) &
\rTo^{\subset}  & D^{m}_{x,y}\cap
Cat(s)_1(x,y)  & \rTo^{\subset} & Cat(p'_{m+1}\cap s)_1(x,y) \\ 
\dTo^{\subset} & & \dTo^{\subset} & & \dTo^{\widetilde{\subset}}  & & \dTo^{\subset} \\
Cat(p'_m)_1(x,y) & \rTo^{\subset} & Cat(a(p'_m))_1(x,y) & \rTo^{\subset}  & D^{m}_{x,y} & \rTo^{\subset} &
Cat(p'_{m+1})_1(x,y) \\
\end{diagram}
Les flèches $Cat(p'_m\cap
s)_1(x,y)\rightarrow Cat(p'_m)_1(x,y)$ et $D^m_{x,y}\cap
Cat(s)_1(x,y)\rightarrow D^m_{x,y}$ alternant, on obtient que la colimite pour
$m\geq n$ des $Cat(p'_m\cap
s)_1(x,y)\rightarrow Cat(p'_m)_1(x,y)$, qui n'est autre que $Cat(p'\cap
s)_1(x,y)\rightarrow Cat(p')_1(x,y)$, est égale à la colimite des $D^m_{x,y}\cap
Cat(s)_1(x,y)\rightarrow D^m_{x,y}$. En utilisant la
proposition~\ref{eqcolim} dans $\mathcal{C}$, on obtient que la colimite des
cofibrations triviales $D^m_{x,y}\cap
Cat(s)_1(x,y)\rightarrow D^m_{x,y}$ est une \eq~faible. Or comme cette
dernière est égale à $Cat(p'\cap
s)_1(x,y)\rightarrow Cat(p')_1(x,y)$, ceci montre que $Cat(p'\cap
s)_1(x,y)\rightarrow Cat(p')_1(x,y)$ est une \eq~faible et donc que $Cat(s\cap
p')\rightarrow Cat(p')$ est pleinement fidèle.\\
CQFD.\\

Nous avons donc démontré le résultat qui permet d'extraire de toute
cofibration triviale une cofibration triviale petite, c'est-à-dire un
élément de $\mathcal{J}$. Nous avons donc maintenant en notre possession
l'outil majeur permettant de vérifier les hypothèses du lemme~\ref{ast2} qui
montrera l'engendrement des cofibrations triviales par $\mathcal{J}$.

\subsection{Engendrement des cofibrations triviales par $\mathcal{J}$}

Afin de démontrer la dernière hypothèse du lemme de
reconnaissance~\ref{reco} de la structure de \cmfcof, nous devons vérifier les
hypothèses du lemme~\ref{ast2} qui permettra de montrer que $\mathcal{J}$
engendre les cofibrations triviales. Or nous venons de montrer l'argument de
base qui permet de vérifier les hypothèses du lemme~\ref{ast2}. Aussi
va-t-on en récolter les hypothèses, qui serviront pour toute cette partie.

\begin{hyp}\label{Jcoftriv}
Soit $(\mathcal{C},\mathcal{F}_1,\mathcal{F}_2)$ une donnée de Segal facile
vérifiant les propriétés suivantes :
\item 1) la \cmf~sous-jacente $\mathcal{C}$ est une \cmf~simpliciale
engendrée par monomorphismes,
\item 2) les cofibrations génératrices de la \cmfcof~$\mathcal{C}$ ont leurs
sources et buts $\alpha$-petits, où l'on note $\alpha$ le plus petit cardinal régulier strictement plus
grand que celui pour lequel les
familles $\mathcal{F}_1$ et $\mathcal{F}_2$ sont petites,
\item 3) les cofibrations génératrices de $\mathcal{C}$ ont leurs buts connexes,
\item 4) les réunions de deux
sous-objets sont exactement les sommes amalgamées de ces sous-objets au-dessus
de leur intersection,
\item 5) pour tout monomorphisme $j:A\rightarrow B$ de $\mathcal{C}$, la classe contenant un représentant de chaque classe d'isomorphismes d'objets $C$ par lequel $j$ se
factorise dans la sous-catégorie des monomorphismes est un ensemble,
\item 6) les colimites séquentielles transfinies de la sous-catégorie des
monomorphismes de $\mathcal{C}$ existent et sont les colimites séquentielles
transfinies de $\mathcal{C}$,
\item 7) il n'existe aucun morphisme de but l'ensemble vide qui ne soit pas un
isomorphisme,
\item 8) les colimites séquentielles transfinies
d'\obcs~le long de cofibrations sont des \obcs,
\item 9) il existe une \cat~$\bar{J}$ munie de deux objets notés $0$ et
$1$ et ayant les propriétés suivantes :
\item a) pour
toute \cat~facile $A$ et pour tout couple d'objets équivalents $a$ et $b$ de $A$, il
existe un morphisme de $\bar{J}$ vers $A$ envoyant $0$ sur $a$ et $1$ sur $b$,
\item b) pour toute \cat~$A$ et pour tout morphisme de $\bar{J}$ vers $A$,
les images de $0$ et $1$ par ce morphisme sont des objets équivalents dans
$A$,
\item c) notons $\bar{L}$, la sous-\cat~pleine de $\bar{J}$ d'objet $0$, alors
il existe un morphisme $p:\bar{J}\rightarrow \bar{L}$ envoyant $0$ et $1$ sur $0$ qui soit
une \eq~de \cats,
\item 10) il existe dans $\mathcal{C-PC}$ une
\precat~$\bar{J}^{pre}$ $\alpha$-petite dont la catégorisation est $\bar{J}$
et vérifiant également la propriété a) de $\bar{J}$.
\end{hyp}

Sous ces hypothèses, nous savons donc extraire de toute cofibration triviale
un élément de $\mathcal{J}$, ce qui entraîne que les flèches
$\mathcal{J}$-injectives se relèvent aussi par rapport aux cofibrations
triviales quelconques.

\begin{lem}
Supposons vraies les hypothèses~\ref{Jcoftriv}.\\
Les flèches $\mathcal{J}$-injectives sont exactement les flèches ayant la
\prd~par rapport aux cofibrations triviales.
\end{lem}
{\it Preuve :}\\
Tout d'abord comme $\mathcal{J}$ est un sous-ensemble de la classe des
cofibrations triviales, les flèches ayant la \prd~par rapport aux cofibrations
triviales ont en particulier cette propriété vis-à-vis des éléments de
$\mathcal{J}$ donc sont $\mathcal{J}$-injectives. Il reste donc à montrer
l'autre sens.\\
Considérons maintenant le diagramme commutatif suivant :
\begin{diagram}
A & \rTo^a & X \\
\dTo^{j} & & \dTo_p \\
B & \rTo_b & Y \\
\end{diagram}
avec $j$ une cofibration triviale et $p$ une flèche $\mathcal{J}$-injective.
Si $j$ est un isomorphisme, le diagramme admet un relèvement naturel. Sinon il
existe une sous-\precat~$D$ $\alpha'$-petite de $B$ non incluse dans $j(A)$ mais
d'intersection non vide avec $j(A)$, où $\alpha'$ est le cardinal de la
définition de $\mathcal{J}$. Par la proposition
précédente, il existe une sous-\precat~$\widetilde{\alpha'}$-petite $D'$ comprise entre
$D$ et $B$ tel que $A\cap D'\rightarrow D'$ soit une sous-cofibration triviale
de $j$ à source et but
$\widetilde{\alpha'}$-petits. En particulier $A\cap D'\rightarrow D'$ est un
élément de $\mathcal{J}$. On a donc le diagramme suivant :
\begin{diagram}
A\cap D' & \rTo^{\subset} & A & \rTo^a & X \\
\dTo^{\in \mathcal{J}} & & \dTo^{j} & & \dTo_p \\
D' & \rTo_{\subset} & B & \rTo_b & Y \\
\end{diagram}
Comme $p$ est $\mathcal{J}$-injective, le diagramme extérieur admet un
relèvement $\theta$.\\ 
Considérons maintenant la réunion $A\cup D'$. Par hypothèse sur $D$,
cette réunion n'est pas isomorphe à $A$. D'après la proposition
précédente, $j$ se factorise à travers $A\cup D'$ en deux cofibrations
triviales dont la première est l'inclusion canonique $i_A:A\rightarrow A\cup
D'$ et la seconde $j':A\cup D'\rightarrow B$. Avec ce qui précède on
obtient le diagramme commutatif suivant :
\begin{diagram}
A\cap D' & \rTo^{\subset} & A &  & \\
\dTo^{\in \mathcal{J}} & & \dTo_{i_A} & \rdTo(2,4)^{a} &  \\
D' & \rTo_{i_D'} & A\cup D' &  &  \\
 & \rdTo(4,2)_{\theta} & & \rdTo~{\exists ! \phi} & \\
 & & & & X \\
\end{diagram}
Comme la partie extérieure du grand diagramme commute par propriété de
$\theta$ et que $A\cup D'$ est une somme amalgamée, il existe un unique
morphisme $\phi$ de $A\cup D'$ vers $X$ tel que $\phi\circ i_A=a$ et $\phi\circ
i_D=\theta$. De là il vient les égalités suivantes :
$$p\circ\phi \circ i_A=p\circ a=b\circ j=b\circ j'\circ i_A$$
$$p\circ\phi \circ i_D'=p\circ\theta=b\circ j'\circ i_D' $$
Comme $p\circ\phi$ et $b\circ j'$ ont même valeur sur les injections
canoniques, par propriété universelle de la somme amalgamée $A\cup D'$, on
obtient que $p\circ\phi=b\circ j'$. Ainsi $A\cup D'$ réalise un relèvement
partiel du diagramme $(j,a,b,p)$ avec $A\cup D'$ non isomorphe à $A$. En appliquant le lemme~\ref{ast2}, il vient que
$(j,a,b,p)$ possède un relèvement global. Ceci montre que les morphismes
$\mathcal{J}$-injectifs ont la \prd~par rapport aux cofibrations triviales.\\
CQFD.\\

Comme l'on vient de montrer que les morphismes $\mathcal{J}$-injectifs sont
exactement ceux qui se relèvent par rapport aux cofibrations triviales, on
obtient donc que les cofibrations triviales sont des $\mathcal{J}$-cofibrations,
c'est-à-dire des morphismes ayant la \prg~par rapport aux morphismes
$\mathcal{J}$-injectifs. Comme avec l'hypothèse 3) du lemme de
reconnaissance~\ref{reco}, on a montré que les $\mathcal{J}$-cofibrations sont
des cofibrations triviales, on a donc bien que $\mathcal{J}$ engendre les
cofibrations triviales.

\begin{cor}\label{p5}
Supposons vraies les hypothèses~\ref{Jcoftriv}.\\
$\mathcal{J}$ engendre les cofibrations triviales.
\end{cor}
{\it Preuve :}\\
Le lemme précédent montre que les cofibrations triviales ont la \prg~par
rapport aux flèches $\mathcal{J}$-injectives. Donc par définition de
$\mathcal{J}$-cofibrations comme flèches ayant la \prg~par rapport aux
flèches $\mathcal{J}$-injectives, les cofibrations triviales sont des
$\mathcal{J}$-cofibra\-tions. En outre le corollaire~\ref{p3} montre que les
$\mathcal{J}$-cofibrations sont des cofibrations triviales. Donc les
cofibrations triviales sont exactement les $\mathcal{J}$-cofibrations, ce qui
montre que $\mathcal{J}$ engendre les cofibrations triviales.\\
CQFD.\\

Avec ce résultat d'engendrement des cofibrations triviales par $\mathcal{J}$,
qui mérite donc bien son nom d'ensemble des cofibrations triviales
génératrices, nous venons de vérifier la dernière hypothèse du
lemme~\ref{reco} de reconnaissance de la structure de \cmfcof. Nous pouvons donc
maintenant donner le théorème central qui établit que la catégorie
$\mathcal{C-PC}$ des \precats~avec les monomorphismes et les \eqs~faibles de
\precats~forme bien une \cmf.

\newpage

\section{Théorème central}

Comme la vérification des six hypothèses du lemme~\ref{reco} de
reconnaissance de la structure de \cmfcof~a donné lieu à un bon nombre
d'hypothèses, nous allons dans un premier temps les rassembler dans une ébauche de
théorème sur la \cmfcof~des \precats, puis les retravailler afin d'obtenir une meilleure formulation de
ce dernier.

\subsection{Ebauche du théorème central}

En regroupant tous les résultats précédents, on obtient l'ébauche du
théo\-rème central suivante :

\begin{theo}[-central (ébauche)]
Soit $(\mathcal{C},\mathcal{F}_1,\mathcal{F}_2)$ une donnée de Segal proto-facile
vérifiant les propriétés suivantes :
\item A) la catégorie sous-jacente $\mathcal{C}$ est une \cmfcof~simpliciale dont tous les
objets sont cofibrants,
\item B) il n'existe aucun morphisme de but l'ensemble vide qui ne soit pas un
isomorphisme,
\item C) les réunions de deux
sous-objets sont exactement les sommes amalgamées de ces sous-objets au-dessus
de leur intersection,
\item D) les objets discrets sont stables par limites dans
$\mathcal{C}$,
\item E) les colimites séquentielles transfinies
d'\obcs~le long de cofibrations sont des \obcs,
\item F) le foncteur $\tau_0$ est tel que, pour tout \obc~$C$ de
$\mathcal{C}$, $\tau_0(C)$ est un quotient de l'ensemble des morphismes dans
$\mathcal{C}$ de l'objet final vers $C$,
\item G) il existe un \obc~contractile $K$ muni de deux
morphismes de l'objet final vers $K$ notés 0 et 1, et ayant la propriété
suivante : pour tout couple $(f,g)$ de morphismes de l'objet final vers un \obc~$C$ tel que leurs images par $\tau_0$ soient égales, il existe un
morphisme de $K$ vers $C$ envoyant 0 sur $f$ et 1 sur $g$,
\item H) les $\mathcal{F}_1\cup\mathcal{F}_2\cup\{\emptyset\rightarrow *\}$-cofibrations sont des monomorphismes,
\item I) les sources et buts des flèches de la famille $\mathcal{F}_1$, ainsi
que les buts de la famille $\mathcal{F}_2$, sont
connexes et non isomorphes à l'ensemble vide,
\item J) il existe un cardinal transfini régulier strictement supérieur à
$\aleph_0$, que l'on notera $\alpha$, pour lequel tout morphisme de
$\mathcal{F}_1$ et $\mathcal{F}_2$ a sa source et son but $\alpha$-petits au
sens~\ref{alphapetit},
\item K) la famille $\mathcal{F}_1$ est incluse dans la classe des cofibrations
triviales de la \cmf~$\mathcal{C}$ et la famille $\mathcal{F}_2$ dans celle des
cofibrations,
\item L) les \eqcs~d'\obcs~sont exactement les \eqs~faibles de la
\cmf~$\mathcal{C}$ entre \obcs,
\item M) le produit fibré d'une équivalence faible entre
\obcs~le long d'une fibration est une équivalence faible,
\item N) les cofibrations et les \eqs~faibles sont stables par produit fibré au-dessus d'un objet
discret,
\item O) les cofibrations de $\mathcal{C}$ sont les monomorphismes,
\item P) les isomorphismes ne se factorisent qu'en isomorphismes dans la
sous-catégorie des cofibrations de $\mathcal{C}$,
\item Q) les colimites séquentielles transfinies de la sous-catégorie des
cofibrations de $\mathcal{C}$ existent et sont les colimites séquentielles
transfinies de $\mathcal{C}$,
\item R) pour toute cofibration $j:A\rightarrow B$ de $\mathcal{C}$, la classe contenant un représentant de chaque classe d'isomorphismes d'objets $C$ par lequel $j$ se
factorise dans la sous-catégorie des cofibrations est un ensemble,
\item S) les cofibrations génératrices de la \cmfcof~$\mathcal{C}$ ont leurs buts connexes,
\item T) les cofibrations génératrices de $\mathcal{C}$ ont leurs
sources et buts $\alpha$-petits,
\item U) il existe une \cat~$\bar{J}$ munie de deux objets notés $0$ et
$1$ et ayant les propriétés suivantes :
\item a) pour
toute \cat~facile $A$ et pour tout couple d'objets équivalents $a$ et $b$ de $A$, il
existe un morphisme de $\bar{J}$ vers $A$ envoyant $0$ sur $a$ et $1$ sur $b$,
\item b) pour toute \cat~$A$ et pour tout morphisme de $\bar{J}$ vers $A$,
les images de $0$ et $1$ par ce morphisme sont des objets équivalents dans
$A$,
\item c) notons $\bar{L}$, la sous-\cat~pleine de $\bar{J}$ d'objet $0$, alors
il existe un morphisme $p:\bar{J}\rightarrow \bar{L}$ envoyant $0$ et $1$ sur $0$ qui soit
une \eq~de \cats,
\item V) il existe dans $\mathcal{C-PC}$ une
\precat~$\bar{J}^{pre}$ $\alpha$-petite dont la catégorisation est $\bar{J}$
et vérifiant également la propriété a) de $\bar{J}$.\\
\\
Alors $\mathcal{C-PC}$ admet une structure de \cmfcof~avec pour \eqs~faibles les morphismes dont la
catégorisation est une \eq~de \cats, pour $\mathcal{I}$ l'ensemble constitué de $\emptyset\rightarrow *$ et des
flèches de type $Attach_n(g):\Delta[n]\Theta X\coprod_{\partial\Delta[n]\Theta X}\partial\Delta[n]\Theta
Y\rightarrow \Delta[n]\Theta Y $ pour $n>0$ et $g:X\rightarrow Y$ décrivant
l'ensemble des cofibrations génératrices de
$\mathcal{C}$, et pour
$\mathcal{J}$ l'ensemble des classes d'\eq~de morphismes à la fois cofibrations niveau par niveau et
\eqs~faibles dont les sources et buts sont petits par rapport à un cardinal
régulier bien choisi.
\end{theo}
{\it Preuve :}\\
D'après le lemme~\ref{reco} de reconnaissance de \cmfcof, il y a cinq axiomes
à montrer.\\ 
Tout d'abord il faut que $\mathcal{C-PC}$ soit complète et
cocomplète, ce qui est le cas d'après la proposition~\ref{precatclos}.\\
Le premier axiome est
la stabilité des \eqs~faibles par rétract donnée par le lemme~\ref{eqfre} ainsi
que l'axiome "trois pour deux" pour les \eqs~fai\-bles donné quant à lui par
le lemme~\ref{3pour2}.\\ 
Le deuxième axiome est que $\mathcal{I}$ et $\mathcal{J}$
vérifient l'argument du petit objet, ce qui est montré par les lemmes~\ref{I}
et~\ref{J}.\\ 
Le troisième axiome est que toute $\mathcal{J}$-cofibration est à la fois
une $\mathcal{I}$-cofibration et une \eq~faible, ce qui nous est donné par le
corollaire~\ref{p3}.\\
Le quatrième axiome est que tout morphisme $\mathcal{I}$-injectif est à la
fois $\mathcal{J}$-injectif et une \eq~faible, ce que montre le lemme~\ref{p4}.\\
Enfin le corollaire~\ref{p5} montre l'une des deux propriétés équivalentes du dernier axiome, à
savoir que les morphismes à la fois $\mathcal{I}$-cofibrations et \eqs~faibles
sont des $\mathcal{J}$-cofibrations.\\
CQFD.\\

Cependant il convient aussi de remarquer que certaines hypothèses permettent
d'en éliminer d'autres. En effet, l'hypothèse O) demande que les
cofibrations de $\mathcal{C}$ soient des monomorphismes. De ce fait l'hypothèse P) ne porte plus que sur les monomorphismes et elle est donc trivialement vraie. Il en est de même pour la première partie de l'hypothèse N). De plus l'hypothèse B) entraîne que la flèche $\emptyset\rightarrow *$ est un monomorphisme. L'hypothèse O) modifie aussi l'énoncé de l'hypothèse K), ce qui entraîne que les familles $\mathcal{F}_1$ et $\mathcal{F}_2$ sont des familles de monomorphismes. Comme A) postule que $\mathcal{C}$ est une \cmf, les cofibrations, qui par l'hypothèse O) sont exactement les monomorphismes, sont stables rétracts, sommes amalgamées et colimites séquentielles transfinies. Ainsi les $\mathcal{F}_1\cup\mathcal{F}_2\cup\{\emptyset\rightarrow *\}$-cofibrations sont des monomorphismes, ce qui rend triviale l'hypothèse H).
 En outre,
si dans l'hypothèse A) l'on demande directement que la \cmfcof~soit une
\cmfmono, alors l'hypothèse O) devient redondante. Par ailleurs, le fait
d'identifier les cofibrations et les monomorphismes changent la nature des
hypothèses Q) et R) qui ne portent donc plus sur la structure de \cmf~mais directement
sur la structure de catégorie de $\mathcal{C}$. Ainsi on peut trier ou modifier les
hypothèses ci-dessus obtenues par regroupement des résultats précédents
et donner de nouvelles définitions plus pertinentes au vu de ces
hypothèses.

\subsection{Remaniement des hypothèses du théorème central et nouvelles
définitions}

Tout d'abord, regroupons les hypothèses portant directement sur la caté\-gorie
de base $\mathcal{C}$. Comme on considère une donnée de Segal
proto-facile sur $\mathcal{C}$, il faut que $\mathcal{C}$ soit \discret~et possède
les produits fibrés de deux objets au-dessus d'un objet discret. Les
hypothèses B), C) et D) portent également sur $\mathcal{C}$ et les deux
dernières supposent que $\mathcal{C}$ est complète et cocomplète. Enfin
comme on l'a vu ci-dessus, le fait pour le théorème central d'identifier
cofibrations et monomorphismes entraîne que les hypothèses Q) et R) portent elles aussi
sur $\mathcal{C}$. Avec toutes ces hypothèses ainsi qu'avec la définition de
catégorie \discret, on peut définir une nouvelle notion de catégorie
\discret~directement utilisable pour le théorème central.

\begin{defin}\index{catégorie!bien discrétisante}
Une catégorie $\mathcal{C}$ est une catégorie bien \discret~si elle
vérifie les propriétés suivantes :
\item -elle est complète et cocomplète,
\item -pour tout morphisme $Z\rightarrow X$ vers un objet discret, on a 
$$ Z=\coprod_{x\in X}Z(x),$$
\item -$Enssj\circ Discret$ est naturellement isomorphe à l'identité de la
catégorie des ensembles,
\item -les objets discrets sont stables par limites dans
$\mathcal{C}$,
\item -il n'existe aucun morphisme de but l'ensemble vide qui ne soit pas un
isomorphisme,
\item -les réunions de deux sous-objets sont exactement les sommes
amalgamées de ces sous-objets au-dessus de leur intersection,
\item -les colimites séquentielles transfinies de la sous-catégorie des
monomorphismes de $\mathcal{C}$ existent et sont les colimites séquentielles
transfinies de $\mathcal{C}$,
\item -pour tout monomorphisme $j:A\rightarrow B$ de $\mathcal{C}$, la classe contenant un représentant de chaque classe d'isomorphismes d'objets $C$ par lequel $j$ se
factorise dans la sous-catégorie des monomorphismes est un ensemble.
\end{defin}

Regardons maintenant les hypothèses concernant uniquement la donnée de Segal. Ce sont
les hypothèses E) et F). Parmi ces hypothèses, l'hypothèse F) donne
une description concrète du foncteur $\tau_0$, aussi va-t-on donner la
définition de donnée de Segal concrète à une donnée de Segal basée
sur une catégorie bien \discret~et vérifiant les hypothèses E) et F).

\begin{defin}\index{donnée de Segal!concrète}
Une donnée de Segal $(\mathcal{C},\mathcal{C}_c,\mathcal{C}_{eq},\tau_0)$ est
dite concrète si :
\item -la catégorie sous-jacente $\mathcal{C}$ est bien \discret,
\item -les colimites séquentielles transfinies d'\obcs~le long de
monomorphismes sont des \obcs,
\item -le foncteur $\tau_0$ est tel que, pour tout \obc~$C$ de $\mathcal(C)$,
$\tau_0(C)$ est un quotient de l'ensemble des morphismes dans $\mathcal{C}$
de l'objet final vers $C$.
\end{defin}

Comme les hypothèses H), I) et J) sont déjà contenues dans la notion de donnée de Segal
pré-facile, on va obtenir la notion de donnée de
Segal pré-facilitée qui ne sera autre qu'une donnée de Segal concrète
 et pré-facile.

\begin{defin}\index{donnée de Segal!pré-facilitée}
Une donnée de Segal pré-facile $(\mathcal{C},\mathcal{F}_1,\mathcal{F}_2)$ est une
donnée de Segal facilitée si la donnée de Segal sous-jacente est concrète,
\end{defin}
  
On remarque que les hypothèses G), K), L) et M) sont des hypothèses
demandées dans la définition de donnée de Segal facile. L'hypothèse N)
n'y est qu'à moitié par sa partie concernant les \eqs~faibles. Quant aux
autres hypothèses de la donnée de Segal facile, elles sont déjà prises en
charge dans les définitions précédentes. Ainsi une donnée de
Segal facilitée sera simplement une donnée de Segal facile et
préfacilitée ou encore une donnée de Segal facile et concrète.

\begin{defin}\index{donnée de Segal!facilitée}
Une donnée de Segal facile $(\mathcal{C},\mathcal{F}_1,\mathcal{F}_2)$ est
dite donnée de Segal facilitée si la donnée de Segal pré-facile
sous-jacente est une donnée de Segal pré-facilitée.\\
Ou ce qui revient au même :\\
Une donnée de Segal facile $(\mathcal{C},\mathcal{F}_1,\mathcal{F}_2)$ est une
donnée de Segal facilitée si la donnée de Segal sous-jacente est concrète,
\end{defin}  

L'hypothèse A) demande à la \cmf~sous-jacente $\mathcal{C}$ d'une donnée
de Segal facilitée d'être une \cmfcof~simpliciale dont tous les objets sont
cofibrants et l'hypothèse O) demande que les cofibrations soient les
monomorphismes. On peut donc regrouper ces deux hypothèses en une, comme on
l'avait déjà fait remarquer plus haut, en demandant à la \cmf~$\mathcal{C}$
d'être une \cmf~simpliciale engendrée par monomorphismes. En revanche, on n'a plus à
demander que ses objets soient tous cofibrants, car ceci est déjà une
hypothèse de la donnée de Segal facile.\\

Si l'on regarde les hypothèses concernant la \cmfmono, on trouve les
hypothèses S) et T) qui portent sur les cofibrations génératrices de la
\cmfmono. Comme le fait de demander à la \cmf~de base $\mathcal{C}$ d'être
engendrée par monomorphisme se justifie en ceci que cela permet de simplifier la
démonstration de la structure de \cmf~sur les \precats~et leurs \eqs~faibles,
nous allons définir une notion de donnée de Segal facilitante (pour la
structure de \cmf) qui prendra en charge ces hypothèses d'engendrement par
monomorphismes.

\begin{defin}\index{donnée de Segal!facilitante}
Soit $(\mathcal{C},\mathcal{F}_1,\mathcal{F}_2)$ une donnée de Segal
facilitée, on dit qu'elle est facilitante si :
\item -la \cmf~sous-jacente $\mathcal{C}$ est \cmf\\ simpliciale engendrée par
monomorphismes,
\item -les cofibrations génératrices de la \cmfcof~$\mathcal{C}$ ont leurs buts connexes,
\item -les cofibrations génératrices de $\mathcal{C}$ ont leurs
sources et buts $\alpha$-petits.
\end{defin} 

Toutes les hypothèses ont été prises en charge par cette notion de
donnée de Segal facilitante à l'exception des deux dernières hypothèses
U) et V) qui concernent l'existence d'un bon intervalle $\bar{J}$ pour
représenter l'\eq~d'objets dans les \cats. Aussi va-t-on définir la notion
pour une donnée de Segal facilitante d'être munie d'un bon intervalle pour
les \cats.

\begin{defin}\index{bon intervalle pour les \cats}
Soit $(\mathcal{C},\mathcal{F}_1,\mathcal{F}_2)$ une donnée de Segal
facilitante, on dit qu'elle admet un bon intervalle pour les \cats~si :
\item -il existe une \cat~$\bar{J}$ munie de deux objets notés $0$ et
$1$ et ayant les propriétés suivantes :
\item a) pour
toute \cat~facile $A$ et pour tout couple d'objets équivalents $a$ et $b$ de $A$, il
existe un morphisme de $\bar{J}$ vers $A$ envoyant $0$ sur $a$ et $1$ sur $b$,
\item b) pour toute \cat~$A$ et pour tout morphisme de $\bar{J}$ vers $A$,
les images de $0$ et $1$ par ce morphisme sont des objets équivalents dans
$A$,
\item c) notons $\bar{L}$, la sous-\cat~pleine de $\bar{J}$ d'objet $0$, alors
il existe un morphisme $p:\bar{J}\rightarrow \bar{L}$ envoyant $0$ et $1$ sur $0$ qui soit
une \eq~de \cats,
\item -il existe dans $\mathcal{C-PC}$ une
\precat~$\bar{J}^{pre}$ $\alpha$-petite dont la catégorisation est $\bar{J}$
et vérifiant également la propriété a) de $\bar{J}$.
\end{defin}

Avec toutes ces définitions, nous pouvons donner un énoncé plus simple du
théorème central, puisqu'il suffit d'être une donnée de Segal
facilitante munie d'un bon intervalle pour les \cats~afin de vérifier toutes
les hypothèses de l'ébauche du théorème central.

\newpage

\subsection{Enoncé du théorème central}

Voici donc le nouvel énoncé du théorème central qui est bien plus simple
que celui de l'ébauche du fait qu'il prenne en compte les nouvelles
définitions dégagées précédemment de l'analyse de la foule des
hypothèses de l'ébauche du théorème central. 

\begin{theo}[central]\index{catégorie de modèles fermée!des \precat}
Soit $(\mathcal{C},\mathcal{F}_1,\mathcal{F}_2)$ une donnée de Segal
facilitante admettant un bon intervalle pour les \cats. 
Alors $\mathcal{C-PC}$ admet une structure de \cmfmono~avec pour \eqs~faibles
les morphismes dont la catégorisation est une \eq~de \cats, les ensembles
générateurs des cofibrations et des cofibrations triviales étant
respectivement les ensembles $\mathcal{I}$ et $\mathcal{J}$.
\end{theo}
 
La preuve est bien évidemment celle de l'ébauche du théorème central.
Bien que cet énoncé soit plus agréable à l'oeil que celui de
l'ébauche, il est cependant extrêmement dense car derrière la
notion de donnée de Segal facilitante se cache une multitude d'hypothèses
que l'ébauche du théorème centrale, ô combien rebutante, a néanmoins
l'avantage d'exposer clairement. Après le balayage des hypothèses triviales des sections précédentes, nous sommes en mesure de donner un meilleur énoncé de notre ébauche. Aussi c'est en fait le nouvel énoncé de l'ébauche du théorème
central qui nous servira par la suite pour vérifier qu'une donnée de Segal
proto-facile vérifie les hypothèses du théorème central.\\

\begin{theo}[-central (ébauche 2)]
Soit $(\mathcal{C},\mathcal{F}_1,\mathcal{F}_2)$ une donnée de Segal proto-facile
vérifiant les propriétés suivantes :
\item A') la catégorie sous-jacente $\mathcal{C}$ est une \cmfmono~simpliciale,
\item B') il n'existe aucun morphisme de but l'ensemble vide qui ne soit pas un
isomorphisme,
\item C') les réunions de deux
sous-objets sont exactement les sommes amalgamées de ces sous-objets au-dessus
de leur intersection,
\item D') les objets discrets sont stables par limites dans
$\mathcal{C}$,
\item E') les colimites séquentielles transfinies
d'\obcs~le long de monomorphismes sont des \obcs,
\item F') le foncteur $\tau_0$ est tel que, pour tout \obc~$C$ de
$\mathcal{C}$, $\tau_0(C)$ est un quotient de l'ensemble des morphismes dans
$\mathcal{C}$ de l'objet final vers $C$,
\item G') il existe un \obc~contractile $K$ muni de deux
morphismes de l'objet final vers $K$ notés 0 et 1, et ayant la propriété
suivante : pour tout couple $(f,g)$ de morphismes de l'objet final vers un \obc~$C$ tel que leurs images par $\tau_0$ soient égales, il existe un
morphisme de $K$ vers $C$ envoyant 0 sur $f$ et 1 sur $g$,
\item H') les sources et buts des flèches de la famille $\mathcal{F}_1$, ainsi
que les buts de la famille $\mathcal{F}_2$, sont
connexes et non isomorphes à l'ensemble vide,
\item I') il existe un cardinal transfini régulier strictement supérieur à
$\aleph_0$, que l'on notera $\alpha$, pour lequel tout morphisme de
$\mathcal{F}_1$ et $\mathcal{F}_2$ a sa source et son but $\alpha$-petits au
sens~\ref{alphapetit},
\item J') la famille $\mathcal{F}_1$ est incluse dans la classe des cofibrations
triviales de la \cmfmono~$\mathcal{C}$ et la famille $\mathcal{F}_2$ dans celle des
cofibrations,
\item K') les \eqcs~d'\obcs~sont exactement les \eqs~faibles de la
\cmf~$\mathcal{C}$ entre \obcs,
\item L') le produit fibré d'une équivalence faible entre
\obcs~le long d'une fibration est une équivalence faible,
\item M') les \eqs~faibles sont stables par produit fibré au-dessus d'un objet
discret,
\item N') les colimites séquentielles transfinies de la sous-catégorie des
monomorphismess de $\mathcal{C}$ existent et sont les colimites séquentielles
transfinies de $\mathcal{C}$,
\item O') pour tout monomorphisme $j:A\rightarrow B$ de $\mathcal{C}$, la classe contenant un représentant de chaque classe d'isomorphismes d'objets $C$ par lequel $j$ se
factorise dans la sous-catégorie des monomorphismes est un ensemble,
\item P') les cofibrations génératrices de la \cmfcof~$\mathcal{C}$ ont leurs buts connexes,
\item Q') les cofibrations génératrices de $\mathcal{C}$ ont leurs
sources et buts $\alpha$-petits,
\item R') il existe une \cat~$\bar{J}$ munie de deux objets notés $0$ et
$1$ et ayant les propriétés suivantes :
\item a) pour
toute \cat~facile $A$ et pour tout couple d'objets équivalents $a$ et $b$ de $A$, il
existe un morphisme de $\bar{J}$ vers $A$ envoyant $0$ sur $a$ et $1$ sur $b$,
\item b) pour toute \cat~$A$ et pour tout morphisme de $\bar{J}$ vers $A$,
les images de $0$ et $1$ par ce morphisme sont des objets équivalents dans
$A$,
\item c) notons $\bar{L}$, la sous-\cat~pleine de $\bar{J}$ d'objet $0$, alors
il existe un morphisme $p:\bar{J}\rightarrow \bar{L}$ envoyant $0$ et $1$ sur $0$ qui soit
une \eq~de \cats,
\item S') il existe dans $\mathcal{C-PC}$ une
\precat~$\bar{J}^{pre}$ $\alpha$-petite dont la catégorisation est $\bar{J}$
et vérifiant également la propriété a) de $\bar{J}$.\\
\\
Alors $\mathcal{C-PC}$ admet une structure de \cmfcof~avec pour \eqs~faibles les morphismes dont la
catégorisation est une \eq~de \cats, pour $\mathcal{I}$ l'ensemble constitué de $\emptyset\rightarrow *$ et des
flèches de type $Attach_n(g):\Delta[n]\Theta X\coprod_{\partial\Delta[n]\Theta X}\partial\Delta[n]\Theta
Y\rightarrow \Delta[n]\Theta Y $ pour $n>0$ et $g:X\rightarrow Y$ décrivant
l'ensemble des cofibrations génératrices de
$\mathcal{C}$, et pour
$\mathcal{J}$ l'ensemble des classes d'\eq~de morphismes à la fois cofibrations niveau par niveau et
\eqs~faibles dont les sources et buts sont petits par rapport à un cardinal
régulier bien choisi.
\end{theo}

Maintenant que nous avons obtenu une belle structure de \cmf~sur les \precats, nous
allons voir un exemple de \precats~ayant cette structure par vérification des
hypothèses du théorème central.

\chapter{Catégories de Segal}

\newpage

Comme on l'a vu dans l'exemple~\ref{hypenssimp1}, la
catégorie des ensembles simpliciaux, avec pour \obcs~les ensembles
simpliciaux quelconques, pour \eqcs~d'\obcs~les équivalences faibles
d'ensembles simpliciaux (i.e. les morphismes dont les réalisations
géométriques induisent des isomorphes sur
les groupes d'homotopie) et pour $\tau_0$ la composée des foncteurs
réalisation géométrique
et composante connexe, forme une donnée de Segal. Les \precats~pour cette
donnée de Segal sont les précatégories de Segal, les \cats~sont les catégories de Segal et les \eqs~de \cats~sont
les \eqs~de catégories de Segal, ces trois notions étant celles définies
dans \cite{h-s}.\\

En outre avec pour $\mathcal{F}_1$ la famille vide et pour $\mathcal{F}_2$ la
famille des inclusions des bords des simplexes standards dans ces derniers, la
catégorie des ensembles simpliciaux forme une donnée de Segal proto-facile, ce que
l'on a vu dans l'exemple~\ref{hypenssimp2}. Pour cette donnée de Segal proto-facile,
les \obcs~faciles sont les ensembles simpliciaux quelconques et les
\eqcs~faciles sont les \eqs~faibles d'ensembles simpliciaux qui sont aussi des
fibrations de Kan.\\

Nous allons montrer que cette donnée de Segal proto-facile vérifie les
hypothèses du théorème central afin d'obtenir une structure de \cmf~sur
les précatégorie de Segal.

\newpage

\section{Vérification des premières hypothèses du théorème central}

Afin de montrer que le théorème central s'applique au cas des
précatégories de Segal, nous devons montrer que les dix-neuf hypothèses
de l'ébauche 2 du théorème central sont vérifiées par la donnée de
Segal proto-facile de
l'exemple~\ref{hypenssimp2}. Nous commencerons par montrer les dix-sept premières
hypothèses, ce qui nous permettra par la suite d'utiliser tous les résultats
de cette thèse ne faisant pas intervenir l'intervalle $\bar{J}$.

\begin{lem}\label{seghyp1}
La donnée de Segal proto-facile sur les ensembles simpliciaux définies dans
l'exemple~\ref{hypenssimp2} vérifie les hypothèses A') à Q') de l'ébauche 2 du théorème
central.
\end{lem}
{\it Preuve :}\\
L'hypothèse A') est un résultat classique sur les ensembles simpliciaux, voir
en particulier \cite{h}. En
effet ceux-ci forment une \cmf~simpliciale dont tous les objets sont cofibrants
avec pour cofibrations les monomorphismes et pour \eqs~faibles les morphismes
dont la réalisation géométrique induit des isomorphismes sur les groupes
d'homotopie, or comme ceux-ci ne sont autres que les \eqcs~d'\obcs, ceci montre
l'hypothèse K'). De plus cette \cmf~est engendrée de manière cofibrante avec
pour cofibrations génératrices les inclusions des bords des simplexes
standards dans ces derniers et pour cofibrations triviales génératrices les
inclusions des cornes des simplexes standards dans ces derniers.\\
\\
Comme la catégorie des ensembles simpliciaux est une catégorie de
préfaisceaux d'ensembles dont les limites et colimites sont calculées niveau
par niveau, les hypothèses B'), C') et D'), étant vérifiées par la catégorie des
ensembles, sont donc vraies pour les ensembles simpliciaux. En outre comme on
l'a vu ci-dessus, les cofibrations sont les monomorphismes. Or ceux-ci sont en
fait les injections niveau par niveau, donc l'hypothèse N'), étant vérifiée
avec les injections ensemblistes, est aussi vraie avec les cofibrations
d'ensembles simpliciaux.\\
\\
L'hypothèse E') est triviale car les \obcs~sont les ensembles simpliciaux
quelconques et que la catégorie des ensembles simpliciaux est cocomplète.\\
L'hypothèse F') découle de la définition du foncteur composante connexe
$\pi_0$ et du fait que les points de la réalisation géométrique d'un
ensemble simplicial sont exactement les morphismes de $\Delta[0]$ vers cet
ensemble simplicial.\\
L'hypothèse G') est vérifiée facilement si on prend pour $K$ l'intervalle
$\bar{I}$, i.e. le nerf du groupoïde à deux objets et un unique
isomorphisme entre ceux-ci.\\
\\
Pour les hypothèses  H'), I') et J'), la partie concernant la famille
$\mathcal{F}_1$ est triviale car $\mathcal{F}_1$ est vide. La famille
$\mathcal{F}_2$ est formée par les inclusions des bords des simplexes
standards dans ces derniers. Comme les cofibrations sont les monomorphismes, la
famille $\mathcal{F}_2$ est bien incluse dans la classe des cofibrations, ce qui
montre la partie de l'hypothèse J') concernant $\mathcal{F}_2$. Les buts des flèches de $\mathcal{F}_2$ sont les simplexes
standards. Or ceux-ci sont connexes au sens de l'exemple~\ref{enssimpconnex} et
non vides d'après l'exemple~\ref{enssimpconnex2}, ce qui montre la partie de
l'hypothèse H') concernant $\mathcal{F}_2$. Enfin les sources et buts des
flèches de $\mathcal{F}_2$, étant les simplexes standards et leurs bords,
sont $\alpha$-petits pour $\alpha=2^{\aleph_0}$ d'après
l'exemple~\ref{enssimppet}, ce qui montre la partie de
l'hypothèse I') concernant $\mathcal{F}_2$.\\
\\
Comme les cofibrations génératrices ne sont autres que les inclusions des
bords des simplexes standards dans ces derniers, l'ensemble des cofibrations
génératrices est donc $\mathcal{F}_2$. Ainsi comme $\mathcal{F}_2$ vérifie
les hypothèses H') et J'), alors l'ensemble des cofibrations génératrices
vérifie les hypothèses P') et Q').\\
\\
Pour montrer l'hypothèse L'), considérons le diagramme cartésien suivant :
\begin{diagram}
A & \rTo & B \\
\dTo & & \dTo^{f}_{\sim} \\
C & \rOnto^{g} & D \\
\end{diagram}
Appliquons le foncteur réalisation géométrique à ce diagramme. Comme le
foncteur réalisation géométrique préserve les produits fibrés, on
obtiendra encore un diagramme cartésien dans la catégorie des espaces
topologiques. Comme $f$ est une \eq~faible d'ensembles simpliciaux, sa
réalisation géométrique est une \eq~d'homotopie faible. Comme $g$ est une
fibration de la \cmf~des ensembles simpliciaux, sa réalisation géométrique est une
fibration de la \cmf~des espaces topologiques. Comme dans la \cmf~des espaces
topologiques tous les objets sont fibrants, cette \cmf~est propre à droite.
Ainsi le produit fibré de l'\eq~d'homotopie faible $|f|$ entre objets fibrants
le long de la fibration $|g|$ est une \eq~d'homotopie faible. Or ce produit
fibré n'est autre que la réalisation géométrique de $A\rightarrow C$ par
préservation du produit fibré par réalisation géométrique, ce qui montre
que $A\rightarrow C$ est bien une \eq~faible et donc que l'hypothèse L') est
vérifiée.\\
\\
Le foncteur réalisation géométrique et les foncteurs des groupes
d'homotopie $\pi_i$ préservent les produits fibrés. Ainsi la réalisation
géométrique du produit fibré de deux \eqs~faibles n'est autre que le
produit fibré de deux \eqs~d'homotopie faibles. En appliquant à ce dernier
les foncteurs $\pi_i$, on obtient des produits fibrés de deux isomorphismes
qui sont donc eux-mêmes des isomorphismes. Ainsi la réalisation
géométrique du produit fibré de deux \eqs~faibles induit des isomorphismes
sur les groupes d'homotopie et donc le produit fibré des deux \eqs~faibles est
une \eq~faible. Ceci montre l'hypothèse M').\\
\\
Pour terminer les vérifications d'hypothèses, il reste juste à montrer
que l'hypothèse O') est vraie et ceci n'est autre que le lemme~\ref{ast1}, car
comme on l'a fait remarquer plus haut la catégorie des ensembles simpliciaux
est une catégorie de préfaisceaux d'ensembles et que la catégorie des
ensembles vérifie l'hypothèse O').\\
CQFD.\\

Il ne nous reste donc plus qu'à montrer l'existence de l'intervalle $\bar{J}$ pour
les précatégories de Segal pour pouvoir appliquer le théorème central au
cas des ensembles simpliciaux. Toutefois l'existence d'un tel intervalle
nécessite la notion de groupoïde de Segal ainsi que l'\eq~des théories
homotopiques de ces groupoïdes de Segal et des espaces topologiques.

\newpage 

\section{Groupoïde de Segal}

Pour montrer l'existence d'un bon intervalle $\bar{J}$ pour les \cats~de Segal,
nous allons introduire, plus précisément rappeler, la notion de groupoïde
de Segal. En effet l'intervalle $\bar{J}$ a pour vocation de représenter la
notion d'\eq~d'objets dans une catégorie de Segal. Comme deux objets sont
équivalents s'il existe un morphisme de l'un à l'autre qui est inversible
à homotopie près, on voit apparaître le lien entre $\bar{J}$ et la
notion de groupoïde de Segal : l'intervalle $\bar{J}$ devra être un
groupoïde de Segal. Aussi avant de construire l'intervalle $\bar{J}$ et de
montrer que c'est un groupoïde de Segal, nous allons rappeler la
définition de groupoïde de Segal et passer en revue les principales
propriétés de ces derniers.

\subsection{Groupes d'homotopie d'un groupoïde de Segal}

Commençons tout d'abord par rappeler la notion de groupoïde de Segal définie dans
\cite{h-s}.

\begin{defin}\index{groupoïde!de Segal}
Une catégorie de Segal $A$ est un groupoïde de Segal si $\tau_1(A)$ est le
nerf d'un groupoïde au sens des catégories. Les morphismes de groupoïdes de Segal sont les morphismes de précatégories de Segal et les
\eqs~de groupoïdes de Segal sont les équivalences de catégories
de Segal.
\end{defin}

Si la notion d'\eq~pour les groupoïdes de Segal est la même que pour les
catégories de Segal, il existe cependant une caractérisation de l'\eq~entre
groupoïdes de Segal par l'intermédiaire de la notion de groupes
d'homotopie d'un groupoïde de Segal.

\begin{defin}\index{groupes d'homotopie d'un groupoïde!de Segal}
Soit $A$ un groupoïde de Segal, on définit ses groupes d'homotopie de la
manière suivante :
\item - $\pi_0(A)=\tau_0(A)$,
\item - pour tout objet $a$ de $A$, $\pi_1(A,a)=\pi_0(|A_1(a,a)|)$
\item - pour $i>1$ et pour tout objet $a$ de $A$,
$\pi_i(A,a)=\pi_{i-1}(|A_1(a,a)|,Id_a)$.
\end{defin}

Ces groupes d'homotopie de groupoïdes de Segal sont bien fonctoriels comme
composées de foncteurs et, pour $i>1$, ils sont bien munis d'une structure de
groupes. Afin de montrer que ces groupes d'homotopie caractérisent les
\eqs~entre groupoïdes de Segal, nous allons introduire la notion de
groupoïde de Segal simplement connexe pour laquelle une telle
caractérisation a lieu.

\begin{defin}\index{simple connexité des groupoïdes de Segal}
Un groupoïde de Segal est simplement connexe si son $\pi_0$ et tous ses $\pi_1$ sont
triviaux.
\end{defin}

\begin{lem}\label{eqgp1}
Soit $f:A\rightarrow B$ un morphisme de groupoïdes de Segal.\\
Si $A$ et $B$ sont simplement connexes, alors $f$ est
une \eq~de catégories de Segal si et seulement si $\pi_0(f)$ est une bijection
d'ensembles et que, pour tout $i>0$ et pour tout objet $a$ de $A$, $\pi_i(f,a)$
est un isomorphisme.
\end{lem}
{\it Preuve :}\\
On a déjà vu que $f$ est essentiellement surjective si et seulement
si $\tau_0(f)$ est une bijection ensembliste. Or $\tau_0(f)$ n'est autre que $\pi_0(f)$.
Ainsi l'essentielle surjectivité est caractérisée en terme de $\pi_0$.\\
\\
Par définition, $f$ est pleinement fidèle si et seulement si $f_1(x,y)$ est
une \eq~faible d'ensembles simpliciaux pour tout couple $(x,y)$ d'objets de $A$.
Or $f_1(x,y)$ est une \eq~faible d'ensembles simpliciaux si et seulement si
$\pi_0(|f_1(x,y)|)$ est une bijection et, pour
tout $j>0$ et pour tout objet $g$ de $A_1(x,y)$, $\pi_j(|f_1(x,y)|,g)$ est
un morphisme de groupes. Ainsi on obtient que si $f$ est pleinement fidèle,
alors, pour tout objet $a$ de $A$, $f_1(a,a)$ est une \eq~faible donc
$\pi_0(|f_1(a,a)|)$ est une bijection et, pour tout
$j>0$, les $\pi_j(|f_1(a,a),Id_a|)$ sont des isomorphismes de groupes,
c'est-à-dire que, pour tout $i>0$, les $\pi_i(f_1(a,a),a)$ sont des
isomorphismes.\\
\\
Montrons que si pour tout $i>0$ et pour tout objet $a$ de $A$, $\pi_i(f,a)$
est un isomorphisme et que $\pi_0(f)$ est bijective, alors $f$ est pleinement fidèle. Pour cela
montrons que, pour tout couple $(x,y)$ d'objets de $A$, $f_1(x,y)$ est une
\eq~faible. Deux cas se présentent.\\ 
\\
Si $A_1(x,y)$ est vide, alors
nécessairement $B_1(f(x),f(y))$ est vide. Sinon il existe un morphisme entre
$f(x)$ et $f(y)$ et, comme $B$ est un groupoïde de Segal, ce morphisme est une
\eq~entre $f(x)$ et $f(y)$, qui ont donc même image dans $\tau_0(B)$. Or
$\tau_0(f)$ est une bijection par hypothèse, donc $x$ et $y$ ont même image
dans $\tau_0(A)$ ce qui contredit la vacuité de $A_1(x,y)$.\\
\\
Si $A_1(x,y)$ est non vide, alors il existe un objet $g$ de $A_1(x,y)$ qui est
une \eq~entre $x$ et $y$ car $A$ est un groupoïde de Segal. En utilisant la
démonstration du lemme~\ref{3pour2c}, ceci entraîne que dans la catégorie
homotopique des ensembles simpliciaux les remplacements fibrants $A^f_1(x,y)$ et
$A^f_1(x,x)$ sont équivalents. De même, par l'\eq~$f(g)$ dans $B$ entre $f(x)$
et $f(y)$, on obtient que $B^f_1(f(x),f(y))$ et $B^f_1(f(x),f(x))$ sont
équivalents. Et par l'axio\-me "trois pour deux" dans la \cmf~des ensembles
simpliciaux, il vient que le remplacement
fibrant de $f_1(x,y)$ est une \eq~faible si et seulement si celui de $f_1(x,x)$
l'est. Or tout ensemble simplicial est équivalent à son remplacement
fibrant, d'où toujours par l'axiome "trois pour deux", les remplacements
fibrants de $f_1(x,y)$ et de $f_1(x,x)$ sont des \eqs~faibles si et seulement si
$f_1(x,y)$ et $f_1(x,x)$ le sont. Par transitivité de l'\eq~logique, on a donc obtenu que $f_1(x,y)$ est une
\eq~faible si et seulement si $f_1(x,x)$ l'est.\\ 
\\
Ainsi pour montrer que $f$ est pleinement fidèle, on en est réduit à
montrer que $f_1(x,x)$ est une \eq~faible pour tout objet $x$ de $A$. Pour cela
il faut montrer que $|f_1(x,x)|$ est une \eq~d'homotopie faible, i.e. que
$\pi_0(|f_1(x,x)|)$ est une bijection, ce qui est vrai car par hypothèse
$\pi_1(f,x)=\pi_0(|f_1(x,x)|)$ est une bijection, et que, pour
tout objet $h$ de $A_1(x,x)$ et pour tout $j>0$, les $\pi_j(|f_1(x,x)|,h)$
sont des isomorphismes de groupes. On sait que c'est vrai quand $h$ vaut
l'identité de $x$, car par hypothèse les $\pi_j(f,x)=\pi_j(|f_1(x,x)|,Id_x)$
sont des isomorphismes de groupes. Or $A$ est simplement connexe donc il existe un chemin
$\gamma$ de $Id_x$ vers $h$ qui induit un isomorphisme de
$\pi_j(|A_1(x,x),Id_x|)$ vers $\pi_j(|A_1(x,x),h|)$, pour $j>0$. L'image de ce chemin par
$|f_1(x,x)|$ induit un isomorphisme de $\pi_j(|B_1(f(x),f(x)),Id_{f(x)}|)$ vers
$\pi_j(|B_1(f(x),f(x)),f(h)|)$, pour $j>0$. On obtient alors le diagramme commutatif suivant
:
\begin{diagram}
\pi_j(|A_1(x,x),Id_x|) & \rTo^{\pi_{j+1}(f,x)} & \pi_j(|B_1(f(x),f(x)),Id_{f(x)}|)\\
\dTo^{\gamma*.*\gamma^{-1}} & & \dTo_{|f_1(x,x)|(\gamma})*.*|f_1(x,x)|(\gamma)^{-1} \\
\pi_j(|A_1(x,x),h|) & \rTo_{\pi_j(|f_1(x,x),h|)} & \pi_j(|B_1(f(x),f(x)),f(h)|)\\
\end{diagram}
Par hypothèse, $\pi_{j+1}(f,x)$ est un isomorphisme de groupes et par ce qui
précède les flèches verticales sont aussi des isomorphismes de groupes. On
en déduit donc que $\pi_j(|f_1(x,x),h|)$ est bien un isomorphisme de
groupes, pour $j>0$. Donc $f_1(x,x)$ est une \eq~faible pour tout objet $x$ de
$A$ et, par ce qui précède, ceci entraîne que pour tout couple $(x,y)$
d'objets de $A$, $f_1(x,y)$ est une \eq~faible et par conséquent que $f$ est
pleinement fidèle.\\
CQFD.\\

On remarque en regardant la démonstration de près que pour $f:A\rightarrow B$ morphisme
de groupoïdes de Segal quelconques, on a :\\
- l'essentielle surjectivité de $f$ équivaut à la bijectivité de
$\pi_0(f)$,\\
- si $f$ est une \eq~de catégories de Segal, alors $f$ induit des
isomorphismes sur tous les $\pi_i$.\\
En revanche la réciproque n'est vraie que pour les groupoïdes de Segal
simplement connexes.\\

Après ces quelques rappels, nous allons maintenant passer à la définition de
l'intervalle $\bar{J}$.

\subsection{Proto-groupoïde de Segal}

Pour définir l'intervalle $\bar{J}$, nous allons en fait passer par la
construction de la précatégorie de Segal $\bar{J}^{pre}$ dont la
catégorisation sera justement $\bar{J}$.  Or cette précatégorie de Segal
$\bar{J}^{pre}$ doit représenter elle-aussi l'\eq~d'objets dans les
catégories de Segal, donc, comme on l'a fait remarquer plus haut, elle
représente en fait les morphismes inversibles à homotopie près. Cependant
cette notion n'a de sens que dans les catégories de Segal, ce que
$\bar{J}^{pre}$ n'est pas. C'est la raison pour laquelle nous allons définir
une notion de morphisme quasi-inversibles dans une précatégorie de Segal
ainsi qu'une notion de
proto-groupoïde de Segal qui permettront de comprendre ce qu'est
un élément d'une précatégorie de Segal qui aura vocation à être
inversible dans la catégorisation de cette précatégorie, mais aussi ce qu'est une
précatégorie de Segal dont la catégorisation est un groupoïde de
Segal. 

\begin{defin}\index{proto-groupoïde de Segal}
Une précatégorie de Segal $A$ est un proto-groupoïde de Segal si, pour
tout élément $u$ de $A_1$ qui n'est pas l'arête extérieure d'un
simplexe, il existe un élément $v$ de $A_1$, deux éléments $T_1$ et
$T_2$ de $A_2$ et deux morphismes $\alpha$ et $\beta$ de l'intervalle $\bar{I}$
dans $A_1$ tels que :
\item -$u=\delta_{0,1}(T_1)=\delta_{1,2}(T_2)$ et
$v=\delta_{1,2}(T_1)=\delta_{0,1}(T_2)$,
\item -$\alpha(0)=\delta_{0,2}(T_1)$ et $\alpha(1)=Id_{s(u)}$,
\item -$\beta(0)=\delta_{0,2}(T_2)$ et $\beta(1)=Id_{s(v)}$.\\
\\
On dira alors que l'élément $u$ est proto-inversible.
\end{defin}

Il est clair que tout groupoïde de Segal est un proto-groupoïde de
Segal. Nous allons maintenant montrer que les précatégories de Segal
proto-groupoïdiques sont bien des précatégories de Segal dont la
catégorisation est un groupoïde de Segal. Nous allons montrer ce
résultat en plusieurs étapes. Supposons dans un premier temps que la
catégorisation préserve la notion de proto-groupoïde, alors le
résultat de la catégorisation est une catégorie de Segal proto-groupoïdique. Montrons donc qu'une telle catégorie de Segal est en fait un groupoïde de Segal.

\begin{lem}
Toute catégorie de Segal proto-groupoïdique est un groupoïde de Segal.
\end{lem}
{\it Preuve :}\\
Soit $A$ une catégorie de Segal proto-groupoïdique. Soit $\phi$ un
morphisme de la catégorie $\tau_1(A)$. Deux cas se présentent.\\
\\
Si $\phi$ n'est pas une composée de morphismes de $\tau_1(A)$, alors ses représentants
$f$ dans $A_1$ non plus. Et comme $A$ est un proto-groupoïde, $f$ est
proto-inversible. Donc il existe $v, T_1, T_2, \alpha, \beta$ le rendant
proto-inversible. Or les morphismes $\alpha,\beta:\bar{I}\rightarrow A_1$
signifient par propriété de $\bar{I}$ que $\delta_{0,2}(T_1)\sim Id_{s(f)}$ et
$\delta_{0,2}(T_2)\sim Id_{b(f)}$. Donc dans $\tau_1(A)$, on obtient que
$\bar{v}\circ\bar{f}=Id_{b(f)}$ et $\bar{f}\circ\bar{v}=Id_{s(f)}$, i.e. que
$\phi=\bar{f}$ est inversible.\\
\\
Si $\phi$ est une composée, il existe un $n$-uplet de morphismes 
$(\phi_1,\ldots,\phi_n)$ qui ne sont pas des composés et tels que $\phi$ en soit
la composée.  Comme $A$ est proto-groupoïdique, d'après ce qui
précède, les $\phi_i$ sont des isomorphismes pour $i$ de $1$ à $n$. Donc leur
composée $\phi$ aussi.\\
\\
Ainsi tous les morphismes de $\tau_1(A)$ sont des isomorphismes et donc
$\tau_1(A)$ est un groupoïde.\\
CQFD.\\

Avec ce résultat, il ne reste donc plus qu'à prouver que la catégorisation
préserve la notion de proto-groupoïde pour montrer que la catégorisation
d'un proto-groupoïde de Segal est un groupoïde de Segal.

\begin{lem}
Si $A$ est un proto-groupoïde de Segal, alors $Cat(A)$ est un groupoïde
de Segal.
\end{lem}
{\it Preuve :}\\
Le procédé de catégorisation consiste en une colimite séquentielle transfinie de sommes
amalgamées par des flèches génératrices. Celles-ci sont de deux types.
Comme, pour la donnée de Segal que l'on considère sur les ensembles
simpliciaux, la famille $\mathcal{FG}1$ est vide, regardons de plus près la
famille $\mathcal{FG}2$. Cette dernière est formée par les flèches de type
$Boit_m(f)$ où $f$ appartient à $\mathcal{F}_2$. En outre tous les
morphismes du but de $Boit_m(f)$ vers $A$ (ou  vers l'une des étapes de la
catégorisation) correspondent à des diagrammes du type suivant :
\begin{diagram}
\partial\Delta[n] & \rTo & A_m\\
\dTo & & \dTo \\
\Delta[n] & \rTo & A_1\times_{A_0}\ldots\times_{A_0} A_1\\
\end{diagram}
Or les éléments de $A_1$ apparaissant dans ce diagramme sont les projections
des images dans $A_1\times_{A_0}\ldots\times_{A_0} A_1$ des objets de $\partial\Delta[n]$ et de
$\Delta[n]$. Comme $\partial\Delta[n]$ et $\Delta[n]$ ont les mêmes objets
pour $n>0$, seules les flèches de type $Boit_m(f)$ avec $f:\emptyset\rightarrow *$
modifient (en en rajoutant) des éléments de $A_1$. Plus explicitement,
$Boit_m(\emptyset\rightarrow *)$ n'est autre que $\Upsilon(m)\rightarrow \Delta[m]$
dont le rajout équivaut à donner une composée à un $m$-uplet de morphismes
composables. On en conclut donc que les éléments de $Cat(A)_1$ correspondent
tous à des $n$-uplets d'éléments de $A_1$ non composés d'autres
éléments.\\
\\
Soit $\phi$ un morphisme de $\tau_1(Cat(A))$. Par ce qui précède, c'est la
composée d'un $n$-uplet $(\phi_1,\ldots,\phi_n)$ de morphismes composables de $\tau_1(A)$
qui ne sont pas des composés d'autres morphismes. Or
ces derniers sont inversibles dans $\tau_1(A)$ car $A$ est un proto-groupoïde de Segal. Donc ils sont aussi inversibles dans
$\tau_1(Cat(A))$. Comme toute composée d'inversibles est inversible, $\phi$
est inversible, donc $\tau_1(Cat(A))$ est un groupoïde.\\
CQFD.\\

Avec cette notion de proto-groupoïde de Segal, on peut mieux formuler l'existence et les
propriétés de la précatégorie de Segal $\bar{J}^{pre}$. En outre si l'on
demande que $\bar{J}^{pre}$ soit un proto-groupoïde de Segal, alors, par le
résultat que l'on vient de montrer, sa catégorisation $\bar{J}$ sera bien un
groupoïde de Segal.

\begin{prop}[-définition]
Il existe un proto-groupoïde de Segal $\bar{J}^{pre}$ ayant deux objets
notés $0$ et $1$ et un morphisme noté $u$ de $0$ à $1$ et qui vérifie la propriété
suivante :\\
pour toute précatégorie de Segal $A$ et pour tout élément $f$ de $A_1$
proto-inversible, il existe un morphisme de $\bar{J}^{pre}$ vers $A$ envoyant
$u$ sur $f$.
\end{prop}
{\it Preuve :}\\
Donnons une construction pour $\bar{J}^{pre}$. On commence par prendre le
coproduit de deux $\Delta[0]$ notés $0$ et $1$. On prend ensuite la somme
amalgamée de ce coproduit avec un coproduit de deux $\Delta[1]$ notés $u$ et
$v$ sur un coproduit de deux $\Delta[0]\coprod \Delta[0]$ telle que $s(u)=0=b(v)$
et $b(u)=1=s(v)$. Puis prenons la somme amalgamée de ce qui précède par un
coproduit de deux $\Delta[2]$ notés $T_1$ et $T_2$ sur un coproduit de deux
$\Delta[1]\coprod\Delta[1]$ telle que $\delta_{0,1}(T_1)=u=\delta_{1,2}(T_2)$ et
$\delta_{1,2}(T_1)=v=\delta_{0,1}(T_2)$. Enfin considérons la somme
amalgamée de la construction précédente par un coproduit de deux
$\Delta[1]\Theta \bar{I}$ notés $\alpha$ et $\beta$ sur un coproduit de deux
$\Delta[1]\coprod \Delta[1]$ telle que $\alpha\rightarrow (\delta_{0,2}(T_1),
Id_0)$ et $\beta\rightarrow (\delta_{0,2}(T_2),Id_1)$. Il est facile de voir
que $\bar{J}^{pre}$ ainsi défini est bien un proto-groupoïde de Segal. La
propriété de $\bar{J}^{pre}$ est également évidente car il suffit
d'envoyer $u$ sur $f$ et $v, T_1, T_2, \alpha, \beta$ sur les éléments
rendant $f$ proto-inversible.\\
CQFD.\\

Le proto-groupoïde de Segal $\bar{J}^{pre}$ étant défini, on peut
maintenant énoncer l'existence et les propriétés de sa catégorisation
$\bar{J}$.
 
\begin{prop}[-définition]
Il existe un groupoïde de Segal $\bar{J}$ ayant deux objets équivalents notés $0$ et
$1$ vérifiant la propriété suivante :
pour toute catégorie de Segal facile $A$ et pour tout couple $(a,b)$ d'objets de $A$
équivalents, il existe un morphisme de $\bar{J}$ vers $A$ envoyant $0$ et $1$
sur $a$ et $b$.
\end{prop}  
{\it Preuve :}\\
Posons $\bar{J}=Cat(\bar{J}^{pre})$. Comme $\bar{J}^{pre}$ est proto-groupoïdique par définition, sa catégorisation $\bar{J}$ est bien un groupoïde
de Segal, d'après le lemme précédent. En outre, le fait que $u$ soit un
morphisme proto-inversible entre $0$ et $1$ dans $\bar{J}^{pre}$ assure que $u$
est inversible à homotopie près dans $\bar{J}$, donc que $0$ et $1$ sont
équivalents dans $\bar{J}$.\\
\\
Soit $A$ une catégorie de Segal facile et $a$ et $b$ deux de ses objets ayant la
propriété d'être équivalents. Ceci signifie qu'il existe un morphisme
$f$ de $0$ vers $1$ qui est inversible à homotopie près, en particulier $f$
est proto-inversible. Par définition de $\bar{J}^{pre}$, il existe donc un
morphisme de $\bar{J}^{pre}$ vers $A$ envoyant $u$ sur $f$, donc $0$ et $1$ sur
$a$ et $b$. Comme $A$ est une
catégorie facile, ce morphisme se factorise par $Cat(\bar{J}^{pre})=\bar{J}$.\\
CQFD.\\

On a donc bien défini un intervalle $\bar{J}$ qui reconnaît les objets
équivalents d'une catégorie de Segal facile. Toutefois pour que cet
intervalle vérifie la partie c) de l'hypothèse R') de l'ébauche 2 du théorème central,
il va falloir au préalable montrer qu'une \eq~entre groupoïdes de Segal équivaut à
une \eq~d'homotopie faible entre leurs réalisations géométriques.

\newpage

\section{Théorie homotopique des groupoïdes de Segal}

L'idée pour montrer qu'il existe une \eq~de catégories de Segal entre
l'intervalle $\bar{J}$ et sa sous-catégorie de Segal pleine d'objet $0$ est de
travailler au niveau des réalisations géométriques de $\bar{J}$ et de sa
sous-catégorie de Segal pleine d'objet 0. Evidemment cela n'est réalisable
que si l'on montre qu'une \eq~entre groupoïdes de Segal équivaut à
une \eq~d'homotopie faible entre leurs réalisations géométriques. Afin de
montrer ce résultat, nous allons montrer l'\eq~des théories homotopiques des
groupoïdes de Segal et des espaces topologiques en nous inspirant de
l'\eq~de théories homotopiques entre les ensembles simpliciaux et les espaces
topologiques. Mais commençons tout d'abord par définir la réalisation
géométrique des groupoïdes de Segal.

\subsection{Réalisation géométrique}

Comme les précatégories de Segal sont des ensembles bi-simpliciaux, il
existe une réalisation géométrique naturelle les concernant.

\begin{defin}\index{réalisation géométrique!des précatégories de Segal}
Le foncteur réalisation géométrique $|.|$ des précatégories de Segal
est la composée de l'inclusion des précatégories de Segal dans les
ensembles bi-simpliciaux avec le foncteur réalisation géométrique des ensembles
bi-simpliciaux.
\end{defin}

Comme nous voulons étudier de près la réalisation géométrique de
l'intervalle $\bar{J}$, nous allons chercher un espace topologique simple avec
lequel sa réalisation est homotope. Cependant comme $\bar{J}$ est défini
comme un catégorisé, il est vraiment difficile de calculer sa réalisation,
ce qui n'est pas le cas du proto-groupoïde de Segal $\bar{J}^{pre}$ dont il
est le catégorisé. Commençons donc par trouver à quel espace est
homotope $\bar{J}^{pre}$.

\begin{lem}
$|\bar{J}^{pre}|$ est homotope à une sphère.
\end{lem} 
{\it Preuve :}\\
Comme le foncteur $|.|$
préserve les colimites par composée de foncteurs préser\-vant les colimites,
on obtient que $|\bar{J}^{pre}|$ est la colimite des réalisations
géo\-métriques des éléments dont $\bar{J}^{pre}$ est la colimite, à
savoir deux points $0$ et $1$, deux segments $u$ et $v$, deux triangles $T_1$ et
$T_2$ et deux disques $\alpha$ et $\beta$, réalisations géométriques de
$\Delta[1]\Theta \bar{I}$.\\ 
Montrons que $|\bar{J}^{pre}|$ est isomorphe à un cylindre fermé.
Considérons le cylindre ouvert $\mathbb{R}\times [0,1]/(x,y)\sim (x+1,y)$.
Envoyons alors $0$ sur le point $(0,1)$, $1$ sur le point $(1/2,0)$,
$u$ sur le segment $\{t(0,1)+t'(1/2,0), t+t'=1, 0\leq t,t' \leq 1\}$,
$v$ sur le segment $\{t'(1,1)+t(1/2,0), t+t'=1, 0\leq t,t' \leq 1\}$,
$T_1$ sur le triangle $\{t(0,1)+t'(1/2,0)+t''(1,1), t+t'+t''=1, 0\leq t,t',t'' \leq 1\}$
et $T_2$ sur le triangle $\{t(1/2,0)+t'(1,1)+t''(3/2,0), t+t'+t''=1, 0\leq
t,t',t'' \leq 1\}$.
Il est facile de voir d'une part que cela définit bien une application de
$|\bar{J}^{pre}|$ privé des deux disques $\alpha$ et $\beta$ vers le cylindre
ouvert, et d'autre part que cette application est bijective. Le disque $\alpha$
se colle sur $|\bar{J}^{pre}|$ privé des deux disques de la manière suivante
: le bord du disque se colle sur  le cercle réalisation géométrique de
$\delta_{0,2}(T_1)$. Or ce dernier s'envoie dans le cylindre sur le bord
haut du cylindre ($\{(t,1), 0\leq t \leq 1\}$). Ainsi le disque $\alpha$
correspond au couvercle du haut du cylindre. De même, le disque $\beta$
correspond au couvercle du bas du cylindre. Donc $|\bar{J}^{pre}|$ est bien
isomorphe à un cylindre fermé, i.e. avec ses couvercles, et ainsi $|\bar{J}^{pre}|$ est bien homotope
à une sphère.\\
CQFD.\\

Bien entendu ceci ne répond pas à la question du type d'homotopie de la
réalisation de $\bar{J}$. Néanmoins $\bar{J}$ est la catégorisation de
$\bar{J}^{pre}$ et, comme la réalisation géométrique commute aux colimites, la
réalisation géométrique de $\bar{J}$ n'est autre que la "catégorisation"
de la réalisation géométrique de $\bar{J}^{pre}$. Il faut donc montrer que
cette "catégorisation" des réalisations géométriques conserve
l'homotopie. Plus précisément, il s'agit de faire une colimite séquen\-tielle
transfinie de sommes amalgamées à partir des réalisations géométriques
des flèches génératrices des catégories de Segal faciles. Comme la
catégorie des espaces topologiques est une \cmf, il suffit de montrer que les
réalisations géométriques des flèches génératrices des catégories
de Segal faciles sont des cofibrations triviales d'espaces topologiques pour
montrer que la "catégorisation" des réalisations géométriques conserve
l'homotopie. 

\begin{lem}
La réalisation géométrique d'une flèche de $\mathcal{FG}2$ est une
cofibration triviale dans la \cmf~$\mathcal{TOP}$.
\end{lem}
{\it Preuve :}\\
Tout d'abord on remarque que la réalisation géométrique préserve les
cofibrations comme composée de foncteurs les préservant. Ainsi comme les
flèches de $\mathcal{FG}2$ sont des cofibrations, leurs réalisations
géométriques également. Pour montrer que les réalisations
géométriques des flèches de $\mathcal{FG}2$ sont des \eqs~d'homotopie
faibles, on prendra bien garde de ne pas dire que c'est par préservation des
\eqs~faibles par réalisation géométrique, car un tel résultat sur les
précatégories de Segal n'existe pas. En revanche, on calculera
explicitement ce que valent les réalisations des flèches de $\mathcal{FG}2$
et on observera aisément que ce sont des inclusions d'espaces
topologiques dont la source est un rétract par déformation du but. Ces
calculs sont assez semblables, pour ne pas dire identiques, à ceux de Tamsamani
dans \cite{t}. Ainsi les réalisations des flèches de $\mathcal{FG}2$ sont
des \eqs~d'homotopie faibles.\\
CQFD.\\

Grâce à ce résultat, nous pouvons maintenant montrer que les
réalisations géométriques d'une précatégorie de Segal et de sa
catégorisation sont homotopes.

\begin{lem}
Soit $A$ une précatégorie de Segal. Alors $|can_A|:|A|\rightarrow |Cat(A)|$
est une \eq~d'homotopie faible.
\end{lem}
{\it Preuve :}\\
Comme le foncteur réalisation géométrique préserve les colimites, la
réalisation de $can_A$ est la colimite séquentielle transfinie
commençant avec la réalisation géométrique de $A$ de sommes
amalgamées le long d'un morphisme des réalisa\-tions géométriques des
flèches génératrices. Ici les seules flèches génératrices
appartiennent à $\mathcal{FG}2$ car $\mathcal{FG}1$ est vide. Or, d'après le
lemme précédent les réalisations géométriques de ces flèches sont
des cofibrations triviales dans $\mathcal{TOP}$. Comme $\mathcal{TOP}$ est une
\cmf, les cofibrations triviales sont stables par sommes amalgamées le long
d'un morphisme et par colimite séquentielle transfinie. Donc la réalisation
de $can_A$ est une cofibration triviale, donc en particulier une \eq~d'homotopie
faible.\\
CQFD.\\

Comme nous avons déjà montré que la réalisation géométrique de
$\bar{J}^{pre}$ est homotope à une sphère, il découle de la préservation
de l'homotopie des réalisations géométriques par catégorisation que la
réalisation géométrique de l'intervalle $\bar{J}$ est aussi homotope à
une sphère.

\begin{cor}\label{sph}
$|\bar{J}|$ est faiblement homotope à une sphère.
\end{cor}
{\it Preuve :}\\ 
Comme $|\bar{J}^{pre}|$ est homotope à une sphère et que
$\bar{J}$ n'est autre que $Cat(\bar{J}^{pre})$, alors $|\bar{J}|$ est faiblement
homotope à une sphère par le lemme précédent appliqué à
$\bar{J}^{pre}$.\\
CQFD.\\

Avant de se préoccuper de montrer qu'il existe un morphisme entre
l'intervalle $\bar{J}$ et sa sous-catégorie pleine d'objet $0$ dont la
réalisation géométrique est une \eq~d'homotopie, montrons tout d'abord que
pour prouver qu'un morphisme de groupoïdes de Segal est une \eq~de
catégories de Segal il faut et il suffit de montrer que sa réalisation
géométrique est une \eq~d'homotopie. Comme nous avons un foncteur
réalisation géométrique des précatégories de Segal vers les espaces
topologiques, définissons un foncteur singulier des espaces topologiques vers
les précatégories de Segal.

\subsection{Foncteur singulier des espaces topologiques vers les groupoïdes
de Segal}

Pour définir un foncteur singulier qui à tout espace topologique associe une
précatégorie de Segal, comme pour la réalisation géométrique, on va
remarquer que la catégorie des précatégories de Segal est une
sous-catégorie pleine de celle des ensembles bi-simpliciaux.

\begin{defin}\index{SSing(X)}
Soit $X$ un espace topologique. On définit $SSing(X)$ de la manière suivante
: on considère d'abord le foncteur covariant $\Delta[.]:\Delta\rightarrow
\mathcal{ENSSIMP}$ qui à un entier $n$ associe l'ensemble simplicial
$\Delta[n]$. On prend son produit avec lui-même
$\Delta[.]\times\Delta[.]:\Delta\times\Delta\rightarrow
\mathcal{ENSSIMP}\times\mathcal{ENSSIMP}$ et on le compose avec le bi-foncteur
covariant $*\Theta *:\mathcal{ENSSIMP}\times\mathcal{ENSSIMP}\rightarrow
\mathcal{ENSSIMP-PC}$. Puis on compose cela par le foncteur réalisation
géométrique des précatégories de Segal
$|.|:\mathcal{ENSSIMP-PC}\rightarrow \mathcal{TOP}$ qui est covariant. Enfin on compose le tout
avec le foncteur contravariant $Hom_{\mathcal{TOP}}(.,X):\mathcal{TOP}\rightarrow
\mathcal{ENS}$. Le résultat de toutes ces compositions est un foncteur
contravariant de $\Delta\times\Delta$ vers $\mathcal{ENS}$ que l'on notera
$Sing(X)$.\\
On peut résumer par la formule suivante :
$$ SSing(X)=Hom_{\mathcal{TOP}}(|\Delta[.]\Theta\Delta[.]|,X) $$
\end{defin}

Bien évidemment cette construction est bien fonctorielle et à valeurs dans les
précatégories de Segal.

\begin{lem}
La construction $SSing(.)$ est un foncteur covariant de la caté\-gorie des
espaces topologiques vers celle des précatégories de Segal, appelé
foncteur singulier.
\end{lem}
{\it Preuve :}\\
La fonctorialité de $SSing(.)$ vient de la fonctorialité de
$Hom_{\mathcal{TOP}}(|\Delta[.]\Theta\Delta[.]|,X)$ en $X$. Par construction
$SSing(X)$ est un ensemble bi-simplicial. Pour montrer que c'est une
précatégorie de Segal, il faut montrer que $SSing(X)_0$ est un ensemble
discret. Pour cela il suffit de remarquer que par construction les
$\Delta[0]\Theta\Delta[m]$, pour $m\geq 0$, sont tous égaux entre eux et
valent $\Delta[0]$ et que les applications entre eux induites par la 
catégorie simpliciale $\Delta$ valent toute l'identité de $\Delta_0$. Comme
$Hom_{\mathcal{TOP}}(.,X)$ est un foncteur à $X$ fixé, il vient que tous les
$SSing(X)_{0,m}$ sont égaux entre eux et que les morphismes dus à la
structure simpliciale en le second indice sont tous des identités. Ceci montre
que $SSing(X)_0$ est discret et donc que $SSing(X)$ est une précatégorie de
Segal pour tout espace topologique $X$.\\
CQFD.\\

Comme pour les ensembles simpliciaux le foncteur singulier est adjoint au
foncteur réalisation géométrique, on s'attend à ce que ce soit aussi le
cas ici. L'idée pour montrer cette adjonction est d'utiliser le lemme de Yoneda. Donnons
donc une sorte de lemme de Yoneda pour les précatégories de Segal. 

\begin{lem}
Pour toute précatégorie de Segal $A$ et pour tout couple $(m,n)$ d'entiers,
l'ensemble $A_{m,n}$ est naturellement isomorphe à l'ensemble des morphismes
de précatégories de Segal allant de $\Delta[m]\Theta\Delta[n]$ vers $A$.
\end{lem}
{\it Preuve :}\\
Par propriété de la construction $\Theta$, on a un isomorphisme naturel
entre l'ensemble des morphismes dans $\mathcal{C-PC}$ de $\Delta[m]\Theta X$
vers $A$ envoyant le $m+1$-uplet d'objets $(0,\ldots,m)$ de $\Delta[m]\Theta X$
sur le $m+1$-uplet d'objets $(a_0,\ldots,a_m)$ de $A$ et l'ensemble des
morphismes dans $\mathcal{C}$ de $X$ vers $A_m(a_0,\ldots,a_m)$. Ceci est en
particulier vrai pour $\mathcal{C}=\mathcal{ENSSIMP}$. En outre, par Yoneda,
l'ensemble des morphismes dans $\mathcal{ENSSIMP}$ de $\Delta[n]$ vers
$A_m(a_0,\ldots,a_m)$ est naturellement isomorphe à l'ensemble
$(A_m(a_0,\ldots,a_m))_n$ qui n'est autre que $A_{m,n}(a_0,\ldots,a_m)$, par
bi-simplicialité de $A$. Ceci montre que l'ensemble des morphismes dans
$\mathcal{ENSSIMP-PC}$ de $\Delta[m]\Theta \Delta[n]$
vers $A$ envoyant le $m+1$-uplet d'objets $(0,\ldots,m)$ de $\Delta[m]\Theta X$
sur le $m+1$-uplet d'objets $(a_0,\ldots,a_m)$ de $A$ est naturellement
isomorphe à l'ensemble $A_{m,n}(a_0,\ldots,a_m)$. Comme l'ensemble $A_{m,n}$
est le coproduit des $A_{m,n}(a_0,\ldots,a_m)$ avec $(a_0,\ldots,a_m)$
décrivant le produit $m$ fois de l'ensemble $A_0$ avec lui-même, il vient
que l'ensemble des morphismes dans
$\mathcal{ENSSIMP-PC}$ de $\Delta[m]\Theta \Delta[n]$
vers $A$ est naturellement isomorphe à l'ensemble $A_{m,n}$.\\
CQFD.\\

Comme par ce simili-lemme de Yoneda, un élément de $A_{m,n}$ est
représenté par la précatégorie de Segal $\Delta[m]\Theta \Delta[n]$, il
vient que toute précatégorie de Segal n'est autre qu'une colimite de ces
précatégories élémentaires $\Delta[m]\Theta \Delta[n]$.

\begin{cor}
Toute précatégorie de Segal $A$ est naturellement isomorphe à la colimite
des $\Delta[m]\Theta\Delta[n]$, pour $m,n\in\mathbb{N}$, indexé par les
ensembles de morphismes de précatégories de Segal de $\Delta[m]\Theta\Delta[n]$ vers $A$.
\end{cor}
{\it Preuve :} c'est une conséquence immédiate de l'isomorphisme naturel
du lemme précédent. CQFD.\\

De ce résultat découle l'adjonction entre le foncteur singulier et le
foncteur réalisation géométrique des précatégories de Segal.

\begin{cor}\label{adj}
Le foncteur singulier $SSing$ est adjoint à droite du foncteur réalisation
géométrique.
\end{cor}
{\it Preuve :}\\
Soient $A$ une précatégorie de Segal et $X$ un espace topologique. Par le
corollaire précédent, $A$ est naturellement isomorphe à la colimite des
$\Delta[m]\Theta\Delta[n]$. Donc comme $|.|$ est un foncteur préservant les
colimites, $|A|$ est naturellement isomorphe à la colimite des
$|\Delta[m]\Theta\Delta[n]|$. Par propriété universelle de la colimite, il
vient que $Hom_{\mathcal{TOP}}(|A|,X)$ est naturellement isomorphe à la
colimite des $Hom_{\mathcal{TOP}}(|\Delta[m]\Theta\Delta[n]|,X)$. Or cet
ensemble n'est autre que $SSing(X)_{m,n}$, qui par le lemme précédent est
naturellement isomorphe à l'ensemble des morphismes de précatégories de
Segal de $\Delta[m]\Theta\Delta[n]$ dans $SSing(X)$. Ainsi\\
$Hom_{\mathcal{TOP}}(|A|,X)$ est naturellement isomorphe à la colimite des\\
$Hom_{\mathcal{ENSSIMP-PC}}(\Delta[m]\Theta\Delta[n],SSing(X))$. Par le
corollaire précédent, $A$ est isomorphe à la colimite des
$\Delta[m]\Theta\Delta[n]$. Donc par propriété universelle de la colimite,
il vient que la colimite des
$Hom_{\mathcal{ENSSIMP-PC}}(\Delta[m]\Theta\Delta[n],SSing(X))$ est
naturellement isomorphe à $Hom_{\mathcal{ENSSIMP-PC}}(A,SSing(X))$.\\ On a donc
bien que $Hom_{\mathcal{TOP}}(|A|,X)$ est naturellement isomorphe à\\
$Hom_{\mathcal{ENSSIMP-PC}}(A,SSing(X))$.\\
CQFD.\\

Cependant de même que le foncteur singulier des ensembles simpliciaux tombe
dans les ensembles simpliciaux fibrants de Kan, le foncteur singulier $SSing$
tombe dans les groupoïdes de Segal.

\begin{lem}\label{singgp}
Pour tout espace topologique $X$, la précatégorie de Segal\\ $SSing(X)$ est
un groupoïde de Segal.
\end{lem}
{\it Preuve :}\\
On a déjà vu que $SSing(X)$ est une précatégorie de Segal pour tout espace
topologique $X$. Montrons tout d'abord que cette précatégorie de Segal est
une catégorie de Segal facile, i.e. se relève par rapport aux flèches de
la famille $\mathcal{FG}2$, puisque la famille $\mathcal{FG}1$ est vide.
Considérons alors un diagramme du type suivant :
\begin{diagram}
A & \rTo^{g} & SSing(X)\\
\dTo^{f\in\mathcal{FG}_2} & \ruDotsto~{\exists ?}& \\
B & & \\
\end{diagram}
Par adjonction des foncteurs $SSing$ et $|.|$, on obtient alors le diagramme
équivalent suivant :
\begin{diagram}
|A| & \rTo^{|g|} & X\\
\dTo^{|f|} & \ruDotsto~{\exists ?}& \\
|B| & & \\
\end{diagram}
Comme dans la démonstration de la préservation du type d'homotopie de la
réalisation géométrique par catégorisation, en calculant explicitement
les réalisations géométriques des flèches de $\mathcal{FG}2$, on trouve
que ce sont des inclusions munies de rétractions qui en font des
\eqs~d'homotopie faibles. Soit $r:|B|\rightarrow |A|$ la rétraction de $|f|$, on a donc que
$r\circ |f|=Id_{|A|}$. Alors la composée $|g|\circ r$ qui va de $|B|$ vers $X$
vérifie l'égalité suivante :\\ 
$$(|g|\circ r)\circ |f|=|g|\circ (r\circ |f|)=|g|\circ Id_{|A|}=|g|.$$
Donc $|g|\circ r$ est un relèvement du second diagramme. Comme par adjonction
des foncteurs $SSing$ et $|.|$, les deux diagrammes sont équivalents, alors le
premier diagramme se relève, ce qui montre que $SSing(X)$ se relève par
rapport aux flèches de $\mathcal{FG}_2$. Ainsi $SSing(X)$ est une catégorie
de Segal facile.\\
\\
Montrons maintenant que $SSing(X)$ un groupoïde de Segal, i.e. que\\
$\tau_1(SSing(X))$ est un groupoïde. Soit $u$ un morphisme de
$\tau_1(SSing(X))$ et $f$ un élément de $SSing(X)_1$ représentant la
classe $u$. $f$ est donc une application continue de
$|\Delta[1]\Theta\Delta[0]|$, qui n'est autre que le segment $[0,1]$, vers $X$,
donc $f$ est un chemin de $X$. Considérons alors $f^{-1}$ le chemin inverse de
$f$. Notons $v$ sa classe dans $\tau_1(SSing(X))$. Comme les chemins composés
$f^{-1}\circ f$ et $f\circ f^{-1}$ sont homotopes aux chemins constants $Id_{s(f)}$
et $Id_{b(f)}$, dans $\tau_1(SSing(X))$ ceci se traduit par le fait que $v\circ
u$ et $u\circ v$ sont des identités, donc que $u$ est inversible d'inverse
$v$. On a donc montré que $\tau_1(SSing(X))$ est un groupoïde, donc que
$SSing(X)$ est un groupoïde de Segal.\\
CQFD.\\

Ainsi pour tout espace topologique $X$, $SSing(X)$ est un groupoïde de
Segal. De ce fait, le groupoïde de Segal $SSing(X)$ possède des groupes d'homotopie qu'il serait
intéressant de comparer à ceux de l'espace topologique $X$.

\begin{lem}\label{singhom}
Pour tout espace topologique $X$, les groupes d'homotopie du groupoïde de
Segal $SSing(X)$ sont isomorphes naturellement à ceux de $X$.
\end{lem}
{\it Preuve :}\\
Pour commencer, remarquons que la réalisation géométrique de $\Delta[1]\Theta\Delta[0]$ est le
segment $[0,1]$, qui est l'espace topologique caractérisant les chemins. La
réalisation géométrique de $\Delta[1]\Theta\Delta[1]$ est le carré $[0,1]\times [0,1]$
quotienté par la relation suivante : $(0,y)\sim (0,y')$ et $(1,y)\sim (1,y')$.
Or cet espace topologique caractérise les homotopies de chemins à
extrémités fixes. Pour $n>1$, la
réalisation géométrique de $\Delta[1]\Theta\Delta[n]$ est le cylindre de
base un $n$-simplexe $[0,1]\times [0,1]^n$
quotienté par la relation suivante : $(0,Y)\sim (0,Y')$ et $(1,Y)\sim (1,Y')$.
Cet espace topologique caractérise les compositions de $n$ homotopies entre
$n+1$ chemins ayant les mêmes extrémités.\\
\\
De ces constatations découlent les identifications suivantes :
$SSing(X)_{1,0}$ est l'ensemble des chemins de $X$, $SSing(X)_{1,1}$ est
l'ensemble des homotopies de chemins de $X$ à extrémités fixes et, pour
$n>1$, $SSing(X)_{1,n}$ est l'ensemble des compositions de $n$ homotopies de
chemins de $X$ à extrémités fixes. Or ceci n'est autre que l'ensemble
simplicial $Sing(X^{[0,1]})$, image de l'espace des chemins de $X$ par le
foncteur singulier des ensembles simpliciaux. De même, on obtient que
l'ensemble simplicial $SSing(X)_1(x,x)$ n'est autre que l'ensemble simplicial
$Sing(\Omega(X,x))$, image de l'espace des lacets de $X$ en $x$ par le
foncteur singulier des ensembles simpliciaux.\\
\\
Comme on a vu au lemme précédent que $SSing(X)$ est un groupoïde de
Segal, on peut donc parler de ses groupes d'homotopie. Soit $n>1$, on a par
définition des groupes d'homotopie d'un groupoïde de Segal que
$$\pi_n(SSing(X),x)=\pi_{n-1}(|SSing(X)_1(x,x)|,Id_x).$$ Or on a vu ci-dessus que
$SSing(X)_1(x,x)$ n'est autre que $Sing(\Omega(X,x))$. Comme en outre le couple
de foncteurs $(|.|,Sing)$ forme une paire d'\eqs~de Quillen entre la catégorie des
ensembles simpliciaux et celle des espaces topologiques, la réalisation
géométrique de $Sing(\Omega(X,x))$, donc de $SSing(X)_1(x,x)$, est
naturellement isomorphe à $\Omega(X,x)$. Ainsi\\
$\pi_{n-1}(|SSing(X)_1(x,x)|,Id_x)$ est naturellement isomorphe au groupe\\
$\pi_{n-1}(\Omega(X,x),Id_x)$ qui n'est autre que $\pi_n(X,x)$, par
propriété de $\Omega(X,x)$. On a donc montré que pour tout entier $n>1$,
les groupes d'homotopie $\pi_n(SSing(X),x)$ et $\pi_n(X,x)$ sont naturellement
isomorphes. Et on montre de manière similaire que le groupe d'homotopie
$\pi_1(SSing(X),x)$, qui par définition n'est autre que
$\pi_0(|SSing(X)(1,1)|)$, est naturellement en bijection avec
$\pi_0(\Omega(X,x))$, qui n'est autre que $\pi_1(X,x)$. Ceci au passage munit
$\pi_1(SSing(X),x)$ d'une structure de groupe.\\
\\
Comme on l'a vu précédemment, $SSing(X)_1$ est en fait $Sing(X^{[0,1]})$ et
donc sa réalisation est naturellement isomorphe à $X^{[0,1]}$, car
$(|.|,Sing)$ est une paire d'\eqs~de Quillen. Ainsi $\tau_1(SSing(X))_1$, qui par
définition vaut $\pi_0(|SSing(X)_1|)$, est naturellement en bijection avec
$\pi_0(X^{[0,1]})$. Comme tout élément de $\tau_1(SSing(X))_1$ est une
classe représenté par un chemin de $X$ et que toute classe est inversible
car $SSing(X)$ est un groupoïde de Segal, tout chemin représente dans
$\tau_1(SSing(X))$ une \eq~entre ses extrémités. Ceci montre que si deux
points de $X$ sont dans la même composante connexe, comme alors il existe un
chemin les liant et que ce chemin représente une \eq~entre eux, par ce qui
précède, alors ces deux points ont même classe dans
$\tau_0(SSing(X))=\pi_0(SSing(X))$. Réciproquement, si deux points de $SSing(X)$ sont
équivalents, il existe un isomorphisme dans $\tau_1(SSing(X))$ qui est une
classe dont tous les représentants sont des chemins liant ces deux points qui
sont, par conséquent, dans la même composante connexe. Ainsi on a montré que
le morphisme identique entre les objets de $SSing(X)$ et les points de $X$
induit une bijection entre les ensembles quotients $\pi_0(SSing(X))$ et $\pi_0(X)$.\\
CQFD.\\

Ce résultat d'\eq~de notions homotopiques entre $X$ et $SSing(X)$ entraîne
le fait que toute \eq~d'homotopie faible entre espaces topologiques donne lieu
par le foncteur singulier $SSing$ à une \eq~de catégories de Segal entre
groupoïdes de Segal.

\begin{cor}
Soit $f:X\rightarrow Y$ une application continue d'espaces topologiques
simplement connexes.
Si $f$ est une \eq~d'homotopie faible alors $SSing(f)$ est une \eq~de
catégories de Segal.
\end{cor}
{\it Preuve :}\\
Par le lemme précédent, les groupes d'homotopie des groupoïdes de Segal
$SSing(X)$ et $SSing(Y)$ sont naturellement isomorphes à ceux de $X$ et $Y$.
Comme $f$ est une \eq~d'homotopie faible, i.e. induit des isomorphismes sur les
groupes d'homotopie de $X$ et $Y$, ce résultat entraîne que $SSing(f)$ induit des
isomorphismes sur les groupes d'homotopie des groupoïdes de Segal $SSing(X)$
et $SSing(Y)$. Toutefois comme on l'a vu au lemme~\ref{eqgp1}, cela ne suffit
pas à montrer que $SSing(f)$ est une \eq~de catégories de Segal, sauf si
$SSing(X)$ et $SSing(Y)$ sont simplement connexes. Or ceci est vrai car $X$ et
$Y$ sont simplement connexes par hypothèse et que ces derniers ont mêmes groupes
d'homotopie que $SSing(X)$ et $SSing(Y)$.\\
CQFD.\\

Ici l'hypothèse de simple connexité est nécessaire car, dans ce cas, un
morphisme de groupoïde de Segal induisant des isomorphismes sur les groupes
d'homotopie est une \eq~de catégories de Segal.\\

Grâce à ce résultat, nous avons fait un premier pas dans la
démonstration de l'\eq~entre \eqs~de groupoïdes de Segal et
\eqs~d'homotopie de leurs réalisations géométriques. Afin d'aller plus
avant sur cette voie, nous avons besoin d'introduire la notion de catégories
topologiques, car Segal a montré sur les groupoïdes topologiques un
résultat fondamental pour ces histoires d'\eq~de théories homotopiques (voir
\cite{t}).

\subsection{Groupoïde topologique}

Pour définir la notion de catégorie topologique, nous devons tout d'abord nous donner une
donnée de Segal sur les espaces topologiques. Comme nous voulons que les
catégories topologiques ne soient pas trop éloignées des caté\-gories de
Segal, nous allons nous inspirer de la donnée de Segal définissant les
catégories de Segal pour définir les éléments constitutifs de la
donnée de Segal qui va définir les catégories topologiques.

\begin{prop}[-définition]\index{précatégorie topologique}\index{catégorie!topologique}
La catégorie $\mathcal{TOP}$ munie des notions suivantes forme une donnée de
Segal :
\item -un espace topologique régal est un espace topologique quelconque,
\item -une \eqc~d'espaces topologiques régaux est une \eq~d'homotopie faible,
\item -le foncteur de troncature $\tau_0$ des espaces topologiques régaux est
le foncteur composante connexe $\pi_0$,\\
\\
Pour cette donnée de Segal, les \precats~seront appelées les
précatégories topologiques, les \cats~les catégories topologiques et les
\eqs~de \cats~les \eqs~de catégories topologiques.
\end{prop}
{\it Preuve :} les propriétés de donnée de Segal se vérifient très
facilement.\\

De même que nous avons défini la notion de groupoïde de Segal pour les
catégories de Segal, nous allons pareillement définir une notion de
groupoïde topologique pour les catégories topologiques.

\begin{defin}\index{groupoïde!topologique}
Une catégorie topologique $A$ est un groupoïde topologique si $\tau_1(A)$ est le
nerf d'un groupoïde au sens des catégories. Les morphismes de groupoïdes topologiques sont les morphismes de précatégories topologiques et les
\eqs~de groupoïdes topologiques sont les équivalences de catégories
topologiques.
\end{defin}

Comme le foncteur réalisation géométrique des ensembles simpliciaux
permet de passer de la catégorie $\mathcal{ENSSIMP}$ à la catégorie
$\mathcal{TOP}$ tout en pré\-servant la notion de donnée de Segal, nous
allons voir comment le foncteur induit des précatégories de Segal aux
précatégories topologiques se comporte vis-à-vis des notions de
catégorie, de groupoïde et d'\eq.

\begin{lem}\label{segtop}
Notons $|.|_{simp}$ le foncteur réalisation géométrique des ensembles
simpliciaux. Ce foncteur induit un foncteur, que l'on notera identiquement, de
la catégorie des précatégories de Segal vers la catégorie des
précatégories topologiques préservant les notions de catégorie, d'\eq~de
catégories, de groupoïde et de $\tau_0$. Plus précisément, on a :
\item 1) $A$ est une catégorie de Segal si et seulement si $|A|_{simp}$ est une catégorie
topologique,
\item 2) pour toute catégorie de Segal $A$, on a que
$\tau_1(A)=\tau_1(|A|_{simp})$,
en particulier $\tau_0(A)=\tau_0(|A|_{simp})$,
\item 3) soit $f$ est un morphisme de catégories de Segal, $f$ est
une \eq~de catégories de Segal si et seulement si $|f|_{simp}$ est une
\eq~de catégories topologiques,
\item 4) soit $A$ une catégorie de Segal, $A$ est un groupoïde de Segal si
et seulement si $|A|_{simp}$ est un groupoïde topologique.
\end{lem}
{\it Preuve :}\\
1) Soit $A$ une précatégorie de Segal. La précatégorie topologique $|A|_{simp}$
n'est autre que la composée du préfaisceau $A$ avec le foncteur $|.|_{simp}$.
Comme le foncteur réalisation géométrique des ensembles simpliciaux
préserve les produits fibrés, les morphismes de Segal de $|A|_{simp}$ ne
sont autres que les réalisations géométriques de ceux de $A$. Or ces
derniers sont des \eqs~faibles d'ensembles simpliciaux par définition si et
seulement si leurs
réalisations géométriques sont des \eqs~d'homotopie faibles. Ainsi les
morphismes de Segal de $A$ sont des \eqs~faibles d'ensembles simpliciaux si et
seulement si les
morphismes de Segal de $|A|_{simp}$ sont des \eqs~d'homotopie faibles, ce qui
signifie que $A$ est une catégorie de Segal si et seulement si $|A|_{simp}$
est une catégorie topologique. Ceci montre la première partie
du lemme.\\
\\
2) Par définition, la catégorie $\tau_1(A)$ est obtenue par composition du
préfaisceau $A$ par le foncteur $\tau_0$ des ensembles simpliciaux. Or ce
dernier n'est autre que la composée du foncteur réalisation géométrique
des ensembles simpliciaux avec le foncteur composante connexe des espaces
topologiques. On a donc l'égalité ci-dessous :
$$\tau_1(A)=\pi_0\circ |.|_{simp} \circ A=\pi_0\circ |A|_{simp}=\pi_0(|A|_{simp}).$$ 
D'un autre côté, pour une catégorie topologique $B$, la catégorie
$\tau_1(B)$ est obtenue par composition du préfaisceau $B$ par le foncteur
$\tau_0$ des espaces topologiques qui n'est autre que le foncteur composante
connexe $\pi_0$. Par ce qui précède, comme $A$ est une catégorie
de Segal, $|A|_{simp}$ est une catégorie topologique, donc la catégorie
$\tau_1(|A|_{simp})$ est la composée du préfaisceau $|A|_{simp}$ par le
foncteur $\pi_0$. On a donc l'égalité suivante :
$$\tau_1(|A|_{simp})=\pi_0 \circ |A|_{simp}=\pi_0(|A|_{simp}).$$
En comparant les égalités, il vient que $\tau_1(A)$ et $\tau_1(|A|_{simp})$
sont égaux par définition même des $\tau_1$. Comme en outre $\tau_0(A)$
est l'ensemble des classes d'isomorphismes d'objets de la catégorie $\tau_1(A)$ et que
$\tau_0(|A|_{simp})$ est l'ensemble des classes d'isomorphismes d'objets de la
catégorie $\tau_1(|A|_{simp})$, l'égalité des deux catégories $\tau_1(A)$ et $\tau_1(|A|_{simp})$
entraîne celle des ensembles $\tau_0(A)$ et $\tau_0(|A|_{simp})$, ce qui
montre la deuxième partie du lemme.\\
\\
4) En outre dire qu'une catégorie de Segal $A$ est un groupoïde de Segal
c'est dire que $\tau_1(A)$ est le nerf d'un groupoïde, et dire que la
catégorie topologique $|A|_{simp}$ est un groupoïde topologique c'est dire
que $\tau_1(|A|_{simp})$ est un nerf de groupoïde. Or l'égalité des
catégories $\tau_1(A)$ et $\tau_1(|A|_{simp})$ entraîne que $A$ est un groupoïde de Segal si et seulement si
$|A|_{simp}$ est un groupoïde topologique, ce qui montre la quatrième
partie du lemme.\\
\\
3) Soit $f$ un morphisme de catégories de Segal. L'égalité entre les $\tau_0$
des catégories de Segal et de leurs catégories topologiques réalisées
entraîne l'égalité entre $\tau_0(f)$ et $\tau_0(|f|_{simp})$. Ainsi on
obtient que $f$ est essentiellement surjective si et seulement si $|f|_{simp}$
l'est. En outre comme la réalisation géométrique des ensembles simpliciaux
préserve les produits fibrés, il vient que, pour tout couple d'objets $(x,y)$
de la source de $f$, $(|f|_{simp})_1(x,y)$ n'est autre
que la réalisation géométrique de $f_1(x,y)$. Ainsi on obtient que
$f_1(x,y)$ est une \eq~faible d'ensembles simpliciaux si et seulement si
$(|f|_{simp})_1(x,y)$ est une \eq~d'homotopie faible, ce qui signifie que $f$
est pleinement fidèle si et seulement si $|f|_{simp}$ l'est. On a donc bien
montré que $f$ est une \eq~de catégories de Segal si et seulement si
$|f|_{simp}$ est une \eq~de catégories topologiques, ce qui montre la
troisième partie du lemme.\\
CQFD.\\

Nous venons donc de montrer que la notion d'\eq~de catégories de Segal est
équivalente à celle d'\eq~de catégories topologiques. Ainsi le problème
de l'\eq~entre \eqs~de groupoïdes de Segal et \eqs~d'homotopie faibles de
leurs réalisations géométriques devient un problème d'\eq~entre \eqs~de
groupoïdes topologiques et \eqs~d'homotopie faibles de
leurs réalisations géométriques. Nous allons donc maintenant montrer que
la réalisation géométrique d'une \eq~entre groupoïdes topologiques est
une \eq~d'homotopie faible. Pour cela, définissons d'abord une
réalisation géométrique des précatégories topologiques.

\begin{defin}\index{réalisation géométrique!des précatégories
topologiques}
Comme les précatégories topologiques peuvent être vues comme des ensembles simpliciaux dont les
niveaux admettent une structure topologique (compatible avec la structure
simpliciale), on définit le foncteur $|.|_{top}$ de réalisation
géométrique des précatégories topologiques de la manière suivante.\\
Soit $A$ une précatégorie topologique, la réalisation géométrique de
l'ensemble simplicial sous-jacent à $A$ est de la forme suivante :
$$|A|_{simp}=\Big(\coprod_{n\geq 0}\Delta[n]\times A_n\Big)/\sim $$
où la relation d'\eq~est induite de la structure simpliciale et où les $A_n$
sont considérés comme des ensembles, i.e. des espaces topologiques discrets.
Alors on pose : 
$$|A|_{top}=\Big(\coprod_{n\geq 0}\Delta[n]\times A_n\Big)/\sim $$
La formule est la même que pour $|A|_{simp}$ mais cette fois les $A_n$ gardent
leurs topologies propres. Ces topologies étant compatibles avec la structure
simpliciale qui engendre la relation d'\eq~$\sim$, la topologie quotient de
$|A|_{top}$ a bien un sens.
\end{defin}

La réalisation géométrique des précatégories topologiques définie,
nous pouvons montrer que les réalisations géométriques d'\eqs~de
catégories topologiques entre groupoïdes topologiques sont des
\eqs~d'homotopie faibles.

\begin{prop}
Le foncteur $|.|_{top}$ vérifie les deux propriétés suivantes :
\item - si $f$ est une \eq~de catégories topologiques entre groupoïdes
topologiques alors $|f|_{top}$ est une \eq~d'homotopie faible,
\item - pour tout groupoïde topologique $\Phi$, on a un isomorphisme naturel
entre $\tau_0(\Phi)$ et $\pi_0(|\Phi|_{top})$.
\end{prop}
{\it Preuve :}\\
Pour démontrer cette proposition, nous nous servirons du théorème de
Segal cité dans \cite{t} qui s'énonce ainsi :\\
\\
\emph{Si $\Phi$ est un groupoïde topologique, alors :\\
\\
- pour tout objet $x,y$ de $\Phi$, le morphisme naturel $\Phi_1(x,y)\rightarrow
Chemin_{x,y}(|\Phi|_{top})$ est une \eq~d'homotopie faible,\\
\\
- l'injection naturelle $\Phi_0\rightarrow |\Phi|_{top}$ induit une bijection de
$\tau_0(\Phi)$ avec $\tau_0(|\Phi|_{top})$, qui n'est autre que
$\pi_0(|\Phi|_{top})$.}\\
\\
La deuxième partie du théorème de Segal nous donne directement la seconde
partie de la proposition. Il ne reste donc plus qu'à montrer la première
partie de la proposition.\\
\\
Soit $f:\Phi\rightarrow \Psi$ une \eq~de catégories topologiques entre
groupoïdes topologiques. Considérons le diagramme commutatif suivant :
\begin{diagram}
\tau_0(\Phi) & \rTo^{\tau_0(f)} & \tau_0(\Psi) \\
\dTo^{\cong} & & \dTo_{\cong} \\
\pi_0(|\Phi|_{top}) & \rTo_{\pi_0(|f|_{top})} & \pi_0(|\Psi|_{top}) \\
\end{diagram}
Comme $f$ est une \eq~de catégories topologiques, $f$ est essentiellement
surjective, i.e. $\tau_0(f)$ est une bijection ensembliste. Comme, en outre,  $\Phi$ et
$\Psi$ sont des groupoïdes topologiques, on peut leur appliquer le
théorème de Segal dont la seconde partie nous donne la bijectivité des
deux flèches verticales. D'où $|f|_{top}$ induit une bijection sur les
$\pi_0$.\\
\\
Considérons maintenant le diagramme commutatif suivant :
\begin{diagram}
\Phi_1(x,y) & \rTo^{f_1(x,y)} & \Psi_1(f(x),f(y)) \\
\dTo^{\sim} & & \dTo_{\sim} \\
Chemin_{x,y}(|\Phi|_{top}) & \rTo & Chemin_{f(x),f(y)}(|\Psi|_{top}) \\
\end{diagram} 
Comme $f$ est une \eq~de catégories topologiques, $f$ est pleinement fidèle,
i.e. $f_1(x,y)$ est une \eq~d'homotopie faible. En outre comme $\Phi$ et $\Psi$
sont des groupoïdes topologiques, on peut leur appliquer le théorème de
Segal dont la première partie nous donne le fait que les deux flèches
verticales sont des \eqs~d'homotopie faibles. Comme les \eqs~d'homotopie faibles
vérifient l'axiome "trois pour deux", le morphisme entre espaces de chemins
est aussi une \eq~d'homotopie faible. En particulier si $x=y$, on obtient une
\eq~d'homotopie faible entre les espaces de lacets $\Omega(|\Phi|_{top},x)$ et
$\Omega(|\Psi|_{top},f(x))$, i.e. $|f|_{top}$ induit des isomorphismes des
$\pi_j(\Omega(|\Phi|_{top},x), Id_x)$
vers les $\pi_j(\Omega(|\Psi|_{top},f(x)), Id_{f(x)})$ pour $j\geq 0$. Par
ailleurs par propriétés des espaces de lacets, on a les égalités
suivantes :
$$\pi_j(\Omega(|\Phi|_{top},x), Id_x)=\pi_{j+1}(|\Phi|_{top},x),$$  
$$\pi_j(\Omega(|\Psi|_{top},f(x)), Id_{f(x)})=\pi_{j+1}(|\Psi|_{top},f(x)).$$
Ainsi on a obtenu que $|f|_{top}$ induit des isomorphismes de groupes sur les $\pi_i$ pour
$i>0$. Comme on a vu précédemment que $|f|_{top}$ induit une bijection sur
les $\pi_0$, on a bien que $|f|_{top}$ est une \eq~d'homotopie faible, ce qui
montre la première partie de la proposition.\\
CQFD.\\

Comme les notions d'\eq~de catégories topologiques et d'\eq~de catégories
de Segal sont équivalentes, on obtient donc que tout morphisme de groupoïdes de Segal qui est une \eq~de catégories de Segal se réalise
géométriquement en une \eq~d'homotopie faible.

\begin{cor}\label{eqgpreal}
Soit $f:A\rightarrow B$ un morphisme de groupoïdes de Segal.
Si $f$ est une \eq~de catégories de Segal alors sa réalisation
géométrique $|f|$ est une \eq~d'homotopie faible.
\end{cor}
{\it Preuve :}\\
Soit $f$ une \eq~de catégories de Segal entre groupoïdes de Segal, alors
par le lemme~\ref{segtop}, le morphisme de précatégories topologiques
$|f|_{simp}$
est une \eq~de catégories topologiques entre groupoïdes topologiques. Par
la proposition précédente, il vient que $||f|_{simp}|_{top}$ est une
\eq~d'homotopie faible. Or le foncteur réalisation géométrique des
ensembles bi-simpliciaux et le foncteur composé $|.|_{top}\circ |.|_{simp}$ sont
naturellement isomorphes. Donc la réalisation géométrique $|f|$ de $f$ est
bien une \eq~d'homotopie faible.\\
CQFD.\\ 
 
On obtient également comme corollaire de la proposition précédente un
résultat très utile pour la suite qui compare la troncature $\tau_0$ d'un
groupoïde de Segal avec le groupe d'homotopie $\pi_0$ de sa réalisation. 

\begin{cor}
Pour tout groupoïde de Segal $A$, on a un isomorphisme naturel entre
$\tau_0(A)$ et $\pi_0(|A|)$.
\end{cor}
{\it Preuve :}\\
Soit $A$ un groupoïde de Segal. Par le lemme~\ref{segtop}, $|A|_{simp}$ est
un groupoïde topologique et le foncteur $|.|_{simp}$ préserve les $\tau_0$.
Ainsi $\tau_0(A)$ est égal à $\tau_0(|A|_{simp})$ qui par la proposition
précédente est naturellement isomorphe à $\pi_0(||A|_{simp}|_{top})$.
Comme $|A|$ est homéomorphe à $||A|_{simp}|_{top}$, on a bien que
$\tau_0(A)$ est naturellement isomorphe à $\pi_0(|A|)$.\\
CQFD.\\

Nous avons désormais en main tous les atouts nous permettant de montrer qu'un
morphisme de groupoïdes de Segal est une \eq~de catégories de Segal si et
seulement si sa réalisation géométrique est une \eq~d'homotopie faible. La
démonstration de ce résultat fait l'objet de la partie suivante.

\subsection{Caractérisation par réalisation géométrique des
équivalences de groupoïdes de Segal}

Afin de montrer l'\eq~entre \eqs~de groupoïdes de Segal et \eqs~d'homotopie
faibles de leurs réalisations géométriques, nous allons démontrer toute
une série de lemmes sur la préservation de la notion d'homotopie par
réalisation géométrique comme par foncteur singulier. Débutons par le
fait que tout groupoïde de Segal est équivalent au groupoïde de Segal
résultant de l'application du foncteur singulier à sa réalisation
géométrique.

\begin{lem}\label{adjunit}
Pour tout groupoïde de Segal $A$, le morphisme naturel $A\rightarrow
SSing(|A|)$, donné par l'unité de l'adjonction des foncteurs $SSing$ et
$|.|$, est une \eq~de catégories de Segal
\end{lem}
{\it Preuve :}\\
Notons $\eta$ le morphisme naturel de $A$ vers $SSing(|A|)$. Montrons que c'est
une \eq~de catégories de Segal. Calculons tout d'abord $\tau_0(SSing(|A|))$.
Comme, d'après le lemme~\ref{singgp}, $SSing(|A|)$ est un groupoïde de Segal, par définition de son
$\pi_0$, on a que $\pi_0(SSing(|A|))$ n'est autre que $\tau_0(SSing(|A|))$. En
outre, on a par le lemme~\ref{singhom} que les groupes d'homotopie de $SSing(|A|)$ sont naturellement
isomorphes à ceux de $|A|$. Ainsi $\pi_0(SSing(|A|))$ est naturellement
isomorphe à $\pi_0(|A|)$ et donc $\tau_0(\eta)$ n'est autre que le morphisme
naturel de $\tau_0(A)$ dans $\pi_0(|A|)$ qui est une bijection par le corollaire
précédent appliqué au groupoïde de Segal $A$. Ceci montre qu' $\eta$
est essentiellement surjective.\\ 
\\
Montrons qu'$\eta$ est aussi pleinement fidèle,
i.e. que, pour tout couple $(x,y)$ d'objets de $A$, $\eta_1(x,y)$ est une
\eq~faible d'ensembles simpliciaux, donc que $|\eta_1(x,y)|_{simp}$ est une
\eq~d'homotopie faible. D'après la démonstration du lemme~\ref{singhom}, on
a que $SSing(|A|)_1(x,y)$ n'est autre que l'ensemble simplicial
$Sing(Chemin_{x,y}(|A|))$. Ainsi $|\eta_1(x,y)|_{simp}$ est un morphisme de
$|A_1(x,y)|_{simp}$ dans $|Sing(Chemin_{x,y}(|A|))|_{simp}$. Or les foncteurs
$Sing$ et $|.|_{simp}$ réalisent une \eq~de Quillen entre les ensembles
simpliciaux et les espaces topologiques. En particulier, le morphisme naturel de
$|Sing(Chemin_{x,y}(|A|))|_{simp}$ dans $Chemin_{x,y}(|A|)$ induit par
l'adjonction est une \eq~d'homotopie faible. Précomposée par
$|\eta_1(x,y)|_{simp}$ on obtient le morphisme naturel de $|A_1(x,y)|_{simp}$ dans
$Chemin_{x,y}(|A|)$. Comme la réalisation géométrique des ensembles
simpliciaux commute aux produits fibrés, on a que $|A_1(x,y)|_{simp}$ n'est
autre que $(|A|_{simp})_1(x,y)$ et comme $|A|$ et $||A|_{simp}|_{top}$ sont
homéomorphes, on obtient que la composée de $|\eta_1(x,y)|_{simp}$ par le
morphisme d'adjonction n'est autre que le morphisme naturel de\\
$(|A|_{simp})_1(x,y)$ vers $Chemin_{x,y}(||A|_{simp}|_{top})$. Or $A$ est un
groupoïde de Segal donc par le lemme~\ref{segtop}, $|A|_{simp}$ est un
groupoïde topologique. On peut donc lui appliquer le théorème de Segal
que l'on a rappelé au début de la démonstration de la proposition précédente, ce qui nous donne
que le morphisme naturel de $(|A|_{simp})_1(x,y)$ vers
$Chemin_{x,y}(||A|_{simp}|_{top})$ est une \eq~d'homotopie faible. En outre ce
morphisme n'est autre que la composée de $|\eta_1(x,y)|_{simp}$ par le
morphisme d'adjonction qui est une \eq~d'homotopie faible, donc, comme les
\eqs~d'homotopie faibles vérifient l'axiome "trois pour deux",
$|\eta_1(x,y)|_{simp}$ est aussi une \eq~d'homotopie faible. Ainsi $\eta_1(x,y)$
est bien une \eq~faible d'ensembles simpliciaux, ce qui montre la pleine
fidélité d'$\eta$. Comme on a déjà vu qu'$\eta$ est essentiellement
surjective, on a bien montré qu'$\eta$ est une \eq~de catégories de Segal.\\
CQFD.\\

Ce résultat a un corollaire très intéressant qui montre que les groupes
d'homotopie des groupoïdes de Segal sont isomorphes à ceux de leurs
réalisations géométriques.

\begin{cor}\label{gphom}
Les groupes d'homotopie d'un groupoïde de Segal sont naturellement
isomorphes à ceux de sa réalisation géométrique.
\end{cor}
{\it Preuve :}\\
Par le lemme précédent, pour tout groupoïde de Segal $A$, le morphisme
naturel d'adjonction $A\rightarrow SSing(|A|)$ est une \eq~de catégories de
Segal. Or $A$ est un groupoïde de Segal par hypothèse et $SSing(|A|)$ est un
groupoïde de Segal par le lemme~\ref{singgp}, donc $A\rightarrow SSing(|A|)$
est une \eq~de catégories de Segal entre groupoïdes de Segal. Par la
remarque qui suit le lemme~\ref{eqgp1}, $A\rightarrow SSing(|A|)$ induit des
isomorphismes au niveau des groupes d'homotopie des groupoïdes de Segal. En
outre, d'après le lemme~\ref{singhom}, les groupes d'homotopie de $SSing(|A|)$
sont naturellement isomorphes à ceux de $|A|$. Ainsi les groupes d'homotopie du
groupoïde de Segal $A$ sont isomorphes à ceux de $SSing(|A|)$, car
l'unité de l'adjonction est une \eq~de catégorie de Segal entre groupoïdes de Segal, et les groupes d'homotopie de $SSing(|A|)$ sont naturellement
isomorphes à ceux de $|A|$ par le lemme~\ref{singhom}. Par composition des
isomorphismes naturels, on obtient bien que les groupes d'homotopie du
groupoïde de Segal $A$ sont naturellement isomorphes à ceux de $|A|$.\\
CQFD.\\

On remarque au passage que la bijection naturelle des $\pi_1$ des groupoïdes de Segal avec
ceux de leurs réalisations géométriques munissent les $\pi_1$ des
groupoïdes de Segal d'une structure de groupe.\\

On en arrive alors au résultat tant attendu de l'\eq~entre les notions
d'\eq~de groupoïdes de Segal et d'\eq~d'homotopie faible de leurs
réalisations géométriques.

\begin{cor}\label{eqgp2}
Si $f:A\rightarrow B$ est un morphisme de groupoïdes de Segal simplement
connexes, alors les énoncés suivants sont équivalents :
\item i) $f:A\rightarrow B$ est une \eq~de catégories de Segal,
\item ii) $|f|:|A|\rightarrow |B|$ est une \eq~d'homotopie faible. 
\end{cor}
{\it Preuve :}\\
On a déjà montré au corollaire~\ref{eqgpreal} que i) implique ii). Montrons
donc la réciproque. Soit $f:A\rightarrow B$ un morphisme de groupoïdes de
Segal simplement connexes dont la réalisation géométrique est une \eq~d'homotopie faible.
Montrons qu'alors $f$ est une \eq~de catégories de Segal. Comme par
hypothèse, $|f|$ est une \eq~d'homotopie faible, elle induit des isomorphismes
entre les groupes d'homotopie de $|A|$ et de $|B|$. Or l'on a, d'après le corollaire
précédent, que, comme $A$ et $B$ sont des groupoïdes de Segal, leurs
groupes d'homotopie sont naturellement isomorphes à ceux de leurs
réalisations géométriques. Ainsi $f$ induit des isomorphismes sur les
groupes d'homotopie des groupoïdes de Segal $A$ et $B$. Or par hypothèse
$A$ et $B$ sont simplement connexes donc on peut appliquer le lemme~\ref{eqgp1}
à $f$, ce qui nous donne que $f$ est bien une \eq~de catégories de Segal. On
a donc montré que ii) implique i).\\
CQFD.\\

Encore une fois on constate que l'hypothèse de simple connexité intervient dans la
démonstration de cette équivalence. Avant de clore cette partie, montrons le
résultat dual à celui du lemme~\ref{adjunit}.

\begin{cor}
Pour tout espace topologique $X$, le morphisme naturel\\ $|SSing(X)|\rightarrow
X$, donné par la co-unité de l'adjonction des foncteurs $SSing$ et $|.|$,
est une \eq~d'homotopie faible.
\end{cor}
{\it Preuve :}\\
Notons $\epsilon$ la co-unité de l'adjonction. Montrons que c'est une
\eq~d'homotopie faible. Considérons $SSing(\epsilon)$, c'est un morphisme de
groupoïdes de Segal de $SSing(|SSing(X)|)$ vers $SSing(X)$ d'après le
lemme~\ref{singgp}. Par ailleurs l'unité de l'adjonction
$\eta:SSing(X)\rightarrow SSing(|SSing(X)|)$ appliquée au
groupoïde de Segal $SSing(X)$ est une \eq~de catégories de Segal par le
lemme~\ref{adjunit}. Or la composée $SSing(\epsilon)\circ \eta$ n'est autre que
l'identité de $SSing(X)$, car les foncteurs $SSing$ et $|.|$ sont adjoints par
le corollaire~\ref{adj}. Comme les hypothèses de A') à Q') de l'ébauche 2 du théorème
central sont vérifiées par la donnée de Segal définissant les
catégories de Segal, l'axiome "trois pour deux" est vrai pour les \eqs~de
catégories de Segal, ce qui entraîne que $SSing(\epsilon)$ est une \eq~de
catégories de Segal. Comme en outre c'est une \eq~de catégories de Segal
entre groupoïdes de Segal, elle induit des isomorphismes entre les groupes
d'homotopie de $SSing(|SSing(X)|)$ et de $SSing(X)$ d'après la remarque
suivant le lemme~\ref{eqgp1}. Or par le
lemme~\ref{singhom}, les groupes d'homotopie de $SSing(|SSing(X)|)$ et de $SSing(X)$ sont naturellement
isomorphes à ceux de $|SSing(X)|$ et de $X$. Ainsi $\epsilon$ induit des
isomorphismes entre les groupes d'homotopie de $|SSing(X)|$ et de $X$, ce qui fait de
$\epsilon$ une \eq~d'homotopie faible.\\
CQFD.\\  

Nous savons maintenant que, pour montrer qu'un morphisme de groupoïdes de
Segal est une \eq~de catégories de Segal, il faut et il suffit de montrer que
sa réalisation géométrique est une \eq~d'homotopie faible. Nous pouvons
donc utiliser ce résultat pour montrer que l'intervalle $\bar{J}$ est équivalent
à sa sous-catégorie de Segal pleine d'objet $0$, ce qui terminera la
vérification des hypothèses du théorème central et ainsi montrera que
les précatégories de Segal avec leurs \eqs~faibles forment une \cmf.

\newpage

\section{Catégorie de modèles fermée sur les pré\-catégories de Segal}

Afin de pouvoir appliquer l'\eq~entre \eqs~de groupoïdes de Segal et
\eqs~d'homotopie faibles de leurs réalisations géométriques à
l'intervalle $\bar{J}$, il faut tout d'abord s'assurer que ce dernier est bien simplement
connexe.

\begin{lem}
Le groupoïde de Segal $\bar{J}$ est simplement connexe.
\end{lem}
{\it Preuve :}\\
Par définition, $\bar{J}$ est un groupoïde de Segal.
D'après le corollaire~\ref{gphom}, les groupes d'homotopie des groupoïdes
de Segal sont naturellement isomorphes à ceux de leurs réalisations
géométriques. Or la réalisation géométrique de $\bar{J}$ est
faiblement homotope à une sphère par le corollaire~\ref{sph}. Comme la
sphère est un espace topologique simplement connexe, alors la réalisation
géométrique de $\bar{J}$ aussi. Donc le $\pi_0$ et les $\pi_1$ de la
réalisation géométrique de $\bar{J}$ sont triviaux et, comme le $\pi_0$ et
les $\pi_1$ de $\bar{J}$ leur sont naturellement isomorphes, on a bien que
$\bar{J}$ a son $\pi_0$ et ses $\pi_1$ triviaux, ce qui en fait un groupoïde
de Segal simplement connexe.\\
CQFD.\\

Nous pouvons donc bien appliquer l'\eq~des notions d'\eqs\\ à l'intervalle
$\bar{J}$ et montrer au niveau de sa réalisation géométrique qu'il est
bien en \eq~avec sa sous-catégorie de Segal pleine d'objet $0$.

\begin{lem}
Il existe un morphisme $p:\bar{J}\rightarrow \bar{J}$ envoyant $0$ et $1$ sur
$0$ qui a la propriété d'être une \eq~de catégories de Segal.
\end{lem}
{\it Preuve :}\\
Tout d'abord définissons $p$. Pour ce faire, définissons un morphisme
$p^{pre}$ de
$\bar{J}^{pre}$ dans $\bar{J}$. Comme $\bar{J}$ est une catégorie facile marquée,
car c'est une catégorisation d'une précatégorie de Segal, nous ne parlerons
que des éléments marqués quand on parlera de "composées". L'idée est de dire que $0$ est équivalent
à lui-même dans le groupoïde de Segal $\bar{J}$ par l'intermédiaire de
$w$ et de son inverse $w$, où on note $w$ l'élément $\delta_{0,2}(T_1)$.  La
composée $w^2$ est homotope de deux façons à $Id_0$. 
On va considérer deux homotopies distinctes entre $w^2$ et $Id_0$. La première
est simplement $\alpha *\alpha$. La seconde est la composée des homotopies
suivantes : on part de l'homotopie entre $w^2$ et $v\circ u \circ v \circ u$
puis on fait celle entre $v\circ u \circ v \circ u$ et $v\circ x\circ u$, où
$x$ dénote $\delta_{0,2}(T_2)$, puis on fait $Id_v *\beta * Id_u$ qui aboutit
à $v\circ Id_1\circ u$, puis on fait celle entre $v\circ Id_1\circ u$ et $w$
et on termine par $\alpha$. Pour simplifier nous noterons $\alpha\circ \beta$
cette homotopie.
 Ainsi le morphisme
$p^{pre}:\bar{J}^{pre}\rightarrow \bar{J}$ est défini ainsi : $0$ et $1$
s'envoient sur $0$, $u$ et $v$ sur $w$, $T_1$ et $T_2$ sur le triangle
$(w,w,w^2)$, $\alpha$ sur $\alpha*\alpha$ et $
\beta$ sur $\alpha \circ \beta$. Ceci définit bien un morphisme
$p^{pre}$ de $\bar{J}^{pre}$ vers $\bar{J}$. Comme $\bar{J}$ est une catégorie
facile marquée, par propriété universelle de la catégorisation, ce
morphisme se factorise de manière unique à travers $can_{\bar{J}^{pre}}$.
Notons alors $p$ l'unique morphisme de $\bar{J}$ dans lui-même qui
précomposé par $can_{\bar{J}^{pre}}$ redonne $p^{pre}$.\\
\\
Montrons que $p$ ainsi défini est une \eq~de catégories de Segal. Comme,
par le lemme précédent, $p$
est un morphisme entre groupoïdes de Segal simplement connexes,
d'après le corollaire~\ref{eqgp2}, il faut et il suffit
de montrer que $|p|$ est une \eq~d'homotopie faible pour montrer que $p$ est une
\eq~de catégories de Segal. Comme $|p^{pre}|$ se
factorise en l'équivalence $|can_{\bar{J}^{pre}}|$ suivie de $|p|$, par
l'axiome "trois pour deux" dans la \cmf~$\mathcal{TOP}$, montrer que $|p|$ est
une \eq~d'homotopie faible équivaut à montrer que $|p^{pre}|$ en est une.\\
\\
Pour cela, on va construire un objet intermédiaire entre $\bar{J}^{pre}$ et
$\bar{J}$. Pour faciliter la lecture, au lieu de dire par quelle flèche
génératrice je fais la somme amalgamée, je dirai quelle composée ou
quelle homotopie je veux rajouter.  On part de $\bar{J}^{pre}$ et on veut lui
rajouter dans l'ordre : le morphisme composé $w^2$, l'homotopie composée
$\alpha *\alpha$, les morphismes composés $v\circ u \circ v \circ u$ et
$v\circ x\circ u$, les homotopies entre $v\circ u$ et $w$, entre $u\circ v$ et
$x$, puis entre $v\circ u \circ v\circ u$ et $w^2$, entre $v\circ u\circ v\circ
u$ et $v\circ x \circ u$, puis l'homotopie composée $Id_v *\beta * Id_u$, puis
l'homotopie composée $\alpha\circ \beta$. Notons $\bar{J'}$ le résultat de
ce sous-plan d'addition fini de $Cat(\bar{J}^{pre})$.\\ 
\\
On remarque aisément que
$p^{pre}$ se factorise à travers l'inclusion de $\bar{J'}$ dans $\bar{J}$. En
outre comme $\bar{J'}$ est construit à partir de $\bar{J}^{pre}$ par des
sommes amalgamées avec des flèches génératrices, dont on a vu que les
réalisations géométriques sont des cofibrations triviales, et que la
réalisation géométrique respecte les colimites, il vient que l'inclusion
de $|\bar{J}^{pre}|$ dans $|\bar{J'}|$ est une cofibration triviale, comme suite
de sommes amalgamées le long d'un morphisme de cofibrations triviales dans la
\cmf~$\mathcal{TOP}$. Ceci montre que $|\bar{J'}|$ est faiblement homotope à
la sphère.\\ 
\\
Regardons maintenant de près l'image de la sphère $|\bar{J}^{pre}|$ dans
$|\bar{J'}|$. Les deux triangles $T_1$ et $T_2$ s'envoyent sur le triangle
$(w,w,w^2)$, le disque $\alpha$ sur le disque $\alpha *\alpha$ et le disque
$\beta$ sur le disque $\alpha\circ \beta$. Donc $|\bar{J}^{pre}|$ s'envoie sur
la sphère de grand cercle $w^2$ et de calottes $\alpha * \alpha$ et
$\alpha\circ\beta$ munie d'une collerette contractile qu'est $(w,w,w^2)$, et ceci en
envoyant calottes sur calottes. Ainsi la co-restriction de $|p^{pre}|$ à
$|\bar{J'}|$ est homotope à l'identité de la sphère.\\
\\
Comme la réalisation géométrique préserve les colimites et que
$\bar{J'}$ est une partie de $Cat(\bar{J}^{pre})$ qui n'est autre que
$\bar{J}$, l'inclusion de $|\bar{J'}|$ dans $|\bar{J}|$ est une colimite
séquentielle transfinie de sommes amalgamées le long d'un morphisme des
réalisations géométriques des  flèches génératrices. Or on a vu, dans
la démonstration du fait que $|can_A|$ est \eq~d'homotopie faible, que les
réalisations géométriques des flèches génératrices sont des
cofibrations triviales. Comme $\mathcal{TOP}$ est une \cmf, les cofibrations
triviales sont stables par colimite séquentielle transfinie et par somme
amalgamée le long d'un morphisme, d'où l'on tire que l'inclusion de
$|\bar{J'}|$ dans $|\bar{J}|$ est une cofibration triviale.
Par l'axiome
"trois pour deux" dans $\mathcal{TOP}$, la composée $|p^{pre}|$ est une
\eq~d'homotopie faible, et comme $|can_{\bar{J}^{pre}}|$ est une \eq~d'homotopie
faible, toujours par
"trois pour deux", il vient que $|p|$ est une \eq~d'homotopie faible. On conclut
que $p$ est une
\eq~de catégories de Segal entre groupoïdes de Segal simplement connexes par le
corollaire~\ref{eqgp2}.\\
CQFD. \\
 
Comme nous avions montré auparavant que l'intervalle $\bar{J}$ vérifie les
parties a) et b) de l'hypothèse R') de l'ébauche 2 du théorème central et
que $\bar{J}^{pre}$ vérifie l'hypothèse S'), ce dernier résultat finit de
montrer que la donnée de Segal proto-facile définissant les catégories de
Segal admet une bon intervalle pour les catégories de Segal.

\begin{cor}
La donnée de Segal facile sur les ensembles simpliciaux définie dans
l'exemple~\ref{hypenssimp2} vérifie les hypothèses R') et S') de l'ébauche 2
du théorème central.
\end{cor}
{\it Preuve :}\\
Montrons d'abord l'hypothèse R'). La partie a) de l'hypothèse R') est
vérifiée par définition de $\bar{J}$. Toujours par définition de
$\bar{J}$, ses deux objets $0$ et $1$ sont équivalents. Donc leurs images par
un morphisme de catégories de Segal quelconque sont encore équivalentes, ce
qui montre la partie b) de l'hypothèse R').\\ 
\\
Notons maintenant $\bar{L}$ la
sous-catégorie pleine de $\bar{J}$ d'objet unique $0$. L'inclusion de
$\bar{L}$ dans $\bar{J}$ est une équivalence de catégories
de Segal car d'une part, comme $1$ est équivalent à $0$ dans $\bar{J}$,
l'inclusion est essentiellement surjective, et d'autre part, comme $\bar{L}$
est une sous-catégorie pleine de $\bar{J}$, l'inclusion est pleinement
fidèle. Remarquons que l'\eq~$p$ du lemme précédent se factorise par
l'inclusion de $\bar{L}$ dans $\bar{J}$. Par l'axiome "trois pour deux" pour les
précatégories de Segal, qui est vrai car la donnée de Segal de
l'exemple~\ref{hypenssimp2} vérifie les hypothèses A') à Q') de l'ébauche 2
du théorème central, il vient que le morphisme $\bar{J}\rightarrow \bar{L}$ est une \eq~de
catégories de Segal, ce qui montre la partie c) de l'hypothèse R').\\
\\
Par définition de $\bar{J}^{pre}$, $\bar{J}^{pre}$ vérifie la partie a) de
l'hypothèse R'). De plus, par construction, $\bar{J}$ est la catégorisation de
$\bar{J}^{pre}$. Pour montrer que l'hypothèse S') est vérifiée, il ne reste
plus qu'à montrer que $\bar{J}^{pre}$ est $\alpha$-petit. Par construction,
$\bar{J}^{pre}$ est une colimite constituée de deux $\Delta[0]$, de deux
$\Delta[1]$, de deux $\Delta[2]$ et de deux $\Delta[1]\Theta \bar{I}$. Or tous
ces éléments sont dénombrables donc $\alpha$-petits car dans la
démonstration du lemme~\ref{seghyp1} on a pris $\alpha=2^{\aleph_0}>\aleph_0$.
Comme colimite finie de précatégories de Segal $\alpha$-petites,
$\bar{J}^{pre}$ est $\alpha$-petit, ce qui termine la vérification de
l'hypothèse S').\\
CQFD.\\

Ceci termine la vérification des dix-neuf hypothèses de l'ébauche 2 du
théorème central par la donnée de Segal sur les ensembles simpliciaux qui
engendre la notion de catégorie de Segal. Ainsi le théorème central peut
s'appliquer à cette donnée de Segal, ce qui montre que la catégorie des
précatégories de Segal avec leurs \eqs~faibles est une \cmf.

\begin{theo}
La catégorie des précatégories de Segal au sens \cite{h-s} admet une
structure de \cmf~avec pour cofibrations les monomorphismes et pour \eqs~faibles
les morphismes de précatégories de Segal dont la catégorisation est une
\eq~de catégories de Segal au sens \cite{h-s}.
\end{theo}
{\it Preuve :}\\
D'après le lemme~\ref{seghyp1} et le corollaire précédent, la donnée de
Segal proto-facile sur les ensembles simpliciaux engendrant la notion de catégorie de Segal
(celle de l'exemple~\ref{hypenssimp2})
vérifie les hypothèses A') à S') de l'ébauche 2 du théorème central. On peut donc lui
appliquer le théorème central, ce qui nous donne le résultat.\\
CQFD.\\

Ce résultat non seulement permet de démontrer l'existence de la \cmf~pour les
catégories de Segal décrite dans \cite{h-s} mais sert également
d'application d'une part de la théorie des catégories enrichies faibles développée
dans cette thèse et d'autre part au théorème central munissant ces
catégories enrichies faibles d'une \cmf~traduisant leurs théories
homotopiques.\\

En bilan de tout ce long et fastidieux exposé, nous pouvons dire qu'il reste
encore du travail à faire dans ce domaine. Tout d'abord, il serait
intéressant de montrer que la \cmf~obtenue avec le théorème central en
vérifie aussi les hypothèses. En effet, dans ce cas, nous pourrions
directement munir les $n$-catégories enrichies faibles d'une structure de
\cmf. Ensuite il serait bon de montrer d'autres exemples d'applications du
théorème central et (mais peut-être cela revient-il au même) de
redémontrer le théorème central avec des hypothèses moins contraignantes.
Si en effet nous pouvions nous débarrasser de l'hypothèse demandant aux
cofibrations d'être les monomorphismes, nous pourrions alors appliquer le
théorème au cas des $n$-catégories de Tamsamani, ce qui redémontrerait
l'existence de la \cmf~sur ces catégories faibles démontrée dans \cite{s}.
Enfin il serait également intéressant de voir si l'on ne peut pas construire
une théorie encore plus générale de catégories faibles enrichies en
regardant ce qui se passe si l'on essaie de prendre une autre catégorie que la
catégorie simpliciale $\Delta$ comme catégorie de base pour la définition des catégories
enrichies faibles.\\

Je conclurai ce travail tout d'abord en remerciant cordialement les
courageux lecteurs qui auraient eu la témérité de lire tout ou partie de
ma thèse. Mais également je consacrerai mon dernier chapitre à analyser les apports de mes quatre années de thèse en termes de compétences professionnelles. 

\newpage
....

\chapter{Projet de thèse et compétences}

\newpage

Quatre ans de ma vie ont été consacrés presque entièrement à un seul projet : une thèse de Mathématiques Pures intitulée ``Catégories enrichies faibles.'' Je vais tenter d'expliquer en quelques mots le sujet de cette thèse qui a mobilisé le plus gros de mon énergie durant une si longue période.\\

Deux théories principales fondent actuellement les Mathématiques. Si la première nommée théorie des ensembles consiste en l'étude des objets mathé\-mati\-ques pris individuellement, la seconde plus récente, dite théorie des caté\-gories, vise à modéliser les structures mathématiques. En effet la catégorie correspondant à une structure donnée est constituée par les objets mathématiques munis de cette structure et par les liens entre ces objets respectant la structure. En particulier, la structure de catégorie donne lieu à la catégorie des catégories dans laquelle les liens entre catégories forment eux aussi une catégorie. C'est un premier exemple de catégorie enrichie (ici sur les catégories). Plus généralement, une catégorie enrichie est une catégorie dont les liens forment une structure. L'adjectif faible rajouté au composé `` catégorie enrichie'' signifie seulement que les liens de la catégorie considérée ne sont pas toujours composables. Le sujet de la thèse est donc d'une part de mettre en place une théorie des catégories enrichies faibles qui permette de définir une structure commune aux différents exemples rencontrés dans ce domaine et d'autre part de munir ces catégories enrichies faibles d'une structure permettant de faire des calculs dans cette théorie.\\

Si ce sujet a bien un intérêt au niveau de la recherche mathématique comme je vais le montrer plus loin, il est cependant très intéressant de se demander ce que moi j'ai bien pu retirer de cette expérience relativement longue. En particulier, une question se pose : en quoi ces quatre années consacrées à une thèse de Mathématiques Pures ont-elles pu m'être profitables au niveau de l'insertion professionnelle. Afin de faire un retour sur mes quatre ans de thèse et d'étudier plus en détail ces points importants, je vais tout d'abord exposer les motivations qui ont conduit à cette thèse. Puis je vais examiner les principales étapes du déroulement de la thèse et je terminerai par une mise au point sur les compétences acquises durant ces quatre années aussi bien dans le travail de la thèse même que dans les autres travaux que j'ai eu à faire durant cette période.

\newpage

\section{En amont de la thèse ...}

Pour commencer, je vais analyser comment je me suis retrouvé à traiter un tel sujet. Pour tenter de dénouer les différentes causes qui ont abouti à la situation actuelle, je vais tout d'abord expliquer mes motivations personnelles concernant la thèse. Puis j'aborderai les motivations de l'équipe qui m'a accueilli.\\

Pour comprendre ce qui m'a poussé à faire une thèse de Mathématiques Pures concernant les catégories, il me faut remonter à la maîtrise. Ce fut en effet l'année dans mon cursus scolaire où j'entendis pour la première fois parler de la théorie des catégories. Mon goût pour la Logique et les fondements des Mathématiques, qui jusqu'alors furent mes principales motivations pour poursuivre des études supérieures en Mathématiques Pures, me firent tout de suite accrocher à cette théorie toute nouvelle pour moi et possédant une puissance à mes yeux infiniment plus grande que la théorie des ensembles que j'avais étudiée jusqu'alors. La matière nommée topologie algébrique, dans laquelle j'avais rencontré les catégories pour la première fois et qui était approfondie en DEA, me fit préférer la poursuite des études de Mathématiques à la préparation des concours de l'enseignement public qui constituait pour moi l'entrée dans la vie active.\\

C'est véritablement en DEA que le domaine de ma thèse allait se préciser. En effet en DEA le troisième trimestre correspond à la rédaction d'un mémoire. Ma nouvelle marotte qu'était pour lors la théorie des catégories me poussa à chercher un directeur de mémoire qui me permettrait de travailler dans ce domaine. Après quelques difficultés, je finis par trouver André Hirschowitz qui me proposa de travailler sur les articles de Carlos Simpson, le créateur d'une théorie de catégories faibles d'ordre supérieur. Le sujet du mémoire fut donc de lire et de rédiger une petite partie d'une pré-publication de Simpson. Mes capacités de compréhension d'articles mathématiques liées à ma capacité de vulgarisation des théories mathématiques amenèrent André Hirschowitz à me proposer de faire ma thèse sur les pré-publications de Simpson. Outre l'intérêt de cette proposition, ce fut la perspective d'obtenir le plus haut diplôme du cursus universitaire et le titre de docteur ès sciences qui fut déterminante dans mon choix de continuer mes études par une thèse de Mathématiques Pures plutôt que de me lancer dans la préparation des concours d'entrée à l'éducation nationale, même si cette dernière option constituait alors mon projet professionnel.\\

Considérons maintenant les motivations d'André Hirschowitz qui l'ont poussé à me proposer de faire une thèse sur les théories de Simpson. Tout d'abord il faut préciser qu'André Hirschowitz travaillait souvent en collaboration avec Carlos Simpson qui était encore à l'époque à Toulouse mais qui rejoignit l'équipe d'André Hirschowitz à Nice l'année d'après. L'année suivante, un de ses anciens élèves vint lui aussi à Nice, faisant ainsi de l'équipe niçoise une équipe à la pointe de la recherche sur les catégories faibles d'ordre supérieur, les champs supérieurs et leurs applications. Il est à noter à ce propos que l'équipe niçoise a organisé la première conférence internationale sur les catégories faibles d'ordre supérieur avec leurs applications à la géométrie algébrique.\\

 Toutefois il existe aux USA, au Canada, au Royaume-Uni ou en Australie d'autres théories sur le même thème qui sont en intense développement avec des techniques et des résultats plus ou moins avancés selon les équipes. En outre, si l'équipe niçoise est l'une des plus avancées dans ce domaine, tous ses résultats et ses théories sont basés sur un théorème fondamental de Carlos Simpson non publié dans une revue mathématique et dont la preuve contient des passages qui n'ont pas fait l'objet d'une clarification. Aussi l'enjeu de ma thèse est de mettre au point la preuve de ce théorème et d'exploiter les idées de cette preuve pour formuler une théorie plus générale permettant une meilleure compréhension de ces théories de catégories faibles d'ordre supérieur. L'impact principal qui en est attendu par l'équipe niçoise est d'assurer par la publication de ma thèse un statut stable au théorème de Simpson et ainsi de permettre d'utiliser ses méthodes et d'avoir un accès plus compréhensible au travail de pionnier de Simpson qui, s'il reste non publié, se serait de lui-même condamné à disparaître dans l'oubli.\\

Il est donc clair que la conjonction simultanée de mon intérêt à faire une thèse dans le domaine des catégories avec celui d'André Hirschowitz de publier des travaux reprenant et exploitant les idées originales mais non publiées de la théorie des catégories faibles d'ordre supérieur de Simpson a été le point de départ d'un travail de quatre ans dont il est intéressant de connaître les différentes étapes afin de mieux comprendre ce que j'en ai retiré.

\newpage

\section{Au cours de la thèse ...}

Je vais maintenant analyser le déroulement de la thèse. Pour débuter cette analyse, je commencerai par évaluer les risques et les avantages que comportait la thèse en elle-même. Puis je m'attarderai sur l'encadrement et le financement de la thèse avant de terminer par un exposé des principales étapes de la thèse.\\

Un des aspects importants et récurrents de ces quatre années passées à faire une thèse a été la gestion des risques inhérents à ce travail. Aussi est-il intéressant de se demander en quoi la thèse était-elle risquée et quels avantages comportait-elle qui purent m'éviter d'abandonner malgré de fréquentes crises de découragement. Bien évidemment déjà en DEA, mes camarades et moi-même avions été avertis de la difficulté intrinsèque aux thèses de Mathématiques Pures qui réside en une très grande abstraction. En revanche, ma thèse concernant essentiellement la démonstration généralisée d'un théorème à la preuve non formatée, le risque principal qui en découle est de pouvoir ne rien trouver. En effet contrairement à beaucoup d'autres thèses mathématiques, la mienne n'est pas une compilation, d'au plus une centaine de pages, de divers résultats plus au moins liés entre eux mais au contraire un immense édifice de trois cents pages dont le point culminant est la démonstration du théorème de Simpson. Si la preuve est fausse, incomplète ou si l'on n'arrive pas à prouver un des résultats intermédiaires, tout s'écroule. Le fait que cette thèse n'aboutisse qu'à un tout ou rien a été l'un des aspects de cette thèse des plus démoralisants. En outre l'un des corollaires de cette spécificité de ma thèse est que la thèse ne peut donner lieu à aucune autre publication qu'elle-même une fois achevée, ce qui est un handicap dans la recherche d'un poste de maître de conférence.\\

Heureusement la thèse comportait aussi des avantages qui m'ont permis de surmonter les graves moments de dépression. Tout d'abord, en elle-même la thèse concerne un domaine de la recherche mathématique récent en pleine ébullition et très porteur d'avenir mathématique. Ce domaine, comme on l'a vu précédemment, correspondait à un domaine qui m'intéressait. Même s'il a fallu faire la thèse pour prendre la véritable mesure de ce domaine, la théorie des catégories s'est avérée  aussi passionnante que je l'espérais. En outre, il est important de remarquer que pendant ces quatre années de thèse, j'ai eu l'opportunité d'enseigner en DEUG première année. En effet, premier de ma promotion à l'issue du DEA, j'ai obtenu le poste de moniteur, c'est-à-dire d'apprenti enseignant-chercheur ayant une charge d'enseignement correspondant à un tiers de la charge d'enseignement d'un enseignant-chercheur titulaire. De plus après mes trois années de monitorat, j'ai obtenu un poste d'ATER d'un an, c'est-à-dire un poste d'enseignant-chercheur temporaire avec une demi-charge d'enseignement. Or ces enseignements m'ont véritablement permis d'acquérir de l'assurance tout au long de ces quatre années. Face à la déstabilisation morale et psychique due à la thèse, le fait de pouvoir maîtriser un savoir et de l'enseigner à des étudiants à peine plus jeunes que moi-même a été très réconfortant, très intéressant et très enrichissant. En effet l'enseignement m'a toujours attiré. De plus j'ai eu la chance de pouvoir enseigner dans des équipes pédagogiques soudées qui avaient une véritable réflexion pédagogique et me permettaient souvent d'appliquer mes propres méthodes. Enfin à travers l'enseignement du raisonnement scientifique, j'ai pu enseigner mes domaines de prédilection, à savoir la Logique et les fondements des Mathématiques, ce qui me permettait d'ailleurs de faire un retour sur mon propre travail de thèse dont le domaine touche par certains endroits à ces domaines. Dans le même ordre d'idées, lors de groupes de travail au laboratoire de Mathématiques, on me demandait souvent d'exposer des sujets ayant trait à la théorie des catégories, ce qui me permettait d'approfondir mes connaissances dans ce domaine et m'a souvent aidé à certains endroits de la thèse tout en me redonnant confiance en moi-même.\\

Si les avantages apportés par la thèse ont pu contrebalancer l'incertitude liée aux risques inhérents à mon sujet de thèse, l'encadrement a lui aussi joué un grand rôle tout au long de la thèse. Comme je l'ai dit précédemment, plus j'avançais dans ma thèse et plus l'équipe d'André Hirschowitz a accueilli de spécialistes de la théorie de Simpson. Aussi puis-je dire que j'ai été très bien entouré : Simpson lui-même a été très souvent disponible, son ancien élève Toen n'a pas hésité à venir nous aider lorsque nous bloquions à certains endroits de la preuve du théorème de Simpson. Pour en revenir au déroulement de la thèse, c'est André Hirschowitz qui, m'ayant proposé une thèse sur la théorie de Simpson, a dirigé effectivement ma thèse même s'il n'est pas mon directeur officiel pour cette thèse. André Hirschowitz a toujours été très disponible, en particulier lorsque la thèse posait de graves problèmes même si, je dois l'avouer, j'ai toujours préféré trouver les solutions aux problèmes par moi-même plutôt que de lui demander de l'aide, sauf bien entendu dans les rares cas de blocage total sur certaines questions. Ceci s'est traduit par des réunions avec André Hirschowitz assez irrégulières et assez espacées dans le temps, phénomène amplifié par l'éloignement géographique de nos deux bureaux. Quant à Jean-Michel Lemaire qui m'a fait l'honneur de diriger officiellement ma thèse, je lui dois également beaucoup car, même s'il a suivi d'assez loin ma thèse dont le domaine ne correspond pas à son domaine de recherche, il s'y est toutefois beaucoup intéressé, m'a soutenu lors de moments difficiles et a apporté sur ma thèse un point de vue extérieur au sujet qui m'a été très enrichissant.\\

Au vu de cet encadrement, il est naturel de se poser la question du coût de la thèse. Tout d'abord, considérons mon propre salaire. Pendant mes trois premières années, étant moniteur, j'ai reçu un double salaire : 6.100 F net par mois pour la thèse et 1.800 F net par mois pour les enseignements. Lors de ma quatrième année, j'ai été payé en tant qu'A.T.E.R. (Attaché Temporaire à l'Enseignement et à la Recherche) 7.400 F net par mois, ce salaire comprenant à la fois thèse et enseignement. Ajoutons à cela mes déplacements à des colloques et à des écoles d'été qui se montent à 8.500 F. Au total (avec les enseignements) mes quatre années sont revenues à environ 375.000 F net. André Hirschowitz m'a consacré un vingtième de son temps durant ces quatre années, ce qui se chiffre à environ 92.000 F net. Comme le temps passé par mon directeur de thèse, Simpson et les autres chercheurs de Nice m'ayant aidé est négligeable par rapport au temps passé avec André Hirschowitz et par ailleurs souvent inévaluable, on peut dire qu'au total la thèse aura coûté environ un demi-million de francs français ($\cong 470.000 F$). Il est à noter que ce coût ne prend pas en compte la fraction du budget d'infrastructure, de la bibliothèque et du réseau informatique du laboratoire, d'une part car cette fraction est difficilement évaluable et d'autre part car le budget du laboratoire ne m'a pas été communiqué.\\

Enfin je vais terminer cette partie par une brève description des principales étapes de la thèse. La première phase d'une durée d'environ un an fut surtout une phase d'apprentissage et d'approfondissement de la théorie de Simpson avec des essais de définition de nouvelles catégories faibles d'ordre supérieur. A la fin de cette période, l'idée de paramétrer les définitions déjà existantes de catégories faibles d'ordre supérieur m'obligea à revenir au théorème central de Simpson et à le comprendre en détail. Lors de cette lecture approfondie de ce théorème, je découvris une erreur dans la preuve, erreur qui fut très aisément rectifiée par Carlos Simpson et qui amena André Hirschowitz à me proposer de reprendre la preuve du théorème afin d'en exploiter les idées dans un cadre plus général. Ceci ouvre la deuxième phase consistant en la relecture et la compréhension du théorème de Simpson qui dura environ un an. La preuve comprise dans ses grandes lignes, la troisième phase d'aussi un an fut celle de la tentative de généralisation du théorème, ce qui nous amena aux principaux problèmes rencontrés pendant cette thèse et aussi à la période qui faillit voir mon abandon. Après beaucoup d'énergie et de ténacité mais aussi grâce à l'important soutien d'André Hirschowitz, cette période se termina avec l'éclaircissement des derniers points obscurs dans la preuve de Simpson. Commença alors la quatrième phase de la thèse qui fut celle de la rédaction d'un théorème dont la preuve exploite les idées de la preuve de Simpson mais dans un cadre plus général et qui donna lieu à la création de la théorie des catégories enrichies faibles. Cette période d'un an fut aussi celle de la rédaction de la thèse. Au final, l'un des principaux problèmes rencontrés dans la preuve de Simpson laisse l'énoncé du théorème général de ma thèse incomplet dans le sens où la manière dont on a résolu le problème n'a pas été suffisamment bonne pour que mon théorème s'applique dans tous les cas que l'on voulait regarder. Ceci laisse encore du travail pour améliorer la thèse et c'est l'une des motivations qui m'ont poussé à demander une seconde année d'ATER.\\

Après avoir expliqué en quelques lignes les principales caractéristiques du déroulement de mes quatre années de thèse, je vais essayer de déterminer en quoi ces années m'ont été bénéfiques, en particulier pour mon avenir professionnel.

\newpage

\section{En aval de la thèse ...}

Dans cette dernière partie, je vais analyser les retombées de ces quatre années de thèse essentiellement en termes de compétences acquises durant la thèse. J'essayerai au fur et à mesure de jauger l'impact de ces acquis sur mon avenir professionnel. Je commencerai par développer les compétences mises en oeuvre pour la thèse, puis celle pour l'enseignement, enfin j'exposerai quelques conséquences de mon initiative personnelle durant ces quatre années.\\

Bien évidemment, les premières retombées de la thèse auxquelles on pense naturellement s'expriment en terme de connaissances mathématiques : théories des ensembles, des catégories, des catégories faibles d'ordre supérieur, des faisceaux, des champs, des topoi, de l'homotopie et des catégories de modèles fermées. Ces savoirs sont en outre interdisciplinaires car certains ont trait aux fondements des Mathématiques et par là à la Logique, à l'Épistémologie et donc à la Philosophie. Mais ils sont aussi transversaux car ces mêmes savoirs ont à la fois une importance pour la recherche mathématique mais aussi pour l'enseignement de la méthodologie en première année de DEUG et permettent un retour didactique sur l'enseignement des Mathématiques en DEUG.\\

Comme cela est clair dans les nouvelles maquettes de DEUG, le savoir mathématique n'est pas seulement confiné dans un amas de connaissances mais il comporte également toute une méthodologie à la fois interdisciplinaire, comme la lecture suivie et l'analyse approfondie des articles, l'exposition d'une théorie déjà connue ou d'un travail de recherche, mais aussi une méthodologie caractéristique à la matière comme la rigueur mathématique, la démonstration d'un résultat mathématique et la mise en place d'une théorie globale généralisant des théories particulières existantes. Mais la principale compétence méthodologique que ma thèse nécessitait est de savoir faire une construction d'un degré de complexité énorme pour aboutir à un unique résultat (300 pages pour démontrer un seul théorème !), ce qui pose des problè\-mes de cohérence monstrueux et demande de bonnes prises de recul, entre autres. Cette dernière compétence est un atout assez rare qui me sera sûrement utile dans ma vie professionnelle ainsi d'ailleurs que les méthodes de base interdisciplinaires, et ce quel que soit mon devenir professionnel. Je dois rajouter à cela la faculté que j'ai acquise durant ces quatre années à gérer mais aussi à surmonter les déprimes profondes et courantes intrinsèques au travail d'un chercheur en Mathématiques Pures, déprimes que j'ai analysées dans la partie précédente. Cette ténacité face à de fréquents moments de démoralisation est, je le crois, une qualité assez recherchée dans la vie active.\\

Outre le travail de thèse, les enseignements menés durant ces quatre années ont eux aussi exigé des compétences particulières. Comme je l'ai déjà évoqué, j'ai eu l'occasion de participer à divers enseignements dans lesquels j'ai eu une part active au niveau de la pédagogie à adopter. J'ai également participé à la conduite de nouvelles méthodes pédagogiques qui ont nécessité une grande adaptabilité aux situations nouvelles ou inattendues, une capacité d'intégration à une équipe pédagogique ainsi qu'un apprentissage du travail en équipe, une remise en question permanente due au caractère interactif de ces méthodes et enfin un énorme investissement. La mise en place de ces méthodes m'a également amené à conduire une réflexion didactique et épistémologique et à échanger avec des philosophes et des didacticiens. La réussite de ces méthodes m'a aussi mené à devoir exposer lors de conférences sur l'enseignement et à écrire un article sur mon enseignement de la méthodologie en première année de DEUG.\\

 Si tout ceci m'a donné une expérience à la fois unique et riche au niveau de l'enseignement, cette expérience liée d'une part à mes échecs répétés aux concours de l'enseignement secondaire et d'autre part à ce que vivent mes proches dans l'éducation nationale m'a plutôt convaincu de ne pas enseigner en-dessous du DEUG. Aussi mon projet professionnel initial (depuis la licence) d'intégrer l'éducation nationale s'est assez vite vu transmué en un nouveau projet : celui d'intégrer la fonction publique. Bien évidemment, la question qu'il est alors naturel de poser est pourquoi ne pas envisager le concours de maître de conférence. La réponse est aisée : c'est tout simplement parce que ce concours est un Nirvâna auquel tendent toutes les bonnes âmes alors que seules celles qui ont reçu la grâce pourront, non sans peine, y accéder (pour leur bonheur ?).\\

A l'image de mes enseignements, ces quatre années ont montré que je suis capable de m'impliquer énormément mais aussi de faire preuve d'initiative : outre mon implication dans la thèse et les innovations mathématiques que j'ai apportées grâce à une bonne interaction avec l'équipe d'André Hirschowitz, en plus de mes initiatives pédagogiques dont j'ai déjà parlé, j'ai créé et animé divers groupes de travail soit avec des thésards mathématiciens sur les champs, les axiomatiques des ensembles, les cardinaux transfinis et les catégories enrichies faibles, soit avec des philosophes et des logiciens sur la théorie des catégories et des topoi et sur les logiques déviantes (intuitionnistes, modales, para-consistantes ...) Je pense que ma capacité d'intégration dans diverses équipes en même temps ainsi que ma capacité à organiser des groupes de travail intra et inter-équipe seront des atouts importants pour ma future vie professionnelle.\\

Une autre grande preuve de mes capacités d'investissement, d'initiative mais aussi d'adaptabilité à de nouvelles situations est ma participation à la coordination de la première année de DEUG MASS ainsi que la responsabilité du module de méthodologie (UMED2). Ceci m'a obligé à apprendre à créer un emploi du temps, à gérer des salles, des horaires de modules et des horaires d'enseignants (avec leurs contraintes personnelles...) mais aussi à gérer les effectifs d'étudiants répartis en plusieurs groupes. En tant que responsable de module, j'ai dû aussi m'occuper de la mise en place d'un programme compatible avec les impératifs de la maquette, créer une équipe pédagogique, organiser les examens, gérer les notes et participer aux délibérations. Je pense que tout cet investissement purement bénévole de ma part (car non payé ni demandé par les statuts de moniteur ou d'ATER) est très important pour mon avenir professionnel que je vois en ce moment plutôt dans l'administration publique.\\

Il est donc aisé de constater au vu de ce qui précède que ces quatre années ont été très riches et très remplies et m'ont permis d'acquérir bon nombre de compétences dont certaines seront, je pense, très utiles pour mon avenir professionnel. Que peut-on donc tirer comme conclusion générale concernant cette thèse ?\\

\newpage

\section{En guise de conclusion ...}

Pour résumer ce qui précède, je peux dire que la thèse m'a apporté une bonne culture mathématique mais aussi un avant-goût du travail administratif, du travail en équipe et des exposés à des colloques. J'ai pu aussi me découvrir des qualités jusqu'alors insoupçonnées telles l'initiative, l'intégration et l'adaptabilité qui se sont exprimées surtout dans les enseignements et les travaux d'ordre administratif, ce qui est assez encourageant pour mon avenir professionnel. \\

En conclusion, je dirai que la thèse m'a apporté un certain savoir-faire qui sera très certainement utile pour ce à quoi je me destine, à savoir la fonction publique, mais aussi un savoir et un diplôme que j'espère pouvoir valoriser dans ma future vie professionnelle. Il n'est peut-être pas inutile de rappeler que je n'ai pas du tout fait ma thèse en vue d'un quelconque projet professionnel, comme je l'ai expliqué dans la première partie. Néanmoins, j'ai pu retirer de ces quatre années un maximum d'expériences et d'enseignements personnels et surtout le sujet de ma thèse, sur lequel je n'avais pas d'idées réelles avant de ne l'aborder en thèse, s'est trouvé à la hauteur de mes espérances. Aussi, malgré de gros moments de ``déprime'', je peux réellement dire que je ne regrette absolument pas d'avoir fait cette thèse qui m'a permis véritablement de vivre une expérience somme toute courte mais si dense et, de toute manière, unique dans ma vie.

\newpage
 ...

\printindex 

\newpage
...

\newpage

\begin{center}
{\bf Résumé :}  
\end{center}

Cette thèse est consacrée à la démonstration d'un théorème montrant l'existence d'une structure de \cmf~concernant les catégories faiblement enrichies. Il faut au préalable définir les notions de catégories faiblement enrichies et d'\eq~de catégories faiblement enrichies de telle manière que ces notions recouvrent diverses notions déjà existantes de catégories faibles d'ordre supérieur telles les catégories de Segal, les n-catégories de Tamsamani et les n-catégories strictes. Afin de démontrer notre théorème, nous devons mettre au point une théorie de plans d'addition de cellules sur le modèle de l'argument du petit objet à la Quillen. Nous terminons ce travail en montrant que notre théorème recouvre le cas des catégories de Segal. Ce dernier résultat nécessite de montrer une adjonction ``groupoïde fondamental-réalisation géométrique'' entre les groupoïdes de Segal et les espaces topologiques.

\begin{center}
{\bf English title:}\\ {\bf Weak enriched categories.} 
\end{center}

\begin{center}
{\bf Abstract :}  
\end{center}

This thesis is devoted to the proof of a theorem showing the existence of a closed model category structure for weakly enriched categories. It requires first of all the definitions of weakly enriched categories and equivalences of weakly enriched categories such that these definitions recover some existing notions of higher order weak categories, for example Segal categories, Tamsamani n-categories and strict n-categories. In order to prove our theorem, we elaborate a theory of plans for cell addition following the approach of the small object argument {\it à la} Quillen. We conclude this work with the proof that our theorem recovers the case of Segal categories. This last result requires a fundamental groupoid-geometric realization adjunction between Segal groupoids and topological spaces.\\

{\bf Discipline :} Mathématiques.\\

{\bf Mots-clefs :} \cmf, catégorie de Segal, catégorie enrichie, groupoïde fondamental, n-catégorie, réalisation géométrique.\\
\\
Laboratoire de Mathématiques J.A. Dieudonné U.M.R. 6621\\
Université de Nice Sophia-Antipolis Parc Valrose 06108 Nice Cedex 2


\begin{thebibliography}{99}
\bibitem{h} P.Hirschhorn. \emph{Localization of Model Categories}, preprint.
\bibitem{h-s} A.Hirschowitz, C.Simpson. \emph{Descente pour les $n$-champs},
preprint.
\bibitem{hv} M.Hovey. \emph{Model categories}, Mathematical Surveys and Monomographs, no. 63, American Mathematical Society, 1999.
\bibitem{q} D.G.Quillen. \emph{Homotopical algebra}, Lectures Notes in Mathematics, vol. 43, Springer-Verlag, Berlin, 1967.
\bibitem{s} C.Simpson. \emph{A closed model structure for
n-categories, internal Hom, n-stacks and generalized Seifert-Van Kampen}, preprint.
\bibitem{t} Z. Tamsamani. \emph{Sur des notions de $n$-catégorie et de
$n$-groupoïde non-stricts via des ensembles multi-simpliciaux}. K-theory, 1999. 
\end{thebibliography}
\end{document}